\documentclass[11pt]{amsproc}

\usepackage{amsmath}
\usepackage{amssymb}
\usepackage{epsfig}
\usepackage{psfrag}

\makeindex

\newtheorem{theorem}{\bf Theorem}[section]
\newtheorem{lemma}[theorem]{\bf Lemma}
\newtheorem{proposition}[theorem]{\bf Proposition}
\newtheorem{corollary}[theorem]{\bf Corollary}

\theoremstyle{definition}
\newtheorem{definition}[theorem]{\bf Definition}
\newtheorem{remark}[theorem]{\bf Remark}
\newtheorem{claim}{\bf Claim}
\newtheorem{assertion}{\bf Assertion}

\newenvironment{pf}{{\medskip \bf\em Proof.}}{{\hfill$\square$}\medskip}

\numberwithin{section}{part} \numberwithin{equation}{section}

\def\Ric{{\rm Ric}\,}
\def\Vol{{\rm Vol}\,}
\def\inj{{\rm inj}\,}
\def\tr{{\rm tr}\,}
\def\({\left(}
\def\){\right)}
\def\be{\begin{equation}}
\def\ee{\end{equation}}
\def\nn{\nonumber}
\def\<{\langle}
\def\>{\rangle}
\def\endproof{{\hfill $\square$}\medskip}

\begin{document}

\title[Hamilton-Perelman's Proof]{{\LARGE Hamilton-Perelman's Proof of the
Poincar\'{e} Conjecture\vspace*{3mm}\\
and the Geometrization Conjecture*}}

\author[H.-D. Cao]{Huai-Dong Cao}   
\address{Department of Mathematics\\
Lehigh University\\
Bethlehem, PA 18015} \email{huc2@lehigh.edu}

\author[X.-P. Zhu]{Xi-Ping Zhu}
\address{Department of Mathematics\\
Zhongshan University\\
Guangzhou 510275, \hspace*{5mm} P.R. China}
\email{stszxp@zsu.edu.cn}

\thanks{*This is a revised version of the article by the same authors
that originally appeared in Asian J. Math., {\bf 10}({\bf 2})
(2006), 165--492.}

\begin{abstract}

In this paper, we provide an essentially self-contained and
detailed account of the fundamental works of Hamilton and the
recent breakthrough of Perelman on the Ricci flow and their
application to the geometrization of three-manifolds. In
particular, we give a detailed exposition of a complete proof of
the Poincar\'e conjecture due to Hamilton and Perelman.

\end{abstract}

\maketitle

\newpage
\tableofcontents


\newpage
\part*{{\Large Introduction}}

This is a revised version of the article that originally appeared
in Asian J. Math., {\bf 10}({\bf 2}) (2006), 165--492. In July, we
received complaint from Kleiner and Lott for the lack of
attribution of their work caused by our oversight in the prior
version; we have apologized and acknowledged their contribution in
an erratum to appear in December 2006 issue of the Asian Journal
of Mathematics. In this revision, we have also tried to amend
other possible oversights by updating the references and
attributions. In the meantime, we have changed the title and
modified the abstract in order to better reflect our view that the
full credit of proving the Poincare's conjecture goes to Hamilton
and Perelman. We regret all the oversights occurred in the prior
version and hope that this revision will right these inattentions.
Other than these modifications mentioned, the paper remains
unchanged.

In this paper, we shall present the Hamilton-Perelman theory of
Ricci flow. Based on it, we shall give a detailed account of
complete proofs of the Poincar\'e conjecture and the
geometrization conjecture of Thurston. While the results presented
here are based on the accumulated works of many geometric
analysts, the complete proof of the Poincar\'e conjecture is due
to Hamilton and Perelman.

An important problem in differential geometry is to find a
canonical metric on a given manifold. In turn, the existence of a
canonical metric often has profound topological implications. A
good example is the classical uniformization theorem in two
dimensions which, on one hand, provides a complete topological
classification for compact surfaces, and on the other hand shows
that every compact surface has a canonical geometric structure: a
metric of constant curvature.

How to formulate and generalize this two-dimensional result to
three and higher dimensional manifolds has been one of the most
important and challenging topics in modern mathematics. In 1977,
W. Thurston \cite{Th82}, based on ideas about Riemann surfaces,
Haken's work and Mostow's rigidity theorem, etc, formulated a
geometrization conjecture for three-manifolds which, roughly
speaking, states that every compact orientable three-manifold has
a canonical decomposition into pieces, each of which admits a
canonical geometric structure. In particular, Thurston's
conjecture contains, as a special case, the Poincar\'e conjecture:
A closed three-manifold with trivial fundamental group is
necessarily homeomorphic to the 3-sphere $\mathbb{S}^3$. In the
past thirty years, many mathematicians have contributed to the
understanding of this conjecture of Thurston. While Thurston's
theory is based on beautiful combination of techniques from
geometry and topology, there has been a powerful development of
geometric analysis in the past thirty years, lead by S.-T. Yau, R.
Schoen, C. Taubes, K. Uhlenbeck, and S. Donaldson, on the
construction of canonical geometric structures based on nonlinear
PDEs (see, e.g., Yau's survey papers \cite{Y78, Y06}). Such
canonical geometric structures include K\"ahler-Einstein metrics,
constant scalar curvature metrics, and self-dual metrics, among
others. However, the most important contribution for geometric
analysis on three-manifolds is due to Hamilton.

In 1982, Hamilton \cite{Ha82} introduced the Ricci flow
$$\frac{\partial g_{ij}}{\partial t}=-2R_{ij}$$
to study compact three-manifolds with positive Ricci curvature.
The Ricci flow, which evolves a Riemannian metric by its Ricci
curvature, is a natural analogue of the heat equation for metrics.
As a consequence, the curvature tensors evolve by a system of
diffusion equations which tends to distribute the curvature
uniformly over the manifold. Hence, one expects that the initial
metric should be improved and evolve into a canonical metric,
thereby leading to a better understanding of the topology of the
underlying manifold. In the celebrated paper \cite{Ha82}, Hamilton
showed that on a compact three-manifold with an initial metric
having positive Ricci curvature, the Ricci flow converges, after
rescaling to keep constant volume, to a metric of positive
constant sectional curvature, proving the manifold is
diffeomorphic to the three-sphere $\mathbb{S}^3$ or a quotient of
the three-sphere $\mathbb{S}^3$ by a linear group of isometries.
Shortly after, Yau suggested that the Ricci flow should be the
best way to prove the structure theorem for general
three-manifolds. In the past two decades, Hamilton proved many
important and remarkable theorems for the Ricci flow, and laid the
foundation for the program to approach the Poincar\'e conjecture
and Thurston's geometrization conjecture via the Ricci flow.

The basic idea of Hamilton's program can be briefly described as
follows. For any given compact three-manifold, one endows it with
an arbitrary (but can be suitably normalized by scaling) initial
Riemannian metric on the manifold and then studies the behavior of
the solution to the Ricci flow. If the Ricci flow develops
singularities, then one tries to find out the structures of
singularities so that one can perform (geometric) surgery by
cutting off the singularities, and then continue the Ricci flow
after the surgery. If the Ricci flow develops singularities again,
one repeats the process of performing surgery and continuing the
Ricci flow. If one can prove there are only a finite number of
surgeries during any finite time interval and if the long-time
behavior of solutions of the Ricci flow with surgery is well
understood, then one would recognize the topological structure of
the initial manifold.

Thus Hamilton's program, when carried out successfully, will give
a proof of the Poincar\'e conjecture and Thurston's geometrization
conjecture. However, there were obstacles, most notably the
verification of the so called ``Little Loop Lemma" conjectured by
Hamilton \cite{Ha95F} (see also \cite{CCCY}) which is a certain
local injectivity radius estimate, and the verification of the
discreteness of surgery times. In the fall of 2002 and the spring
of 2003, Perelman \cite{P1, P2} brought in fresh new ideas to
figure out important steps to overcome the main obstacles that
remained in the program of Hamilton. (Indeed, in page 3 of
\cite{P1}, Perelman said ``the implementation of Hamilton program
would imply the geometrization conjecture for closed
three-manifolds" and ``In this paper we carry out some details of
Hamilton program".) Perelman's breakthrough on the Ricci flow
excited the entire mathematics community. His work has since been
examined to see whether the proof of the Poincar\'e conjecture and
geometrization program, based on the combination of Hamilton's
fundamental ideas and Perelman's new ideas, holds together. The
present paper grew out of such an effort.

Now we describe the three main parts of Hamilton's program in more
detail.

\vskip 0.1cm \noindent {\bf (i) Determine the structures of
singularities}

Given any compact three-manifold $M$ with an arbitrary Riemannian
metric, one evolves the metric by the Ricci flow. Then, as
Hamilton showed in \cite{Ha82}, the solution $g(t)$ to the Ricci
flow exists for a short time and is unique (also see Theorem
1.2.1). In fact, Hamilton \cite{Ha82} showed that the solution
$g(t)$ will exist on a maximal time interval $[0, T)$, where
either $T=\infty$, or $0<T<\infty$ and the curvature becomes
unbounded as $t$ tends to $T$. We call such a solution $g(t)$ a
maximal solution of the Ricci flow. If $T<\infty$ and the
curvature becomes unbounded as $t$ tends to $T$, we say the
maximal solution develops singularities as $t$ tends to $T$ and
$T$ is the singular time.

In the early 1990s, Hamilton systematically developed methods to
understand the structure of singularities. In \cite{Ha93}, based
on suggestion by Yau, he proved the fundamental Li-Yau \cite{LY}
type differential Harnack estimate (the Li-Yau-Hamilton estimate)
for the Ricci flow with nonnegative curvature operator in all
dimensions. With the help of Shi's interior derivative estimate
\cite{Sh89}, he \cite{Ha95} established a compactness theorem for
smooth solutions to the Ricci flow with uniformly bounded
curvatures and uniformly bounded injectivity radii at the marked
points. By imposing an injectivity radius condition, he rescaled
the solution to show that each singularity is asymptotic to one of
the three types of singularity models \cite{Ha95F}. In \cite
{Ha95F} he discovered (also independently by Ivey \cite{Iv}) an
amazing curvature pinching estimate for the Ricci flow on
three-manifolds. This pinching estimate implies that any
three-dimensional singularity model must have nonnegative
curvature. Thus in dimension three, one only needs to obtain a
complete classification for nonnegatively curved singularity
models.

For Type I singularities in dimension three, Hamilton \cite{Ha95F}
established an isoperimetric ratio estimate to verify the
injectivity radius condition and obtained spherical or necklike
structures for any Type I singularity model. Based on the
Li-Yau-Hamilton estimate, he showed that any Type II singularity
model with nonnegative curvature is either a steady Ricci soliton
with positive sectional curvature or the product of the so called
cigar soliton with the real line \cite{Ha93E}. (Characterization
for nonnegatively curved Type III models was obtained in
\cite{CZ00}.) Furthermore, he developed a dimension reduction
argument to understand the geometry of steady Ricci solitons
\cite{Ha95F}. In the three-dimensional case, he showed that each
steady Ricci soliton with positive curvature has some necklike
structure. Hence Hamilton had basically obtained a canonical
neighborhood structure at points where the curvature is comparable
to the maximal curvature for solutions to the three-dimensional
Ricci flow.

However two obstacles remained: (a) the verification of the
imposed injectivity radius condition in general; and (b) the
possibility of forming a singularity modelled on the product of
the cigar soliton with a real line which could not be removed by
surgery. The recent spectacular work of Perelman \cite{P1} removed
these obstacles by establishing a local injectivity radius
estimate, which is valid for the Ricci flow on compact manifolds
in all dimensions. More precisely, Perelman proved two versions of
``no local collapsing" property (Theorem 3.3.3 and Theorem 3.3.2),
one with an entropy functional he introduced in \cite{P1}, which
is monotone under the Ricci flow, and the other with a space-time
distance function obtained by path integral, analogous to what
Li-Yau did in \cite{LY}, which gives rise to a monotone
volume-type (called reduced volume by Perelman) estimate. By
combining Perelman's no local collapsing theorem I${'}$ (Theorem
3.3.3) with the injectivity radius estimate in Theorem 4.2.2, one
immediately obtains the desired injectivity radius estimate, or
the Little Loop Lemma (Theorem 4.2.4) conjectured by Hamilton.

Furthermore, Perelman \cite{P1} developed a refined rescaling
argument (by considering local limits and weak limits in
Alexandrov spaces) for singularities of the Ricci flow on
three-manifolds to obtain a uniform and global version of the
canonical neighborhood structure theorem. We would like to point
out that our proof of the singularity structure theorem (Theorem
7.1.1) is different from that of Perelman in two aspects: (1) in
Step 2 of the proof, we only prove a weaker version of Perelman's
Claim 2 in section 12.1 of \cite{P1}, namely finite distance
implies finite curvature. Our treatment, with some modifications,
follows the notes of Kleiner-Lott \cite{KL} in June 2003 on
Perelman's first paper \cite{P1}; (2) in Step 4 of the proof, we
give a new approach to extend the limit backward in time to an
ancient solution.

\vskip 0.1cm \noindent {\bf (ii) Geometric surgeries and the
discreteness of surgery times}

After obtaining the canonical neighborhoods (consisting of
spherical, necklike and caplike regions) for the singularities,
one would like to perform geometric surgery and then continue the
Ricci flow. In \cite{Ha97}, Hamilton initiated such a surgery
procedure for the Ricci flow on four-manifolds with positive
isotropic curvature and presented a concrete method for performing
the geometric surgery. His surgery procedures can be roughly
described as follows: cutting the neck-like regions, gluing back
caps, and removing the spherical regions. As will be seen in
Section 7.3 of this paper, Hamilton's geometric surgery method
also works for the Ricci flow on compact orientable
three-manifolds.

Now an important challenge is to prevent surgery times from
accumulating and make sure one performs only a finite number of
surgeries on each finite time interval. The problem is that, when
one performs the surgeries with a given accuracy at each surgery
time, it is possible that the errors may add up to a certain
amount which could cause the surgery times to accumulate. To
prevent this from happening, as time goes on, successive surgeries
must be performed with increasing accuracy. In \cite{P2}, Perelman
introduced some brilliant ideas which allow one to find ``fine"
necks, glue ``fine" caps, and use rescaling to prove that the
surgery times are discrete.

When using the rescaling argument for surgically modified
solutions of the Ricci flow, one encounters the difficulty of how
to apply Hamilton's compactness theorem (Theorem 4.1.5), which
works only for smooth solutions. The idea to overcome this
difficulty consists of two parts. The first part, due to Perelman
\cite{P2}, is to choose the cutoff radius in neck-like regions
small enough to push the surgical regions far away in space. The
second part, due to the authors and Chen-Zhu \cite{CZ05F}, is to
show that the surgically modified solutions are smooth on some
uniform (small) time intervals (on compact subsets) so that
Hamilton's compactness theorem can still be applied. To do so, we
establish three time-extension results (see Assertions 1-3 of the
Step 2 in the proof of Proposition 7.4.1). Perhaps, this second
part is more crucial. Without it, Shi's interior derivative
estimate (Theorem 1.4.2) may not applicable, and hence one cannot
be certain that Hamilton's compactness theorem holds when only
having the uniform $C^{0}$ bound on curvatures. We remark that in
our proof of this second part, as can be seen in the proof of
Proposition 7.4.1, we require a deep comprehension of the
prolongation of the gluing ``fine" caps for which we will use the
recent uniqueness theorem of Bing-Long Chen and the second author
\cite{CZ05U} for solutions of the Ricci flow on noncompact
manifolds.

Once surgeries are known to be discrete in time, one can complete
the classification, started by Schoen-Yau \cite{ScY79m, ScY79},
for compact orientable three-manifolds with positive scalar
curvature. More importantly, for simply connected three-manifolds,
if one can show that solutions to the Ricci flow with surgery
become extinct in finite time, then the Poincar\'e conjecture
would follow. Such a finite extinction time result was proposed by
Perelman \cite{P3}, and a proof also appears in Colding-Minicozzi
\cite{CM}. Thus, the combination of Theorem 7.4.3 (i) and the
finite extinction time result provides a complete proof to the
Poincar\'e conjecture.

\vskip 0.1cm \noindent {\bf (iii) The long-time behavior of
surgically modified solutions.}

To approach the structure theorem for general three-manifolds, one
still needs to analyze the long-time behavior of surgically
modified solutions to the Ricci flow. In \cite{Ha99}, Hamilton
studied the long time behavior of the Ricci flow on compact
three-manifolds for a special class of (smooth) solutions, the so
called nonsingular solutions. These are the solutions that, after
rescaling to keep constant volume, have (uniformly) bounded
curvature for all time. Hamilton \cite{Ha99} proved that any
three-dimensional nonsingular solution either collapses or
subsequently converges to a metric of constant curvature on the
compact manifold or, at large time, admits a thick-thin
decomposition where the thick part consists of a finite number of
hyperbolic pieces and the thin part collapses. Moreover, by
adapting Schoen-Yau's minimal surface arguments in \cite{ScY79}
and using a result of Meeks-Yau \cite{MY}, Hamilton showed that
the boundary of hyperbolic pieces are incompressible tori.
Consequently, when combined with the collapsing results of
Cheeger-Gromov \cite{ChG86, ChG90}, this shows that any
nonsingular solution to the Ricci flow is geometrizable in the
sense of Thurston \cite{Th82}. Even though the nonsingular
assumption seems very restrictive and there are few conditions
known so far which can guarantee a solution to be nonsingular,
nevertheless the ideas and arguments of Hamilton's work
\cite{Ha99} are extremely important.

In \cite{P2}, Perelman modified Hamilton's arguments to analyze
the long-time behavior of arbitrary smooth solutions to the Ricci
flow and solutions with surgery to the Ricci flow in dimension
three. Perelman also argued that the proof of Thurston's
geometrization conjecture could be based on a thick-thin
decomposition, but he could only show the thin part will only have
a (local) lower bound on the sectional curvature. For the thick
part, Perelman \cite{P2} established a crucial elliptic type
estimate, which allowed him to conclude that the thick part
consists of hyperbolic pieces. For the thin part, he announced in
\cite{P2} a new collapsing result which states that if a
three-manifold collapses with (local) lower bound on the sectional
curvature, then it is a graph manifold. Assuming this new
collapsing result, Perelman \cite{P2} claimed that the solutions
to the Ricci flow with surgery have the same long-time behavior as
nonsingular solutions in Hamilton's work, a conclusion which would
imply a proof of Thurston's geometrization conjecture. Although
the proof of this new collapsing result promised by Perelman in
\cite{P2} is still not available in literature, Shioya-Yamaguchi
\cite{ShY} has published a proof of the collapsing result in the
special case when the manifold is closed. In the last section of
this paper (see Theorem 7.7.1), we will provide a proof of
Thurston's geometrization conjecture by only using
Shioya-Yamaguchi's collapsing result. In particular, this gives
another proof of the Poincar\'e conjecture.

We would like to point out that Perelman \cite{P2} did not quite
give an explicit statement of the thick-thin decomposition for
surgical solutions. When we were trying to write down an explicit
statement, we needed to add a restriction on the relation between
the accuracy parameter $\varepsilon$ and the collapsing parameter
$w$. Nevertheless, we are still able to obtain a weaker version of
the thick-thin decomposition (Theorem 7.6.3) that is sufficient to
deduce the geometrization result.

In this paper, we shall give detailed exposition of what we
outlined above, especially of Perelman's work in his second paper
\cite{P2} in which many key ideas of the proofs are sketched or
outlined but complete details of the proofs are often missing.
Since our paper is aimed at both graduate students and researchers
who intend to learn Hamilton's Ricci flow and to understand the
Hamilton-Perelman theory and its application to the geometrization
of three-manifolds, we have made the paper essentially
self-contained so that the proof of the geometrization is readily
accessible to those who are familiar with basics of Riemannian
geometry and elliptic and parabolic partial differential
equations. In particular, our presentation of Hamilton's works in
general follows very closely his original papers, which are
beautifully written. We hope the readers will find this helpful
and convenient and we thank Professor Hamilton for his generosity.

The reader may find many original papers appeared before
Perelman's preprints, particularly those of Hamilton's on the
Ricci flow, in the book ``Collected Papers on Ricci Flow"
\cite{CCCY} in which the editors provided numerous helpful
footnotes. For introductory materials to the Hamilton-Perelman
theory of Ricci flow, we also refer the reader to the recent book
by B. Chow and D. Knopf \cite{CK04} and the forthcoming book by B.
Chow, P. Lu and L. Ni \cite{CLN}. We remark that there have also
appeared several sets of notes on Perelman's work which cover part
of the materials that are needed for the geometrization program,
especially the notes by Kleiner and Lott \cite{KL} from which we
have benefited for the Step 2 of the proof of Theorem 7.1.1.
(After the original version of our paper went to press, Kleiner
and Lott had put up their latest notes \cite{KL2} on Perelman's
papers on the arXiv on May 25, 2006. Also Morgan and Tian posted
their manuscript \cite{MT} on the arXiv on July 25, 2006.) There
also have appeared several survey articles, such as Cao-Chow
\cite{CC99}, Milnor \cite{Mil}, Anderson \cite{An} and Morgan
\cite{Mor}, on the geometrization of three-manifolds via the Ricci
flow.

We are very grateful to Professor S.-T. Yau, who suggested us to
write this paper based on our notes, for introducing us to the
wonderland of the Ricci flow. His vision and strong belief in the
Ricci flow encouraged us to persevere. We also thank him for his
many suggestions and constant encouragement. Without him, it would
be impossible for us to finish this paper. We are enormously
indebted to Professor Richard Hamilton for creating the Ricci flow
and developing the entire program to approach the geometrization
of three-manifolds. His work on the Ricci flow and other geometric
flows has influenced on virtually everyone in the field. The first
author especially would like to thank Professor Hamilton for
teaching him so much about the subject over the past twenty years,
and for his constant encouragement and friendship.

We are indebted to Dr. Bing-Long Chen, who contributed a great
deal in the process of writing this paper. We benefited a lot from
constant discussions with him on the subjects of geometric flows
and geometric analysis. He also contributed many ideas in various
proofs in the paper. We would like to thank Ms. Huiling Gu, a Ph.D
student of the second author, for spending many months of going
through the entire paper and checking the proofs. Without both of
them, it would take much longer time for us to finish this paper.

We also would like to thank Ben Andrews, Simon Brendle, Shu-Cheng
Chang, Ben Chow, Sun-Chin Chu, Panagiota Daskalopoulos, Klaus
Ecker, Gerhard Huisken, Tom Ilmanen, Dan Knopf, Peter Li, Peng Lu,
Lei Ni, Natasa Sesum, Carlo Sinestrari, Luen-fai Tam, Jiaping
Wang, Mu-Tao Wang, and Brian White from whom the authors have
benefited a lot on the subject of geometric flows through
discussions or collaborations.

\enlargethispage{5mm} The first author would like to express his
gratitude to the John Simon Guggenheim Memorial Foundation, the
National Science Foundation (grants DMS-0354621 and DMS-0506084),
and the Outstanding Overseas Young Scholar Fund of Chinese
National Science Foundation for their support for the research in
this paper. He also would like to thank Tsinghua University in
Beijing for its hospitality and support while he was working
there. The second author wishes to thank his wife, Danlin Liu, for
her understanding and support over all these years. The second
author is also indebted to the National Science Foundation of
China for the support in his work on geometric flows, some of
which has been incorporated in this paper. The last part of the
work in this paper was done and the material in Chapter 3, Chapter
6 and Chapter 7 was presented while the second author was visiting
the Harvard Mathematics Department in the fall semester of 2005
and the early spring semester of 2006. He wants to especially
thank Professor Shing-Tung Yau, Professor Cliff Taubes and
Professor Daniel W. Stroock for the enlightening comments and
encouragement during the lectures. Also he gratefully acknowledges
the hospitality and the financial support of Harvard University.
\newpage


\part{{\Large Evolution Equations}}



In this chapter, we collect some basic material on Hamilton's
Ricci flow. In Section 1.1, we introduce the Ricci flow equation
of Hamilton \cite{Ha82} and examine some special solutions. The
short time existence and uniqueness theorem of the Ricci flow on a
compact manifold \cite{Ha82} (also cf. \cite{De}) is presented in
Section 1.2. Evolution equations of curvatures under the Ricci
flow \cite{Ha82, Ha86} are described in Section 1.3. In Section
1.4, we recall Shi's local derivative estimate \cite{Sh89}, which
plays an important role in the Ricci flow. Perelman's two
functionals in \cite{P1} and their monotonicity properties are
discussed in Section 1.5.

\section{The Ricci Flow}

In this section, we introduce Hamilton's Ricci flow and examine
some special solutions. The presentation follows closely
\cite{Ha82, Ha95F}.

Let $M$ be an $n$-dimensional complete Riemannian manifold with
the Riemannian metric $g_{ij}$.  The Levi-Civita connection is
given by the Christoffel symbols
$$
\Gamma^k_{ij}=\frac{1}{2}g^{kl}\left( \frac{\partial
g_{jl}}{\partial x^i}+\frac{\partial g_{il}}{\partial
x^j}-\frac{\partial g_{ij}}{\partial x^l}\right)
$$
where $g^{ij}$ is the inverse of $g_{ij}$. The summation
convention of summing over repeated indices is used here and
throughout the book. The Riemannian curvature tensor is given by
$$
R^k_{ijl}=\frac{\partial \Gamma^k_{jl}}{\partial x^i}
-\frac{\partial \Gamma^k_{il}}{\partial
x^j}+\Gamma^k_{ip}\Gamma^p_{jl}-\Gamma^k_{jp}\Gamma^p_{il}.
$$
We lower the index to the third position, so that
$$
R_{ijkl}=g_{kp}R^p_{ijl}.
$$
The curvature tensor $R_{ijkl}$ is anti-symmetric in the pairs
$i,\; j$ and $k,\; l$ and symmetric in their interchange:
$$
R_{ijkl}=-R_{jikl}=-R_{ijlk}=R_{klij}.
$$
Also the first Bianchi identity holds \be
R_{ijkl}+R_{jkil}+R_{kijl}=0.              
\ee The Ricci tensor is the contraction
$$
R_{ik}=g^{jl}R_{ijkl},
$$
and the scalar curvature is
$$
R=g^{ij}R_{ij}.
$$

We denote the covariant derivative of a vector field
$v=v^j\frac{\partial}{\partial x^j}$ by
$$
\nabla_iv^j=\frac{\partial v^j}{\partial x^i}+\Gamma^j_{ik}v^k
$$
and of a $1$-form by
$$
\nabla_iv_j=\frac{\partial v_j}{\partial x^i}-\Gamma^k_{ij}v_k.
$$
These definitions extend uniquely to tensors so as to preserve the
product rule and contractions. For the exchange of two covariant
derivatives, we have
\begin{align}
\nabla_i\nabla_jv^l-\nabla_j\nabla_iv^l&=R^l_{ijk}v^k,\\
\nabla_i\nabla_jv_k-\nabla_j\nabla_iv_k&=R_{ijkl}g^{lm}v_m,
\end{align}
and similar formulas for more complicated tensors.
 The second Bianchi identity is given by
\be \nabla_mR_{ijkl}+\nabla_iR_{jmkl}+\nabla_jR_{mikl}=0.
\ee

For any tensor $T=T^i_{jk}$ we define its length by
$$
|T^i_{jk}|^2=g_{il}g^{jm}g^{kp}T^i_{jk}T^l_{mp},
$$
and we define its Laplacian by
$$
\Delta T^i_{jk}=g^{pq}\nabla_p\nabla_qT^i_{jk},
$$
the trace of the second iterated covariant derivatives. Similar
definitions hold for more general tensors.

The {\bf Ricci flow}\index{Ricci flow} of Hamilton \cite{Ha82} is
the evolution equation \be
\frac{\partial g_{ij}}{\partial t}=-2R_{ij} 
\ee for a family of Riemannian metrics $g_{ij}(t)$ on $M$. It is a
nonlinear system of second order partial differential equations on
metrics.

In order to get a feel for the Ricci flow (1.1.5) we first present
some examples of specific solutions (cf. Section 2 of
\cite{Ha95F}). \vskip 0.2cm

(1) Einstein metrics \vskip 0.1cm

A Riemannian metric $g_{ij}$ is called {\bf
Einstein}\index{Einstein! metric} if
$$
R_{ij}=\lambda g_{ij}
$$
for some constant $\lambda$. A smooth manifold $M$ with an
Einstein metric is called an {\bf Einstein
manifold}\index{Einstein! manifold}.

If the initial metric is Ricci flat, so that $R_{ij}=0$, then
clearly the metric does not change under (1.1.5). Hence any Ricci
flat metric is a stationary solution of the Ricci flow. This
happens, for example, on a flat torus or on any $K3$-surface with
a Calabi-Yau metric.

If the initial metric is Einstein with positive scalar curvature,
then the metric will shrink under the Ricci flow by a
time-dependent factor. Indeed, since the initial metric is
Einstein, we have
$$
R_{ij}(x,0)=\lambda g_{ij}(x,0), \quad \forall x\in M
$$
and some $\lambda>0$. Let
$$
g_{ij}(x,t)=\rho^2(t)g_{ij}(x,0).
$$
{}From the definition of the Ricci tensor, one sees that
$$
R_{ij}(x,t)=R_{ij}(x,0)=\lambda g_{ij}(x,0).
$$
Thus the equation (1.1.5) corresponds to
$$
\frac{\partial (\rho^2(t) g_{ij}(x,0))}{\partial t}=-2\lambda
g_{ij}(x,0).
$$
This gives the ODE \be
\frac{d\rho}{dt}=-\frac{\lambda}{\rho}, 
\ee whose solution is given by
$$
\rho^2(t)=1-2\lambda t.
$$
Thus the evolving metric $g_{ij}(x,t)$ shrinks homothetically to a
point as $t\rightarrow T=1/{2\lambda}$. Note that as $t\rightarrow
T$, the scalar curvature becomes infinite like $1/(T-t)$.

By contrast, if the initial metric is an Einstein metric of
negative scalar curvature, the metric will expand homothetically
for all times. Indeed if
$$
R_{ij}(x,0)=-\lambda g_{ij}(x,0)
$$
with $\lambda >0$ and
$$
g_{ij}(x,t)=\rho^2(t)g_{ij}(x,0).
$$
Then $\rho(t)$ satisfies the ODE \be
\frac{d\rho}{dt}=\frac{\lambda}{\rho}, 
\ee with the solution
$$
\rho^2(t)=1+2\lambda t.
$$
Hence the evolving metric $g_{ij}(x,t)=\rho^2(t)g_{ij}(x,0)$
exists and expands homothetically for all times, and the curvature
will fall back to zero like $-1/t$. Note that now the evolving
metric $g_{ij}(x,t)$ only goes back in time to $-1/{2\lambda}$,
when the metric explodes out of a single point in a ``big bang".
\vskip 0.2cm

(2) Ricci Solitons \vskip 0.1cm

The concept of a Ricci soliton is introduced by Hamilton
\cite{Ha88}. We will call a solution to an evolution equation
which moves under a one-parameter subgroup of the symmetry group
of the equation a {\bf steady soliton}\index{soliton!steady}. The
symmetry group of the Ricci flow contains the full diffeomorphism
group. Thus a solution to the Ricci flow (1.1.5) which moves by a
one-parameter group of diffeomorphisms $\varphi_t$ is called a
\textbf{steady Ricci soliton}\index{Ricci soliton!steady }.

If $\varphi_t$ is a one-parameter group of diffeomorphisms
generated by a vector field $V$ on $M$, then the Ricci soliton is
given by \be
g_{ij}(x,t)=\varphi^*_tg_{ij}(x,0) 
\ee which implies that the Ricci term $-2 Ric$ on the RHS of
(1.1.5) is equal to the Lie derivative $\mathcal L_V g$ of the
evolving metric $g$. In particular, the initial metric
$g_{ij}(x,0)$ satisfies the following steady Ricci soliton
equation \be
2R_{ij}+g_{ik}\nabla_j V^k+g_{jk}\nabla_i V^k=0 .
\ee If the vector field $V$ is the gradient of a function $f$ then
the soliton is called a {\bf steady gradient Ricci
soliton}\index{Ricci soliton!gradient }. Thus \be
R_{ij}+\nabla_i\nabla_jf=0,
\quad \text{or}\quad Ric+{\nabla}^2f=0, 
\ee is the steady gradient Ricci soliton equation.

Conversely, it is clear that a metric $g_{ij}$ satisfying (1.1.10)
generates a steady gradient Ricci soliton $g_{ij}(t)$ given by
(1.1.8). For this reason we also often call such a metric $g_{ij}$
a steady gradient Ricci soliton and do not necessarily distinguish
it with the solution $g_{ij} (t)$ it generates.

More generally, we can consider a solution to the Ricci flow
(1.1.5) which moves by diffeomorphisms and also shrinks or expands
by a (time-dependent) factor at the same time. Such a solution is
called a homothetically {\bf shrinking}\index{Ricci soliton!
shrinking} or homothetically {\bf expanding}\index{Ricci soliton!
expanding} \textbf{Ricci soliton}. The equation for a homothetic
Ricci soliton is \be 2R_{ij}+g_{ik}\nabla_j V^k+g_{jk}\nabla_i
V^k-2\lambda g_{ij}=0,
\ee or for a homothetic gradient Ricci soliton, \be
R_{ij}+\nabla_i\nabla_jf-\lambda g_{ij}=0, 
\ee where $\lambda$ is the homothetic constant. For $\lambda>0$
the soliton is shrinking, for $\lambda<0$ it is expanding. The
case $\lambda=0$ is a steady Ricci soliton, the case $V=0$ (or $f$
being a constant function) is an Einstein metric. Thus Ricci
solitons can be considered as natural extensions of Einstein
metrics. In fact, the following result states that there are no
nontrivial gradient steady or expanding Ricci solitons on any
compact manifold.

We remark that if the underlying manifold $M$ is a complex
manifold and the initial metric is K\"ahler, then it is well known
(see, e.g., \cite {Ha95, Cao85}) that the solution metric to the
Ricci flow (1.1.5) remains K\"ahler. For this reason, the Ricci
flow on a K\"ahler manifold is called the {\bf K\"ahler-Ricci
flow}\index{K\"ahler-Ricci flow}. A (steady, or shrinking, or
expanding) Ricci soliton to the K\"ahler-Ricci flow is called a
({\bf steady}, or {\bf shrinking}, or {\bf expanding} repectively)
{\bf K\"ahler-Ricci soliton}. \index{K\"ahler-Ricci soliton!
steady} \index{K\"ahler-Ricci soliton! shrinking}
\index{K\"ahler-Ricci soliton! expanding} The following result,
essentially due to Hamilton (cf. Theorem 20.1 of \cite{Ha95F}, or
Proposition 5.20 of \cite{CK04}), says there are no nontrivial
compact steady and expanding Ricci solitons.

\begin{proposition}
On a compact n-dimensional manifold $M$, a gradient steady or
expanding Ricci soliton is necessarily an Einstein metric.
\end{proposition}

\begin{pf}
We shall only prove the steady case and leave the expanding case
as an exercise. Our argument here follows that of Hamilton
\cite{Ha95F}.

Let $g_{ij}$ be a complete steady gradient Ricci soliton on a
manifold $M$ so that
$$
R_{ij}+\nabla_i\nabla_jf=0.
$$

Taking the trace, we get \be
R+\Delta f=0. 
\ee Also, taking the covariant derivatives of the Ricci soliton
equation, we have
$$
\nabla_i\nabla_j\nabla_kf-\nabla_j\nabla_i\nabla_kf
=\nabla_jR_{ik}-\nabla_iR_{jk}.
$$
On the other hand, by using the commutating formula (1.1.3), we
otain
$$
\nabla_i\nabla_j\nabla_kf-\nabla_j\nabla_i\nabla_kf
=R_{ijkl}\nabla_lf.
$$
Thus
$$
\nabla_iR_{jk}-\nabla_jR_{ik}+R_{ijkl}\nabla_lf=0.
$$
Taking the trace on $j$ and $k$, and using the contracted second
Bianchi identity \be
\nabla_j R_{ij}=\frac {1}{2}\nabla_i R, 
\ee we get
$$
\nabla_iR-2R_{ij}\nabla_jf=0.
$$
Then
$$
\nabla_i(|\nabla f|^2+R)=2\nabla_jf(\nabla_i\nabla_jf+R_{ij})=0.
$$
Therefore \be
R+|\nabla f|^2=C 
\ee for some constant $C$.

Taking the difference of (1.1.13) and (1.1.15), we get \be
\Delta f -|\nabla f|^2=-C. 
\ee We claim $C=0$ when $M$ is compact. Indeed, this follows
either from \be 0=-\int_M \Delta(e^{-f}) dV
=\int_M (\Delta f-|\nabla f|^2)e^{-f}dV, 
\ee or from considering (1.1.16) at both the maximum point and
minimum point of $f$. Then, by integrating (1.1.16) we obtain
$$
\int_M |\nabla f|^2 dV=0.
$$
Therefore $f$ is a constant and $g_{ij}$ is Ricci flat.
\end{pf}

\begin{remark}
By contrast, there do exist nontrivial compact gradient shrinking
Ricci solitons (see Koiso \cite{Ko}, Cao \cite{Cao94} and Wang-Zhu
\cite{WZ} ). Also, there exist complete noncompact steady gradient
Ricci solitons that are not Ricci flat. In two dimensions Hamilton
\cite{Ha88} wrote down the first such example on $\mathbb R^2$,
called the {\bf cigar soliton}\index{soliton!cigar}, where the
metric is given by \be
ds^2=\frac {dx^2 +dy^2} {1+x^2+y^2}, 
\ee and the vector field is radial, given by $V=-\partial/\partial
r=-(x\partial/\partial x+y\partial/\partial y).$ This metric has
positive curvature and is asymptotic to a cylinder of finite
circumference $2\pi$ at $\infty$. Higher dimensional examples were
found by Robert Bryant \cite{BR} on $\mathbb R^n$ in the
Riemannian case, and by the first author \cite{Cao94} on $\mathbb
C^n$ in the K\"ahler case. These examples are complete,
rotationally symmetric, of positive curvature and found by solving
certain nonlinear ODE (system). Noncompact expanding solitons were
also constructed by the first author \cite{Cao94}. More recently,
Feldman, Ilmanen and Knopf \cite{FIM} constructed new examples of
noncompact shrinking and expanding K\"ahler-Ricci solitons.
\end{remark}

\section{Short-time Existence and Uniqueness}

In this section we establish the short-time existence and
uniqueness result \cite{Ha82, De} for the Ricci flow (1.1.5) on a
compact $n$-dimensional manifold $M$. Our presentation follows
closely Hamilton \cite{Ha82} and DeTurck \cite{De}.

We will see that the Ricci flow is a system of second order
nonlinear weakly parabolic partial differential equations. In
fact, as observed by Hamilton \cite{Ha82}, the degeneracy of the
system is caused by the diffeomorphism group of $M$ which acts as
the gauge group of the Ricci flow. For any diffeomorphism
$\varphi$ of $M$, we have $\Ric(\varphi^*(g))=\varphi^*(\Ric(g))$.
Thus, if $g(t)$ is a solution to the Ricci flow (1.1.5), so is
$\varphi^*(g(t))$. Because the Ricci flow (1.1.5) is only weakly
parabolic, even the existence and uniqueness result on a compact
manifold does not follow from standard PDE theory. The short-time
existence and uniqueness result in the compact case is first
proved by Hamilton \cite{Ha82} using the Nash-Moser implicit
function theorem. Shortly after Denis De Turck \cite{De} gave a
much simpler proof using the gauge fixing idea which we will
present here. In the noncompact case, the short-time existence was
established by Shi \cite{Sh89} in 1989, but the uniqueness result
has been proved only very recently by Bing-Long Chen and the
second author. These results will be presented at the end of this
section.

First, let us follow Hamilton \cite{Ha82} to examine the equation
of the Ricci flow and see it is weakly parabolic. Let $M$ be a
compact $n$-dimensional Riemannian manifold. The Ricci flow
equation is a second order nonlinear partial differential system
\be
\frac{\partial}{\partial t}g_{ij}=E(g_{ij}),
\ee for a family of Riemannian metrics $g_{ij}(\cdot,t)$ on $M$,
where
\begin{align*}
E(g_{ij})&  =  -2R_{ij}\\
&  =  -2\left(\frac{\partial}{\partial x^k}
\Gamma^k_{ij}-\frac{\partial}{\partial x^i}\Gamma^k_{kj}
+\Gamma^k_{kp}\Gamma^p_{ij}-\Gamma^k_{ip}\Gamma^p_{kj}\right)\\
&  = \frac{\partial}{\partial x^i}
\left\{g^{kl}\frac{\partial}{\partial x^j}g_{kl}\right\}
-\frac{\partial}{\partial x^k}
\left\{g^{kl}\left(\frac{\partial}{\partial x^i}g_{jl}
+\frac{\partial}{\partial x^j}g_{il}
-\frac{\partial}{\partial x^l}g_{ij}\right)\right\}\\
&\quad  +2\Gamma^k_{ip}\Gamma^p_{kj}-2\Gamma^k_{kp}\Gamma^p_{ij}.
\end{align*}
The linearization of this system is
$$
\frac{\partial \tilde{g}_{ij}}{\partial t}=DE(g_{ij})
\tilde{g}_{ij}
$$
where $\tilde{g}_{ij}$ is the variation in $g_{ij}$ and $DE$ is
the derivative of $E$ given by
\begin{align*}
DE(g_{ij})\tilde{g}_{ij}&  = g^{kl}\left\{\frac{\partial^2
\tilde{g}_{kl}}{\partial x^i \partial x^j}-\frac{\partial^2
\tilde{g}_{jl}}{\partial x^i \partial x^k}-\frac{\partial^2
\tilde{g}_{il}}{\partial x^j \partial x^k}+\frac{\partial^2
\tilde{g}_{ij}}{\partial x^k \partial x^l}\right\}\\
& \quad + \text{(lower order terms).}
\end{align*}
We now compute the symbol of $DE$. This is to take the highest
order derivatives and replace $\frac{\partial}{\partial x^i}$ by
the Fourier transform variable $\zeta_i$. The symbol of the linear
differential operator $DE(g_{ij})$ in the direction
$\zeta=(\zeta_1 , \ldots,\zeta_n )$ is
$$
\sigma DE(g_{ij})(\zeta)\tilde{g}_{ij} =g^{kl}(\zeta_i \zeta_j
\tilde{g}_{kl} + \zeta_k \zeta_l \tilde{g}_{ij} -\zeta_i \zeta_k
\tilde{g}_{jl} -\zeta_j \zeta_k \tilde{g}_{il}).
$$
To see what the symbol does, we can always assume $\zeta$ has
length 1 and choose coordinates at a point such that
$$
      \left\{
       \begin{array}{lll}
       \displaystyle
g_{ij}=\delta_{ij},
          \\[4mm]
      \displaystyle
\zeta=(1,0,\ldots,0).
       \end{array}
      \right.
$$
Then
\begin{align*}
(\sigma DE(g_{ij})(\zeta))(\tilde{g}_{ij}) &  =
\tilde{g}_{ij}+\delta_{i1}\delta_{j1}(\tilde{g}_{11}
+\cdots+\tilde{g}_{nn})\\
&\quad -\delta_{i1}\tilde{g}_{1j}-\delta_{j1}\tilde{g}_{1i},
\end{align*}
i.e.,
\begin{align*}
[\sigma DE(g_{ij})(\zeta)(\tilde{g}_{ij})]_{11}
&  = \tilde{g}_{22}+\cdots+\tilde{g}_{nn} ,\\
[\sigma DE(g_{ij})(\zeta)(\tilde{g}_{ij})]_{1k} &  = 0 ,
\;\;\;\;\;\; \mbox{ if }\; k\neq1,\\
[\sigma DE(g_{ij})(\zeta)(\tilde{g}_{ij})]_{kl} & =
\tilde{g}_{kl}, \;\;\;\;\mbox{ if }\; k\neq1 ,l\neq1.
\end{align*}
In particular
$$
(\tilde{g}_{ij})=\left(
\begin{array}{cccc}\ast&  \ast&  \cdots&  \ast\\
\ast&  0&  \cdots&  0\\
\vdots&  \vdots&  \ddots&  \vdots\\
\ast&  0&  \cdots&  0
\end{array}
\right)
$$
are zero eigenvectors of the symbol. The presence of the zero
eigenvalue shows that the system cannot be strictly parabolic.

Next, instead of considering the system (1.2.1) (or the Ricci flow
equation (1.1.5)), we follow a trick of DeTurck \cite{De} to
consider a modified evolution equation, which is equivalent to the
Ricci flow up to diffeomorphisms and turns out to be strictly
parabolic so that the standard theory of parabolic equations
applies.

Suppose $\hat{g}_{ij}(x,t)$ is a solution of the Ricci flow
(1.1.5), and $\varphi_t: M\rightarrow M$ is a family of
diffeomorphisms of $M$. Let
$$
g_{ij}(x,t)=\varphi^\ast_t\hat{g}_{ij}(x,t)
$$
be the pull-back metrics. We now want to find the evolution
equation for the metrics $g_{ij}(x,t)$.

Denote by
$$
y(x,t)=\varphi_t(x)=\{y^1(x,t),y^2(x,t),\ldots,y^n(x,t)\}
$$
in local coordinates. Then \be g_{ij}(x,t) =\frac{\partial
y^\alpha}{\partial x^i} \frac{\partial y^\beta}{\partial
x^j}\hat{g}_{\alpha\beta}(y,t),
\ee and
\begin{align*}
\frac{\partial}{\partial t}g_{ij}(x,t) & =\frac{\partial}{\partial
t}\left[\frac{\partial y^\alpha}{\partial x^i}\frac{\partial
y^\beta}{\partial x^j}
\hat{g}_{\alpha\beta}(y,t)\right]\\
&  = \frac{\partial y^\alpha}{\partial x^i} \frac{\partial
y^\beta}{\partial x^j}\frac{\partial}{\partial t}
\hat{g}_{\alpha\beta}(y,t)+\frac{\partial}{\partial x^i}
\left(\frac{\partial y^\alpha}{\partial t}\right)
\frac{\partial y^\beta}{\partial x^j}\hat{g}_{\alpha\beta}(y,t)\\
&\quad +\frac{\partial y^\alpha}{\partial x^i}
\frac{\partial}{\partial x^j}\left(\frac{\partial
y^\beta}{\partial t}\right) \hat{g}_{\alpha\beta}(y,t).
\end{align*}
Let us choose a normal coordinate $\{x^i\}$ around a fixed point
$p\in M$ such that $\frac{\partial g_{ij}}{\partial x^k}=0$ at
$p$. Since
$$
\frac{\partial}{\partial t}\hat{g}_{\alpha\beta}(y,t)
=-2\hat{R}_{\alpha\beta}(y,t)+\frac{\partial
\hat{g}_{\alpha\beta}}{\partial y^\gamma}\frac{\partial
y^\gamma}{\partial t} ,
$$
we have in the normal coordinate,
\begin{align*}
&\frac{\partial}{\partial t}g_{ij}(x,t) \\
&  =-2\frac{\partial y^\alpha}{\partial x^i}\frac{\partial
y^\beta}{\partial x^j}\hat{R}_{\alpha\beta}(y,t)+\frac{\partial
y^\alpha}{\partial x^i}\frac{\partial y^\beta}{\partial
x^j}\frac{\partial \hat{g}_{\alpha\beta}}{\partial
y^\gamma}\frac{\partial y^\gamma}{\partial t}\\
&\quad +\frac{\partial}{\partial x^i}\left(\frac{\partial
y^\alpha}{\partial t}\right)\frac{\partial y^\beta}{\partial
x^j}\hat{g}_{\alpha\beta}(y,t)+\frac{\partial}{\partial
x^j}\left(\frac{\partial y^\beta}{\partial t}\right)\frac{\partial
y^\alpha}{\partial x^i}\hat{g}_{\alpha\beta}(y,t)\\
&  = -2R_{ij}(x,t)+\frac{\partial y^\alpha}{\partial
x^i}\frac{\partial y^\beta}{\partial x^j}\frac{\partial
\hat{g}_{\alpha\beta}}{\partial y^\gamma}\frac{\partial
y^\gamma}{\partial t}+\frac{\partial}{\partial
x^i}\left(\frac{\partial y^\alpha}{\partial
t}\right)\frac{\partial x^k}{\partial
y^\alpha}g_{jk} \\
&\quad+\frac{\partial}{\partial x^j}\left(\frac{\partial
y^\beta}{\partial t}\right)\frac{\partial x^k}{\partial
y^\beta}g_{ik}\\
&  = -2R_{ij}(x,t)+\frac{\partial y^\alpha}{\partial
x^i}\frac{\partial y^\beta}{\partial x^j}\frac{\partial
\hat{g}_{\alpha\beta}}{\partial y^\gamma}\frac{\partial
y^\gamma}{\partial t}+\frac{\partial}{\partial
x^i}\left(\frac{\partial y^\alpha}{\partial t}\frac{\partial
x^k}{\partial
y^\alpha}g_{jk}\right) \\
&\quad+\frac{\partial}{\partial x^j}\left(\frac{\partial
y^\beta}{\partial t}\frac{\partial x^k}{\partial
y^\beta}g_{ik}\right)
-\frac{\partial y^\alpha}{\partial t}\frac{\partial}{\partial
x^i}\left(\frac{\partial x^k}{\partial
y^\alpha}\right)g_{jk}-\frac{\partial y^\beta}{\partial
t}\frac{\partial}{\partial x^j}\left(\frac{\partial x^k}{\partial
y^\beta}\right)g_{ik}.
\end{align*}
The second term on the RHS gives, in the normal coordinate,
\begin{align*}
\frac{\partial y^\alpha}{\partial x^i}\frac{\partial
y^\beta}{\partial x^j}\frac{\partial y^\gamma}{\partial
t}\frac{\partial \hat{g}_{\alpha\beta}}{\partial y^\gamma} &
=\frac{\partial y^\alpha}{\partial x^i}\frac{\partial
y^\beta}{\partial x^j}\frac{\partial y^\gamma}{\partial
t}g_{kl}\frac{\partial}{\partial y^\gamma} \(\frac{\partial
x^k}{\partial y^\alpha}
\frac{\partial x^l}{\partial y^\beta}\)\\
&  = \frac{\partial y^\alpha}{\partial x^i}\frac{\partial
y^\gamma}{\partial t}\frac{\partial}{\partial
y^\gamma}\(\frac{\partial x^k}{\partial
y^\alpha}\)g_{jk}+\frac{\partial y^\beta}{\partial
x^j}\frac{\partial y^\gamma}{\partial t}\frac{\partial}{\partial
y^\gamma}\(\frac{\partial x^k}{\partial y^\beta}\)g_{ik}\\
&  = \frac{\partial y^\alpha}{\partial t}\frac{\partial^2
x^k}{\partial y^\alpha\partial y^\beta}\frac{\partial
y^\beta}{\partial x^i}g_{jk}+\frac{\partial y^\beta}{\partial
t}\frac{\partial^2 x^k}{\partial y^\alpha\partial
y^\beta}\frac{\partial y^\alpha}{\partial x^j}g_{ik}\\
&  = \frac{\partial y^\alpha}{\partial t}\frac{\partial}{\partial
x^i}\(\frac{\partial x^k}{\partial
y^\alpha}\)g_{jk}+\frac{\partial y^\beta}{\partial
t}\frac{\partial}{\partial x^j}\(\frac{\partial x^k}{\partial
y^\beta}\)g_{ik}.
\end{align*}
So we get
\begin{align}
&\frac{\partial}{\partial t}g_{ij}(x,t) \\
&=-2R_{ij}(x,t)+\nabla_i\left(\frac{\partial y^\alpha}{\partial
t}\frac{\partial x^k}{\partial
y^\alpha}g_{jk}\right)+\nabla_j\left(\frac{\partial
y^\beta}{\partial
t}\frac{\partial x^k}{\partial y^\beta}g_{ik}\right).\nn
\end{align}
If we define $y(x,t)=\varphi_t(x)$ by the equations \be
      \left\{
       \begin{array}{lll}
  \frac{\partial y^\alpha}{\partial t}=\frac{\partial y^\alpha}{\partial
x^k}g^{jl}(\Gamma^k_{jl}-\stackrel{o}{\Gamma}^k_{jl}),
          \\[4mm]
  y^\alpha(x,0)=x^\alpha,
       \end{array}
    \right. 
\ee and
$V_i=g_{ik}g^{jl}(\Gamma^k_{jl}-\stackrel{o}{\Gamma}^k_{jl})$, we
get the following evolution equation for the pull-back metric \be
      \left\{
       \begin{array}{lll}
  \frac{\partial}{\partial t}g_{ij}(x,t)
=-2R_{ij}(x,t)+\nabla_iV_j+\nabla_jV_i, \\[4mm]
  g_{ij}(x,0)= \stackrel{o}{g}_{ij}(x),
       \end{array}
    \right. 
\ee where $\stackrel{o}{g}_{ij}(x)$ is the initial metric and
$\stackrel{o}{\Gamma}^k_{jl}$ is the connection of the initial
metric.

\begin{lemma} [{DeTurck \cite{De}}]
The modified evolution equation $(1.2.5)$ is a strictly parabolic
system.
\end{lemma}

\begin{pf}
The RHS of the equation (1.2.5) is given by
\begin{align*}
&  -2R_{ij}(x,t)+\nabla_iV_j+\nabla_jV_i\\
&= \frac{\partial}{\partial x^i}\left\{g^{kl}\frac{\partial
g_{kl}}{\partial x^j}\right\}-\frac{\partial}{\partial
x^k}\left\{g^{kl}\(\frac{\partial g_{jl}}{\partial
x^i}+\frac{\partial
g_{il}}{\partial x^j}-\frac{\partial g_{ij}}{\partial x^l}\)\right\}\\
&\quad  +g_{jk}g^{pq}\frac{\partial}{\partial
x^i}\left\{\frac{1}{2}g^{kl}\(\frac{\partial g_{pl}}{\partial
x^q}+\frac{\partial g_{ql}}{\partial
x^p}-\frac{\partial g_{pq}}{\partial x^l}\)\right\}\\
&\quad  +g_{ik}g^{pq}\frac{\partial}{\partial
x^j}\left\{\frac{1}{2}g^{kl}\(\frac{\partial g_{pl}}{\partial
x^q}+\frac{\partial g_{ql}}{\partial
x^p}-\frac{\partial g_{pq}}{\partial x^l}\)\right\}\\
&\quad  +\text{(lower order terms)}\\
&= g^{kl}\left\{\frac{\partial^2 g_{kl}}{\partial x^i\partial
x^j}-\frac{\partial^2 g_{jl}}{\partial x^i\partial
x^k}-\frac{\partial^2 g_{il}}{\partial x^j\partial
x^k}+\frac{\partial^2 g_{ij}}{\partial x^k\partial
x^l}\right\}\\
&\quad  +\frac{1}{2}g^{pq}\left\{\frac{\partial^2 g_{pj}}{\partial
x^i\partial x^q}+\frac{\partial^2 g_{qj}}{\partial x^i\partial
x^p}-\frac{\partial^2
g_{pq}}{\partial x^i\partial x^j}\right\}\\
&\quad  +\frac{1}{2}g^{pq}\left\{\frac{\partial^2 g_{p
i}}{\partial x^j\partial x^q}+\frac{\partial^2 g_{q i}}{\partial
x^j\partial x^p}-\frac{\partial^2
g_{pq}}{\partial x^i\partial x^j}\right\}\\
&\quad +\text{(lower order terms)}\\
&= g^{kl}\frac{\partial^2 g_{ij}}{\partial x^k\partial x^l}
+\text{(lower order terms)}.
\end{align*}
Thus its symbol is $(g^{kl}\zeta_k\zeta_l)\tilde{g}_{ij}$. Hence
the equation in (1.2.5) is strictly parabolic.
\end{pf}

Now since the equation (1.2.5) is strictly parabolic and the
manifold $M$ is compact,  it follows from the standard theory of
parabolic equations (see for example \cite{LSU}) that (1.2.5) has
a solution for a short time.  From the solution of (1.2.5) we can
obtain a solution of the Ricci flow from (1.2.4) and (1.2.2). This
shows existence. Now we argue the uniqueness of the solution.
Since
$$
\Gamma^k_{jl} = \frac{\partial y^{\alpha}}{\partial x^j}
\frac{\partial y^{\beta}}{\partial x^l} \frac{\partial
x^k}{\partial y^{\gamma}}\hat{\Gamma}^{\gamma}_{\alpha \beta} +
\frac{\partial x^k}{\partial y^{\alpha}} \frac{\partial ^2
y^{\alpha}}{\partial x^j \partial x^l},
$$
the initial value problem (1.2.4) can be written as \be
\begin{cases}
\frac{\partial y^\alpha}{\partial t} =g^{jl}\(\frac{\partial ^2
y^\alpha}{\partial x^j \partial
x^l}-\stackrel{o}{\Gamma}^k_{jl}\frac{\partial
y^{\alpha}}{\partial x^k}+\hat{\Gamma}^{\alpha}_{\gamma
\beta}\frac{\partial y^{\beta}}{\partial x^j}\frac{\partial
y^{\gamma}}{\partial x^l}\),\\[4mm]
  y^\alpha(x,0)=x^\alpha.
       \end{cases} 
\ee This is clearly a strictly parabolic system. For any two
solutions $\hat{g}^{(1)}_{ij}(\cdot,t)$ and
$\hat{g}^{(2)}_{ij}(\cdot,t)$ of the Ricci flow (1.1.5) with the
same initial data, we can solve the initial value problem (1.2.6)
(or equivalently, (1.2.4)) to get two families $\varphi ^{(1)}_t$
and $\varphi ^{(2)}_t$ of diffeomorphisms of $M$. Thus we get two
solutions, ${g}^{(1)}_{ij}(\cdot,t) = (\varphi
^{(1)}_t)^*\hat{g}^{(1)}_{ij}(\cdot,t)$ and
${g}^{(2)}_{ij}(\cdot,t) = (\varphi
^{(2)}_t)^*\hat{g}^{(2)}_{ij}(\cdot,t)$, to the modified evolution
equation (1.2.5) with the same initial metric. The uniqueness
result for the strictly parabolic equation (1.2.5) implies that
$g^{(1)}_{ij} = g^{(2)}_{ij}$. Then by equation (1.2.4) and the
standard uniqueness result of ODE systems, the corresponding
solutions $\varphi ^{(1)}_t$ and $\varphi ^{(2)}_t$ of (1.2.4) (or
equivalently (1.2.6)) must agree. Consequently the metrics
$\hat{g}^{(1)}_{ij}$ and $\hat{g}^{(2)}_{ij}$ must agree also.
Thus we have proved the following result.

\begin{theorem}[{Hamilton \cite{Ha82}, De Turck \cite{De}}]
Let $(M,\, g_{ij}(x))$ be a compact Riemannian manifold. Then
there exists a constant $T>0$ such that the initial value problem
$$
\begin{cases}
  \frac{\partial}{\partial t}g_{ij}(x,t)=-2R_{ij}(x,t)
          \\[4mm]
  g_{ij}(x,0)=g_{ij}(x)
       \end{cases}
$$
has a unique smooth solution $g_{ij}(x,t)$ on $M\times [0,T)$.
\end{theorem}

The case of a noncompact manifold is much more complicated and
involves a huge amount of techniques from the theory of partial
differential equations. Here we will only state the existence and
uniqueness results and refer the reader to the cited references
for the proofs.

The following existence result was obtained by Shi in \cite{Sh89}
published in 1989.

\begin{theorem}[{Shi \cite{Sh89}}]
Let $(M,\, g_{ij}(x))$ be a complete noncompact Riemannian
manifold of dimension $n$ with bounded curvature. Then there
exists a constant $T>0$ such that the initial value problem
$$
\begin{cases}
  \frac{\partial}{\partial t}g_{ij}(x,t)=-2R_{ij}(x,t)
          \\[4mm]
  g_{ij}(x,0)=g_{ij}(x)
       \end{cases}
$$
has a smooth solution $g_{ij}(x,t)$ on $M\times [0,T]$ with
uniformly bounded curvature.
\end{theorem}

The Ricci flow is a heat type equation. It is well-known that the
uniqueness of a heat equation on a complete noncompact manifold is
not always held if there are no further restrictions on the growth
of the solutions. For example, the heat equation on Euclidean
space with zero initial data has a nontrivial solution which grows
faster than $\exp(a|x|^2)$ for any $a>0$ whenever $t>0$. This
implies that even for the standard linear heat equation on
Euclidean space, in order to ensure the uniqueness one can only
allow the solution to grow at most as $\exp(C|x|^2)$ for some
constant $C>0$. Note that on a K\"ahler manifold, the Ricci
curvature is given by $R_{\alpha\bar\beta} = - \frac{\partial
^2}{\partial z^{\alpha}\partial \bar{z}^{\beta}}
\log\det(g_{\gamma\bar{\delta}})$. So the reasonable growth rate
for the uniqueness of the Ricci flow to hold is that the solution
has bounded curvature. Thus the following uniqueness result of
Bing-Long Chen and the second author \cite{CZ05U} is essentially
the best one can hope for.

\begin{theorem}[{Chen-Zhu \cite{CZ05U}}]
Let $(M, \hat{g}_{ij})$ be a complete noncompact Riemannian
manifold of dimension $n$ with bounded curvature. Let
${g}_{ij}(x,t)$ and $\bar{g}_{ij}(x,t)$ be two solutions, defined
on $M\times[0,T]$, to the Ricci flow $(1.1.5)$ with $\hat{g}_{ij}$
as initial data and with bounded curvatures. Then
$g_{ij}(x,t)\equiv \bar{g}_{ij}(x,t)$ on $M\times[0,T]$.
\end{theorem}

\section{Evolution of Curvatures}

The Ricci flow is an evolution equation on the metric. The
evolution for the metric implies a nonlinear heat equation for the
Riemannian curvature tensor $R_{ijkl}$ which we now derive. Our
presentation in this section follows closely the original papers
of Hamilton \cite{Ha82, Ha86}.

\begin{proposition}[{Hamilton \cite{Ha82}}]
Under the Ricci flow $(1.1.5),$ the curvature tensor satisfies the
evolution equation
\begin{align*}
\frac{\partial}{\partial t}R_{ijkl}&  = \Delta
R_{ijkl}+2(B_{ijkl}-B_{ijlk}-B_{iljk}+B_{ikjl})\\
& \quad
-g^{pq}(R_{pjkl}R_{qi}+R_{ipkl}R_{qj}+R_{ijpl}R_{qk}+R_{ijkp}R_{ql})
\end{align*}
where $B_{ijkl}=g^{pr}g^{qs}R_{piqj}R_{rksl}$ and $\Delta$ is the
Laplacian with respect to the evolving metric.
\end{proposition}

\begin{pf}
Choose $\{x^1,\ldots,x^m\}$ to be a normal coordinate system at a
fixed point. At this point, we compute
\begin{align*}
\frac{\partial}{\partial t}\Gamma^h_{jl} &
=\frac{1}{2}g^{hm}\left\{\frac{\partial}{\partial
x^j}\(\frac{\partial}{\partial t}g_{lm}\)+\frac{\partial}{\partial
x^l}\(\frac{\partial}{\partial t}g_{jm}\)-\frac{\partial}{\partial
x^m}\(\frac{\partial}{\partial t}g_{jl}\)\right\}\\
&  = \frac{1}{2}g^{hm}(\nabla_j(-2R_{lm})
+\nabla_l(-2R_{jm})-\nabla_m(-2R_{jl})),\\
\frac{\partial}{\partial t}R_{ijl}^h & = \frac{\partial}{\partial
x^i}\(\frac{\partial}{\partial
t}\Gamma^h_{jl}\)-\frac{\partial}{\partial
x^j}\(\frac{\partial}{\partial t}\Gamma^h_{il}\),\\
\frac{\partial}{\partial t}R_{ijkl} &
=g_{hk}\frac{\partial}{\partial t}R^h_{ijl}+\frac{\partial
g_{hk}}{\partial t}R^h_{ijl}.
\end{align*}

\noindent Combining these identities we get
\begin{align*}
\frac{\partial}{\partial t}R_{ijkl} &
=g_{hk}\bigg[\(\frac{1}{2}\nabla_i
[g^{hm}(\nabla_j(-2R_{lm})+\nabla_l(-2R_{jm})-\nabla_m(-2R_{jl}))]\)\\
&\quad -\(\frac{1}{2}\nabla_j[g^{hm}(\nabla_i(-2R_{lm})
+\nabla_l(-2R_{im})-\nabla_m(-2R_{il}))]\)\bigg]\\
& \quad -2R_{hk}R^h_{ijl}\\
&  = \nabla_i\nabla_kR_{jl}-\nabla_i\nabla_lR_{jk}
-\nabla_j\nabla_kR_{il}+\nabla_j\nabla_lR_{ik}\\
&\quad -R_{ijlp}g^{pq}R_{qk}-R_{ijkp}g^{pq}R_{ql}-2R_{ijpl}g^{pq}R_{qk}\\
& = \nabla_i\nabla_kR_{jl}-\nabla_i\nabla_lR_{jk}
-\nabla_j\nabla_kR_{il}+\nabla_j\nabla_lR_{ik}\\
&\quad -g^{pq}(R_{ijkp}R_{ql}+R_{ijpl}R_{qk}).
\end{align*}
Here we have used the exchanging formula (1.1.3).

Now it remains to check the following identity, which is analogous
to the Simon$'$s identity in extrinsic geometry,
\begin{align}
& \Delta R_{ijkl}+2(B_{ijkl}-B_{ijlk}-B_{iljk}+B_{ikjl})\\
&  = \nabla_i\nabla_k R_{jl}-\nabla_i\nabla_l
R_{jk}-\nabla_j\nabla_k R_{il}+\nabla_j\nabla_l R_{ik} \nn\\
&\quad +g^{pq}(R_{pjkl}R_{qi}+R_{ipkl}R_{qj}). \nn
\end{align}
Indeed, from the second Bianchi identity (1.1.4), we have
\begin{align*}
\Delta R_{ijkl}&  = g^{pq}\nabla_p\nabla_q R_{ijkl}\\
&  = g^{pq}\nabla_p\nabla_i
R_{qjkl}-g^{pq}\nabla_p\nabla_jR_{qikl}.
\end{align*}
Let us examine the first term on the RHS. By using the exchanging
formula (1.1.3) and the first Bianchi identity (1.1.1), we have
\begin{align*}
& g^{pq}\nabla_p\nabla_iR_{qjkl}-g^{pq}\nabla_i\nabla_pR_{qjkl}\\
& = g^{pq}g^{mn}(R_{piqm}R_{njkl}+R_{pijm}R_{qnkl}
+R_{pikm}R_{qjnl}+R_{pilm}R_{qjkn})\\
& = R_{im}g^{mn}R_{njkl}+g^{pq}g^{mn}R_{pimj}(R_{qkln}+R_{qlnk})\\
&\quad +g^{pq}g^{mn}R_{pikm}R_{qjnl}+g^{pq}g^{mn}R_{pilm}R_{qjkn}\\
& = R_{im}g^{mn}R_{njkl}-B_{ijkl}+B_{ijlk}-B_{ikjl}+B_{iljk},
\end{align*}
while using the contracted second Bianchi identity \be
g^{pq}\nabla_pR_{qjkl}=\nabla_kR_{jl}-\nabla_lR_{jk},
\ee we have
$$
g^{pq}\nabla_i\nabla_pR_{qjkl}
=\nabla_i\nabla_kR_{jl}-\nabla_i\nabla_lR_{jk}.
$$
Thus
\begin{align*}
& g^{pq}\nabla_p\nabla_iR_{qjkl} \\
&
=\nabla_i\nabla_kR_{jl}-\nabla_i\nabla_lR_{jk}-(B_{ijkl}-B_{ijlk}
-B_{iljk}+B_{ikjl})+g^{pq}R_{pjkl}R_{qi}.
\end{align*}
Therefore we obtain
\begin{align*}
& \Delta R_{ijkl} \\
&  = g^{pq}\nabla_p\nabla_iR_{qjkl}-g^{pq}\nabla_p\nabla_jR_{qikl}\\
&
=\nabla_i\nabla_kR_{jl}-\nabla_i\nabla_lR_{jk}-(B_{ijkl}-B_{ijlk}
-B_{iljk}+B_{ikjl})+g^{pq}R_{pjkl}R_{qi}\\
&\quad
-\nabla_j\nabla_kR_{il}+\nabla_j\nabla_lR_{ik}+(B_{jikl}-B_{jilk}
-B_{jlik}+B_{jkil})-g^{pq}R_{pikl}R_{qj}\\
& = \nabla_i\nabla_kR_{jl}-\nabla_i\nabla_lR_{jk}
-\nabla_j\nabla_kR_{il}+\nabla_j\nabla_lR_{ik}\\[2mm]
&\quad +g^{pq}(R_{pjkl}R_{qi}+R_{ipkl}R_{qj})-2(B_{ijkl}-B_{ijlk}
-B_{iljk}+B_{ikjl})
\end{align*}
as desired, where in the last step we used the symmetries \be
B_{ijkl}=B_{klij}=B_{jilk}.
\ee
\end{pf}

\begin{corollary}[{Hamilton \cite{Ha82}}]
The Ricci curvature satisfies the evolution equation
$$
\frac{\partial}{\partial t}R_{ik} =\Delta
R_{ik}+2g^{pr}g^{qs}R_{piqk}R_{rs}-2g^{pq}R_{pi}R_{qk}.
$$
\end{corollary}

\begin{pf}
\begin{align*}
\frac{\partial}{\partial t}R_{ik} &
=g^{jl}\frac{\partial}{\partial t}R_{ijkl}
+\(\frac{\partial}{\partial t}g^{jl}\)R_{ijkl}\\
& = g^{jl}[\Delta R_{ijkl}+2(B_{ijkl}-B_{ijlk}
-B_{iljk}+B_{ikjl})\\
&\quad -g^{pq}(R_{pjkl}R_{qi}+R_{ipkl}R_{qj}
+R_{ijpl}R_{qk}+R_{ijkp}R_{ql})]\\
&\quad -g^{jp}\(\frac{\partial}{\partial t}g_{pq}\)g^{ql}R_{ijkl}\\
& = \Delta
R_{ik}+2g^{jl}(B_{ijkl}-2B_{ijlk})+2g^{pr}g^{qs}R_{piqk}R_{rs}\\
&\quad -2g^{pq}R_{pk}R_{qi}.
\end{align*}

We claim that $g^{jl}(B_{ijkl}-2B_{ijlk})=0$. Indeed by using the
first Bianchi identity, we have
\begin{align*}
g^{jl}B_{ijkl}&  = g^{jl}g^{pr}g^{qs}R_{piqj}R_{rksl}\\
&  = g^{jl}g^{pr}g^{qs}R_{pqij}R_{rskl}\\
&  = g^{jl}g^{pr}g^{qs}(R_{piqj}-R_{pjqi})(R_{rksl}-R_{rlsk})\\
&  = 2g^{jl}(B_{ijkl}-B_{ijlk})
\end{align*}
as desired.

Thus we obtain
$$
\frac{\partial}{\partial t}R_{ik} =\Delta
R_{ik}+2g^{pr}g^{qs}R_{piqk}R_{rs}-2g^{pq}R_{pi}R_{qk}.
$$
\end{pf}

\begin{corollary}[{Hamilton \cite{Ha82}}]
The scalar curvature satisfies the evolution equation
$$
\frac{\partial R}{\partial t}=\Delta R+2|\Ric|^2.
$$
\end{corollary}

\begin{pf}
\begin{align*}
\frac{\partial R}{\partial t}&  = g^{ik}\frac{\partial
R_{ik}}{\partial t}+\(-g^{ip}\frac{\partial g_{pq}}{\partial
t}g^{qk}\)R_{ik}\\
&  = g^{ik}(\Delta R_{ik}+2g^{pr}g^{qs}R_{piqk}R_{rs}
-2g^{pq}R_{pi}R_{qk})+2R_{pq}R_{ik}g^{ip}g^{qk}\\
&  = \Delta R+2|\Ric|^2.
\end{align*}
\end{pf}

To simplify the evolution equations of curvatures, we will follow
Hamilton \cite{Ha86} and represent the curvature tensors in an
orthonormal frame and evolve the frame so that it remains
orthonormal. More precisely, let us pick an abstract vector bundle
$V$ over $M$ isomorphic to the tangent bundle $TM$. Locally, the
frame $F=\{ F_1,\ldots,F_a,\ldots,F_n\}$ of $V$ is given by
$F_a=F^i_a\frac{\partial}{\partial x^i}$ with the isomorphism
$\{F^i_a\}$. Choose $\{F^i_a\}$ at $t=0$ such that
$F=\{F_1,\ldots,F_a,\ldots,F_n\}$ is an orthonormal frame at
$t=0$, and evolve $\{F_a^i\}$ by the equation
$$
\frac{\partial}{\partial t}F^i_a=g^{ij}R_{jk}F^k_a.
$$
Then the frame $F=\{F_1,\ldots,F_a,\ldots,F_n\}$ will remain
orthonormal for all times since the pull back metric on $V$
$$
h_{ab}=g_{ij}F^i_aF^j_b
$$
remains constant in time. In the following we will use indices
$a,b,\ldots$ on a tensor to denote its components in the evolving
orthonormal frame. In this frame we have the following:
\begin{align*}
R_{abcd}&  = F^i_aF^j_bF^k_cF^l_dR_{ijkl},\\
\Gamma^a_{jb}&  = F^a_i\frac{\partial F_b^i}{\partial
x^j}+\Gamma_{jk}^iF_i^aF_b^k,\quad ((F^a_i)=(F^i_a)^{-1})\\
\nabla_iV^a&  = \frac{\partial}{\partial
x^i}V^a+\Gamma_{ib}^aV^b,\\
\nabla_bV^a&  = F^i_b\nabla_iV^a,
\end{align*}
where $\Gamma_{jb}^a$ is the metric connection of the vector
bundle $V$ with the metric $h_{ab}$. Indeed, by direct
computations,
\begin{align*}
\nabla_iF_b^j&  = \frac{\partial F_b^j}{\partial
x^i}+F_b^k\Gamma_{ik}^j-F_c^j\Gamma_{ib}^c\\
&  = \frac{\partial F_b^j}{\partial x^i}+F_b^k\Gamma_{ik}^j
-F^j_c\(F_k^c\frac{\partial F_b^k}{\partial
x^i}+\Gamma_{ik}^lF^c_lF_b^k\)\\
&  = 0,\\
\nabla_ih_{ab}&  = \nabla_i(g_{ij}F^i_aF^j_b)=0.
\end{align*}
So
$$
\nabla_aV_b=F^i_aF^j_b\nabla_iV^j,\hskip 2.5cm$$ and
\begin{align*}
\Delta R_{abcd}&  = \nabla_l\nabla_lR_{abcd}\\
&  = g^{ij}\nabla_i\nabla_jR_{abcd}\\
&  = g^{ij}F_a^kF_b^lF_c^mF_d^n\nabla_i\nabla_jR_{klmn}.
\end{align*}

In an orthonormal frame $F=\{F_1,\ldots,F_a,\ldots,F_n\}$, the
evolution equations of curvature tensors become
\begin{align}
\frac{\partial}{\partial t}R_{abcd}&=\Delta
R_{abcd}+2(B_{abcd}-B_{abdc}-B_{adbc}+B_{acbd}) \\  
\frac{\partial}{\partial t}R_{ab}&=\Delta
R_{ab}+2R_{acbd}R_{cd}   \\
\frac{\partial}{\partial t}R&=\Delta R+2|\Ric|^2  
\end{align}
where $B_{abcd}=R_{aebf}R_{cedf}$.

Equation (1.3.4) is a reaction-diffusion equation. We can
understand the quadratic terms of this equation better if we think
of the curvature tensor $R_{abcd}$ as a symmetric bilinear form on
the two-forms $\Lambda^2(V)$ given by the formula
$$
Rm(\varphi,\psi)=R_{abcd}\varphi_{ab}\psi_{cd},\qquad \text{for
}\; \varphi,\psi\in\Lambda^2(V).
$$

A two-form $\varphi\in\Lambda^2(V)$ can be regarded as an element
of the Lie algebra $so(n)$ (i.e. the skew-symmetric matrix
$(\varphi _{ab})_{n\times n}$), where the metric on $\Lambda^2(V)$
is given by
$$
\<\varphi,\psi\>=\varphi_{ab}\psi_{ab}
$$
and the Lie bracket is given by
$$
[\varphi,\psi]_{ab}=\varphi_{ac}\psi_{bc}-\psi_{ac}\varphi_{bc}.
$$

Choose an orthonormal basis of $\Lambda^2(V)$
$$
\Phi=\{\varphi^1,\ldots,\varphi^\alpha,\ldots,\varphi^{\frac{n(n-1)}{2}}\}
$$
where $\varphi^\alpha=\{\varphi_{ab}^\alpha\}$. The Lie bracket is
given by
$$
[\varphi^\alpha,\varphi^\beta]=C^{\alpha\beta}_\gamma\varphi^\gamma,
$$
where
$C^{\alpha\beta\gamma}=C^{\alpha\beta}_\sigma\delta^{\sigma\gamma}
=\<[\varphi^\alpha,\varphi^\beta],\varphi^\gamma\>$ are the Lie
structure constants.

Write
$R_{abcd}=M_{\alpha\beta}\varphi_{ab}^\alpha\varphi^\beta_{cd}.$
We now claim that the first part of the quadratic terms in (1.3.4)
is given by \be
2(B_{abcd}-B_{abdc})=M_{\alpha\gamma}M_{\beta\gamma}
\varphi_{ab}^\alpha\varphi^\beta_{cd}.  
\ee Indeed, by the first Bianchi identity,
\begin{align*}
B_{abcd}-B_{abdc}&  = R_{aebf}R_{cedf}-R_{aebf}R_{decf}\\
&  = R_{aebf}(-R_{cefd}-R_{cfde})\\
&  = R_{aebf}R_{cdef}.
\end{align*}
On the other hand,
\begin{align*}
R_{aebf}R_{cdef}&  = (-R_{abfe}-R_{afeb})R_{cdef}\\
&  = R_{abef}R_{cdef}-R_{afeb}R_{cdef}\\
&  = R_{abef}R_{cdef}-R_{afbe}R_{cdfe}
\end{align*}
which implies $R_{aebf}R_{cdef}=\frac{1}{2}R_{abef}R_{cdef}$. Thus
we obtain
$$
2(B_{abcd}-B_{abdc}) =R_{abef}R_{cdef}
=M_{\alpha\gamma}M_{\beta\gamma}\varphi_{ab}^\alpha\varphi^\beta_{cd}.
$$

We next consider the last part of the quadratic terms:
\begin{align*}
& 2(B_{acbd}-B_{adbc}) \\
&  = 2(R_{aecf}R_{bedf}-R_{aedf}R_{becf})\\
&  = 2(M_{\gamma\delta}\varphi_{ae}^\gamma\varphi_{cf}^\delta
M_{\eta\theta}\varphi_{be}^\eta\varphi^\theta_{df}
-M_{\gamma\theta}\varphi_{ae}^\gamma\varphi_{df}^\theta
M_{\eta\delta}\varphi_{be}^\eta\varphi_{cf}^\delta)\\
& =2[M_{\gamma\delta}(\varphi_{ae}^\eta\varphi_{be}^\gamma
+C_\alpha^{\gamma\eta}\varphi_{ab}^\alpha)\varphi_{cf}^\delta
M_{\eta\theta}\varphi^\theta_{df}
-M_{\eta\theta}\varphi^\eta_{ae}\varphi^\theta_{df}
M_{\gamma\delta}\varphi^\gamma_{be}\varphi^\delta_{cf}]\\
&  = 2M_{\gamma\delta}\varphi_{cf}^\delta
M_{\eta\theta}\varphi^\theta_{df}C_\alpha^{\gamma\eta}\varphi_{ab}^\alpha.
\end{align*}
But
\begin{align*}
& M_{\gamma\delta}\varphi_{cf}^\delta
M_{\eta\theta}\varphi^\theta_{df}C_\alpha^{\gamma\eta}
\varphi_{ab}^\alpha\\
&
=M_{\gamma\delta}M_{\eta\theta}C_\alpha^{\gamma\eta}\varphi_{ab}^\alpha
[\varphi^\theta_{cf}\varphi^\delta_{df}
+C_\beta^{\delta\theta}\varphi^\beta_{cd}]\\
&
=-M_{\eta\theta}M_{\gamma\delta}C_\alpha^{\gamma\eta}\varphi_{ab}^\alpha
\varphi^\delta_{cf}\varphi^\theta_{df}
+M_{\gamma\delta}M_{\eta\theta}C_\alpha^{\gamma\eta}
C_\beta^{\delta\theta}\varphi^\alpha_{ab}\varphi^\beta_{cd}\\
&
=-M_{\eta\theta}M_{\gamma\delta}C_\alpha^{\gamma\eta}\varphi_{ab}^\alpha
\varphi^\delta_{cf}\varphi^\theta_{df}
+(C_\alpha^{\gamma\eta}C_\beta^{\delta\theta}
M_{\gamma\delta}M_{\eta\theta})\varphi^\alpha_{ab}\varphi^\beta_{cd}
\end{align*}
which implies
$$
M_{\gamma\delta}\varphi_{cf}^\delta
M_{\eta\theta}\varphi^\theta_{df}C_\alpha^{\gamma\eta}\varphi_{ab}^\alpha
=\frac{1}{2}(C_\alpha^{\gamma\eta}C_\beta^{\delta\theta}
M_{\gamma\delta}M_{\eta\theta})\varphi^\alpha_{ab}\varphi^\beta_{cd}.
$$
Then we have \be 2(B_{acbd}-B_{adbc})=
(C_\alpha^{\gamma\eta}C_\beta^{\delta\theta}
M_{\gamma\delta}M_{\eta\theta})\varphi^\alpha_{ab}
\varphi^\beta_{cd}. 
\ee

Therefore, combining (1.3.7) and (1.3.8), we can reformulate the
curvature evolution equation (1.3.4) as follows.

\begin{proposition}[{Hamilton \cite{Ha86}}]
Let
$R_{abcd}=M_{\alpha\beta}\varphi_{ab}^\alpha\varphi^\beta_{cd}.$
Then under the Ricci flow $(1.1.5),$ $M_{\alpha\beta}$ satisfies
the evolution equation \be \frac{\partial
M_{\alpha\beta}}{\partial t} =\Delta
M_{\alpha\beta}+M_{\alpha\beta}^2+M_{\alpha\beta}^\# 
\ee where $M_{\alpha\beta}^2=M_{\alpha\gamma}M_{\beta\gamma}$ is
the operator square and
$M_{\alpha\beta}^\#=(C_\alpha^{\gamma\eta}C_\beta^{\delta\theta}
M_{\gamma\delta}M_{\eta\theta})$ is the Lie algebra square.
\end{proposition}

Let us now consider the operator $M_{\alpha\beta}^\#$ in
dimensions 3 and 4 in more detail.

In dimension $3$, let $\omega_1,\omega_2,\omega_3$ be a positively
oriented orthonormal basis for one-forms. Then
$$
\varphi^1=\sqrt{2}\omega_1\wedge\omega_2,\qquad\varphi^2
=\sqrt{2}\omega_2\wedge\omega_3,
\qquad\varphi^3=\sqrt{2}\omega_3\wedge\omega_1
$$
form an orthonormal basis for two-forms $\Lambda^2$. Write
$\varphi ^{\alpha} = \{\varphi ^{\alpha}_{ab}\}, \alpha =1,2,3,$
as
\begin{gather*}
(\varphi^1_{ab})= \begin{pmatrix}
0& \frac{\sqrt{2}}{2}& 0\\
-\frac{\sqrt{2}}{2}& 0&  0\\
0&  0&  0\end{pmatrix},\qquad (\varphi^2_{ab})=\begin{pmatrix} 0&
0& 0\\0& 0& \frac{\sqrt{2}}{2}\\0& -\frac{\sqrt{2}}{2}&
0\end{pmatrix},\\
(\varphi^3_{ab})=\begin{pmatrix} 0& 0& -\frac{\sqrt{2}}{2}\\0& 0&
0\\\frac{\sqrt{2}}{2}& 0& 0\end{pmatrix},
\end{gather*}
then
$$\arraycolsep=1.5pt\begin{array}{rcl}
[\varphi^1,\varphi^2]&  =&  \left(\begin{array}{ccc} 0&
\frac{\sqrt{2}}{2}& 0\\
-\frac{\sqrt{2}}{2}& 0& 0\\0& 0&
0\end{array}\right)\left(\begin{array}{ccc} 0& 0& 0\\0& 0&
-\frac{\sqrt{2}}{2}\\0& \frac{\sqrt{2}}{2}&
0\end{array}\right)-\left(\begin{array}{ccc} 0& 0& 0\\0& 0&
\frac{\sqrt{2}}{2}\\0& -\frac{\sqrt{2}}{2}&
0\end{array}\right)\left(\begin{array}{ccc}
0&  -\frac{\sqrt{2}}{2}&  0\\
\frac{\sqrt{2}}{2}&  0&  0\\
0&  0&  0\end{array}\right)\\[7mm]
&  =&  \left(\begin{array}{ccc}
0&  0&  -\frac{1}{2}\\0&  0&  0\\
\frac{1}{2}&  0&  0\end{array}\right)\\[7mm]
&  =& {\displaystyle \frac{\sqrt{2}}{2}\varphi^3.}
\end{array}
$$
So
$C^{123}=\<[\varphi^1,\varphi^2],\varphi^3\>=\frac{\sqrt{2}}{2}$,
in particular
$$
C^{\alpha\beta\gamma}=\begin{cases} \pm\frac{\sqrt{2}}{2},&
\text{ if }\;
\alpha\neq\beta\neq\gamma,\\
0,&  \text{ otherwise}.
\end{cases}
$$
Hence the matrix $M^\#=(M^{\#}_{\alpha\beta})$ is just the adjoint
matrix of $M=(M_{\alpha\beta})$: \be
M^\#=\det M\cdot {}^tM^{-1}. 
\ee

In dimension 4, we can use the Hodge star operator $\*$ to
decompose the space of two-forms $\Lambda^2$ as
$$
\Lambda^2={\Lambda^2_+}{\oplus}{\Lambda^2_-}
$$
where $\Lambda^2_+$ (resp. $\Lambda^2_-$) is the eigenspace of the
star operator with eigenvalue $+1$ (resp. $-1$). Let $\omega_1,
\omega_2, \omega_3, \omega_4$ be a positively oriented orthonormal
basis for one-forms. A basis for $\Lambda^2_+$ is then given by
$$
\varphi^1=\omega_1\wedge\omega_2 +
\omega_3\wedge\omega_4,\quad\varphi^2 =\omega_1\wedge\omega_3 +
\omega_4\wedge\omega_2, \quad\varphi^3=\omega_1\wedge\omega_4 +
\omega_2\wedge\omega_3,
$$
while a basis for $\Lambda^2_-$ is given by
$$
\psi^1=\omega_1\wedge\omega_2 - \omega_3\wedge\omega_4,\quad\psi^2
=\omega_1\wedge\omega_3 - \omega_4\wedge\omega_2,
\quad\psi^3=\omega_1\wedge\omega_4 - \omega_2\wedge\omega_3.
$$
In particular, $\{\varphi^1, \varphi^2, \varphi^3, \psi^1, \psi^2,
\psi^3 \}$ forms an orthonormal basis for the space of two-forms
$\Lambda^2$. By using this basis we obtain a block decomposition
of the curvature operator matrix $M$ as
$$
M=(M_{\alpha\beta})=\left(\begin{array}{cc} A&  B\\ {}^tB&
C\end{array}\right).
$$
Here $A, B$ and $C$ are $3\times 3$ matrices with $A$ and $C$
being symmetric. Then we can write each element of the basis as a
skew-symmetric $4\times 4$ matrix and compute as above to get \be
M^\#=(M_{\alpha\beta}^\#)=2\left(\begin{array}{cc}
A^\#&  B^\#\\ {}^tB^\#&  C^\#\end{array}\right), 
\ee where $A^\#,B^\#,C^\#$ are the adjoint of $3\times3$
submatrices as before.

For later applications in Chapter 5, we now give some computations
for the entries of the matrices $A$, $C$ and $B$ as follows. First
for the matrices $A$ and $C$, we have
$$
A_{11} = Rm(\varphi^1,\varphi^1) = R_{1212} + R_{3434} + 2R_{1234}
$$
$$
A_{22} = Rm(\varphi^2,\varphi^2) = R_{1313} + R_{4242} + 2R_{1342}
$$
$$
A_{33} = Rm(\varphi^3,\varphi^3) = R_{1414} + R_{2323} + 2R_{1423}
$$
and
$$
C_{11} = Rm(\psi^1,\psi^1) = R_{1212} + R_{3434} -2R_{1234}
$$
$$
C_{22} = Rm(\psi^2,\psi^2) = R_{1313} + R_{4242} -2R_{1342}
$$
$$
C_{33} = Rm(\psi^3,\psi^3) = R_{1414} + R_{2323} -2R_{1423}.
$$
By the Bianchi identity
$$
R_{1234} + R_{1342} + R_{1423} =0,
$$
so we have
$$
trA = trC = \frac{1}{2}R.
$$
Next for the entries of the matrix $B$, we have
$$
B_{11} = Rm(\varphi^1,\psi^1) = R_{1212} - R_{3434}
$$
$$
B_{22} = Rm(\varphi^2,\psi^2) = R_{1313} - R_{4242}
$$
$$
B_{33} = Rm(\varphi^3,\psi^3) = R_{1414} - R_{2323}
$$
and
$$
B_{12} = Rm(\varphi^1,\psi^2) = R_{1213} + R_{3413} -R_{1242} -
R_{3442} \ \ \mbox{ etc. }
$$
Thus the entries of $B$ can be written as
$$
B_{11} = \frac{1}{2}(R_{11} + R_{22} -R_{33} -R_{44})
$$
$$
B_{22} = \frac{1}{2}(R_{11} + R_{33} -R_{44} -R_{22})
$$
$$
B_{33} = \frac{1}{2}(R_{11} + R_{44} -R_{22} -R_{33})
$$
and
$$
B_{12} = R_{23} - R_{14} \ \ \mbox{ etc. }
$$
If we choose the frame $\{ \omega_1, \omega_2, \omega_3, \omega_4
\}$ so that the Ricci tensor is diagonal, then the matrix $B$ is
also diagonal. In particular, the matrix $B$ is identically zero
when the four-manifold is Einstein.

\section{Derivative Estimates}

In the previous section we have seen that the curvatures satisfy
nonlinear heat equations with quadratic growth terms. The
parabolic nature will give us a bound on the derivatives of the
curvatures at any time $t>0$ in terms of a bound of the
curvatures. In this section, we derive Shi's local derivative
estimate \cite{Sh89}. Our presentation follows Hamilton
\cite{Ha95F}.

We begin with the global version of the derivative estimate of Shi
\cite{Sh89}.

\begin{theorem}[{Shi \cite{Sh89}}]
There exist constants $C_m$, $m=1,2,\ldots,$ such that if the
curvature of a complete solution to Ricci flow is bounded by
$$
|R_{ijkl}|\leq M
$$
up to time t with $0<t\leq \frac{1}{M}$, then the covariant
derivative of the curvature is bounded by
$$
|\nabla R_{ijkl}|\leq C_1M/\sqrt{t}
$$
and the $m^{th}$ covariant derivative of the curvature is bounded
by
$$
|\nabla^mR_{ijkl}|\leq C_mM/t^{\frac{m}{2}}.
$$
Here the norms are taken with respect to the evolving metric.
\end{theorem}

\begin{pf}
We shall only give the proof for the compact case.  The noncompact
case can be deduced from the next local derivative estimate
theorem. Let us denote the curvature tensor by $Rm$ and denote by
$A*B$ any tensor product of two tensors $A$ and $B$ when we do not
need the precise expression. We have from Proposition 1.3.1 that
\be
\frac{\partial}{\partial t}Rm=\Delta Rm+Rm*Rm.  
\ee Since
\begin{align*}
\frac{\partial}{\partial t}\Gamma^i_{jk} &  =
\frac{1}{2}g^{il}\left\{\nabla_j\frac{\partial g_{kl}}{\partial t}
+\nabla_k\frac{\partial g_{jl}}{\partial t}
-\nabla_l\frac{\partial g_{jk}}{\partial t}\right\}\\
&  = \nabla Rm ,
\end{align*}
it follows that \be \frac{\partial}{\partial t}(\nabla
Rm)=\Delta(\nabla Rm)
+Rm*(\nabla Rm).  
\ee Thus
\begin{align*}
\frac{\partial}{\partial t}|Rm|^2
&\leq\Delta|Rm|^2-2|\nabla Rm|^2+C|Rm|^3,\\
\frac{\partial}{\partial t}|\nabla Rm|^2 &\leq\Delta|\nabla
Rm|^2-2|\nabla^2 Rm|^2+C|Rm|\cdot|\nabla Rm|^2,
\end{align*}
for some constant $C$ depending only on the dimension $n$.

Let $A>$0 be a constant (to be determined) and set
$$
F=t|\nabla Rm|^2+A|Rm|^2.
$$
We compute
\begin{align*}
\frac{\partial F}{\partial t} &  = |\nabla
Rm|^2+t\frac{\partial}{\partial t}|\nabla Rm|^2
+A\frac{\partial}{\partial t}|Rm|^2\\
& \leq \Delta(t|\nabla Rm|^2+A|Rm|^2)+|\nabla Rm|^2(1+tC|Rm|-2A)
+CA|Rm|^3.
\end{align*}
Taking $A\geq C+1,$ we get
$$
\frac{\partial F}{\partial t}\leq \Delta F+\bar{C}M^3
$$
for some constant $\bar{C}$ depending only on the dimension $n$.
We then obtain
$$
F\leq F(0)+\bar{C}M^3t\leq (A+\bar{C})M^2,
$$
and then
$$
|\nabla Rm|^2\leq(A+\bar{C})M^2/t.
$$

The general case follows in the same way. If we have bounds
$$
|\nabla^k Rm|\leq C_kM/t^\frac{k}{2},
$$
we know from (1.4.1) and (1.4.2) that
$$
\frac{\partial}{\partial t}|\nabla^kRm|^2
\leq\Delta|\nabla^kRm|^2-2|\nabla^{k+1}Rm|^2+\frac{CM^3}{t^k},
$$
and
$$
\frac{\partial}{\partial t}|\nabla^{k+1}Rm|^2
\leq\Delta|\nabla^{k+1}Rm|^2-2|\nabla^{k+2}Rm|^2
+CM|\nabla^{k+1}Rm|^2+\frac{CM^3}{t^{k+1}}.
$$
Let $A_k > 0$ be a constant (to be determined) and set
$$
F_k=t^{k+2}|\nabla^{k+1}Rm|^2+A_k t^{k+1}|\nabla^kRm|^2.
$$
Then
\begin{align*}
\frac{\partial }{\partial t}F_k &  =(k+2)t^{k+1}|\nabla ^{k+1}
Rm|^2
+t^{k+2}\frac{\partial}{\partial t}|\nabla ^{k+1} Rm|^2\\
&\quad  +A_k (k+1) t^k |\nabla ^{k} Rm|^2
+ A_k t^{k+1}\frac{\partial}{\partial t}|\nabla ^{k}Rm|^2\\
& \leq (k+2)t^{k+1}|\nabla ^{k+1} Rm|^2 \\
&\quad+t^{k+2}\left[\Delta|\nabla^{k+1}Rm|^2 -2|\nabla^{k+2}Rm|^2
+CM|\nabla^{k+1}Rm|^2+\frac{CM^3}{t^{k+1}} \right]\\
&\quad  +A_k (k+1) t^k |\nabla ^{k} Rm|^2  \\
&\quad+ A_k t^{k+1}\left[\Delta|\nabla^kRm|^2
-2|\nabla^{k+1}Rm|^2+\frac{CM^3}{t^k} \right]\\
& \leq \Delta F_k + C_{k+1}M^2
\end{align*}
for some positive constant $C_{k+1}$, by choosing $A_k$ large
enough. This implies that
$$
|\nabla^{k+1}Rm|\leq\frac{C_{k+1}M}{t^\frac{k+1}{2}}.
$$
\end{pf}

The above derivative estimate is a somewhat standard Bernstein
estimate in PDEs. By using a cutoff argument, we will derive the
following local version in \cite{Sh89}, which is called
\textbf{Shi's derivative estimate}\index{Shi's derivative
estimate}. The proof is adapted from Hamilton \cite{Ha95F}.

\begin{theorem}[{Shi \cite{Sh89}}]
There exist positive constants $\theta, C_k, k=1, 2, \ldots,$
depending only on the dimension with the following property.
Suppose that the curvature of a solution to the Ricci flow is
bounded
$$
|Rm|\leq M, \ \ \ \ \mbox{on} \
U\times\left[0,\frac{\theta}{M}\right]
$$
where $U$ is an open set of the manifold. Assume that the closed
ball ${B_0(p,r)}$, centered at $p$ of radius $r$ with respect to
the metric at $t=0$, is contained in $U$ and the time $t\leq
{\theta}/{M}$. Then we can estimate the covariant derivatives of
the curvature at $(p,t)$ by
$$
|\nabla Rm(p,t)|^2\leq C_1M^2\(\frac{1}{r^2}+\frac{1}{t}+M\),
$$
and the $k^{th}$ covariant derivative of the curvature at $(p,t)$
by
$$
|\nabla ^k Rm(p,t)|^2\leq
C_kM^2\(\frac{1}{r^{2k}}+\frac{1}{t^k}+M^k\).
$$
\end{theorem}

\begin{pf}
Without loss of generality, we may assume $r\leq\theta/\sqrt{M}$
and the exponential map at $p$ at time $t=0$ is injective on the
ball of radius $r$ (by passing to a local cover if necessary, and
pulling back the local solution of the Ricci flow to the ball of
radius $r$ in the tangent space at $p$ at time $t=0$).

Recall
\begin{align*}
\frac{\partial}{\partial t}|Rm|^2
&\leq\Delta|Rm|^2-2|\nabla Rm|^2+C|Rm|^3, \\
\frac{\partial}{\partial t}|\nabla Rm|^2 &\leq\Delta|\nabla
Rm|^2-2|\nabla^2 Rm|^2+C|Rm|\cdot|\nabla Rm|^2.
\end{align*}
Define
$$
S=(BM^2+|Rm|^2)|\nabla Rm|^2
$$
where $B$ is a positive constant to be determined. By choosing
$B\geq {C^2}/{4}$ and using the Cauchy inequality, we have
\begin{align*}
\frac{\partial}{\partial t}S
&  \leq \Delta S-2BM^2|\nabla^2Rm|^2-2|\nabla Rm|^4 \\
&\quad+CM|\nabla Rm|^2\cdot|\nabla^2Rm|+CBM^3|\nabla Rm|^2\\
&  \leq \Delta S-|\nabla Rm|^4+CB^2M^6\\
&  \leq \Delta S-\frac{S^2}{(B+1)^2M^4}+CB^2M^6.
\end{align*}
If we take
$$
F=b(BM^2+|Rm|^2)|\nabla Rm|^2/M^4=bS/M^4,
$$
and $b\leq \min\{1/(B+1)^2,1/CB^2\}$, we get \be
\frac{\partial F}{\partial t}\leq\Delta F-F^2+M^2.
\ee

We now want to choose a cutoff function $\varphi$ with the support
in the ball $B_0(p,r)$ such that at $t=0$,
$$
\varphi(p)=r,\ \ \ 0\leq\varphi\leq Ar,
$$
and
$$ |\nabla\varphi|\leq A,\ \ \
|\nabla^2\varphi|\leq\frac{A}{r}
$$
for some positive constant $A$ depending only on the dimension.
Indeed, let $g:(-\infty,+\infty)\rightarrow[0,+\infty)$ be a
smooth, nonnegative function satisfying
$$
g(u)= \begin{cases}
         1, & u\in(-\frac{1}{2},\frac{1}{2}),\\
         0, & \text{outside }\; (-1,1).
       \end{cases}
$$
Set
$$
\varphi=rg\(\frac{s^2}{r^2}\),
$$
where $s$ is the geodesic distance function from $p$ with respect
to the metric at $t=0$. Then
$$
\nabla\varphi=\frac{1}{r}g'\(\frac{s^2}{r^2}\)\cdot 2s\nabla s
$$
and hence $$|\nabla\varphi|\leq 2C_1.$$ Also,
$$
\nabla^2\varphi=\frac{1}{r}g''\(\frac{s^2}{r^2}\)\frac{1}{r^2}4s^2\nabla
s\cdot\nabla s+\frac{1}{r}g'\(\frac{s^2}{r^2}\)2\nabla s\cdot
\nabla s+\frac{1}{r}g'\(\frac{s^2}{r^2}\)\cdot 2s\nabla^2 s.
$$
Thus, by using the standard Hessian comparison,
\begin{align*}
|\nabla^2\varphi|&  \leq \frac{C_1}{r}+\frac{C_1}{r}|s\nabla^2s|\\
&  \leq \frac{C_1}{r}\(1+s\(\frac{C_2}{s}+\sqrt{M}\)\)\\
&  \leq \frac{C_3}{r}.
\end{align*}
Here $C_1, C_2$ and $C_3$ are positive constants depending only on
the dimension.

Now extend $\varphi$ to $U\times[0,\frac{\theta}{M}]$ by letting
$\varphi$ to be zero outside $B_0(p,r)$ and independent of time.
Introduce the barrier function \be
H=\frac{(12+4\sqrt{n})A^2}{\varphi^2}+\frac{1}{t}+M 
\ee which is defined and smooth on the set
$\{\varphi>0\}\times(0,T]$.

As the metric evolves, we will still have $0\leq\varphi\leq Ar$
(since $\varphi$ is independent of time $t$); but
$|\nabla\varphi|^2$ and $\varphi|\nabla^2\varphi|$ may increase.
By continuity it will be a while before they double.

\begin{claim}
As long as
$$
|\nabla\varphi|^2\leq 2A^2,\ \ \varphi|\nabla^2\varphi|\leq 2A^2,
$$
we have
$$
\frac{\partial H}{\partial t}>\Delta H-H^2+M^2.$$
\end{claim}

Indeed, by the definition of $H$, we have
\begin{align*}
H^2 &>\frac{(12+4\sqrt{n})^2A^4}{\varphi^4}+\frac{1}{t^2}+M^2,\\
\frac{\partial H}{\partial t}&=-\frac{1}{t^2},
\end{align*}
and
\begin{align*}
\Delta H&  = (12+4\sqrt{n})A^2\Delta \(\frac{1}{\varphi^2}\)\\
&  = (12+4\sqrt{n})A^2\(\frac{6|\nabla\varphi|^2
-2\varphi\Delta\varphi}{\varphi^4}\)\\
&  \leq (12+4\sqrt{n})A^2\(\frac{12A^2+4\sqrt{n}A^2}{\varphi^4}\)\\
&  = \frac{(12+4\sqrt{n})^2A^4}{\varphi^4}.
\end{align*}
Therefore,
$$
H^2>\Delta H-\frac{\partial H}{\partial t}+M^2.
$$

\begin{claim}
If the constant $\theta>0$ is small enough compared to $b$, $B$
and $A$, then we have the following property: as long as $r\leq
{\theta}/{\sqrt{M}}$, $t\leq {\theta}/{M}$ and $F\leq H$, we will
have
$$
|\nabla \varphi|^2\le 2A^2\quad \text{ and }\quad
\varphi|\nabla^2\varphi|\le 2A^2.
$$
\end{claim}

Indeed, by considering the evolution of $\nabla\varphi$, we have
\begin{align*}
\frac{\partial}{\partial t}\nabla_a\varphi
& =\frac{\partial}{\partial t}(F^i_a\nabla_i\varphi)\\
&  = F^i_a\nabla_i\(\frac{\partial \varphi}{\partial t}\)
+\nabla_i\varphi R_k^iF_a^k\\
&  = R_{ab}\nabla_b\varphi
\end{align*}
which implies
$$
\frac{\partial}{\partial t}|\nabla \varphi|^2\le CM|\nabla
\varphi|^2,
$$
and then
$$|\nabla \varphi|^2\le A^2e^{CMt}\le 2A^2,
$$
provided $t\le \theta/M$ with $\theta\le \log 2/C$.

By considering the evolution of $\nabla^2\varphi$, we have
\begin{align*}
\frac{\partial}{\partial t}(\nabla_a\nabla_b\varphi)
& =\frac{\partial}{\partial t}(F_a^iF_b^j\nabla_i\nabla_j\varphi)\\
&  = \frac{\partial}{\partial t}\(F_a^iF_b^j
\(\frac{\partial^2\varphi}{\partial x^i\partial x^j}
-\Gamma_{ij}^k\frac{\partial\varphi}{\partial x^k}\)\)\\
&  = \nabla_a\nabla_b\(\frac{\partial\varphi}{\partial t}\)
+R_{ac}\nabla_b\nabla_c\varphi+R_{bc}\nabla_a\nabla_c\varphi\\
&\quad
+(\nabla_cR_{ab}-\nabla_aR_{bc}-\nabla_bR_{ac})\nabla_c\varphi
\end{align*}
which implies
\be \frac{\partial}{\partial t}|\nabla^2\varphi| \le
C|Rm|\cdot|\nabla^2\varphi|
+C|\nabla Rm|\cdot|\nabla \varphi|.
\ee By assumption $F\le H$, we have \be |\nabla Rm|^2\le
\frac{2M^2}{bB}\(\frac{(12+4\sqrt{n})A^2}{\varphi^2}+\frac{1}{t}\),\quad
\mbox{for}\quad t\le \theta/M. 
\ee Thus by noting $\varphi$ independent of $t$ and $\varphi\le
Ar$, we get from (1.4.5) and (1.4.6) that
$$
\frac{\partial}{\partial t}(\varphi|\nabla^2\varphi|) \le
CM\(\varphi|\nabla^2\varphi|+1+\frac{r}{\sqrt{t}}\)
$$
which implies
\begin{align*}
\varphi|\nabla^2\varphi| &  \le
e^{CMt}\left[(\varphi|\nabla^2\varphi|)|_{t=0}
+CM\int_0^t\(1+\frac{r}{\sqrt{t}}\)dt\right]\\
&  \le e^{CMt}\left[A^2+CM(t+2r\sqrt{t})\right]\\
&  \le 2A^2
\end{align*}
provided $r\le\theta/\sqrt{M}$, and $t\le\theta/M$ with $\theta$
small enough. Therefore we have obtained Claim 2.

The combination of Claim 1 and Claim 2 gives us
$$
\frac{\partial H}{\partial t}>\Delta H-H^2+M^2
$$
as long as $r\le\theta/\sqrt{M},\quad t\le\theta/M$ and $F\le H$.
And (1.4.3) tells us
$$
\frac{\partial F}{\partial t}\le\Delta F-F^2+M^2.
$$
Then the standard maximum principle immediately gives the estimate
$$
|\nabla Rm|^2\le CM^2\(\frac{1}{\varphi^2}+\frac{1}{t}+M\)\qquad
\mbox{on }\;\{\varphi>0\} \times\Big(0,\frac{\theta}{M}\Big],
$$
which implies the first order derivative estimate.

The higher order derivative estimates can be obtained in the same
way by induction. Suppose we have the bounds
$$
|\nabla ^k Rm|^2\leq C_kM^2\(\frac{1}{r^{2k}}+\frac{1}{t^k}+M^k\).
$$
As before, by (1.4.1) and (1.4.2), we have
$$
\frac{\partial}{\partial t}|\nabla^kRm|^2
\leq\Delta|\nabla^kRm|^2-2|\nabla^{k+1}Rm|^2
+CM^3\(\frac{1}{r^{2k}}+\frac{1}{t^k}+M^k\),
$$
and
\begin{align*}
\frac{\partial}{\partial t}|\nabla^{k+1}Rm|^2
&\leq\Delta|\nabla^{k+1}Rm|^2-2|\nabla^{k+2}Rm|^2 \\
&\quad+CM|\nabla^{k+1}Rm|^2
+CM^3\(\frac{1}{r^{2(k+1)}}+\frac{1}{t^{k+1}}+M^{k+1}\).
\end{align*}
Here and in the following we denote by $C$ various positive
constants depending only on $C_k$ and the dimension.

Define
$$
S_k = \left[B_kM^2\(\frac{1}{r^{2k}}+\frac{1}{t^k}+M^k\) +
|\nabla^kRm|^2\right] \cdot |\nabla^{k+1}Rm|^2
$$
where $B_k$ is a positive constant to be determined. By choosing
$B_k$ large enough and Cauchy inequality, we have
{\allowdisplaybreaks
\begin{align*}
\frac{\partial}{\partial t}S_k &  \leq \bigg[-\frac{k}{t^{k+1}}
+\Delta|\nabla^kRm|^2-2|\nabla^{k+1}Rm|^2 \\
&\qquad+CM^3\(\frac{1}{r^{2k}}+\frac{1}{t^k}+M^k\)\bigg]
\cdot|\nabla^{k+1}Rm|^2\\
&\quad+\left[B_kM^2\(\frac{1}{r^{2k}}+\frac{1}{t^k}+M^k\)
+|\nabla^{k}Rm|^2\right] \\
&\qquad\cdot\bigg[\Delta|\nabla^{k+1}Rm|^2
-2|\nabla^{k+2}Rm|^2 +CM|\nabla^{k+1}Rm|^2 \\
&\qquad+CM^3\(\frac{1}{r^{2(k+1)}}+\frac{1}{t^{k+1}}+M^{k+1}\)\bigg]\\
&  \leq \Delta S_k +
8|\nabla^{k}Rm|\cdot|\nabla^{k+1}Rm|^2\cdot|\nabla^{k+2}Rm| -
\frac{k}{t^{k+1}}|\nabla^{k+1}Rm|^2 \\
&\quad- 2|\nabla^{k+1}Rm|^4
+ CM^3|\nabla^{k+1}Rm|^2\(\frac{1}{r^{2k}}+\frac{1}{t^k}+M^k\) \\
&\quad-2|\nabla^{k+2}Rm|^2
\left[B_kM^2\(\frac{1}{r^{2k}}+\frac{1}{t^k}+M^k\)
+|\nabla^kRm|^2\right]\\
&\quad +CM|\nabla^{k+1}Rm|^2\left[B_kM^2\(\frac{1}{r^{2k}}
+\frac{1}{t^k}+M^k\)+|\nabla^kRm|^2\right]\\
&\quad+CM^3\(\frac{1}{r^{2(k+1)}}+\frac{1}{t^{k+1}}
+M^{k+1}\)\\
&\qquad\cdot\left[B_kM^2\(\frac{1}{r^{2k}}+\frac{1}{t^k}+M^k\)
+ |\nabla^kRm|^2\right]\\
& \leq \Delta S_k - |\nabla^{k+1}Rm|^4 +
CB_k^2M^6\(\frac{1}{r^{4k}}+\frac{1}{t^{2k}}+M^{2k}\) \\
&\quad +CB_kM^5\(\frac{1}{r^{2(2k+1)}}+\frac{1}{t^{2k+1}}+M^{2k+1}\)\\
&  \leq \Delta S_k - |\nabla^{k+1}Rm|^4 +
CB_k^2M^5\(\frac{1}{r^{2(2k+1)}}+\frac{1}{t^{2k+1}}+M^{2k+1}\)\\
&  \leq \Delta S_k -
\frac{S_k}{(B+1)^2M^4\(\frac{1}{r^{2k}}+\frac{1}{t^k}+M^k\)^2} \\
&\quad
+CB_k^2M^5\(\frac{1}{r^{2(2k+1)}}+\frac{1}{t^{2k+1}}+M^{2k+1}\).
\end{align*}
} Let $u = {1}/{r^2} + {1}/{t} +M $ and set $F_k = bS_k/u^k$. Then
\begin{align*}
\frac{\partial F_k}{\partial t} &  \leq \Delta F_k
-\frac{F_k^2}{b(B_k+1)^2M^4u^k}
+ bCB_k^2M^5u^{k+1} + kF_ku\\
&  \leq \Delta F_k -\frac{F_k^2}{2b(B_k+1)^2M^4u^k} +
b(C+2k^2)(B_k+1)^2M^4u^{k+2}.
\end{align*}
By choosing $b \leq 1/(2(C+2k^2)(B_k+1)^2M^4)$, we get
$$
\frac{\partial F_k}{\partial t}\leq\Delta F_k - \frac{1}{u^k}F_k^2
+ u^{k+2}.
$$
Introduce
$$
H_k = 5(k+1)(2(k+1)+1+\sqrt{n})A^2\varphi^{-2(k+1)} + Lt^{-(k+1)}
+ M^{k+1},
$$
where $L \geq k+2$. Then by using Claim 1 and Claim 2, we have
$$
\frac{\partial H_k}{\partial t} = -(k+1)Lt^{-(k+2)},
$$
$$
\Delta H_k \leq 20(k+1)^2(2(k+1)+1+\sqrt{n})A^4\varphi^{-2(k+2)}
$$
and
$$
H_k^2 > 25(k+1)^2(2(k+1)+1+\sqrt{n})A^4\varphi^{-4(k+1)} +
L^2t^{-2(k+1)} + M^{2(k+1)}.
$$
These imply
$$
\frac{\partial H_k}{\partial t} >\Delta H_k - \frac{1}{u^k}H_k^2 +
u^{k+2}.
$$
Then the maximum principle immediately gives the estimate
$$
F_k \leq H_k.
$$
In particular,
\begin{align*}
&\frac{b}{u^k}B_kM^2\(\frac{1}{r^{2k}}
+\frac{1}{t^k}+M^k\)\cdot|\nabla^{k+1}Rm|^2 \\
& \leq 5(k+1)\(2(k+1)+1+\sqrt{n}\)A^2\varphi^{-2(k+1)} +
Lt^{-(k+1)} + M^{k+1}.
\end{align*}
So by the definition of $u$ and the choosing of $b$, we obtain the
desired estimate
$$
|\nabla ^{k+1} Rm|^2 \leq
C_{k+1}M^2\(\frac{1}{r^{2(k+1)}}+\frac{1}{t^{k+1}}+M^{k+1}\).
$$

Therefore we have completed the proof of the theorem.
\end{pf}

\section{Variational Structure and Dynamic Property}

In this section, we introduce two functionals of Perelman
\cite{P1}, $\mathcal F$ and $\mathcal W$, and discuss their
relations with the Ricci flow. Our presentation here follows
sections 1.1-1.3 of Perelman \cite{P1}.

It was not known whether the Ricci flow is a gradient flow until
Perelman \cite{P1} showed that the Ricci flow is, in a certain
sense, the gradient flow of the functional $\mathcal F$. If we
consider the Ricci flow as a dynamical system on the space of
Riemannian metrics, then these two functionals are of Lyapunov
type for this dynamical system. Obviously, Ricci flat metrics are
fixed points of the dynamical system. When we consider the space
of Riemannian metrics modulo diffeomorphism and scaling, fixed
points of the Ricci flow dynamical system correspond to steady, or
shrinking, or expanding Ricci solitons. The following concept
corresponds to a periodic orbit.

\begin{definition}
A metric $g_{ij}(t)$ evolving by the Ricci flow is called a {\bf
breather}\index{breather} if for {\bf some} $t_1<t_2$ and
$\alpha>0$ the metrics $\alpha g_{ij}(t_1)$ and $g_{ij}(t_2)$
differ only by a diffeomorphism; the case
$\alpha=1,\;\alpha<1,\;\alpha>1$ correspond to {\bf steady,
shrinking and expanding
breathers}\index{breather!steady}\index{breather!shrinking}\index{breather!expanding},
respectively.
\end{definition}

Clearly,\; (steady,\; shrinking\; or\; expanding)\; Ricci\;
solitons\; are\; trivial breathers for which the metrics
$g_{ij}(t_1)$ and $g_{ij}(t_2)$ differ only by diffeomorphism and
scaling for every pair $t_1$ and $t_2$.

We always assume $M$ is a compact $n$-dimensional manifold in this
section. Let us first consider the functional \be
{\mathcal{F}} (g_{ij},f)=\int_M(R+|\nabla f|^2)e^{-f}dV 
\ee of Perelman \cite{P1} defined on the space of Riemannian
metrics, and smooth functions on $M$. Here $R$ is the scalar
curvature of $g_{ij}$.

\begin{lemma}[{Perelman \cite{P1}}]
If $\delta g_{ij}=v_{ij}$ and $\delta f=h$ are variations of
$g_{ij}$ and $f$ respectively, then the first variation of
${\mathcal{F}}$ is given by
$$
\delta{\mathcal{F}}(v_{ij},h) =\!\!\int_M\left[-v_{ij}
(R_{ij}+\nabla_i\nabla_jf)+\(\frac{v}{2}-h\)(2\Delta f-|\nabla
f|^2+R)\right]e^{-f}dV
$$
where $v=g^{ij}v_{ij}$.
\end{lemma}

\begin{pf}
In any normal coordinates at a fixed point, we have
\begin{align*}
\delta R^h_{ijl}&  = \frac{\partial}{\partial
x^i}(\delta\Gamma^h_{jl})-\frac{\partial}{\partial
x^j}(\delta\Gamma^h_{il})\\
&  = \frac{\partial}{\partial x^i}\left[\frac{1}{2}g^{hm}(\nabla_j
v_{lm}+\nabla_lv_{jm}-\nabla_mv_{jl})\right]\\
&\quad -\frac{\partial}{\partial
x^j}\left[\frac{1}{2}g^{hm}(\nabla_i
v_{lm}+\nabla_lv_{im}-\nabla_mv_{il})\right],\\
\delta R_{jl}&  = \frac{\partial}{\partial
x^i}\left[\frac{1}{2}g^{im}(\nabla_j
v_{lm}+\nabla_lv_{jm}-\nabla_mv_{jl})\right]\\
&\quad -\frac{\partial}{\partial
x^j}\left[\frac{1}{2}g^{im}(\nabla_i
v_{lm}+\nabla_lv_{im}-\nabla_mv_{il})\right]\\
&  = \frac{1}{2}\frac{\partial}{\partial
x^i}[\nabla_jv^i_l+\nabla_lv^i_j-\nabla^iv_{jl}]
-\frac{1}{2}\frac{\partial}{\partial x^j}[\nabla_lv],\\
\delta R&  = \delta(g^{jl}R_{jl})\\
&  = -v_{jl}R_{jl}+g^{jl}\delta R_{jl}\\
&  = -v_{jl}R_{jl}+\frac{1}{2}\frac{\partial}{\partial x^i}
[\nabla^lv^i_l+\nabla_lv^{il}-\nabla^iv]
-\frac{1}{2}\frac{\partial}{\partial x^j}[\nabla^jv]\\
&  = -v_{jl}R_{jl}+\nabla_i\nabla_lv_{il}-\Delta v.
\end{align*}
Thus \be \delta R(v_{ij})=-\Delta
v+\nabla_i\nabla_jv_{ij}-v_{ij}R_{ij}.  
\ee The first variation of the functional
${\mathcal{F}}(g_{ij},f)$ is
\begin{align}
& \delta\(\int_M(R+|\nabla f|^2)e^{-f}dV\)\\
&  = \int_M\bigg([\delta
R(v_{ij})+\delta(g^{ij}\nabla_if\nabla_jf)]e^{-f}dV \nn\\
&\qquad+(R+|\nabla f|^2)
\left[-he^{-f}dV+e^{-f}\frac{v}{2}dV\right]\bigg) \nn\\
&  = \int_M\bigg[-\Delta
v+\nabla_i\nabla_jv_{ij}-R_{ij}v_{ij}-v_{ij}\nabla_if\nabla_jf \nn\\
&\qquad +2\<\nabla f,\nabla h\>
+(R+|\nabla f|^2)\(\frac{v}{2}-h\)\bigg]e^{-f}dV.\nn  
\end{align}

On the other hand,
\begin{align*}
\int_M(\nabla_i\nabla_jv_{ij}-v_{ij}\nabla_if\nabla_jf)e^{-f}dV &=
\int_M(\nabla_if\nabla_jv_{ij}
-v_{ij}\nabla_if\nabla_jf)e^{-f}dV\\
&  = -\int_M(\nabla_i\nabla_jf)v_{ij}e^{-f}dV,\\
\int_M\!2\<\nabla f,\nabla h\>e^{-f}dV &  = -2\!\int_M\!h\Delta
fe^{-f}dV+2\!\int_M\!|\nabla f|^2he^{-f}dV,
\end{align*}
and
\begin{align*}
\int_M(-\Delta v)e^{-f}dV
& = -\int_M\<\nabla f,\nabla v\>e^{-f}dV\\
& = \int_Mv\Delta fe^{-f}dV-\int_M|\nabla f|^2ve^{-f}dV.
\end{align*}
Plugging these identities into (1.5.3) the first variation formula
follows.
\end{pf}

Now let us study the functional ${\mathcal{F}}$ when the metric
evolves under the Ricci flow and the function evolves by a
backward heat equation.

\begin{proposition}[{Perelman \cite{P1}}]
Let $g_{ij}(t)$ and $f(t)$ evolve according to the coupled flow
$$\arraycolsep=1.5pt\left\{\begin{array}{rcl}
\frac{\partial g_{ij}}{\partial t}&  =&  -2R_{ij},\\[2mm]
\frac{\partial f}{\partial t}&  =&  -\Delta f+|\nabla f|^2-R.
\end{array}\right.
$$
Then
$$
\frac{d}{dt}{\mathcal{F}}(g_{ij}(t),f(t))
=2\int_M|R_{ij}+\nabla_i\nabla_j f|^2e^{-f}dV
$$
and $\int_Me^{-f}dV$ is constant. In particular
${\mathcal{F}}(g_{ij}(t),f(t))$ is nondecreasing in time and the
monotonicity is strict unless we are on a steady gradient soliton.
\end{proposition}

\begin{pf}
Under the coupled flow and using the first variation formula in
Lemma 1.5.2, we have
\begin{align*}
&\frac{d}{dt}{\mathcal{F}}(g_{ij}(t),f(t)) \\
& = \int_M\bigg[-(-2R_{ij})(R_{ij}+\nabla_i\nabla_jf)\\
&\qquad +\(\frac{1}{2}(-2R) -\frac{\partial f}{\partial
t}\)(2\Delta
f-|\nabla f|^2+R)\bigg]e^{-f}dV\\
&  = \int_M[2R_{ij}(R_{ij}+\nabla_i\nabla_jf)+(\Delta f-|\nabla
f|^2)(2\Delta f-|\nabla f|^2+R)]e^{-f}dV.
\end{align*}
Now
\begin{align*}
& \int_M(\Delta f-|\nabla f|^2)(2\Delta f-|\nabla
f|^2)e^{-f}dV\\
&  = \int_M-\nabla_if\nabla_i(2\Delta f-|\nabla f|^2)e^{-f}dV\\
& =\int_M-\nabla_if(2\nabla_j(\nabla_i\nabla_jf)-2R_{ij}\nabla_jf
-2\<\nabla f,\nabla_i\nabla f\>)e^{-f}dV\\
& =\!-2\!\int_M[(\nabla_if\nabla_jf\!-\!\nabla_i\nabla_jf)
\nabla_i\nabla_jf\!-\!R_{ij}\nabla_if\nabla_jf\!-\!\<\nabla
f,\nabla_i\nabla f\>\nabla_if]e^{-f}dV\\
&  = 2\int_M[|\nabla_i\nabla_jf|^2+R_{ij}\nabla_i
f\nabla_jf]e^{-f}dV,
\end{align*}
and
\begin{align*}
&\int_M(\Delta f-|\nabla f|^2)Re^{-f}dV \\
&  = \int_M-\nabla_if\nabla_iRe^{-f}dV\\
& =2\int_M\nabla_i\nabla_jfR_{ij}e^{-f}dV
-2\int_M\nabla_if\nabla_jfR_{ij}e^{-f}dV.
\end{align*}
Here we have used the contracted second Bianchi identity.
Therefore we obtain
\begin{align*}
& \frac{d}{dt}{\mathcal{F}}(g_{ij}(t), f(t)) \\
& = \int_M[2R_{ij}(R_{ij}+\nabla_i\nabla_jf)
+2(\nabla_i\nabla_jf)(\nabla_i\nabla_jf+R_{ij})]e^{-f}dV\\
& = 2\int_M|R_{ij}+\nabla_i\nabla_jf|^2e^{-f}dV.
\end{align*}
It remains to show $\int_Me^{-f}dV$ is a constant. Note that the
volume element $dV=\sqrt{{\det}g_{ij}}\ dx$ evolves under the
Ricci flow by
\begin{align}
\frac{\partial}{\partial t}dV&  = \frac{\partial}{\partial
t}(\sqrt{\det g_{ij}}) dx\\
&  = \frac {1}{2}\(\frac{\partial }{\partial t}
\log(\det g_{ij})\) dV \nn\\
&  =\frac{1}{2}(g^{ij}\frac{\partial}{\partial t}g_{ij})dV \nn\\
&  = -RdV.  \nn  
\end{align}
Hence
\begin{align}
\frac{\partial}{\partial t}\(e^{-f} dV\)
& = e^{-f}\(-\frac{\partial f}{\partial t}-R\)dV \\
&  = (\Delta f-|\nabla f|^2)e^{-f}dV \nn\\
&  = -\Delta(e^{-f})dV. \nn
\end{align}
It then follows that
$$
\frac{d}{dt}\int_Me^{-f}dV=-\int_M\Delta(e^{-f})dV=0.
$$
This finishes the proof of the proposition.
\end{pf}

Next we define the associated energy \be \lambda(g_{ij})
=\inf\left\{{\mathcal{F}}(g_{ij},f)\ |\  f\in C^{\infty}(M),
\int_Me^{-f}dV=1\right\}. 
\ee

If we set $u=e^{-f/2}$, then the functional ${\mathcal{F}}$ can be
expressed in terms of $u$ as
$$
{\mathcal{F}}=\int_M(R u^2+4|\nabla u|^2)dV,
$$
and the constraint $\int_Me^{-f}dV=1$ becomes $\int_M u^2 dV=1$.
Therefore $\lambda (g_{ij})$ is just the first eigenvalue of the
operator $-4\Delta +R$.  Let $u_0>0$ be a first eigenfunction of
the operator $-4\Delta +R$ satisfying
$$
-4\Delta u_0+Ru_0=\lambda(g_{ij})u_0.
$$
The $f_0=-2\log u_0$ is a minimizer:
$$
\lambda(g_{ij})={\mathcal{F}}(g_{ij},f_0).
$$
Note that $f_0$ satisfies the equation \be
-2\Delta f_0+|\nabla f_0|^2-R=-\lambda(g_{ij}). 
\ee

Observe that the evolution equation
$$
\frac{\partial f}{\partial t}=-\Delta f+|\nabla f|^2-R
$$
can be rewritten as the following linear equation
$$
\frac{\partial}{\partial t}(e^{-f})=-\Delta(e^{-f})+R(e^{-f}).
$$
Thus we can always solve the evolution equation for $f$ backwards
in time. Suppose at $t=t_0$, the infimum $\lambda(g_{ij})$ is
achieved by some function $f_0$ with $\int_M e^{-f_0}dV=1$. We
solve the backward heat equation
$$
\arraycolsep=1.5pt\left\{\begin{array}{l}
\frac{\partial f}{\partial t}=-\Delta f+|\nabla f|^2-R\\[2mm]
f|_{t=t_0}=f_0
\end{array}\right.
$$
to obtain a solution $f(t)$ for $t\le t_0$ which satisfies
$\int_Me^{-f}dV=1$. It then follows from Proposition 1.5.3 that
$$
\lambda(g_{ij}(t))\le{\mathcal{F}}(g_{ij}(t),f(t))\le{\mathcal{F}}
(g_{ij}(t_0),f(t_0))=\lambda(g_{ij}(t_0)).
$$
Also note $\lambda(g_{ij})$ is invariant under diffeomorphism.
Thus we have proved

\begin{corollary} [{Perelman \cite{P1}}]\
\begin{itemize}
\item[{\rm (i)}] $\lambda(g_{ij}(t))$ is nondecreasing along the
Ricci flow and the monotonicity is strict unless we are on a
steady gradient soliton; \item[{\rm (ii)}] A steady breather is
necessarily a steady gradient soliton.
\end{itemize}
\end{corollary}

To deal with the expanding case we consider a scale invariant
version
$$
\bar{\lambda}(g_{ij})=\lambda(g_{ij})V^{\frac{2}{n}}(g_{ij}).
$$
Here $V=Vol(g_{ij})$ denotes the volume of $M$ with respect to the
metric $g_{ij}$. \pagebreak

\begin{corollary} [{Perelman \cite{P1}}]\
\begin{itemize}
\item[{\rm (i)}] $\bar{\lambda}(g_{ij})$ is nondecreasing along
the Ricci flow whenever it is nonpositive; moreover, the
monotonicity is strict unless we are on a gradient expanding
soliton; \item[{\rm (ii)}] An expanding breather is necessarily an
expanding gradient soliton.
\end{itemize}
\end{corollary}

\begin{pf}
Let $f_0$ be a minimizer of $\lambda(g_{ij}(t))$ at $t=t_0$ and
solve the backward heat equation
$$
\frac{\partial f}{\partial t}=-\Delta f+|\nabla f|^2-R
$$
to obtain $f(t)$, $t\le t_0$, with $\int_Me^{-f(t)}dV=1$. We
compute the derivative (understood in the barrier sense) at
$t=t_0$,
\begin{align*}
& \frac{d}{dt}\bar{\lambda}(g_{ij}(t)) \\
&  \ge\frac{d}{dt}({\mathcal{F}}(g_{ij}(t),f(t)) \cdot
V^{\frac{2}{n}}(g_{ij}(t)))\\
&  = V^{\frac{2}{n}}\int_M2|R_{ij}+\nabla_i\nabla_jf|^2e^{-f}dV\\
&\quad
+\frac{2}{n}V^{\frac{2-n}{n}}\int_M(-R)dV\cdot\int_M(R+|\nabla
f|^2)e^{-f}dV\\
& =2V^{\frac{2}{n}}\bigg[\int_M\bigg|R_{ij}
+\nabla_i\nabla_jf-\frac{1}{n}(R+\Delta f)g_{ij}\bigg|^2e^{-f}dV\\
&\quad +\frac{1}{n}\int_M\!(R\!+\!\Delta f)^2e^{-f}dV\!
+\!\frac{1}{n}\(\!-\!\int_M\!(R\!+\!|\nabla f|^2)e^{-f}dV\)
\(\frac{1}{V}\int_MRdV\)\bigg],
\end{align*}
where we have used the formula (1.5.4) in the computation of
$dV/dt$.

Suppose $\lambda(g_{ij}(t_0))\le 0$, then the last term on the RHS
is given by,
\begin{align*}
& \frac{1}{n}\(-\int_M(R+|\nabla
f|^2\)e^{-f}dV)\(\frac{1}{V}\int_MRdV\)\\
&  \ge \frac{1}{n}\(-\int_M(R+|\nabla f|^2)e^{-f}dV\)
\(\int_M(R+|\nabla f|^2)e^{-f}dV\)\\
&  = -\frac{1}{n}\(\int_M(R+\Delta f)e^{-f}dV\)^2.
\end{align*}
Thus at $t=t_0$,
\begin{align}
& \frac{d}{dt}\bar{\lambda}(g_{ij}(t)) \\
&  \ge 2V^{\frac{2}{n}}
\bigg[\int_M|R_{ij}+\nabla_i\nabla_jf-\frac{1}{n}(R+\Delta
f)g_{ij}|^2e^{-f}dV \nn\\
&\quad +\frac{1}{n}\(\int_M(R+\Delta f)^2e^{-f}dV
-\(\int_M(R+\Delta f)e^{-f}dV\)^2\)\bigg]\ge 0 \nn  
\end{align}
by the Cauchy-Schwarz inequality. Thus we have proved statement
(i).

We note that on an expanding breather on $[t_1,t_2]$ with $\alpha
g_{ij}(t_1)$ and $g_{ij}(t_2)$ differ only by a diffeomorphism for
some $\alpha>1$, it would necessary have
$$
\frac{dV}{dt}>0,\; \text{ for some }\;t\in[t_1,t_2].
$$
On the other hand, for every $t$,
$$
-\frac{d}{dt}\log V=\frac{1}{V}\int_MRdV\ge \lambda(g_{ij}(t))
$$
by the definition of $\lambda(g_{ij}(t))$. It follows that on an
expanding breather on $[t_1,t_2]$,
$$
\bar{\lambda}(g_{ij}(t))
=\lambda(g_{ij}(t))V^\frac{2}{n}(g_{ij}(t))<0
$$
for some $t\in[t_1,t_2]$. Then by using statement (i), it implies
$$
\bar{\lambda}(g_{ij}(t_1))<\bar{\lambda}(g_{ij}(t_2))
$$ unless we
are on an expanding gradient soliton. We also note that
$\bar{\lambda}(g_{ij}(t))$ is invariant under diffeomorphism and
scaling which implies
$$
\bar{\lambda}(g_{ij}(t_1))=\bar{\lambda}(g_{ij}(t_2)).
$$
Therefore the breather must be an expanding gradient soliton.
\end{pf}

In particular part (ii) of Corollaries 1.5.4 and 1.5.5  imply that
all compact steady or expanding Ricci solitons are gradient ones.
Combining this fact with Proposition (1.1.1), we immediately get

\begin{proposition} [{Perelman \cite{P1}}]
On a compact manifold, a steady or expanding breather is necessarily
an Einstein metric.
\end{proposition}

In order to handle the shrinking case, we introduce the following
important functional, also due to Perelman \cite{P1}, \be
\mathcal{W}(g_{ij},f,\tau) =\int_M [\tau (R+|\nabla
f|^2)+f-n](4\pi
\tau)^{-\frac{n}{2}}e^{-f}dV 
\ee where $g_{ij}$ is a Riemannian metric, $f$ is a smooth
function on $M$, and $\tau$ is a positive scale parameter. Clearly
the functional $\mathcal{W}$ is invariant under simultaneous
scaling of $\tau$ and $g_{ij}$ (or equivalently the parabolic
scaling), and invariant under diffeomorphism. Namely, for any
positive number $a$ and any diffeomorphism $\varphi$ \be
\mathcal{W}(a\varphi^*g_{ij},\varphi^*f,a
\tau)=\mathcal{W}(g_{ij},f,\tau).
\ee

Similar to Lemma 1.5.2, we have the following first variation
formula for $\mathcal{W}$.

\begin{lemma}[{Perelman \cite{P1}}]
If $v_{ij}=\delta g_{ij},\; h=\delta f,\;\mbox{ and }\;
\eta=\delta\tau$, then
\begin{align*}
& \delta \mathcal{W}(v_{ij},h,\eta)\\
&  = \int_M-\tau v_{ij}\(R_{ij}+\nabla_i\nabla_jf
-\frac{1}{2\tau}g_{ij}\)(4\pi\tau)^{-\frac{n}{2}}e^{-f}dV\\
&\quad +\int_M\(\frac{v}{2}-h-\frac{n}{2\tau}\eta\)[\tau(R+2\Delta
f
-|\nabla f|^2)+f-n-1](4\pi\tau)^{-\frac{n}{2}}e^{-f}dV\\
&\quad +\int_M \eta\(R+|\nabla f|^2
-\frac{n}{2\tau}\)(4\pi\tau)^{-\frac{n}{2}}e^{-f}dV.
\end{align*}
Here $v=g^{ij}v_{ij}$ as before.
\end{lemma}

\begin{pf}
Arguing as in the proof of Lemma 1.5.2, the first variation of the
functional $\mathcal{W}$ can be computed as follows,
\begin{align*}
&\delta W(v_{ij},h,\eta)  \\
&  = \int_M[\eta(R+|\nabla f|^2)+\tau(-\Delta
v+\nabla_i\nabla_jv_{ij}-R_{ij}v_{ij}-v_{ij}\nabla_if\nabla_jf \\
&\quad +2\<\nabla f,\nabla h\>)+h](4\pi\tau)^{-\frac{n}{2}}e^{-f}dV\\
&\quad +\int_M\left[(\tau(R+|\nabla
f|^2)+f-n)\(-\frac{n}{2}\frac{\eta}{\tau}+\frac{v}{2}-h\)\right]
(4\pi\tau)^{-\frac{n}{2}}e^{-f}dV\\
& = \int_M[\eta(R+|\nabla f|^2)+h](4\pi\tau)^{-\frac{n}{2}}e^{-f}dV\\
&\quad +\int_M[-\tau v_{ij}(R_{ij}+\nabla_i\nabla_j f)
+\tau(v-2h)(\Delta f-|\nabla f|^2)](4\pi\tau)^{-\frac{n}{2}}e^{-f}dV\\
&\quad +\int_M\left[(\tau(R+|\nabla
f|^2)+f-n)\(-\frac{n}{2}\frac{\eta}{\tau}+\frac{v}{2}-h\)\right]
(4\pi\tau)^{-\frac{n}{2}}e^{-f}dV\\
& = -\int_M\tau
v_{ij}\(R_{ij}+\nabla_i\nabla_jf-\frac{1}{2\tau}g_{ij}\)
(4\pi\tau)^{-\frac{n}{2}}e^{-f}dV\\
&\quad +\int_M\(\frac{v}{2}-h-\frac{n}{2\tau}\eta\)
[\tau(R+|\nabla f|^2) \\
&\qquad+f-n+2\tau(\Delta f-|\nabla f|^2)]
(4\pi\tau)^{-\frac{n}{2}}e^{-f}dV\\
&\quad +\int_M\left[\eta\(R+|\nabla
f|^2-\frac{n}{2\tau}\)+\(h-\frac{v}{2}+\frac{n}{2\tau}\eta\)\right]
(4\pi\tau)^{-\frac{n}{2}}e^{-f}dV\\
& = \int_M-\tau
v_{ij}\(R_{ij}+\nabla_i\nabla_jf-\frac{1}{2\tau}g_{ij}\)
(4\pi\tau)^{-\frac{n}{2}}e^{-f}dV\\
&\quad +\int_M\(\frac{v}{2}-h-\frac{n}{2\tau}\eta\)[\tau(R+2\Delta
f
-|\nabla f|^2)+f-n-1](4\pi\tau)^{-\frac{n}{2}}e^{-f}dV\\
&\quad +\int_M \eta\(R+|\nabla f|^2
-\frac{n}{2\tau}\)(4\pi\tau)^{-\frac{n}{2}}e^{-f}dV.
\end{align*}
\end{pf}

The following result is analogous to Proposition 1.5.3.

\begin{proposition}[{Perelman \cite{P1}}]
If $g_{ij}(t), f(t)$ and $\tau(t)$ evolve according to the system
$$
      \left\{
       \begin{array}{lll}
       \displaystyle
\frac{\partial g_{ij}}{\partial t}=-2R_{ij},
          \\[4mm]
      \displaystyle
\frac{\partial f}{\partial t}=-\Delta f+|\nabla
f|^2-R+\frac{n}{2\tau},
          \\[4mm]
      \displaystyle
\frac{\partial \tau}{\partial t}=-1,
       \end{array}
      \right.
$$
then we have the identity
$$
\frac{d} {d t} \mathcal{W}(g_{ij}(t),f(t),\tau(t))=\int_M
2\tau\left|R_{ij}+\nabla_i\nabla_jf-\frac
{1}{2\tau}g_{ij}\right|^2 (4\pi\tau)^{-\frac{n}{2}} e^{-f}dV
$$
and $\int_M(4\pi\tau)^{-\frac{n}{2}}e^{-f}dV$ is constant. In
particular $\mathcal{W}(g_{ij}(t),f(t),\tau(t))$ is nondecreasing
in time and the monotonicity is strict unless we are on a
shrinking gradient soliton.
\end{proposition}

\begin{pf}
Using Lemma 1.5.7, we have
\begin{align}
&\frac{d} {d t} \mathcal{W}(g_{ij}(t),f(t),\tau(t)) \\
&  = \int_M 2\tau R_{ij}\(R_{ij}+\nabla_i\nabla_jf
-\frac{1}{2\tau}g_{ij}\)(4\pi\tau)^{-\frac{n}{2}}e^{-f}dV \nn\\
&\quad +\int_M(\Delta f-|\nabla f|^2)[\tau(R+2\Delta f-|\nabla
f|^2)+f](4\pi\tau)^{-\frac{n}{2}}e^{-f}dV \nn\\
&\quad -\int_M\(R+|\nabla f|^2-\frac{n}{2\tau}\)
(4\pi\tau)^{-\frac{n}{2}}e^{-f}dV. \nn  
\end{align}
Here we have used the fact that $\int_M(\Delta f-|\nabla
f|^2)e^{-f}dV=0.$

The second term on the RHS of (1.5.11) is {\allowdisplaybreaks
\begin{align*}
& \int_M(\Delta f-|\nabla f|^2)[\tau(R+2\Delta f-|\nabla f|^2)+f]
(4\pi\tau)^{-\frac{n}{2}}e^{-f}dV\\
& = \int_M(\Delta f-|\nabla f|^2)(2\tau\Delta f-\tau|\nabla
f|^2)(4\pi\tau)^{-\frac{n}{2}}e^{-f}dV \\
&\quad-\int_M|\nabla f|^2(4\pi\tau)^{-\frac{n}{2}}e^{-f}dV
 +\tau\int_M(-\nabla_if)(\nabla_iR)
(4\pi\tau)^{-\frac{n}{2}}e^{-f}dV\\
& = \tau\int_M(-\nabla_if)(\nabla_i(2\Delta f-|\nabla
f|^2))(4\pi\tau)^{-\frac{n}{2}}e^{-f}dV \\
&\quad-\int_M\Delta f(4\pi\tau)^{-\frac{n}{2}}e^{-f}dV
-2\tau\int_M\nabla_if\nabla_jR_{ij}(4\pi\tau)^{-\frac{n}{2}}
e^{-f}dV\\
& = -2\tau \int_M(\nabla_if)(\nabla_i\Delta f-\<\nabla
f,\nabla_i\nabla f\>)(4\pi\tau)^{-\frac{n}{2}}e^{-f}dV \\
&\quad+2\tau\int_M[(\nabla_i\nabla_jf)R_{ij}
-\nabla_if\nabla_jfR_{ij}](4\pi\tau)^{-\frac{n}{2}}e^{-f}dV\\
&\quad+2\tau\int_M\(-\frac{1}{2\tau}g_{ij}\)(\nabla_i\nabla_jf)
(4\pi\tau)^{-\frac{n}{2}}e^{-f}dV\\
& =-2\tau\int_M[(\nabla_if\nabla_jf-\nabla_i\nabla_jf)
\nabla_i\nabla_jf-R_{ij}\nabla_if\nabla_jf\\
&\qquad-\nabla_i\nabla_j
f\nabla_if\nabla_jf](4\pi\tau)^{-\frac{n}{2}}e^{-f}dV\\
&\quad +2\tau\int_M[(\nabla_i\nabla_jf)R_{ij}
-\nabla_if\nabla_jfR_{ij}](4\pi\tau)^{-\frac{n}{2}}e^{-f}dV\\
&\quad +2\tau\int_M\(-\frac{1}{2\tau}g_{ij}\)(\nabla_i\nabla_jf)
(4\pi\tau)^{-\frac{n}{2}}e^{-f}dV\\
& =2\tau\int_M(\nabla_i\nabla_jf)\(\nabla_i\nabla_jf+R_{ij}
-\frac{1}{2\tau}g_{ij}\)(4\pi\tau)^{-\frac{n}{2}}e^{-f}dV.
\end{align*}
} Also the third term on the RHS of (1.5.11) is
\begin{align*}
& \int_M-\(R+|\nabla f|^2-\frac{n}{2\tau}\)
(4\pi\tau)^{-\frac{n}{2}}e^{-f}dV\\
&  = \int_M-\(R+\Delta f-\frac{n}{2\tau}\)
(4\pi\tau)^{-\frac{n}{2}}e^{-f}dV\\
& =2\tau\int_M\(\frac{-1}{2\tau}g_{ij}\)\(R_{ij}
+\nabla_i\nabla_jf-\frac{1}{2\tau}g_{ij}\)
(4\pi\tau)^{-\frac{n}{2}}e^{-f}dV.
\end{align*}
Therefore, by combining the above identities, we obtain
$$
\frac{d} {d t} \mathcal{W}(g_{ij}(t),f(t),\tau(t))
=2\tau\int_M\left|R_{ij}+\nabla_i\nabla_jf-\frac{1}{2\tau}g_{ij}\right|^2
(4\pi\tau)^{-\frac{n}{2}}e^{-f}dV.
$$

Finally, by using the computations in (1.5.5) and the evolution
equations of $f$ and $\tau$, we have
\begin{align*}
\frac{\partial}{\partial t}\((4\pi\tau)^{-\frac{n}{2}} e^{-f}dV\)
&= (4\pi\tau)^{-\frac{n}{2}}\left[\frac{\partial}{\partial t}(
e^{-f}dV)+\frac {n}{2\tau} e^{-f}dV\right]\\
&= -(4\pi\tau)^{-\frac{n}{2}}\Delta (e^{-f})dV.
\end{align*}
Hence
$$
\frac{d}{d t}\int_M(4\pi\tau)^{-\frac{n}{2}}e^{-f}dV
=-(4\pi\tau)^{-\frac{n}{2}}\int_M\Delta (e^{-f})dV=0.
$$
\end{pf}

Now we set \be
\mu(g_{ij},\tau)=\inf\left\{\mathcal{W}(g_{ij},f,\tau)\ |\  f\in
C^\infty(M),
\frac{1}{(4\pi\tau)^{n/2}}\int_M e^{-f}dV=1\right\}  
\ee and
$$
\nu(g_{ij})=\inf\left\{\mathcal W(g,f,\tau)\ |\  f\in C^\infty(M),
\tau>0, \frac{1}{(4\pi\tau)^{n/2}}\int e^{-f}dV=1\right\}.
$$
Note that if we let $u=e^{-f/2}$, then the functional
$\mathcal{W}$ can be expressed as
$$
\mathcal{W}(g_{ij},f,\tau) =\int_M[\tau(Ru^2+4|\nabla
u|^2)-u^2\log u^2-nu^2] (4\pi\tau)^{-\frac{n}{2}}dV
$$
and the constraint $\int_M(4\pi\tau)^{-\frac{n}{2}}e^{-f}dV=1$
becomes $\int_Mu^2(4\pi\tau)^{-\frac{n}{2}}dV=1.$ Thus
$\mu(g_{ij},\tau)$ corresponds to the best constant of a
logarithmic Sobolev inequality. Since the nonquadratic term is
subcritical (in view of Sobolev exponent), it is rather
straightforward to show that
$$
\inf\!\bigg\{\!\int_M\![\tau(4|\nabla u|^2\!+Ru^2) -u^2\log
u^2\!-nu^2](4\pi\tau)^{-\frac{n}{2}}dV \Big|
\int_Mu^2(4\pi\tau)^{-\frac{n}{2}}dV\!=\!1\!\bigg\}
$$
is achieved by some nonnegative function $u\in H^1(M)$ which
satisfies the Euler-Lagrange equation
$$\tau (-4\Delta u+Ru)-2u\log u-nu=\mu(g_{ij},\tau)u
.$$ One can further show that $u$ is positive (see \cite{Ro}).
Then the standard regularity theory of elliptic PDEs shows that
$u$ is smooth. We refer the reader to Rothaus \cite{Ro} for more
details. It follows that $\mu(g_{ij},\tau)$ is achieved by a
minimizer $f$ satisfying the nonlinear equation \be
\tau (2\Delta f-|\nabla f|^2+R)+f-n=\mu(g_{ij},\tau).
\ee

\begin{corollary} [{Perelman \cite{P1}}]\
\begin{itemize}
\item[(i)] $\mu(g_{ij}(t),\tau-t)$ is nondecreasing along the
Ricci flow; moveover, the monotonicity is strict unless we are on
a shrinking gradient soliton; \item[(ii)] A shrinking breather is
necessarily a shrinking gradient soliton.
\end{itemize}
\end{corollary}

\begin{pf}
Fix any time $t_0$, let $f_0$ be a minimizer of
$\mu(g_{ij}(t_0),\tau-t_0).$ Note that the backward heat equation
$$
\frac{\partial f}{\partial t}=-\Delta f+|\nabla
f|^2-R+\frac{n}{2\tau}
$$
is equivalent to the linear equation
$$
\frac{\partial }{\partial
t}((4\pi\tau)^{-\frac{n}{2}}e^{-f})=-\Delta
((4\pi\tau)^{-\frac{n}{2}}e^{-f})+R((4\pi\tau)^{-\frac{n}{2}}e^{-f}).
$$

Thus we can solve the backward heat equation of $f$ with
$f|_{t=t_0}=f_0$ to obtain $f(t)$, $t\leq t_0$, with
$\int_M(4\pi\tau)^{-\frac{n}{2}}e^{-f(t)}dV=1.$ It then follows
from Proposition 1.5.8 that
\begin{align*}
\mu(g_{ij}(t),\tau-t)
& \leq \mathcal{W}(g_{ij}(t),f(t),\tau-t)\\
& \leq \mathcal{W}(g_{ij}(t_0),f(t_0),\tau-t_0)\\
& = \mu(g_{ij}(t_0),\tau-t_0)
\end{align*}
for $t\leq t_0$ and the second inequality is strict unless we are
on a shrinking gradient soliton. This proves statement (i).

Consider a shrinking breather on $[t_1,t_2]$ with $\alpha
g_{ij}(t_1)$ and $g_{ij}(t_2)$ differ only by a diffeomorphism for
some $\alpha < 1.$ Recall that the functional $\mathcal{W}$ is
invariant under simultaneous scaling of $\tau$ and $g_{ij}$ and
invariant under diffeomorphism. Then for $\tau >0$ to be
determined,
$$
\mu(g_{ij}(t_1),\tau-t_1)=\mu(\alpha
g_{ij}(t_1),\alpha(\tau-t_1))=\mu(g_{ij}(t_2),\alpha(\tau-t_1))
$$
and by the monotonicity of $\mu(g_{ij}(t),\tau-t),$
$$
\mu(g_{ij}(t_1),\tau-t_1)\leq\mu(g_{ij}(t_2),\tau-t_2).
$$
Now take $\tau>0$ such that
$$
\alpha(\tau-t_1)=\tau-t_2,
$$
i.e.,
$$
\tau=\frac{t_2-\alpha t_1}{1-\alpha}.
$$

\noindent This shows the equality holds in the monotonicity of
$\mu(g_{ij}(t),\tau-t).$ So the shrinking breather must be a
shrinking gradient soliton.
\end{pf}

Finally, we remark that Hamilton, Ilmanen and the first author
\cite{Cao04} have obtained the second variation formulas for both
$\lambda$-energy and $\nu$-energy. We refer the reader to their
paper \cite{Cao04}  for more details and related stability
questions.

\newpage
\part{{\Large Maximum Principle and Li-Yau-Hamilton Inequalities}}

\bigskip
The maximum principle is a fundamental tool in the study of
parabolic equations in general. In this chapter, we present
various maximum principles for tensors developed by Hamilton in
the Ricci flow. As an immediate consequence, the Ricci flow
preserves the nonnegativity of the curvature operator. We also
present two crucial estimates in the Ricci flow: the Hamilton-Ivey
curvature pinching estimate \cite{Ha95F, Iv} when dimension $n=3$,
and the Li-Yau-Hamilton estimate \cite{Ha88, Ha93} from which one
obtains the Harnack inequality for the evolved scalar curvature
via a Li-Yau \cite{LY} path integral. Most of the presentation in
Sections 2.1-2.5 follows closely Hamilton \cite{Ha82, Ha86, Ha88,
Ha93, Ha95F}, and some parts of Section 2.3 also follows Chow-Lu
\cite{Cl04}. Finally, in Section 2.6 we describe Perelman's Li-Yau
type estimate for solutions to the conjugate heat equation and
show how Li-Yau type path integral leads to a space-time distance
function (i.e., what Perelman \cite{P1} called the reduced
distance).

\section{Preserving Positive Curvature}

Let $M$ be an $n$-dimensional complete manifold. Consider a family
of smooth metrics $g_{ij}(t)$ evolving by the Ricci flow with
uniformly bounded curvature for $t\in[0,T]$ with $T<+\infty$.
Denote by $d_t(x,y)$ the distance between two points $x,y\in M$
with respect to the metric $g_{ij}(t)$. First we need the
following useful fact (cf. \cite{ScY})

\begin{lemma}
There exists a smooth function $f$ on $M$ such that $f\ge 1$
everywhere, $f(x)\rightarrow +\infty$ as $d_0(x,x_0)\rightarrow
+\infty$ $($for some fixed $x_0\in M),$
$$
|\nabla f|_{g_{ij}(t)}\le C \quad \text{and}\quad |\nabla^2
f|_{g_{ij}(t)}\le C
$$
on $M \times [0,T]$ for some positive constant $C$.
\end{lemma}

\begin{pf}
Let $\varphi(v)$ be a smooth function on $\mathbb{R}^n$ which is
nonnegative, rotationally symmetric and has compact support in a
small ball centered at the origin with
$\int_{\mathbb{R}^n}\varphi(v)dv=1$.

For each $x\in M$, set
$$
f(x)=\int_{\mathbb{R}^n}\varphi(v)(d_0(x_0,exp_x(v))+1)dv,
$$
where the integral is taken over the tangent space $T_xM$ at $x$
which we have identified with $\mathbb{R}^n$. If the size of the
support of $\varphi(v)$ is small compared to the maximum
curvature, then it is well known that this defines a smooth
function $f$ on $M$ with $f(x)\rightarrow +\infty$ as
$d_0(x,x_0)\rightarrow +\infty$, while the bounds on the first and
second covariant derivatives of $f$ with respect to the metric
$g_{ij}(\cdot,0)$ follow from the Hessian comparison theorem. Thus
it remains to show these bounds hold with respect to the evolving
metric ${g_{ij}(t)}$.

We compute, using the frame $\{F_a^i\nabla_if\}$ introduced in
Section 1.3,
$$
\frac{\partial}{\partial t}\nabla_a f =\frac{\partial}{\partial
t}(F_a^i\nabla_if) =R_{ab}\nabla_b f.
$$
Hence
$$
|\nabla f|\le C_1\cdot e^{C_2 t},
$$
where $C_1, C_2$ are some positive constants depending only on the
dimension. Also
\begin{align*}
\frac{\partial}{\partial t}(\nabla_a\nabla_bf) &
=\frac{\partial}{\partial t} \(F_a^iF_b^j\(\frac{\partial^2
f}{\partial x^i\partial x^j}
-\Gamma_{ij}^k\frac{\partial f}{\partial x^k}\)\)\\
& =R_{ac}\nabla_b\nabla_cf+R_{bc}\nabla_a\nabla_cf
+(\nabla_cR_{ab}-\nabla_aR_{bc}-\nabla_bR_{ac})\nabla_cf.
\end{align*}
Then by Shi's derivative estimate (Theorem 1.4.1), we have
$$
\frac{\partial}{\partial t}|\nabla^2 f|\le C_3|\nabla^2 f|
+\frac{C_3}{\sqrt{t}},
$$
which implies
$$
|\nabla^2f|_{g_{ij}(t)}\le e^{C_3 t}\(|\nabla^2f|_{g_{ij}(0)}
+\int_0^t\frac{C_3}{\sqrt{\tau}}e^{-C_3\tau}d\tau\)
$$
for some positive constants $C_3$ depending only on the dimension
and the curvature bound.
\end{pf}

We now use the weak maximum principle to derive the following result
(cf. \cite{Ha82} and \cite{Sh89}).

\begin{proposition}
If the scalar curvature $R$ of the solution $g_{ij}(t), 0\le t\le
T$, to the Ricci flow is nonnegative at $t=0$, then it remains so
on $0\le t\le T$.
\end{proposition}

\begin{pf}
Let $f$ be the function constructed in Lemma 2.1.1 and recall
$$
\frac{\partial R}{\partial t}=\Delta R+2|\Ric|^2.
$$
For any small constant $\varepsilon>0$ and large constant $A>0$,
we have
\begin{align*}
\frac{\partial}{\partial t}(R+\varepsilon e^{At}f)
&  = \frac{\partial R}{\partial t}+\varepsilon Ae^{At}f\\
&  = \Delta(R+\varepsilon e^{At}f)+2|Ric|^2+\varepsilon
e^{At}(Af-\Delta f)\\
& > \Delta(R+\varepsilon e^{At}f)
\end{align*}
by choosing $A$ large enough.

We claim that
$$
R+\varepsilon e^{At}f>0\qquad \text{on}\quad M\times[0,T].
$$
Suppose not, then there exist a first time $t_0>0$ and a point
$x_0\in M$ such that
\begin{align*}
(R+\varepsilon e^{At}f)(x_0,t_0)&=0,\\
\nabla(R+\varepsilon e^{At}f)(x_0,t_0)&=0,\\
\Delta(R+\varepsilon e^{At}f)(x_0,t_0)&\ge 0,\\
\text{and}\quad \frac{\partial}{\partial t}(R+\varepsilon
e^{At}f)(x_0,t_0)&\le 0.
\end{align*}
Then
$$
0\ge \frac{\partial}{\partial t}(R+\varepsilon
e^{At}f)(x_0,t_0)>\Delta(R+\varepsilon e^{At}f)(x_0,t_0)\ge 0,
$$
which is a contradiction. So we have proved that
$$
R+\varepsilon e^{At}f>0\qquad \text{on}\quad M\times[0,T].
$$
Letting $\varepsilon\rightarrow 0$, we get
$$
R\ge 0\qquad \text{on}\quad M\times[0,T].
$$
This finishes the proof of the proposition.
\end{pf}

Next we derive a maximum principle of Hamilton for tensors. Let
$M$ be a complete manifold with a metric $g=\{g_{ij}\}$, $V$ a
vector bundle over $M$ with a metric $h=\{h_{\alpha\beta}\}$ and a
connection $\nabla=\{\Gamma^\alpha_{i\beta}\}$ compatible with
$h$, and suppose $h$ is fixed but $g$ and $\nabla$ may vary
smoothly with time $t$. Let $\Gamma(V)$ be the vector space of
$C^\infty$ sections of $V$. The Laplacian $\Delta$ acting on a
section $\sigma\in\Gamma(V)$ is defined by
$$
\Delta\sigma=g^{ij}\nabla_i\nabla_j\sigma.
$$
Let $M_{\alpha\beta}$ be a symmetric bilinear form on $V$. We say
$M_{\alpha\beta}\ge 0$ if $M_{\alpha\beta}v^\alpha v^\beta\ge 0$
for all vectors $v=\{v^\alpha\}$. Assume
$N_{\alpha\beta}={\mathcal{P}}(M_{\alpha\beta},h_{\alpha\beta})$
is a polynomial in $M_{\alpha\beta}$ formed by contracting
products of $M_{\alpha\beta}$ with itself using the metric
$h=\{h_{\alpha\beta}\}$. Assume that the tensor $M_{\alpha\beta}$
is uniformly bounded in space-time and let $g_{ij}$ evolve by the
Ricci flow with bounded curvature.

\begin{lemma} [{Hamilton \cite{Ha82}}]
Suppose that on $0\le t\le T$,
$$
\frac{\partial}{\partial t}M_{\alpha\beta} =\Delta
M_{\alpha\beta}+u^i\nabla_i M_{\alpha\beta}+N_{\alpha\beta}
$$
where $u^i(t)$ is a time-dependent vector field on $M$ with
uniform bound and
$N_{\alpha\beta}={\mathcal{P}}(M_{\alpha\beta},h_{\alpha\beta})$
satisfies
$$
N_{\alpha\beta}v^\alpha v^\beta\ge 0 \quad \text{whenever} \quad
M_{\alpha\beta}v^\beta=0.
$$
If $M_{\alpha\beta}\ge 0$ at $t=0$, then it remains so on $0\le
t\le T$.
\end{lemma}

\begin{pf}
Set
$$
\tilde{M}_{\alpha\beta} =M_{\alpha\beta}+\varepsilon e^{At}f
h_{\alpha\beta},
$$
where $A>0$ is a suitably large constant (to be chosen later) and
$f$ is the function constructed in Lemma 2.1.1.

We claim that $\tilde{M}_{\alpha\beta}>0$ on $M\times[0,T]$ for
every $\varepsilon>0$. If not, then for some $\varepsilon>0$,
there will be a first time $t_0>0$ where $\tilde{M}_{\alpha\beta}$
acquires a null vector $v^\alpha$ of unit length at some point
$x_0\in M$. At $(x_0,t_0)$,
\begin{align*}
N_{\alpha\beta}v^\alpha v^\beta & \geq N_{\alpha\beta}v^\alpha
v^\beta
-\tilde{N}_{\alpha\beta}v^\alpha v^\beta\\
& \ge -C\varepsilon e^{At_0}f(x_0),
\end{align*}
where
$\tilde{N}_{\alpha\beta}=\mathcal{P}(\tilde{M}_{\alpha\beta},
h_{\alpha\beta})$, and $C$ is a positive constant (depending on
the bound of $M_{\alpha\beta}$, but independent of $A$).

Let us extend $v^\alpha$ to a local vector field in a neighborhood
of $x_0$ by parallel translating $v^{\alpha}$ along geodesics
(with respect to the metric $g_{ij}(t_0)$) emanating radially out
of $x_0$, with $v^{\alpha}$ independent of $t$. Then, at
$(x_0,t_0)$, we have
\begin{align*}
\frac{\partial}{\partial t}(\tilde{M}_{\alpha\beta}v^\alpha
v^\beta)&  \le 0 ,\\
\nabla (\tilde{M}_{\alpha\beta}v^\alpha v^\beta)&  = 0, \\
 \text{and} \quad \Delta(\tilde{M}_{\alpha\beta}v^\alpha
v^\beta)&  \ge 0.
\end{align*}
But
\begin{align*}
0\ge\frac{\partial}{\partial t}(\tilde{M}_{\alpha\beta}v^\alpha
v^\beta) &  = \frac{\partial}{\partial t}(M_{\alpha\beta}v^\alpha
v^\beta+\varepsilon e^{At}f),\\
&  = \Delta(\tilde{M}_{\alpha\beta}v^\alpha
v^\beta)-\Delta(\varepsilon
e^{At}f)+u^i\nabla_i(\tilde{M}_{\alpha\beta}v^\alpha
v^\beta)\\
&\quad -u^i\nabla_i(\varepsilon e^{At}f)+N_{\alpha\beta}v^\alpha
v^\beta+\varepsilon A e^{At_0}f(x_0)\\
&  \ge -C\varepsilon e^{At_0}f(x_0)+\varepsilon A e^{At_0}f(x_0)>0
\end{align*}
when $A$ is chosen sufficiently large. This is a contradiction.
\end{pf}

By applying Lemma 2.1.3 to the evolution equation
$$
\frac{\partial}{\partial t}M_{\alpha\beta} =\Delta
M_{\alpha\beta}+M_{\alpha\beta}^2+M_{\alpha\beta}^\#
$$
of the curvature operator $M_{\alpha\beta}$, we immediately obtain
the following important result.

\begin{proposition}[{Hamilton \cite{Ha86}}]
Nonnegativity of the curvature operator $M_{\alpha\beta}$ is
preserved by the Ricci flow.
\end{proposition}   

In the K\"ahler case, the nonnegativity of the holomorpic
bisectional curvature is preserved under the K\"ahler-Ricci flow.
This result is proved by Bando \cite{Bando} for complex dimension
$n=3$ and by Mok \cite{Mo88} for general dimension $n$ when the
manifold is compact, and by Shi \cite{Sh90} when the manifold is
noncompact.

\begin{proposition}
Under the K\"ahler-Ricci flow if the initial metric has positive
$($nonnegative$)$ holomorphic bisectional curvature then the
evolved metric also has positive $($nonnegative$)$ holomorphic
bisectional curvature.
\end{proposition}   

\section{Strong Maximum Principle}

Let $\Omega$ be a bounded, connected open set of a complete
$n$-dimensional manifold $M$, and let $g_{ij}(x,t)$ be a smooth
solution to the Ricci flow on $\Omega \times [0,T]$. Consider a
vector bundle $V$ over $\Omega$ with a fixed metric
$h_{\alpha\beta}$ (independent of time), and a connection $\nabla
= \{\Gamma^\alpha_{i\beta}\}$ which is compatible with
$h_{\alpha\beta}$ and may vary with time $t$. Let $\Gamma(V)$ be
the vector space of $C^\infty$ sections of $V$ over $\Omega$. The
Laplacian $\Delta$ acting on a section $\sigma\in\Gamma(V)$ is
defined by
$$
\Delta\sigma=g^{ij}(x,t)\nabla_i \nabla_j\sigma.
$$

Consider a family of smooth symmetric bilinear forms
$M_{\alpha\beta}$ evolving by \be \frac{\partial}{\partial
t}M_{\alpha\beta}=\Delta M_{\alpha\beta}+N_{\alpha\beta},\ \
\text{ on }\
\Omega\times[0,T], 
\ee where $N_{\alpha\beta}=P(M_{\alpha\beta},h_{\alpha\beta})$ is
a polynomial in $M_{\alpha\beta}$ formed by contracting products
of $M_{\alpha\beta}$ with itself using the metric
$h_{\alpha\beta}$ and satisfies
$$
N_{\alpha\beta}\geq 0, \; \mbox{ whenever }\; M_{\alpha\beta}\geq
0.
$$
The following result, due to Hamilton \cite{Ha86}, shows that the
solution of (2.2.1) satisfies a strong maximum principle.

\begin{theorem}[Hamilton's strong maximum principle
\index{Hamilton's strong maximum principle}\!\!] Let
$M_{\alpha\beta}$ be a smooth solution of the equation $(2.2.1)$.
Suppose $M_{\alpha\beta}\geq0$ on $\Omega \times [0,T]$. Then
there exists a positive constant $0<\delta \leq T$ such that on
$\Omega \times (0,\delta)$, the rank of $M_{\alpha\beta}$ is
constant, and the null space of $M_{\alpha\beta}$ is invariant
under parallel translation and invariant in time and also lies in
the null space of $N_{\alpha\beta}$.
\end{theorem}

\begin{pf}
Set
$$
l=\max_{x\in \Omega}\{\text{rank of } M_{\alpha\beta}(x,0)\}.
$$
Then we can find a nonnegative smooth function $\rho(x)$, which is
positive somewhere and has compact support in $\Omega$, so that at
every point $x \in \Omega$,
$$
\sum^{n-l+1}_{i=1}M_{\alpha\beta}(x,0)v^{\alpha}_{i}v^{\beta}_{i}
\geq \rho(x)
$$
for any $(n-l+1)$ orthogonal unit vectors $\{ v_1, \ldots,
v_{n-l+1} \}$ at $x$.

Let us evolve $\rho(x)$ by the heat equation
$$
\frac{\partial}{\partial t}\rho=\Delta\rho
$$
with the Dirichlet condition $\rho|_{\partial \Omega}=0$ to get a
smooth function $\rho(x,t)$ defined on $\Omega \times [0,T]$. By
the standard strong maximum principle, we know that $\rho(x,t)$ is
positive everywhere in $\Omega$ for all $t\in (0,T]$.

For every $\varepsilon>0$, we claim that at every point $(x,t) \in
\Omega \times [0,T]$, there holds
$$
\sum^{n-l+1}_{i=1}M_{\alpha\beta}(x,t)v^{\alpha}_{i}v^{\beta}_{i}
 + \varepsilon e^t > \rho(x,t)
$$
for any $(n-l+1)$ orthogonal unit vectors $\{ v_1, \ldots,
v_{n-l+1} \}$ at $x$.

We argue by contradiction. Suppose not, then for some
$\varepsilon>0$, there will be a first time $t_0 > 0$ and some
$(n-l+1)$ orthogonal unit vectors $\{ v_1, \ldots, v_{n-l+1} \}$
at some point $x_0 \in \Omega$ so that
$$
\sum^{n-l+1}_{i=1}M_{\alpha\beta}(x_0,t_0)v^{\alpha}_{i}v^{\beta}_{i}
 + \varepsilon e^{t_0} = \rho(x_0,t_0)
$$

Let us extend each $v_i$ ($i=1, \ldots, n-l+1$) to a local vector
field, independent of $t$, in a neighborhood of $x_0$ by parallel
translation along geodesics (with respect to the metric
$g_{ij}(t_0)$) emanating radially out of $x_0$. Clearly $\{v_1,
\ldots, v_{n-l+1}\}$ remain orthogonal unit vectors in the
neighborhood.  Then, at $(x_0,t_0)$, we have
\begin{align*}
\frac{\partial}{\partial t}
\(\sum^{n-l+1}_{i=1}M_{\alpha\beta}v^{\alpha}_{i}v^{\beta}_{i} +
\varepsilon e^t -\rho\)
& \le 0 ,\\
\text{and} \quad
\Delta\(\sum^{n-l+1}_{i=1}M_{\alpha\beta}v^{\alpha}_{i}v^{\beta}_{i}
 + \varepsilon e^t -\rho\)&  \ge 0.
\end{align*}
But, since $N_{\alpha\beta}\ge 0$ by our assumption, we have
\begin{align*}
0&  \ge \frac{\partial}{\partial t}
\(\sum^{n-l+1}_{i=1}M_{\alpha\beta}v^{\alpha}_{i}v^{\beta}_{i}
+ \varepsilon e^t -\rho\)\\
&  = \sum^{n-l+1}_{i=1}(\Delta
M_{\alpha\beta}+N_{\alpha\beta})v^{\alpha}_{i}v^{\beta}_{i}
 + \varepsilon e^t -\Delta \rho\\
& \geq\sum^{n-l+1}_{i=1}\Delta(M_{\alpha\beta}
v^{\alpha}_{i}v^{\beta}_{i})+ \varepsilon e^t -\Delta \rho\\
&  =\sum^{n-l+1}_{i=1}\Delta(M_{\alpha\beta}v^{\alpha}_{i}
v^{\beta}_{i}+\varepsilon e^t -\rho)+ \varepsilon e^t \\
& \ge \varepsilon e^t >0.
\end{align*}
This is a contradiction. Thus by letting $\varepsilon\rightarrow
0$, we prove that
$$
\sum^{n-l+1}_{i=1}M_{\alpha\beta}(x,t)v^{\alpha}_{i}v^{\beta}_{i}
 \geq \rho(x,t)
$$
for any $(n-l+1)$ orthogonal unit vectors $\{ v_1, \ldots,
v_{n-l+1} \}$ at $x \in \Omega$ and $t \in [0,T]$. Hence
$M_{\alpha\beta}$ has at least rank $l$ everywhere in the open set
$\Omega$ for all $t \in (0,T]$. Therefore we can find a positive
constant $\delta (\leq T)$ such that the rank $M_{\alpha\beta}$ is
constant over $\Omega \times (0,\delta)$.

Next we proceed to analyze the null space of $M_{\alpha\beta}$.
Let $v$ be any smooth section of $V$ in the null of
$M_{\alpha\beta}$ on $0<t<\delta$. Then
\begin{align*}
0&  = \frac{\partial}{\partial t}(M_{\alpha\beta}v^\alpha
v^\beta)\\
&  = \(\frac{\partial}{\partial t}M_{\alpha\beta}\)v^\alpha
v^\beta+2M_{\alpha\beta}v^\alpha\frac{\partial v^\beta}{\partial t}\\
&  = \(\frac{\partial}{\partial t}M_{\alpha\beta}\)v^\alpha
v^\beta,
\end{align*}
and
\begin{align*}
0&  = \Delta(M_{\alpha\beta}v^\alpha v^\beta)\\
&  = (\Delta M_{\alpha\beta})v^\alpha v^\beta
+4g^{kl}\nabla_kM_{\alpha\beta}\cdot v^\alpha\nabla_lv^\beta\\
&\quad +2M_{\alpha\beta}g^{kl}\nabla_kv^\alpha\cdot\nabla_lv^\beta
+2M_{\alpha\beta}v^\alpha\Delta v^\beta\\
&  = (\Delta M_{\alpha\beta})v^\alpha
v^\beta+4g^{kl}\nabla_kM_{\alpha\beta}\cdot
v^\alpha\nabla_lv^\beta
+2M_{\alpha\beta}g^{kl}\nabla_kv^\alpha\cdot\nabla_lv^\beta.
\end{align*}
By noting that
$$
0=\nabla_k(M_{\alpha\beta}v^\beta)
=(\nabla_kM_{\alpha\beta})v^\alpha+M_{\alpha\beta}\nabla_kv^\alpha
$$
and using the evolution equation (2.2.1), we get
$$
N_{\alpha\beta}v^\alpha v^\beta
+2M_{\alpha\beta}g^{kl}\nabla_kv^\alpha\cdot\nabla_lv^\beta=0.
$$
Since $M_{\alpha\beta}\ge 0$ and $N_{\alpha\beta}\ge 0$, we must
have
$$
v\in \text{null}\; (N_{\alpha\beta})\quad \text{and} \quad
\nabla_i v\in \text{null}\; (M_{\alpha\beta}),\quad \mbox{for all
}\; i.
$$
The first inclusion shows that $\text{null}\;
(M_{\alpha\beta})\subset \text{null}\; (N_{\alpha\beta}),$ and the
second inclusion shows that $\text{null}\; (M_{\alpha\beta})$ is
invariant under parallel translation.

To see $\text{null}\; (M_{\alpha\beta})$ is also invariant in
time, we first note that
$$
\Delta v=\nabla^i(\nabla_i v)\in \text{null}\; (M_{\alpha\beta})
$$
and then
$$
g^{kl}\nabla_kM_{\alpha\beta}\cdot\nabla_lv^\alpha
=g^{kl}\nabla_k(M_{\alpha\beta}\nabla_lv^\alpha)
-M_{\alpha\beta}\Delta v^\alpha=0.
$$
Thus we have
\begin{align*}
0& = \Delta(M_{\alpha\beta}v^\alpha)\\
& = (\Delta M_{\alpha\beta})v^\alpha
+2g^{kl}\nabla_kM_{\alpha\beta}\cdot \nabla_lv^\alpha
+M_{\alpha\beta}\Delta v^\alpha\\
&  = (\Delta M_{\alpha\beta})v^\alpha,
\end{align*}
and hence
\begin{align*}
0&  = \frac{\partial}{\partial t}(M_{\alpha\beta}v^\alpha)\\
&  = (\Delta M_{\alpha\beta}+N_{\alpha\beta})v^\alpha
+M_{\alpha\beta}\frac{\partial v^\alpha}{\partial t}\\
&  = M_{\alpha\beta}\frac{\partial v^\alpha}{\partial t}.
\end{align*}
This shows that
$$
\frac{\partial v}{\partial t}\in \text{null}\; (M_{\alpha\beta}),
$$
so the null space of $M_{\alpha\beta}$ is invariant in time.
\end{pf}

We now apply Hamilton's strong maximum principle to the evolution
equation of the curvature operator $M_{\alpha\beta}$. Recall
$$
\frac{\partial M_{\alpha\beta}}{\partial t} =\Delta
M_{\alpha\beta}+M_{\alpha\beta}^2+M_{\alpha\beta}^{\#}
$$
where $M_{\alpha\beta}^{\#}
=C_\alpha^{\xi\gamma}C_\beta^{\eta\theta}M_{\xi\eta}M_{\gamma\theta}$.
Suppose we have a solution to the Ricci flow with nonnegative
curvature operator. Then by Theorem 2.2.1, the null space of the
curvature operator $M_{\alpha\beta}$ of the solution has constant
rank and is invariant in time and under parallel translation over
some time interval $0<t<\delta$ . Moreover the null space of
$M_{\alpha\beta}$ must also lie in the null space of
$M_{\alpha\beta}^\#$.

Denote by $(n-k)$ the rank of $M_{\alpha\beta}$ on $0<t<\delta$.
Let us diagonalize $M_{\alpha\beta}$ so that $M_{\alpha\alpha}=0$
if $\alpha\le k$ and $M_{\alpha\alpha}>0$ if $\alpha>k$. Then we
have $M_{\alpha\alpha}^\#=0$ also for $\alpha\le k$ from the
evolution equation of $M_{\alpha\alpha}$. Since
$$
0=M_{\alpha\alpha}^\#
=C_\alpha^{\xi\gamma}C_\alpha^{\eta\theta}M_{\xi\eta}M_{\gamma\theta},
$$
it follows that
\begin{align*}
C_\alpha^{\xi\gamma}&  = \<v^\alpha,[v^\xi,v^\gamma]\>\\
&  = 0, \quad \text{if}\quad \alpha\le k\; \text{ and }\;
\xi,\gamma> k.
\end{align*}
This says that the image of $M_{\alpha\beta}$ is a Lie subalgebra
(in fact it is the subalgebra of the restricted holonomy group by
using the Ambrose-Singer holonomy theorem \cite{As53}). This
proves the following result.

\begin{theorem}[{Hamilton \cite{Ha86}}]
Suppose the curvature operator $M_{\alpha\beta}$ of the initial
metric is nonnegative. Then, under the Ricci flow, for some
interval $0<t<\delta$ the image of $M_{\alpha\beta}$ is a Lie
subalgebra of $so(n)$ which has constant rank and is invariant
under parallel translation and invariant in time.
\end{theorem}

\section{Advanced Maximum Principle for Tensors}

In this section we present Hamilton's advanced maximum principle
\cite{Ha86} for tensors which generalizes Lemma 2.1.3 and shows
how a tensor evolving by a nonlinear heat equation may be
controlled by a system of ODEs. Our presentation in this section
follows closely the papers of Hamilton \cite{Ha86} and Chow-Lu
\cite{Cl04}. An important application of the advanced maximum
principle is the Hamilton-Ivey curvature pinching estimate for the
Ricci flow on three-manifolds given in the next section. More
applications will be given in Chapter 5.

Let $M$ be a complete manifold equipped with a one-parameter
family of Riemannian metrics $g_{ij}(t)$, $0\le t\le T$, with
$T<+\infty$. Let $V\rightarrow M$ be a vector bundle with a
time-independent bundle metric $h_{ab}$ and $\Gamma(V)$ be the
vector space of $C^\infty$ sections of $V$. Let
$$
\nabla_t:\quad \Gamma(V)\rightarrow\Gamma(V\otimes T^*M), \quad
t\in[0,T]
$$
be a smooth family of time-dependent connections compatible with
$h_{ab}$, i.e.
$$
(\nabla_t)_ih_{ab}\stackrel{\Delta}{=}
(\nabla_t)_{\frac{\partial}{\partial x^i}}h_{ab}=0,
$$
for any local coordinate $\{\frac{\partial}{\partial
x^1},\ldots,\frac{\partial}{\partial x^n}\}.$ The Laplacian
$\Delta_t$ acting on a section $\sigma\in\Gamma(V)$ is defined by
$$
\Delta_t\sigma=g^{ij}(x,t)(\nabla_t)_i(\nabla_t)_j\sigma.
$$
For the application to the Ricci flow, we will always assume that
the metrics $g_{ij}(\cdot,t)$ evolve by the Ricci flow. Since $M$
may be noncompact, we assume that, for the sake of simplicity, the
curvature of $g_{ij}(t)$ is uniformly bounded on $M\times [0,T]$.

Let $N: V\times [0,T]\rightarrow V$ be a fiber preserving map,
i.e., $N(x,\sigma,t)$ is a  time-dependent vector field defined on
the bundle $V$ and tangent to the fibers. We assume that
$N(x,\sigma,t)$ is continuous in $x,t$ and satisfies
$$
|N(x,\sigma_1,t) - N(x,\sigma_2,t)| \leq C_B |\sigma_1 -\sigma_2|
$$
for all $ x \in M$, $t\in [0,T]$ and $|\sigma_1|\leq B,
|\sigma_2|\leq B$, where $C_B$ is a positive constant depending
only on $B$. Then we can form the nonlinear heat equation
\begin{equation*}
\frac{\partial}{\partial t}\sigma(x,t)=
\Delta_t\sigma(x,t)+u^i(\nabla_t)_i\sigma(x,t)+N(x,\sigma(x,t),t)
\tag{PDE}
\end{equation*}
where $u^i=u^i(t)$ is a time-dependent vector field on $M$ which
is uniformly bounded on $M\times[0,T]$. Let $K$ be a closed subset
of $V$. One important question is under what conditions will
solutions of the PDE which start in $K$ remain in $K$. To answer
this question, Hamilton \cite{Ha86} imposed the following two
conditions on $K$:
\begin{itemize}
\item[(H1)] $K$ is invariant under parallel translation defined by
the connection $\nabla_t$ for each $t\in[0,T]$; \item[(H2)] in
each fiber $V_x$, the set $K_x\stackrel{\Delta}{=}V_x\cap K$ is
closed and convex.
\end{itemize}

\noindent Then one can judge the behavior of the PDE by comparing
to that of the following ODE
\begin{equation*}
\frac{d\sigma_x}{dt}=N(x,\sigma_x,t)\tag{ODE}
\end{equation*}
for $\sigma_x=\sigma_x(t)$ in each fiber $V_x$. The following
version of Hamilton's advanced maximum principle is from Chow-Lu
\cite{Cl04}

\begin{theorem}[{Hamilton's advanced maximum
principle \index{Hamilton's advanced maximum principle}
\cite{Ha86}}] Let $K$ be a closed subset of $V$ satisfying the
hypothesis {\rm (H1)} and {\rm (H2)}. Suppose that for any $x\in
M$ and any initial time $t_0\in[0,T)$, any solution $\sigma_x(t)$
of the {\rm (ODE)} which starts in $K_x$ at $t_0$ will remain in
$K_x$ for all later times. Then for any initial time $t_0\in
[0,T)$ the solution $\sigma(x,t)$ of the {\rm (PDE)} will remain
in $K$ for all later times provided $\sigma(x,t)$ starts in $K$ at
time $t_0$ and $\sigma(x,t)$ is uniformly bounded with respect to
the bundle metric $h_{ab}$ on $M\times[t_0,T]$.
\end{theorem}

We remark that Lemma 2.1.3 is a special case of the above theorem
where $V$ is given by a symmetric tensor product of a vector
bundle and $K$ corresponds to the convex set consisting of all
nonnegative symmetric bilinear forms. We also remark that Hamilton
\cite{Ha86} established the above theorem for a general evolving
metric $g_{ij}(x,t)$ which does not necessarily satisfy the Ricci
flow.

Before proving Theorem 2.3.1, we need to establish three lemmas in
\cite{Ha86}. Let $\varphi:[a,b]\rightarrow \mathbb{R}$ be a
Lipschitz function. We consider $\frac{d\varphi}{dt}(t)$ at
$t\in[a,b)$ in the sense of limsup of the forward difference
quotients, i.e.,
$$
\frac{d\varphi}{dt}(t) =\limsup_{h\rightarrow
0^+}\frac{\varphi(t+h)-\varphi(t)}{h}.
$$

\begin{lemma} [{Hamilton \cite{Ha86}}]
Suppose $\varphi: [a,b]\rightarrow \mathbb{R}$ is Lipschitz
continuous and suppose for some constant $C<+\infty$,
\begin{align*}
&  \frac{d}{dt}\varphi(t)\le C\varphi(t),\quad whenever\ \
\varphi(t)\ge 0\;\ on\ \;[a,b),\\[2mm]
\text{and}\hskip 1cm &  \varphi(a)\le 0.
\end{align*}
Then $\varphi(t)\le 0$ on $[a,b]$.
\end{lemma}

\begin{pf}
By replacing $\varphi$ by $e^{-Ct}\varphi$, we may assume
\begin{align*}
&  \frac{d}{dt}\varphi(t)\le 0,\quad \text{whenever}\quad
\varphi(t)\ge 0\;\mbox{ on }\;[a,b),\\
\mbox{and }\hskip 1cm &  \varphi(a)\le 0.
\end{align*}
For arbitrary $\varepsilon>0$, we shall show $\varphi(t)\le
\varepsilon (t-a)$ on $[a,b]$. Clearly we may assume
$\varphi(a)=0$. Since
$$
\limsup_{h\rightarrow 0^+}\frac{\varphi(a+h)-\varphi(a)}{h}\le 0,
$$
there must be some interval $a\le t<\delta$ on which
$\varphi(t)\le \varepsilon (t-a)$.

Let $a\le t<c$ be the largest interval with $c\le b$ such that
$\varphi(t)\le \varepsilon (t-a)$ on $[a,c)$. Then by continuity
$\varphi(t)\le \varepsilon (t-a)$ on the closed interval $[a,c]$.
We claim that $c=b$. Suppose not, then we can find $\delta>0$ such
that $\varphi(t)\le \varepsilon (t-a)$ on $[a,c+\delta]$ since
$$
\limsup_{h\rightarrow 0^+}\frac{\varphi(c+h)-\varphi(c)}{h}\le 0.
$$
This contradicts the choice of the largest interval $[a,c)$.
Therefore, since $\varepsilon>0$ can be arbitrary small, we have
proved $\varphi(t)\le 0$ on $[a,b]$.
\end{pf}

The second lemma below is a general principle on the derivative of
a sup-function which will bridge solutions between ODEs and PDEs.
Let $X$ be a complete smooth manifold and $Y$ be a compact subset
of $X$. Let $\psi(x,t)$ be a smooth function on $X\times[a,b]$ and
let $\varphi (t)=\sup\{\psi(y,t)\ |\  y\in Y\}$. Then it is clear
that $\varphi(t)$ is Lipschitz continuous. We have the following
useful estimate on its derivative.

\begin{lemma} [{Hamilton \cite{Ha86}}]
$$
\frac{d}{dt}\varphi(t) \le \sup\left\{\frac{\partial
\psi}{\partial t}(y,t)\ |\ y\in Y \; \mbox{ satisfies }\;
\psi(y,t)=\varphi(t)\right\}.
$$
\end{lemma}

\begin{pf}
Choose a sequence of times $\{t_j\}$ decreasing to $t$ for which
$$
\lim_{t_j\rightarrow t}\frac{\varphi(t_j)-\varphi(t)}{t_j-t}
=\frac{d\varphi(t)}{dt}.
$$
Since $Y$ is compact, we can choose $y_j\in Y$ with
$\varphi(t_j)=\psi(y_j,t_j)$. By passing to a subsequence, we can
assume $y_j\rightarrow y$ for some $y\in Y$. By continuity, we
have $\varphi(t)=\psi(y,t)$. It follows that $\psi(y_j,t)\le
\psi(y,t)$, and then
\begin{align*}
\varphi(t_j)-\varphi(t)&  \le \psi(y_j,t_j)-\psi(y_j,t)\\[2mm]
&  = \frac{\partial}{\partial t}\psi(y_j,\tilde{t}_j)\cdot(t_j-t)
\end{align*}
for some $\tilde{t}_j\in[t,t_j]$ by the mean value theorem. Thus
we have
$$
\lim_{t_j\rightarrow t}\frac{\varphi(t_j)-\varphi(t)}{t_j-t}
\le\frac{\partial}{\partial t}\psi(y,t).
$$
This proves the result.
\end{pf}

We remark that the above two lemmas are somewhat standard facts in
the theory of PDEs and we have implicitly used them in the
previous sections when we apply the maximum principle. The third
lemma gives a characterization of when a system of ODEs preserve
closed convex sets in Euclidean space. We will use the version
given in \cite{Cl04}. Let $Z\subset \mathbb R^n$ be a closed
convex subset. We define the {\bf tangent cone}\index{tangent
cone} $T_\varphi Z$ to the closed convex set $Z$ at a point
$\varphi\in
\partial Z$ as the smallest closed convex cone with vertex at
$\varphi$ which contains $Z$.

\begin{lemma} [{Hamilton \cite{Ha86}}]
Let $U\subset R^n$ be an open set and $Z\subset U$ be a closed
convex subset. Consider the {\rm ODE} \be
\frac{d\varphi}{dt}=N(\varphi,t)
\ee where $N:\ U\times[0,T]\rightarrow R^n$ is continuous and
Lipschitz in $\varphi$. Then the following two statements are
equivalent.
\begin{itemize}
\item[(i)] For any initial time $t_0\in[0,T]$, any solution of the
{\rm ODE} $(2.3.1)$ which starts in Z at $t_0$ will remain in Z
for all later times; \item[(ii)] $\varphi+N(\varphi,t)\in
T_\varphi Z$ for all $\varphi\in\partial Z$ and $t\in[0,T)$.
\end{itemize}
\end{lemma}

\begin{pf}
We say that a linear function $l$ on $\mathbb R^n$ is a
\textbf{support function}\index{support function} for $Z$ at
$\varphi\in\partial Z$ and write $l\in S_\varphi Z$ if $|l|=1$ and
$l(\varphi)\ge l(\eta)$ for all $\eta \in Z$. Then
$\varphi+N(\varphi,t)\in T_\varphi Z$ if and only if
$l(N(\varphi,t))\le 0$ for all $l\in S_\varphi Z$. Suppose
$l(N(\varphi,t))> 0$ for some $\varphi\in
\partial Z$ and some $l\in S_\varphi Z.$ Then
$$
\frac{d}{dt}l(\varphi)=l\(\frac{d\varphi}{dt}\)=l(N(\varphi,t))>0,
$$
so $l(\varphi)$ is strictly increasing and the solution
$\varphi(t)$ of the ODE (2.3.1) cannot remain in $Z$.

To see the converse, first note that we may assume $Z$ is compact.
This is because we can modify the vector field $N(\varphi,t)$ by
multiplying a cutoff function which is everywhere nonnegative,
equals one on a large ball and equals zero on the complement of a
larger ball. The paths of solutions of the ODE are unchanged
inside the first large ball, so we can intersect $Z$ with the
second ball to make $Z$ convex and compact. If there were a
counterexample before the modification there would still be one
after as we chose the first ball large enough.

Let $s(\varphi)$ be the distance from $\varphi$ to $Z$ in
$\mathbb{R}^n$. Clearly $s(\varphi)=0$ if $\varphi\in Z$. Then
$$
s(\varphi)=\sup\{l(\varphi-\eta)\ |\ \eta\in\partial Z\; \text{
and }\; l\in S_\eta Z\}.
$$
The sup is taken over a compact subset of $\mathbb R^n\times
\mathbb R^n$. Hence by Lemma 2.3.3
$$
\frac{d}{dt}s(\varphi)\le \sup\{l(N(\varphi,t))\ |\
\eta\in\partial Z, l\in S_\eta Z\; \text{ and }\;
s(\varphi)=l(\varphi-\eta)\}.
$$
It is clear that the sup on the RHS of the above inequality can be
takeen only when $\eta$ is the unique closest point in $Z$ to
$\varphi$ and $l$ is the linear function of length one with
gradient in the direction of $\varphi-\eta$. Since $N(\varphi,t)$
is Lipschitz in $\varphi$ and continuous in $t$, we have
$$
|N(\varphi,t)-N(\eta,t)|\le C|\varphi-\eta|
$$
for some constant $C$ and all $\varphi$ and $\eta$ in the compact
set $Z$.

By hypothesis (ii), $$l(N(\eta,t))\le 0,$$ and for the unique
$\eta$, the closest point in $Z$ to $\varphi$,
$$
|\varphi-\eta|=s(\varphi).
$$
Thus
\begin{align*}
\frac{d}{dt}s(\varphi) &  \le \sup \left\{\begin{array}{c}
l(N(\eta,t))+ |l(N(\varphi,t))-l(N(\eta,t))| \quad |\quad
\eta\in\partial Z,\vspace*{3mm}\\
 l\in S_\eta Z,\quad \text{and}\quad
s(\varphi)=l(\varphi-\eta)
\end{array}\right\} \\
&  \le Cs(\varphi).
\end{align*}
Since $s(\varphi)=0$ to start at $t_0$, it follows from Lemma
2.3.2 that $s(\varphi)=0$ for $t\in[t_0,T]$. This proves the
lemma.
\end{pf}

We are now ready to prove Theorem 2.3.1.

\medskip
{\bf\em Proof of Theorem} {\bf 2.3.1.} \ Since the solution
$\sigma(x,t)$ of the (PDE) is uniformly bounded with respect to
the bundle metric $h_{ab}$ on $M\times [t_0,T]$ by hypothesis, we
may assume that $K$ is contained in a tubular neighborhood $V(r)$
of the zero section in $V$ whose intersection with each fiber
$V_x$ is a ball of radius $r$ around the origin measured by the
bundle metric $h_{ab}$ for some large $r>0$.

Recall that $g_{ij}(\cdot,t), t\in[0,T]$, is a smooth solution to
the Ricci flow with uniformly bounded curvature on $M\times
[0,T]$. From Lemma 2.1.1, we have a smooth function $f$ such that
$f\geq 1$ everywhere, $f(x)\rightarrow +\infty$ as
$d_0(x,x_0)\rightarrow +\infty$ for some fixed point $x_0\in M$,
and the first and second covariant derivatives with respect to the
metrics $g_{ij}(\cdot,t)$ are uniformly bounded on $M\times
[0,T]$. Using the metric $h_{ab}$ in each fiber $V_x$ and writing
$|\varphi- \eta|$ for the distance between $\varphi \in V_x$ and
$\eta \in V_x$, we set
$$
s(t)=\sup_{x\in M}\{ \inf\{|\sigma(x,t)-\eta|\  | \ \eta\in
K_x\stackrel{\Delta}{=} K\cap V_x\}-\epsilon e^{At}f(x)\}
$$
where $\epsilon$ is an arbitrarily small positive number and $A$
is a positive constant to be determined. We rewrite the function
$s(t)$ as
$$
s(t)=\sup\{l(\sigma(x,t)-\eta)-\epsilon e^{At}f(x)\ |\  x\in M,
\eta\in \partial K_x  \mbox{ and }  l\in S_\eta K_x\}.
$$
By the construction of the function $f$, we see that the sup is
taken in a compact subset of $M\times V \times V^*$ for all $t$.
Then by Lemma 2.3.3, \be \frac{ds(t)}{dt}\leq \sup\left\{
\frac{\partial}{\partial
t}l(\sigma(x,t)-\eta)-\epsilon Ae^{At}f(x)\right\} 
\ee where the sup is over all $x\in M, \eta\in\partial K_x$ and $
l \in S_\eta K_x$ such that
$$
l(\sigma(x,t)-\eta)-\epsilon e^{At}f(x)=s(t);
$$
in particular we have $|\sigma(x,t)-\eta|=l(\sigma(x,t)-\eta)$,
where $\eta$ is the unique closest point in $K_x$ to
$\sigma(x,t)$, and $l$ is the linear function of length one on the
fiber $V_x$ with gradient in the direction of $\eta$ to
$\sigma(x,t)$. We compute at these $(x,\eta,l)$,
\begin{align}
& \frac{\partial}{\partial t}l(\sigma(x,t)-\eta)-\epsilon A
e^{At}f(x)\\
&  =l\(\frac{\partial \sigma(x,t)}{\partial t}\)
-\epsilon A e^{At}f(x)\nn\\
&  =l(\Delta_t \sigma(x,t))
+l(u^i(x,t)(\nabla_t)_i\sigma(x,t))+l(N(x,\sigma(x,t),t))-\epsilon
A e^{At}f(x).\nn
\end{align} 
By the assumption and Lemma 2.3.4 we have $\eta+N(x,\eta,t)\in
T_\eta K_x$. Hence, for those $(x,\eta,l)$, $l(N(x,\eta,t))\leq 0$
and then
\begin{align}
& l(N(x,\sigma(x,t),t))\\
&  \leq l(N(x,\eta,t))+|N(x,\sigma(x,t),t)-N(x,\eta,t)| \nn\\
&  \leq C|\sigma(x,t)-\eta|=C(s(t)+\epsilon e^{At}f(x)) \nn
\end{align} 
for some positive constant $C$ by the assumption that
$N(x,\sigma,t)$ is Lipschitz in $\sigma$ and the fact that the sup
is taken on a compact set. Thus the combination of
(2.3.2)--(2.3.4) gives \be \frac{ds(t)}{dt}\leq l(\Delta_t
\sigma(x,t)) +l(u^i(x,t)(\nabla_t)_i\sigma(x,t))+Cs(t)+\epsilon
(C-A)e^{At}f(x)
\ee for those $x\in M, \eta \in \partial K_x$ and $l\in  S_\eta
K_x$ such that $ l(\sigma(x,t)-\eta)-\epsilon e^{At}f(x)=s(t)$.

Next we estimate the first two terms of (2.3.5). As we extend a
vector in a bundle from a point $x$ by parallel translation along
geodesics emanating radially out of $x$, we will get a smooth
section of the bundle in some small neighborhood of $x$ such that
all the symmetrized covariant derivatives at $x$ are zero. Now let
us extend $\eta\in V_x$ and $l\in V^*_x$ in this manner.  Clearly,
we continue to have $|l|(\cdot)=1$. Since $K$ is invariant under
parallel translations, we continue to have $\eta (\cdot)\in
\partial K$ and $l(\cdot)$ as a support function for $K$ at $\eta
(\cdot)$. Therefore
$$
l(\sigma(\cdot,t)-\eta(\cdot))-\epsilon e^{At}f(\cdot)\leq s(t)
$$
in the neighborhood. It follows that the function
$l(\sigma(\cdot,t)-\eta(\cdot))-\epsilon e^{At}f(\cdot)$ has a
local maximum at $x$, so at $x$
\begin{align*}
(\nabla _t)_i(l(\sigma(x,t)-\eta)-\epsilon e^{At}f(x))&=0,\\
\text{and }\; \Delta_t(l(\sigma(x,t)-\eta)-\epsilon e^{At}f(x))&
\leq 0.
\end{align*}
Hence at $x$
\begin{align*}
l((\nabla _t)_i\sigma(x,t))-\epsilon e^{At}(\nabla_t)_if(x)&=0,\\
\text{and }\; l(\Delta_t \sigma(x,t))-\epsilon e^{At}\Delta_t
f(x)&\leq 0.
\end{align*}
Therefore by combining with (2.3.5), we have
\begin{align*}
\frac{d}{dt}s(t)&  \leq Cs(t)+\epsilon
(\Delta_tf(x)+u^i(\nabla_t)_if(x)+(C-A)f(x))e^{At} \\
&  \leq Cs(t)
\end{align*}
for $A>0$ large enough, since $f(x)\geq 1$ and the first and
second covariant derivatives of $f$ are uniformly bounded on
$M\times[0,T]$. So by applying Lemma 2.3.2 and the arbitrariness
of $\epsilon$, we have completed the proof of Theorem 2.3.1.
\endproof

Finally, we would like to state a useful generalization of Theorem
2.3.1 by Chow and Lu in \cite{Cl04} which allows the set $K$ to
depend on time. One can consult the paper \cite{Cl04} for the
proof.

\begin{theorem}[{Chow and Lu \cite{Cl04}}]
Let $K(t)\subset V$, $t\in [0,T]$ be closed subsets which satisfy
the following hypotheses
\begin{itemize}
\item[(H3)] $K(t)$ is invariant under parallel translation defined
by the connection $\nabla_t$ for each $t\in [0,T]$; \item[(H4)] in
each fiber $V_x$, the set $K_x(t)\stackrel{\Delta}{=}K(t)\cap V_x$
is nonempty, closed and convex for each $t\in [0,T]$; \item[(H5)]
the space-time track $\bigcup\limits_{t\in[0,T]}(\partial
K(t)\times\{t\})$ is a closed subset of $V\times[0,T]$.
\end{itemize}
Suppose that, for any $x\in M$ and any initial time $t_0\in[0,T)$,
and for any solution $\sigma_x(t)$ of the {\rm (ODE)} which starts
in $K_x(t_0)$, the solution $\sigma_x(t)$ will remain in $K_x(t)$
for all later times. Then for any initial time $t_0\in[0,T)$ the
solution $\sigma(x,t)$ of the {\rm (PDE)} will remain in $K(t)$
for all later times if $\sigma(x,t)$ starts in $K(t_0)$ at time
$t_0$ and the solution $\sigma(x,t)$ is uniformly bounded with
respect to the bundle metric $h_{ab}$ on $M\times[t_0,T]$.
\end{theorem}

\section{Hamilton-Ivey Curvature Pinching Estimate}

The Hamilton-Ivey curvature pinching estimate \cite{Ha95F, Iv}
roughly says that if a solution to the Ricci flow on a
three-manifold becomes singular (i.e., the curvature goes to
infinity) as time $t$ approaches the maximal time $T$, then the
most negative sectional curvature will be small compared to the
most positive sectional curvature. This pinching estimate plays a
crucial role in analyzing the formation of singularities in the
Ricci flow on three-manifolds. The proof here is based on the
argument in Hamilton \cite{Ha95F}. The estimate was later improved
by Hamilton \cite{Ha99} which will be presented in Section 5.3
(see Theorem 5.3.2).

Consider a complete solution to the Ricci flow
$$
\frac{\partial}{\partial t}g_{ij}=-2R_{ij}
$$
on a complete three-manifold with bounded curvature in space for
each time $t\geq 0$. Recall from Section 1.3 that the evolution
equation of the curvature operator $M_{\alpha\beta}$ is given by
\be \frac{\partial}{\partial t}M_{\alpha\beta}=\Delta
M_{\alpha\beta}+M_{\alpha\beta}^2+M_{\alpha\beta}^\# 
\ee where $M_{\alpha\beta}^2$ is the operator square
$$
M_{\alpha\beta}^2=M_{\alpha\gamma}M_{\beta\gamma}
$$
and $M_{\alpha\beta}^\#$ is the Lie algebra $so(n)$ square
$$
M_{\alpha\beta}^\#=C_\alpha^{\gamma\zeta}
C_\beta^{\eta\theta}M_{\gamma\eta}M_{\zeta\theta}.
$$
In dimension $n=3$, we know that $M_{\alpha\beta}^\#$ is the
adjoint matrix of $M_{\alpha\beta}$. If we diagonalize
$M_{\alpha\beta}$ with eigenvalues $\lambda \ge \mu \ge \nu$ so
that
$$(M_{\alpha\beta})=\left(
\begin{array}{ccc}\lambda&  \ &  \ \\
\ &  \mu &  \ \\
\ &  \ &  \nu
\end{array}
\right),
$$
then $M_{\alpha\beta}^2$ and $M_{\alpha\beta}^\#$ are also
diagonal, with
$$
(M_{\alpha\beta}^2)=\left(
\begin{array}{ccc}\lambda^2&  \ &  \ \\
\ &  \mu^2 &  \ \\
\ &  \ &  \nu^2
\end{array}
\right) \ \ \text{and} \ \ (M_{\alpha\beta}^\#)=\left(
\begin{array}{ccc}\mu\nu&  \ &  \ \\
\ &  \lambda\nu &  \ \\
\ &  \ &  \lambda\mu
\end{array}
\right).
$$

Thus the ODE corresponding to PDE (2.4.1) for $M_{\alpha\beta}$
(in the space of $3\times 3 $ matrices) is given by the following
system \be
      \left\{
       \begin{array}{lll}
  \frac{d}{dt}\lambda=\lambda^2+\mu\nu,\\[4mm]
  \frac{d}{dt}\mu=\mu^2+\lambda\nu,\\[4mm]
  \frac{d}{dt}\nu=\nu^2+\lambda\mu.
       \end{array}
    \right.
\ee

Let $P$ be the principal bundle of the manifold and form the
associated bundle $V=P\times_GE$, where $G=O(3)$ and $E$ is the
vector space of symmetric bilinear forms on $so(3)$. The curvature
operator $M_{\alpha\beta}$ is a smooth section of $V=P\times_GE$.
According to Theorem 2.3.1, any closed convex set of curvature
operator matrices $M_{\alpha\beta}$ which is $O(3)$-invariant (and
hence invariant under parallel translation) and preserved by ODE
(2.4.2) is also preserved by the Ricci flow.

We are now ready to state and prove the \textbf{Hamilton-Ivey
pinching estimate} \index{Hamilton-Ivey pinching estimate}.

\begin{theorem}[{Hamilton \cite{Ha95F}, Ivey \cite{Iv}}]
Suppose we have a solution to the Ricci flow on a three-manifold
which is complete with bounded curvature for each $t\geq0$. Assume
at $t=0$ the eigenvalues $\lambda\geq\mu\geq\nu$ of the curvature
operator at each point are bounded below by $\nu\geq-1$. The
scalar curvature $R=\lambda+\mu+\nu$ is their sum. Then at all
points and all times $t\geq0$ we have the pinching estimate
$$
R\geq (-\nu)[\log(-\nu)-3],
$$
whenever $\nu<0$.
\end{theorem}

\begin{pf} The proof is taken from Hamilton \cite{Ha95F}.
Consider the function
$$
y=f(x)=x(\log x-3)
$$
defined on $e^2\leq x<+\infty$. It is easy to check that $f$ is
increasing and convex  with range $-e^2\leq y<+\infty$. Let
$f^{-1}(y)=x$ be the inverse function, which is also increasing
but concave and satisfies \be
\lim_{y\rightarrow\infty}\frac{f^{-1}(y)}{y}=0 
\ee Consider also the set $K$ of matrices $M_{\alpha\beta}$
defined by the inequalities \be K:
   \left  \{
     \begin{array}{lll}
         \lambda+\mu+\nu\geq-3,\\[4mm]
         \nu+f^{-1}(\lambda+\mu+\nu)\geq0.
     \end{array}
     \right. 
\ee By Theorem 2.3.1 and the assumptions in Theorem 2.4.1 at
$t=0$, we only need to check that the set $K$ defined above is
closed, convex and preserved by the ODE (2.4.2).

Clearly $K$ is closed because $f^{-1}$ is continuous.
$\lambda+\mu+\nu$ is just the trace function of $3\times 3$
matrices which is a linear function. Hence the first inequality in
(2.4.4) defines a linear half-space, which is convex. The function
$\nu$ is the least eigenvalue function, which is concave. Also
note that $f^{-1}$ is concave. Thus the second inequality in
(2.4.4) defines a convex set as well. Therefore we proved $K$ is
closed and convex.

Under the ODE (2.4.2)
\begin{align*}
    \frac{d}{dt}(\lambda+\mu+\nu)
&  = \lambda^2+\mu^2+\nu^2+\lambda\mu+\lambda\nu+\mu\nu\\
&  = \frac{1}{2}[(\lambda+\mu)^2+(\lambda+\nu)^2+(\mu+\nu)^2]\\
&  \geq 0.
\end{align*}
Thus the first inequality in (2.4.4) is preserved by the ODE.

The second inequality in (2.4.4) can be written as
$$
\lambda+\mu+\nu\geq f(-\nu),\ \ \ \text{whenever}\ \nu\leq-e^2,
$$
which becomes \be \lambda+\mu\geq (-\nu)[\log(-\nu)-2],\ \ \
\text{whenever}\
\nu\leq-e^2.
\ee To show the inequality is preserved we only need to look at
points on the boundary of the set. If
$\nu+f^{-1}(\lambda+\mu+\nu)=0$ then
$\nu=-f^{-1}(\lambda+\mu+\nu)\le -e^2$ since $f^{-1}(y)\ge e^2$.
Hence the RHS of (2.4.5) is nonnegative. We thus have
$\lambda\geq0$ because $\lambda\geq\mu$. But $\mu$ may have either
sign. We split our consideration into two cases:

\medskip
{\it Case} (i): $\mu\geq0$.

\smallskip
We need to verify
$$
\frac{d\lambda}{dt}+\frac{d\mu}{dt}\geq(\log(-\nu)-1)\frac{d(-\nu)}{dt}
$$
when $\lambda+\mu=(-\nu)[\log(-\nu)-2]$. Solving for
$$
\log(-\nu)-2=\frac{\lambda+\mu}{(-\nu)}
$$
and substituting above, we must show
$$
\lambda^2+\mu\nu+\mu^2+\lambda\nu
\geq\(\frac{\lambda+\mu}{(-\nu)}+1\)(-\nu^2-\lambda\mu)
$$
which is equivalent to
$$
(\lambda^2+\mu^2)(-\nu)+\lambda\mu(\lambda+\mu+(-\nu))+(-\nu)^3\geq
0.
$$
Since $\lambda,\mu$ and $(-\nu)$ are all nonnegative we are done
in the first case.

\medskip
{\it Case} (ii): $\mu<0$.

\smallskip
We need to verify
$$
\frac{d\lambda}{dt}\geq\frac{d(-\mu)}{dt}
+(\log(-\nu)-1)\frac{d(-\nu)}{dt}
$$
when $\lambda=(-\mu)+(-\nu)[\log(-\nu)-2]$. Solving for
$$
\log(-\nu)-2=\frac{\lambda-(-\mu)}{(-\nu)}
$$
and substituting above, we need to show
$$
\lambda^2+\mu\nu\geq-\mu^2-\lambda\nu
+\(\frac{\lambda-(-\mu)}{(-\nu)}+1\)(-\nu^2-\lambda\mu)
$$
or
$$
\lambda^2+(-\mu)(-\nu)\geq\lambda(-\nu)-(-\mu)^2
+\(\frac{\lambda-(-\mu)}{(-\nu)}+1\)(\lambda(-\mu)-(-\nu)^2)
$$
which reduces to
$$
\lambda^2(-\nu)+\lambda(-\mu)^2+(-\mu)^2(-\nu)+(-\nu)^3
\geq\lambda^2(-\mu)+\lambda(-\mu)(-\nu)
$$
or equivalently
$$
(\lambda^2-\lambda(-\mu)+(-\mu)^2)((-\nu)-(-\mu))
+(-\mu)^3+(-\nu)^3\geq0.
$$
Since $\lambda^2-\lambda(-\mu)+(-\mu)^2\geq0\ \ \text{and}\ \
(-\nu)-(-\mu)\geq0$ we are also done in the second case.

Therefore the proof is completed.
\end{pf}

\section{Li-Yau-Hamilton Estimates}

In \cite{LY}, Li-Yau developed a fundamental gradient estimate,
now called Li-Yau estimate, for positive solutions to the heat
equation on a complete Riemannian manifold with nonnegative Ricci
curvature. They used it to derive the Harnack inequality for such
solutions by path integration. Then based on the suggestion of
Yau, Hamilton \cite{Ha88} developed a similar estimate for the
scalar curvature of solutions to the Ricci flow on a Riemann
surface with positive curvature, and later obtained a matrix
version of the Li-Yau estimate \cite{Ha93} for solutions to the
Ricci flow with positive curvature operator in all dimensions.
This matrix version of the Li-Yau estimate is the
\textbf{Li-Yau-Hamilton estimate}\index{Li-Yau-Hamilton estimate},
which we will present in this section. Most of the presentation
follows the original papers of Hamilton \cite{Ha88, Ha93, Ha93E}.

We have seen that in the Ricci flow the curvature tensor satisfies
a nonlinear heat equation, and the nonnegativity of the curvature
operator is preserved by the Ricci flow. Roughly speaking, the
Li-Yau-Hamilton estimate says the nonnegativity of a certain
combination of the derivatives of the curvature up to second order
is also preserved by the Ricci flow. This estimate plays a central
role in the analysis of formation of singularities and the
application of the Ricci flow to three-manifold topology.

Let us begin by describing the Li-Yau estimate \cite{LY} for
positive solutions to the heat equation on a complete Riemannian
manifold with nonnegative Ricci curvature.

\begin{theorem}[{Li-Yau \cite{LY}}]
Let $(M, g_{ij})$ be an $n$-dimensional complete Riemannian
manifold with nonnegative Ricci curvature. Let $u(x,t)$ be any
positive solution to the heat equation
$$
\frac{\partial u}{\partial t}=\Delta u \qquad {\text{on}} \ \
M\times [0,\infty).
$$
Then we have \be \frac{\partial u}{\partial t}-\frac{|\nabla
u|^2}{u}+\frac{n}{2t}u\ge 0 \qquad \text{on } \ \ M\times
(0,\infty).
\ee
\end{theorem}

We remark that, as observed by Hamilton (cf. \cite{Ha93}), one can
in fact prove that for any vector field $V^i$ on $M$, \be
\frac{\partial u}{\partial t}+2\nabla u\cdot V+u|V|^2
+\frac{n}{2t}u\ge 0. 
\ee If we take the optimal vector field $V=-\nabla u/u$, we
recover the inequality (2.5.1).

Now we consider the Ricci flow on a Riemann surface. Since in
dimension two the Ricci curvature is given by
$$
R_{ij}=\frac{1}{2} Rg_{ij},
$$
the Ricci flow (1.1.5) becomes \be
\frac{\partial g_{ij}}{\partial t}=-R g_{ij}.
\ee

Now let $g_{ij}(x,t)$ be a complete solution of the Ricci flow
(2.5.3) on a Riemann surface $M$ and $0\le t<T$. Then the scalar
curvature $ R(x,t)$ evolves by the semilinear equation
$$
\frac{\partial R}{\partial t}=\triangle R+R^{2}
$$
on $M\times[0,T)$. Suppose the scalar curvature of the initial
metric is bounded, nonnegative everywhere and positive somewhere.
Then it follows from Proposition 2.1.2 that the scalar curvature
$R(x,t)$ of the evolving metric remains nonnegative. Moreover,
from the standard strong maximum principle (which works in each
local coordinate neighborhood), the scalar curvature is positive
everywhere for $t>0$. In \cite{Ha88}, Hamilton obtained the
following Li-Yau estimate for the scalar curvature $R(x,t)$.

\begin{theorem}[{Hamilton \cite{Ha88}}]
Let $g_{ij}(x,t)$ be a complete solution of the Ricci flow on a
surface $M$. Assume the scalar curvature of the initial metric is
bounded, nonnegative everywhere and positive somewhere. Then the
scalar curvature $R(x,t)$ satisfies the Li-Yau estimate \be
\frac{\partial R}{\partial t} - \frac{|\nabla R|^2}{R} +
\frac{R}{t} \geq 0.
\ee
\end{theorem}

\begin{pf}
By the above discussion, we know $R(x,t)>0$ for $t>0$. If we set
$$
L=\log R(x,t)\quad \text{for}\quad t>0,
$$
then
\begin{align*}
\frac{\partial}{\partial t}L&  = \frac{1}{R}(\triangle
R+R^{2})\\
&  = \triangle L+|\nabla L|^{2}+R
\end{align*}
and (2.5.4) is equivalent to
$$
\frac{\partial L}{\partial t}-|\nabla L|^{2}+\frac{1}{t}
=\triangle L +R +\frac{1}{t}\geq0.
$$

Following Li-Yau \cite{LY} in the linear heat equation case, we
consider the quantity \be Q =\frac{\partial L}{\partial t}-|\nabla
L|^{2}=\triangle L+R.
\ee Then by a direct computation,
\begin{align*}
\frac{\partial Q}{\partial t}
&  = \frac{\partial}{\partial t}(\triangle L+R)\\
&  = \triangle\(\frac{\partial L}{\partial t}\)+R\triangle
L+\frac{\partial R}{\partial t}\\
&  = \triangle Q+2\nabla L\cdot\nabla Q+2|\nabla ^{2}L|^{2}
+2R(\triangle L)+R^{2}\\
&  \geq \triangle Q+2\nabla L\cdot\nabla Q+Q^{2}.
\end{align*}
So we get
$$
\frac{\partial}{\partial t}
\(Q+\frac{1}{t}\)\geq\triangle\(Q+\frac{1}{t}\)+2\nabla
L\cdot\nabla \(Q+\frac{1}{t}\)+\(Q-\frac{1}{t}\)\(Q+\frac{1}{t}\).
$$
Hence by a similar maximum principle argument as in the proof of
Lemma 2.1.3, we obtain
$$
Q+\frac{1}{t}\ge 0.
$$
This proves the theorem.
\end{pf}

As an immediate consequence, we obtain the following Harnack
inequality for the scalar curvature $R$ by taking the Li-Yau type
path integral as in \cite{LY}.

\begin{corollary}[{Hamilton \cite{Ha88}}]
Let $g_{ij}(x,t)$ be a complete solution of the Ricci flow on a
surface with bounded and nonnegative scalar curvature. Then for
any points $x_{1},x_{2}\in M$, and $0<t_{1}<t_{2}$, we have
$$
R(x_2,t_2)
\geq\frac{t_1}{t_2}e^{-d_{t_1}{(x_1,x_2)}^2/{4(t_2-t_1)}}R(x_1,t_1).
$$
\end{corollary}

\begin{pf}
Take the geodesic path $\gamma(\tau)$, $\tau\in[t_1,t_2],$ from
$x_1$ to $ x_2$ at time $t_1$ with constant velocity
$d_{t_1}(x_1,x_2)/(t_2-t_1).$ Consider the space-time path
$\eta(\tau)=(\gamma(\tau),\tau)$, $\tau\in[t_1,t_2]$. We compute
\begin{align*}
\log\frac{R(x_2,t_2)}{R(x_1,t_1)}
&  = \int^{t_2}_{t_1}\frac{d}{d\tau}L(\gamma(\tau),\tau)d\tau\\
&  = \int^{t_2}_{t_1}\frac{1}{R}\(\frac{\partial R}{\partial \tau}
+\nabla R \cdot \frac{d\gamma}{d\tau}\)d\tau\\
&  \geq \int^{t_2}_{t_1}\(\frac{\partial L}{\partial\tau} -|\nabla
L|^2_{g_{ij}(\tau)}
-\frac{1}{4}\left|\frac{d\gamma}{d\tau}\right|^2_{g_{ij}(\tau)}\)d\tau.
\end{align*}
Then by Theorem 2.5.2 and the fact that the metric is shrinking
(since the scalar curvature is nonnegative), we have
\begin{align*}
\log\frac{R(x_2,t_2)}{R(x_1,t_1)} &
\geq\int^{t_2}_{t_1}\(-\frac{1}{\tau}-
\frac{1}{4}\left|\frac{d\gamma}{d\tau}\right|^2_{g_{ij}(\tau)}\)d\tau \\
&  = \log\frac{t_1}{t_2}-\frac{d_{t_1}(x_1,x_2)^2}{4(t_2-t_1)}
\end{align*}
After exponentiating above, we obtain the desired Harnack
inequality.
\end{pf}

To prove a similar inequality as (2.5.4) for the scalar curvature
of solutions to the Ricci flow in higher dimensions is not so
simple. First of all, we will need to require nonnegativity of the
curvature operator (which we know is preserved under the Ricci
flow).  Secondly, one does not get inequality (2.5.4) directly,
but rather indirectly as the trace of certain matrix estimate. The
key ingredient in formulating this matrix version is to derive
some identities from the soliton solutions and prove an elliptic
inequality based on these quantities. Hamilton found such a
general principle which was based on the idea of Li-Yau \cite{LY}
when an identity is checked on the heat kernel before an
inequality was found. To illustrate this point, let us first
examine the heat equation case. Consider the heat kernel
$$
u(x,t)=(4\pi t)^{-n/2} e^{-|x|^2/4t}
$$
for the standard heat equation on $\mathbb R^n$ which can be
considered as an expanding soliton solution.

Differentiating the function $u$, we get \be \nabla_j u=-u
\frac{x_j}{2t} \; \text{ or }\;
\nabla_j u+uV_j=0,
\ee where
$$
V_j=\frac{x_j}{2t}=-\frac{\nabla_j u}{u}.
$$

Differentiating (2.5.6), we have \be \nabla_i\nabla_j
u+\nabla_iuV_j+\frac{u}{2t}\delta_{ij}=0.
\ee To make the expression in (2.5.7) symmetric in $i, j$, we
multiply $V_i$ to (2.5.6) and add to (2.5.7) and obtain \be
\nabla_i\nabla_j u+\nabla_iuV_j+\nabla_juV_i+uV_iV_j
+\frac{u}{2t}\delta_{ij}=0. 
\ee Taking the trace in (2.5.8) and using the equation ${\partial
u}/{\partial t}=\Delta u$, we arrive at
$$
\frac{\partial u}{\partial t} +2\nabla u\cdot V+u|V|^2
+\frac{n}{2t} u=0,
$$
which shows that the Li-Yau inequality (2.5.1) becomes an equality
on our expanding soliton solution $u$! Moreover, we even have the
matrix identity (2.5.8).

Based on the above observation and using a similar process,
Hamilton found a matrix quantity, which vanishes on expanding
gradient Ricci solitons and is nonnegative for any solution to the
Ricci flow with nonnegative curvature operator. Now we describe
the process of finding the Li-Yau-Hamilton quadratic for the Ricci
flow in arbitrary dimension.

Consider a homothetically expanding gradient soliton $g$, we have
\be
R_{ab}+\frac{1}{2t}g_{ab}=\nabla_aV_b 
\ee in the orthonormal frame coordinate chosen as in Section 1.3.
Here $V_b=\nabla_bf$ for some function $f$. Differentiating
(2.5.9) and commuting give the first order relations
\begin{align}
\nabla_aR_{bc}-\nabla_bR_{ac}
&  = \nabla_a\nabla_bV_c-\nabla_b\nabla_aV_c\\
&  = R_{abcd}V_d,\nn
\end{align}   
and differentiating again, we get
\begin{align*}
\nabla_a\nabla_bR_{cd}-\nabla_a\nabla_cR_{bd}
&  = \nabla_a(R_{bcde}V_e)\\
&  = \nabla_aR_{bcde}V_e+R_{bcde}\nabla_aV_e\\
&  = \nabla_aR_{bcde}V_e+R_{ae}R_{bcde}+\frac{1}{2t}R_{bcda}.
\end{align*}
We further take the trace of this on $a$ and $b$ to get
$$
\Delta R_{cd}-\nabla_a\nabla_cR_{ad}-R_{ae}R_{acde}
+\frac{1}{2t}R_{cd}-\nabla_aR_{acde}V_e=0,
$$
and then by commuting the derivatives and second Bianchi identity,
$$
\Delta R_{cd}-\frac{1}{2}\nabla_c\nabla_dR
+2R_{cade}R_{ae}-R_{ce}R_{de}+\frac{1}{2t}R_{cd}+(\nabla
_eR_{cd}-\nabla_dR_{ce})V_e=0.
$$
Let us define
\begin{align*}
M_{ab}&  = \Delta R_{ab}-\frac{1}{2}\nabla_a\nabla_bR
+2R_{acbd}R_{cd}-R_{ac}R_{bc}+\frac{1}{2t}R_{ab},\\
P_{abc}&  = \nabla_aR_{bc}-\nabla_bR_{ac}.
\end{align*}
Then \be
M_{ab}+P_{cba}V_c=0,
\ee We rewrite (2.5.10) as
$$
P_{abc}=R_{abcd}V_d
$$
and then \be
P_{cab}V_c+R_{acbd}V_cV_d=0. 
\ee Adding (2.5.11) and (2.5.12) we have
$$
M_{ab}+(P_{cab}+P_{cba})V_c+R_{acbd}V_cV_d=0
$$
and then
$$
M_{ab}W_aW_b+(P_{cab}+P_{cba})W_aW_bV_c+R_{acbd}W_aV_cW_bV_d=0.
$$
If we write
$$
U_{ab}=\frac{1}{2}(V_aW_b-V_bW_a)=V\wedge W,
$$
then the above identity can be rearranged as \be
Q\stackrel{\Delta}{=}M_{ab}W_aW_b+2P_{abc}U_{ab}W_c
+R_{abcd}U_{ab}U_{cd}=0.
\ee This is the \textbf{Li-Yau-Hamilton quadratic}
\index{Li-Yau-Hamilton quadratic} we look for. Note that the proof
of the Li-Yau-Hamilton estimate below does not depend on the
existence of such an expanding gradient Ricci soliton. It is only
used as inspiration.

Now we are ready to state the remarkable \textbf{Li-Yau-Hamilton
estimate} \index{Li-Yau-Hamilton estimate} for the Ricci flow.

\begin{theorem}[{Hamilton \cite{Ha93}}]
Let $g_{ij}(x,t)$ be a complete solution with bounded curvature to
the Ricci flow on a manifold $M$ for $t$ in some time interval
$(0,T)$ and suppose the curvature operator of $g_{ij}(x,t)$ is
nonnegative. Then for any one-form $W_a$ and any two-form $U_{ab}$
we have
$$
M_{ab}W_aW_b+2P_{abc}U_{ab}W_c+R_{abcd}U_{ab}U_{cd}\geq 0
$$
on $M\times (0,T)$.
\end{theorem}

The proof of this theorem requires some rather intense
calculations. Here we only give a sketch of the proof. For more
details, we refer the reader to Hamilton's original paper
\cite{Ha93}.

\medskip
{\bf\em Sketch of the Proof.} \  Let $g_{ij}(x,t)$ be the complete
solution with bounded and nonnegative curvature operator. Recall
that in the orthonormal frame coordinate system, the curvatures
evolve by
$$
\left\{
\begin{array}{lll}
\frac{\partial}{\partial t}R_{abcd}
=\Delta R_{abcd}+2(B_{abcd}-B_{abdc}-B_{adbc} +B_{acbd}),\\[4mm]
\frac{\partial}{\partial t}R_{ab}
=\Delta R_{ab}+2R_{acbd}R_{cd},\\[4mm]
\frac{\partial}{\partial t}R=\Delta R+ 2|Ric|^2,
       \end{array}
    \right.
$$
where $B_{abcd}=R_{aebf}R_{cedf}.$

By a long but straightforward computation from these evolution
equations, one can get
$$
\(\frac{\partial}{\partial t}-\Delta\)P_{abc}=2R_{adbe}P_{dec}
+2R_{adce}P_{dbe}+2R_{bdce}P_{ade}-2R_{de}\nabla_dR_{abce}
$$
and
\begin{align*}
\(\frac{\partial}{\partial t}-\Delta\)M_{ab}
&  =2R_{acbd}M_{cd}+2R_{cd}(\nabla_cP_{dab}+\nabla_cP_{dba}) \\
&\quad +2P_{acd}P_{bcd}-4P_{acd}P_{bdc}+2R_{cd}R_{ce}R_{adbe}
-\frac{1}{2t^2}R_{ab}.
\end{align*}
Now consider
$$
Q\stackrel{\Delta}{=}M_{ab}W_aW_b+2P_{abc}U_{ab}W_c
+R_{abcd}U_{ab}U_{cd}.
$$
At a point where \be \(\frac{\partial}{\partial
t}-\Delta\)W_a=\frac{1}{t}W_a, \ \
 \(\frac{\partial}{\partial t}-\Delta\)U_{ab}=0,
\ee and \be \nabla_aW_b=0,\ \
\nabla_aU_{bc}=\frac{1}{2}(R_{ab}W_c-R_{ac}W_b)
+\frac{1}{4t}(g_{ab}W_c-g_{ac}W_b),
\ee we have
\begin{align}
\(\frac{\partial}{\partial t}-\Delta\)Q
&  = 2R_{acbd}M_{cd}W_aW_b-2P_{acd}P_{bdc}W_aW_b \\
&\quad  + 8R_{adce}P_{dbe}U_{ab}W_c+4R_{aecf}R_{bedf}U_{ab}U_{cd}\nn\\
&\quad  +
(P_{abc}W_c+R_{abcd}U_{cd})(P_{abe}W_e+R_{abef}U_{ef}).\nn
\end{align} 

For simplicity we assume the manifold is compact and the curvature
operator is strictly positive. (For the general case we shall mess
the formula up a bit to sneak in the term $\epsilon e^{At}f$, as
done in Lemma 2.1.3). Suppose not; then there will be a first time
when the quantity $Q$ is zero, and a point where this happens, and
a choice of $U$ and $W$ giving the null eigenvectors. We can
extend $U$ and $W$ any way we like in space and time and still
have $Q\geq$0, up to the critical time. In particular we can make
the first derivatives in space and time to be anything we like, so
we can extend first in space to make (2.5.15) hold at that point.
And then, knowing $\Delta W_a$ and $\Delta U_{ab}$, we can extend
in time to make (2.5.14) hold at that point and that moment. Thus
we have (2.5.16) at the point.

In the RHS of (2.5.16) the quadratic term
$$
(P_{abc}W_c+R_{abcd}U_{cd})(P_{abe}W_e+R_{abef}U_{ef})
$$
is clearly nonnegative. By similar argument as in the proof of
Lemma 2.1.3, to get a contradiction we only need to show the
remaining part in the RHS of (2.5.16) is also nonnegative.

A nonnegative quadratic form can always be written as a sum of
squares of linear forms. This is equivalent to diagonalizing a
symmetric matrix and writing each nonnegative eigenvalue as a
square. Write
$$
Q=\sum_{k}(X^k_aW_a+Y^k_{ab}U_{ab})^2,\ \ \ \ \(1\leq k\leq
n+\frac{n(n-1)}{2}\).
$$
This makes
$$
M_{ab}=\sum_kX^k_aX^k_b, \ \ \ P_{abc}=\sum_kY^k_{ab}X^k_c
$$
and
$$
R_{abcd}=\sum_kY^k_{ab}Y^k_{cd}.
$$
It is then easy to compute {\allowdisplaybreaks
\begin{align*}
& 2R_{acbd}M_{cd}W_aW_b-2P_{acd}P_{bdc}W_aW_b+8R_{adce}P_{dbe}U_{ab}W_e\\
&\quad+4R_{aecf}R_{bedf}U_{ab}U_{cd}\\
&  = 2\(\sum_kY^k_{ac}Y^k_{bd}\)\(\sum_lX^l_{a}Y^l_{c}\)W_aW_b \\
&\quad-2\(\sum_kY^k_{ac}X^k_{d}\)\(\sum_lY^l_{bd}X^l_{c}\)W_aW_b\\
&\quad +8\(\sum_kY^k_{ad}Y^k_{ce}\)\(\sum_lY^l_{db}X^l_{e}\)U_{ab}W_c \\
&\quad+4\(\sum_kY^k_{ae}Y^k_{cf}\)\(\sum_lY^l_{be}Y^l_{df}\)U_{ab}U_{cd}\\
&  = \sum_{k,l}(Y^k_{ac}X^l_cW_a-Y^l_{ac}X^k_cW_a
-2Y^k_{ac}Y^l_{bc}U_{ab})^2\\
&  \geq 0.
\end{align*}
} This says that the remaining part in the RHS of (2.5.16) is also
nonnegative. Therefore we have completed the sketch of the proof.
\endproof

By taking $U_{ab}=\frac{1}{2}(V_aW_b-V_bW_a)$ and tracing over
$W_a$, we immediately get \vskip 0.2cm

\begin{corollary}[{Hamilton \cite{Ha93}}]
For any one-form $V_a$ we have
$$
\frac{\partial R}{\partial t}+\frac{R}{t} +2\nabla_aR\cdot
V_a+2R_{ab}V_aV_b\geq 0.
$$
\end{corollary}

In particular by taking $V\equiv 0$, we see that the function
$tR(x,t)$ is pointwise nondecreasing in time. By combining this
property with the local derivative estimate of curvature, we have
the following elliptic type estimate.

\begin{corollary}[{Hamilton \cite{Ha95F}}]
Suppose we have a solution to the Ricci flow for $t>0$ which is
complete with bounded curvature, and has nonnegative curvature
operator. Suppose also that at some time $t>0$ we have the scalar
curvature $R\leq M$ for some constant M in the ball of radius $r$
around some point $p$. Then for $k=1,2,\ldots$, the $k^{th}$ order
derivatives of the curvature at $p$ at the time $t$ satisfy a
bound
$$
|\nabla^k Rm(p,t)|^2\leq
C_kM^2\(\frac{1}{r^{2k}}+\frac{1}{t^k}+M^k\)
$$
for some constant $C_k$ depending only on the dimension and $k$.
\end{corollary}

\begin{pf}
Since $tR$ is nondecreasing in time, we get a bound $R\leq 2M$ in
the given region for times between $t/2$ and $t$. The nonnegative
curvature hypothesis tells us the metric is shrinking. So we can
apply the local derivative estimate in Theorem 1.4.2 to deduce the
result.
\end{pf}

By a similar argument as in Corollary 2.5.3, one readily has the
following Harnack inequality.

\begin{corollary} [{Hamilton \cite{Ha93}}]
Let $g_{ij}(x,t)$ be a complete solution of the Ricci flow on a
manifold with bounded and nonnegative curvature operator, and let
$x_1,x_2\in M,\ 0<t_1<t_2.$ Then the following inequality holds
$$
R(x_2,t_2)\geq\frac{t_1}{t_2}e^{-d_{t_1}(x_1,x_2)^2/2(t_2-t_1)}\cdot
R(x_1,t_1).$$
\end{corollary}

In the above discussion, we assumed that the solution to the Ricci
flow exists on $0\le t<T$, and we derived the Li-Yau-Hamilton
estimate with terms $1/t$ in it. When the solution happens to be
{\bf ancient}\index{solution!ancient}\index{ancient!solution},
i.e., defined on $-\infty<t<T$, Hamilton \cite{Ha93} found an
interesting and simple procedure for getting rid of them. Suppose
we have a solution on $\alpha<t<T$ we can replace $t$ by
$t-\alpha$ in the Li-Yau-Hamilton estimate. If we let $\alpha\to
-\infty$, then the expression $1/(t-\alpha)\to 0$ and disappears!
In particular the trace Li-Yau-Hamilton estimate in Corollary
2.5.5 becomes \be \frac{\partial R}{\partial t}+2\nabla_aR\cdot
V_a
+2R_{ab}V_aV_b\geq 0.
\ee By taking $V=0$, we see that $\frac{\partial R}{\partial
t}\geq 0$. Thus, we have the following

\begin{corollary}[{Hamilton \cite{Ha93}}]
Let $g_{ij}(x,t)$ be a complete ancient solution of the Ricci flow
on $M\times (-\infty, T)$ with bounded and nonnegative curvature
operator, then the scalar curvature $R(x,t)$ is pointwise
nondecreasing in time $t$.
\end{corollary}

Corollary 2.5.8  will be very useful later on when we study
ancient $\kappa$-solutions in Chapter 6, especially combined with
Shi's derivative estimate.

We end this section by stating the Li-Yau-Hamilton estimate for
the K\"ahler-Ricci flow, due to the first author \cite{Cao92},
under the weaker curvature assumption of nonnegative holomorphic
bisectional curvature. Note that the following Li-Yau-Hamilton
estimate in the K\"ahler case is really a Li-Yau-Hamilton estimate
for the Ricci tensor of the evolving metric, so not only can we
derive an estimate on the scalar curvature, which is the trace of
the Ricci curvature, similar to Corollary 2.5.5 but also an
estimate on the determinant of the Ricci curvature as well.

\begin{theorem}[{Cao \cite{Cao92}}]
Let $g_{\alpha\bar \beta}(x,t)$ be a complete solution to the
K\"ahler-Ricci flow on a complex  manifold $M$ with bounded
curvature and nonnegtive bisectional curvature and $0\leq t < T$.
For any point $x\in M$ and any vector $V$ in the holomorphic
tangent space $T^{1,0}_xM$, let
$$
Q_{\alpha\bar \beta}=\frac{\partial}{\partial t}R_{\alpha\bar
\beta}+ R_{\alpha\bar \gamma}R_{\gamma\bar \beta}
+\nabla_{\gamma}R_{\alpha\bar \beta} V^{\gamma}
+\nabla_{\bar\gamma}R_{\alpha\bar
\beta}V^{\bar\gamma}+R_{\alpha\bar \beta\gamma\bar \delta}V^{
\gamma}V^{\bar\delta}+\frac {1}{t}R_{\alpha\bar \beta} .
$$
Then we have
$$
Q_{\alpha\bar \beta}W^{\alpha}W^{\bar \beta}\geq 0
$$
for all $x\in M$, $V, W\in T^{1,0}_xM$, and $t>0$.
\end{theorem}

\begin{corollary}[{Cao \cite{Cao92}}]
Under the assumptions of Theorem $2.5.9,$ we have
\begin{itemize}
\item[(i)] the scalar curvature $R$ satisfies the estimate
$$
\frac{\partial R}{\partial t}-\frac{|\nabla R|^2}{R}
+\frac{R}{t}\geq 0,$$ and \item[(ii)] assuming $R_{\alpha\bar
\beta}>0$, the determinant
$\phi=\det(R_{\alpha\bar\beta})/\det(g_{\alpha\bar\beta})$ of the
Ricci curvature satisfies the estimate
$$
\frac {\partial \phi}{\partial t}-\frac{|\nabla\phi|^2}{n\phi}
+\frac {n\phi}{t}\geq 0
$$
for all $x\in M$ and $t>0$.
\end{itemize}
\end{corollary}

\section{Perelman's Estimate for Conjugate Heat Equations}

In \cite{P1} Perelman obtained a Li-Yau type estimate for
fundamental solutions of the conjugate heat equation, which is a
backward heat equation, when the metric evolves by the Ricci flow.
In this section we shall describe how to get this estimate along
the same line as in the previous section. More importantly, we
shall show how the Li-Yau path integral, when applied to
Perelman's Li-Yau type estimate, leads to an important space-time
distance function introduced by Perelman \cite{P1}. We learned
from Hamilton \cite{HaL} this idea of looking at Perelman's Li-Yau
estimate.

We saw in the previous section that the Li-Yau quantity and the
Li-Yau-Hamilton quantity vanish on expanding solutions. Note that
when we consider a backward heat equation, shrinking solitons can
be viewed as expanding backward in time. So we start by looking at
shrinking gradient Ricci solitons.

Suppose we have a shrinking gradient Ricci soliton $g_{ij}$ with
potential function $f$ on manifold $M$ and $- \infty < t < 0$ so
that, for $\tau = -t$, \be R_{ij} + \nabla_j \nabla_i f -
\frac{1}{2 \tau} g_{ij} = 0.
\ee Then, by taking the trace, we have \be
R + \Delta f - \frac{n}{2 \tau} = 0. 
\ee Also, by similar calculations as in deriving (1.1.15), we get
\be
R + |\nabla f |^2 - \frac{f}{\tau} = C 
\ee where $C$ is a constant which we can set to be zero.

Moreover, observe \be
\frac{\partial f}{\partial t} = | \nabla f |^2 
\ee because $f$ evolves in time with the rate of change given by
the Lie derivative in the direction of $\nabla f$ generating the
one-parameter family of diffeomorphisms.

Combining (2.6.2) with (2.6.4), we see $f$ satisfies the backward
heat equation \be \frac{\partial f}{\partial t} = - \Delta f + |
\nabla f |^2 - R
+ \frac{n}{2 \tau}, 
\ee or equivalently \be \frac{\partial f}{\partial \tau} = \Delta
f - | \nabla f |^2 +
R - \frac{n}{2 \tau}. 
\ee

Recall the Li-Yau-Hamilton quadratic is a certain combination of
the second order space derivative (or first order time
derivative), first order space derivatives and zero orders.
Multiplying (2.6.2) by a factor of $2$ and subtracting (2.6.3)
yields
$$
2 \Delta f - | \nabla f |^2 + R + \frac{1}{\tau} f -
\frac{n}{\tau}=0
$$
valid for our potential function $f$ of the shrinking gradient
Ricci soliton. The quantity on the LHS of the above identity is
precisely the Li-Yau-Hamilton type quadratic found by Perelman
(cf. section 9 of \cite{P1}).

Note that a function $f$ satisfies the backward heat equation
(2.6.6) if and only if the function
$$
u=(4\pi \tau)^{-\frac{n}{2}}e^{- f}
$$
satisfies the so called \textbf{conjugate heat equation}
\index{conjugate heat equation} \be \square^\ast u \triangleq
\frac{\partial u}{\partial \tau} - \Delta
u + Ru = 0. 
\ee

\begin{lemma}[{Perelman \cite{P1}}]
Let $g_{ij}(x,t)$, $0 \leq t < T$, be a complete solution to the
Ricci flow on an $n$-dimensional manifold $M$ and let $u=(4 \pi
\tau)^{-\frac{n}{2}} e^{-f}$ be a solution to the conjugate
equation $(2.6.7)$ with $\tau=T-t$. Set
$$
H=2 \Delta f -| \nabla f |^2 + R + \frac{f-n}{\tau}
$$
and
$$
v = \tau H u = \big(\tau(R+2 \Delta f - | \nabla f|^2)+f-n \big)
u.
$$
Then we have
$$
\frac{\partial H}{\partial \tau} = \Delta H - 2 \nabla f \cdot
\nabla H - \frac{1}{\tau} H - 2 \left| R_{ij} + \nabla_i\nabla_j f
- \frac{1}{2 \tau} g_{ij} \right|^2,
$$
and
$$
\frac{\partial v}{\partial \tau} = \Delta v - Rv - 2 \tau u
\left|R_{ij}+\nabla_i\nabla_jf - \frac{1}{2 \tau} g_{ij}\right|^2.
$$
\end{lemma}

\begin{pf}
By direct computations, we have
\begin{align*}
\frac{\partial }{\partial \tau}H &  =\frac{\partial}{\partial
\tau}\(2\triangle f -|\nabla f|^2 +R + \frac{f-n}{\tau}\)\\
&  = 2\triangle \(\frac{\partial f}{\partial \tau}\) - 2\langle
2R_{ij}, f_{ij}\rangle - 2\left\langle \nabla
f,\nabla\(\frac{\partial
f}{\partial \tau}\)\right\rangle + 2\Ric(\nabla f,\nabla f)\\
&\quad + \frac{\partial }{\partial
\tau}R+\frac{1}{\tau}\frac{\partial
}{\partial \tau}f - \frac{f-n}{\tau^2}\\
& =2\triangle \(\triangle f -|\nabla f|^2 +R -\frac{n}{2\tau}\) -
4\langle
R_{ij}, f_{ij}\rangle + 2\Ric(\nabla f,\nabla f)\\
&\quad - 2\left\langle \nabla f,\nabla\(\triangle f -|\nabla f|^2
+R
-\frac{n}{2\tau}\)\right\rangle - \triangle R -2|R_{ij}|^2\\
&\quad + \frac{1}{\tau}\(\triangle f -|\nabla f|^2 +R
-\frac{n}{2\tau}\)-\frac{f-n}{\tau^2},\\
\nabla H &  = \nabla \(2\triangle f -|\nabla f|^2 +R
+ \frac{f-n}{\tau}\)\\
&  = 2\nabla(\triangle f) - 2\langle \nabla
\nabla_if,\nabla_if\rangle + \nabla R + \frac{1}{\tau} \nabla f,\\
\triangle H &  = \triangle \(2\triangle f -|\nabla f|^2 +R
+ \frac{f-n}{\tau}\)\\
&  = 2\triangle(\triangle f) - \triangle(|\nabla f|^2)+ \triangle
R + \frac{1}{\tau} \triangle f,
\end{align*}
and
\begin{align*}
2\nabla H \cdot \nabla f &  =2\langle 2\nabla(\triangle f) -
2\langle \nabla \nabla_if,\nabla_if\rangle + \nabla R +
\frac{1}{\tau} \nabla f,\nabla f\rangle \\
&  = 2\left[\langle 2\nabla(\triangle f),\nabla f\rangle -2\langle
f_{ij},f_i f_j\rangle + \langle \nabla R,\nabla f\rangle\right] +
\frac{2}{\tau}|\nabla f|^2.
\end{align*}
Thus we get
\begin{align*}
&  \frac{\partial}{\partial \tau}H - \triangle H + 2 \nabla f
\cdot\nabla H + \frac{1}{\tau} H\\
&  =-4\langle R_{ij},f_{ij}\rangle + 2\Ric(\nabla f,\nabla f) -
2|R_{ij}|^2 -\triangle (|\nabla f|^2) + 2\langle \nabla(\triangle
f),\nabla f\rangle \\
&\quad + \frac{2}{\tau}\triangle f + \frac{2}{\tau}R - \frac{n}{2\tau^2}\\
&  = -2 \left[|R_{ij}|^2 + |f_{ij}|^2 + \frac{n}{4\tau^2} +
2\langle
R_{ij},f_{ij}\rangle -\frac{R}{\tau}-\frac{1}{\tau} \triangle f\right]\\
&  = -2\left|R_{ij} + \nabla_i\nabla_j
f-\frac{1}{2\tau}g_{ij}\right|^2,
\end{align*}
and
\begin{align*}
\(\frac{\partial}{\partial \tau} - \triangle + R\)v
& =\(\frac{\partial}{\partial \tau} - \triangle + R\)(\tau H u)\\
&  = \(\frac{\partial}{\partial \tau} - \triangle\)(\tau
H)\cdot u - 2\langle \nabla(\tau H),\nabla u\rangle\\
&  = \(\(\frac{\partial}{\partial \tau} - \triangle\)(\tau
H) - 2\langle \nabla(\tau H),\nabla f\rangle\)u\\
&  = \tau \(\frac{\partial H}{\partial \tau} - \triangle H + 2
\nabla f
\cdot \nabla H + \frac{1}{\tau} H\)u\\
&  = -2\tau u \left|R_{ij} + \nabla_i\nabla_j f
-\frac{1}{2\tau}g_{ij}\right|^2.
\end{align*}
\end{pf}

Note that, since $f$ satisfies the equation (2.6.6), we can
rewrite $H$ as \be H = 2 \frac{\partial f}{\partial \tau} + |
\nabla f |^2 - R +
\frac{1}{\tau} f. 
\ee Then, by Lemma 2.6.1, we have
$$
\frac{\partial}{\partial \tau}(\tau H)= \Delta (\tau H) -2 \nabla
f \cdot \nabla(\tau H) - 2 \tau \left|\Ric+ \nabla^2f - \frac{1}{2
\tau}g\right|^2.
$$
So by the maximum principle, we find $\operatorname{max}(\tau H)$
is nonincreasing as $\tau$ increasing.  When $u$ is chosen to be a
fundamental solution to (2.6.7), one can show that $\lim_{\tau\to
0}\tau H \leq 0$ and hence $H \leq 0$ on $M$ for all $\tau\in
(0,T]$ (see, for example, \cite{Ni06}). Since this fact is not
used in later chapters and will be only used in the rest of the
section to introduce a space-time distance via Li-Yau path
integral, we omit the details of the proof.

\medskip
Once we have Perelman's Li-Yau type estimate $H\leq 0$, we can
apply the Li-Yau path integral as in \cite{LY} to estimate the
above solution $u$ (i.e., a heat kernel estimate for the conjugate
heat equation, see also the earlier work of Cheeger-Yau
\cite{CY}). Let $p, q \in M$ be two points and $\gamma(\tau), \tau
\in [0, \bar\tau ],$ be a curve joining $p$ and $q$, with
$\gamma(0)=p$ and $\gamma(\bar{\tau})=q$. Then along the
space-time path $(\gamma(\tau), \tau)$, $\tau \in [0, \bar\tau]$,
we have
\begin{align*}
\frac{d}{d \tau} \big(2\sqrt{\tau} f (\gamma(\tau), \tau) \big) &=
2 \sqrt{\tau} \( \frac{\partial f}{\partial \tau} + \nabla
f \cdot \dot{\gamma}(\tau)\) + \frac{1}{\sqrt{\tau}}f \\
&  \leq \sqrt{\tau}\big(-|\nabla f|_{g_{ij}(\tau)}^2 + 2 \nabla
f\cdot\dot{\gamma}(\tau) \big)+\sqrt{\tau}R \\
&  = - \sqrt{\tau} |\nabla f -
\dot{\gamma}(\tau)|_{g_{ij}(\tau)}^2 +
\sqrt{\tau}(R+|\dot{\gamma}(\tau)|_{g_{ij}(\tau)}^2) \\
&  \leq \sqrt{\tau}(R + |\dot{\gamma}(\tau)|_{g_{ij}(\tau)}^2)
\end{align*}
where we have used the fact that $H\leq 0$ and the expression for
$H$ in (2.6.8).

Integrating the above inequality from $\tau =0$ to $\tau =
\bar{\tau}$, we obtain
$$
2 \sqrt{\bar{\tau}} f(q, \bar{\tau}) \leq \int^{\bar{\tau}}_0
\sqrt{\tau} (R+| \dot{\gamma}(\tau)|_{g_{ij}(\tau)}^2)d \tau,
$$

\noindent or
$$f(q, \bar{\tau}) \leq \frac{1}{2
\sqrt{\bar{\tau}}}\; \mathcal{L}(\gamma),
$$
where \be \mathcal{L}(\gamma) \triangleq \int^{\bar{\tau}}_0
\sqrt{\tau}(R+|
\dot{\gamma}(\tau)|_{g_{ij}(\tau)}^2)d\tau.  
\ee

Denote by \be l(q, \bar{\tau}) \triangleq \inf_{\gamma} \frac{1}{2
\sqrt{\bar{\tau}}}\; \mathcal{L}(\gamma),
\ee where the $inf$ is taken over all space curves $\gamma (\tau),
0 \leq \tau \leq \bar{\tau}$, joining $p$ and $q$. The space-time
distance function $l(q, \bar{\tau})$ obtained by the above Li-Yau
path integral argument is first introduced by Perelman in
\cite{P1} and is what Perelman calls reduced distance. Since
Perelman pointed out in page 19 of \cite{P1} that ``an even closer
reference in \cite{LY}, where they use `length', associated to a
linear parabolic equation, which is pretty much the same as in our
case", it is natural to call $l(q, \bar{\tau})$ the
\textbf{Li-Yau-Perelman distance}. \index{Li-Yau-Perelman
distance} See Chapter 3 for much more detailed discussions.

Finally, we conclude this section by relating the quantity $H$ (or
$v$) and the $\mathcal W$-functional of Perelman defined in
(1.5.9). Observe that $v$ happens to be the integrand of the
$\mathcal W$-functional,
$$
\mathcal{W}(g_{ij}(t), f, \tau) = \int_M v dV.
$$
Hence, when $M$ is compact,
\begin{align*}
\frac{d}{d \tau}\mathcal{W}
&  =\int_M \(\frac{\partial}{\partial \tau}v+Rv\)dV \\
& =-2\tau\int_M\left|\Ric+\nabla^2f-\frac{1}{2\tau}g\right|^2u d V\\
&  \leq 0,
\end{align*}
or equivalently,
$$
\frac{d} {d t} \mathcal{W}(g_{ij}(t),f(t),\tau(t))=\int_M
2\tau\left|R_{ij}+\nabla_i\nabla_jf-\frac
{1}{2\tau}g_{ij}\right|^2 (4\pi\tau)^{-\frac{n}{2}} e^{-f}dV,
$$
which is the same as stated in Proposition 1.5.8.

\newpage
\part{{\Large Perelman's Reduced Volume}}

\bigskip
In Section 1.5 we introduced the $\mathcal{F}$-functional and the
$\mathcal{W}$-functional of Perelman \cite{P1} and proved their
monotonicity properties under the Ricci flow. In the last section
of the previous chapter we have defined the Li-Yau-Perelman
distance. The main purpose of this chapter is to use the
Li-Yau-Perelman distance to define the Perelman's reduced volume,
which was introduced by Perelman in \cite{P1}, and prove the
monotonicity property of the reduced volume under the Ricci flow.
This new monotonicity formula of Perelman \cite{P1} is more useful
for local considerations, especially when we consider the
formation of singularities in Chapter 6 and work on the Ricci flow
with surgery in Chapter 7. As first applications we will present
two no local collapsing theorems of Perelman \cite{P1} in this
chapter. More applications can be found in Chapter 6 and 7. This
chapter is a detailed exposition of sections 6-8 of Perelman
\cite{P1}.

\section{Riemannian Formalism in Potentially Infinite Dimensions}

In Section 2.6, from an analytic view point, we saw how the Li-Yau
path integral of Perelman's estimate for fundamental solutions to
the conjugate heat equation leads to the Li-Yau-Perelman distance.
In this section we present, from a geometric view point, another
motivation (cf. section 6 of \cite{P1}) how one is lead to the
consideration of the Li-Yau-Perelman distance function, as well as
a reduced volume concept. Interestingly enough, the
Li-Yau-Hamilton quadratic introduced in Section 2.5 appears again
in this geometric consideration.

We consider the Ricci flow
$$
\frac{\partial}{\partial t}g_{ij}=-2R_{ij}
$$
on a manifold $M$ where we assume that $g_{ij}(\cdot,t)$ are
complete and have uniformly bounded curvatures.

Recall from Section 2.5 that the Li-Yau-Hamilton quadratic
introduced in \cite{Ha93} is
$$
Q=M_{ij}W_iW_j+2P_{ijk}U_{ij}W_k+R_{ijkl}U_{ij}U_{kl}
$$
where
\begin{align*}
M_{ij}& = \Delta R_{ij}-\frac{1}{2}\nabla_i\nabla_jR
+2R_{ikjl}R_{kl}-R_{ik}R_{jk}+\frac{1}{2t}R_{ij},\\
P_{ijk}&  = \nabla_iR_{jk}-\nabla_jR_{ik}
\end{align*}
and $U_{ij}$ is any two-form and $W_i$ is any 1-form. Here and
throughout this chapter we do not always bother to raise indices;
repeated indices is short hand for contraction with respect to the
metric.

In \cite{Ha95F}, Hamilton predicted that the Li-Yau-Hamilton
quadratic is some sort of jet extension of positive curvature
operator on some larger space. Such an interpretation of the
Li-Yau-Hamilton quadratic as a curvature operator on the space
$M\times \mathbb{R}^{+}$ was found by Chow and Chu \cite{ChCh95}
where a potentially degenerate Riemannian metric on $M\times
\mathbb{R}^+$ was constructed. The degenerate Riemannian metric on
$M\times \mathbb{R}^+$ is the limit of the following two-parameter
family of Riemannian metrics
$$
g_{N,\delta}(x,t)=g(x,t)+(R(x,t)+\frac{N}{2(t+\delta)})dt^{2}
$$
as $N$ tends to infinity and $\delta$ tends to zero, where
$g(x,t)$ is the solution of the Ricci flow on $M$ and $t\in
\mathbb{R}^{+}$.

To avoid the degeneracy, Perelman \cite{P1} considers the manifold
$\tilde{M}=M\times \mathbb{S}^N\times \mathbb{R}^+$ with the
following metric:
\begin{align*}
\tilde{g}_{ij}&  = g_{ij},\\
\tilde{g}_{\alpha\beta}&  = \tau g_{\alpha\beta},\\
\tilde{g}_{oo}&  = \frac{N}{2\tau}+R,\\
\tilde{g}_{i\alpha} &  = \tilde{g}_{io} =\tilde{g}_{\alpha o}=0,
\end{align*}
where $i,j$ are coordinate indices on $M$; $\alpha,\;\beta$ are
coordinate indices on $\mathbb{S}^N$; and the coordinate $\tau$ on
$R^+$ has index $o$. Let $\tau = T-t$ for some fixed constant $T$.
Then $g_{ij}$ will evolve with $\tau$ by the backward Ricci flow
$\frac{\partial}{\partial \tau}g_{ij}=2R_{ij}$. The metric
$g_{\alpha\beta}$ on $\mathbb{S}^N$ is a metric with constant
sectional curvature $\frac{1}{2N}$.

We remark that the metric $\tilde{g}_{\alpha\beta}$ on
$\mathbb{S}^N$ is chosen so that the product metric
$(\tilde{g}_{ij},\ \tilde{g}_{\alpha\beta})$ on $M\times
\mathbb{S}^N$ evolves by the Ricci flow, while the component
$\tilde{g}_{oo}$ is just the scalar curvature of
$(\tilde{g}_{ij},\ \tilde{g}_{\alpha\beta})$. Thus the metric
$\tilde{g}$ defined on $\tilde{M}=M\times \mathbb{S}^N \times
\mathbb{R}^+$ is exactly a ``regularization" of Chow-Chu's
degenerate metric on $M\times\mathbb{R}^+$.  \vskip 0.2cm

\begin{proposition} [cf. \cite{ChCh95}]
The components of the curvature tensor of the metric $\tilde{g}$
coincide $($modulo $N^{-1})$ with the components of the
Li-Yau-Hamilton quadratic.
\end{proposition}

\begin{pf}
By definition, the Christoffel symbols of the metric $\tilde{g}$
are given by the following list: {\allowdisplaybreaks
\begin{align*}
\tilde{\Gamma}_{ij}^k&=\Gamma_{ij}^k,\\
\tilde{\Gamma}_{i\beta}^k&=0\quad \text{and}\quad
\tilde{\Gamma}_{ij}^\gamma=0,\\
\tilde{\Gamma}_{\alpha\beta}^k&=0\quad \text{and}\quad
\tilde{\Gamma}_{i\beta}^\gamma=0,\\
\tilde{\Gamma}_{io}^k&=g^{kl}R_{li}\quad
\text{and}\quad\tilde{\Gamma}_{ij}^o=-\tilde{g}^{oo}R_{ij},\\
\tilde{\Gamma}_{oo}^k&=-\frac{1}{2}g^{kl} \frac{\partial}{\partial
x^l}R\quad
\text{and}\quad\tilde{\Gamma}_{io}^o=\frac{1}{2}\tilde{g}^{oo}
\frac{\partial}{\partial x^i}R,\\
\tilde{\Gamma}_{i\beta}^o&=0,\;\tilde{\Gamma}_{o\beta}^k=0\quad
\text{and}\quad\tilde{\Gamma}_{oj}^\gamma=0,\\
\tilde{\Gamma}_{\alpha\beta}^\gamma&=\Gamma_{\alpha\beta}^\gamma,\\
\tilde{\Gamma}_{\alpha o}^\gamma
&=\frac{1}{2\tau}\delta_\alpha^\gamma\quad
\text{and}\quad\tilde{\Gamma}_{\alpha\beta}^o
=-\frac{1}{2}\tilde{g}^{oo}g_{\alpha\beta},\\
\tilde{\Gamma}_{oo}^\gamma&=0\quad \text{and}\quad
\tilde{\Gamma}_{o\beta}^o=0,\\
\tilde{\Gamma}_{oo}^o&=\frac{1}{2}\tilde{g}^{oo}
\(-\frac{N}{2\tau^2}+\frac{\partial}{\partial\tau}R\).
\end{align*}
} Fix a point $(p,s,\tau)\in M\times \mathbb{S}^N\times
\mathbb{R}^+$ and choose normal coordinates around $p\in M$ and
normal coordinates around $s\in \mathbb{S}^N$ such that
$\Gamma_{ij}^k(p)=0$ and $\Gamma_{\alpha\beta}^\gamma(s)=0$ for
all $i,j,k$ and $\alpha,\beta,\gamma$. We compute the curvature
tensor $\tilde{R}m$ of the metric $\tilde{g}$ at the point as
follows:
\begin{align*}
\tilde{R}_{ijkl}& =R_{ijkl}+\tilde{\Gamma}_{io}^k
\tilde{\Gamma}_{jl}^o-\tilde{\Gamma}_{jo}^k\tilde{\Gamma}_{il}^o=
R_{ijkl}+O\(\frac{1}{N}\),\\
\tilde{R}_{ijk\delta}&  = 0,\\
\tilde{R}_{ij\gamma\delta}&  = 0\quad \text{and} \quad
\tilde{R}_{i\beta
k\delta}=\tilde{\Gamma}_{io}^k\tilde{\Gamma}_{\beta\delta}^o
-\tilde{\Gamma}_{\beta o}^k\tilde{\Gamma}_{i\delta}^o
=-\frac{1}{2}\tilde{g}^{oo}g_{\beta\delta}g^{kl}R_{li}
=O\(\frac{1}{N}\),\\
\tilde{R}_{i\beta\gamma\delta}&  = 0,\\
\tilde{R}_{ijko}&  = \frac{\partial}{\partial
x^i}R_{jk}-\frac{\partial}{\partial
x^j}R_{ik}+\tilde{\Gamma}_{io}^k\tilde{\Gamma}_{jo}^o-
\tilde{\Gamma}_{jo}^k\tilde{\Gamma}_{io}^o=P_{ijk}+O\(\frac{1}{N}\),\\
\tilde{R}_{ioko}&  = -\frac{1}{2}\frac{\partial^2}{\partial
x^i\partial x^k}R-\frac{\partial}{\partial
\tau}(R_{il}g^{lk})+\tilde{\Gamma}_{io}^k\tilde{\Gamma}_{oo}^o
-\tilde{\Gamma}_{oj}^k\tilde{\Gamma}_{io}^j
-\tilde{\Gamma}_{oo}^k\tilde{\Gamma}_{io}^o\\
&  = -\frac{1}{2}\nabla_i\nabla_kR-\frac{\partial}{\partial
\tau}R_{ik}+2R_{ik}R_{lk}-\frac{1}{2\tau}R_{ik}-R_{ij}R_{jk}
+O\(\frac{1}{N}\)\\
&  = M_{ik}+O(\frac{1}{N}),\\
\tilde{R}_{ij\gamma o}&  = 0\quad
\text{and}\quad \tilde{R}_{i\gamma jo}=0,\\
\tilde{R}_{i\beta\gamma o}&  = -\tau\tilde{\Gamma}_{\beta
o}^\gamma\tilde{\Gamma}_{io}^o=O\(\frac{1}{N}\)\quad
\mbox{and}\quad
\tilde{R}_{io\gamma\delta}=0,\\
\tilde{R}_{io\gamma o}&  = 0,\\
\tilde{R}_{\alpha\beta\gamma o}&  = 0, \\
\tilde{R}_{\alpha o\gamma o}& =
\(\frac{1}{2\tau^2}\delta_\alpha^\gamma +\tilde{\Gamma}_{\alpha
o}^\gamma\tilde{\Gamma}_{oo}^o
-\tilde{\Gamma}_{o\beta}^\gamma\tilde{\Gamma}_{\alpha
o}^\beta\)\tau=O\(\frac{1}{N}\),\\
\tilde{R}_{\alpha\beta\gamma\delta}&  = O\(\frac{1}{N}\).
\end{align*}

Thus the components of the curvature tensor of the metric
$\tilde{g}$ coincide (modulo $N^{-1}$) with the components of the
Li-Yau-Hamilton quadratic.
\end{pf}

The following observation due to Perelman \cite{P1} gives an
important motivation to define Perelman's reduced volume.

\begin{corollary}
All components of the Ricci tensor of $\tilde{g}$ are zero
$($modulo $N^{-1})$.
\end{corollary}

\begin{pf}
{}From the list of the components of the curvature tensor of
$\tilde{g}$ given above, we have
\begin{align*}
\tilde{R}_{ij}& =
\tilde{g}^{kl}\tilde{R}_{ijkl}+\tilde{g}^{\alpha\beta}\tilde{R}_{i\alpha
j\beta}+\tilde{g}^{oo}\tilde{R}_{iojo}\\
& =R_{ij}-\frac{1}{2\tau}g^{\alpha\beta}\tilde{g}^{oo}
g^{\alpha\beta}R_{ij}+\tilde{g}^{oo}
\(M_{ij}-\frac{1}{2\tau}R_{ij}+O\(\frac{1}{N}\)\)\\
&  =R_{ij}-\frac{N}{2\tau}\tilde{g}^{oo}R_{ij}+O\(\frac{1}{N}\)\\
&  =O\(\frac{1}{N}\),\\
\tilde{R}_{i\gamma}&  = \tilde{g}^{kl}\tilde{R}_{ik\gamma
l}+\tilde{g}^{\alpha\beta}\tilde{R}_{i\alpha\gamma\beta}
+\tilde{g}^{oo}\tilde{R}_{io\gamma o}=0,\\
\tilde{R}_{io}& =\tilde{g}^{kl}\tilde{R}_{ikol}
+\tilde{g}^{\alpha\beta}\tilde{R}_{i\alpha
o\beta}+\tilde{g}^{oo}\tilde{R}_{iooo}\\
&  = -g^{kl}P_{ikl}+O\(\frac{1}{N}\),\\
\tilde{R}_{\alpha\beta}&  = \tilde{g}^{kl}\tilde{R}_{\alpha k\beta
l}+\tilde{g}^{\gamma\delta}\tilde{R}_{\alpha\gamma\beta\delta}
+\tilde{g}^{oo}\tilde{R}_{\alpha o \beta o}\\
&  = O\(\frac{1}{N}\),\\
\tilde{R}_{\alpha o}&  = \tilde{g}^{kl}\tilde{R}_{\alpha
kol}+\tilde{g}^{\beta\gamma}\tilde{R}_{\alpha\beta
o\gamma}+\tilde{g}^{oo}\tilde{R}_{\alpha ooo}=0,\\
\tilde{R}_{oo}& =
\tilde{g}^{kl}\tilde{R}_{okol}+\tilde{g}^{\alpha\beta}\tilde{R}_{o\alpha
o\beta}+\tilde{g}^{oo}\tilde{R}_{oooo}\\
&  = g^{kl}\(M_{kl}+O\(\frac{1}{N}\)\)+O\(\frac{1}{N}\).
\end{align*}
Since $\tilde{g}^{oo}$ is of order $N^{-1}$, we see that the norm
of the Ricci tensor is given by
$$
|\widetilde{\Ric}|_{\tilde{g}}=O\(\frac{1}{N}\).
$$
This proves the result.
\end{pf}

We now use the Ricci-flatness of the metric $\tilde{g}$ to
interpret the Bishop-Gromov relative volume comparison theorem
which will motivate another monotonicity formula for the Ricci
flow. The argument in the following will not be rigorous. However
it gives an intuitive picture of what one may expect. Consider a
metric ball in $(\tilde{M},\tilde{g})$ centered at some point
$(p,s,0)\in\tilde{M}$. Note that the metric of the sphere
$\mathbb{S}^N$ at $\tau=0$ degenerates and it shrinks to a point.
Then the shortest geodesic $\gamma(\tau)$ between $(p,s,0)$ and an
arbitrary point $(q,\bar{s},\bar\tau)\in\tilde{M}$ is always
orthogonal to the $\mathbb{S}^N$ fibre. The length of
$\gamma(\tau)$ can be computed as
\begin{align*}
& \int_0^{\bar\tau}\sqrt{\(\frac{N}{2\tau}+R\)
+|\dot{\gamma}(\tau)|^2_{g_{ij}(\tau)}}d\tau\\
&
=\sqrt{2N\bar\tau}+\frac{1}{\sqrt{2N}}\int_0^{\bar\tau}\sqrt{\tau}
(R+|\dot{\gamma}(\tau)|^2_{g_{ij}})d\tau+O(N^{-\frac{3}{2}}).
\end{align*}
Thus a shortest geodesic should minimize
$$
\mathcal{L}(\gamma)=\int_0^{\bar\tau}\sqrt{\tau}
(R+|\dot{\gamma}(\tau)|^2_{g_{ij}})d\tau.
$$

Let $L(q,\bar\tau)$ denote the corresponding minimum. We claim
that a metric sphere $S_{\tilde{M}}(\sqrt{2N\bar\tau})$ in
$\tilde{M}$ of radius $\sqrt{2N\bar\tau}$ centered at $(p,s,0)$ is
$O(N^{-1})$-close to the hypersurface $\{\tau=\bar\tau \}$.
Indeed, if $(x,s',\tau(x))$ lies on the metric sphere
$S_{\tilde{M}}(\sqrt{2N\bar\tau})$, then the distance between
$(x,s',\tau(x))$ and $(p,s,0)$ is
$$
\sqrt{2N\bar\tau}
=\sqrt{2N\tau(x)}+\frac{1}{\sqrt{2N}}L(x,\tau(x))+O\(N^{-\frac{3}{2}}\)
$$
which can be written as
$$
\sqrt{\tau(x)}-\sqrt{\bar\tau}=-\frac{1}{2N}L(x,\tau(x))+O(N^{-2})
=O(N^{-1}).$$ This shows that the metric sphere
$S_{\tilde{M}}(\sqrt{2N\bar\tau})$ is $O(N^{-1})$-close to the
hypersurface $\{\tau=\bar\tau \}$. Note that the metric
$g_{\alpha\beta}$ on $\mathbb{S}^N$ has constant sectional
curvature $\frac{1}{2N}$. Thus
\begin{align*}
& {\rm Vol}\,\(S_{\tilde{M}}\(\sqrt{2N\bar\tau}\)\) \\
& \approx \int_{M}\(\int_{S^N}dV_{\tau(x)g_{\alpha\beta}}\)dV_{g_{ij}}(x)\\
& = \int_M(\tau(x))^{\frac{N}{2}}{\rm Vol}\,(S^N)dV_M\\
& \approx (2N)^{\frac{N}{2}}\omega_N\int_M\(\sqrt{\bar\tau}
-\frac{1}{2N}L(x,\tau(x))+O(N^{-2})\)^NdV_M\\
& \approx (2N)^{\frac{N}{2}}\omega_N\int_M
\(\sqrt{\bar\tau}-\frac{1}{2N}L(x,\bar\tau)+o(N^{-1})\)^NdV_M,
\end{align*}
where $\omega_N$ is the volume of the standard $N$-dimensional
sphere. Now the volume of Euclidean sphere of radius
$\sqrt{2N\bar\tau}$ in $\mathbb{R}^{n+N+1}$ is
$$
{\rm Vol}\,(S_{\mathbb{R}^{n+N+1}}(\sqrt{2N\bar\tau}))
=(2N\bar\tau)^\frac{N+n}{2}\omega_{n+N}.
$$
Thus we have
$$
\frac{{\rm Vol}\,(S_{\tilde{M}}(\sqrt{2N\bar\tau}))}{{\rm
Vol}\,(S_{
\mathbb{R}^{n+N+1}}(\sqrt{2N\bar\tau}))}\approx\text{const}\cdot
N^{-\frac{n}{2}}\cdot \int_M(\bar\tau)^{-\frac{n}{2}}\exp
\left\{{-\frac{1}{2\sqrt{\bar\tau}}L(x,\bar\tau)}\right\}dV_M.
$$
Since the Ricci curvature of $\tilde{M}$ is zero (modulo
$N^{-1}$), the Bishop-Gromov volume comparison theorem then
suggests that the integral
$$
\tilde{V}(\bar\tau)\stackrel{\Delta}{=}\int_M
(4\pi\bar\tau)^{-\frac{n}{2}}\exp
\left\{-\frac{1}{2\sqrt{\bar\tau}}L(x,\bar\tau)\right\}dV_M,
$$
which we will call {\bf Perelman's reduced
volume}\index{Perelman's reduced volume}, should be nonincreasing
in $\bar\tau$. A rigorous proof of this monotonicity property will
be given in the next section. One should note the analog of
reduced volume with the heat kernel and there is a parallel
calculation for the heat kernel of the Shr\"odinger equation in
the paper of Li-Yau \cite{LY}.

\section{Comparison Theorems for Perelman's Reduced Volume}

In this section we will write the Ricci flow in the backward
version
$$
\frac{\partial}{\partial \tau}g_{ij}=2R_{ij}
$$
on a manifold $M$ with $\tau=\tau(t)$ satisfying $d\tau/dt=-1$ (in
practice we often take $\tau=t_0-t$ for some fixed time $t_0$). We
always assume that either $M$ is compact or $g_{ij}(\tau)$ are
complete and have uniformly bounded curvature. To each (smooth)
space curve $\gamma(\tau)$, $0<\tau_1\le\tau\le\tau_2$, in $M$, we
define its {\bf $\mathcal{L}$-length}\index{$\mathcal{L}$-length}
as
$$
\mathcal{L}(\gamma)
=\int_{\tau_1}^{\tau_2}\sqrt{\tau}(R(\gamma(\tau),\tau)+
|\dot{\gamma}(\tau)|^2_{g_{ij}(\tau)})d\tau.
$$

Let $X(\tau)=\dot{\gamma}(\tau)$, and let $Y(\tau)$ be any
(smooth) vector field along $\gamma(\tau)$. First of all, we
compute the first variation formula for $\mathcal{L}$-length (cf.
section 7 of \cite{P1}).

\begin{lemma}[First variation formula]
$$
\delta_Y(\mathcal{L})=2\sqrt{\tau}\<X,Y\>|_{\tau_1}^{\tau_2}
+\int_{\tau_1}^{\tau_2}\sqrt{\tau}\left\<Y,\nabla
R-2\nabla_XX-4\Ric(\cdot,X)-\frac{1}{\tau}X\right\>d\tau
$$
where $\<\cdot,\cdot\>$ denotes the inner product with respect to
the metric $g_{ij}(\tau)$.
\end{lemma}

\begin{pf}
By direct computations,
\begin{align*}
\delta_Y(\mathcal{L}) &
=\int_{\tau_1}^{\tau_2}\sqrt{\tau}(\<\nabla R,Y\>
+2\<X,\nabla_YX\>)d\tau\\
&  = \int_{\tau_1}^{\tau_2}\sqrt{\tau}(\<\nabla R,Y\>
+2\<X,\nabla_XY\>)d\tau\\
&  = \int_{\tau_1}^{\tau_2}\sqrt{\tau}\(\<\nabla R,Y\>
+2\frac{d}{d\tau}\<X,Y\>-2\<\nabla_XX,Y\>-4\Ric(X,Y)\)d\tau\\
& =2\sqrt{\tau}\<X,Y\>|_{\tau_1}^{\tau_2}
+\int_{\tau_1}^{\tau_2}\sqrt{\tau}\left\<Y,\nabla R
-2\nabla_XX-4\Ric(\cdot,X)-\frac{1}{\tau}X\right\>d\tau.
\end{align*}
\end{pf}

A smooth curve $\gamma(\tau)$ in $M$ is called an {\bf
$\mathcal{L}$-geodesic}\index{$\mathcal{L}$-geodesic}\index{$\mathcal{L}$-geodesic!curve}
if it satisfies the following {\bf $\mathcal{L}$-geodesic
equation}\index{$\mathcal{L}$-geodesic!equation} \be
\nabla_XX-\frac{1}{2}\nabla R+\frac{1}{2\tau}X+2\Ric(X,\cdot)=0.
\ee Given any two points $p,q\in M$ and $\tau_2>\tau_1>0$, there
always exists an $\mathcal{L}$-shortest curve (or shortest
$\mathcal{L}$-geodesic) $\gamma(\tau)$:
$[\tau_1,\tau_2]\rightarrow M$ connecting $p$ to $q$ which
satisfies the above $\mathcal{L}$-geodesic equation. Multiplying
the $\mathcal{L}$-geodesic equation (3.2.1) by $\sqrt{\tau}$, we
get
$$
\nabla_X(\sqrt{\tau}X)=\frac{\sqrt{\tau}}{2}\nabla R
-2\sqrt{\tau}\Ric(X,\cdot)\quad \text{on}\quad [\tau_1,\tau_2],
$$
or equivalently
$$
\frac{d}{d\tau}(\sqrt{\tau}X)=\frac{\sqrt{\tau}}{2}\nabla R
-2\Ric(\sqrt{\tau}X,\cdot) \quad \text{on}\quad [\tau_1,\tau_2].
$$
Thus if a continuous curve, defined on $[0,\tau_2]$, satisfies the
$\mathcal{L}$-geodesic equation on every subinterval
$0<\tau_1\le\tau\le\tau_2$, then $\sqrt{\tau_1}X(\tau_1)$ has a
limit as $\tau_1\rightarrow 0^+$. This allows us to extend the
definition of the $\mathcal{L}$-length to include the case
$\tau_1=0$ for all those (continuous) curves $\gamma:
[0,\tau_2]\rightarrow M$ which are smooth on $(0,\tau_2]$ and have
limits
$\lim\limits_{\tau\rightarrow0^+}\sqrt{\tau}\dot{\gamma}(\tau)$.
Clearly, there still exists an $\mathcal{L}$-shortest curve
$\gamma(\tau):[0,\tau_2]\rightarrow M$ connecting arbitrary two
points $p,q\in M$ and satisfying the $\mathcal{L}$-geodesic
equation (3.2.1) on $(0,\tau_2]$. Moreover, for any vector $v \in
T_pM$, we can find an $\mathcal{L}$-geodesic $\gamma(\tau)$
starting at $p$ with
$\lim\limits_{\tau\rightarrow0^+}\sqrt{\tau}\dot{\gamma}(\tau)=v$.

{}From now on, we fix a point $p\in M$ and set $\tau_1=0$. The
{\bf $\mathcal{L}$-distance function}\index{$\mathcal{L}$-distance
function} on the space-time $M\times \mathbb{R}^+$ is denoted by
$L(q,\bar{\tau})$ and defined to be the $\mathcal{L}$-length of
the $\mathcal{L}$-shortest curve $\gamma(\tau)$ connecting $p$ and
$q$ with $0\le\tau\le\bar{\tau}$.

Consider a shortest $\mathcal{L}$-geodesic $\gamma: [0,\bar \tau]
\rightarrow M$ connecting $p$ to $q$. In the computations below we
pretend that $\mathcal{L}$-shortest geodesics between $p$ and $q$
are unique for all pairs $(q,\bar \tau)$; if this is not the case,
the inequalities that we obtain are still valid, by a standard
barrier argument, when understood in the sense of distributions
(see, for example, \cite{ScY}).

The first variation formula in Lemma 3.2.1 implies that
$$
\nabla_YL(q,\bar{\tau}) =\left\<2\sqrt{\bar{\tau}}X(\bar{\tau}),
Y(\bar{\tau})\right\>.
$$
Thus
$$
\nabla L(q,\bar{\tau}) =2\sqrt{\bar{\tau}}X(\bar{\tau}),
$$
and \be |\nabla L|^2 =4\bar{\tau}|X|^2
=-4\bar{\tau}R+4\bar{\tau}(R+|X|^2).
\ee We also compute
\begin{align}
L_{\bar{\tau}}(\gamma(\bar{\tau}),\bar{\tau})& =
\frac{d}{d\tau}L(\gamma(\tau),\tau)|_{\tau=\bar{\tau}}
-\<\nabla L,X\>\\
&  = \sqrt{\bar{\tau}}(R+|X|^2)-2\sqrt{\bar{\tau}}|X|^2 \nn\\
&  = 2\sqrt{\bar{\tau}}R-\sqrt{\bar{\tau}}(R+|X|^2).\nn
\end{align}
To evaluate $R+|X|^2$, we compute by using (3.2.1),
\begin{align*}
& \frac{d}{d\tau}(R(\gamma(\tau),\tau)+|X(\tau)|^2_{g_{ij}(\tau)})\\
& = R_\tau+\<\nabla R,X\>+2\<\nabla_XX,X\>+2{ \Ric}(X,X)\\
& = R_\tau+\frac{1}{\tau}R+2\<\nabla R,X\>-2 {\Ric}(X,X)
-\frac{1}{\tau}(R+|X|^2)\\
&  = -Q(X)-\frac{1}{\tau}(R+|X|^2),
\end{align*}
where
$$
Q(X)=-R_\tau-\frac{R}{\tau}-2\<\nabla R,X\>+2{ \Ric}(X,X)
$$
is the trace Li-Yau-Hamilton quadratic in Corollary 2.5.5. Hence
\begin{align*}
\frac{d}{d\tau}(\tau^\frac{3}{2}(R+|X|^2))|_{\tau=\bar{\tau}} &
=\frac{1}{2}\sqrt{\bar{\tau}}(R+|X|^2)
-\bar{\tau}^\frac{3}{2}Q(X)\\
&
=\frac{1}{2}\frac{d}{d\tau}L(\gamma(\tau),\tau)|_{\tau=\bar{\tau}}
-\bar{\tau}^\frac{3}{2}Q(X).
\end{align*}
Therefore, \be
\bar{\tau}^\frac{3}{2}(R+|X|^2)=\frac{1}{2}L(q,\bar{\tau})-K,
\ee where \be K=\int_0^{\bar{\tau}}\tau^\frac{3}{2}Q(X)d\tau.
\ee Combining (3.2.2) with (3.2.3), we obtain \be |\nabla L|^2
=-4\bar{\tau}R+\frac{2}{\sqrt{\bar{\tau}}}L
-\frac{4}{\sqrt{\bar{\tau}}}K  
\ee and \be L_{\bar{\tau}}
=2\sqrt{\bar{\tau}}R-\frac{1}{2\bar{\tau}}L
+\frac{1}{\bar{\tau}}K.  
\ee

Next we compute the second variation of an $\mathcal{L}$-geodesic
(cf. section 7 of \cite{P1}).

\begin{lemma}[Second variation formula]
For any $\mathcal{L}$-geodesic $\gamma$, we have
\begin{align*}
\delta_Y^2(\mathcal{L}) &
=2\sqrt{\tau}\<\nabla_YY,X\>|_0^{\bar{\tau}}
+\int_0^{\bar{\tau}}\sqrt{\tau}[2|\nabla_XY|^2
+2\<R(Y,X)Y,X\>\\
&\quad +\nabla_Y\nabla_YR+2\nabla_X{\Ric}(Y,Y)-4\nabla_Y
\Ric(Y,X)]d\tau.
\end{align*}
\end{lemma}

\begin{pf}
We compute
\begin{align*}
\delta_Y^2(\mathcal{L})
& = Y\(\int_0^{\bar{\tau}}\sqrt{{\tau}}(Y(R)+2\<\nabla_YX,X\>)d\tau\)\\
& = \int_0^{\bar{\tau}}\sqrt{{\tau}}(Y(Y(R))
+2\<\nabla_Y\nabla_YX,X\>+2|\nabla_YX|^2)d\tau\\
& =\int_0^{\bar{\tau}}\sqrt{{\tau}}(Y(Y(R))
+2\<\nabla_Y\nabla_XY,X\>+2|\nabla_XY|^2)d\tau
\end{align*}
and
\begin{align*}
& 2\<\nabla_Y\nabla_XY,X\>\\
&  = 2\<\nabla_X\nabla_YY,X\>+2\<R(Y,X)Y,X\>\\
&  = 2\frac{d}{d\tau}\<\nabla_YY,X\>-4{ \Ric}(\nabla_YY,X)
-2\<\nabla_YY,\nabla_XX\> \\
&\quad-\(2\left\<\frac{d}{d\tau}\nabla_YY,X\right\>
-2\<\nabla_X\nabla_YY,X\>\)+2\<R(Y,X)Y,X\>\\
&  = 2\frac{d}{d\tau}\<\nabla_YY,X\>
-4{\Ric}(\nabla_YY,X)-2\<\nabla_YY,\nabla_XX\>\\
&\quad -2\left\<Y^iY^j(g^{kl}(\nabla_iR_{lj}+\nabla_jR_{li}
-\nabla_lR_{ij}))\frac{\partial}{\partial x^k},X\right\>+2\<R(Y,X)Y,X\>\\
& = 2\frac{d}{d\tau}\<\nabla_YY,X\>-4{\Ric}(\nabla_YY,X)
-2\<\nabla_YY,\nabla_XX\>-4\nabla_Y{\Ric}(X,Y)\\
&\quad +2\nabla_X{\Ric}(Y,Y)+2\<R(Y,X)Y,X\>,
\end{align*}
where we have used the computation $$\frac{\partial}{\partial
\tau}\Gamma_{ij}^k=g^{kl}(\nabla_iR_{lj}+\nabla_jR_{li}-\nabla_lR_{ij}).$$
Thus by using the $\mathcal{L}$-geodesic equation (3.2.1), we get
\begin{align*}
\delta_Y^2(\mathcal{L}) &
=\int_0^{\bar{\tau}}\sqrt{\tau}\bigg[Y(Y(R))
+2\frac{d}{d\tau}\<\nabla_YY,X\>-4\Ric(\nabla_YY,X)\\
&\quad-2\<\nabla_YY,\nabla_XX\>
-4\nabla_Y\Ric(X,Y)+2\nabla_X \Ric(Y,Y) \\
&\quad+2\<R(Y,X)Y,X\>+2|\nabla_XY|^2\bigg]d\tau\\
&  =
\int_0^{\bar{\tau}}\bigg[2\sqrt{\tau}\frac{d}{d\tau}\<\nabla_YY,X\>
+\frac{1}{\sqrt{\tau}}\<\nabla_YY,X\>\bigg]d\tau\\
&\quad +\int_0^{\bar{\tau}}\sqrt{\tau}[Y(Y(R))
-\<\nabla_YY,\nabla R\>-4\nabla_Y \Ric(X,Y)\\
&\quad +2\nabla_X \Ric(Y,Y)+2\<R(Y,X)Y,X\>+2|\nabla_XY|^2]d\tau\\
&  = 2\sqrt{\tau}\<\nabla_YY,X\>|_0^{\bar{\tau}}
+\int_0^{\bar{\tau}}\sqrt{\tau}[2|\nabla_XY|^2
+2\<R(Y,X)Y,X\> \\
&\quad+\nabla_Y\nabla_YR -4\nabla_Y \Ric(X,Y) +2\nabla_X
\Ric(Y,Y)]d\tau.
\end{align*}

\vspace*{-5mm}
\end{pf}

We now use the above second variation formula to estimate the
Hessian of the $\mathcal{L}$-distance function.

Let $\gamma(\tau):[0,\bar{\tau}]\rightarrow M$ be an
$\mathcal{L}$-shortest curve connecting $p$ and $q$ so that the
$\mathcal{L}$-distance function $L=L(q,\bar\tau)$ is given by the
$\mathcal{L}$-length of $\gamma$. We fix a vector $Y$ at
$\tau=\bar{\tau}$ with $|Y|_{g_{ij}(\bar \tau)}=1$, and extend $Y$
along the $\mathcal{L}$-shortest geodesic $\gamma$ on
$[0,\bar{\tau}]$ by solving the following ODE \be
\nabla_XY=-\Ric(Y,\cdot)+\frac{1}{2\tau}Y.  
\ee This is similar to the usual parallel translation and
multiplication with proportional parameter. Indeed, suppose
$\{Y_1,\ldots,Y_n\}$ is an orthonormal basis at $\tau=\bar{\tau}$
(with respect to the metric $g_{ij}(\bar \tau)$) and extend this
basis along the $\mathcal{L}$-shortest geodesic $\gamma$ by
solving the above ODE (3.2.8). Then
\begin{align*}
\frac{d}{d\tau}\<Y_i,Y_j\>
&  = 2\Ric(Y_i,Y_j)+\<\nabla_XY_i,Y_j\>+\<Y_i,\nabla_XY_j\>\\
&  = \frac{1}{\tau}\<Y_i,Y_j\>
\end{align*}
for all $i,j$. Hence, \be \<Y_i(\tau),Y_j(\tau)\>
=\frac{\tau}{\bar{\tau}}\delta_{ij}
\ee and $\{Y_1(\tau),\ldots,Y_n(\tau)\}$ remains orthogonal on
$[0,\bar{\tau}]$ with $Y_i(0)=0,\;i=1,\ldots,n$.

\begin{proposition} [{Perelman \cite{P1}}]
Given any unit vector $Y$ at any point $q \in M$ with
$\tau=\bar{\tau}$, consider an $\mathcal{L}$-shortest geodesic
$\gamma$ connecting $p$ to $q$ and extend $Y$ along $\gamma$ by
solving the {\rm ODE} $(3.2.8)$. Then the Hessian of the
$\mathcal{L}$-distance function $L$ on $M$ with $\tau=\bar{\tau}$
satisfies
$$
{\rm
Hess}_L(Y,Y)\leq\frac{1}{\sqrt{\bar{\tau}}}-2\sqrt{\bar{\tau}}\Ric(Y,Y)
-\int_0^{\bar{\tau}}\sqrt{{\tau}}Q(X,Y)d\tau
$$
in the sense of distributions, where
\begin{align*}
Q(X,Y)&  = -\nabla_Y\nabla_YR-2\<R(Y,X)Y,X\>
-4\nabla_X \Ric(Y,Y)+4\nabla_Y \Ric(Y,X)\\
&\quad -2\Ric_\tau(Y,Y)+2|\Ric(Y,\cdot)|^2-\frac{1}{\tau}\Ric(Y,Y)
\end{align*}
is the Li-Yau-Hamilton quadratic. Moreover the equality holds if
and only if the vector field $Y(\tau),\tau\in[0,\bar{\tau}]$, is
an {\bf $\mathcal{L}$-Jacobian field} $($i.e., $Y$ is the
derivative of a variation of $\gamma$ by
$\mathcal{L}$-geodesics$)$.\index{$\mathcal{L}$-Jacobian field}
\end{proposition}

\begin{pf}
As said before, we pretend that the shortest
$\mathcal{L}$-geodesics between $p$ and $q$ are unique so that
$L(q,\bar{\tau})$ is smooth. Otherwise, the inequality is still
valid, by a standard barrier argument, when understood in the
sense of distributions (see, for example, \cite{ScY}).

Recall that $\nabla L(q,\bar{\tau})=2\sqrt{\bar{\tau}}X$. Then
$\<\nabla_YY,\nabla L\>=2\sqrt{\bar{\tau}}\<\nabla_YY,X\>.$ We
compute by using Lemma 3.2.2, (3.2.8) and (3.2.9),
\begin{align*}
{\rm Hess}_L(Y,Y)&  = Y(Y(L))(\bar{\tau})-\langle\nabla_YY,\nabla
L\rangle(\bar\tau)\\
& \leq \delta_Y^2({\mathcal{L}})-2\sqrt{\bar{\tau}}
\langle\nabla_YY,X\rangle(\bar\tau)\\
&  = \int_0^{\bar{\tau}}\sqrt{\tau}[2 |\nabla_XY|^2
+2\langle R(Y,X)Y,X\rangle+\nabla_Y\nabla_YR\\
&\quad +2\nabla_X{\Ric}(Y,Y)-4\nabla_Y{\Ric}(Y,X)]d\tau\\
&  = \int_0^{\bar{\tau}}\!\!\sqrt{\tau}\bigg[2
\bigg|-\Ric(Y,\cdot)
+\frac{1}{2\tau}Y\bigg|^2\!+2\langle R(Y,X)Y,X\rangle+\nabla_Y\nabla_YR\\
&\quad +2\nabla_X \Ric(Y,Y)-4\nabla_Y \Ric(Y,X)\bigg]d\tau\\
& =\int_0^{\bar{\tau}}\sqrt{\tau}\bigg[2|\Ric(Y,\cdot)|^2
-\frac{2}{\tau}\Ric(Y,Y)
+\frac{1}{2\tau\bar{\tau}}+2\langle R(Y,X)Y,X\rangle\\
&\quad +\nabla_Y\nabla_YR+2\nabla_X \Ric(Y,Y)-4\nabla_Y
\Ric(Y,X)\bigg]d\tau.
\end{align*}
Since
\begin{align*}
\frac{d}{d\tau}\Ric(Y,Y)
&  =\Ric_\tau(Y,Y)+\nabla_X \Ric(Y,Y)+2\Ric(\nabla_XY,Y)\\
&  = \Ric_\tau(Y,Y)+\nabla_X
\Ric(Y,Y)-2|\Ric(Y,\cdot)|^2+\frac{1}{\tau}\Ric(Y,Y),
\end{align*}
we have
\begin{align*}
&{\rm Hess}_L(Y,Y) \\
& \leq\int_0^{\bar{\tau}}\sqrt{\tau}\bigg[2|\Ric(Y,\cdot)|^2
-\frac{2}{\tau}\Ric(Y,Y)+\frac{1}{2\tau\bar{\tau}}
+2\langle R(Y,X)Y,X\rangle\\
&\quad +\nabla_Y\nabla_YR-4(\nabla_Y \Ric)(X,Y)
-\bigg(2\frac{d}{d\tau}\Ric(Y,Y)-2\Ric_\tau(Y,Y)\\
&\quad +4|\Ric(Y,\cdot)|^2-\frac{2}{\tau}\Ric(Y,Y)\bigg)
+4\nabla_X \Ric(Y,Y)\bigg]d\tau\\
& =-\int_0^{\bar{\tau}}\bigg[2\sqrt{\tau}\frac{d}{d\tau}\Ric(Y,Y)
+\frac{1}{\sqrt{\tau}}\Ric(Y,Y)\bigg]d\tau+\frac{1}{2\bar\tau}
\int_0^{\bar{\tau}}\frac{1}{\sqrt{\tau}}d\tau\\
&\quad +\int_0^{\bar{\tau}}\sqrt{\tau}\bigg[2\langle
R(Y,X)Y,X\rangle
+\nabla_Y\nabla_YR+\frac{1}{\tau}\Ric(Y,Y)\\
&\quad +4(\nabla_X \Ric(Y,Y)-\nabla_Y \Ric(X,Y))+2\Ric_\tau(Y,Y)
-2|\Ric(Y,\cdot)|^2\bigg]d\tau\\
& =\frac{1}{\sqrt{\bar{\tau}}}-2\sqrt{\bar{\tau}}\Ric(Y,Y)
-\int_0^{\bar{\tau}}\sqrt{\tau}Q(X,Y)d\tau.
\end{align*}
This proves the inequality.

As usual, the quadratic form
$$
I(V,V)=\int_0^{\bar{\tau}}\sqrt{\tau}[2 |\nabla_XV|^2 +2\langle
R(V,X)V,X\rangle+\nabla_V\nabla_VR
$$
$$
+2\nabla_X \Ric(V,V)-4\nabla_V \Ric(V,X)]d\tau,
$$
for any vector field $V$ along $\gamma$, is called the index form.
Since $\gamma$ is shortest, the standard Dirichlet principle for
$I(V,V)$ implies that the equality holds if and only if the vector
field $Y$ is the derivative of a variation of $\gamma$ by
$\mathcal{L}$-geodesics.
\end{pf}

\begin{corollary} [{Perelman \cite{P1}}]
We have
$$
\Delta L\leq\frac{n}{\sqrt{\bar{\tau}}}
-2\sqrt{\bar{\tau}}R-\frac{1}{\bar{\tau}}K
$$
in the sense of distribution. Moreover, the equality holds if and
only if we are on a gradient shrinking soliton with
$$
R_{ij}+\frac{1}{2\sqrt{\bar\tau}}\nabla_i\nabla_jL
=\frac{1}{2\bar\tau}g_{ij}.$$
\end{corollary}

\begin{pf}
Choose an orthonormal basis $\{Y_1,\ldots,Y_n\}$ at
$\tau=\bar{\tau}$ and extend them along the shortest
$\mathcal{L}$-geodesic $\gamma$ to get vector fields $Y_i(\tau)$,
$i=1,\ldots,n$, by solving the ODE (3.2.8), with $\langle
Y_i(\tau),Y_j(\tau)\rangle=\frac{\tau}{\bar{\tau}}\delta_{ij}$ on
$[0,\bar{\tau}]$. Taking $Y=Y_i$ in Proposition 3.2.3 and summing
over $i$, we get
\begin{align}
\Delta L& \leq
\frac{n}{\sqrt{\bar{\tau}}}-2\sqrt{\bar{\tau}}R-\sum\limits_{i=1}^n
\int_0^{\bar{\tau}}\sqrt{\tau}Q(X,Y_i)d\tau \\
&  = \frac{n}{\sqrt{\bar{\tau}}}-2\sqrt{\bar{\tau}}R
-\int_0^{\bar{\tau}}\sqrt{\tau}\(\frac{\tau}{\bar{\tau}}\)Q(X)d\tau \nn\\
& =\frac{n}{\sqrt{\bar{\tau}}}-2\sqrt{\bar{\tau}}R
-\frac{1}{\bar{\tau}}K. \nn
\end{align} 
Moreover, by Proposition 3.2.3, the equality in (3.2.10) holds
everywhere if and only if for each $(q,\bar{\tau})$ and any
shortest $\mathcal{L}$-geodesic $\gamma$ on $[0,\bar{\tau}]$
connecting $p$ and $q$, and for any unit vector $Y$ at
$\tau=\bar{\tau}$, the extended vector field $Y(\tau)$ along
$\gamma$ by the ODE (3.2.8) must be an $\mathcal{L}$-Jacobian
field. When $Y_i(\tau),$ $i=1,\ldots,n$ are $\mathcal{L}$-Jacobian
fields along $\gamma$, we have
\begin{align*}
& \frac{d}{d\tau}\langle Y_i(\tau),Y_j(\tau)\rangle \\
&  = 2\Ric(Y_i,Y_j)+\langle\nabla_XY_i,Y_j\rangle
+\langle Y_i,\nabla_XY_j\rangle\\
& =2\Ric(Y_i,Y_j)+\left\langle\nabla_{Y_i}
\(\frac{1}{2\sqrt{\bar{\tau}}}\nabla L\),Y_j\right\rangle
+\left\langle Y_i,\nabla_{Y_j}
\(\frac{1}{2\sqrt{\bar{\tau}}}\nabla L\)\right\rangle\\
&  = 2\Ric(Y_i,Y_j)+\frac{1}{\sqrt{\bar{\tau}}}{\rm
Hess}_L(Y_i,Y_j)
\end{align*}
and then by (3.2.9),
$$
2\Ric(Y_i,Y_j)+\frac{1}{\sqrt{\bar{\tau}}}{\rm Hess}_L(Y_i,Y_j)
=\frac{1}{\bar{\tau}}\delta_{ij},\ \ \ \text{at}\ \
\tau=\bar{\tau}.
$$
Therefore the equality in (3.2.10) holds everywhere if and only if
we are on a gradient shrinking soliton with
$$
R_{ij}+\frac{1}{2\sqrt{\bar\tau}}\nabla_i\nabla_jL
=\frac{1}{2\bar\tau}g_{ij}.
$$
\end{pf}

In summary, from (3.2.6), (3.2.7) and Corollary 3.2.4, we have
$$   \left\{
       \begin{array}{lll}
  \frac{\partial L}{\partial \bar{\tau}}
=2\sqrt{\bar{\tau}}R-\frac{L}{2\bar{\tau}}
  +\frac{K}{\bar{\tau}},\\[4mm]
  |\nabla L|^2=-4\bar{\tau}R+\frac{2}{\sqrt{\bar{\tau}}}L
-\frac{4}{\sqrt{\bar{\tau}}}K,\\[4mm]
  \Delta
  L\leq-2\sqrt{\bar{\tau}}R+\frac{n}{\sqrt{\bar{\tau}}}
-\frac{K}{\bar{\tau}},
       \end{array}
    \right.
$$
in the sense of distributions.

Now the \textbf{Li-Yau-Perelman distance} \index{Li-Yau-Perelman
distance} $l=l(q,\bar{\tau})$ is defined by
$$l
(q,\bar{\tau})=L(q,\bar{\tau})/2\sqrt{\bar{\tau}}.
$$
We thus have the following

\begin{lemma} [{Perelman \cite{P1}}]
For the Li-Yau-Perelman distance $l(q,\bar{\tau})$ defined above,
we have
\begin{align}
\frac{\partial l}{\partial  \bar{\tau}}
&=-\frac{l}{\bar{\tau}}+R+\frac{1}{2\bar{\tau}^{3/2}}K,\\
|\nabla l|^2&=-R+\frac{l}{\bar{\tau}}-\frac{1}{\bar{\tau}^{3/2}}K,\\
\Delta l
&\leq-R+\frac{n}{2\bar{\tau}}-\frac{1}{2\bar{\tau}^{3/2}}K,
\end{align}
in the sense of distributions. Moreover, the equality in
$(3.2.13)$ holds if and only if we are on a gradient shrinking
soliton.
\end{lemma}

As the first consequence, we derive the following upper bound on
the minimum of $l(\cdot,\tau)$ for every $\tau$ which will be
useful in proving the no local collapsing theorem in the next
section.

\begin{corollary} [{Perelman \cite{P1}}]
Let $g_{ij}(\tau)$, $\tau\geq 0$, be a family of metrics evolving
by the Ricci flow $\frac{\partial}{\partial \tau}g_{ij}=2R_{ij}$
on a compact $n$-dimensional manifold $M$. Fix a point $p$ in $M$
and let $l(q,\tau)$ be the Li-Yau-Perelman distance from $(p,0)$.
Then for all $\tau$,
$$
\min\{l(q,\tau)\ |\ \ q\in M\}\leq\frac{n}{2}.
$$
\end{corollary}

\begin{pf}
Let
$$
\bar{L}(q,\tau)=4\tau l(q,\tau).
$$
Then, it follows from (3.2.11) and (3.2.13) that
$$
\frac{\partial \bar{L}}{\partial \tau}=4\tau
R+\frac{2K}{\sqrt{\tau}},
$$
and
$$
\Delta \bar{L}\leq-4\tau R+2n-\frac{2K}{\sqrt{\tau}}.
$$
Hence
$$
\frac{\partial \bar{L}}{\partial \tau}+\Delta \bar{L}\leq 2n.
$$
Thus, by a standard maximum principle argument,
$\min\{\bar{L}(q,\tau)-2n\tau| \ q\in M\}$ is nonincreasing and
therefore $\min\{\bar{L}(q,\tau)|\ \ q\in M\}\leq 2n\tau$.
\end{pf}

As another consequence of Lemma 3.2.5, we obtain
$$
\frac{\partial l}{\partial \bar{\tau}}-\Delta l+|\nabla l|^2-R
+\frac{n}{2\bar{\tau}}\geq 0.
$$
or equivalently
$$
\(\frac{\partial}{\partial \bar{\tau}}-\Delta +R\)
\((4\pi\bar{\tau})^{-\frac{n}{2}}\exp(-l)\)\leq 0.
$$
If $M$ is compact, we define {\bf Perelman's reduced
volume}\index{Perelman's reduced volume} by
$$
\tilde{V}(\tau)=\int_M(4\pi\tau)^{-\frac{n}{2}}
\exp(-l(q,\tau))dV_{\tau}(q),
$$
where $dV_{\tau}$ denotes the volume element with respect to the
metric $g_{ij}(\tau)$. Note that Perelman's reduced volume
resembles the expression in Huisken's monotonicity formula for the
mean curvature flow \cite{Hu90}. It follows, from the above
computation, that
\begin{align*}
& \frac{d}{d\bar{\tau}}\int_M(4\pi\bar{\tau})^{-\frac{n}{2}}
\exp(-l(q,\bar{\tau}))dV_{\bar{\tau}}(q)\\
&  = \int_M\left[\frac{\partial}{\partial \bar{\tau}}
((4\pi\bar{\tau})^{-\frac{n}{2}}\exp(-l(q,\bar{\tau})))
+R(4\pi\bar{\tau})^{-\frac{n}{2}}\exp(-l(q,\bar{\tau}))\right]
dV_{\bar{\tau}}(q)\\
&  \leq \int_M\Delta((4\pi\bar{\tau})^{-\frac{n}{2}}
\exp(-l(q,\bar{\tau})))dV_{\bar{\tau}}(q)\\
&  = 0.
\end{align*}
This says that if $M$ is compact, then Perelman's reduced volume
$\tilde{V}(\tau)$ is nonincreasing in $\tau$; moreover, the
monotonicity is strict unless we are on a gradient shrinking
soliton.

\medskip
In order to define and to obtain the monotonicity of  Perelman's
reduced volume for a complete noncompact manifold, we need to
formulate the monotonicity of Perelman's reduced volume in a local
version. This local version is very important and will play a
crucial role in the analysis of the Ricci flow with surgery in
Chapter 7.

\medskip
We define the {\bf $\mathcal{L}$-exponential map (with parameter
$\bar \tau$)}\index{$\mathcal{L}$-exponential map!with parameter
$\bar \tau$} ${\mathcal{L}}\exp(\bar \tau): T_pM \rightarrow M $
as follows: for any $X\in T_pM$, we set
$$
{\mathcal{L}}\exp_X(\bar{\tau})=\gamma(\bar{\tau})
$$
where $\gamma(\tau)$ is the ${\mathcal{L}}$-geodesic, starting at
$p$ and having $X$ as the limit of $\sqrt{\tau}\dot{\gamma}(\tau)$
as $\tau\rightarrow 0^+$. The associated Jacobian of the
$\mathcal{L}$-exponential map is called {\bf
$\mathcal{L}$-Jacobian}\index{$\mathcal{L}$-Jacobian}. We denote
by $\mathcal{J}(\tau)$ the $\mathcal{L}$-Jacobian  of
${\mathcal{L}}\exp(\tau):$ $T_pM\rightarrow M $. We can now deduce
an estimate for the $\mathcal{L}$-Jacobian as follows.

Let $q={\mathcal{L}}\exp_X(\bar{\tau})$ and $\gamma(\tau)$,
$\tau\in[0,\bar{\tau}]$, be the shortest ${\mathcal{L}}$-geodesic
connecting $p$ and $q$ with
$\sqrt{\tau}\dot{\gamma}(\tau)\rightarrow X$ as $\tau\rightarrow
0^+$. For any vector $v\in T_pM$, we consider the family of
${\mathcal{L}}$-geodesics:
$$
\gamma_s(\tau)={\mathcal{L}}\exp_{(X+sv)}(\tau),  \ \ \ \
0\leq\tau\leq\bar{\tau},\ \  s\in(-\epsilon,\epsilon).
$$
The associated variation vector field $V(\tau)$,
$0\leq\tau\leq\bar{\tau}$, is an $\mathcal{L}$-Jacobian field with
$V(0)=0$ and $V(\tau)=({\mathcal{L}}\exp_X(\tau))_*(v)$.

Let $v_1,\dots,v_n$ be $n$ linearly independent vectors in $T_pM$.
Then
$$
V_i(\tau)=({\mathcal{L}}\exp_X(\tau))_*(v_i),  \ \ \ \ i=1,2,
\ldots,n,
$$
are $n$ $\mathcal{L}$-Jacobian fields along $\gamma(\tau)$,
$\tau\in[0,\bar{\tau}]$. The $\mathcal{L}$-Jacobian
$\mathcal{J}(\tau)$ is given by
$$
{\mathcal{J}(\tau)}=|V_1(\tau)\wedge\dots\wedge
V_n(\tau)|_{g_{ij}(\tau)}/|v_1\wedge\dots\wedge v_n|.
$$

Now for fixed $b\in(0,\bar{\tau})$, we can choose linearly
independent vectors $v_1,\ldots,v_n\in T_pM$ such that $\langle
V_i(b),V_j(b)\rangle_{g_{ij}(b)}=\delta_{ij}$. We compute
\begin{align*}
&\frac{d}{d\tau}\mathcal{J}^2 \\
&  =\frac{2}{|v_1\wedge\cdots\wedge
v_n|^2}\sum\limits_{j=1}^n\langle V_1\wedge\cdots\wedge
\nabla_XV_j\wedge\cdots\wedge V_n,V_1\wedge\cdots\wedge
V_n\rangle_{g_{ij}(\tau)}\\
&\; +\frac{2}{|v_1\wedge\cdots\wedge
v_n|^2}\sum\limits_{j=1}^n\langle V_1\wedge\cdots\wedge
\Ric(V_j,\cdot)\wedge\cdots\wedge V_n,V_1\wedge\cdots\wedge
V_n\rangle_{g_{ij}(\tau)}.
\end{align*}
At $\tau=b$,
$$
\frac{d}{d\tau}\mathcal{J}^2(b)=\frac{2}{|v_1\wedge\cdots\wedge
v_n|^2}\sum\limits_{j=1}^n(\langle\nabla_XV_j,V_j\rangle_{g_{ij}(b)}
+\Ric(V_j,V_j)).
$$
Thus,
\begin{align*}
\frac{d}{d\tau}\log\mathcal{J}(b) &
=\sum\limits_{j=1}^n(\langle\nabla_XV_j,V_j\rangle_{g_{ij}(b)}
+\Ric(V_j,V_j))\\
& =\sum\limits_{j=1}^n\(\left\langle\nabla_{V_j}
\(\frac{1}{2\sqrt{b}}\nabla L\),V_j\right\rangle_{g_{ij}(b)}
+\Ric(V_j,V_j)\)\\
&  = \frac{1}{2\sqrt{b}}\(\sum\limits_{j=1}^n{\rm Hess}_L(V_j,V_j)\)+R\\
&  = \frac{1}{2\sqrt{b}}\Delta L+R.
\end{align*}
Therefore, in view of Corollary 3.2.4, we obtain the following
estimate for $\mathcal{L}$-Jacobian: \be
\frac{d}{d\tau}\log{\mathcal{J}}(\tau)
\leq\frac{n}{2\tau}-\frac{1}{2\tau^{3/2}}K\
\ \ \text{on} \ \ \ [0,\bar{\tau}].  
\ee On the other hand, by the definition of the Li-Yau-Perelman
distance and (3.2.4), we have
\begin{align}
\frac{d}{d\tau} l(\tau)
&  =-\frac{1}{2\tau}l+\frac{1}{2\sqrt{\tau}}\frac{d}{d\tau}L \\
&  =-\frac{1}{2\tau}l+\frac{1}{2\sqrt{\tau}}(\sqrt{\tau}(R+|X|^2)) \nn\\
&  =-\frac{1}{2\tau^{3/2}}K.\nn
\end{align}   
Here and in the following we denote by $l(\tau) =
l(\gamma(\tau),\tau)$.  Now the combination of (3.2.14) and
(3.2.15) implies the following important {\bf Jacobian comparison
theorem}\index{Jacobian comparison theorem} of Perelman \cite{P1}.

\begin{theorem}[Perelman's Jacobian comparison theorem]
Let $g_{ij}(\tau)$ be a family of complete solutions to the Ricci
flow $\frac{\partial}{\partial \tau}g_{ij}=2R_{ij}$ on a manifold
$M$ with bounded curvature. Let $\gamma:[0,\bar\tau]\rightarrow M$
be a shortest $\mathcal{L}$-geodesic starting from a fixed point
$p$. Then {\bf Perelman's reduced volume element}
$$
(4\pi\tau)^{-\frac{n}{2}}\exp(-l(\tau)){\mathcal{J}(\tau)}
$$
is nonincreasing in $\tau$ along $\gamma$.\index{Perelman's
reduced volume! element}
\end{theorem}

We now show how to integrate Perelman's reduced volume element
over $T_pM$ to deduce the following monotonicity result of
Perelman \cite{P1}.

\begin{theorem}[Monotonicity of Perelman's reduced volume]
Let $g_{ij}$ be a family of complete metrics evolving by the Ricci
flow $\frac{\partial}{\partial \tau}g_{ij}=2R_{ij}$ on a manifold
$M$ with bounded curvature. Fix a point $p$ in $M$ and let
$l(q,\tau)$ be the reduced distance from $(p,0)$. Then
\begin{itemize}
\item[(i)] Perelman's reduced volume
$$
\tilde{V}(\tau)
=\int_M(4\pi\tau)^{-\frac{n}{2}}\exp(-l(q,\tau))dV_{\tau}(q)
$$
is finite and nonincreasing in $\tau$; \item[(ii)] the
monotonicity is strict unless we are on a gradient shrinking
soliton.
\end{itemize}
\end{theorem}

\begin{pf}
For any $v\in T_pM$ we can find an $\mathcal{L}$-geodesic
$\gamma(\tau)$, starting at $p$, with
$\lim\limits_{\tau\rightarrow0^+}\sqrt{\tau}\dot{\gamma}(\tau)=v$.
Recall that $\gamma(\tau)$ satisfies the $\mathcal{L}$-geodesic
equation
$$
\nabla_{\dot{\gamma}(\tau)}\dot{\gamma}(\tau)-\frac{1}{2}\nabla
R+\frac{1}{2\tau}\dot{\gamma}(\tau)+2\Ric(\dot{\gamma}(\tau),\cdot)=0.
$$
Multiplying this equation by $\sqrt{\tau}$, we get \be
\frac{d}{d\tau}(\sqrt{\tau}\dot{\gamma})-\frac{1}{2}\sqrt{\tau}\nabla
R+2\Ric(\sqrt{\tau}\dot{\gamma}(\tau),\cdot)=0.
\ee Since the curvature of the metric $g_{ij}(\tau)$ is bounded,
it follows from Shi's derivative estimate (Theorem 1.4.1) that
$|\nabla R|$ is also bounded for small $\tau>0$. Thus by
integrating (3.2.16), we have \be
|\sqrt{\tau}\dot{\gamma}(\tau)-v|\leq C\tau(|v|+1) 
\ee for $\tau$ small enough and for some positive constant $C$
depending only the curvature bound.

Let $v_1,\ldots,v_n$ be $n$ linearly independent vectors in $T_pM$
and let
$$
V_i(\tau)=(\mathcal{L}\exp_v(\tau))_*(v_i)
=\frac{d}{ds}|_{s=0}\mathcal{L}\exp_{(v+sv_i)}(\tau),\ \
i=1,\ldots,n.
$$
The $\mathcal{L}$-Jacobian $\mathcal{J}(\tau)$ is given by
$$
\mathcal{J}(\tau)=|V_1(\tau)\wedge\cdots\wedge
V_n(\tau)|_{g_{ij}(\tau)}/|v_1\wedge\cdots\wedge v_n|
$$
By (3.2.17), we see that
$$
\left|\sqrt{\tau}\frac{d}{d\tau}\mathcal{L}
\exp_{(v+sv_i)}(\tau)-(v+sv_i)\right| \leq C\tau (|v|+|v_i|+1)
$$
for $\tau$ small enough and all $s\in(-\epsilon,\epsilon)$ (for
some $\epsilon>0$ small) and $i=1,\ldots,n$. This implies that
$$
\lim\limits_{\tau\rightarrow0^+}\sqrt{\tau}\dot{V}_i(\tau) =v_i,\
\ \ i=1,\ldots,n,
$$
so we deduce that \be
\lim\limits_{\tau\rightarrow0^+}\tau^{-\frac{n}{2}}\mathcal{J}(\tau)=1.
\ee Meanwhile, by using (3.2.17), we have
\begin{displaymath}
\begin{split}
l(\tau)&  =\frac{1}{2\sqrt{\tau}}\int_0^\tau\sqrt{\tau}
(R+|\dot{\gamma}(\tau)|^2)d\tau\\
&  \rightarrow |v|^2\ \ \text{as}\ \tau\rightarrow0^+.
\end{split}
\end{displaymath}
Thus \be
l(0)=|v|^2.
\ee Combining (3.2.18) and (3.2.19) with Theorem 3.2.7, we get
\begin{displaymath}
\begin{split}
\tilde{V}(\tau)
&  =\int_M(4\pi\tau)^{-\frac{n}{2}}\exp(-l(q,\tau))dV_{\tau}(q) \\
&  \leq\int_{T_pM}(4\pi\tau)^{-\frac{n}{2}}
\exp(-l(\tau))\mathcal{J}(\tau)|_{\tau=0}dv\\
&  =(4\pi)^{-\frac{n}{2}}\int_{\mathbb{R}^n}\exp(-|v|^2)dv\\
&  <+\infty.
\end{split}
\end{displaymath}
This proves that Perelman's reduced volume is always finite and
hence well defined. Now the monotonicity assertion in (i) follows
directly from Theorem 3.2.7.

For the assertion (ii), we note that the equality in (3.2.13)
holds everywhere when the monotonicity of  Perelman's reduced
volume is not strict. Therefore we have completed the proof of the
theorem.
\end{pf}

\section{No Local Collapsing Theorem I}

In this section we apply the monotonicity of Perelman's reduced
volume in Theorem 3.2.8 to prove Perelman's  {\bf no local
collapsing theorem I}\index{no local collapsing theorem I}(cf.
section 7.3 and section 4 of \cite{P1}), which is extremely
important not only because it gives a local injectivity radius
estimate in terms of local curvature bound but also it will
survive the surgeries in Chapter 7.

\begin{definition}
Let $\kappa$, $r$ be two positive constants and let $g_{ij}(t), 0
\leq t <T,$ be a solution to the Ricci flow on an $n$-dimensional
manifold $M$. We call the solution $g_{ij}(t)$
\textbf{$\kappa$-noncollapsed} at $(x_0,t_0)\in M \times [0,T)$ on
the scale $r$ \index{$\kappa$-noncollapsed} if it satisfies the
following property: whenever
$$
|Rm|(x,t)\leq r^{-2}
$$
for all $x \in B_{t_0}(x_0,r)$ and $t \in [t_0-r^2,t_0]$, we have
$$Vol_{t_0}(B_{t_0}(x_0,r)) \geq \kappa r^n.$$ Here $B_{t_0}(x_0,r)$ is
the geodesic ball centered at $x_0\in M $ and of radius $r$ with
respect to the metric $g_{ij}(t_0)$.
\end{definition}

Now we are ready to state the {\bf no local collapsing theorem
I}\index{no local collapsing theorem I} of Perelman \cite{P1}.

\begin{theorem}[No local collapsing theorem I]
Given any metric $g_{ij}$ on an $n$-dimensional compact manifold
$M$. Let $g_{ij}(t)$ be the solution to the Ricci flow on $[0,T)$,
with $T<+\infty$, starting at $g_{ij}$. Then there exist positive
constants $\kappa$ and $\rho_0$ such that for any $t_0\in[0,T)$
and any point $x_0\in M$, the solution $g_{ij}(t)$ is
$\kappa$-noncollapsed at $(x_0,t_0)$ on all scales less than
$\rho_0$.
\end{theorem}

\begin{pf}
We argue by contradiction. Suppose that there are sequences
$p_k\in M$, $t_k\in[0,T)$ and $r_k\rightarrow 0$ such that \be
|Rm|(x,t)\leq r_k^{-2} 
\ee for $x\in B_k=B_{t_k}(p_k,r_k)$ and $t_k-r_k^2\leq t\leq t_k$,
but \be \epsilon_k=r_k^{-1}Vol_{t_k}(B_k)^{\frac{1}{n}}\rightarrow
0 \
\ \text{as} \ k\rightarrow \infty. 
\ee

Without loss of generality, we may assume that $t_k\rightarrow T$
as $k\rightarrow +\infty$. \vskip 0.1cm \noindent Let
$\bar{\tau}(t)=t_k-t$, $p=p_k$ and
$$
\tilde{V_k}(\bar{\tau})=\int_M(4\pi\bar{\tau})^{-\frac{n}{2}}
\exp(-l(q,\bar{\tau}))dV_{t_k-\bar{\tau}}(q),
$$
where $l(q,\bar\tau)$ is the Li-Yau-Perelman distance with respect
to $p=p_k$.

\medskip
{\it Step} 1. \ We first want to show that for $k$ large enough,
$$
\tilde{V_k}(\epsilon_kr_k^2)\leq2\epsilon_k^\frac{n}{2}.
$$

For any $v\in T_pM$ we can find an $\mathcal{L}$-geodesic
$\gamma(\tau)$ starting at $p$ with $\lim\limits_{\tau\rightarrow
0}\sqrt{\tau}\dot{\gamma}(\tau)$ $=v$. Recall that $\gamma(\tau)$
satisfies the equation (3.2.16). It follows from assumption
(3.3.1) and Shi's local derivative estimate (Theorem 1.4.2) that
$|\nabla R|$ has a bound in the order of ${1}/{r_k^3}$ for $t\in
[t_k-\epsilon_kr_k^2,t_k]$. Thus by integrating (3.2.16) we see
that for $\tau\leq\epsilon_kr_k^2$ satisfying the property that
$\gamma(\sigma)\in B_k$ as long as $\sigma<\tau$, there holds \be
|\sqrt{\tau}\dot{\gamma}(\tau)-v|\leq C\epsilon_k(|v|+1)
\ee where $C$ is some positive constant depending only on the
dimension. Here we have implicitly used the fact that the metric
$g_{ij}(t)$ is equivalent for $x\in B_k$ and $t\in
[t_k-\epsilon_kr_k^2,t_k]$. In fact since $\frac{\partial
g_{ij}}{\partial t}=-2R_{ij}$ and $|Rm|\leq r_k^{-2}$ on
$B_k\times [t_k-r_k^2,t_k]$, we have \be
e^{-2\epsilon_k}g_{ij}(x,t_k)\leq g_{ij}(x,t)\leq
e^{2\epsilon_k}g_{ij}(x,t_k), 
\ee for $x\in B_k$ and $t\in [t_k-\epsilon_kr_k^2,t_k]$.\vskip
0.3cm Suppose $v\in T_pM$ with
$|v|\leq\frac{1}{4}\epsilon_k^{-\frac{1}{2}}$. Let $\tau\leq
\epsilon_kr_k^2$ such that $\gamma(\sigma)\in B_k$ as long as
$\sigma<\tau$, where $\gamma$ is the ${\mathcal{L}}$-geodesic
starting at $p$ with $\lim\limits_{\tau\rightarrow
0}\sqrt{\tau}\dot{\gamma}(\tau)=v$. Then, by (3.3.3) and (3.3.4),
for $k$ large enough,
\begin{align*}
d_{t_k}(p_k,\gamma(\tau))
&  \leq \int_0^\tau|\dot{\gamma}(\sigma)|_{g_{ij}(t_k)}d\sigma\\
&  < \frac{1}{2}\epsilon_k^{-\frac{1}{2}}
\int_0^\tau\frac{d\sigma}{\sqrt{\sigma}}\\
&  = \epsilon_k^{-\frac{1}{2}}\sqrt{\tau}\\
&  \leq r_k.
\end{align*}
This shows that for $k$ large enough, \be
{\mathcal{L}}\exp_{\{|v|\leq\frac{1}{4}\epsilon_k^{-1/2}\}}
(\epsilon_kr_k^2)\subset
B_k=B_{t_k}(p_k,r_k). 
\ee We now estimate the integral of $\tilde{V_k}(\epsilon_kr_k^2)$
as follows,
\begin{align}
\tilde{V_k}(\epsilon_kr_k^2) &
=\int_M(4\pi\epsilon_kr_k^2)^{-\frac{n}{2}}
\exp(-l(q,\epsilon_kr_k^2))dV_{t_k-\epsilon_kr_k^2}(q)\\
& =\int\limits_{{\mathcal{L}}\exp_{\{|v|\leq\frac{1}{4}
\epsilon_k^{-1/2}\}}(\epsilon_kr_k^2)}
(4\pi\epsilon_kr_k^2)^{-\frac{n}{2}}
\exp(-l(q,\epsilon_kr_k^2))dV_{t_k-\epsilon_kr_k^2}(q)\nn\\
&\quad +\!\int\limits_{M\setminus{\mathcal{L}}
\exp_{\{|v|\leq\frac{1}{4}\epsilon_k^{-1/2}\}}(\epsilon_kr_k^2)}\!\!
(4\pi\epsilon_kr_k^2)^{-\frac{n}{2}}
\exp(-l(q,\epsilon_kr_k^2))dV_{t_k-\epsilon_kr_k^2}(q).\nn
\end{align}   
We observe that for each $q\in B_k$,
$$
L(q,\epsilon_kr^2_k)
=\int_0^{\epsilon_kr_k^2}\sqrt{\tau}(R+|\dot{\gamma}|^2)d\tau
\geq-C(n)r_k^{-2}(\epsilon_kr_k^2)^{\frac{3}{2}}
=-C(n)\epsilon_k^{\frac{3}{2}}r_k,
$$
hence $l(q,\epsilon_kr^2_k)\geq -C(n)\epsilon_k.$ Thus, the first
term on the RHS of (3.3.6) can be estimated by
\begin{align}
&\int\limits_{\mathcal{L}\exp_{\{|v|\leq\frac{1}{4}
\epsilon_k^{-1/2}\}}(\epsilon_kr^2_k)}
(4\pi\epsilon_kr^2_k)^{-\frac{n}{2}}
\exp(-l(q,\epsilon_kr_k^2))dV_{t_k-\epsilon_kr^2_k}(q)\\
&\leq e^{n\epsilon_k}\int\limits_{B_k}
(4\pi\epsilon_kr^2_k)^{-\frac{n}{2}}
\exp(-l(q,\epsilon_kr_k^2))dV_{t_k}(q)\nn\\
&\leq e^{n\epsilon_k}(4\pi)^{-\frac{n}{2}}\cdot
e^{C(n)\epsilon_k}\cdot
\epsilon_k^{-\frac{n}{2}}\cdot(r_k^{-n}\Vol_{t_k}(B_k))\nn\\
&=e^{(n+C(n))\epsilon_k}(4\pi)^{-\frac{n}{2}}
\cdot\epsilon_k^\frac{n}{2},\nn
\end{align}  
where we have also used (3.3.5) and (3.3.4).

Meanwhile, by using (3.2.18), (3.2.19) and the Jacobian Comparison
Theorem 3.2.7, the second term on the RHS of (3.3.6) can be
estimated as follows
\begin{align}
&\int\limits_{M\setminus\mathcal{L}\exp_{\{|v|\leq\frac{1}{4}
\epsilon_k^{-\frac{1}{2}}\}}(\epsilon_kr^2_k)}
(4\pi\epsilon_kr^2_k)^{-\frac{n}{2}}
\exp(-l(q,\epsilon_kr_k^2))dV_{t_k-\epsilon_kr^2_k}(q)\\
&\leq \int\limits_{\{|v|>\frac{1}{4}\epsilon_k^{-\frac{1}{2}}\}}
(4\pi\tau)^{-\frac{n}{2}}
\exp(-l(\tau))\mathcal{J}(\tau)|_{\tau=0}dv \nn\\
&=(4\pi)^{-\frac{n}{2}}\int\limits_{\{|v|>\frac{1}{4}
\epsilon_k^{-\frac{1}{2}}\}}\exp(-|v|^2)dv \nn\\
&\leq \epsilon_k^{\frac{n}{2}},\nn
\end{align} 
for $k$ sufficiently large. Combining (3.3.6)-(3.3.8), we finish
the proof of Step 1.

\medskip
{\it Step} 2. \ We next want to show
$$
\tilde{V_k}(t_k)=(4\pi
t_k)^{-\frac{n}{2}}\int_M\exp(-l(q,t_k))dV_0(q)>C'
$$
for all $k$, where $C'$ is some positive constant independent of
$k$.

It suffices to show the Li-Yau-Perelman distance $l(\cdot, t_k)$
is uniformly bounded from above on $M$. By Corollary 3.2.6 we know
that the minimum of $l(\cdot,\tau)$ does not exceed $\frac{n}{2}$
for each $\tau>0$. Choose $q_k\in M$ such that the minimum of
$l(\cdot,t_k-\frac{T}{2})$ is attained at $q_k$. We now construct
a path $\gamma: [0,t_k]\rightarrow M$ connecting $p_k$ to any
given point $q\in M$ as follows: the first half path
$\gamma|_{[0,t_k-\frac{T}{2}]}$ connects $p_k$ to $q_k$ so that
$$
l\(q_k,t_k-\frac{T}{2}\) =\frac{1}{2\sqrt{t_k-\frac{T}{2}}}
\int_0^{t_k-\frac{T}{2}}\sqrt{\tau}
(R+|\dot{\gamma}(\tau)|^2)d\tau\leq\frac{n}{2}
$$
and the second half path $\gamma|_{[t_k-\frac{T}{2},t_k]}$ is a
shortest geodesic connecting $q_k$ to $q$ with respect to the
initial metric $g_{ij}(0)$. Then, for any $q\in M^n$,
\begin{align*}
l(q,t_k)&  = \frac{1}{2\sqrt{t_k}}L(q,t_k)\\
& \leq \frac{1}{2\sqrt{t_k}}\(\int_0^{t_k-\frac{T}{2}}
+\int_{t_k-\frac{T}{2}}^{t_k}\)
\sqrt{\tau}(R+|\dot{\gamma}(\tau)|^2)d\tau\\
& \leq \frac{1}{2\sqrt{t_k}}\(n\sqrt{t_k-\frac{T}{2}}
+\int_{t_k-\frac{T}{2}}^{t_k}
\sqrt{\tau}(R+|\dot{\gamma}(\tau)|^2)d\tau\)\\
&  \leq  C
\end{align*}
for some constant $C>0$, since all geometric quantities in
$g_{ij}$ are uniformly bounded when $t\in[0,\frac{T}{2}]$ (or
equivalently, $\tau\in[t_k-\frac{T}{2},t_k]$).

Combining Step 1 with Step 2, and using the monotonicity of
$\tilde{V_k}(\tau)$, we get
$$
C'<\tilde{V_k}(t_k)\leq\tilde{V_k}(\epsilon_kr_k^2) \leq
2\epsilon_k^\frac{n}{2}\to 0
$$
as $k\to \infty$. This gives the desired contradiction. Therefore
we have proved the theorem.
\end{pf}

The above no local collapsing theorem I says that if $|Rm|\leq
r^{-2}$ on the parabolic ball $\{(x,t)\ |\ d_{t_0}(x,x_0)\leq r, \
t_0-r^2\leq t\leq t_0\}$, then the volume of the geodesic ball
$B_{t_0}(x_0,r)$ (with respect to the metric $g_{ij}(t_0)$) is
bounded from below by $\kappa r^n$. In \cite{P1}, Perelman used
the monotonicity of the $\mathcal{W}$-functional (defined by
(1.5.9)) to obtain a stronger version of the no local collapsing
theorem, where the curvature bound assumption on the parabolic
ball is replaced by that on the geodesic ball $B_{t_0}(x_0,r)$.
The following result, called {\bf no local collapsing theorem
I$'$}\index{no local collapsing theorem I$'$}, gives a further
extension where the bound on the curvature tensor is replaced by
the bound on the scalar curvature only.

\begin{theorem}[No local collapsing theorem I$'$]
Suppose $M$ is a compact Riemannian manifold, and $g_{ij}(t)$,
$0\leq t<T<+\infty$, is a solution to the Ricci flow. Then there
exists a positive constant $\kappa$ depending only the initial
metric and $T$ such that for any $(x_0,t_0)\in M\times(0,T)$ if
$$
R(x,t_0)\leq r^{-2},\ \ \ \forall x\in B_{t_0}(x_0,r)
$$
with $0<r\leq\sqrt{T}$, then we have
$$
{\rm Vol}_{t_0}(B_{t_0}(x_0,r))\geq\kappa r^n.
$$
\end{theorem}

\begin{pf}
We will prove the assertion
$$
{\rm Vol}_{t_0}(B_{t_0}(x_0,a))\geq\kappa a^n\leqno{(\ast)_a}
$$
for all $0<a\leq r$. Recall that
$$
\mu(g_{ij},\tau)=\inf\left\{\mathcal{W}(g_{ij},f,\tau)\ \Big|\
\int_M(4\pi\tau)^{-\frac{n}{2}}e^{-f}dV=1\right\}.
$$
Set
$$
\mu_0=\inf\limits_{0\leq\tau\leq2T}\mu(g_{ij}(0),\tau)>-\infty.
$$
By Corollary 1.5.9, we have
\begin{align}
\mu(g_{ij}(t_0),b)&  \geq\mu(g_{ij}(0),t_0+b)\\
                  &  \geq\mu_0 \nn
\end{align}     
for $0<b\leq r^2$. Let $0<\zeta\leq1$ be a positive smooth
function on $\mathbb{R}$ where $\zeta(s)=1$ for
$|s|\leq\frac{1}{2}$, $|\zeta'|^2/\zeta\leq20$ everywhere, and
$\zeta(s)$ is very close to zero for $|s|\geq1$. Define a function
$f$ on $M$ by
$$
(4\pi r^2)^{-\frac{n}{2}}e^{-f(x)}=e^{-c}(4\pi
r^2)^{-\frac{n}{2}}\zeta\(\frac{d_{t_0}(x,x_0)}{r}\),
$$
where the constant $c$ is chosen so that $\int_M(4\pi
r^2)^{-\frac{n}{2}}e^{-f}dV_{t_0}=1$. Then it follows from (3.3.9)
that
\begin{align}
\mathcal{W}(g_{ij}(t_0),f,r^2) &  =\int_M[r^2(|\nabla f|^2
+R)+f-n]
(4\pi r^2)^{-\frac{n}{2}}e^{-f}dV_{t_0}\\
                 &  \geq\mu_0.\nn
\end{align} 
Note that
\begin{displaymath}
\begin{split}
   1&  =\int_M(4\pi
   r^2)^{-\frac{n}{2}}e^{-c}\zeta\(\frac{d_{t_0}(x,x_0)}{r}\)dV_{t_0}\\
    &  \geq\int_{B_{t_0}(x_0,\frac{r}{2})}(4\pi
    r^2)^{-\frac{n}{2}}e^{-c}dV_{t_0}\\
    &  =(4\pi
   r^2)^{-\frac{n}{2}}e^{-c}\Vol_{t_0}\(B_{t_0}\(x_0,\frac{r}{2}\)\).
\end{split}
\end{displaymath}
By combining with (3.3.10) and the scalar curvature bound, we have
\begin{displaymath}
\begin{split}
   c&  \geq-\int_M\(\frac{(\zeta')^2}{\zeta}
-\log\zeta\cdot\zeta\)e^{-c}(4\pi r^2)^{-\frac{n}{2}}
dV_{t_0}+(n-1)+\mu_0\\
    &  \geq-2(20+e^{-1})e^{-c}(4\pi
    r^2)^{-\frac{n}{2}}\Vol_{t_0}(B_{t_0}(x_0,r))+(n-1)+\mu_0\\
    &  \geq-2(20+e^{-1})
      \frac{{\rm Vol}_{t_0}(B_{t_0}(x_0,r))}{\Vol_{t_0}(B_{t_0}
(x_0,\frac{r}{2}))}+(n-1)+\mu_0,
\end{split}
\end{displaymath}
where we used the fact that $\zeta(s)$ is very close to zero for
$|s|\geq1$. Note also that
$$
2\int_{B_{t_0}(x_0,r)}e^{-c}(4\pi r^2)^{-\frac{n}{2}}dV_{t_0}
\geq\int_M(4\pi r^2)^{-\frac{n}{2}}e^{-f}dV_{t_0}=1.
$$
Let us set
$$\kappa=\min\left\{\frac{1}{2}\exp(-2(20+e^{-1})3^{-n}+(n-1)+\mu_0),\
\frac{1}{2}\alpha_n\right\}
$$
where $\alpha_n$ is the volume of the unit ball in $\mathbb{R}^n$.
Then we obtain
\begin{displaymath}
\begin{split}
  {\rm Vol}_{t_0}(B_{t_0}(x_0,r))
   &  \geq\frac{1}{2}e^c(4\pi r^2)^{\frac{n}{2}}\\
   &  \geq\frac{1}{2}(4\pi)^{\frac{n}{2}}
\exp(-2(20+e^{-1})3^{-n}+(n-1)+\mu_0)\cdot  r^n\\
   &  \geq\kappa r^n
\end{split}
\end{displaymath}
provided Vol$_{t_0}(B_{t_0}(x_0,\frac{r}{2}))
\geq3^{-n}Vol_{t_0}(B_{t_0}(x_0,r))$.

Note that the above argument also works for any smaller radius
$a\leq r$. Thus we have proved the following assertion: \be
{\rm Vol}_{t_0}(B_{t_0}(x_0,a))\geq\kappa a^n  
\ee whenever $a\in(0,r]$ and ${\rm
Vol}_{t_0}(B_{t_0}(x_0,\frac{a}{2})) \geq3^{-n}{\rm
Vol}_{t_0}(B_{t_0}(x_0,a))$.

Now we argue by contradiction to prove the assertion $(\ast)_a$
for any $a\in(0,r]$. Suppose $(\ast)_a$ fails for some
$a\in(0,r]$. Then by (3.3.11) we have
\begin{displaymath}
\begin{split}
   \Vol_{t_0}(B_{t_0}(x_0,\frac{a}{2}))
&  <3^{-n}\Vol_{t_0}(B_{t_0}(x_0,a))\\
   &  <3^{-n}\kappa a^n\\
   &  <\kappa\(\frac{a}{2}\)^n.
\end{split}
\end{displaymath}
This says that $(\ast)_{\frac{a}{2}}$ would also fail. By
induction, we deduce that
$$
\Vol_{t_0}\(B_{t_0}\(x_0,\frac{a}{2^k}\)\)<\kappa\(\frac{a}{2^k}\)^n\
\ \ \text{for all } k\geq1.
$$
This is a contradiction since $\lim\limits_{k\rightarrow\infty}
\Vol_{t_0}\(B_{t_0}\(x_0,\frac{a}{2^k}\)\)/\(\frac{a}{2^k}\)^n=\alpha_n$.
\end{pf}

\section{No Local Collapsing Theorem II}

By inspecting the arguments in the previous section, one can see
that if the injectivity radius of the initial metric is uniformly
bounded from below, then the no local collapsing theorem I of
Perelman also holds for complete solutions with bounded curvature
on a complete noncompact manifold. In this section we will use a
cut-off argument to extend the no local collapsing theorem to any
complete solution with bounded curvature. In some sense, the
second no local collapsing theorem of Perelman \cite{P1} gives a
good relative estimate of the volume element for the Ricci flow.

We first need the following useful lemma which contains two
assertions. The first one is a parabolic version of the Laplacian
comparison theorem (where the curvature sign restriction in the
ordinary Laplacian comparison is essentially removed in the Ricci
flow). The second one is a generalization of a result of Hamilton
(Theorem 17.2 in \cite{Ha95F}), where it was derived by an
integral version of Bonnet-Myers' theorem.

\begin{lemma}[{Perelman \cite{P1}}]
Let $g_{ij}(x,t)$ be a solution to the Ricci flow on an
$n$-dimensional manifold $M$ and denote by $d_t(x,x_0)$ the
distance between $x$ and $x_0$ with respect to the metric
$g_{ij}(t)$.
\begin{itemize}
\item[(i)] Suppose $\Ric(\cdot,t_0)\leq(n-1)K$ on
$B_{t_0}(x_0,r_0)$ for some $x_0\in M$ and some positive constants
K and $r_0$. Then the distance function $d(x,t)=d_t(x,x_0)$
satisfies, at $t=t_0$ and outside $B_{t_0}(x_0,r_0)$, the
differential inequality:
$$
\frac{\partial}{\partial t}d-\Delta d
\geq-(n-1)\(\frac{2}{3}Kr_0+r_0^{-1}\).
$$
\item[(ii)] Suppose $\Ric(\cdot,t_0)\leq (n-1)K$ on
$B_{t_0}(x_0,r_0)\bigcup B_{t_0}(x_1,r_0)$ for some $x_0,x_1\in M$
and some positive constants K and $r_0$.  Then, at $t=t_0$,
$$
\frac{d}{dt}d_t(x_0,x_1)\geq-2(n-1)\(\frac{2}{3}Kr_0+r_0^{-1}\).
$$
\end{itemize}
\end{lemma}

\begin{pf}
Let $\gamma: [0,d(x,t_0)]\rightarrow M$ be a shortest normal
geodesic from $x_0$ to $x$ with respect to the metric
$g_{ij}(t_0)$. As usual, we may assume that $x$ and $x_0$ are not
conjugate to each other in the metric $g_{ij}(t_0)$, otherwise we
can understand the differential inequality in the barrier sense.
Let $X=\dot{\gamma}(0)$ and let $\{X,e_1,\ldots,e_{n-1}\}$ be an
orthonormal basis of $T_{x_0}M$. Extend this basis parallel along
$\gamma$ to form a parallel orthonormal basis
$\{X(s),e_1(s),\ldots,e_{n-1}(s)\}$ along $\gamma$.

\medskip
(i) Let $X_i(s),\ i=1,\ldots,n-1$, be the Jacobian fields along
$\gamma$ such that $X_i(0)=0$ and $X_i(d(x,t_0))=e_i(d(x,t_0))$
for $i=1,\ldots,n-1$. Then it is well-known that (see for example
\cite{ScY})
$$
\Delta d_{t_0}(x,x_0)
=\sum_{i=1}^{n-1}\int_0^{d(x,t_0)}(|\dot{X}_i|^2-R(X,X_i,X,X_i))ds
$$
(in Proposition 3.2.3 we actually did this for the more
complicated $\mathcal{L}$-distance function).

Define vector fields $Y_i,\ i=1,\ldots,n-1$, along $\gamma$ as
follows:
$$
Y_i(s)=\begin{cases}
\frac{s}{r_0}e_i(s),&   \text{if }\;s\in[0,r_0],\\
e_i(s),&   \text{if }\;s\in[r_0,d(x,t_0)].\end{cases}
$$
which have the same value as the corresponding Jacobian fields
$X_i(s)$ at the two end points of $\gamma$. Then by using the
standard index comparison theorem (see for example \cite{CE}) we
have
\begin{align*}
\Delta d_{t_0}(x,x_0) &  =
\sum\limits_{i=1}^{n-1}\int_0^{d(x,t_0)}
(|\dot{X}_i|^2-R(X,X_i,X,X_i))ds\\
&  \leq \sum\limits_{i=1}^{n-1}\int_0^{d(x,t_0)}(|\dot{Y}_i|^2
-R(X,Y_i,X,Y_i))ds\\
&  = \int_0^{r_0}\frac{1}{r_0^2}(n-1-s^2\Ric(X,X))ds
+\int_{r_0}^{d(x,t_0)}(-\Ric(X,X))ds\\
&  = -\int_\gamma \Ric(X,X)+\int^{r_0}_0\(\frac{(n-1)}{r_0^2}
+\(1-\frac{s^2}{r_0^2}\)\Ric(X,X)\)ds\\
&  \leq -\int_\gamma \Ric(X,X)+(n-1)\(\frac{2}{3}Kr_0+r_0^{-1}\).
\end{align*}
On the other hand,
\begin{align*}
\frac{\partial}{\partial t}d_t(x,x_0) & = \frac{\partial}{\partial
t}
\int_0^{d(x,t_0)}\sqrt{g_{ij}X^iX^j}ds\\
&  = -\int_\gamma \Ric(X,X)ds.
\end{align*}
Hence we obtain the desired differential inequality. \vskip 3mm

\medskip
(ii) The proof is divided into three cases.

\medskip
{\it Case} (1): $d_{t_0}(x_0,x_1)\geq 2r_0$.

Define vector fields $Y_i,\ i=1,\ldots,n-1$, along $\gamma$ as
follows:
$$
Y_i(s)=\begin{cases}
\frac{s}{r_0}e_i(s),&   \mbox{ if }\ s\in[0,r_0],\\[4mm]
e_i(s),&   \mbox{ if }\ s\in[r_0,d(x_1,t_0)],\\[4mm]
\frac{d(x_1,t_0)-s}{r_0}e_i(s),&   \mbox{ if }\
s\in[d(x_1,t_0)-r_0,d(x_1,t_0)].\end{cases}
$$
Then by the second variation formula, we have
$$
\sum\limits_{i=1}^{n-1}\int_0^{d(x_1,t_0)}R(X,Y_i,X,Y_i)ds
\leq\sum\limits_{i=1}^{n-1}\int_0^{d(x_1,t_0)}|\dot{Y}_i|^2ds,
$$
which implies
\begin{align*}
& \int_0^{r_0}\frac{s^{2}}{r_0^2}\Ric(X,X)ds
+\int_{r_0}^{d(x,t_0)-r_0}\Ric(X,X)ds\\
&\quad+\int_{d(x_1,t_0)-r_0}^{d(x_1,t_0)}
\(\frac{d(x_1,t_0)-s}{r_0}\)^{2}\Ric(X,X)ds \leq
\frac{2(n-1)}{r_0}.
\end{align*}
Thus
\begin{align*}
& \frac{d}{dt}(d_t(x_0,x_1))\\
& \geq-\int_0^{r_0}\(1-\frac{s^{2}}{r_0^2}\)\Ric(X,X)ds \\
&\quad-\int_{d(x_1,t_0)-r_0}^{d(x_1,t_0)}
\(1-\(\frac{d(x_1,t_0)-s}{r_0}\)^{2}\)\Ric(X,X)ds
- \frac{2(n-1)}{r_0}\\
& \geq-2(n-1)\(\frac{2}{3}Kr_0+r_0^{-1}\).
\end{align*}

\smallskip
{\it Case} (2): $\frac{2}{\sqrt{\frac{2K}{3}}}\leq
d_{t_0}(x_0,x_1)\leq2r_0$.

In this case, letting $r_1=\frac{1}{\sqrt{\frac{2K}{3}}}$ and
applying case (1) with $r_0$ replaced by $r_1,$ we get
\begin{align*}
\frac{d}{dt}(d_t(x_0,x_1))
&  \geq -2(n-1)\(\frac{2}{3}Kr_1+r_1^{-1}\)\\
&  \geq -2(n-1)\(\frac{2}{3}Kr_0+r_0^{-1}\).
\end{align*}

\medskip
{\it Case} (3): $ d_{t_0}(x_0,x_1)\leq
\min\Big\{\frac{2}{\sqrt{\frac{2K}{3}}}, 2r_0\Big\}$.

In this case,
$$
\int_0^{d(x_1,t_0)}\Ric(X,X)ds \leq (n-1)K
\frac{2}{\sqrt{\frac{2K}{3}}}=(n-1)\sqrt{6K},
$$
and
$$
2(n-1)\(\frac{2}{3}Kr_0+r_0^{-1}\)\geq (n-1)\sqrt{\frac{32}{3}K}.
$$
This proves the lemma.
\end{pf}

The following result, called the {\bf no local collapsing theorem
II}\index{no local collapsing theorem II}, was obtained by
Perelman in \cite{P1}.


\begin{theorem}[No local collapsing theorem II]
For any $A >0$ there exists $\kappa=\kappa(A)>0$ with the
following property: if $g_{ij}(t)$ is a complete solution to the
Ricci flow on $0\leq t\leq r_0^2$ with bounded curvature and
satifying
$$
|Rm|(x,t)\leq r_0^{-2}\ \ \ on \ B_0(x_0,r_0)\times [0,r_0^2]
$$
and
$$
\Vol_0(B_0(x_0,r_0))\geq A^{-1}r_0^n,
$$
then $g_{ij}(t)$ is $\kappa$-noncollapsed on all scales less than
$r_0$ at every point $(x,r_0^2)$ with $d_{r_0^2}(x,x_0)\leq Ar_0$.
\end{theorem}

\begin{pf}
{}From the evolution equation of the Ricci flow, we know that the
metrics $g_{ij}(\cdot,t)$ are equivalent to each other on
$B_0(x_0,r_0)\times [0,r_0^2]$. Thus, without loss of generality,
we may assume that the curvature of the solution is uniformly
bounded for all $t \in [0,r_0^2]$ and all points in
$B_t(x_0,r_0)$. Fix a point $(x,r_0^2) \in M\times \{r_0^2\}$. By
scaling we may assume $r_0=1$. We may also assume $d_1(x,x_0)=A$.
Let $p=x$, $\bar{\tau}=1-t$, and consider Perelman's reduced
volume
$$
\tilde{V}(\bar{\tau}) =\int_M(4\pi\bar{\tau})^{-\frac{n}{2}}
\exp(-l(q,\bar{\tau}))dV_{1-\bar{\tau}}(q),
$$
where
\begin{multline*}
l(q,\bar{\tau})=\inf\bigg\{\frac{1}{2\sqrt{\bar{\tau}}}
\int_0^{\bar{\tau}}\sqrt{\tau}(R+|\dot{\gamma}|^2)d\tau
\mid\ \gamma: [0,\bar{\tau}]\rightarrow M\\
\mbox{with}\ \gamma(0)=p,\ \gamma(\bar{\tau})=q \bigg\}
\end{multline*}
is the Li-Yau-Perelman distance. We argue by contradiction.
Suppose for some $0<r<1$ we have
$$
|Rm|(y,t)\leq r^{-2}
$$
whenever $y\in B_{1}(x,r)$ and $1-r^2\leq t\leq 1$, but
$\epsilon=r^{-1}\Vol_{1}(B_1(x,r))^{\frac{1}{n}}$ is very small.
Then arguing as in the proof of the no local collapsing theorem I
(Theorem 3.3.2), we see that Perelman's reduced volume
$$
\tilde{V}(\epsilon r^2)\leq 2\epsilon^\frac{n}{2}
$$
On the other hand, from the monotonicity of Perelman's reduced
volume we have
$$
(4\pi)^{-\frac{n}{2}}\int_M\exp(-l(q,1))dV_0(q)
=\tilde{V}(1)\leq\tilde{V}(\epsilon r^2).
$$
Thus once we bound the function $l(q,1)$ over $B_0(x_0,1)$ from
above, we will get the desired contradiction and will prove the
theorem.

For any $q\in B_0(x_0,1)$, exactly as in the proof of the no local
collapsing theorem I, we choose a path $\gamma: [0,1]\rightarrow
M$ with $\gamma(0)=x,\ \gamma(1)=q$, $\gamma(\frac{1}{2})=y\in
B_\frac{1}{2}(x_0,\frac{1}{10})$ and $\gamma (\tau) \in
B_{1-\tau}(x_0,1)$ for $\tau \in [\frac{1}{2},1]$ such that
$$
{\mathcal{L}}(\gamma|_{[0,\frac{1}{2}]})
=2\sqrt{\frac{1}{2}}l\(y,\frac{1}{2}\)\ \ \(=L\(y,\frac{1}{2}\)\).
$$
Now ${\mathcal{L}}(\gamma|_{[\frac{1}{2},1]})
=\int_\frac{1}{2}^1\sqrt{\tau}(R(\gamma(\tau),1-\tau)+
|\dot{\gamma}(\tau)|^2_{g_{ij}(1-\tau)})d\tau$ is bounded from
above by a uniform constant since all geometric quantities in
$g_{ij}$ are uniformly bounded on $\{(y,t)\ |\ t\in[0,{1}/{2}],
y\in B_t(x_0,1)\}$ (where $t\in [0,1/2]$ is equivalent to
$\tau\in[{1}/{2},1]$).  Thus all we need is to estimate the
minimum of $l(\cdot,\frac{1}{2})$, or equivalently
$\bar{L}(\cdot,\frac{1}{2})=4\frac{1}{2}l(\cdot,\frac{1}{2})$, in
the ball $B_\frac{1}{2}(x_0,\frac{1}{10})$.

Recall that $\bar{L}$ satisfies the differential inequality \be
\frac{\partial \bar{L}}{\partial \tau}+\Delta\bar{L}\leq 2n.
\ee We will use this in a maximum principle argument. Let us
define
$$
h(y,t)=\phi(d(y,t)-A(2t-1))\cdot(\bar{L}(y,1-t)+2n+1)
$$
where $d(y,t)=d_t(y,x_0)$, and $\phi$ is a function of one
variable, equal to 1 on $(-\infty,\frac{1}{20})$, and rapidly
increasing to infinity on $(\frac{1}{20},\frac{1}{10})$ in such a
way that: \be
2\frac{(\phi')^2}{\phi}-\phi''\geq(2A+100n)\phi'-C(A)\phi
\ee for some constant $C(A)<+\infty$. The existence of such a
function $\phi$ can be justified as follows: put
$v=\frac{\phi'}{\phi}$, then the condition (3.4.2) for $\phi$ can
be written as
$$
3v^2-v'\geq(2A+100n)v-C(A)
$$
which can be solved for $v$.

Since the scalar curvature $R$ evolves by
$$
\frac{\partial R}{\partial t}=\Delta R+2|Rc|^2\geq\Delta
R+\frac{2}{n}R^2,
$$
we can apply the maximum principle as in Chapter 2 to deduce
$$
R(x,t)\geq-\frac{n}{2t}\ \ \text{for}\ t\in(0,1]\ \text{and} \
x\in M.
$$
Thus for $\bar{\tau}=1-t\in[0,\frac{1}{2}]$,
\begin{align*}
\bar{L}(\cdot,\bar{\tau}) &
=2\sqrt{\bar{\tau}}\int_0^{\bar{\tau}}
\sqrt{\tau}(R+|\dot{\gamma}|^2)d\tau\\
&  \geq 2\sqrt{\bar{\tau}}\int_0^{\bar{\tau}}
\sqrt{\tau}\(-\frac{n}{2(1-\tau)}\)d\tau\\
&  \geq 2\sqrt{\bar{\tau}}\int_0^{\bar{\tau}}
\sqrt{\tau}(-n)d\tau\\
&  > -2n.
\end{align*}
That is \be \bar{L}(\cdot,1-t)+2n+1\geq1,\ \ \text{for}\
t\in\left[\frac{1}{2},1\right]. 
\ee Clearly $\min\limits_{y\in M}h(y,\frac{1}{2})$ is achieved by
some $y\in B_{\frac{1}{2}}(x_0,\frac{1}{10})$ and \be
\min\limits_{y\in M}h(y,1)\leq h(x,1)=2n+1.  
\ee We compute
\begin{align*}
\(\frac{\partial}{\partial t}-\Delta\)h
& =\(\frac{\partial}{\partial t}-\Delta\)\phi\cdot(\bar{L}(y,1-t)+2n+1)\\
&\quad+\phi\cdot\(\frac{\partial}{\partial
t}-\Delta\)\bar{L}(y,1-t)
 -2\langle\nabla\phi,\nabla\bar{L}(y,1-t)\rangle\\
&  =\(\phi'\left[\(\frac{\partial}{\partial t}-\Delta\)d-2A\right]
-\phi''|\nabla d|^2\)\cdot(\bar{L}+2n+1)\\
&\quad +\phi\cdot\(-\frac{\partial}{\partial \tau}-\Delta\)
\bar{L}-2\langle\nabla\phi,\nabla\bar{L}\rangle\\
&  \geq \(\phi'\left[\(\frac{\partial}{\partial
t}-\Delta\)d-2A\right]-\phi''\)\cdot(\bar{L}+2n+1)\\
&\quad-2n\phi-2\langle\nabla\phi,\nabla\bar{L}\rangle
\end{align*}
by using (3.4.1). At a minimizing point of $h$ we have
$$
\frac{\nabla\phi}{\phi}=-\frac{\nabla\bar{L}}{(\bar{L}+2n+1)}.
$$
Hence
$$
-2\langle\nabla\phi,\nabla\bar{L}\rangle
=2\frac{|\nabla\phi|^2}{\phi}(\bar{L}+2n+1)
=2\frac{(\phi')^2}{\phi}(\bar{L}+2n+1).$$ Then at the minimizing
point of $h$, we compute
\begin{align*}
\(\frac{\partial}{\partial t}-\Delta\)h& \geq
\(\phi'\left[\(\frac{\partial}{\partial
t}-\Delta\)d-2A\right]-\phi''\)
\cdot(\bar{L}+2n+1) \\
&\quad-2n\phi+2\frac{(\phi')^2}{\phi}(\bar{L}+2n+1)\\
&  \geq \(\phi'\left[\(\frac{\partial}{\partial t}
-\Delta\)d-2A\right]-\phi''\)
\cdot(\bar{L}+2n+1) \\
&\quad-2nh+2\frac{(\phi')^2}{\phi}(\bar{L}+2n+1)
\end{align*}
for $t\in[\frac{1}{2},1]$ and
$$
\Delta h\geq0.
$$
Let us denote by $h_{\min}(t)=\min\limits_{y\in M}h(y,t)$. By
applying Lemma 3.4.1(i) to the set where $\phi'\neq0$, we further
obtain
\begin{align*}
\frac{d}{dt}h_{\min} &  \geq
(\bar{L}+2n+1)\cdot\left[\phi'(-100n-2A)
-\phi''+2\frac{(\phi')^2}{\phi}\right]-2nh_{\min}\\
&  \geq -(2n+C(A))h_{\min},\quad \text{for }\;
t\in[\frac{1}{2},1].
\end{align*}
This implies that $h_{\min}(t)$ cannot decrease too fast. By
combining (3.4.3) and (3.4.4) we get the required estimate for the
minimum $\bar{L}(\cdot,\frac{1}{2})$ in the ball
$B_{\frac{1}{2}}(x_0,\frac{1}{10})$.

Therefore we have completed the proof of the theorem.
\end{pf}

\newpage
\part{{\Large Formation of Singularities}}

\bigskip
Let $g_{ij}(x,t)$ be a solution to the Ricci flow on $M\times
[0,T)$ and suppose $[0,T)$, $T\le \infty$, is the maximal time
interval. If $T<+\infty$, then the short time existence theorems
in Section 1.2 tells us the curvature of the solution becomes
unbounded as $t\rightarrow T$ (cf. Theorem 8.1 of \cite{Ha95F}).
We then say the solution \textbf{develops a
singularity}\index{solution develops a singularity} as
$t\rightarrow T$. As in the minimal surface theory and harmonic
map theory, one usually tries to understand the structure of a
singularity of the Ricci flow by rescaling the solution (or blow
up) to obtain a sequence of solutions to the Ricci flow with
uniformly bounded curvature on compact subsets and looking at its
limit.

The main purpose of this chapter is to present the convergence
theorem of Hamilton \cite{Ha95} for a sequence of solutions to the
Ricci flow with uniform bounded curvature on compact subsets, and
to use the convergence theorem to give a rough classification in
\cite{Ha95F} for singularities of solutions to the Ricci flow.
Further studies on the structures of singularities of the Ricci
flow will be given in Chapter 6 and 7.

\section{Cheeger Type Compactness}

In this section, we establish Hamilton's compactness theorem
(Theorem 4.1.5) for solutions to the Ricci flow. The presentation
is based on Hamilton \cite{Ha95}.

We begin with the concept of $C_{loc}^{\infty}$ convergence of
tensors on a given manifold $M$. Let ${T_i}$ be a sequence of
tensors on $M$. We say that \mbox{\boldmath{$T_i$}} {\bf converges
to a tensor} \mbox{\boldmath{$T$}}\index{$T_i$ converges to a
tensor $T$} in the $C_{loc}^{\infty}$ topology if we can find a
covering $\{(U_s,\varphi_s)\}$, $\varphi_s: U_s\rightarrow
\mathbb{R}^n$, of $C^\infty$ coordinate charts so that for every
compact set $K\subset M$, the components of $T_i$ converge in the
$C^\infty$ topology to the components of $T$ in the intersections
of $K$ with these coordinate charts, considered as functions on
$\varphi_s(U_s)\subset \mathbb{R}^n$. Consider a Riemannian
manifold $(M,g)$. A {\bf marking}\index{marking} on $M$ is a
choice of a point $p\in M$ which we call the {\bf
origin}\index{origin}. We will refer to such a triple $(M,g,p)$ as
a {\bf marked Riemannian manifold}.\index{marked Riemannian
manifold}

\begin{definition}
Let $(M_k,g_k,p_k)$ be a sequence of marked complete Riemannian
manifolds, with metrics $g_k$ and marked points $p_k\in M_k$. Let
$B(p_k, s_k)\subset M_k$ denote the geodesic ball centered at
$p_k\in M_k$ and of radius $s_k$ ($0 < s_k \leq +\infty$). We say
a sequence of marked geodesic balls $(B(p_k,s_k),g_k, p_k)$ with
$s_k \rightarrow s_{\infty}(\leq +\infty)$ {\bf converges} in the
$C_{loc}^{\infty}$ topology {\bf to a marked} (maybe noncomplete)
{\bf manifold}\index{converges to a marked manifold}
$(B_{\infty},g_{\infty},p_{\infty})$, which is an open geodesic
ball centered at $p_{\infty}\in B_{\infty}$ and of radius
$s_{\infty}$ with respect to the metric $g_{\infty}$, if we can
find a sequence of exhausting open sets $U_k$ in $B_{\infty}$
containing $p_{\infty}$ and a sequence of diffeomorphisms $f_k$ of
the sets $U_k$ in $B_{\infty}$ to open sets $V_k$ in $B(p_k,s_k)
\subset M_k$ mapping $p_{\infty}$ to $p_k$ such that the pull-back
metrics $\tilde{g}_k=(f_k)^*g_k$ converge in $C^\infty$ topology
to $g_{\infty}$ on every compact subset of $B_{\infty}$.
\end{definition}

We remark that this concept of $C_{\rm loc}^\infty$-convergence of
a sequence of marked manifolds $(M_k,g_k,p_k)$ is not the same as
that of $C_{\rm loc}^\infty$-convergence of metric tensors on a
given manifold, even when we are considering the sequence of
Riemannian metric $g_k$ on the same space $M$. This is because one
can have a sequence of diffeomorphisms $f_k: M\rightarrow M$ such
that $(f_k)^*g_k$ converges in $C^\infty_{\rm loc}$ topology while
$g_k$ itself does not converge.

There have been a lot of work in Riemannian geometry on the
convergence of a sequence of compact manifolds with bounded
curvature, diameter and injectivity radius (see for example Gromov
\cite{Gro}, Peters \cite{Pe}, and Greene and Wu \cite{Gre}). The
following theorem, which is a slight generalization of Hamilton's
convergence theorem \cite{Ha95}, modifies these results in three
aspects: the first one is to allow noncompact limits and then to
avoid any diameter bound; the second one is to avoid having to
assume a uniform lower bound for the injectivity radius over the
whole manifold, a hypothesis which is much harder to satisfy in
applications; the last one is to avoid a uniform curvature bound
over the whole manifold so that we can take a local limit.

\begin{theorem}[{Hamilton\; \cite{Ha95}}]
Let\; $(M_k,\,g_k,\,p_k)$\; be\; a\; sequence\; of marked complete
Riemannian manifolds of dimension $n$. Consider a sequence of
geodesic balls $B(p_k,s_k) \subset M_k$ of radius $s_k$ $(0<s_k\le
\infty)$, with $s_k \rightarrow s_{\infty}(\le \infty)$, around
the base point $p_k$ of $M_k$ in the metric $g_k$. Suppose
\begin{itemize}
\item[(a)] for every radius $r<s_{\infty}$ and every integer $l\ge
0$ there exists a constant $B_{l,r}$, independent of $k$, and
positive integer $k(r,l) < +\infty$ such that as $k \geq k(r,l)$,
the curvature tensors $Rm(g_k)$ of the metrics $g_k$ and their
$l^{\mbox{th}}$-covariant derivatives satisfy the bounds
$$
|\nabla^lRm(g_k)|\leq B_{l,r}
$$
on the balls $B(p_k, r)$ of radius $r$ around $p_k$ in the metrics
$g_k$; and \item[(b)] there exists a constant $\delta>0$
independent of $k$ such that the injectivity radii
inj$\,(M_k,p_k)$ of $M_k$ at $p_k$ in the metric $g_k$ satisfy the
bound
$$
{\rm inj}\,(M_k,p_k)\geq\delta.$$
\end{itemize}
Then there exists a subsequence of the marked geodesic balls
$(B(p_k,\,s_k),\,g_k,$ $p_k)$ which converges to a marked geodesic
ball $(B(p_{\infty},s_{\infty}),g_{\infty},p_{\infty})$ in
$C^\infty_{\rm loc}$ topology. Moreover the limit is complete if
$s_{\infty} = +\infty$.

\end{theorem}

\begin{pf}
In \cite{Ha95} (see also Theorem 16.1 of \cite{Ha95F}), Hamilton
proved the above convergence theorem for the case $s_{\infty} =
+\infty$. Thus it remains to prove the case of $s_{\infty} <
+\infty$. Our proof here follows, with some slight modifications,
the argument of Hamilton \cite{Ha95}. In fact, except Step 1, the
proof is essentially taken from Hamilton \cite{Ha95}. Suppose we
are given a sequence of geodesic balls $(B(p_k,s_k),g_k, p_k)
\subset (M_k,g_k,p_k)$, with $s_k \rightarrow s_{\infty} (<
+\infty)$, satisfying the assumptions of Theorem 4.1.2. We will
split the proof into three steps.

\medskip
{\it Step} 1: \ Picking the subsequence.

\smallskip
By the local injectivity radius estimate (4.2.2) in Corollary
4.2.3 of the next section, we can find a positive decreasing $C^1$
function $\rho (r)$, $0 \leq r < s_{\infty}$, independent of $k$
such that
\begin{align}
\rho(r) &< \frac{1}{100}(s_{\infty} - r), \\  
0 &\geq \rho'(r) \geq -\frac{1}{1000}, 
\end{align}
and a sequence of positive constants $\varepsilon_k \rightarrow 0$
so that the injectivity radius at any point $x\in B(p_k,s_k)$ with
$r_k = d_k(x,p_k) < s_{\infty} - \varepsilon_k$ is bounded from
below by \be
{\rm inj}\, (M_k, x)\ge 500\rho (r_k(x)),  
\ee where $r_k(x)=d_k(x, p_k)$ is the distance from $x$ to $p_k$
in the metric $g_k$ of $M_k$. We define
$$
\tilde\rho(r)=\rho(r+20\rho(r)),\quad
\tilde{\tilde\rho}(r)=\tilde\rho(r+20\tilde\rho(r)).
$$
By (4.1.2) we know that both $\tilde{\rho}(r)$ and
$\tilde{\tilde{\rho}}(r)$ are nonincreasing positive functions on
$[0,s_{\infty})$.

In each $B(p_k,s_{\infty})$ we choose inductively a sequence of
points $x_k^{\alpha}$ for $\alpha=0,1,2,\ldots$ in the following
way. First we let $x_k^0=p_k$. Once $x_k^{\alpha}$ are chosen for
$\alpha=0,1,2,\ldots, \sigma$, we pick $x_k^{\sigma+1}$ closest to
$p_k$ so that $r_k^{\sigma+1}=r_k(x_k^{\sigma+1})$ is as small as
possible, subject to the requirement that the open ball
$B(x_k^{\sigma+1}, \tilde{\tilde\rho}_k^{\sigma+1})$ around
$x_k^{\sigma+1}$ of radius $\tilde{\tilde\rho}_k^{\sigma+1}$ is
disjoint from the balls $B(x_k^{\alpha},
\tilde{\tilde\rho}_k^{\alpha})$ for $\alpha=0,1,2,\ldots, \sigma$,
where
$\tilde{\tilde\rho}_k^{\alpha}=\tilde{\tilde\rho}(r_k^{\alpha})$
and $r_k^{\alpha}=r_k(x_k^{\alpha})$. In particular, the open
balls $B(x_k^{\alpha}, \tilde{\tilde\rho}_k^{\alpha}),
\alpha=0,1,2,\ldots$, are all disjoint. We claim the balls
$B(x_k^{\alpha}, 2\tilde{\tilde\rho}_k^{\alpha})$ cover
$B(p_k,s_{\infty}-\varepsilon_k)$ and moreover for any $ r$, $0< r
<s_{\infty} -\varepsilon_k$, we can find $\lambda ( r)$
independent of $k$ such that for $k$ large enough, the geodesic
balls $B(x_k^{\alpha}, 2\tilde{\tilde\rho}_k^{\alpha})$ for
$\alpha\le\lambda ( r)$ cover the ball $B(p_k,r)$.

To see this, let $x\in B(p_k,s_{\infty}-\varepsilon_k)$ and let
$r(x)$ be the distance from $x$ to $p_k$ and let
$\tilde{\tilde\rho}=\tilde{\tilde\rho}(r(x))$. Consider those
$\alpha$ with $r_k^{\alpha}\le r(x) < s_{\infty} - \varepsilon_k$.
Then
$$
\tilde{\tilde\rho}\le \tilde{\tilde\rho}_k^{\alpha}.
$$
Now the given point $x$ must lie in one of the balls
$B(x_k^{\alpha}, 2\tilde{\tilde\rho}_k^{\alpha})$. If not, we
could choose the next point in the sequence of $x_k^{\beta}$ to be
$x$ instead, for since
$\tilde{\tilde\rho}_k^{\alpha}+\tilde{\tilde\rho}\le
2\tilde{\tilde\rho}_k^{\alpha}$ the ball $B(x,
\tilde{\tilde\rho})$ would miss $B(x_k^{\alpha},
\tilde{\tilde\rho}_k^{\alpha})$ with $r_k^\alpha\le r(x)$. But
this is a contradiction. Moreover for any $ r$, $0< r <s_{\infty}
-\varepsilon_k$, using the curvature bound and the injectivity
radius bound, each ball $B(x_k^{\alpha},
\tilde{\tilde\rho}_k^{\alpha})$ with $r_k^{\alpha}\le r$ has
volume at least $\epsilon(r)\tilde{\tilde\rho}^n$ where
$\epsilon(r)>0$ is some constant depending on $r$ but independent
of $k$. Now these balls are all disjoint and contained in the ball
$B(p_k, (r+s_{\infty})/2)$. On the other hand, for large enough
$k$, we can estimate the volume of this ball from above, again
using the curvature bound, by a positive function of $r$ that is
independent of $k$. Thus there is a $k'(r) > 0$ such that for each
$k \geq k'(r)$, there holds \be \# \{ \alpha \ |\ r_k^{\alpha}\le
r \} \leq \lambda(r)
\ee for some positive constant $\lambda(r)$ depending only on $r$,
and the geodesic balls $B(x_k^{\alpha},
2\tilde{\tilde\rho}_k^{\alpha})$ for $\alpha\le\lambda ( r)$ cover
the ball $B(p_k,r)$.

By the way, since
\begin{align*}
r_k^\alpha &  \leq r_k^{\alpha -1} +
\tilde{\tilde{\rho}}_k^{\alpha -1} +\tilde{\tilde{\rho}}_k^{\alpha} \\
&  \leq r_k^{\alpha -1} + 2\tilde{\tilde{\rho}}_k^{\alpha -1},
\end{align*}
and by (4.1.1)
$$
\tilde{\tilde{\rho}}_k^{\alpha -1} \leq \frac{1}{100}(s_{\infty}
-r_k^{\alpha -1}),
$$
we get by induction
\begin{align*}
r_k^\alpha & \leq \frac{49}{50}r_k^{\alpha -1} +\frac{1}{50}s_{\infty} \\
&  \leq \(\frac{49}{50}\)^{\alpha}r_k^{0}
+\frac{1}{50}\(1+\frac{49}{50}+ \cdots
+ \(\frac{49}{50}\)^{\alpha-1}\)s_{\infty} \\
&  = \(1-\(\frac{49}{50}\)^{\alpha}\)s_{\infty}.
\end{align*}
So for each $\alpha$, with $\alpha \leq \lambda(r)$
$(r<s_{\infty})$, there holds \be r_k^\alpha \leq
\(1-\(\frac{49}{50}\)^{\lambda(r)}\)s_{\infty}
\ee for all $k$. And by passing to a subsequence (using a
diagonalization argument) we may assume that $r_k^{\alpha}$
converges to some $r^{\alpha}$ for each ${\alpha}$. Then
$\tilde{\tilde\rho}_k^{\alpha}$ (respectively
$\tilde\rho_k^{\alpha}, \rho_k^{\alpha}$) converges to
$\tilde{\tilde\rho}^{\alpha}=\tilde{\tilde\rho}(r^{\alpha})$
(respectively $\tilde\rho^{\alpha}=\tilde\rho(r^{\alpha}),
\rho^{\alpha}=\rho(r^{\alpha})$).

Hence for all $\alpha$ we can find $k({\alpha})$ such that
$$
\frac{1}{2}\tilde{\tilde\rho}^{\alpha}
\le\tilde{\tilde\rho}_k^{\alpha} \le 2\tilde{\tilde\rho}^{\alpha}
$$
$$
\frac{1}{2}\tilde\rho^{\alpha}\le \tilde\rho_k^{\alpha}\le
2\tilde\rho^{\alpha} \quad \mbox{and} \quad
\frac{1}{2}\rho^{\alpha}\le \rho_k^{\alpha}\le 2\rho^{\alpha}
$$
whenever $k\ge k(\alpha)$. Thus for all $\alpha$,
$\tilde{\tilde\rho}_k^{\alpha}$ and $\tilde{\tilde\rho}^{\alpha}$
are comparable when $k$ is large enough so we can work with balls
of a uniform size, and the same is true for
$\tilde\rho_k^{\alpha}$ and $\tilde\rho^{\alpha}$, and
$\rho_k^{\alpha}$ and $\rho^{\alpha}$. Let $\hat
B_k^{\alpha}=B(x_k^{\alpha}, 4\tilde{\tilde\rho}^{\alpha})$, then
$\tilde{\tilde\rho}_k^{\alpha}\le 2\tilde{\tilde\rho}^{\alpha}$
and $B(x_k^{\alpha}, 2\tilde{\tilde\rho}_k^{\alpha})\subset
B(x_k^{\alpha}, 4\tilde{\tilde\rho}^{\alpha})=\hat B_k^{\alpha}$.
So for every $r$ if we let $k(r)=\max\{k(\alpha) \  |\ \alpha\le
\lambda(r)\}$ then when $k\ge k(r)$, the balls $\hat B_k^{\alpha}$
for $\alpha\le \lambda(r)$ cover the ball $B(p_k, r)$ as well.
Suppose that $\hat B_k^{\alpha}$ and $\hat B_k^{\beta}$ meet for
$k\ge k({\alpha})$ and $k\ge k({\beta})$, and suppose
$r_k^{\beta}\le r_k^{\alpha}$. Then, by the triangle inequality,
we must have
$$
r_k^{\alpha}\le r_k^{\beta}
+4\tilde{\tilde\rho}^{\alpha}+4\tilde{\tilde\rho}^{\beta}\le
r_k^{\beta} +8\tilde{\tilde\rho}^{\beta}<r_k^{\beta}
+16\tilde\rho_k^{\beta}.
$$
This then implies
$$
\tilde{\tilde\rho}_k^{\beta}=\tilde{\tilde\rho}(r_k^{\beta})=
\tilde\rho(r_k^{\beta}+20\tilde\rho(r_k^{\beta}))
<\tilde\rho(r_k^{\alpha})=\tilde\rho_k^{\alpha}
$$
and hence
$$
\tilde{\tilde\rho}^{\beta}\le 4\tilde\rho^{\alpha}.
$$
Therefore $\hat B_k^{\beta}\subset B(x_k^{\alpha},
36\tilde\rho^{\alpha})$ whenever $\hat B_k^{\alpha}$ and $\hat
B_k^{\beta}$ meet and $k\ge \max\{ k({\alpha}),$ $k({\beta})\}.$

Next we define the balls $B_k^{\alpha}=B(x_k^{\alpha},
5\tilde{\tilde\rho}^{\alpha})$ and $\tilde
B_k^{\alpha}=B(x_k^{\alpha}, \tilde{\tilde\rho}^{\alpha}/2)$. Note
that $\tilde B_k^{\alpha}$ are disjoint since  $\tilde
B_k^{\alpha}\subset B(x_k^{\alpha},
\tilde{\tilde{\rho}}^{\alpha}_k)$. Since $\hat B_k^{\alpha}\subset
B_k^{\alpha}$, the balls $B_k^{\alpha}$ cover $B(p_k, r)$ for
$\alpha\le \lambda(r)$ as before. If $B_k^{\alpha}$ and
$B_k^{\beta}$ meet for $k\ge k({\alpha})$ and $k\ge k({\beta})$
and $r_k^{\beta}\le r_k^{\alpha}$, then by the triangle inequality
we get
$$
r_k^{\alpha}\le r_k^{\beta}
+10\tilde{\tilde\rho}^{\beta}<r_k^{\beta} +20\tilde\rho_k^{\beta},
$$
and hence
$$
\tilde{\tilde\rho}^{\beta}\le 4\tilde\rho^{\alpha}
$$
again. Similarly,
$$
\tilde\rho_k^{\beta}=\tilde\rho(r_k^{\beta})=
\rho(r_k^{\beta}+20\rho(r_k^{\beta}))
<\rho(r_k^{\alpha})=\rho_k^{\alpha}.
$$
This makes
$$
\tilde\rho^{\beta}\le 4\rho^{\alpha}.
$$
Now any point in $B_k^{\beta}$ has distance at most
$$
5\tilde{\tilde\rho}^{\alpha}+5\tilde{\tilde\rho}^{\beta}
+5\tilde{\tilde\rho}^{\beta}\le 45\tilde\rho^{\alpha}
$$
from $x_k^{\alpha}$, so $B_k^{\beta}\subset B(x_k^{\alpha},
45\tilde\rho^{\alpha})$.  Likewise, whenever $B_k^{\alpha}$ and
$B_k^{\beta}$ meet for $k\ge k({\alpha})$ and $k\ge k({\beta})$,
any point in the larger ball $B(x_k^{\beta},
45\tilde\rho^{\beta})$ has distance at most
$$
5\tilde{\tilde\rho}^{\alpha}+5\tilde{\tilde\rho}^{\beta}
+45\tilde\rho^{\beta}\le 205\rho^{\alpha}
$$
from $x_k^{\alpha}$ and hence $B(x_k^{\beta},
45\tilde\rho^{\beta})\subset B(x_k^{\alpha}, 205\rho^{\alpha})$.
Now we define ${\bar B}_k^{\alpha}=B(x_k^{\alpha},
45\tilde\rho^{\alpha})$ and ${\bar{\bar
B}}_k^{\alpha}=B(x_k^{\alpha}, 205\rho^{\alpha})$. Then the above
discussion says that whenever $B_k^{\alpha}$ and $B_k^{\beta}$
meet for $k\ge k({\alpha})$ and $k\ge k({\beta})$, we have \be
B_k^{\beta}\subset {\bar B}_k^{\alpha} \quad \mbox{and}\quad {\bar
B}_k^{\beta}\subset {\bar{\bar
B}}_k^{\alpha}. 
\ee Note that ${\bar{\bar B}}_k^{\alpha}$ is still a nice embedded
ball since, by (4.1.3), $205\rho^{\alpha}\le
410\rho_k^{\alpha}<inj(M_k, x_k^{\alpha}). $

We claim there exist positive numbers $N(r)$ and $k''(r)$ such
that for any given $\alpha$ with $r^{\alpha}<r$, as $k \geq
k''(r)$, there holds \be \# \{\beta \ |\ B_k^{\alpha} \cap
B_k^{\beta} \neq \phi \}
\leq N(r). 
\ee

Indeed, if $B_k^{\alpha}$ meets $B_k^{\beta}$ then there is a
positive $k''(\alpha)$ such that as $k \geq k''(\alpha)$,
\begin{align*}
r_k^\beta
&  \leq r_k^{\alpha} + 10\tilde{\tilde{\rho}}_k^{\alpha }\\
&  \leq r + 20\rho(r)\\
&  \leq r + \frac{1}{5}(s_{\infty}-r),
\end{align*}
where we used (4.1.2) in the third inequality. Set
$$
k''(r)=\max\{ k''(\alpha), k'(r)\ |\ \alpha \leq \lambda(r)\}
$$ and
$$
N(r)= \lambda\(r + \frac{1}{5}(s_{\infty}-r)\).
$$
Then by combining with (4.1.4), these give the desired estimate
(4.1.7)

Next we observe that by passing to another subsequence we can
guarantee that for any pair $\alpha$ and $\beta$ we can find a
number $k(\alpha,\beta)$ such that if $k\ge k(\alpha,\beta)$ then
either $B_k^{\alpha}$ always meets $B_k^{\beta}$ or it never does.

Hence by setting
\begin{multline*}
\bar{k}(r) = \max \bigg\{ k(\alpha,\beta), k(\alpha), k(\beta),
k''(r) \ |\
\alpha \leq\lambda(r)\; \mbox{ and}\\
\beta \leq \lambda\(r +\frac{1}{5}(s_{\infty}-r)\) \bigg\},
\end{multline*}
we have shown the following results: for every $r < s_{\infty}$,
if $k\ge \bar{k}(r)$, we have
\begin{itemize}
\item[(i)] the ball $B(p_k, r)$ in $M_k$ is covered by the balls
$B_k^{\alpha}$ for $\alpha\le \lambda(r)$, \item[(ii)] whenever
$B_k^{\alpha}$ and $B_k^{\beta}$ meet for $\alpha \leq
\lambda(r)$, we have
$$
B_k^{\beta}\subset {\bar B}_k^{\alpha} \quad \mbox{and}\quad {\bar
B}_k^{\beta}\subset {\bar{\bar B}}_k^{\alpha},
$$
\item[(iii)] for each $\alpha \leq \lambda(r)$, there no more than
$N(r)$ balls ever meet $B_k^{\alpha}$, and \item[(iv)] for any
$\alpha \le \lambda(r)$ and any $\beta$, either $B_k^{\alpha}$
meets $B_k^{\beta}$ for all $k\ge \bar{k}(r)$ or none for all
$k\ge \bar{k}(r)$.
\end{itemize}

Now we let $\hat E^{\alpha}, E^{\alpha}$, $\bar E^{\alpha}$, and
${\bar{\bar E}}^{\alpha}$ be the balls of radii
$4\tilde{\tilde\rho}^{\alpha}, 5\tilde{\tilde\rho}^{\alpha},
45\tilde\rho^{\alpha}$, and $205\rho^{\alpha}$ around the origin
in Euclidean space $\mathbb R^n$. At each point $x_k^{\alpha}\in
B(p_k, s_k)$ we define coordinate charts
$H_k^{\alpha}:E^{\alpha}\to B^{\alpha}_k$ as the composition of a
linear isometry of $\mathbb R^n$ to the tangent space
$T_{x_k^{\alpha}}M_k$ with the exponential map
$\exp_{x_k^{\alpha}}$ at $x_k^{\alpha}$. We also get maps $\bar
H_k^{\alpha}:\bar E^{\alpha}\to \bar B^{\alpha}_k$ and ${\bar{\bar
H}}_k^{\alpha}:{\bar{\bar E}}^{\alpha}\to {\bar{\bar
B}}^{\alpha}_k$ in the same way. Note that (4.1.3) implies that
these maps are all well defined. We denote by $g_k^{\alpha}$ (and
$\bar g_k^{\alpha}$ and ${\bar{\bar g}}_k^{\alpha}$) the
pull-backs of the metric $g_k$ by $H_k^{\alpha}$ (and $\bar
H_k^{\alpha}$ and ${\bar{\bar H}}_k^{\alpha}$). We also consider
the coordinate transition functions $J_k^{\alpha\beta}:
E^{\beta}\to \bar E^{\alpha}$ and $\bar J_k^{\alpha\beta}: \bar
E^{\beta}\to {\bar{\bar E}}^{\alpha}$ defined by
$$
J_k^{\alpha\beta}=(\bar H_k^{\alpha})^{-1}H_k^{\beta} \quad
{\text{and}} \quad \bar J_k^{\alpha\beta}=({\bar{\bar
H}}_k^{\alpha})^{-1}\bar H_k^{\beta}.
$$
Clearly $\bar J_k^{\alpha\beta} J_k^{\beta\alpha}=I$. Moreover
$J_k^{\alpha\beta}$ is an isometry from $g_k^{\beta}$ to $\bar
g_k^{\alpha}$ and $\bar J_k^{\alpha\beta}$ from $\bar g_k^{\beta}$
to ${\bar{\bar g}}_k^{\alpha}$.

Now for each fixed $\alpha$, the metrics $g_k^{\alpha}$ are in
geodesic coordinates and have their curvatures and their covariant
derivatives uniformly bounded.

\vskip 0.2cm {\bf Claim 1.} By passing to another subsequence we
can guarantee that for each $\alpha$ (and indeed all $\alpha$ by
diagonalization) the metrics $g_k^{\alpha}$ (or $\bar
g_k^{\alpha}$ or ${\bar{\bar g}}_k^{\alpha}$) converge uniformly
with their derivatives to a smooth metric $g^{\alpha}$ (or $\bar
g^{\alpha}$ or ${\bar{\bar g}}^{\alpha}$) on $E^{\alpha}$ (or
$\bar E^{\alpha}$ or ${\bar{\bar E}}^{\alpha}$) which is also in
geodesic coordinates. \vskip 0.2cm

Look now at any pair $\alpha, \beta$ for which the balls
$B_k^{\alpha}$ and $B_k^{\beta}$ always meet for large $k$, and
thus the maps $J_k^{\alpha\beta}$ (and $\bar J_k^{\alpha\beta}$
and $J_k^{\beta\alpha}$ and $\bar J_k^{\beta \alpha}$) are always
defined for large $k$.

\vskip 0.2cm {\bf Claim 2.} The isometries $J_k^{\alpha\beta}$
(and $\bar J_k^{\alpha\beta}$ and $J_k^{\beta\alpha}$ and $\bar
J_k^{\beta\alpha}$) always have a convergent subsequence. \vskip
0.2cm

So by passing to another subsequence we may assume
$J_k^{\alpha\beta}\to J^{\alpha\beta}$ (and $\bar
J_k^{\alpha\beta}\to \bar J^{\alpha\beta}$ and
$J_k^{\beta\alpha}\to J^{\beta\alpha}$ and $\bar
J_k^{\beta\alpha}\to \bar J^{\beta\alpha}$). The limit maps
$J^{\alpha\beta}: E^{\beta}\to \bar E^{\alpha}$ and $\bar
J^{\alpha\beta}: \bar E^{\beta}\to {\bar{\bar E}}^{\alpha}$ are
isometries in the limit metrics $g^{\beta}$ and $g^{\alpha}$.
Moreover
$$
J^{\alpha\beta}\bar J^{\beta\alpha}=I.
$$
We are now done picking subsequences, except we still owe the
reader the proofs of Claim 1 and Claim 2.

\medskip
{\it Step} 2: Finding local diffeomorphisms which are approximate
isometries.

\smallskip
Take the subsequence $(B(p_k, s_k), g_k, p_k)$ chosen in Step 1
above. We claim that for every $r < s_{\infty}$ and every
$(\epsilon_1,\epsilon_2,\ldots,\epsilon_p)$, and for all $k$ and
$l$ sufficiently large in comparison, we can find a diffeomorphism
$F_{kl}$ of a neighborhood of the ball $B(p_k, r)\subset B(p_k,
s_k)$ into an open set in $B(p_l, s_l)$ which is an
$(\epsilon_1,\epsilon_2,\ldots,\epsilon_p)$ approximate isometry
in the sense that
$$
|^{t}\nabla F_{kl} \nabla F_{kl}-I|<\epsilon_1
$$
and
$$
|\nabla^2F_{kl}|<\epsilon_2,\ldots, |\nabla^pF_{kl}|<\epsilon_p
$$
where $\nabla^pF_{kl}$ is the $p^{th}$ covariant derivative of
$F_{kl}$.

The idea (following Peters \cite{Pe} or Greene and Wu \cite{Gre})
of proving the claim is to define the map $F^\alpha_{k
\ell}=H^{\alpha}_l\circ(H^{\alpha}_k)^{-1}$ (or $\bar{F}^\alpha_{k
\ell}=\bar{H}^{\alpha}_l\circ(\bar{H}^{\alpha}_k)^{-1}$, resp.)
from $B^\alpha_k$ to $B^\alpha_\ell$ (or $\bar{B}^\alpha_k$ to
$\bar{B}^\alpha_\ell$, resp.) for $k$ and $\ell$ large compared to
$\alpha$ so as to be the identity map on ${E}^\alpha$ (or
$\bar{E}^\alpha$, resp.) in the coordinate charts $H^\alpha_k$ and
$H^\alpha_\ell$ (or $\bar{H}^\alpha_k$ and $\bar{H}^\alpha_\ell$,
resp.), and then to define $F_{k \ell}$ on a neighborhood of
$B(p_k,r)$ for $k, \ell \geq \bar{k}(r)$ be averaging the maps
$\bar{F}^\beta_{k \ell}$ for $\beta \leq \
\lambda(r+\frac{1}{5}(s_{\infty}-r))$. To describe the averaging
process on $B^\alpha_k$ with $\alpha \leq \lambda(r)$ we only need
to consider those $B^\beta_k$ which meet $B^\alpha_k$; there are
never more than $N(r)$ of them and each $\beta \leq
\lambda(r+\frac{1}{5}(s_{\infty}-r))$, and they are the same for
$k$ and $\ell$ when $k, \ell \geq \bar{k}(r)$. The averaging
process is defined by taking $F_{k \ell}(x)$ to be the center of
mass of the $\bar{F}^\beta_{k \ell}(x)$ for $x \in B^\alpha_k$
averaging over those $\beta$ where $B^\beta_k$ meets $B^\alpha_k$
using weights $\mu^\beta_k(x)$ defined by a partition of unity.
The center of mass of the points $y^\beta=F^\beta_{k \ell}(x)$
with weights $\mu^\beta$ is defined to be the point $y$ such that
$$
\exp_yV^\beta = y^\beta \quad \text{and} \quad \sum \mu^\beta
V^\beta=0.
$$
When the points $y^\beta$ are all close and the weights
$\mu^\beta$ satisfy $0 \leq \mu^\beta \leq 1$ then there will be a
unique solution $y$ close to $y^\beta$ which depends smoothly on
the $y^\beta$ and the $\mu^\beta$ (see \cite{Gre} for the
details). The point $y$ is found by the inverse function theorem,
which also provides bounds on all the derivatives of $y$ as a
function of the $y^\beta$ and the $\mu^\beta$.

Since $B^\alpha_k \subseteq \bar{B}^\beta_k$ and
$\bar{B}^\beta_\ell \subseteq \bar{\bar{B}}^{\alpha}_\ell$, the
map $\bar{F}^\beta_{k
\ell}=\bar{H}^{\beta}_l\circ(\bar{H}^{\beta}_k)^{-1}$ can be
represented in local coordinates by the map
$$
P^{\alpha \beta}_{k \ell}: E^\alpha \rightarrow
\bar{\bar{E}}^\alpha
$$
defined by
$$
P^{\alpha \beta}_{k \ell} = \bar{J}^{\alpha \beta}_\ell \circ
J^{\beta \alpha}_k.
$$
Since $J^{\beta \alpha}_k \rightarrow J^{\beta \alpha}$ as $k
\rightarrow \infty$ and $\bar{J}^{\alpha \beta}_\ell \rightarrow
\bar{J}^{\alpha \beta}$ as $\ell \rightarrow \infty$ and
$\bar{J}^{\alpha \beta} \circ J^{\beta \alpha} = I$, we see that
the maps $P^{\alpha \beta}_{k \ell} \rightarrow I$ as $k, \ell
\rightarrow \infty$ for each choice of $\alpha$ and $\beta$. The
weights $\mu^\beta_k$ are defined in the following way. We pick
for each $\beta$ a smooth function $\psi^\beta$ which equals 1 on
$\hat E^\beta$ and equals 0 outside $E^\beta$. We then transfer
$\psi^\beta$  to a function $\psi^\beta_k$ on $M_k$ by the
coordinate map $\bar{\bar{H}}^\beta_k$ (i.e. $\psi^\beta_k =
\psi^\beta\circ(\bar{\bar{H}}^\beta_k)^{-1}$). Then let
$$
\mu^\beta_k = \psi^\beta_k \Big/ \sum_\gamma \psi^\gamma_k
$$
as usual. In the coordinate chart $E^\alpha$ the function
$\psi^\beta_k$ looks like the composition of $J^{\beta \alpha}_k$
with $\psi^\beta$. Call this function
$$
\psi^{\alpha \beta}_k = \psi^\beta \circ J^{\beta \alpha}_k.
$$
Then as $k \rightarrow \infty, \psi^{\alpha \beta}_k \rightarrow
\psi^{\alpha \beta}$ where
$$
\psi^{\alpha \beta} = \psi^\beta \circ J^{\beta \alpha}.
$$
In the coordinate chart $E^\alpha$ the function $\mu^\beta_k$
looks like
$$
\mu_k^{\alpha \beta} = \psi_k^{\alpha \beta} \Big/ \sum_\gamma
\psi_k^{\alpha \gamma}
$$
and $\mu^{\alpha \beta}_k \rightarrow \mu^{\alpha \beta}$ as $k
\rightarrow \infty$ where
$$
\mu^{\alpha \beta} = \psi^{\alpha \beta} \Big/ \sum_\gamma
\psi^{\alpha \gamma}.
$$
Since the sets $\hat{B}^\alpha_k$ cover $B(p_k,r)$, it follows
that $\sum_\gamma \psi^\gamma_k \geq 1$ on this set and by
combining with (4.1.5) and (4.1.7) there is no problem bounding
all these functions and their derivatives. There is a small
problem in that we want to guarantee that the averaged map still
takes $p_k$ to $p_\ell$. This is true at least for the map $F^0_{k
\ell}$. Therefore it will suffice to guarantee that
$\mu^\alpha_k=0$ in a neighborhood of $p_k$ if $\alpha \neq 0$.
This happens if the same is true for $\psi^\alpha_k$. If not, we
can always replace $\psi^\alpha_k$ by $\tilde{\psi}^\alpha_k=(1-
\psi^0_k)\psi^\alpha_k$ which still leaves $\tilde{\psi}^\alpha_k
\geq \frac{1}{2} \psi^\alpha_k$ or $\psi^0_k \geq \frac{1}{2}$
everywhere, and this is sufficient to make $\sum_\gamma
\tilde{\psi}^\gamma_k \geq \frac{1}{2}$ everywhere.

Now in the local coordinate $E^\alpha$ we are averaging maps
$P^{\alpha \beta}_{k \ell}$ which converge to the identity with
respect to weights $\mu^{\alpha \beta}_k$ which converge. It
follows that the averaged map converges to the identity in these
coordinates. Thus $F_{k \ell}$ can be made to be an $(\epsilon_1,
\epsilon_2, \ldots, \epsilon_p)$ approximate isometry on
$B(p_k,r)$ when $k$ and $\ell$ are suitably large. At least the
estimates
$$
|^t\nabla F_{k \ell} \cdot \nabla F_{k\ell}-I| < \epsilon_1
$$
and $|\nabla^2F_{k \ell}| < \epsilon_2, \ldots, |\nabla^pF_{k
\ell}| < \epsilon_p$ on $B(p_k,r)$ follow from the local
coordinates. We still need to check that $F_{k \ell}$ is a
diffeomorphism on a neighborhood of $B(p_k,r)$.

This, however, follows quickly enough from the fact that we also
get a map $F_{\ell k}$ on a slightly larger ball $B(p_\ell,
r^\prime)$ which contains the image of $F_{k \ell}$ on $B(p_k,r)$
if we take $r^\prime=(1+ \epsilon_1)r$, and $F_{\ell k}$ also
satisfies the above estimates. Also $F_{k \ell}$ and $F_{\ell k}$
fix the markings, so the composition $F_{\ell k} \circ F_{k \ell}$
satisfies the same sort of estimates and fixes the origin $p_k$.

Since the maps $P^{\alpha \beta}_{k \ell}$  and $P^{\alpha
\beta}_{\ell k}$ converge to the identity as $k, \ell$ tend to
infinity, $F_{\ell k} \circ F_{k \ell}$ must be very close to the
identity on $B(p_k,r)$. It follows that $F_{k \ell}$ is
invertible.  This finishes the proof of the claim and the Step 2.

\medskip
{\it Step} 3: Constructing the limit geodesic ball $(B_{\infty},
g_{\infty}, p_{\infty})$.

\smallskip
We now know the geodesic balls $(B(p_k,s_k),g_k,p_k)$ are nearly
isometric for large $k$. We are now going to construct the limit
$B_{\infty}$.  For a sequence of positive numbers $r_j \nearrow
s_{\infty}$ with each $r_j < s_j$, we choose the numbers
$(\epsilon_1(r_j), \ldots, \epsilon_j(r_j))$ so small that when we
choose $k(r_j)$ large in comparison and find the maps
$F_{k(r_j),k(r_{j+1})}$ constructed above on neighborhoods of
$B(p_{k(r_j)},r_j)$, in $M_{k(r_j)}$ into $M_{k(r_{j+1})}$ the
image always lies in $B(p_{k(r_{j+1})}, r_{j+1})$ and the
composition of $F_{k(r_j),k(r_{j+1})}$ with
$F_{k(r_{j+1}),k(r_{j+2})}$ and $\cdots$ and
$F_{k(r_{s-1}),k(r_s)}$ for any $s> j$ is still an $(\eta_1(r_j),
\ldots, \eta_j(r_j))$ isometry for any choice of $\eta_i(r_j)$,
say $\eta_i(r_j)=1/j$ for $1 \leq i \leq j$.  Now we simplify the
notation by writing $M_j$ in place of $M_{k(r_j)}$ and $F_j$ in
place of $F_{k(r_j),k(r_{j+1})}$. Then
$$
F_j:B(p_j,r_j) \rightarrow B(p_{j+1}, r_{j+1})
$$
is a diffeomorphism map from $B(p_j,r_j)$ into $B(p_{j+1},
r_{j+1})$, and the composition
$$
F_{s-1} \circ \cdots \circ F_j: B(p_j,j) \rightarrow B(p_s,s)
$$
is always an $(\eta_1(r_j), \ldots, \eta_j(r_j))$ approximate
isometry.

We now construct the limit $B_{\infty}$ as a topological space by
identifying the balls $B(p_j,r_j)$ with each other using the
homeomorphisms $F_j$. Given any two points $x$ and $y$ in
$B_{\infty}$, we have $x \in B(p_j,r_j)$ and $y \in B(p_s,r_s)$
for some $j$ and $s$. If $j \leq s$ then $x \in B(p_s,r_s)$ also,
by identification. A set in $B_{\infty}$ is open if and only if it
intersects each $B(p_j,r_j)$ in an open set. Then choosing
disjoint neighborhoods of $x$ and $y$ in $B(p_s,r_s)$ gives
disjoint neighborhoods of $x$ and $y$ in $B_{\infty}$. Thus
$B_{\infty}$ is a Hausdorff space.

Any smooth chart on $B(p_j,r_j)$ also gives a smooth chart on
$B(p_s,r_s)$ for all $s > j$. The union of all such charts gives a
smooth atlas on $B_{\infty}$. It is fairly easy to see the metrics
$g_j$ on $B(p_j,r_j)$, converge to a smooth metric $g_{\infty}$ on
$B_{\infty}$ uniformly together with all derivatives on compact
sets. For since the $F_{s-1} \circ \cdots \circ F_j$ are very good
approximate isometries, the $g_j$ are very close to each other,
and hence form a Cauchy sequence (together with their derivatives,
in the sense that the covariant derivatives of $g_j$ with respect
to $g_s$ are very small when $j$ and $s$ are both large). One
checks in the usual way that such a Cauchy sequence converges.

The origins $p_j$ are identified with each other, and hence with
an origin $p_{\infty}$ in $B_{\infty}$. Now it is the inverses of
the maps identifying $B(p_j,r_j)$ with open subsets of
$B_{\infty}$ that provide the diffeomorphisms of (relatively
compact) open sets in $B_{\infty}$ into the geodesic balls
$B(p_j,s_j) \subset M_j$ such that the pull-backs of the metrics
$g_j$ converge to $g_{\infty}$. This completes the proof of Step
3.

\medskip
Now it remains to prove both Claim 1 and Claim 2 in Step 1.

\medskip
{\bf\em Proof of Claim} {\bf 1.} \ It suffices to show the
following general result:

\medskip
{\it There exists a constant $c>0$ depending only on the
dimension, and constants $C_q$ depending only on the dimension and
$q$ and bounds $B_j$ on the curvature and its derivatives for $j
\leq q$ where $|D^jRm| \leq B_j$, so that for any metric $g_{k
\ell}$ in geodesic coordinates in the ball $|x| \leq r \leq c/
\sqrt{B_0}$, we have
$$
\frac{1}{2} I_{k \ell} \leq g_{k \ell} \leq 2I_{k \ell}
$$
and
$$
\Big| \frac{\partial}{\partial x^{j_1}} \cdots
\frac{\partial}{\partial x^{j_q}} g_{k \ell} \Big| \leq C_q,
$$
where $I_{k\ell}$ is the Euclidean metric.}

\medskip
Suppose we are given a metric $g_{ij}(x)dx^idx^j$ in geodesic
coordinates in the ball $|x| \leq r \leq c/ \sqrt{B_0}$ as in
Claim 1. Then by definition every line through the origin is a
geodesic (parametrized proportional to arc length) and
$g_{ij}=I_{ij}$ at the origin. Also, the Gauss Lemma says that the
metric $g_{ij}$ is in geodesic coordinates if and only if
$g_{ij}x^i=I_{ij}x^i$. Note in particular that in geodesic
coordinates
$$
|x|^2=g_{ij}x^ix^j=I_{ij}x^ix^j
$$
is unambiguously defined. Also, in geodesic coordinates we have
$\Gamma^k_{ij}(0)=0$, and all the first derivatives for $g_{jk}$
vanish at the origin.

Introduce the symmetric tensor
$$
A_{ij} = \frac{1}{2} x^k \frac{\partial}{\partial x^k} g_{ij}.
$$
Since we have $g_{jk}x^k=I_{jk}x^k$, we get
$$
x^k \frac{\partial}{\partial x^i}g_{jk} = I_{ij}-g_{ij}=x^k
\frac{\partial}{\partial x^j}g_{ik}
$$
and hence from the formula for $\Gamma^i_{jk}$
$$
x^j \Gamma^i_{jk}=g^{i \ell}A_{k \ell}.
$$
Hence $A_{k \ell}x^k=0$. Let $D_i$ be the covariant derivative
with respect to the metric $g_{ij}$. Then
$$
D_ix^k=I^k_i+ \Gamma^k_{ij}x^j=I^k_i+g^{k \ell}A_{i \ell}.
$$
Introduce the potential function
$$
P=|x|^2/2 = \frac{1}{2}g_{ij}x^ix^j.
$$
We can use the formulas above to compute
$$
D_iP=g_{ij}x^j.
$$
Also we get
$$
D_iD_jP=g_{ij}+A_{ij}.
$$
The defining equation for $P$ gives
$$
g^{ij}D_iPD_jP=2P.
$$
If we take the covariant derivative of this equation we get
$$
g^{k \ell}D_jD_kPD_{\ell}P=D_jP
$$
which is equivalent to $A_{jk}x^k=0$. But if we take the covariant
derivative again we get
$$
g^{k \ell}D_iD_jD_kPD_\ell P + g^{k \ell}D_jD_kPD_iD_\ell P = D_i
D_j P.
$$
Now switching derivatives
$$
D_iD_jD_kP=D_iD_kD_jP=D_kD_iD_jP+R_{ikj \ell}g^{\ell m}D_mP
$$
and if we use this and $D_iD_jP=g_{ij}+A_{ij}$ and $g^{k
\ell}D_{\ell}P=x^k$ we find that
$$
x^kD_kA_{ij}+A_{ij}+g^{k \ell} A_{ik}A_{j \ell} + R_{ikj \ell}x^k
x^\ell=0.
$$
{}From our assumed curvature bounds we can take $|R_{ijk \ell}|
\leq B_0$. Then we get the following estimate:
$$
|x^kD_kA_{ij}+A_{ij}| \leq C|A_{ij}|^2+CB_0 r^2
$$
on the ball $|x| \leq r$ for some constant $C$ depending only on
the dimension.

We now show how to use the maximum principle on such equations.
First of all, by a maximum principle argument, it is easy to show
that if $f$ is a function on a ball $|x| \leq r$ and $\lambda>0$
is a constant, then
$$
\lambda \operatorname{sup}|f| \leq \operatorname{sup} \Big| x^k
\frac{\partial f}{\partial x^k} + \lambda f \Big|.
$$
For any tensor $T=\{T_{i \cdots j}\}$ and any constant $\lambda
>0$, setting $f=|T|^2$ in the above inequality, we have \be
\lambda \sup|T| \leq \sup|x^kD_kT+
\lambda T|. 
\ee Applying this to the tensor $A_{ij}$ we get
\[
\sup_{|x| \leq r} |A_{ij}| \leq C \sup_{|x| \leq r}
|A_{ij}|^2+CB_0r^2
\]
for some constant depending only on the dimension.

It is fairly elementary to see that there exist constants $c>0$
and $C_0 < \infty$ such that if the metric $g_{ij}$ is in geodesic
coordinates with $|R_{ijk \ell}| \leq B_0$ in the ball of radius
$r \leq c /\sqrt{B_0}$ then
$$
|A_{ij}| \leq C_0 B_0r^2.
$$
Indeed, since the derivatives of $g_{ij}$ vanish at the origin, so
does $A_{ij}$. Hence the estimate holds near the origin. But the
inequality
$$
\sup_{|x| \leq r} |A_{ij}| \leq C \sup_{|x| \leq r}
|A_{ij}|^2+CB_0r^2
$$
says that $|A_{ij}|$ avoids an interval when $c$ is chosen small.
In fact the inequality
$$
X \leq CX^2+D
$$
is equivalent to
$$
|2CX-1| \geq \sqrt{1-4CD}
$$
which makes $X$ avoid an interval if $4CD<1$. (Hence in our case
we need to choose $c$ with $4C^2c^2<1.)$ Then if $X$ is on the
side containing $0$ we get
$$
X \leq \frac{1- \sqrt{1 - 4CD}}{2C} \leq 2D.
$$
This gives $|A_{ij}| \leq C_0B_0r^2$ with $C_0=2C$.

\smallskip
We can also derive bounds on all the covariant derivatives of $P$
in terms of bounds on the covariant derivatives of the curvature.
To simplify the notation, we let
$$
D^qP= \{D_{j_1}D_{j_2} \cdots D_{j_q}P \}
$$
denote the $q^{\operatorname{th}}$ covariant derivative, and in
estimating $D^qP$ we will lump all the lower order terms into a
general slush term $\Phi^q$ which will be a polynomial in
$D^1P,D^2P, \ldots, D^{q-1}P$ and $Rm, D^1Rm, \ldots, D^{q-2}Rm$.
We already have estimates on a ball of radius $r$
$$
P \leq r^2/2
$$
$$
|D^1P|\leq r
$$
$$
|A_{ij}| \leq C_0 B_0r^2
$$
and since $D_iD_jP=g_{ij}+A_{ij}$ and $r \leq c /\sqrt{B_0}$ if we
choose $c$ small we can make
$$
|A_{ij}|\leq 1/2,
$$
and we get
$$
|D^2P| \leq C_2
$$
for some constant $C_2$ depending only on the dimension.

\smallskip
Start with the equation $g^{ij}D_iPD_jP=2P$ and apply repeated
covariant derivatives. Observe that we get an equation which
starts out
$$
g^{ij}D_iPD^qD_jP+ \cdots =0
$$
where the omitted terms only contain derivatives $D^qP$ and lower.
If we switch two derivatives in a term $D^{q+1}P$ or lower, we get
a term which is a product of a covariant derivative of $Rm$ of
order at most $q-2$ (since the two closest to $P$ commute) and a
covariant derivative of $P$ of order at most $q-1$; such a term
can be lumped in with the slush term $\Phi^q$. Therefore up to
terms in $\Phi^q$ we can regard the derivatives as commuting. Then
paying attention to the derivatives in $D^1P$ we get an equation

\smallskip\noindent
$g^{ij}D_iPD_jD_{k_1} \cdots D_{k_q}P+g^{ij}D_iD_{k_1}PD_jD_{k_2}
\cdots D_{k_q}P$

\smallskip\noindent
$+g^{ij}D_iD_{k_2}PD_jD_{k_1}D_{k_3} \cdots D_{k_q}P + \cdots +
g^{ij}D_iD_{k_q}PD_jD_{k_1} \cdots D_{k_{q-1}}P$

\smallskip\noindent
$ =D_{k_1} \cdots D_{k_q}P+ \Phi^q$.

\medskip
Recalling that $D_iD_jP=g_{ij}+A_{ij}$ we can rewrite this as
\begin{align*}
\Phi^q &  = g^{ij}D_iPD_jD_{k_1} \cdots D_{k_q}+(q-1)D_{k_1}
\cdots
D_{k_q}P \\
&\quad  + g^{ij}A_{ik_1}D_jD_{k_2} \cdots D_{k_q}P+ \cdots +
g^{ij}A_{ik_q}D_jD_{k_1} \cdots D_{k_{q-1}}P.
\end{align*}
Estimating the product of tensors in the usual way gives
$$
|x^iD_iD^qP+(q-1)D^qP| \leq q|A||D^qP|+|\Phi^q|.
$$
Applying the inequality $\lambda \sup|T| \leq \sup|x^kD_kT+
\lambda T|$ with $T=D^qP$ gives
$$
(q-1) \sup|D^qP| \leq\sup(q|A||D^qP|+| \Phi^q|).
$$
Now we can make $|A| \leq 1/2$ by making $r \leq c/ \sqrt{B_0}$
with $c$ small; it is important here that $c$ is independent of
$q!$ Then we get
$$
(q-2) \operatorname{sup}|D^qP| \leq 2 \operatorname{sup} | \Phi^q|
$$
which is a good estimate for $q \geq 3$. The term $\Phi^q$ is
estimated inductively from the terms $D^{q-1}P$ and $D^{q-2}Rm$
and lower. This proves that there exist constants $C_q$ for $q
\geq 3$ depending only on $q$ and the dimension and on $|D^jRm|$
for $j \leq q-2$ such that
$$
|D^qP| \leq C_q
$$
on the ball $r \leq c/ \sqrt{B_0}.$

Now we turn our attention to estimating the Euclidean metric
$I_{jk}$ and its covariant derivatives with respect to $g_{jk}$.
We will need the following elementary fact: suppose that $f$ is a
function on a ball $|x| \leq r$ with $f(0)=0$ and
$$
\Big| x^i \frac{\partial f}{\partial x^i} \Big| \leq C |x|^2
$$
for some constant $C$. Then \be
|f| \leq C|x|^2  
\ee for the same constant $C$. As a consequence, if $T = \{T_{j
\cdots k }\}$ is a tensor which vanishes at the origin and if
$$
|x^iD_iT| \leq C|x|^2
$$
on a ball $|x| \leq r$ then $|T| \leq  C|x|^2$ with the same
constant $C$. (Simply apply the inequality (4.1.9) to the function
$f=|T|$. In case this is not smooth, we can use
$f=\sqrt{|T|^2+\epsilon^2}-\epsilon$ and then let $\epsilon\to
0$.)

Our application will be to the tensor $I_{jk}$ which gives the
Euclidean metric as a tensor in geodesic coordinates. We have
$$
D_iI_{jk} = - \Gamma^p_{ij}I_{pk} - \Gamma^p_{ik}I_{pj}
$$
and since
$$
x^i \Gamma^p_{ij} = g^{pq}A_{jq}
$$
we get the equation
$$
x^iD_iI_{jk} = -g^{pq}A_{jp}I_{kq}-g^{pq}A_{kp}I_{jq}
$$.

\noindent We already have $|A_{jk}| \leq C_0B_0|x|^2$ for $|x|
\leq r \leq c/ \sqrt{B_0}$. The tensor $I_{jk}$ doesn't vanish at
the origin, but the tensor
$$
h_{jk}=I_{jk}-g_{jk}
$$
does. We can then use
$$
x^iD_ih_{jk}=-g^{pq}A_{jp}h_{kq}-g^{pq}A_{kq}h_{jq}-2A_{jk}.
$$
Suppose $M(s) = \sup_{|x| \leq s}|h_{jk}|$. Then
$$
|x^iD_ih_{jk}| \leq 2[1+M(s)]C_0B_0|x|^2
$$
and we get
$$
|h_{jk}| \leq 2 [1+M(s)]C_0B_0|x|^2
$$
on $|x| \leq s$. This makes
$$
M(s) \leq 2[1+M(s)]C_0B_0s^2.
$$
Then for $s \leq r \leq c/ \sqrt{B_0}$ with $c$ small compared to
$C_0$ we get $2C_0B_0s^2 \leq 1/2$ and $M(s) \leq 4C_0B_0s^2$.
Thus
$$
|I_{jk}-g_{jk}|=|h_{jk}| \leq 4C_0B_0|x|^2
$$
for $|x| \leq r\leq c/ \sqrt{B_0}$, and hence for $c$ small enough
$$
\frac{1}{2} g_{jk} \leq I_{jk} \leq 2g_{jk}.
$$
Thus the metrics are comparable. Note that this estimate only
needs $r$ small compared to $B_0$ and does not need any bounds on
the derivatives of the curvature.

Now to obtain bounds on the covariant derivative of the Eucliden
metric $I_{k \ell}$ with respect to the Riemannian metric $g_{k
\ell}$ we want to start with the equation
$$
x^iD_iI_{k \ell}+g^{mn}A_{km}I_{\ell n}+g^{mn}A_{\ell m}I_{kn}=0
$$
and apply $q$ covariant derivatives $D_{j_1} \cdots D_{j_q}$. Each
time we do this we must interchange $D_j$ and $x^iD_i$, and since
this produces a term which helps we should look at it closely. If
we write $R_{ji}=[D_j,D_i]$ for the commutator, this operator on
tensors involves the curvature but no derivatives. Since
$$
D_jx^i=I^i_j+g^{im}A_{jm}
$$
we can compute
$$
[D_j, x^iD_i]=D_j+g^{im}A_{jm}D_i+x^iR_{ji}
$$
and the term $D_j$ in the commutator helps, while $A_{jm}$ can be
kept small and $R_{ji}$ is zero order. It follows that we get an
equation of the form
\begin{align*}
0 &  = x^iD_iD_{j_1} \cdots D_{j_q}I_{k \ell} + qD_{j_1} \cdots
D_{j_q}I_{k \ell} \\
&\quad  + \sum^q_{h=1} g^{im}A_{j_hm}D_{j_1} \cdots
D_{j_{h-1}}D_iD_{j_{h+1}} \cdots D_{j_q}I_{k \ell} \\
&\quad  + g^{mn}A_{km}D_{j_1} \cdots D_{j_q}I_{\ell n} + g^{mn}
A_{\ell m} D_{j_1} \cdots D_{j_q}I_{kn}+ \Psi^q,
\end{align*}
where the slush term $\Psi^q$ is a polynomial in derivatives of
$I_{k \ell}$ of degree no more than $q-1$ and derivatives of $P$
of degree no more than $q+2$ (remember $x^i = g^{ij}D_jP$ and
$A_{ij} = D_iD_jP-g_{ij})$ and derivatives of the curvature $Rm$
of degree no more than $q-1$. We now estimate
$$
D^qI_{k \ell} = \{D_{j_1} \cdots D_{j_q}I_{k \ell}\}
$$
by induction on $q$ using (4.1.8) with $\lambda = q$. Noticing a
total of $q+2$ terms contracting $A_{ij}$ with a derivative of
$I_{k \ell}$ of degree $q$, we get the estimate
$$
q \sup |D^qI_{k \ell}| \leq (q+2) \operatorname{sup}|A|
\operatorname{sup} |D^q I_{k \ell} | + \operatorname{sup} |
\Psi^q|.
$$

\noindent and everything works. This proves that there exists a
constant $c>0$ depending only on the dimension, and constants
$C_q$ depending only on the dimension and $q$ and bounds $B_j$ on
the curvature and its derivatives for $j \leq q$ where $|D^jRm|
\leq B_j$, so that for any metric $g_{k \ell}$ in geodesic
coordinates in the ball $|x| \leq r \leq c/ \sqrt{B_0}$ the
Euclidean metric $I_{k \ell}$ satisfies
$$
\frac{1}{2} g_{k \ell}\leq I_{k \ell} \leq 2g_{k \ell}
$$
and the covariant derivatives of $I_{k \ell}$ with respect to
$g_{k \ell}$ satisfy
$$
|D_{j_1} \cdots D_{j_q}I_{k \ell}| \leq C_q.
$$

The difference between a covariant derivative and an ordinary
derivative is given by the connection
$$
-\Gamma^p_{ij}I_{pk} - \Gamma^p_{ik} I_{pj}
$$
to get
$$
\Gamma^k_{ij} = \frac{1}{2} I^{k \ell}(D_\ell I_{ij} - D_iI_{j
\ell} - D_j I_{i \ell}).
$$
This gives us bounds on $\Gamma^k_{ij}$. We then obtain bounds on
the first derivatives of $g_{ij}$ from
$$
\frac{\partial}{\partial x^i} g_{jk} = g_{k \ell} \Gamma^\ell_{ij}
+ g_{j \ell} \Gamma^\ell_{ik}.
$$
Always proceeding inductively on the order of the derivative, we
now get bounds on covariant derivatives of $\Gamma^k_{ij}$ from
the covariant derivatives of $I_{pk}$ and bounds of the ordinary
derivatives of $\Gamma^k_{ij}$ by relating the to the covariant
derivatives using the $\Gamma^k_{ij}$, and bounds on the ordinary
derivatives of the $g_{jk}$ from bounds on the ordinary
derivatives of the $\Gamma^\ell_{ij}$. Consequently, we have
estimates
$$
\frac{1}{2} I_{k \ell} \leq g_{k \ell} \leq 2I_{k \ell}
$$
and
$$
\Big| \frac{\partial}{\partial x^{j_1}} \cdots
\frac{\partial}{\partial x^{j_q}} g_{k \ell} \Big| \leq
\tilde{C}_q
$$
for similar constants $\tilde{C}_q$.

Therefore we have finished the proof of Claim 1.

\medskip
{\bf\em Proof of Claim} {\bf 2:} \  We need to show how to
estimate the derivatives of an isometry. We will prove that if
$y=F(x)$ is an isometry from a ball in Euclidean space with a
metric $g_{ij}dx^idx^j$ to a ball in Euclidean space with a metric
$h_{kl}dy^kdy^l$. Then we can bound all of the derivatives of $y$
with respect to $x$ in terms of bounds on $g_{ij}$ and its
derivatives with respect to $x$ and bound on $h_{kl}$ and its
derivatives with respect to $y$. This would imply Claim 2.

Since $y=F(x)$ is an isometry we have the equation
$$
h_{pq} \frac{\partial y^p}{\partial x^j} \frac{\partial
y^q}{\partial x^k} = g_{jk}.
$$

\noindent Using bounds $g_{jk} \leq CI_{jk}$ and $h_{pq} \geq
cI_{pq}$ comparing to the Euclidean metric, we easily get
estimates
$$
\Big| \frac{\partial y^p}{\partial x^j} \Big| \leq C.
$$

\noindent Now if we differentiate the equation with respect to
$x^i$ we get
$$
h_{pq} \frac{\partial^2y^p}{\partial x^i \partial x^j}
\frac{\partial y^q}{\partial x^k} + h_{pq} \frac{\partial
y^p}{\partial x^j} \frac{\partial^2y^q}{\partial x^i \partial x^k}
= \frac{\partial g_{jk}}{\partial x^i} - \frac{\partial
h_{pq}}{\partial y^r} \frac{\partial y^r}{\partial x^i}
\frac{\partial y^p}{\partial x^j} \frac{\partial y^q}{\partial
x^k}.
$$

\noindent Now let
$$
T_{ijk} = h_{pq} \frac{\partial y^p}{\partial x^i}
\frac{\partial^2 y^q}{\partial x^j \partial x^k}
$$

\noindent and let
$$
U_{ijk} = \frac{\partial g_{jk}}{\partial x^i} - \frac{\partial
h_{pq}}{\partial y^r} \frac{\partial y^r}{\partial x^i}
\frac{\partial y^p}{\partial x^j} \frac{\partial y^q}{\partial
x^k}.
$$

\noindent Then the above equation says
$$
T_{kij}+T_{jik}=U_{ijk}.
$$

\noindent Using the obvious symmetries $T_{ijk}=T_{ikj}$ and
$U_{ijk} = U_{ikj}$ we can solve this in the usual way to obtain
$$
T_{ijk} = \frac{1}{2}(U_{jik}+U_{kij} - U_{ijk}).
$$

\noindent We can recover the second derivatives of $y$ with
respect to $x$ from the formula
$$
\frac{\partial^2 y^p}{\partial x^i \partial x^j} = g^{k
\ell}T_{kij} \frac{\partial y^p}{\partial x^\ell}.
$$

\noindent Combining these gives an explicit formula giving
$\partial^2y^p/\partial x^i \partial x^j$ as a function of
$g^{ij}, h_{pq}, \partial g_{jk}/\partial x^i, \partial
h_{pq}/\partial y^r$, and $ \partial y^p/\partial y^i$. This gives
bounds
$$
\Big| \frac{\partial^2 y^p}{\partial y^i \partial y^j}\Big| \leq C
$$

\noindent and bounds on all higher derivatives follow by
differentiating the formula and using induction. This completes
the proof of Claim 2 and hence the proof of Theorem 4.1.2.
\end{pf}

We now want to show how to use this convergence result on
solutions to the Ricci flow. Let us first state the definition for
the convergence of evolving manifolds.

\begin{definition}
Let $(M_k,g_k(t),p_k)$ be a sequence of evolving marked complete
Riemannian manifolds, with the evolving metrics $g_k(t)$ over a
fixed time interval $t\in(A,\Omega]$, $A < 0 \leq \Omega$, and
with the marked points $p_k\in M_k$. We say a sequence of evolving
marked $(B_0(p_k,s_k),g_k(t),p_k)$ over $t \in (A,\Omega]$, where
$B_0(p_k,s_k)$ are geodesic balls of $(M_k,g_k(0))$ centered at
$p_k$ with the radii $s_k \rightarrow s_{\infty}(\leq +\infty)$,
{\bf converges} in the $C_{\rm loc}^{\infty}$ topology {\bf to an
evolving marked} (maybe noncomplete) {\bf
manifold}\index{converges to an evolving marked manifold}
$(B_{\infty},g_{\infty}(t),p_{\infty})$ over $t\in(A,\Omega]$,
where, at the time $t=0$, $B_{\infty}$ is a geodesic open ball
centered at $p_{\infty}\in B_{\infty}$ with the radius
$s_{\infty}$, if we can find a sequence of exhausting open sets
$U_k$ in $B_{\infty}$ containing $p_{\infty}$ and a sequence of
diffeomorphisms $f_k$ of the sets $U_k$ in $B_{\infty}$ to open
sets $V_k$ in $B(p_k,s_k) \subset M_k$ mapping $p_{\infty}$ to
$p_k$ such that the pull-back metrics
$\tilde{g}_k(t)=(f_k)^*g_k(t)$ converge in $C^\infty$ topology to
$g_{\infty}(t)$ on every compact subset of
$B_{\infty}\times(A,\Omega]$.
\end{definition}

Now we fix a time interval $A<t \leq \Omega$ with $-\infty < A<0$
and $0 \leq \Omega <+\infty$. Consider a sequence of marked
evolving complete manifolds $(M_k,g_k(t),p_k),\ t\in(A,\Omega]$,
with each $g_k(t)$, $k=1, 2, \ldots,$ being a solution of the
Ricci flow
$$
\frac{\partial}{\partial t}g_k(t)=-2\Ric_k(t)
$$
on $B_0(p_k,s_k) \times (A,\Omega]$, where $Ric_k$ is the Ricci
curvature tensor of $g_k$, and  $B_0(p_k,s_k)$ is the geodesic
ball of $(M_k,g_k(0))$ centered at $p_k$ with the radii $s_k
\rightarrow s_{\infty}(\leq +\infty)$.

Assume that for each $r < s_{\infty}$ there are positive constants
$C(r)$ and $k(r)$ such that the curvatures of $g_k(t)$ satisfy the
bound
$$
|Rm(g_k)| \leq C(r)
$$
on $B_0(p_k,r) \times (A,\Omega]$ for all $k \geq k(r).$ We also
assume that $(M_k,g_k(t),p_k)$, $k=1, 2, \ldots,$ have a uniform
injectivity radius bound at the origins $p_k$ at $t=0$. By Shi's
derivatives estimate (Theorem 1.4.1), the above assumption of
uniform bound of the curvatures on the geodesic balls $B_0(p_k,r)$
($r < s_{\infty}$) implies the uniform bounds on all the
derivatives of the curvatures at $t=0$ on the geodesic balls
$B_0(p_k,r)$ ($r < s_{\infty}$). Then by Theorem 4.1.2 we can find
a subsequence of marked evolving manifolds, still denoted by
$(M_k,g_k(t),p_k)$ with $t\in(A,\Omega]$, so that the geodesic
balls $(B_0(p_k,s_k),g_k(0),p_k)$ converge in the $C^\infty_{loc}$
topology to a geodesic ball
$(B_{\infty}(p_{\infty},s_{\infty}),g_{\infty}(0),p_{\infty})$.
>From now on, we consider this subsequence of marked evolving
manifolds.  By Definition 4.1.1, we have a sequence of (relatively
compact) exhausting covering $\{U_k\}$ of
$B_{\infty}(p_{\infty},s_{\infty})$ containing $p_{\infty}$ and a
sequence of diffeomorphisms $f_k$ of the sets $U_k$ in
$B_{\infty}(p_{\infty},s_{\infty})$ to open sets $V_k$ in
$B_{0}(p_{k},s_{k})$ mapping $p_{\infty}$ to $p_k$ such that the
pull-back metrics at $t=0$
$$
\tilde{g}_k(0)=(f_k)^*g_k(0) \stackrel{C_{\rm
loc}^\infty}{\longrightarrow}g_{\infty}(0),\ \quad \text{ as }\;
k\rightarrow+\infty,\; \text{ on }\;
B_{\infty}(p_{\infty},s_{\infty}).
$$
However, the pull-back metrics $\tilde{g}_k(t)=(f_k)^*g_k(t)$ are
also defined at all times $A<t \leq \Omega$ (although
$g_{\infty}(t)$ is not yet). We also have uniform bounds on the
curvature of the pull-back metrics $\tilde{g}_k(t)$ and all their
derivatives, by Shi's derivative estimates (Theorem 1.4.1), on
every compact subset of $B_{\infty}(p_{\infty},s_{\infty}) \times
(A,\Omega]$. What we claim next is that we can find uniform bounds
on all the covariant derivatives of the $\tilde{g}_k$ taken with
respect to the fixed metric $g_{\infty}(0)$.

\begin{lemma} [{Hamilton \cite{Ha95}}]
Let $(M,g)$ be a Riemannian manifold, $K$ a compact subset of $M$,
and $\tilde{g}_k(t)$ a collection of solutions to Ricci flow
defined on neighborhoods of $K\times[\alpha,\beta]$ with
$[\alpha,\beta]$ containing $0$.  Suppose that for each $l\geq0$,
\begin{itemize}
\item[(a)] $C_0^{-1}g\leq\tilde{g}_k(0)\leq C_0g,\ \ on\ K,\ for\
all\ k,$ \item[(b)] $|\nabla^l\tilde{g}_k(0)|\leq C_l,\ \ on\ K,\
for\ all\ k,$ \item[(c)]
$|\tilde{\nabla}^l_kRm(\tilde{g}_k)|_k\leq C'_l,\ on\
K\times[\alpha,\beta],\ for\ all\ k,$
\end{itemize}
for some positive constants $C_l,\ C'_l,\ l=0,1,\ldots,$
independent of $k$, where $Rm(\tilde{g}_k)$ are the curvature
tensors of the metrics $\tilde{g}_k(t)$, $\tilde{\nabla}_k$ denote
covariant derivative with respect to $\tilde{g}_k(t)$, $|\cdot|_k$
are the length of a tensor with respect to $\tilde{g}_k(t)$, and
$|\cdot|$ is the length with respect to $g$. Then the metrics
$\tilde{g}_k(t)$ satisfy
$$
\tilde{C_0}^{-1}g\leq\tilde{g}_k(t) \leq \tilde{C_0}g,\ \ on\
K\times[\alpha,\beta]
$$
and
$$
|\nabla^l\tilde{g}_k|\leq\tilde{C_l},\ on\ K\times[\alpha,\beta],\
\ l=1,2,\ldots,
$$
for all $k$, where $\tilde{C_l},\ l=0,1,\ldots,$ are positive
constants independent of $k$.
\end{lemma}

\begin{pf}
First by using the equation
$$
\frac{\partial}{\partial t}\tilde{g}_k=-2\tilde{\Ric}_k
$$
and the assumption (c) we immediately get \be
\tilde{C_0}^{-1}g\leq\tilde{g}_k(t)
\leq \tilde{C_0}g, \; \mbox{ on }\; K\times[\alpha,\beta] 
\ee for some positive constant $\tilde{C_0}$ independent of $k$.

Next we want to bound $\nabla\tilde{g}_k$. The difference of the
connection $\tilde{\Gamma}_k$ of $\tilde{g}_k$ and the connection
$\Gamma$ of $g$ is a tensor. Taking $\Gamma$ to be fixed in time,
we get
\begin{align*}
\frac{\partial}{\partial t}(\tilde{\Gamma}_k-\Gamma)& =
\frac{\partial}{\partial
t}\(\frac{1}{2}(\tilde{g}_k)^{\gamma\delta}\left[\frac{\partial}{\partial
x^\alpha}(\tilde{g}_k)_{\delta\beta}+\frac{\partial}{\partial
x^\beta}(\tilde{g}_k)_{\delta\alpha}-\frac{\partial}{\partial
x^\delta}(\tilde{g}_k)_{\alpha\beta}\right]\)\\
& = \frac{1}{2}(\tilde{g}_k)^{\gamma\delta}
\Big[(\tilde{\nabla}_k)_\alpha(-2(\tilde{\Ric}_k)_{\beta\delta})
+(\tilde{\nabla}_k)_\beta(-2(\tilde{\Ric}_k)_{\alpha\delta}) \\
&\qquad-(\tilde{\nabla}_k)_\delta(-2(\tilde{\Ric}_k)_{\alpha\beta})\Big]
\end{align*}
and then by the assumption (c) and (4.1.10),
$$
\left|\frac{\partial}{\partial
t}(\tilde{\Gamma}_k-\Gamma)\right|\leq C,\quad \text{for all }\;
k.
$$
Note also that at a normal coordinate of the metric $g$ at a fixed
point and at the time $t=0$,
\begin{align}
(\tilde{\Gamma}_k)^\gamma_{\alpha\beta}-\Gamma^\gamma_{\alpha\beta}
&=
\frac{1}{2}(\tilde{g}_k)^{\gamma\delta}\(\frac{\partial}{\partial
x^\alpha}(\tilde{g}_k)_{\delta\beta}+\frac{\partial}{\partial
x^\beta}(\tilde{g}_k)_{\delta\alpha}-\frac{\partial}{\partial
x^\delta}(\tilde{g}_k)_{\alpha\beta}\)\\
& =\frac{1}{2}(\tilde{g}_k)^{\gamma\delta}
(\nabla_\alpha(\tilde{g}_k)_{\delta\beta}
+\nabla_\beta(\tilde{g}_k)_{\delta\alpha}
-\nabla_\delta(\tilde{g}_k)_{\alpha\beta}),\nn
\end{align}  
thus by the assumption (b) and (4.1.10),
$$
|\tilde{\Gamma}_k(0)-\Gamma|\leq C, \; \text{ for all }\; k.
$$
Integrating over time we deduce that \be
|\tilde{\Gamma}_k-\Gamma|\leq C, \; \text{ on }\;
K\times[\alpha,\beta],\; \text{ for all }\; k.  
\ee By using the assumption (c) and (4.1.10) again, we have
\begin{align*}
\left|\frac{\partial}{\partial t}(\nabla\tilde{g}_k)\right|
&  = |-2\nabla\tilde{\Ric}_k|\\
&  = |-2\tilde{\nabla}_k\tilde{\Ric}_k+(\tilde{\Gamma}_k
-\Gamma)*\tilde{\Ric}_k|\\
&  \leq C,\; \text{ for all }\; k.
\end{align*}
Hence by combining with the assumption (b) we get bounds \be
|\nabla\tilde{g}_k|\leq\tilde{C}_1,\; \text{ on }\;
K\times[\alpha,\beta],
\ee for some positive constant $\tilde{C}_1$ independent of $k$.

Further we want to bound $\nabla^2\tilde{g}_k$. Again regarding
$\nabla$ as fixed in time, we see
$$
\frac{\partial}{\partial
t}(\nabla^2\tilde{g}_k)=-2\nabla^2(\tilde{R}ic_k).
$$
Write
\begin{align*}
\nabla^2\tilde{\Ric}_k& =
(\nabla-\tilde{\nabla}_k)(\nabla\tilde{\Ric}_k)+
\tilde{\nabla}_k(\nabla-\tilde{\nabla}_k)\tilde{\Ric}_k
+\tilde{\nabla}^2_k\tilde{\Ric}_k\\
&  = (\Gamma-\tilde{\Gamma}_k)*\nabla\tilde{\Ric}_k+
\tilde{\nabla}_k((\Gamma-\tilde{\Gamma}_k)*\tilde{\Ric}_k)
+\tilde{\nabla}^2_k\tilde{\Ric}_k\\
&
=(\Gamma-\tilde{\Gamma}_k)*[(\nabla-\tilde{\nabla}_k)\tilde{\Ric}_k
+\tilde{\nabla}_k\tilde{\Ric}_k] \\
&\quad+
\tilde{\nabla}_k(\tilde{g}_k^{-1}*\nabla\tilde{g}_k*\tilde{\Ric}_k)
+\tilde{\nabla}^2_k\tilde{\Ric}_k\\
&
=(\Gamma-\tilde{\Gamma}_k)*[(\Gamma-\tilde{\Gamma}_k)*\tilde{\Ric}_k
+\tilde{\nabla}_k\tilde{\Ric}_k] \\
&\quad+
\tilde{\nabla}_k(\tilde{g}_k^{-1}*\nabla\tilde{g}_k*\tilde{\Ric}_k)
+\tilde{\nabla}^2_k\tilde{\Ric}_k
\end{align*}
where we have used (4.1.11). Then by the assumption (c), (4.1.10),
(4.1.12) and (4.1.13) we have
\begin{align*}
|\frac{\partial}{\partial t}\nabla^2\tilde{g}_k|
&  \leq  C+C\cdot|\tilde{\nabla}_k\nabla\tilde{g}_k|\\
&  = C+C\cdot|\nabla^2\tilde{g}_k
+(\tilde{\Gamma}_k-\Gamma)*\nabla\tilde{g}_k|\\
&  \leq C+C|\nabla^2\tilde{g}_k|.
\end{align*}
Hence by combining with the assumption (b) we get
$$
|\nabla^2\tilde{g}_k|\leq\tilde{C}_2,\; \mbox{ on }\;
K\times[\alpha,\beta],
$$
for some positive constant $\tilde{C}_2$ independent of $k$.

The bounds on the higher derivatives can be derived by the same
argument. Therefore we have completed the proof of the lemma.
\end{pf}

We now apply the lemma to the pull-back metrics
$\tilde{g}_k(t)=(f_k)^*g_k(t)$ on
$B_{\infty}(p_{\infty},s_{\infty}) \times (A,\Omega]$. Since the
metrics $\tilde{g}_k(0)$ have uniform bounds on their curvature
and all derivatives of their curvature on every compact set of
$B_{\infty}(p_{\infty},s_{\infty})$ and converge to the metric
$g_{\infty}(0)$ in $C^\infty_{loc}$ topology, the assumptions (a)
and (b) are certainly held for every compact subset $K \subset
B_{\infty}(p_{\infty},s_{\infty})$ with $g=g_{\infty}(0)$. For
every compact subinterval $[\alpha,\beta] \subset (A,\Omega]$, we
have already seen from Shi's derivative estimates (Theorem 1.4.1)
that the assumption (c) is also held on $K \times [\alpha,\beta]$.
Then all of the $\nabla^l\tilde{g}_k$ are uniformly bounded with
respect to the fixed metric $g = g_{\infty}(0)$ on every compact
set of $B_{\infty}(p_{\infty},s_{\infty})\times(A,\Omega]$. By
using the classical Arzela-Ascoli theorem, we can find a
subsequence which converges uniformly together with all its
derivatives on every compact subset of
$B_{\infty}(p_{\infty},s_{\infty})\times(A,\Omega]$. The limit
metric will agree with that obtained previously at $t=0$, where we
know its convergence already. The limit $g_{\infty}(t),\
t\in(A,\Omega]$, is now clearly itself a solution of the Ricci
flow. Thus we obtain the following Cheeger type compactness
theorem to the Ricci flow, which is essentially obtained by
Hamilton in \cite{Ha95} and is called \textbf{Hamilton's
compactness theorem}\index{Hamilton's compactness theorem}.


\begin{theorem}[Hamilton's compactness theorem]
Let $(M_k,g_k(t),p_k),$ $t\in(A,\Omega]$ with $A<0 \leq \Omega$,
be a sequence of evolving marked complete Riemannian manifolds.
Consider a sequence of geodesic balls $B_0(p_k,s_k) \subset M_k$
of radii $s_k (0<s_k \leq +{\infty})$, with $s_k \rightarrow
s_{\infty}$ $(\leq +{\infty})$, around the base points $p_k$ in
the metrics $g_k(0)$. Suppose each $g_k(t)$ is a solution to the
Ricci flow on $B_0(p_k,s_k) \times (A,\Omega]$. Suppose also
\begin{itemize}
\item[(i)] for every radius $r<s_{\infty}$ there exist positive
constants $C(r)$ and $k(r)$ independent of $k$ such that the
curvature tensors $Rm(g_k)$ of the evolving metrics $g_k(t)$
satisfy the bound
$$
|Rm(g_k)|\leq C(r),
$$
on $B_0(p_k,r) \times (A,\Omega]$ for all $k \geq k(r)$, and
\item[(ii)] there exists a constant $\delta > 0$ such that the
injectivity radii of $M_k$ at $p_k$ in the metric $g_k(0)$ satisfy
the bound
$$
{\rm inj}\,(M_k,p_k,g_k(0))\geq\delta>0,
$$
for all $k = 1, 2, \ldots$.
\end{itemize}
Then there exists a subsequence of evolving marked
$(B_0(p_k,s_k),g_k(t),p_k)$ over $t \in (A,\Omega]$ which converge
in $C^\infty_{loc}$ topology to a solution
$(B_{\infty},g_{\infty}(t),p_{\infty})$ over $t \in (A,\Omega]$ to
the Ricci flow, where, at the time $t=0$, $B_{\infty}$ is a
geodesic open ball centered at $p_{\infty}\in B_{\infty}$ with the
radius $s_{\infty}$. Moreover the limiting solution is complete if
$s_{\infty}=+\infty$.
\end{theorem}

\section{Injectivity Radius Estimates}

We will use rescaling arguments to understand the formation of
singularities and long-time behaviors of the Ricci flow. In view
of the compactness property obtained in the previous section, on
one hand one needs to control the bounds on the curvature, and on
the other hand one needs to control the lower bounds of the
injectivity radius. In applications we usually rescale the
solution so that the (rescaled) curvatures become uniformly
bounded on compact subsets and leave the injectivity radii of the
(rescaled) solutions to be estimated in terms of curvatures. In
this section we will review a number of such injectivity radius
estimates in Riemannian geometry. The combination of these
injectivity estimates with Perelman's no local collapsing theorem
I$'$ yields the well-known little loop lemma to the Ricci flow
which was conjectured by Hamilton in \cite{Ha95F}.

Let $M$ be a Riemannian manifold. Recall that the {\bf injectivity
radius}\index{injectivity radius} at a point $p\in M$ is defined
by
$$
{\rm inj}\,(M,p)=\sup\{r>0\ |\  \exp_p: B(O,r)(\subset T_pM)\to M
\; \text{ is injective}\},
$$
and the {\bf injectivity radius of} ${\mathbf M}$ is
$$
{\rm inj}\,(M)=\inf\{ {\rm inj}\,(M,p)\ |\  p\in M\}.
$$
We begin with a basic lemma due to Klingenberg (cf. Corollary 5.7
of Cheeger \& Ebin \cite{CE}).

\medskip
{\bf Klingenberg's Lemma.}\index{Klingenberg's lemma} \emph{ \ Let
$M$ be a complete Riemannian manifold and let $p\in M$. Let
$l_M(p)$ denote the minimal length of a nontrivial geodesic loop
starting and ending at $p$ $($maybe not smooth at $p)$. Then the
injectivity radius of $M$ at $p$ satisfies the inequality
$$
{\rm inj}\,(M,p)\geq
\min\left\{\frac{\pi}{\sqrt{K_{\max}}},\frac{1}{2}l_M(p)\right\}
$$
where $K_{\max}$ denotes the supermum of the sectional curvature
on $M$ and we understand ${\pi}/{\sqrt{K_{\max}}}$ to be positive
infinity if $K_{\max}\leq0$.}

\medskip
Based on this lemma and a second variation argument, Klingenberg
proved that the injectivity radius of an even-dimensional,
compact, simply connected Riemannian manifold of positive
sectional curvature is bounded from below by
$\pi/\sqrt{K_{\max}}$. For odd-dimensional, compact, simply
connected Riemannian manifold of positive sectional curvature, the
same injectivity radius estimates was also proved by Klingenberg
under an additional assumption that the sectional curvature is
strictly $\frac{1}{4}$-pinched (cf. Theorem 5.9 and 5.10 of
\cite{CE}). We also remark that in dimension 7, there exists a
sequence of simply connected, homogeneous Einstein spaces whose
sectional curvatures are positive and uniformly bounded from above
but their injectivity radii converge to zero. (See \cite{AW}.)

The next result, due to Gromoll and Meyer \cite{GM}, shows that
for complete noncompact Riemannian manifold with positive
sectional curvature, the above injectivity radius estimate
actually holds without any restriction on dimensions. Since the
result and proof were not explicitly given in \cite{GM}, we
include a proof here.


\begin{theorem}[The Gromoll-Meyer injectivity radius
estimate]\index{Gromoll-Meyer injectivity radius estimate} Let $M$
be a complete, noncompact Riemannian manifold with positive
sectional curvature. Then the injectivity radius of $M$ satisfies
the following estimate
$$
{\rm inj}\,(M)\geq\frac{\pi}{\sqrt{K_{\max}}}.
$$
\end{theorem}

\begin{pf}
Let $O$ be an arbitrary fixed point in $M$. We need to show that
the injectivity radius at $O$ is not less than
$\pi/\sqrt{K_{\max}}$. We argue by contradiction.  Suppose not,
then by Klingenberg's lemma there exists a closed geodesic loop
$\gamma$ on $M$ starting and ending at $O$ (may be not smooth at
$O$).

Since $M$ has positive sectional curvature, we know from the work
of Gromoll-Meyer \cite{GM} (also cf. Proposition 8.5 of \cite{CE})
that there exists a compact totally convex subset $C$ of $M$
containing the geodesic loop $\gamma$. Among all geodesic loops
starting and ending at the same point and lying entirely in the
compact totally convex set $C$ there will be a shortest one. Call
it $\gamma_0$, and suppose $\gamma_0$ starts and ends at a point
we call $p_0$.

First we claim that $\gamma_0$ must be also smooth at the point
$p_0$. Indeed by the curvature bound and implicit function
theorem, there will be a geodesic loop $\tilde{\gamma}$ close to
$\gamma_0$ starting and ending at any point $\tilde{p}$ close to
$p_0$. Let $\tilde{p}$ be along $\gamma_0$. Then by total
convexity of the set $C$, $\tilde{\gamma}$ also lies entirely in
$C$. If $\gamma_0$ makes an angle different from $\pi$ at $p_0$,
the first variation formula will imply that $\tilde{\gamma}$ is
shorter than $\gamma_0$. This contradicts with the choice of the
geodesic loop $\gamma_0$ being the shortest.

Now let $L:[0,+\infty)\rightarrow M$ be a ray emanating from
$p_0$. Choose $r>0$ large enough and set $q=L(r)$. Consider the
distance between $q$ and the geodesic loop $\gamma_0$. It is clear
that the distance can be realized by a geodesic $\beta$ connecting
the point $q$ to a point $p$ on $\gamma_0$.

Let $X$ be the unit tangent vector of the geodesic loop $\gamma_0$
at $p$. Clearly $X$ is orthogonal to the tangent vector of $\beta$
at $p$. We then translate the vector $X$ along the geodesic
$\beta$ to get a parallel vector field $X(t),\ 0\leq t\leq r$. By
using this vector field we can form a variation fixing one
endpoint $q$ and the other on $\gamma_0$ such that the variational
vector field is $(1-\frac{t}{r})X(t)$. The second variation of the
arclength of this family of curves is given by
\begin{align*}
&  I\(\(1-\frac{t}{r}\)X(t),\(1-\frac{t}{r}\)X(t)\)\\
&  = \int_0^r\bigg[\left|\nabla_{\frac{\partial}{\partial
t}}\(\(1-\frac{t}{r}\)X(t)\)\right|^2 \\
&\qquad-R\(\frac{\partial}{\partial t},\(1-\frac{t}{r}\)X(t),
\frac{\partial}{\partial t},
\(1-\frac{t}{r}\)X(t)\)\bigg]dt\\
&
=\frac{1}{r}-\int_0^r\(1-\frac{t}{r}\)^2R\(\frac{\partial}{\partial
t},X(t), \frac{\partial}{\partial t},X(t)\)dt\\
&  < 0
\end{align*}
when $r$ is sufficiently large, since the sectional curvature of
$M$ is strictly positive everywhere. This contradicts with the
fact that $\beta$ is the shortest geodesic connecting the point
$q$ to the shortest geodesic loop $\gamma_0$. Thus we have proved
the injectivity radius estimate.
\end{pf}

In contrast to the above injectivity radius estimates, the
following well-known injectivity radius estimate of Cheeger (cf.
Theorem 5.8 of \cite{CE}) does not impose the restriction on the
sign of the sectional curvature.

\medskip
{\bf Cheeger's Lemma.} \emph{ \ Let $M$ be an $n$-dimensional
compact Riemannian manifold with the sectional curvature
$|K_M|\leq\lambda$, the diameter $d(M)\leq D$, and the volume
$\Vol(M)\geq v > 0$. Then, we have
$$
{\rm inj}\,(M)\ge C_n(\lambda,D,v)
$$
for some positive constant $C_n(\lambda,D,v)$ depending only on
$\lambda,D,v$ and the dimension $n$.}\index{Cheeger's lemma}

\medskip

For general manifolds, there are localized versions of the above
Cheeger's Lemma. In 1981, Cheng-Li-Yau \cite{CLY} first obtained
the important estimate of local injectivity radius under the
normalization of the injectivity radius at any fixed base point.
Their result is what Hamilton needed to prove his compactness
result in \cite{Ha95}. Here, we also need the following version of
the local injectivity radius estimate, which appeared in the 1982
paper of Cheeger-Gromov-Taylor \cite{CGT}, in which the
normalization is in terms of local volume of a ball. (According to
Yau, the argument between local volume and local injectivity
radius was, however, initiated by him during a conversation with
Gromov in 1975 in an explanation of his paper \cite{Y76} on
proving complete manifolds with positive Ricci curvature have
infinite volume.)

\begin{theorem}
Let $B(x_0, 4r_0)$, $0<r_0<\infty$, be a geodesic ball in an
$n$-dimensional complete Riemannian manifold $(M,g)$ such that the
sectional curvature $K$ of the metric $g$ on $B(x_0,4r_0)$
satisfies the bounds
$$
\lambda\leq K\leq\Lambda
$$
for some constants $\lambda$ and $\Lambda$. Then for any positive
constant $r\le r_0$ $($we will also require $r\le
\pi/(4\sqrt{\Lambda})$ if $\Lambda>0)$ the injectivity radius of
$M$ at $x_0$ can be bounded from below by
$$
\mbox{\rm inj}(M, x_0)\geq r\cdot\frac{\Vol(B(x_0,
r))}{\Vol(B(x_0,r))+V^n_{\lambda}(2r)},
$$
where $V_{\lambda}^n(2r)$ denotes the volume of a geodesic ball of
radius $2r$ in the $n$-dimensional simply connected space form
$M_{\lambda}$ with constant sectional curvature $\lambda$.
\end{theorem}

\begin{pf} The following proof is essentially from Cheeger-Gromov-Taylor \cite{CGT}
(cf. also \cite{AM}).

It is well known (cf. Lemma 5.6 of \cite{CE}) that
$$
\mbox{inj}(M, x_0)=\min \left\{\mbox{conjugate radius of}\ x_0, \
\frac {1}{2}l_M(x_0)\right\}
$$
where $l_M(x_0)$ denotes the length of the shortest (nontrivial)
closed geodesic starting and ending at $x_0$. Since by assumption
$r\le \pi/(4\sqrt{\Lambda})$ if $\Lambda>0$, the conjugate radius
of $x_0$ is at least $4r$. Thus it suffices to show \be
l_M(x_0)\geq
2r\cdot\frac{\Vol(B(x_0,r))}{\Vol(B(x_0,r))+V^n_{\lambda}(2r)}.
\ee The idea of Cheeger-Gromov-Taylor \cite{CGT} for proving this
inequality, as indicated in \cite{AM}, is to compare the geometry of
the ball $B(x_0,4r)$ $\subseteq B(x_0, 4r_0)\subset M$ with the
geometry of its lifting $\tilde B_{4r}\subset T_{x_0}(M)$, via the
exponential map $\mbox{exp}_{x_0}$, equipped with the pull-back
metric $\tilde g=\mbox{exp}^{*}_{x_0} g$. Thus $\mbox{exp}_{x_0}:
\tilde B_{4r} \to B(x_0,4r)$ is a length-preserving local
diffeomorphism.

Let $\tilde x_0, \tilde x_1, \ldots, \tilde x_N$ be the preimages
of $x_0$ in $\tilde B_r\subset \tilde B_{4r}$ with $\tilde x_0=0$.
Clearly they one-to-one correspond to the geodesic loops
$\gamma_0, \gamma_1, \ldots, \gamma_N$ at $x_0$ of length less
than $r$, where $\gamma_0$ is the trivial loop. Now for each point
$\tilde x_i$ there exists exactly one isometric immersion
$\varphi_i : \tilde B_r \to \tilde B_{4r}$ mapping $0$ to $\tilde
x_i$ and such that $\mbox{exp}_{x_0} \varphi_i=\mbox{exp}_{x_0}$.

Without loss of generality, we may assume $\gamma_1$ is the
shortest nontrivial geodesic loop at $x_0$. By analyzing short
homotopies, one finds that $\varphi_i(\tilde x)\neq
\varphi_j(\tilde x)$ for all $\tilde x\in \tilde B_r$ and $0\le
i<j\le N$. This fact has two consequences:

\medskip
(a) \ $N\ge 2m$, where $m=[r/l_M(x_0)]$. To see this, we first
observe that the points $\varphi_1^k(0), -m\le k\le m$, are
preimages of $x_0$ in $\tilde B_r$ because $\varphi_1$ is an
isometric immersion satisfying $\mbox{exp}_{x_0}
\varphi_1=\mbox{exp}_{x_0}$. Moreover we claim they are distinct.
For otherwise $\varphi_1$ would act as a permutation on the set
$\{\varphi_1^k(0)\ |\ -m\le k\le m\}$. Since the induced metric
$\tilde g$ at each point in $\tilde B_r$ has the injectivity
radius at least $2r$, it follows from the Whitehead theorem (see
for example \cite{CE}) that $\tilde B_r$ is geodesically convex.
Then there would exist the unique center of mass $\tilde y \in
\tilde{B_r}$. But then $\tilde y=\varphi_0(\tilde
y)=\varphi_1(\tilde y)$, a contradiction.

\medskip\noindent
\;\;\;(b) Each point in $B(x_0,r)$ has at least $N\!+\!1$
preimages in $\Omega\!=\!\cup_{i=0}^N B(\tilde x_i, r)$ $\subset
\tilde B_{2r}$. Hence by the Bishop volume comparison,
$$
(N+1)\mbox{Vol}\,(B(x_0,r)) \leq \mbox{Vol}_{\tilde g}(\Omega)
\leq \mbox{Vol}_{\tilde g}(\tilde B_{2r})\le V^n_{\lambda}(2r).
$$

Now the inequality (4.2.1) follows by combining the fact $N\ge
2[r/l_M(x_0)]$ with the above volume estimate.
\end{pf}

For our later applications, we now consider in a complete
Riemannian manifold $M$ a geodesic ball $B(p_0, s_0)$ $(0<s_0\le
\infty)$ with the property that there exists a positive increasing
function $\Lambda : [0, s_0)\to [0, \infty)$ such that for any
$0<s<s_0$ the sectional curvature $K$ on the ball $B(p_0, s)$ of
radius $s$ around $p_0$ satisfies the bound
$$
|K|\le \Lambda(s).
$$
Using Theorem 4.2.2, we can control the injectivity radius at any
point $p\in B(p_0, s_0)$ in terms a positive constant that depends
only on the dimension $n$, the injectivity radius at the base
point $p_0$, the function $\Lambda$ and the distance $d(p_0, p)$
from $p$ to $p_0$. We now proceed to derive such an estimate. The
geometric insight of the following argument belongs to Yau
\cite{Y76} where he obtained a lower bound estimate for volume by
comparing various geodesic balls.  Indeed, it is a finite version
of Yau's Busemann function argument which gives the information on
comparing geodesic balls with centers far apart.

For any point $p\in B(p_0, s_0)$ with $d(p_0, p)=s$, set
$r_0=(s_0-s)/4$ (we define $r_0$=1 if $s_0=\infty$). Define the
set $S$ to be the union of minimal geodesic segments that connect
$p$ to each point in $B(p_0, r_0)$. Now any point $q\in S$ has
distance at most
$$
r_0+r_0+s=s+2r_0
$$
from $p_0$ and hence $S\subseteq B(p_0, s+2r_0)$. For any $0<r\le
\min\{\pi/4\sqrt{\Lambda(s+2r_0)},$ $r_0\}$, we denote by $\alpha
(p, r)$ the sector $S\cap B(p,r)$ of radius $r$ and by $\alpha (p,
s+r_0)=S \cap B(p,s+r_0)$. Let $\alpha_{-\Lambda(s+2r_0)}(r_0)$
(resp. $\alpha_{-\Lambda(s+2r_0)} (s+r_0)$) be a corresponding
sector of the same ``angles" with radius $r_0$ (resp. $s+r_0$) in
the $n$-dimensional simply connected space form with constant
sectional curvature $-\Lambda(s+2r_0)$. Since $B(p_0, r_0)\subset
S\subset\alpha (p,s+r_0)$ and $\alpha (p, r)\subset B(p,r)$, the
Bishop-Gromov volume comparison theorem implies that
\begin{align*}
\frac{\Vol(B(p_0,r_0))}{\Vol(B(p,r))}
&\le \frac{\Vol(\alpha (p,s+r_0))}{\Vol(\alpha (p,r))} \\
&\le\frac{\Vol(\alpha_{-\Lambda(s+2r_0)}(s+r_0))}{\Vol(\alpha_
{-\Lambda(s+2r_0)}(r))}=\frac{V^n_{-\Lambda(s+2r_0)}(s+r_0)}{V^n_
{-\Lambda(s+2r_0)}(r)}.
\end{align*}
Combining this inequality with the local injectivity radius
estimate in Theorem 4.2.2, we get
\begin{align*}
&{\rm inj}\,(M,p) \\
&\geq r \frac{V^{n}_{-\Lambda(s+2r_0)}(r)\cdot
\Vol(B(p_0,r_0))}{V^{n}_{-\Lambda(s+2r_0)}(r)\Vol(B(p_0,r_0))+
V_{-\Lambda(s+2r_0)}^n(2r)V^{n}_{-\Lambda(s+2r_0)}(s+2r_0)}.
\end{align*}
Thus, we have proved the following (cf. \cite{CLY}, \cite{CGT} or
\cite{AM})

\begin{corollary}
Suppose $B(p_0, s_0)$ $(0<s_0\le \infty)$ is a geodesic ball in an
$n$-dimensional complete Riemannian manifold $M$ having the
property that for any $0<s<s_0$ the sectional curvature $K$ on
$B(p_0, s)$ satisfies the bound
$$
|K|\le \Lambda(s)
$$
for\; some\; positive\; increasing\; function\; $\Lambda$\;
defined\; on\; $[0,\,s_0)$.\; Then\; for\; any\; point\; $p\;\in\;
B(p_0,\, s_0)$\; with\; $d(p_0,\, p)\;=\;s$\; and\; any\;
positive\; number $r\le \min\{\pi/4\sqrt{\Lambda(s+2r_0)},r_0\}$
with $r_0=(s_0-s)/4$, the injectivity radius of $M$ at $p$ is
bounded below by
\begin{align*}
& {\rm inj}\,(M,p) \\
&\geq r \frac{V^{n}_{-\Lambda(s+2r_0)}(r)\cdot
\Vol(B(p_0,r_0))}{V^{n}_{-\Lambda(s+2r_0)}(r)\Vol(B(p_0,r_0))+
V_{-\Lambda(s+2r_0)}^n(2r)V^{n}_{-\Lambda(s+2r_0)}(s+2r_0)}.
\end{align*}
In particular, we have \be
\inj(M,p)\geq \rho_{n, \delta, \Lambda}(s) 
\ee where $\delta>0$ is a lower bound of the injectivity radius
$\inj(M, p_0)$ at the origin $p_0$ and $\rho_{n,\delta, \Lambda}:
[0, s_0) \to \mathbb{R^+}$ is a positive decreasing function that
depends only on the dimension $n$, the lower bound $\delta$ of the
injectivity radius $\inj(M,p_0)$, and the function $\Lambda$.
\end{corollary}

We remark that in the above discussion if $s_0=\infty$ then we can
apply the standard Bishop relative volume comparison theorem to
geodesic balls directly. Indeed, for any $p\in M$ and any positive
constants $r$ and $r_0$, we have $B(p_0, r_0)\subseteq B(p,
\hat{r})$ with
$\hat{r}\stackrel{\Delta}{=}\max\{r,r_0+d(p_0,p)\}$. Suppose in
addition the curvature $K$ on $M$ is uniformly bounded by
$\lambda\le K\le\Lambda$ for some constants $\lambda$ and
$\Lambda$, then the Bishop volume comparison theorem implies that
$$
\frac{\Vol(B(p_0,r_0))}{\Vol(B(p,r))}\le
\frac{\Vol(B(p,\hat{r}))}{\Vol(B(p,r))}\le
\frac{V_{\lambda}(\hat{r})}{V_{\lambda}(r)}.$$ Hence \be
\inj(M,p)\geq r\frac{V_{\lambda}^n(r)\cdot
\Vol(B(p_0,r_0))}{V_{\lambda}^n(r)\Vol(B(p_0,r_0))
+V_{\lambda}^n(2r)V_{\lambda}^n(\hat{r})}. 
\ee So we see that the injectivity radius $inj(M,p)$ at $p$ falls
off at worst exponentially as the distance $d(p_0,p)$ goes to
infinity. In other words, \be \inj(M,p)\geq
\frac{c}{\sqrt{B}}(\delta \sqrt{B})^n
e^{-C\sqrt{B}d(p,p_0)} 
\ee where $B$ is an upper bound on the absolute value of the
sectional curvature, $\delta$ is a lower bound on the injectivity
radius at $p_0$ with $\delta<c/\sqrt{B}$, and $c>0$ and $C<+\infty
$ are positive constants depending only on the dimension $n$.

\smallskip
Finally, by combining Theorem 4.2.2 with Perelman's no local
collapsing Theorem I${'}$ (Theorem 3.3.3) we immediately obtain the
following important (due to Perelman \cite{P1}) \textbf{Little Loop
Lemma} \index{Little Loop Lemma} conjectured by Hamilton
\cite{Ha95F}.

\begin{theorem}[Little Loop Lemma]
Let $g_{ij}(t)$, $0\le t<T<+\infty$, be a solution of the Ricci
flow on a compact manifold $M$. Then there exists a constant
$\rho>0$ having the following property: if at a point $x_0\in M$
and a time $t_0\in [0,T)$,
$$ |Rm|(\cdot,t_0)\leq r^{-2} \quad
\text{on}\ B_{t_0}(x_0,r)
$$
for some $r\leq \sqrt{T}$, then the injectivity radius of $M$ with
respect to the metric $g_{ij}(t_0)$ at $x_0$ is bounded from below
by
$$
\inj(M,x_0,g_{ij}(t_0))\geq \rho r.
$$
\end{theorem}

\section{Limiting Singularity Models}

In \cite{Ha95F}, Hamilton classified singularities of the Ricci
flow into three types and showed each type has a corresponding
singularity model. The main purpose of this section is to discuss
these results of Hamilton. Most of the presentation follows
Hamilton \cite{Ha95F}.

Consider a solution $g_{ij}(x,t)$ of the Ricci flow on $M\times
[0,T)$, $T\le +\infty$, where either $M$ is compact or at each
time $t$ the metric $g_{ij}(\cdot,t)$ is complete and has bounded
curvature. We say that $g_{ij}(x,t)$ is a {\bf maximal solution}
\index{maximal solution} if either $T=+\infty$ or $T<+\infty$ and
$|Rm|$ is unbounded as $t\to T$.

Denote by
$$
K_{\max}(t)=\sup_{x\in M}|Rm(x,t)|_{g_{ij}(t)}.
$$

\begin{definition}
We say that $\{(x_k,t_k) \in M\times [0,T) \}$, $k=1, 2, \ldots$,
is a sequence of {\bf (almost) maximum points}\index{(almost)
maximum points} if there exist positive constants $c_1$ and
$\alpha \in (0,1]$ such that
$$
|Rm(x_k,t_k)|\geq c_1 K_{\max}(t), \quad t\in
[t_k-\frac{\alpha}{K_{\max}(t_k)},t_k]
$$
for all $k$.
\end{definition}

\begin{definition}
We say that the solution satisfies {\bf injectivity radius
condition}\index{injectivity radius!condition} if for any sequence
of (almost) maximum points $\{(x_k,t_k)\}$, there exists a
constant $c_2>0$ independent of $k$ such that
$$
\inj(M,x_k,g_{ij}(t_k))\geq \frac{c_2}{\sqrt{K_{\max}(t_k)}} \quad
\text{for all }\;   k.
$$
\end{definition}

Clearly, by the Little Loop Lemma, a maximal solution on a compact
manifold with the maximal time $T<+\infty$ always satisfies the
injectivity radius condition. Also by the Gromoll-Meyer
injectivity radius estimate, a solution on a complete noncompact
manifold with positive sectional curvature also satisfies the
injectivity radius condition.

According to Hamilton \cite{Ha95F}, we can classify maximal
solutions into three types; every maximal solution is clearly of
one and only one of the following three types:
$$
\mbox{\textbf{Type I:}} \ \ \ T<+\infty\ \  \mbox{and}\ \
\sup_{t\in [0,T)}(T-t)K_{\max}(t)<+\infty;\
$$
\index{type I}
$$
\ \ \ \mbox{\textbf{Type II:}} \ \ \ \text{(a)}\  T<+\infty\ \
\mbox{but}\ \ \sup_{t\in [0,T)}(T-t)K_{\max}(t)=+\infty;
$$
\index{type II!(a)}
$$
 \hskip 2.3cm \text{(b)}\  T=+\infty\ \  \mbox{but}\ \
\sup_{t\in [0,T)}tK_{\max}(t)=+\infty;\ \ \ \ \ \
$$
\index{type II!(b)}
$$
\mbox{\textbf{Type III:}} \ \ \ \text{(a)}\  T=+\infty,\ \ \ \
\sup_{t\in [0,T)}tK_{\max}(t)<+\infty,\ \ \mbox{and} \ \
$$
$$ \limsup_{t\rightarrow
+\infty}tK_{\max}(t)>0;
$$
\index{type III!(a)}
$$
\hskip 2.3cm \text{(b)}\  T=+\infty,\ \ \ \  \sup_{t\in
[0,T)}tK_{\max}(t)<+\infty,\ \ \mbox{and} \ \
$$
$$ \limsup_{t\rightarrow
+\infty}tK_{\max}(t)=0;
$$
\index{type III!(b)}

Note that Type III (b) is not compatible with the injectivity
radius condition unless it is a trivial flat solution. (Indeed
under the Ricci flow the length of a curve $\gamma$ connecting two
points $x_0$, $x_1\in M$ evolves by
\begin{align*}
\frac{d}{dt}L_t(\gamma)
&  = \int_{\gamma}-\Ric(\dot{\gamma},\dot{\gamma})ds\\
&  \leq C(n)K_{\max}(t)\cdot L_t(\gamma)\\
&  \leq \frac{\epsilon}{t}L_t(\gamma),\; \text{ as $t$ large
enough,}
\end{align*}
for arbitrarily fixed $\epsilon>0$. Thus when we are considering
the Ricci flow on a compact manifold, the diameter of the evolving
manifold grows at most as $t^{\epsilon}$. But the curvature of the
evolving manifold decays faster than $t^{-1}$. This says, as
choosing $\epsilon>0$ small enough,
$$
{\rm diam}_t(M)^2\cdot |Rm(\cdot,t)|\rightarrow0,\; \text{ as }\;
t\rightarrow+\infty.
$$
Then it is well-known from Cheeger-Gromov \cite{Gro78} that the
manifold is a nilmanifold and the injectivity radius condition can
not be satisfied as $t$ large enough. When we are considering the
Ricci flow on a complete noncompact manifold with nonnegative
curvature operator or on a complete noncompact K\"ahler manifold
with nonnegative holomorphic bisectional curvature,
Li-Yau-Hamilton inequalities imply that $tR(x,t)$ is increasing in
time $t$. Then Type III(b) occurs only when the solution is a
trivial flat metric.)

For each type of maximal solution, Hamilton \cite{Ha95F} defined a
corresponding type of limiting singularity model.

\begin{definition}
A solution $g_{ij}(x,t)$ to the Ricci flow on the manifold $M$,
where either $M$ is compact or at each time $t$ the metric
$g_{ij}(\cdot,t)$ is complete and has bounded curvature, is called
a {\bf singularity model}\index{singularity model} if it is not
flat and of one of the following three types:

\vskip 0.1cm\noindent \textbf{Type I}\index{type I}: The solution
exists for $t\in(-\infty,\Omega)$ for some constant $\Omega$ with
$0<\Omega<+\infty$ and
$$
|Rm|\leq\Omega/(\Omega-t)
$$
everywhere with equality somewhere at $t=0$;

\vskip 0.1cm\noindent \textbf{Type II}\index{type II}: The
solution exists for $t\in(-\infty,+\infty)$  and
$$
|Rm|\leq 1 \ \ \ \ \ \ \ \ \
$$
everywhere with equality somewhere at $t=0$;

\vskip 0.1cm\noindent \textbf{Type III}\index{type III}: The
solution exists for $t\in(-A,+\infty)$ for some constant $A$ with
$0<A<+\infty$ and
$$
\ \ |Rm|\leq A/(A+t)
$$
everywhere with equality somewhere at $t=0$.
\end{definition}

\begin{theorem} [{Hamilton \cite{Ha95F}}]
For any maximal solution to the Ricci flow which satisfies the
injectivity radius condition and is of Type {\rm I, II(a), (b),}
or {\rm III(a),} there exists a sequence of dilations of the
solution along (almost) maximum points which converges in the
$C_{\rm loc}^{\infty}$ topology to a singularity model of the
corresponding type.
\end{theorem}

\begin{pf} \

\smallskip
{\bf Type I:}\index{type I} \ We consider a maximal solution
$g_{ij}(x,t)$ on $M\times [0,T)$ with $T<+\infty$ and
$$
\Omega\stackrel{\Delta}{=}\limsup_{t\rightarrow T}
(T-t)K_{\max}(t)<+\infty.
$$

First we note that $\Omega>0$. Indeed by the evolution equation of
curvature,
$$
\frac{d}{dt}K_{\max}(t)\leq{\rm Const}\, \cdot K_{\max}^2(t).
$$
This implies that
$$
K_{\max}(t)\cdot (T-t)\geq {\rm Const}\, >0,
$$
because
$$
\limsup_{t\rightarrow T}K_{\max}(t)=+\infty.
$$
Thus $\Omega$ must be positive.

Choose a sequence of points $x_k$ and times $t_k$ such that
$t_k\rightarrow T$ and
$$
\lim_{k\rightarrow \infty}(T-t_k)|Rm(x_k,t_k)|=\Omega.
$$
Denote by
$$
\epsilon_k=\frac{1}{\sqrt{|Rm(x_k,t_k)|}}.
$$
We translate in time so that $t_k$ becomes 0, dilate in space by
the factor $\epsilon_k$ and dilate in time by $\epsilon_k^2$ to
get
$$
\tilde{g}_{ij}^{(k)}(\cdot,\tilde{t})
=\epsilon_k^{-2}g_{ij}(\cdot,t_k+\epsilon_k^2\tilde{t}),\ \ \ \
\tilde{t}\in [-t_k/\epsilon_k^2, (T-t_k)/\epsilon_k^2).
$$
Then
\begin{align*}
\frac{\partial}{\partial
\tilde{t}}\tilde{g}_{ij}^{(k)}(\cdot,\tilde{t})& =
\epsilon_k^{-2}\frac{\partial
}{\partial t}g_{ij}(\cdot,t)\cdot \epsilon_k^2 \\
&  = -2R_{ij}(\cdot,t_k+\epsilon_k^2\tilde{t}) \\
&  = -2\tilde{R}_{ij}^{(k)}(\cdot,\tilde{t}),
\end{align*}
where $\tilde{R}_{ij}^{(k)}$ is the Ricci curvature of the metric
$\tilde{g}_{ij}^{(k)}$. So $\tilde{g}_{ij}^{(k)}(\cdot,\tilde{t})$
is still a solution to the Ricci flow which exists on the time
interval $[-t_k/\epsilon_k^2,(T-t_k)/\epsilon_k^2)$, where
$$
t_k/\epsilon_k^2=t_k|Rm(x_k,t_k)|\to +\infty
$$
and
$$
(T-t_k)/\epsilon_k^2=(T-t_k)|Rm(x_k,t_k)|\to \Omega.
$$

For any $\epsilon>0$ we can find a time $\tau <T $ such that for
$t\in [\tau,T)$,
$$
|Rm|\leq (\Omega+\epsilon)/(T-t)
$$
by the assumption. Then for $\tilde{t}\in
[(\tau-t_k)/\epsilon_k^2,(T-t_k)/\epsilon_k^2)$, the curvature of
$\tilde{g}_{ij}^{(k)}(\cdot,\tilde{t})$ is bounded by
\begin{align*}
|\tilde{R}m^{(k)}|&  = \epsilon_k^2|Rm| \\
&  \leq(\Omega+\epsilon)/((T-t)|Rm(x_k,t_k)|) \\
& =(\Omega+\epsilon)/((T-t_k)|Rm(x_k,t_k)|+(t_k-t)|Rm(x_k,t_k)|)\\
&  \rightarrow (\Omega+\epsilon)/(\Omega-\tilde{t}), \qquad
\text{as } \; k\rightarrow +\infty.
\end{align*}
This implies that $\{(x_k,t_k)\}$ is a sequence of (almost)
maximum points. And then by the injectivity radius condition and
Hamilton's compactness theorem 4.1.5, there exists a subsequence
of the metrics $\tilde{g}_{ij}^{(k)}(\tilde {t})$ which converges
in the $C_{loc}^{\infty}$ topology to a limit metric
$\tilde{g}_{ij}^{(\infty)}(\tilde {t})$ on a limiting manifold
$\tilde{M}$ with $\tilde{t}\in(-\infty,\Omega)$ such that
$\tilde{g}_{ij}^{(\infty)}(\tilde {t})$ is a complete solution of
the Ricci flow and its curvature satisfies the bound
$$
|\tilde{R}m^{(\infty)}|\leq \Omega/(\Omega-\tilde{t})
$$
everywhere on $\tilde{M}\times (-\infty,\Omega)$ with the equality
somewhere at $\tilde{t}=0$.

\medskip
{\bf Type II(a):}\index{type II!(a)} \ We consider a maximal
solution $g_{ij}(x,t)$ on $M\times [0,T)$ with
$$
T<+\infty \quad\text{and} \quad \limsup_{t\rightarrow
T}(T-t)K_{\max}(t)=+\infty.
$$

Let $T_k<T<+\infty$ with $T_k\rightarrow T$, and $\gamma_k\nearrow
1$, as $k\rightarrow +\infty$. Pick points $x_k$ and times $t_k$
such that, as $k\rightarrow +\infty$,
$$
(T_k-t_k)|Rm(x_k,t_k)|\geq \gamma_k \sup_{x\in M, t\leq
T_k}(T_k-t)|Rm(x,t)|\rightarrow +\infty.
$$
Again denote by
$$
\epsilon_k=\frac{1}{\sqrt{|Rm(x_k,t_k)|}}
$$
and dilate the solution as before to get
$$
\tilde{g}_{ij}^{(k)}(\cdot,\tilde{t})
=\epsilon_k^{-2}g_{ij}(\cdot,t_k+\epsilon_k^2\tilde{t}),\ \ \ \
\tilde{t}\in [-t_k/\epsilon_k^2,(T_k-t_k)/\epsilon_k^2),
$$
which is still a solution to the Ricci flow and satisfies the
curvature bound
\begin{align*}
|\tilde{R}m^{(k)}|&  = \epsilon_k^2|Rm| \\
&  \leq\frac{1}{\gamma_k}\cdot \frac{(T_k-t_k)}{(T_k-t)}\\
& =\frac{1}{\gamma_k}\frac{(T_k-t_k)|Rm(x_k,t_k)|}{[(T_k-t_k)|
Rm(x_k,t_k)|-\tilde{t}]}\quad \mbox{ for } \; \tilde{t}\in\
\bigg[-\frac{t_k}{\epsilon_k^2},\frac{(T_k-t_k)}{\epsilon_k^2}\bigg),
\end{align*}
since $t=t_k+\epsilon_k^2\tilde{t}$ and
$\epsilon_k=1/\sqrt{|Rm(x_k,t_k)|}$. Hence $\{(x_k,t_k)\}$ is a
sequence of (almost) maximum points. And then as before, by
applying Hamilton's compactness theorem 4.1.5, there exists a
subsequence of the metrics $\tilde{g}_{ij}^{(k)}(\tilde{t})$ which
converges in the $C_{loc}^{\infty}$ topology to a limit
$\tilde{g}_{ij}^{(\infty)}(\tilde {t})$ on a limiting manifold
$\tilde{M}$ and $\tilde{t}\in(-\infty,+\infty)$ such that
$\tilde{g}_{ij}^{(\infty)}(\tilde {t})$ is a complete solution of
the Ricci flow and its curvature satisfies
$$
|\tilde{R}m^{(\infty)}|\leq 1
$$
everywhere on $\tilde{M}\times (-\infty,+\infty)$ and the equality
holds somewhere at $\tilde{t}=0$.

\medskip
{\bf Type II(b):}\index{type II!(b)}\ \ We consider a maximal
solution $g_{ij}(x,t)$ on $M\times [0,T)$ with
$$
T=+\infty\quad \text{and} \quad \limsup_{t\rightarrow
T}tK_{\max}(t)=+\infty.
$$

Again let $T_k\rightarrow T=+\infty$, and $\gamma_k\nearrow 1$, as
$k\rightarrow +\infty$. Pick $x_k$ and $t_k$ such that
$$
t_k(T_k-t_k)|Rm(x_k,t_k)|\geq \gamma_k \sup_{x\in M, t\leq
T_k}t(T_k-t)|Rm(x,t)|.
$$
Define
$$
\tilde{g}_{ij}^{(k)}(\cdot,\tilde{t})
=\epsilon_k^{-2}g_{ij}(\cdot,t_k+\epsilon_k^2\tilde{t}),\ \ \ \
\tilde{t}\in [-t_k/\epsilon_k^2,(T_k-t_k)/\epsilon_k^2),
$$
where $\epsilon_k=1/\sqrt{|Rm(x_k,t_k)|}$.\vskip 0.1cm \noindent
Since
\begin{align*}
t_k(T_k-t_k)|Rm(x_k,t_k)|&  \geq \gamma_k \sup\limits_{x\in M,
t\leq T_k}t(T_k-t)|Rm(x,t)|\\
&  \geq\gamma_k \sup\limits_{x\in M, t\leq T_k/2}t(T_k-t)|Rm(x,t)|\\
&  \geq\frac{T_k}{2}\gamma_k \sup\limits_{x\in M, t\leq
T_k/2}t|Rm(x,t)|,
\end{align*}
we have
$$
\frac{t_k}{\epsilon_k^2}=t_k|Rm(x_k,t_k)|\geq
\frac{\gamma_k}{2}\(\frac{T_k}{T_k-t_k}\)\sup_{x\in M, t\leq
T_k/2}t|Rm(x,t)| \to +\infty,
$$
and
$$
\frac{(T_k-t_k)}{\epsilon_k^2}=(T_k-t_k)|Rm(x_k,t_k)|\geq
\frac{\gamma_k}{2}\(\frac{T_k}{t_k}\)\sup_{x\in M, t\leq
T_k/2}t|Rm(x,t)| \to +\infty,
$$
as $k \rightarrow +\infty$. As before, we also have
$$
\frac{\partial}{\partial
\tilde{t}}\tilde{g}_{ij}^{(k)}(\cdot,\tilde{t})
=-2\tilde{R}_{ij}^{(k)}(\cdot,\tilde{t})
$$
and {\small
\begin{align*}
& |\tilde{R}m^{(k)}| \\
&  =\epsilon_k^2|Rm| \\
&  \leq \frac{1}{\gamma_k} \cdot \frac{t_k(T_k-t_k)}{t(T_k-t)}\\
&  =\frac{1}{\gamma_k}\cdot \frac{t_k(T_k-t_k)|Rm(x_k,t_k)|}
{(t_k+\epsilon_k^2\tilde{t})
[(T_k-t_k)-\epsilon_k^2\tilde{t}]\cdot |Rm(x_k,t_k)|}\\
&  =\frac{1}{\gamma_k}\cdot \frac{t_k(T_k-t_k)|Rm(x_k,t_k)|}
{(t_k+\epsilon_k^2\tilde{t})[(T_k-t_k)|Rm(x_k,t_k)|-\tilde{t}]}\\
&  =\frac{t_k(T_k-t_k)|Rm(x_k,t_k)|}{\gamma_k
(\!1\!+\!\tilde{t}/(t_k|Rm(x_k,t_k)|))[t_k(T_k\!-\!t_k)|Rm(x_k,t_k)|]
(1\!-\!\tilde{t}/((T_k\!-\!t_k)|Rm(x_k,t_k)|))}\\
&  \rightarrow 1,\quad \text{as } \; k\rightarrow+\infty.
\end{align*}
} Hence $\{(x_k,t_k)\}$ is again a sequence of (almost) maximum
points. As before, there exists a subsequence of the metrics
$\tilde{g}_{ij}^{(k)}(\tilde {t})$ which converges in the
$C_{loc}^{\infty}$ topology to a limit
$\tilde{g}_{ij}^{(\infty)}(\tilde {t})$ on a limiting manifold
$\tilde{M}$ and $\tilde{t}\in(-\infty,+\infty)$ such that
$\tilde{g}_{ij}^{(\infty)}(\tilde {t})$ is a complete solution of
the Ricci flow and its curvature satisfies
$$
|\tilde{R}m^{(\infty)}|\leq 1
$$
everywhere on $\tilde{M}\times (-\infty,+\infty)$ with the
equality somewhere at $\tilde{t}=0$.

\medskip
{\bf Type III(a):}\index{type III!(a)} \ We consider a maximal
solution $g_{ij}(x,t)$ on $M\times [0,T)$ with $T=+\infty$ and
$$
\limsup\limits_{t\rightarrow T}tK_{\max}(t)=A\in (0,+\infty).
$$
Choose a sequence of $x_k$ and $t_k$ such that
$t_k\rightarrow+\infty$ and
$$
\lim_{k\rightarrow \infty}t_k|Rm(x_k,t_k)|=A.
$$
Set $\epsilon_k=1/\sqrt{|Rm(x_k,t_k)|}$ and dilate the solution as
before to get
$$
\tilde{g}_{ij}^{(k)}(\cdot,\tilde{t})
=\epsilon_k^{-2}g_{ij}(\cdot,t_k+\epsilon_k^2\tilde{t}),\ \ \ \
\tilde{t}\in [-t_k/\epsilon_k^2,+\infty)
$$
which is still a solution to the Ricci flow. Also, for arbitrarily
fixed $\epsilon>0$, there exists a sufficiently large positive
constant $\tau$ such that for $t\in [\tau,+\infty)$,
\begin{align*}
|\tilde{R}m^{(k)}|&  = \epsilon_k^2|Rm| \\
&  \leq \epsilon_k^2\(\frac{A+\epsilon}{t}\)\\
&  = \epsilon_k^2 \(\frac{A+\epsilon}{t_k+\epsilon_k^2\tilde{t}}\)\\
&  = (A+\epsilon)/(t_k|Rm(x_k,t_k)|+\tilde{t}),\quad \text{for }\;
\tilde{t}\in [(\tau-t_k)/\epsilon_k^2,+\infty).
\end{align*}
Note that
$$
(A+\epsilon)/(t_k|Rm(x_k,t_k)|+\tilde{t})\rightarrow
(A+\epsilon)/(A+\tilde{t}), \quad \text{as } \; k\rightarrow
+\infty
$$
and
$$
(\tau-t_k)/\epsilon_k^2\rightarrow -A,\quad \text{as }\;
k\rightarrow +\infty.
$$
Hence $\{(x_k,t_k)\}$ is a sequence of (almost) maximum points.
And then as before, there exists a subsequence of the metrics
$\tilde{g}_{ij}^{(k)}(\tilde {t})$ which converges in the
$C_{loc}^{\infty}$ topology to a limit
$\tilde{g}_{ij}^{(\infty)}(\tilde {t})$ on a limiting manifold
$\tilde{M}$ and $\tilde{t}\in(-A,+\infty)$ such that
$\tilde{g}_{ij}^{(\infty)}(\tilde {t})$ is a complete solution of
the Ricci flow and its curvature satisfies
$$
|\tilde{R}m^{(\infty)}| \leq A/(A+\tilde{t})
$$
everywhere on $\tilde{M}\times (-A,+\infty)$ with  the equality
somewhere at $\tilde{t}=0$.
\end{pf}

In the case of manifolds with nonnegative curvature operator, or
K\"ahler metrics with nonnegative holomorphic bisectional curvature,
we can bound the Riemannian curvature by the scalar curvature $R$ up
to a constant factor depending only on the dimension. Then we can
slightly modify the statements in the previous theorem as follows


\begin{corollary} [{Hamilton \cite{Ha95F}}]
For any complete maximal solution to the Ricci flow satisfying the
injectivity radius condition and with bounded and nonnegative
curvature operator on a Riemannian manifold, or on a K\"ahler
manifold with bounded and nonnegative holomorphic bisectional
curvature, there exists a sequence of dilations of the solution
along $($almost$\,)$ maximum points which converges to a singular
model.

\vskip 0.1cm \noindent For Type {\rm I} solutions:\  the limit
model exists for $t\in (-\infty,\Omega)$ with $0<\Omega<+\infty$
and has
$$
R\leq \Omega/(\Omega-t)
$$
everywhere with equality somewhere at $t=0$.

\vskip 0.1cm \noindent For Type {\rm II} solutions:\ the limit
model exists for $t\in (-\infty, +\infty)$  and has
$$
R\leq 1
$$
everywhere with equality somewhere at $t=0$.

\vskip 0.1cm \noindent For Type {\rm III} solutions:\ the limit
model exists for $t\in(-A, +\infty)$ with $0<A<+\infty$ and has
$$
R\leq A/(A+t)
$$
everywhere with equality somewhere at $t=0$.
\end{corollary}

A natural and important question is to understand each of the three
types of singularity models. The following results obtained by
Hamilton \cite{Ha93E} and Chen-Zhu \cite{CZ00} characterize the Type
II and Type III singularity models with nonnegative curvature
operator and positive Ricci curvature respectively.

\begin{theorem} \
\begin{itemize}
\item[(i)] \ {\rm (Hamilton \cite{Ha93E})} \ Any Type {\rm II}
singularity model with nonnegative curvature operator and positive
Ricci curvature to the Ricci flow on a manifold $M$ must be a
$($steady$)$ Ricci soliton. \item[(ii)] \ {\rm (Chen-Zhu
\cite{CZ00})} \ Any Type {\rm III} singularity model with
nonnegative curvature operator and positive Ricci curvature on a
manifold $M$ must be a homothetically expanding Ricci soliton.
\end{itemize}
\end{theorem}

\begin{pf}
We only give the proof of (ii), since the proof of (i) is similar
and easier.

After a shift of the time variable, we may assume the Type III
singularity model is defined on $0<t<+\infty$ and $tR$ assumes its
maximum in space-time.

Recall from the Li-Yau-Hamilton inequality (Theorem 2.5.4) that
for any vectors $V^i$ and $W^i$, \be
M_{ij}W^iW^j+(P_{kij}+P_{kji})V^kW^iW^j+R_{ikjl}W^iW^jV^kV^l
\geq 0, 
\ee where
$$
M_{ij}=\Delta R_{ij}-\frac{1}{2}\nabla_i\nabla_jR+2R_{ipjq}R^{pq}
-g^{pq}R_{ip}R_{jq}+\frac{1}{2t}R_{ij}
$$
and
$$
P_{ijk}=\nabla_iR_{jk}-\nabla_jR_{ik}.
$$
Take the trace on $W$ to get \be
Q\stackrel{\Delta}{=}\frac{\partial R}{\partial
t}+\frac{R}{t}+2\nabla_iR\cdot V^i+2R_{ij}V^iV^j\geq0
\ee for any vector $V^i$. Let us choose $V$ to be the vector field
minimizing $Q$, i.e., \be
V^i=-\frac{1}{2}(\Ric^{-1})^{ik}\nabla_kR, 
\ee where $(Ric^{-1})^{ik}$ is the inverse of the Ricci tensor
$R_{ij}$. Substitute this vector field $V$ into $Q$ to get a smooth
function $\tilde{Q}$. The calculations of Chow-Hamilton in the proof
of Theorem 6.1 of \cite{CH} give, \be \frac{\partial}{\partial
t}\tilde{Q}
\ge \Delta \tilde{Q}-\frac{2}{t}\tilde{Q}. 
\ee

Suppose $tR$ assumes its maximum at $(x_0,t_0)$ with $t_0>0$, then
$$
\frac{\partial R}{\partial t}+\frac{R}{t}=0,\qquad \text{at}\quad
(x_0,t_0).
$$
This implies that the quantity
$$
Q=\frac{\partial R}{\partial t}+\frac{R}{t} +2\nabla_iR\cdot
V^i+2R_{ij}V^iV^j
$$
vanishes in the direction $V=0$ at $(x_0,t_0)$. We claim that for
any earlier time $t<t_0$ and any point $x\in M$, there is a vector
$V\in T_xM$ such that $Q=0$.

We argue by contradiction. Suppose not, then there is $\bar{x}\in
M$ and $0<\bar{t}<t_0$ such that $\tilde{Q}$ is positive at $x =
\bar{x}$ and $t=\bar{t}$. We can find a nonnegative smooth
function $\rho$ on $M$ with support in a neighborhood of $\bar{x}$
so that $\rho(\bar{x})>0$ and
$$
\tilde{Q}\geq \frac{\rho}{\bar{t}^2},
$$
at $t=\bar{t}$. Let $\rho$ evolve by the heat equation
$$
\frac{\partial\rho}{\partial t}=\Delta \rho.
$$
It then follows from the standard strong maximum principle that
$\rho>0$ everywhere for any $t>\bar{t}$. From (4.3.4) we see that
$$
\frac{\partial}{\partial t}\(\tilde{Q}-\frac{\rho}{t^2}\)
\ge\Delta\(\tilde{Q}-\frac{\rho}{t^2}\)-\frac{2}{t}\(\tilde{Q}
-\frac{\rho}{t^2}\)
$$
Then by the maximum principle as in Chapter 2, we get
$$
\tilde{Q}\ge \frac{\rho}{t^2}>0,\qquad \text{for all }\; t\ge
\bar{t}.
$$
This gives a contradiction with the fact $Q=0$ for $V=0$ at
$(x_0,t_0)$. We thus prove the claim.

Consider each time $t<t_0$. The null vector field of $Q$ satisfies
the equation \be
\nabla_iR+2R_{ij}V^j=0,
\ee by the first variation of $Q$ in $V$. Since $R_{ij}$ is
positive, we see that such a null vector field is unique and
varies smoothly in space-time.

Substituting (4.3.5) into the expression of $Q$, we have \be
\frac{\partial R}{\partial t}+\frac{R}{t}
+\nabla_iR\cdot V^i=0 
\ee Denote by
$$
Q_{ij}=M_{ij}+(P_{kij}+P_{kji})V^k+R_{ikjl}V^kV^l.
$$

{}From (4.3.1) we see that $Q_{ij}$ is nonnegative definite with
its trace $Q=0$ for such a null vector $V$. It follows that
$$
Q_{ij}=M_{ij}+(P_{kij}+P_{kji})V^k+R_{ikjl}V^kV^l=0.
$$
Again from the first variation of $Q_{ij}$ in $V$, we see that \be
(P_{kij}+P_{kji})+(R_{ikjl}+R_{jkil})V^l=0, 
\ee and hence \be
M_{ij}-R_{ikjl}V^kV^l=0. 
\ee

Applying the heat operator to (4.3.5) and (4.3.6) we get
\begin{align}
0&  = \(\frac{\partial}{\partial
t}-\Delta\)(\nabla_iR+2R_{ij}V^j)\\
&  = 2R_{ij}\(\frac{\partial}{\partial
t}-\Delta\)V^j+\(\frac{\partial}{\partial
t}-\Delta\)(\nabla_iR) \nn\\
&\quad+2V^j\(\frac{\partial}{\partial
t}-\Delta\)R_{ij}-4\nabla_kR_{ij}\nabla^kV^j, \nn
\end{align}  
and
\begin{align}
0&  = \(\frac{\partial}{\partial t}-\Delta\)\(\frac{\partial
R}{\partial t}+\frac{R}{t}+\nabla_iR\cdot V^i\)\\
&  = \nabla_iR\(\frac{\partial}{\partial
t}-\Delta\)V^i+V^i\(\frac{\partial}{\partial
t}-\Delta\)(\nabla_iR)-2\nabla_k\nabla_iR\cdot \nabla^kV^i \nn\\
&\quad +\(\frac{\partial}{\partial t}-\Delta\)\(\frac{\partial
R}{\partial t}+\frac{R}{t}\). \nn
\end{align}   

Multiplying (4.3.9) by $V^i$, summing over $i$ and adding
(4.3.10), as well as using the evolution equations on curvature,
we get
\begin{align}
0&  =
2V^i(2\nabla_i(|Rc|^2)-R_{il}\nabla^lR)+2V^iV^j(2R_{piqj}R^{pq}
-2g^{pq}R_{pi}R_{qj})\\
&\quad -4\nabla_kR_{ij}\cdot\nabla^kV^j\cdot V^i-2\nabla_k\nabla_i
R\cdot\nabla^kV^i+4R^{ij}\nabla_i\nabla_jR \nn\\
&\quad +4g^{kl}g^{mn}g^{pq}R_{km}R_{np}R_{ql}
+4R_{ijkl}R^{ik}V^jV^l-\frac{R}{t^2}. \nn
\end{align}  
{}From (4.3.5), we have the following equalities \be
\begin{cases}
-2V^iR_{il}\nabla^lR-4V^iV^jg^{pq}R_{pi}R_{qj}=0,\\
-4\nabla_kR_{ij}\cdot\nabla^kV^j\cdot V^i-2\nabla_k\nabla_i
R\cdot\nabla^kV^i=4R_{ij}\nabla_kV^i\cdot\nabla^kV^j,\\
\nabla_i\nabla_jR=-2\nabla_iR_{jl}\cdot V^l-2R_{jl}\nabla_iV^l.
\end{cases}   
\ee Substituting (4.3.12) into (4.3.11), we obtain
\begin{align*}
& 8R^{ij}(\nabla_kR_{ij}\cdot
V^k+R_{ikjl}V^kV^l-\nabla_iR_{jl}\cdot
V^l-R_{jl}\nabla_iV^l)\\
&\qquad+4R_{ij}\nabla_kV^i\cdot\nabla^kV^j
+4g^{kl}g^{mn}g^{pq}R_{km}R_{np}R_{ql}-\frac{R}{t^2}=0.
\end{align*}
By using (4.3.7), we know
$$
R^{ij}(\nabla_kR_{ij}\cdot V^k +R_{ikjl}V^kV^l-\nabla_iR_{jl}\cdot
V^l)=0.
$$
Then we have \be
-8R^{ij}R_{jl}\nabla_iV^l+4R_{ij}\nabla_kV^i\cdot\nabla^kV^j
+4g^{kl}g^{mn}g^{pq}R_{km}R_{np}R_{ql}-\frac{R}{t^2}=0.
\ee By taking the trace in the last equality in (4.3.12) and using
(4.3.6) and the evolution equation of the scalar curvature, we can
get \be
R^{ij}(R_{ij}+\frac{g_{ij}}{2t}-\nabla_iV_j)=0. 
\ee Finally by combining (4.3.13) and (4.3.14), we deduce
$$
4R^{ij}g^{kl}\(R_{ik}+\frac{g_{ik}}{2t}-\nabla_kV_i\)
\(R_{jk}+\frac{g_{jk}}{2t}-\nabla_kV_j\)=0.
$$
Since $R_{ij}$ is positive definite, we get \be
\nabla_iV_j=R_{ij}+\frac{g_{ij}}{2t},\qquad
\mbox{for all }\;i,j.  
\ee This means that $g_{ij}(t)$ is a homothetically expanding
Ricci soliton.
\end{pf}

\begin{remark}
Recall from Section 1.5 that any compact steady Ricci soliton or
expanding Ricci soliton must be Einstein. If the manifold $M$ in
Theorem 4.3.6 is noncompact and simply connected, then the steady
(or expanding) Ricci soliton must be a steady (or expanding)
gradient Ricci soliton. For example, we know that $\nabla_iV_j$ is
symmetric from (4.3.15).  Also, by the simply connectedness of $M$
there exists a function $F$ such that
$$
\nabla_i\nabla_jF=\nabla_iV_j,\qquad \text{on }\;M.
$$
So
$$
R_{ij}=\nabla_i\nabla_jF-\frac{g_{ij}}{2t},\qquad \text{on }\;M
$$
This means that $g_{ij}$ is an expanding gradient Ricci soliton.
\end{remark}

In the K\"ahler case, we have the following results for Type II
and Type III singularity models with nonnegative holomorphic
bisectional curvature obtained by the first author in
\cite{Cao97}.

\begin{theorem}[{Cao \cite{Cao97}}] \
\begin{itemize}
\item[(i)] Any Type {\rm II} singularity model on a K\"ahler
manifold with nonnegative holomorphic bisectional curvature and
positive Ricci curvature must be a steady K\"ahler-Ricci soliton.
\item[(ii)] Any Type {\rm III} singularity model on a K\"ahler
manifold with nonnegative holomorphic bisectional curvature and
positive Ricci curvature must be an expanding K\"ahler-Ricci
soliton.
\end{itemize}
\end{theorem}

To conclude this section, we state a result of Sesum \cite{Se} on
compact Type I singularity models. Recall that Perelman's
functional $\mathcal{W}$, introduced in Section 1.5, is given by
$$
\mathcal{W}(g,f,\tau)=\int_M(4\pi\tau)^{-\frac{n}{2}}[\tau(|\nabla
f|^2+R)+f-n]e^{-f}dV_g
$$
with the function $f$ satisfying the constraint
$$
\int_M(4\pi\tau)^{-\frac{n}{2}}e^{-f}dV_g=1.
$$
And recall from Corollary 1.5.9 that
$$
\mu(g(t))=\inf\left\{\mathcal{W}(g(t),f,T-t)|
\int_M(4\pi(T-t))^{-\frac{n}{2}}e^{-f}dV_{g(t)}=1\right\}
$$
is strictly increasing along the Ricci flow unless we are on a
gradient shrinking soliton. If one can show that $\mu(g(t))$ is
uniformly bounded from above and the minimizing functions
$f=f(\cdot,t)$ have a limit as $t\rightarrow T$, then the
rescaling limit model will be a shrinking gradient soliton. As
shown by Natasa Sesum in \cite{Se}, Type I assumption guarantees
the boundedness of $\mu(g(t))$, while the compactness assumption
of the rescaling limit guarantees the existence of the limit for
the minimizing functions $f(\cdot,t)$. Therefore we have

\begin{theorem}[{Sesum \cite{Se}}]
Let $(M,g_{ij}(t))$ be a Type {\rm I} singularity model obtained
as a rescaling limit of a Type {\rm I} maximal solution. Suppose
$M$ is compact. Then $(M,g_{ij}(t))$ must be a gradient shrinking
Ricci soliton.
\end{theorem}

It seems that the assumption on the compactness of the rescaling
limit is superfluous. We conjecture that any noncompact Type I
limit is also a gradient shrinking soliton.

\section{Ricci Solitons}

We now follow Hamilton \cite{Ha95F} to examine the structure of a
steady Ricci soliton that we get as a Type II limit . The material
is from section 20 of Hamilton \cite{Ha95F} and Hamilton
\cite{Ha88}.

\begin{lemma} [{Hamilton \cite{Ha95F}}]
Suppose we have a complete gradient steady Ricci soliton $g_{ij}$
with bounded curvature so that
$$
R_{ij}=\nabla_i\nabla_jF
$$
for some function $F$ on $M$. Assume the Ricci curvature is
positive and the scalar curvature $R$ attains its maximum
$R_{\max}$ at a point $x_0\in M$. Then \be
|\nabla F|^2+R=R_{\max} 
\ee everywhere on $M$, and furthermore $F$ is convex and attains
its minimum at $x_0$.
\end{lemma}

\begin{pf}
Recall that, from (1.1.15) and noting our $F$ here is $-f$ there,
the steady gradient Ricci soliton has the property
$$
|\nabla F|^2+R=C_0
$$
for some constant $C_0$. Clearly, $C_0\ge R_{\max}$.

If $C_0=R_{\max}$, then $\nabla F=0$ at the point $x_0$. Since
$\nabla_i\nabla_jF=R_{ij}>0,$ we see that $F$ is convex and $F$
attains its minimum at $x_0$.

If $C_0>R_{\max}$, consider a gradient path of $F$ in a local
coordinate neighborhood through $x_0=(x_0^1,\ldots,x_0^n):$
$$
\left\{\arraycolsep=1.5pt\begin{array}{l} x^i=x^i(u),\qquad
u\in(-\varepsilon,\varepsilon),\;i=1,\ldots,n\\[2mm]
x_0^i=x^i(0),
\end{array}\right.
$$
and
$$
\frac{dx^i}{du} =g^{ij}\nabla_jF,\quad
u\in(-\varepsilon,\varepsilon).
$$
Now $|\nabla F|^2=C_0-R\ge C_0-R_{\max}>0$ everywhere, while
$|\nabla F|^2$ is smallest at $x=x_0$ since $R$ is largest there.
But we compute
\begin{align*}
\frac{d}{du}|\nabla
F|^2&  = 2g^{jl}\(\frac{d}{du}\nabla_jF\)\nabla_lF\\
&  = 2g^{ik}g^{jl}\nabla_i\nabla_jF\cdot\nabla_kF\nabla_lF\\
&  = 2g^{ik}g^{jl}R_{ij}\nabla_kF\nabla_lF\\
&  > 0
\end{align*}
since $R_{ij}>0$ and $|\nabla F|^2>0$. Then $|\nabla F|^2$ is not
smallest at $x_0$, and we have a contradiction. 
\end{pf}

We remark that when we are considering a complete expanding
gradient Ricci soliton on $M$ with positive Ricci curvature and
$$
R_{ij}+\rho g_{ij}=\nabla_i\nabla_jF
$$
for some constant $\rho>0$ and some function $F$, the above
argument gives
$$
|\nabla F|^2+R-2\rho F=C
$$
for some positive constant $C$. Moreover the function $F$ is an
exhausting and convex function. In particular, such an expanding
gradient Ricci soliton is diffeomorphic to the Euclidean space
$\mathbb R^n$.

Let us introduce a geometric invariant as follows. Let $O$ be a
fixed point in a Riemannian manifold $M$, $s$ the distance to the
fixed point $O$, and $R$ the scalar curvature. We define the {\bf
asymptotic scalar curvature ratio}\index{asymptotic scalar
curvature ratio}
$$
A=\limsup_{s\rightarrow +\infty}Rs^2.
$$
Clearly the definition is independent of the choice of the fixed
point $O$ and invariant under dilation. This concept is particular
useful on manifolds with positive sectional curvature. The first
type of gap theorem was obtained by Mok-Siu-Yau \cite{MoSY} in
understanding the hypothesis of the paper of Siu-Yau \cite{SiY}.
Yau (see \cite{GW82}) suggested that this should be a general
phenomenon. This was later conformed by Greene-Wu \cite{GW82, GW},
Eschenberg-Shrader-Strake \cite{ESS} and Drees \cite{Dr} where
they show that any complete noncompact $n$-dimensional (except
$n=4$ or 8) Riemannian manifold of positive sectional curvature
must have $A>0$. Similar results on complete noncompact
K$\rm\ddot{a}$hler manifolds of positive holomorphic bisectional
curvature were obtained by Chen-Zhu \cite{CZ03} and Ni-Tam
\cite{NT}.

\begin{theorem}[{Hamilton \cite{Ha95F}}]
For a complete noncompact steady gradient Ricci soliton with
bounded curvature and positive sectional curvature of dimension
$n\ge3$ where the scalar curvature assume its maximum at a point
$O\in M$, the asymptotic scalar curvature ratio is infinite, i.e.,
$$
A=\limsup_{s\rightarrow +\infty}Rs^2=+\infty
$$
where s is the distance to the point $O$.
\end{theorem}

\begin{pf}
The solution to the Ricci flow corresponding to the soliton exists
for $-\infty<t<+\infty$ and is obtained by flowing along the
gradient of a potential function $F$ of the soliton. We argue by
contradiction. Suppose $Rs^2\le C$. We will show that the limit
$$
\bar{g}_{ij}(x)=\lim_{t\rightarrow-\infty}g_{ij}(x,t)
$$
exists for $x\neq O$ on the manifold $M$ and is a complete flat
metric on $M\setminus\{O\}.$ Since the sectional curvature of $M$
is positive everywhere, it follows from Cheeger-Gromoll \cite{CG}
that $M$ is diffeomorphic to $\mathbb{R}^n$. Thus
$M\setminus\{O\}$ is diffeomorphic to
$\mathbb{S}^{n-1}\times\mathbb{R}$. But for $n\ge 3$ there is no
flat metric on $\mathbb{S}^{n-1}\times\mathbb{R}$, and this will
finish the proof.

To see the limit metric exists, we note that $R\rightarrow 0$ as
$s\rightarrow+\infty$, so $|\nabla F|^2\rightarrow R_{\max}$ as
$s\rightarrow +\infty$ by (4.4.1). The function $F$ itself can be
taken to evolve with time, using the definition
$$
\frac{\partial F}{\partial t} =\nabla_iF\cdot\frac{\partial
x^i}{\partial t} =-|\nabla F|^2 =\Delta F-R_{\max}
$$
which pulls $F$ back by the flow along the gradient of $F$. Then
we continue to have $\nabla_i\nabla_jF=R_{ij}$ for all time, and
$|\nabla F|^2\rightarrow R_{\max}$ as $s\rightarrow +\infty$ for
each time.

When we go backward in time, this is equivalent to flowing
outwards along the gradient of $F$, and our speed approaches
$\sqrt{R_{\max}}$. So, starting outside of any neighborhood of $O$
we have
$$
\frac{s}{|t|}=\frac{d_t(\cdot,O)}{|t|}\rightarrow
\sqrt{R_{\max}},\quad \text{as }\;t\rightarrow-\infty
$$
and \be R(\cdot,t)\le\frac{C}{R_{\max}\cdot|t|^2},\quad
\text{as}\;|t|\;\text{large enough}.  
\ee Hence for $|t|$ sufficiently large,
\begin{align*}
0&  \ge -2R_{ij}\\
&  = \frac{\partial}{\partial t}g_{ij}\\
&  \ge -2Rg_{ij}\\
&  \ge -\frac{2C}{R_{\max}\cdot|t|^2}g_{ij}
\end{align*}
which implies that for any tangent vector $V$,
$$
0\le \frac{d}{d|t|}(\log(g_{ij}(t)V^iV^j))
\le\frac{2C}{R_{\max}\cdot|t|^2}.
$$
These two inequalities show that $g_{ij}(t)V^iV^j$ has a limit
$\bar{g}_{ij}V^iV^j$ as $t\rightarrow-\infty$.

Since the metrics are all essentially the same, it always takes an
infinite length to get out to the infinity. This shows the limit
$\bar{g}_{ij}$ is complete at the infinity. One the other hand,
any point $P$ other than $O$ will eventually be arbitrarily far
from $O$, so the limit metric $\bar{g}_{ij}$ is also complete away
from $O$ in $M\setminus\{O\}$. Using Shi's derivative estimates in
Chapter 1, it follows that $g_{ij}(\cdot,t)$ converges in the
$C_{\rm loc}^\infty$ topology to a complete smooth limit metric
$\bar{g}_{ij}$ as $t\rightarrow-\infty$, and the limit metric is
flat by (4.4.2).
\end{pf}

The above argument actually shows that \be
\limsup_{s\rightarrow+\infty}Rs^{1+\varepsilon}=+\infty
\ee for arbitrarily small $\varepsilon>0$ and for any complete
gradient Ricci soliton with bounded and positive sectional
curvature of dimension $n\ge3$ where the scalar curvature assumes
its maximum at a fixed point $O$.

Finally we conclude this section with the important result of
Hamilton on the uniqueness of complete Ricci soliton on
two-dimensional Riemannian manifolds.

\begin{theorem}[cf. Theorem 10.1 of {Hamilton \cite{Ha88}}]
The only complete steady Ricci soliton on a two-dimensional
manifold with bounded curvature which assumes its maximum $1$ at
an origin is the ``cigar" soliton on the plane $\mathbb{R}^2$ with
the metric
$$
ds^2=\frac{dx^2+dy^2}{1+x^2+y^2}.
$$
\end{theorem}

\begin{pf}
Recall that the scalar curvature evolves by
$$
\frac{\partial R}{\partial t}=\Delta R+R^2
$$
on a two-dimensional manifold $M$. Denote by
$R_{\min}(t)=\inf\{R(x,t)\ |\ x\in M\}.$ We see from the maximum
principle (see for example Chapter 2) that $R_{\min}(t)$ is
strictly increasing whenever $R_{\min}(t)\neq 0$, for
$-\infty<t<+\infty$. This shows that the curvature of a steady
Ricci soliton on a two-dimensional manifold $M$ must be
nonnegative and $R_{\min}(t) = 0$ for all $t \in
(-\infty,+\infty)$. Further by the strong maximum principle we see
that the curvature is actually positive everywhere. In particular,
the manifold must be noncompact. So the manifold $M$ is
diffiomorphic to $\mathbb{R}^2$ and the Ricci soliton must be a
gradient soliton. Let $F$ be a potential function of the gradient
Ricci soliton. Then, by definition, we have
$$
\nabla_iV_j+\nabla_jV_i=Rg_{ij}
$$
with $V_i=\nabla_iF$. This says that the vector field $V$ must be
conformal. In complex coordinate a conformal vector field is
holomorphic. Hence $V$ is locally given by
$V(z)\frac{\partial}{\partial z}$ for a holomorphic function
$V(z)$. At a zero of $V$ there will be a power series expansion
$$
V(z)=az^p+\cdots,\;(a\neq 0)
$$
and if $p>1$ the vector field will have closed orbits in any
neighborhood of the zero. Now the vector field is gradient and a
gradient flow cannot have a closed orbit. Hence $V(z)$ has only
simple zeros. By Lemma 4.4.1, we know that $F$ is strictly convex
with the only critical point being the minima, chosen to be the
origin of $\mathbb{R}^2$. So the holomorphic vector field $V$ must
be
$$
V(z)\frac{\partial}{\partial z} =cz\frac{\partial}{\partial
z},\qquad \text{for}\;z\in C,
$$
for some complex number $c$.

We now claim that $c$ is real. Let us write the metric as
$$
ds^2=g(x,y)(dx^2+dy^2)
$$
with $z=x+\sqrt{-1}y.$ Then $\nabla F=cz\frac{\partial}{\partial
z}$ means that if $c=a+\sqrt{-1}b$, then
$$
\frac{\partial F}{\partial x}=(ax-by)g,\quad \frac{\partial
F}{\partial y}=(bx+ay)g.
$$
Taking the mixed partial derivatives $\frac{\partial^2 F}{\partial
x\partial y}$ and equating them at the origin $x=y=0$ gives $b=0$,
so $c$ is real.

Let
$$\left\{\arraycolsep=1.5pt\begin{array}{l} x=e^u\cos
v,\qquad -\infty<u<+\infty,\\[2mm]
y=e^u\sin v,\qquad 0\le v\le 2\pi.
\end{array}\right.
$$
Write
\begin{align*}
ds^2& = g(x,y)(dx^2+dy^2)\\
&  = g(e^u\cos v,e^u\sin v)e^{2u}(du^2+dv^2)\\
&  \stackrel{\Delta}{=} g(u,v)(du^2+dv^2).
\end{align*}
Then we get the equations
$$
\frac{\partial F}{\partial u}=ag,\quad \frac{\partial F}{\partial
v}=0
$$
since the gradient of $F$ is just $a\frac{\partial}{\partial u}$
for a real constant $a$. The second equation shows that $F=F(u)$
is a function of $u$ only, then the first equation shows that
$g=g(u)$ is also a function of $u$ only. Then we can write the
metric as
\begin{align}
ds^2&  = g(u)(du^2+dv^2)\\
&  = g(u)e^{-2u}(dx^2+dy^2).\nn
\end{align}
This implies that $e^{-2u}g(u)$ must be a smooth function of
$x^2+y^2=e^{2u}$. So as $u\rightarrow -\infty$, \be
g(u)=b_1e^{2u}+b_2(e^{2u})^2+\cdots,
\ee with $b_1>0$.

The curvature of the metric is given by
$$
R=-\frac{1}{g}\(\frac{g'}{g}\)'
$$
where $(\cdot)'$ is the derivative with respect to $u$. Note that
the soliton is by translation in $u$ with velocity $c$. Hence
$g=g(u+ct)$ satisfies
$$
\frac{\partial g}{\partial t}=-Rg
$$
which becomes
$$
cg'=\(\frac{g'}{g}\)'.
$$
Thus by (4.4.5),
$$
\frac{g'}{g}=cg+2
$$
and then by integrating
$$
e^{2u}\(\frac{1}{g}\)=-\frac{c}{2}e^{2u}+b_1
$$
i.e.,
$$
g(u)=\frac{e^{2u}}{b_1-\frac{c}{2}e^{2u}}.
$$
In particular, we have $c<0$ since the Ricci soliton is not flat.
Therefore
$$
ds^2=g(u)e^{-2u}(dx^2+dy^2)
=\frac{dx^2+dy^2}{\alpha_1+\alpha_2(x^2+y^2)}
$$
for some constants $\alpha_1,\alpha_2>0$. By the normalization
condition that the curvature attains its maximum $1$ at the
origin, we conclude that
$$
ds^2=\frac{dx^2+dy^2}{1+(x^2+y^2)}.
$$
\end{pf}

\newpage
\part{{\Large Long Time Behaviors}}

\bigskip
Let $M$ be a complete manifold of dimension $n$. Consider a
solution of the Ricci flow $g_{ij}(x,t)$ on $M$ and on a maximal
time interval $[0,T)$. When $M$ is compact, we usually consider
the \textbf{normalized Ricci flow} \index{normalized Ricci flow}
$$
\frac{\partial g_{ij}}{\partial t}=\frac{2}{n}rg_{ij}-2R_{ij},
$$
where $r=\int_MRdV/\int_MdV$ is the average scalar curvature. The
factor $r$ serves to normalize the Ricci flow so that the volume
is constant. To see this we observe that $dV=\sqrt{\det
g_{ij}}\;dx$ and then
$$
\frac{\partial}{\partial t}\log\sqrt{\det g_{ij}}
=\frac{1}{2}g^{ij}\frac{\partial}{\partial t}g_{ij}=r-R,
$$
$$
\frac{d}{dt}\int_MdV=\int_M(r-R)dV=0.
$$
As observed by Hamilton \cite{Ha82}, the Ricci flow and the
normalized Ricci flow differ only by a change of scale in space
and a change of parametrization in time. Indeed, we first assume
that $g_{ij}(t)$ evolves by the (unnormalized) Ricci flow and
choose the normalization factor $\psi=\psi(t)$ so that
$\tilde{g}_{ij}=\psi g_{ij}$, and $\int_Md\tilde{\mu}=1$. Next we
choose a new time scale $\tilde{t}=\int\psi(t)dt$. Then for the
normalized metric $\tilde{g}_{ij}$ we have
$$
\tilde{R}_{ij}=R_{ij}, \tilde{R}=\frac{1}{\psi}R,
\tilde{r}=\frac{1}{\psi}r.
$$
Because $\int_Md\tilde{V}=1$, we see that
$\int_MdV=\psi^{-\frac{n}{2}}$. Then
\begin{align*}
\frac{d}{dt}\log\psi&  =\(-\frac{2}{n}\)\frac{d}{dt}\log\int_MdV\\
&  = \(-\frac{2}{n}\)\frac{\int_M\frac{\partial}{\partial
t}\sqrt{\det g_{ij}}\;dx}{\int_MdV}\\
&  = \frac{2}{n}r,
\end{align*}
since $\frac{\partial}{\partial t}g_{ij}=-2R_{ij}$ for the Ricci
flow. Hence it follows that
\begin{align*}
\frac{\partial}{\partial \tilde{t}}\tilde{g}_{ij}& =
\frac{\partial}{\partial t}g_{ij}+\(\frac{d}{dt}\log\psi\)g_{ij}\\
&  = \frac{2}{n}\tilde{r}\tilde{g}_{ij}-2\tilde{R}_{ij}.
\end{align*}
Thus studying the behavior of the Ricci flow near the maximal time
is equivalent to studying the long-time behavior of the normalized
Ricci flow.

In this chapter we will discuss long-time behavior of the normalized
Ricci flow for the following special cases: (1) compact
two-manifolds (cf. Hamilton \cite{Ha88} and Chow \cite{Chow}); (2)
compact three-manifolds with nonnegative Ricci curvature (cf.
Hamilton \cite{Ha82}); (3) compact four-manifolds with nonnegative
curvature operator (cf. Hamilton \cite{Ha86}); and (4) compact
three-manifolds with uniformly bounded normalized curvature (cf.
Hamilton \cite{Ha99}).

\section{The Ricci Flow on Two-manifolds}

Let $M$ be a compact surface, we will discuss in this section the
evolution of a Riemannian metric $g_{ij}$ under the normalized
Ricci flow. Most of the presentation in this section follows
Hamilton \cite{Ha88}, as well as Chow \cite{Chow}. We also refer
the reader to chapter 5 of Chow-Knopf's book \cite{CK04} for an
excellent description of the subject.

On a surface, the Ricci curvature is given by
$$
R_{ij}=\frac{1}{2}Rg_{ij}
$$
so the normalized Ricci flow equation becomes \be
\frac{\partial}{\partial t}g_{ij}=(r-R)g_{ij}. 
\ee Recall the Gauss-Bonnet formula says
$$
\int_MRdV=4\pi\chi(M),
$$
where $\chi(M)$ is the Euler characteristic number of $M$. Thus
the average scalar curvature $r=4\pi\chi(M)/\int_MdV$ is constant
in time.

To obtain the evolution equation of the normalized curvature, we
recall a simple principle in \cite {Ha82} for converting from the
unnormalized to the normalized evolution equation on an
$n$-dimensional manifold. Let $P$ and $Q$ be two expressions
formed from the metric and curvature tensors, and let $\tilde{P}$
and $\tilde{Q}$ be the corresponding expressions for the
normalized Ricci flow. Since they differ by dilations, they differ
by a power of the normalized factor $\psi=\psi(t)$. We say $P$ has
\textbf{degree}\index{degree} $k$ if $\tilde{P}=\psi^kP$. Thus
$g_{ij}$ has degree 1, $R_{ij}$ has degree 0, $R$ has degree $-1$.


\begin{lemma} [{Hamilton \cite{Ha82}}]
Suppose $P$ satisfies
$$
\frac{\partial P}{\partial t}=\Delta P+Q
$$
for the unnormalized Ricci flow, and $P$ has degree $k$. Then $Q$
has degree $k-1$, and for the normalized Ricci flow,
$$
\frac{\partial \tilde{P}}{\partial \tilde{t}}=\tilde{\Delta}
\tilde{P}+\tilde{Q}+\frac{2}{n}k\tilde{r}\tilde{P}.
$$
\end{lemma}

\begin{pf}
We first see $Q$ has degree $k-1$ since $\partial
\tilde{t}/\partial{t}=\psi$ and $\Delta=\psi\tilde{\Delta}$. Then
$$
\psi\frac{\partial}{\partial \tilde{t}}(\psi^{-k}\tilde{P})
=\psi\tilde{\Delta}(\psi^{-k}\tilde{P})+\psi^{-k+1}\tilde{Q}
$$
which implies
\begin{align*}
\frac{\partial \tilde{P}}{\partial \tilde{t}}&  =\tilde{\Delta}
\tilde{P}+\tilde{Q}+\frac{k}{\psi}\frac{\partial \psi}{\partial
\tilde{t}}\tilde{P}\\
&  = \tilde{\Delta}
\tilde{P}+\tilde{Q}+\frac{2}{n}k\tilde{r}\tilde{P}
\end{align*}
since $\frac{\partial}{\partial\tilde{t}}\log\psi
=(\frac{\partial}{\partial t}\log\psi)\psi^{-1}
=\frac{2}{n}\tilde{r}.$
\end{pf}

We now come back to the normalized Ricci flow (5.1.1) on a compact
surface. By applying the above lemma to the evolution equation of
unnormalized scalar curvature, we have \be
\frac{\partial R}{\partial t}=\Delta R+R^2-rR 
\ee for the normalized scalar curvature $R$. As a direct
consequence, by using the maximum principle, both nonnegative
scalar curvature and nonpositive scalar curvature are preserved
for the normalized Ricci flow on surfaces.

Let us introduce a potential function $\varphi$ as in the
K\"ahler-Ricci flow (see for example \cite{Cao85}). Since $R-r$
has mean value zero on a compact surface, there exists a unique
function $\varphi$, with mean value zero, such that \be
\Delta \varphi=R-r.     
\ee Differentiating (5.1.3) in time, we have
\begin{align*}
\frac{\partial }{\partial t}R
&  =\frac{\partial }{\partial t}(\Delta\varphi)\\
&  =(R-r)\Delta\varphi+g^{ij}\frac{\partial}{\partial
t}\(\frac{\partial^2\varphi}{\partial x^i\partial x^j}
-\Gamma^k_{ij}\frac{\partial \varphi}{\partial x^k}\)\\[4mm]
&  =(R-r)\Delta \varphi +\Delta\(\frac{\partial \varphi}{\partial
t}\).
\end{align*}
Combining with the equation (5.1.2), we get
$$
\Delta\(\frac{\partial \varphi}{\partial t}\)
=\Delta(\Delta\varphi)+r\Delta\varphi
$$
which implies that \be \frac{\partial \varphi}{\partial
t}=\Delta\varphi+r\varphi-b(t)
\ee for some function $b(t)$ of time only. Since $\int_M\varphi
dV=0$ for all $t$, we have
\begin{align*}
0=\frac{d}{dt}\int_M\varphi d\mu&  =\int_M(\Delta
\varphi+r\varphi-b(t))d\mu+\int_M\varphi(r-R)d\mu \\
&  = -b(t)\int_Md\mu+\int_M|\nabla \varphi|^2d\mu.
\end{align*}
Thus the function $b(t)$ is given by
$$
b(t)=\frac{\int_M|\nabla\varphi|^2d\mu}{\int_Md\mu}.
$$

Define a function $h$ by
$$
h=\Delta\varphi+|\nabla\varphi|^2=(R-r)+|\nabla\varphi|^2,
$$
and set
$$
M_{ij}=\nabla_i\nabla_{j}\; \varphi-\frac{1}{2}\Delta\;\varphi
g_{ij}
$$
to be the traceless part of $\nabla_i\nabla_j\; \varphi$.

\begin{lemma} [{Hamilton \cite{Ha88}}]
The function $h$ satisfies the evolution equation \be
\frac{\partial h}{\partial t}=\Delta h-2|M_{ij}|^2+rh.
\ee
\end{lemma}

\begin{pf}
Under the normalized Ricci flow,
\begin{align*}
\frac{\partial}{\partial t}|\nabla \varphi|^2& =
\(\frac{\partial}{\partial t}g^{ij}\)
\nabla_i\varphi\nabla_j\varphi+2g^{ij}\(\frac{\partial}{\partial
t}\nabla_i\varphi\)(\nabla_j\varphi)\\
&  =(R-r)|\nabla\varphi|^2+2g^{ij}\nabla_i
(\Delta\varphi+r\varphi-b(t))\nabla_j\varphi\\
&  =(R+r)|\nabla\varphi|^2+2g^{ij}
(\Delta\nabla_i\varphi-R_{ik}\nabla_k\varphi)\nabla_j\varphi\\
&  =(R+r)|\nabla\varphi|^2+\Delta|\nabla\varphi|^2
-2|\nabla^2\varphi|^2-2g^{ij}R_{ik}\nabla_k\varphi\nabla_j\varphi\\
&
=\Delta|\nabla\varphi|^2-2|\nabla^2\varphi|^2+r|\nabla\varphi|^2,
\end{align*}
where $R_{ik}=\frac{1}{2}Rg_{ik}$ on a surface.

On the other hand we may rewrite the evolution equation (5.1.2) as
$$
\frac{\partial}{\partial t}(R-r)=\Delta(R-r)+(\Delta
\varphi)^2+r(R-r).
$$
Then the combination of above two equations yields
\begin{align*}
\frac{\partial}{\partial t}h&  =\Delta
h-2(|\nabla^2\varphi|^2-\frac{1}{2}(\Delta\varphi)^2)+rh\\
&  =\Delta h-2|M_{ij}|^2+rh
\end{align*}
as desired.
\end{pf}

As a direct consequence of the evolution equation (5.1.5) and the
maximum principle, we have \be
R\leq C_1e^{rt}+r 
\ee for some positive constant $C_1$ depending only on the initial
metric.

On the other hand, it follows from (5.1.2) that
$R_{\min}(t)=\min_{x\in M}R(x,t)$ satisfies
$$
\frac{d}{dt}R_{\min}\geq R_{\min}(R_{\min}-r)\geq 0
$$
whenever $R_{\min}\leq 0$. This says that \be
R_{\min}(t)\geq-C_2, \quad\text{ for all }\; t>0 
\ee for some positive constant $C_2$ depending only on the initial
metric.

Thus the combination of (5.1.6) and (5.1.7) implies the following
long time existence result.

\begin{proposition} [{Hamilton \cite{Ha88}}]
For any initial metric on a compact surface, the normalized Ricci
flow $(5.1.1)$ has a solution for all time.
\end{proposition}

To investigate the long-time behavior of the solution, let us now
divide the discussion into three cases: $\chi(M)<0$; $\chi(M)=0$;
and $\chi(M)>0$.

\medskip
{\it Case} (1): $\chi(M)<0$ (i.e., $r<0$).

\smallskip
{}From the evolution equation (5.1.2), we have
\begin{align*}
\frac{d}{dt}R_{\min}&  \geq R_{\min}(R_{\min}-r) \\
&  \geq r(R_{\min}-r), \; \text{ on } \; M\times [0,+\infty)
\end{align*}
which implies that
$$
R-r\geq-\tilde{C_1}e^{rt}, \; \text{ on } \; M\times[0,+\infty)
$$
for some positive constant $\tilde{C_1}$ depending only on the
initial metric. Thus by combining with (5.1.6) we have \be
-\tilde{C_1}e^{rt}\leq R-r \leq C_1e^{rt}, \; \text{ on } \;
M\times
[0,+\infty). 
\ee

\begin{theorem}[{Hamilton \cite{Ha88}}]
On a compact surface with $\chi(M)<0$, for any initial metric the
solution of the normalized Ricci flow $(5.1.1)$ exists for all
time and converges in the $C^{\infty}$ topology to a metric with
negative constant curvature.
\end{theorem}

\begin{pf}
The estimate (5.1.8) shows that the scalar curvature $R(x,t)$
converges exponentially to the negative constant $r$ as
$t\rightarrow+\infty$.

Fix a tangent vector $v\in T_xM$ at a point $x\in M$ and let
$|v|^2_t=g_{ij}(x,t)v^iv^j$. Then we have
\begin{align*}
\frac{d}{dt}|v|_t^2
&  =\(\frac{\partial}{\partial t}g_{ij}(x,t)\)v^iv^j\\
&  = (r-R)|v|_t^2
\end{align*}
which implies
$$
\left|\frac{d}{dt}\log|v|_t^2\right| \leq Ce^{rt}, \; \text{ for
all } \; t>0
$$
for some positive constant $C$ depending only on the initial
metric (by using (5.1.8)). Thus $|v|_t^2$ converges uniformly to a
continuous function $|v|_{\infty}^2$ as $t\rightarrow +\infty$ and
$|v|_{\infty}^2\neq 0$ if $v\neq 0$. Since the parallelogram law
continues to hold to the limit, the limiting norm $|v|_{\infty}^2$
comes from an inner product $g_{ij}(\infty)$. This says, the
metrics $g_{ij}(t)$ are all equivalent and as $t\rightarrow
+\infty$, the metric $g_{ij}(t)$ converges uniformly to a
positive-definite metric tensor $g_{ij}(\infty)$ which is
continuous and equivalent to the initial metric.

By the virtue of Shi's derivative estimates of the unnormalized
Ricci flow in Section 1.4, we see that all derivatives and higher
order derivatives of the curvature of the solution $g_{ij}$ of the
normalized flow are uniformly bounded on $M\times [0,+\infty)$.
This shows that the limiting metric $g_{ij}(\infty)$ is a smooth
metric with negative constant curvature and the solution
$g_{ij}(t)$ converges to the limiting metric $g_{ij}(\infty)$ in
the $C^{\infty}$ topology as $t\rightarrow +\infty$.
\end{pf}

\smallskip
{\it Case} (2): $\chi(M)=0$, i.e., $r=0$.

\smallskip
{}The following argument of dealing with the case $\chi(M)=0$ is
adapted from Chow-Knopf's book (cf. section 5.6 of \cite{CK04}).
From (5.1.6) and (5.1.7) we know that the curvature remains
bounded above and below. To get the convergence, we consider the
potential function $\varphi$ of (5.1.3) again. The evolution of
$\varphi$ is given by (5.1.4). We renormalize the function
$\varphi$ by
$$
\tilde{\varphi}(x,t) =\varphi(x,t)+\int b(t)dt, \qquad \text{on
}\;M\times[0,+\infty).
$$
Then, since $r=0$, $\tilde{\varphi}$ evolves by \be \frac{\partial
\tilde{\varphi}}{\partial t} =\Delta\tilde{\varphi},\qquad
\text{on}
\;M\times[0,+\infty).  
\ee {}From the proof of Lemma 5.1.2, we get \be \frac{\partial
}{\partial t}|\nabla \tilde{\varphi}|^2= \Delta|\nabla
\tilde{\varphi}|^2-2|\nabla^2\tilde{\varphi}|^2. 
\ee Clearly, we have \be \frac{\partial }{\partial
t}\tilde{\varphi}^2
=\Delta\tilde{\varphi}^2-2|\nabla\tilde{\varphi}|^2. 
\ee Thus it follows that
$$
\frac{\partial }{\partial t}(t|\nabla\tilde{\varphi}|^2
+\tilde{\varphi}^2)\le\Delta(t|\nabla \tilde{\varphi}|^2+
\tilde{\varphi}^2).
$$
Hence by applying the maximum principle, there exists a positive
constant $C_3$ depending only on the initial metric such that \be
|\nabla\tilde{\varphi}|^2(x,t)\le \frac{C_3}{1+t},
\qquad \text{on}\;M\times[0,+\infty). 
\ee In the following we will use this decay estimate to obtain a
decay estimate for the scalar curvature.

By the evolution equations (5.1.2) and (5.1.10), we have
\begin{align*}
\frac{\partial }{\partial t}(R+2|\nabla\tilde{\varphi}|^2) &
=\Delta(R+2|\nabla\tilde{\varphi}|^2)
+R^2-4|\nabla^2\tilde{\varphi}|^2\\
&  \le \Delta(R+2|\nabla\tilde{\varphi}|^2)-R^2
\end{align*}
since $R^2=(\Delta\tilde{\varphi})^2\le
2|\nabla^2\tilde{\varphi}|^2$. Thus by using (5.1.12), we have
\begin{align*}
& \frac{\partial }{\partial
t}[t(R+2|\nabla\tilde{\varphi}|^2)]\\
&  \le \Delta[t(R+2|\nabla\tilde{\varphi}|^2)]
-tR^2+R+2|\nabla\tilde{\varphi}|^2\\
& \leq\Delta[t(R+2|\nabla\tilde{\varphi}|^2)]
-t(R+2|\nabla\tilde{\varphi}|^2)^2+(1+4t|\nabla
\tilde{\varphi}|^2)(R+2|\nabla\tilde{\varphi}|^2)\\
&  \le \Delta[t(R+2|\nabla\tilde{\varphi}|^2)]
-[t(R+2|\nabla\tilde{\varphi}|^2)
-(1+4C_3)](R+2|\nabla\tilde{\varphi}|^2)\\
&  \le \Delta[t(R+2|\nabla\tilde{\varphi}|^2)]
\end{align*}
wherever $t(R+2|\nabla\tilde{\varphi}|^2)\ge (1+4C_3)$. Hence by
the maximum principle, there holds \be
R+2|\nabla\tilde{\varphi}|^2\le\frac{C_4}{1+t},\qquad
\text{on}\;M\times[0,+\infty)  
\ee for some positive constant $C_4$ depending only on the initial
metric.

On the other hand, the scalar curvature satisfies
$$
\frac{\partial R}{\partial t} =\Delta R+R^2,\qquad \text{on
}\;M\times[0,+\infty).
$$
It is not hard to see that \be
R\ge\frac{R_{\min}(0)}{1-R_{\min}(0)t},\qquad
\text{on }\;M\times[0,+\infty),  
\ee by using the maximum principle. So we obtain the decay
estimate for the scalar curvature \be
|R(x,t)|\le\frac{C_5}{1+t},\qquad \text{on }\;M\times[0,+\infty),
\ee for some positive constant $C_5$ depending only on the initial
metric.

\begin{theorem}[{Hamilton \cite{Ha88}}]
On a compact surface with $\chi(M)=0$, for any initial metric the
solution of the normalized Ricci flow $(5.1.1)$ exists for all
time and converges in $C^\infty$ topology to a flat metric.
\end{theorem}

\begin{pf}
Since $\frac{\partial \tilde{\varphi}}{\partial
t}=\Delta\tilde{\varphi}$, it follows from the maximum principle
that
$$
|\tilde\varphi(x,t)|\le C_6,\qquad \text{on }\; M\times
[0,+\infty)
$$
for some positive constant $C_6$ depending only on the initial
metric. Recall that $\Delta\tilde{\varphi}=R$. We thus obtain for
any tangent vector $v\in T_xM$ at a point $x\in M$,
\begin{align*}
\frac{d}{dt}|v|^2_t&
= \(\frac{\partial }{\partial t}g_{ij}(x,t)\)v^iv^j\\
&  = -R(x,t)|v|^2_t
\end{align*}
and then
\begin{align*}
\left|\log\frac{|v|^2_t}{|v|^2_0}\right|
&  = \left|\int_0^t\frac{d}{dt}\log\right|v_t|^2dt|\\
&  = \left|\int_0^tR(x,t)dt\right|\\
&  = \left|\tilde{\varphi}(x,t)-\tilde{\varphi}(x,0)\right|\\
&  \le 2C_6,
\end{align*}
for all $x\in M$ and $t\in[0,+\infty).$ This shows that the
solution $g_{ij}(t)$ of the normalized Ricci flow are all
equivalent. This gives us control of the diameter and injectivity
radius.

As before, by Shi's derivative estimates of the unnormalized Ricci
flow, all derivatives and higher order derivatives of the
curvature of the solution $g_{ij}$ of the normalized Ricci flow
(5.1.1) are uniformly bounded on $M\times[0,+\infty)$. By the
virtue of Hamilton's compactness theorem (Theorem 4.1.5) we see
that the solution $g_{ij}(t)$ subsequentially converges in
$C^\infty$ topology. The decay estimate (5.1.15) implies that each
limit must be a flat metric on $M$. Clearly, we will finish the
proof if we can show that limit is unique.

Note that the solution $g_{ij}(t)$ is changing conformally under
the Ricci flow (5.1.1) on surfaces. Thus each limit must be
conformal to the initial metric, denoted by $\bar{g}_{ij}$. Let us
denote $g_{ij}(\infty)=e^u\bar{g}_{ij}$ to be a limiting metric.
Since $g_{ij}(\infty)$ is flat, it is easy to compute
$$
0=e^{-u}(\bar{R}-\bar{\Delta}u),\qquad \mbox{on }\;M,
$$
where $\bar{R}$ is the curvature of $\bar{g}_{ij}$ and
$\bar{\Delta}$ is the Laplacian in the metric $\bar{g}_{ij}$. The
solution of Poission equation
$$\bar{\Delta}u=\bar{R},\qquad \mbox{on }\;M$$
is unique up to constant. Moreover the constant must be also
uniquely determined since the area of the solution of the
normalized Ricci flow (5.1.1) is constant in time. So the limit is
unique and we complete the proof of Theorem 5.1.5.
\end{pf}

\medskip
{\it Case} (3): $\chi(M)>0$, i.e., $r>0$.

\smallskip
This is the most difficult case. The first proof is due to Ben
Chow\cite{Chow}, based on the important work of Hamilton
\cite{Ha88}.  By now there exist several proofs (cf.
Bartz-Struwe-Ye \cite{BSY} and Struwe \cite{St}, etc). But, in
contrast to the previous two cases, none of the proofs depends
only on maximum principle type argument. In fact, all the proofs
rely on some kind of combination of the maximum principle argument
and certain integral estimate of the curvature.

In the pioneer work \cite{Ha88}, Hamilton considered the important
special case when the initial metric is of positive scalar
curvature. He introduced an integral quantity
$$
E=\int_M R \log R \;dV,
$$
which he called {\bf entropy}, for the (normalized) Ricci flow on
a surface $M$ with positive curvature, and showed that the entropy
is monotone decreasing under the flow. By combining this entropy
estimate with the Harnack inequality for the curvature (Corollary
2.5.3), Hamilton obtained the uniform bound on the curvature of
the normalized Ricci flow on $M$ with positive curvature.
Furthermore, he showed that the evolving metric converges to a
shrinking Ricci soliton on $M$ and that the shrinking Ricci
soliton must be a round metric on the $2$-sphere $\mathbb S^2$.
Subsequently, Chow \cite{Chow} extended Hamilton's work to the
general case when the curvature may change signs. More precisely,
he proved that given any initial metric on a compact surface $M$
with $\chi(M)>0$, the evolving metric under the (normalized) Ricci
flow will have positive curvature after a finite time. Hence, when
combined with Hamilton's result, B. Chow's result implies that the
solution under the normalized flow on $M$ converges to the round
metric on $\mathbb S^2$.

In the following we will basically follow the arguments of
Hamilton \cite{Ha88} and Chow \cite{Chow}, except when we prove
the uniform bound of the evolving scalar curvature we will present
a new argument using the Harnack inequality of Chow \cite{Chow}
and Perelman's no local collapsing theorem I$'$ (as was done in
the joint work of Bing-Long Chen and the authors \cite{CCZ} where
they considered the K\"ahler-Ricci flow of nonnegative holomorphic
bisectional curvature).

Given any initial metric on $M$ with $\chi(M)>0$, we consider the
solution $g_{ij}(t)$ of the normalized Ricci flow (5.1.1). Recall
that the (scalar) curvature $R$ satisfies the evolution equation
$$
\frac{\partial }{\partial t}R=\Delta R+R^2-rR.
$$
The corresponding ODE is \be
\frac{ds}{dt}=s^2-rs.  
\ee

Let us choose $c>1$ and close to 1 so that $r/(1-c)<\min_{x\in
M}R(x,0)$. It is clear that the function $s(t)=r/(1-ce^{rt})<0$ is
a solution of the ODE (5.1.16) with $s(0)<\min\limits_{x\in
M}R(x,0)$. Then the difference of $R$ and $s$ evolves by \be
\frac{\partial }{\partial t}(R-s)=\Delta(R-s)+(R-r+s)(R-s).
\ee Since $\min_{x\in M}R(x,0)-s(0)>0$, the maximum principle
implies that $R-s>0$ for all times.

First, we need the Harnack inequality obtained by B. Chow
\cite{Chow}, which is an extension of Theorem 2.5.2, for the
normalized Ricci flow whose curvature may change signs.

Consider the quantity
$$
L=\log(R-s).
$$
It is easy to compute
$$
\frac{\partial L}{\partial t}=\Delta L+|\nabla L|^2+R-r+s.
$$
Set
$$
Q=\frac{\partial L}{\partial t}-|\nabla L|^2-s=\Delta L+R-r.
$$
By a direct computation and using the estimate (5.1.8), we have
\begin{align*}
\frac{\partial }{\partial t}Q&  = \Delta\(\frac{\partial
L}{\partial t}\)
+(R-r)\Delta L+\frac{\partial R}{\partial t}\\
&  =\Delta Q+2|\nabla^2L|^2+2\langle\nabla L,\nabla(\Delta
L)\rangle+R|\nabla L|^2\\
&\quad +(R-r)\Delta L+\Delta R+R(R-r)\\
&  =\Delta Q+2|\nabla^2L|^2+2\langle\nabla L,\nabla
Q\rangle+2(R-r)\Delta L+(R-r)^2\\
&\quad +(r-s)\Delta L+s|\nabla L|^2+r(R-r)\\
&  = \Delta Q+2\langle\nabla L,\nabla Q\rangle
+2|\nabla^2L|^2+2(R-r)\Delta L+(R-r)^2\\
&\quad +(r-s)Q+s|\nabla L|^2+s(R-r)\\
&  \ge \Delta Q+2\langle \nabla L,\nabla
Q\rangle+Q^2+(r-s)Q+s|\nabla L|^2-C.
\end{align*}
Here and below $C$ is denoted by various positive constants
depending only on the initial metric.

In order to control the bad term $s|\nabla L|^2$, we consider
\begin{align*}
\frac{\partial }{\partial t}(sL)
&  =\Delta(sL)+s|\nabla L|^2+s(R-r+s)+s(s-r)L\\
&  \ge \Delta(sL)+2\langle\nabla L,\nabla(sL)\rangle -s|\nabla
L|^2-C
\end{align*}
by using the estimate (5.1.8) again. Thus
\begin{align*}
\frac{\partial }{\partial t}(Q+sL) &
\ge\Delta(Q+sL)+2\langle\nabla
L,\nabla(Q+sL)\rangle+Q^2+(r-s)Q-C\\
&  \ge \Delta(Q+sL)+2\langle\nabla L,\nabla(Q+sL)\rangle
+\frac{1}{2}[(Q+sL)^2-C^2],
\end{align*}
since $sL$ is bounded by (5.1.8). This, by the maximum principle,
implies that
$$
Q\ge -C,\qquad \text{for\;all }\;t\in[0,+\infty).
$$
Then for any two points $x_1,x_2\in M$ and two times $t_2>t_1\ge
0$, and a path $\gamma:[t_1,t_2]\rightarrow M$ connecting $x_1$ to
$x_2$, we have
\begin{align*}
L(x_2,t_2)-L(x_1,t_1)
&  = \int_{t_1}^{t_2}\frac{d}{dt}L(\gamma(t),t)dt\\
&  = \int_{t_1}^{t_2}\(\frac{\partial L}{\partial t}+\langle\nabla
L,\dot{\gamma}\rangle\)dt\\
&  \ge -\frac{1}{4}\Delta-C(t_2-t_1)
\end{align*}
where
\begin{align*}
\Delta&  = \Delta(x_1,t_1;x_2,t_2)\\
&  =
\inf\left\{\int_{t_1}^{t_2}|\dot{\gamma}(t)|^2_{g_{ij}(t)}dt\, |\,
\gamma:[t_1,t_2]\rightarrow
M\;\text{with}\;\gamma(t_1)=x_1,\gamma(t_2)=x_2\right\}.
\end{align*}
Thus we have proved the following Harnack inequality of B. Chow.

\begin{lemma}[{Chow \cite{Chow}}]
There exists a positive constant C depending only on the initial
metric such that for any $x_1$, $x_2\in M$ and $t_2>t_1\ge 0$,
$$
R(x_1,t_1)-s(t_1) \le
e^{\frac{\Delta}{4}+C(t_2-t_1)}(R(x_2,t_2)-s(t_2))
$$
where
$$
\Delta=\inf\left\{\int_{t_1}^{t_2}|\dot{\gamma}(t)|^2_tdt\ |\
\gamma:[t_1,t_2]\rightarrow M\;with
\;\gamma(t_1)=x_1,\gamma(t_2)=x_2\right\}.
$$
\end{lemma}

We now state and prove the following uniform bound estimate for
the curvature, a consequence of the results of Hamilton
\cite{Ha88} and Chow \cite{Chow}. We remark the special case when
the scalar curvature $R>0$ is first proved by Hamilton
\cite{Ha88}. As we mentioned before, the proof here is adapted
from \cite{CCZ}.


\begin{proposition} [{cf. Lemma 5.74 and Lemma 5.76 of \cite{CK04}}]
Let $(M,g_{ij}(t))$ be a solution of the normalized Ricci flow on
a compact surface with $\chi(M)>0$. Then there exist a time
$t_0>0$ and a positive constant $C$ such that the estimate
$$
C^{-1}\le R(x,t)\le C
$$
holds for all $x\in M$ and $t\in [t_0,+\infty).$
\end{proposition}

\begin{pf}
Recall that
$$
R(x,t)\ge s(t) =\frac{r}{1-ce^{rt}},\qquad \text{on}\; M\times
[0,+\infty).
$$
For any $\varepsilon\in(0,r)$, there exists a large enough $t_0>0$
such that \be R(x,t)\ge-\varepsilon^2,\qquad \text{on}\;
M\times[t_0,+\infty).
\ee Let $t$ be any fixed time with $t\ge t_0+1$. Obviously there
is some point $x_0\in M$ such that $R(x_0,t+1)=r$.

Consider the geodesic ball $B_t(x_0,1)$, centered at $x_0$ and
radius 1 with respect to the metric at the fixed time $t$. For any
point $x\in B_t(x_0,1)$, we choose a geodesic $\gamma$:
$[t,t+1]\rightarrow M$ connecting $x$ and $x_0$ with respect to
the metric at the fixed time $t$. Since
$$
\frac{\partial }{\partial t}g_{ij}=(r-R)g_{ij} \le 2rg_{ij}\qquad
\text{on }\;M\times[t_0,+\infty),
$$
we have
$$
\int_t^{t+1}|\dot{\gamma}(\tau)|^2_\tau d\tau\le
e^{2r}\int_t^{t+1} |\dot{\gamma}(\tau)|^2_t d\tau\leq e^{2r}.
$$
Then by Lemma 5.1.6, we have
\begin{align}
R(x,t)& \leq s(t)+\exp\left\{\frac{1}{4}e^{2r}+C\right\}
\cdot(R(x_0,t+1)-s(t+1))\\
& \leq C_1,\qquad \text{as }\; x\in B_t(x_0,1),\nn
\end{align} 
for some positive constant $C_1$ depending only on the initial
metric. Note the the corresponding unnormalized Ricci flow in this
case has finite maximal time since its volume decreases at a fixed
rate $-4\pi \chi (M) < 0$. Hence the no local collapsing theorem
I$'$ (Theorem 3.3.3) implies that the volume of $B_t(x_0,1)$ with
respect to the metric at the fixed time $t$ is bounded from below
by \be
\Vol_t(B_t(x_0,1))\geq C_2 
\ee for some positive constant $C_2$ depending only on the initial
metric.

We now want to bound the diameter of $(M,g_{ij}(t))$ from above.
The following argument is analogous to Yau in \cite{Y76} where he
got a lower bound for the volume of geodesic balls of a complete
Riemannian manifold with nonnegative Ricci curvature. Without loss
of generality, we may assume that the diameter of $(M,g_{ij}(t))$
is at least 3. Choose a point $x_1\in M$ such that the distance
$d_t(x_0,x_1)$ between $x_1$ and $x_0$ with respect to the metric
at the fixed time $t$ is at least a half of the diameter of
$(M,g_{ij}(t))$. By (5.1.18), the standard Laplacian comparison
theorem (c.f. \cite{ScY}) implies
$$
\Delta \rho^2=2\rho\Delta\rho+2\leq2(1+\varepsilon\rho)+2
$$
in the sense of distribution, where $\rho$ is the distance
function from $x_1$ (with respect to the metric $g_{ij}(t))$. That
is, for any $\varphi\in C_0^\infty(M)$, $\varphi\geq0$, we have
\be -\int_M\nabla\rho^2\cdot\nabla\varphi
\leq\int_M[2(1+\varepsilon\rho)+2]\varphi. 
\ee Since $C_0^\infty(M)$ functions can be approximated by
Lipschitz functions in the above inequality, we can set
$\varphi(x)=\psi(\rho(x))$, $x\in M$, where $\psi(s)$ is given by
$$
\psi(s)=\begin{cases}
1,& 0\leq s\leq d_t(x_0,x_1)-1,\\
\psi'(s)=-\frac{1}{2},& d_t(x_0,x_1)-1\leq s\leq
d_t(x_0,x_1)+1,\\
0,& s\geq d_t(x_0,x_1)+1.
\end{cases}
$$
Thus, by using (5.1.20), the left hand side of (5.1.21) is
\begin{displaymath}
\begin{split}
&-\int_M\nabla\rho^2\cdot\nabla\varphi \\
&  =\int_{B_t(x_1,d_t(x_0,x_1)+1)
\setminus B_t(x_1,d_t(x_0,x_1)-1)}\rho\\
&  \geq(d_t(x_0,x_1)-1)\Vol_t(B_t(x_1,d_t(x_0,x_1)+1)
\setminus B_t(x_1,d_t(x_0,x_1)-1))\\
&  \geq(d_t(x_0,x_1)-1)\Vol_t(B_t(x_0,1))\\
&  \geq(d_t(x_0,x_1)-1)C_2,
\end{split}
\end{displaymath}
 and the right hand side of (5.1.21) is
\begin{displaymath}
\begin{split}
\int_M[2(1+\varepsilon\rho)+2]\varphi
&  \leq\int_{B_t(x_1,d_t(x_0,x_1)+1)}[2(1+\varepsilon\rho)+2]\\
&  \leq [2(1+\varepsilon d_t(x_0,x_1))+4]
\Vol_t(B_t(x_1,d_t(x_0,x_1)+1))\\
&  \leq [2(1+\varepsilon d_t(x_0,x_1))+4]A
\end{split}
\end{displaymath}
where $A$ is the area of $M$ with respect to the initial metric.
Here we have used the fact that the area of solution of the
normalized Ricci flow is constant in time. Hence
$$
C_2(d_t(x_0,x_1)-1)\leq[2(1+\varepsilon d_t(x_0,x_1))+4] A,
$$
which implies, by choosing $\varepsilon>0$ small enough,
$$
d_t(x_0,x_1)\leq C_3
$$
for some positive constant $C_3$ depending only on the initial
metric. Therefore, the diameter of $(M,g_{ij}(t))$ is uniformly
bounded above by \be
{\rm diam}\,(M,g_{ij}(t))\leq 2C_3  
\ee for all $t\in[t_0,+\infty)$.

We then argue, as in deriving (5.1.19), by applying Lemma 5.1.6
again to obtain
$$
R(x,t)\leq C_4,\ \ \ \text{on }\; M\times[t_0,+\infty)
$$
for some positive constant $C_4$ depending only on the initial
metric.

It remains to prove a positive lower bound estimate of the
curvature. First, we note that the function $s(t)\to 0$ as
$t\rightarrow+\infty$, and the average scalar curvature of the
solution equals to $r$, a positive constant. Thus the Harnack
inequality in Lemma 5.1.6 and the diameter estimate (5.1.22) imply
a positive lower bound for the curvature. Therefore we have
completed the proof of Proposition 5.1.7.
\end{pf}

Next we consider long-time convergence of the normalized flow.

Recall that the trace-free part of the Hessian of the potential
$\varphi$ of the curvature is the tensor $M_{ij}$ defined by
$$
M_{ij}=\nabla_i\nabla_j\varphi-\frac{1}{2}\Delta\varphi\cdot
g_{ij},
$$
where by (5.1.3),
$$
\Delta\varphi=R-r.
$$

\begin{lemma} [{Hamilton \cite{Ha88}}]
We have \be \frac{\partial}{\partial t}|M_{ij}|^2
=\Delta|M_{ij}|^2-2|\nabla_k M_{ij}|^2-2R|M_{ij}|^2,
\; \text{ on }\; M\times[0,+\infty). 
\ee
\end{lemma}

\begin{pf} This follows from a standard computation
(e.g., cf. Editors' note on p. 217 of \cite{CCCY}).

First we note the time-derivative of the Levi-Civita connection is
\begin{displaymath}
\begin{split}
\frac{\partial}{\partial t}\Gamma^k_{ij} &
=\frac{1}{2}g^{kl}\(\nabla_i\frac{\partial}{\partial t}g_{jl}
+\nabla_j\frac{\partial}{\partial t}g_{il}
-\nabla_l\frac{\partial}{\partial t}g_{ij}\)\\
&  =\frac{1}{2}\(-\nabla_iR\cdot\delta_j^k
-\nabla_jR\cdot\delta_i^k+\nabla^kR\cdot g_{ij}\).
\end{split}
\end{displaymath}
By using this and (5.1.4), we have
\begin{displaymath}
\begin{split}
  \frac{\partial}{\partial t}M_{ij}
& = \nabla_i\nabla_j\(\frac{\partial \varphi}{\partial t}\)
    -\(\frac{\partial}{\partial t}\Gamma^k_{ij}\)\nabla_k\varphi
    -\frac{1}{2}\frac{\partial}{\partial t}[(R-r)g_{ij}]\\
& = \nabla_i\nabla_j\Delta\varphi
+\frac{1}{2}(\nabla_iR\cdot\nabla_j\varphi
    +\nabla_jR\cdot\nabla_i\varphi
-\langle\nabla R,\nabla\varphi\rangle g_{ij})\\
& \quad -\frac{1}{2}\Delta R\cdot g_{ij}+rM_{ij}.
\end{split}
\end{displaymath}
Since on a surface,
$$
R_{ijkl}=\frac{1}{2}R(g_{il}g_{jk}-g_{ik}g_{jl}),
$$
we have
\begin{displaymath}
\begin{split}
& \nabla_i\nabla_j\Delta\varphi \\
& = \nabla_i\nabla_k\nabla_j\nabla^k\varphi
-\nabla_i(R_{jl}\nabla^l\varphi)\\
& = \nabla_k\nabla_i\nabla_j\nabla^k\varphi
-R^l_{ikj}\nabla_l\nabla^k\varphi-R_{il}\nabla_j\nabla^l\varphi
    -R_{jl}\nabla_i\nabla^l\varphi
-\nabla_iR_{jl}\nabla^l\varphi\\
& = \Delta\nabla_i\nabla_j\varphi
-\nabla^k(R^l_{ikj}\nabla_l\varphi)-R^l_{ikj}\nabla_l\nabla^k\varphi\\
&\quad  -R_{il}\nabla_j\nabla^l\varphi
-R_{jl}\nabla_i\nabla^l\varphi-\nabla_iR_{jl}\nabla^l\varphi\\
& = \Delta\nabla_i\nabla_j\varphi
-\frac{1}{2}(\nabla_iR\cdot\nabla_j\varphi+\nabla_jR\cdot\nabla_i\varphi
    -\langle\nabla R,\nabla \varphi\rangle g_{ij})\\
&\quad  -2R\(\nabla_i\nabla_j\varphi-\frac{1}{2}\Delta\varphi\cdot
    g_{ij}\).
\end{split}
\end{displaymath}
Combining these identities, we get
\begin{displaymath}
\begin{split}
\frac{\partial}{\partial t}M_{ij} & =
\Delta\nabla_i\nabla_j\varphi-\frac{1}{2}\Delta R\cdot g_{ij}
+(r-2R)M_{ij}\\
& = \Delta\(\nabla_i\nabla_j\varphi-\frac{1}{2}(R-r)g_{ij}\)
+(r-2R)M_{ij}.
\end{split}
\end{displaymath}
Thus the evolution $M_{ij}$ is given by \be \frac{\partial
M_{ij}}{\partial t}=\Delta M_{ij}+(r-2R)M_{ij}.
\ee Now the lemma follows from (5.1.24) and a straightforward
computation.
\end{pf}

Proposition 5.1.7 tells us that the curvature $R$ of the solution
to the normalized Ricci flow is uniformly bounded from below by a
positive constant for $t$ large. Thus we can apply the maximum
principle to the equation (5.1.23) in Lemma 5.1.8 to obtain the
following estimate.

\begin{proposition}  [{Hamilton \cite{Ha88} and Chow \cite{Chow}}]
Let $(M,g_{ij}(t))$ be a solution of the normalized Ricci flow on
a compact surface with $\chi(M)>0$. Then there exist positive
constants $c$ and $C$ depending only on the initial metric such
that
$$
|M_{ij}|^2\leq Ce^{-ct},\quad \text{ on }\; M\times[0,+\infty).
$$
\end{proposition}

Now we consider a modification of the normalized Ricci flow.
Consider the equation \be \frac{\partial}{\partial
t}g_{ij}=2M_{ij}
=(r-R)g_{ij}+2\nabla_i\nabla_j\varphi. 
\ee As we saw in Section 1.3, the solution of this modified flow
differs from that of the normalized Ricci flow only by a one
parameter family of diffeomorphisms generated by the gradient
vector field of the potential function $\varphi$. Since the
quantity $|M_{ij}|^2$ is invariant under diffeomorphisms, the
estimate $|M_{ij}|^2\leq Ce^{-ct}$ also holds for the solution of
the modified flow (5.1.25). This exponential decay estimate then
implies the solution $g_{ij}(x,t)$ of the modified flow (5.1.25)
converges exponentially to a continuous metric $g_{ij}(\infty)$ as
$t\rightarrow+\infty$. Furthermore, by the virtue of Hamilton's
compactness theorem (Theorem 4.1.5) we see that the solution
$g_{ij}(x,t)$ of the modified flow actually converges
exponentially in $C^\infty$ topology to $g_{ij}(\infty)$. Moreover
the limiting metric $g_{ij}(\infty)$ satisfies
$$
M_{ij}=(r-R)g_{ij}+2\nabla_i\nabla_j\varphi=0,\ \ \mbox{ on }\; M.
$$
That is, the limiting metric is a shrinking gradient Ricci soliton
on the surface $M$.

The next result was first obtained by Hamilton in \cite{Ha88}. The
following simplified proof by using the Kazdan-Warner identity has
been widely known to experts in the field (e.g., cf. Proposition
5.21 of \cite{CK04}).

\begin{proposition} [{Hamilton \cite{Ha88}}]
On a compact surface there are no shrinking Ricci solitons other
than constant curvature.
\end{proposition}

\begin{pf}
By definition, a shrinking Ricci soliton on a compact surface $M$
is given by \be
\nabla_iX_j+\nabla_jX_i=(R-r)g_{ij} 
\ee for some vector field $X={X_j}$. By contracting the above
equation by $Rg^{-1}$, we have
$$
2R(R-r)=2R\; \text{div}\, X,
$$
and hence
$$
\int_M(R-r)^2 dV=\int_MR(R-r)dV=\int_M R\; \text{div}\, X dV.
$$
Since $X$ is a conformal vector field (by the Ricci soliton
equation (5.1.26)), by integrating by parts and applying the
Kazdan-Warner identity \cite{KW}, we obtain
$$
\int_M(R-r)^2 dV=-\int_M\langle\nabla R,X\rangle dV=0.
$$
Hence $R\equiv r$, and the lemma is proved.
\end{pf}

Now back to the solution of the modified flow (5.1.25). We have
seen the curvature converges exponentially to its limiting value
in the $C^\infty$ topology. But since there are no nontrivial
soliton on $M$, we must have $R$ converging exponentially to the
constant value $r$ in the $C^\infty$ topology. This then implies
that the unmodified flow (5.1.1) will converge to a metric of
positive constant curvature in the $C^\infty$ topology.

In conclusion, we have proved the following main theorem of this
section.

\begin{theorem}[{Hamilton \cite{Ha88} and Chow \cite{Chow}}]
On a compact surface with $\chi(M)>0$, for any initial metric, the
solution of the normalized Ricci flow $(5.1.1)$ exists for all
time, and converges in the $C^\infty$ topology to a metric with
positive constant curvature.
\end{theorem}

\section{Differentiable Sphere Theorems in 3-D and 4-D}

An important problem in Riemannian geometry is to understand the
influence of curvatures, in particular the sign of curvatures, on
the topology of underlying manifolds. Classical results of this
type include sphere theorem and its refinements stated below
(e.g., cf. Theorem 6.1, Theorem 7.16, and Theorem 6.6 of
\cite{CE}). In this section we shall use the long-time behavior of
the Ricci flow on positively curved manifolds to establish
Hamilton's differentiable sphere theorems in dimensions three and
four. Our presentation is based on Hamilton \cite{Ha82, Ha86}.

Let us first recall some classical sphere theorems. Given a
Remannian manifold $M$, we denote by $K_M$ the sectional curvature
of $M$.

\medskip
{\bf Classical Sphere Theorems} (cf. \cite{CE}).\index{classical
sphere theorems} \ \emph{ Let $M$ be a complete, simply connected
$n$-dimensional manifold.}
\begin{itemize}
\item[(i)] \emph{If $\frac{1}{4}<K_M\le 1$, then $M$ is
homeomorphic to the $n$-sphere $\mathbb{S}^n$.} \item[(ii)]
\emph{There exists a positive constant $\delta\in(\frac{1}{4},1)$
such that if $\delta<K_M\le 1$, then $M$ is diffeomorphic to the
$n$-sphere $\mathbb{S}^n$.}
\end{itemize}

\medskip
Result (ii) is called the differentiable sphere theorem. If we
relax the assumptions on the strict lower bound in (i), then we
have the following rigidity result.

\medskip
{\bf Berger's Rigidity Theorem (cf. \cite{CE}).}\index{Berger's
rigidity theorem} \emph{ \ Let $M$ be a complete, simply connected
$n$-dimensional manifold with $\frac{1}{4}\le K_M\leq 1$. Then
either $M$ is homeomorphic to $\mathbb{S}^n$ or $M$ is isometric
to a symmetric space.}

\medskip
We remark that it follows from the classification of symmetric
spaces (see for example \cite{He}) that the only simply connected
symmetric spaces with positive curvature are $\mathbb{S}^n$,
$\mathbb{CP}^\frac{n}{2}$, $\mathbb{QP}^\frac{n}{4}$, and the
Cayley plane.

In early and mid 80's respectively, Hamilton \cite{Ha82, Ha86}
used the Ricci flow to prove the following differential sphere
theorems.

\begin{theorem}[{Hamilton \cite{Ha82}}]
A compact three-manifold with positive Ricci curvature must be
diffeomorphic to the three-sphere $\mathbb{S}^3$ or a quotient of
it by a finite group of fixed point free isometries in the
standard metric.
\end{theorem}

\begin{theorem}[{Hamilton \cite{Ha86}}]
A compact four-manifold with positive curvature operator is
diffeomorphic to the four-sphere $\mathbb{S}^4$ or the real
projective space $\mathbb{RP}^4$.
\end{theorem}

Note that in above two theorems, we only assume curvatures to be
strictly positive, but not any strong pinching conditions as in
the classical sphere theorems. In fact, one of the important
special features discovered by Hamilton is that if the initial
metric has positive curvature, then the metric will get rounder
and rounder as it evolves under the Ricci flow, at least in
dimension three and four, so any small initial pinching will get
improved. Indeed, the pinching estimate is a key step in proving
both Theorem 5.2.1 and 5.2.2.

The following results are concerned with compact three-manifolds
or four-manifolds with weakly positive curvatures.

\begin{theorem}[{Hamilton \cite{Ha86}}] \
\begin{itemize}
\item[(i)] A compact three-manifold with nonnegative Ricci
curvature is diffeomorphic to $\mathbb{S}^3$, or a quotient of one
of the spaces $\mathbb{S}^3$ or $\mathbb{S}^2\times\mathbb{R}^1$
or $\mathbb{R}^3$ by a group of fixed point free isometries in the
standard metrics. \item[(ii)] A compact four-manifold with
nonnegative curvature operator is diffeomorphic to $\mathbb{S}^4$
or $\mathbb{CP}^2$ or $\mathbb{S}^2\times \mathbb{S}^2$, or a
quotient of one of the spaces $\mathbb{S}^4$ or $\mathbb{CP}^2$ or
$\mathbb{S}^3\times\mathbb{R}^1$ or $\mathbb{S}^2\times
\mathbb{S}^2$ or $\mathbb{S}^2\times \mathbb{R}^2$ or
$\mathbb{R}^4$ by a group of fixed point free isometries in the
standard metrics.
\end{itemize}
\end{theorem}

The rest of the section will be devoted to prove Theorems
5.2.1-5.2.3 and the presentation follows Hamilton \cite{Ha82,
Ha86} (also cf. \cite{Ha95F}).

Recall that the curvature operator $M_{\alpha\beta}$ evolves by
\be \frac{\partial}{\partial t}M_{\alpha\beta}=\Delta
M_{\alpha\beta}+M^2_{\alpha\beta}+M_{\alpha\beta}^\#.
\ee where (see Section 1.3 and Section 2.4) $M_{\alpha\beta}^2$ is
the operator square
$$
M_{\alpha\beta}^2=M_{\alpha\gamma}M_{\beta\gamma}
$$
and $M_{\alpha\beta}^\#$ is the Lie algebra $so(n)$ square
$$
M_{\alpha\beta}^\#=C_\alpha^{\gamma\zeta}
C_\beta^{\eta\theta}M_{\gamma\eta}M_{\zeta\theta}.
$$

We begin with the curvature pinching estimates of the Ricci flow
in three dimensions. In dimension $n=3$, we know that
$M_{\alpha\beta}^\#$ is the adjoint matrix of $M_{\alpha\beta}$.
If we diagonalize $M_{\alpha\beta}$ with eigenvalues $\lambda \ge
\mu \ge \nu$ so that
$$
(M_{\alpha\beta})=\left(
\begin{array}{ccc}\lambda&  \ &  \ \\
\ &  \mu &  \ \\
\ &  \ &  \nu
\end{array}
\right),
$$
then $M_{\alpha\beta}^2$ and $M_{\alpha\beta}^\#$ are also
diagonal, with
$$
(M_{\alpha\beta}^2)=\left(
\begin{array}{ccc}\lambda^2&  \ &  \ \\
\ &  \mu^2 &  \ \\
\ &  \ &  \nu^2
\end{array}
\right) \ \ \text{and} \ \ (M_{\alpha\beta}^\#)=\left(
\begin{array}{ccc}\mu\nu&  \ &  \ \\
\ &  \lambda\nu &  \ \\
\ &  \ &  \lambda\mu
\end{array}
\right),
$$
and the ODE corresponding to PDE (5.2.1) is then given by the
system \be
      \left\{
       \begin{array}{lll}
  \frac{d}{dt}\lambda=\lambda^2+\mu\nu,\\[4mm]
  \frac{d}{dt}\mu=\mu^2+\lambda\nu,\\[4mm]
  \frac{d}{dt}\nu=\nu^2+\lambda\mu.
       \end{array}
    \right.
\ee

\begin{lemma} [{Hamilton \cite{Ha82, Ha95F}}]
For any $\varepsilon\in[0,\frac{1}{3}]$, the pinching condition
$$
R_{ij}\geq0\ \ \ and\ \ \ R_{ij}\geq\varepsilon Rg_{ij}
$$
is preserved by the Ricci flow.
\end{lemma}

\begin{pf}
If we diagonalize the $3\times 3$ curvature operator matrix
$M_{\alpha\beta}$ with eigenvalues $\lambda \ge \mu \ge \nu$, then
nonnegative sectional curvature corresponds to $\nu\geq0$ and
nonnegative Ricci curvature corresponds to the inequality
$\mu+\nu\geq0$. Also, the scalar curvature $R=\lambda+\mu+\nu$. So
we need to show
$$
\mu+\nu\geq 0 \quad \text{and} \quad \mu+\nu\geq\delta \lambda, \;
\text{ with }\; \delta=2\varepsilon/(1-2\varepsilon),
$$
are preserved by the Ricci flow. By Hamilton's advanced maximum
principle (Theorem 2.3.1), it suffices to show that the closed
convex set
$$
K=\{M_{\alpha\beta}\ |\ \mu+\nu\geq0\; \text{ and }\;
\mu+\nu\geq\delta\lambda\}
$$
is preserved by the ODE system (5.2.2).

Now suppose we have diagonalized $M_{\alpha\beta}$ with
eigenvalues $\lambda \ge \mu \ge \nu$ at $t=0$, then both
$M^2_{\alpha\beta}$ and $M_{\alpha\beta}^\#$ are diagonal, so the
matrix $M_{\alpha\beta}$ remains diagonal for $t>0$. Moreover,
since
$$
\frac{d}{dt}(\mu-\nu)=(\mu-\nu)(\mu+\nu-\lambda),
$$
it is clear that $\mu\geq \nu$ for $t>0$ also. Similarly, we have
$\lambda\geq\mu$ for $t>0$. Hence the inequalities $\lambda \ge
\mu \ge \nu$ persist. This says that the solutions of the ODE
system (5.2.2) agree with the original choice for the eigenvalues
of the curvature operator.

The condition $\mu+\nu\geq 0$ is clearly preserved by the ODE,
because
$$
\frac{d}{dt}(\mu+\nu)=\mu^2+\nu^2+\lambda(\mu+\nu)\geq 0.
$$

It remains to check
$$
\frac{d}{dt}(\mu+\nu)\geq\delta\frac{d}{dt}\lambda
$$
or
$$
\mu^2+\lambda\nu+\nu^2+\lambda\mu\geq\delta(\lambda^2+\mu\nu)
$$
on the boundary where
$$
\mu+\nu=\delta \lambda\geq0.
$$
In fact, since
$$
(\lambda-\nu)\mu^2+(\lambda-\mu)\nu^2\geq0,
$$
we have
$$
\lambda(\mu^2+\nu^2)\geq(\mu+\nu)\mu\nu.
$$
Hence
\begin{displaymath}
\begin{split}
    \mu^2+\mu\lambda+\nu^2+\nu\lambda
    &  \geq\(\frac{\mu+\nu}{\lambda}\)(\lambda^2+\mu\nu)\\[3mm]
    &  =\delta(\lambda^2+\mu\nu).
\end{split}
\end{displaymath}
So we get the desired pinching estimate.
\end{pf}

\begin{proposition} [cf. Hamilton \cite{Ha86} or \cite{Ha95F}]
Suppose that the initial metric of the solution to the Ricci flow
on $M^3\times [0,T)$ has positive Ricci curvature. Then for any
$\varepsilon>0$ we can find $C_\varepsilon<+\infty$ such that
$$
\left|R_{ij}-\frac{1}{3}Rg_{ij}\right|\leq\varepsilon
R+C_\varepsilon
$$
for all subsequent $t\in [0, T)$.
\end{proposition}

\begin{pf}
Again we consider the ODE system (5.2.2). Let $M_{\alpha\beta}$ be
diagonalized with eigenvalues $\lambda \ge \mu \ge \nu$ at $t=0$.
We saw in the proof of Lemma 5.2.4 the inequalities $\lambda \ge
\mu \ge \nu$ persist for $t>0$. We only need to show that there
are positive constants $\delta$ and $C$ such that the closed
convex set
$$
K=\{\ M_{\alpha\beta}\ |\ \lambda-\nu\leq
C(\lambda+\mu+\nu)^{1-\delta}\}
$$
is preserved by the ODE.

We compute
$$
\frac{d}{dt}(\lambda-\nu)=(\lambda-\nu)(\lambda+\nu-\mu)
$$
and
\begin{displaymath}
\begin{split}
   \frac{d}{dt}(\lambda+\mu+\nu)
& = (\lambda+\mu+\nu)(\lambda+\nu-\mu)+\mu^2\\[2mm]
&\quad  +\mu(\mu+\nu)+\lambda(\mu-\nu)\\[2mm]
& \geq (\lambda+\mu+\nu)(\lambda+\nu-\mu)+\mu^2.
\end{split}
\end{displaymath}
Thus, without loss of generality, we may assume $\lambda-\nu>0$
and get
$$
\frac{d}{dt}\log(\lambda-\nu)=\lambda+\nu-\mu
$$
and
$$
\frac{d}{dt}\log(\lambda+\mu+\nu)\geq
\lambda+\nu-\mu+\frac{\mu^2}{\lambda+\mu+\nu}.
$$
By Lemma 5.2.4, there exists a positive constant $C$ depending
only on the initial metric such that
$$
\lambda\leq \lambda+\mu\leq C(\mu+\nu)\leq 2C\mu,
$$
$$
\lambda+\nu-\mu\leq \lambda+\mu+\nu\leq 6C\mu,
$$
and hence with $\epsilon={1}/{36C^2}$,
$$
\frac{d}{dt}\log(\lambda+\mu+\nu)\geq(1+\epsilon)(\lambda+\nu-\mu).
$$
Therefore with $(1-\delta)={1}/{(1+\epsilon)}$,
$$
\frac{d}{dt}\log((\lambda-\nu)/(\lambda+\mu+\nu)^{1-\delta})\leq0.
$$
This proves the proposition.
\end{pf}

We now are ready to prove Theorem 5.2.1.

\medskip
{\bf\em Proof of Theorem} {\bf 5.2.1.} \ Let $M$ be a compact
three-manifold with positive Ricci curvature and let the metric
evolve by the Ricci flow. By Lemma 5.2.4 we know that there exists
a positive constant $\beta>0$ such that
$$
R_{ij}\ge \beta Rg_{ij}
$$
for all $t\ge 0$ as long as the solution exists. The scalar
curvature evolves by
\begin{align*}
\frac{\partial R}{\partial t}&  = \Delta R+2|R_{ij}|^2\\
&  \ge \Delta R+\frac{2}{3}R^2,
\end{align*}
which implies, by the maximum principle, that the scalar curvature
remains positive and tends to $+\infty$ in finite time.

We now use a blow up argument as in Section 4.3 to get the
following gradient estimate.

\medskip
{\bf Claim.} \ For any $\varepsilon>0$, there exists a positive
constant $C_\varepsilon<+\infty$ such that for any time $\tau\ge
0$, we have
$$
\max_{t\le\tau}\max_{x\in M}|\nabla Rm(x,t)| \le \varepsilon
\max_{t\le\tau}\max_{x\in M}|Rm(x,t)|^\frac{3}{2}+C_\varepsilon.
$$

\medskip
We argue by contradiction. Suppose the above gradient estimate
fails for some fixed $\varepsilon_0>0$. Pick a sequence
$C_j\rightarrow+\infty$, and pick points $x_j\in M$ and times
$\tau_j$ such that
$$
|\nabla Rm(x_j,\tau_j)|\ge\varepsilon_0\max_{t\le\tau_j}
\max_{x\in M}|Rm(x,t)|^\frac{3}{2}+C_j,\quad j=1,2,\ldots.
$$
Choose $x_j$ to be the origin, and pull the metric back to a small
ball on the tangent space $T_{x_j}M$ of radius $r_j$ proportional
to the reciprocal of the square root of the maximum curvature up
to time $\tau_j$ (i.e., $\max_{t\le\tau_j}\max_{x\in M}|Rm(x,$
$t)|$). Clearly the maximum curvatures go to infinity by Shi's
derivative estimate of curvature (Theorem 1.4.1). Dilate the
metrics so that the maximum curvature
$$
\max_{t\le\tau_j}\max_{x\in M}|Rm(x,t)|
$$
becomes 1 and translate time so that $\tau_j$ becomes the time 0.
By Theorem 4.1.5, we can take a (local) limit. The limit metric
satisfies
$$
|\nabla Rm(0,0)|\ge \varepsilon_0>0.
$$
However the pinching estimate in Proposition 5.2.5 tells us the
limit metric has
$$
R_{ij}-\frac{1}{3}Rg_{ij}\equiv 0.
$$
By using the contracted second Bianchi identity
$$
\frac{1}{2}\nabla_iR =\nabla^jR_{ij}
=\nabla^j\(R_{ij}-\frac{1}{3}Rg_{ij}\)+\frac{1}{3}\nabla_iR,
$$
we get
$$
\nabla_iR\equiv0\quad \text{and then}\quad \nabla_iR_{jk}\equiv0.
$$
For a three-manifold, this in turn implies
$$
\nabla Rm=0
$$
which is a contradiction. Hence we have proved the gradient
estimate claimed.

We can now show that the solution to the Ricci flow becomes round
as the time $t$ tends to the maximal time $T$. We have seen that
the scalar curvature goes to infinity in finite time. Pick a
sequence of points $x_j\in M$ and times $\tau_j$ where the
curvature at $x_j$ is as large as it has been anywhere for $0\le
t\le\tau_j$ and $\tau_j$ tends to the maximal time. Since $|\nabla
Rm|$ is very small compared to $|Rm(x_j,\tau_j)|$ by the above
gradient estimate and $|R_{ij}-\frac{1}{3}Rg_{ij}|$ is also very
small compared to $|Rm(x_j,\tau_j)|$ by Proposition 5.2.5, the
curvature is nearly constant and positive in a large ball around
$x_j$ at the time $\tau_j$. But then the Bonnet-Myers' theorem
tells us this is the whole manifold. For $j$ large enough, the
sectional curvature of the solution at the time $\tau_j$ is
sufficiently pinched. Then it follows from the Klingenberg
injectivity radius estimate (see Section 4.2) that the injectivity
radius of the metric at time $\tau_j$ is bounded from below by
$c/\sqrt{|Rm(x_j,\tau_j)|}$ for some positive constant $c$
independent of $j$. Dilate the metrics so that the maximum
curvature $|Rm(x_j,\tau_j)|=\max_{t\leq\tau_j}\max_{x \in M}
|Rm(x,t)|$ becomes $1$ and shift the time $\tau_j$ to the new time
$0$. Then we can apply Hamilton's compactness theorem (Theorem
4.1.5) to take a limit. By the pinching estimate in Proposition
5.2.5, we know that the limit has positive constant curvature
which is either the round $\mathbb{S}^3$ or a metric quotient of
the round $\mathbb{S}^3$. Consequently, the compact three-manifold
$M$ is diffeomorphic to the round $\mathbb{S}^3$ or a metric
quotient of the round $\mathbb{S}^3$.
\endproof

Next we consider the pinching estimates of the Ricci flow on a
compact four-manifold $M$ with positive curvature operator. The
arguments are taken Hamilton \cite{Ha86}.

In dimension 4, we saw in Section 1.3 when we decompose
orthogonally $\Lambda^2=\Lambda^2_+\oplus\Lambda^2_-$ into the
eigenspaces of Hodge star with eigenvalue $\pm1$, we have a block
decomposition of $M_{\alpha\beta}$ as
$$
M_{\alpha\beta}=\left(\begin{array}{cc} A&  B\\ {}^tB&
C\end{array}\right)
$$
and then
$$
M_{\alpha\beta}^\#=2\left(\begin{array}{cc} A^\#&  B^\#\\
{}^tB^\#&  C^\#\end{array}\right)
$$
where $A^\#$, $B^\#$, $C^\#$ are the adjoints of $3\times3$
submatrices as before.

Thus the ODE
$$
\frac{d}{dt}M_{\alpha\beta}=M_{\alpha\beta}^2+M_{\alpha\beta}^\#
$$
corresponding to the PDE (5.2.1) breaks up into the system of
three equations \be
   \left\{
   \begin{array}{lll}
    \frac {d }{d t}A=A^2+B {^tB}+2A^\#,
         \\[3mm]
    \frac {d }{d t}B=AB+BC+2B^\#,
         \\[3mm]
    \frac {d }{d t}C=C^2+ {^tB}B+2C^\#.
\end{array}
\right. 
\ee As shown in Section 1.3, by the Bianchi identity, we know that
$\tr A=\tr C.$ For the symmetric matrices $A$ and $C$, we can
choose an orthonormal basis $x_1,x_2,x_3$ of $\Lambda^2_+$ such
that
$$
 A=
{ \left(
  \begin{array}{ccc}
    a_1&   0 &   0\\[4mm]
    0 &   a_2 &   0\\[4mm]
    0 &   0 &   a_3
\end{array}
\right),}
$$
and an orthonormal basis $z_1,z_2,z_3$ of $\Lambda^2_-$ such that
$$
 C=
{ \left(
  \begin{array}{ccc}
    c_1&   0 &   0\\[4mm]
    0 &   c_2 &   0\\[4mm]
    0 &   0 &   c_3
\end{array}
\right).}
$$
For matrix $B$, we can choose orthonormal basis
$y_1^+,y_2^+,y_3^+$ of $\Lambda^2_+$ and $y_1^-,y_2^-,$ $y_3^-$ of
$\Lambda^2_-$ such that
$$
 B=
{ \left(
  \begin{array}{ccc}
    b_1&   0 &   0\\[4mm]
    0 &   b_2 &   0\\[4mm]
    0 &   0 &   b_3
\end{array}
\right).}
$$
with $0\leq b_1\leq b_2\leq b_3$. We may also arrange the
eigenvalues of $A$ and $C$ as $a_1\leq a_2\leq a_3$ and $c_1\leq
c_2\leq c_3$. In view of the advanced maximum principle Theorem
2.3.1, we only need to establish the pinching estimates for the
ODE (5.2.3).

Note that
\begin{align*}
a_1&=\inf\{A(x,x)\ |\  x\in \Lambda_+^2 \ \ \text{and} \ \ |x|=1\},\\
a_3&=\sup\{A(x,x)\ |\  x\in \Lambda_+^2 \ \ \text{and} \ \ |x|=1\},\\
c_1&=\inf\{C(z,z)\ |\  z\in \Lambda_-^2 \ \ \text{and} \ \ |z|=1\},\\
c_3&=\sup\{C(z,z)\ |\  z\in \Lambda_-^2 \ \ \text{and} \ \
|z|=1\}.
\end{align*}
We can compute their derivatives by Lemma 2.3.3 as follows: \be
   \left\{
   \begin{array}{lll}
    \frac {d }{d t}a_1\geq a_1^2+b_1^2+2a_2a_3,
         \\[3mm]
    \frac {d }{d t}a_3\leq a_3^2+b_3^2+2a_1a_2,
         \\[3mm]
    \frac {d }{d t}c_1\geq c_1^2+b_1^2+2c_2c_3,
         \\[3mm]
    \frac {d }{d t}c_3\leq c_3^2+b_3^2+2c_1c_2.
\end{array}
\right. 
\ee We shall make the pinching estimates by using the functions
$b_2+b_3$ and $a-2b+c$, where $a=a_1+a_2+a_3=c=c_1+c_2+c_3$ and
$b=b_1+b_2+b_3$. Since
\begin{align*}
b_2+b_3&  =B(y_2^+,y_2^-)+B(y_3^+,y_3^-) \\
&  = \sup\{B(y^+,y^-)+B(\tilde{y}^+,\tilde{y}^-)\, | \,
y^+,\tilde{y}^+\in \Lambda_+^2 \ \ \text{with} \ \
|y^+|=|\tilde{y}^+|=1,\\
&\qquad y^+ \bot \tilde{y}^+, \;\text{ and } \; \
y^-,\tilde{y}^-\in \Lambda_-^2 \;\text{with } \;
|y^-|=|\tilde{y}^-|=1,  y^- \bot \tilde{y}^- \},
\end{align*}
We compute by Lemma 2.3.3,
\begin{align}
\frac{d}{dt}(b_2+b_3)& \leq
\frac{d}{dt}B(y_2^+,y_2^-)+\frac{d}{dt}B(y_3^+,y_3^-)\\
&  = AB(y_2^+,y_2^-)+BC(y_2^+,y_2^-)+2B^\#(y_2^+,y_2^-) \nn\\
&\quad +AB(y_3^+,y_3^-)+BC(y_3^+,y_3^-)+2B^\#(y_3^+,y_3^-) \nn\\
&  = b_2A(y_2^+,y_2^+)+b_2C(y_2^-,y_2^-)+2b_1b_3 \nn\\
&\quad +b_3A(y_3^+,y_3^+)+b_3C(y_3^-,y_3^-)+2b_1b_2 \nn\\
&  \leq a_2b_2+a_3b_3+b_2c_2+b_3c_3+2b_1b_2+2b_1b_3,\nn
\end{align} 
where we used the facts that $A(y_2^+,y_2^+)+A(y_3^+,y_3^+)\leq
a_2+a_3$ and $C(y_2^-,y_2^-)+C(y_3^-,y_3^-)\leq c_2+c_3$.

Note also that the function $a=trA=c=trC$ is linear, and the
function $b$ is given by
\begin{align*}
b&  = B(y_1^+,y_1^-)+B(y_2^+,y_2^-)+B(y_3^+,y_3^-) \\
& =\sup\Big\{B(Ty_1^+,\tilde{T}y_1^-)+B(Ty_2^+,\tilde{T}y_2^-)
+B(Ty_3^+,\tilde{T}y_3^-)\ |
\ T,\tilde{T} \; \text{ are} \\
&\qquad \text{othogonal transformations of}\ \Lambda^2_+ \; \text{
and }\; \Lambda^2_-\; \text{ respectively}\Big\}.
\end{align*}
Indeed,
\begin{align*}
& B(Ty_1^+,\tilde{T}y_1^-)+B(Ty_2^+,\tilde{T}y_2^-)
+B(Ty_3^+,\tilde{T}y_3^-)\\
&  = B(y_1^+,T^{-1}\tilde{T}(y_1^-))
+B(y_2^+,T^{-1}\tilde{T}(y_2^-))+B(y_3^+,T^{-1}\tilde{T}(y_3^-))\\
&  = b_1t_{11}+b_2t_{22}+b_3t_{33}
\end{align*}
where $t_{11},t_{22},t_{33}$ are diagonal elements of the
orthogonal matrix $T^{-1}\tilde{T}$ with $t_{11},t_{22},t_{33}\leq
1$. Thus by using Lemma 2.3.3 again, we compute
\begin{align*}
\frac{d}{dt}(a-2b+c)&  \geq
{\rm tr}\,\(\frac{d}{dt}A-2\frac{d}{dt}B+\frac{d}{dt}C\)\\
&  ={\rm tr}\,((A-B)^2+(C-B)^2+2(A^\#-2B^\#+C^\#))
\end{align*}
evaluated in those coordinates where $B$ is diagonal as above.
Recalling the definition of Lie algebra product
$$
P\#Q=\frac{1}{2}\varepsilon_{\alpha\beta\gamma}
\varepsilon_{\zeta\eta\theta}P_{\beta\eta}Q_{\gamma\theta}
$$
with $\varepsilon_{\alpha\beta\gamma}$ being the permutation
tensor, we see that the Lie algebra product $\#$ gives a symmetric
bilinear operation on matrices, and then
\begin{align*}
& {\rm tr}\,(2(A^\#-2B^\#+C^\#)) \\
&  = {\rm tr}\,((A-C)^\#+(A+2B+C)\#(A-2B+C))\\
&  = -\frac{1}{2}{\rm tr}\,(A-C)^2+\frac{1}{2}({\rm tr}\,(A-C))^2\\
&\quad + {\rm tr}\,((A+2B+C)\#(A-2B+C))\\
&  = -\frac{1}{2}{\rm tr}\,(A-C)^2+{\rm tr}\,((A+2B+C)\#(A-2B+C))
\end{align*}
by the Bianchi identity. It is easy to check that
$$
{\rm tr}\,(A-B)^2+\tr(C-B)^2-\frac{1}{2}\tr(A-C)^2
=\frac{1}{2}\tr(A-2B+C)^2\geq 0.
$$
Thus we obtain
$$
\frac{d}{dt}(a-2b+c)\geq \tr((A+2B+C)\#(A-2B+C))
$$
Since $M_{\alpha\beta}\geq 0$ and
$$
M_{\alpha\beta}= { \left(
  \begin{array}{cc}
    A &   B\\[4mm]
    ^tB &   C
   \end{array}
\right),}
$$
we see that $A+2B+C\geq 0$ and $A-2B+C\geq 0$, by applying
$M_{\alpha\beta}$ to the vectors $(x,x)$ and $(x,-x)$. It is then
not hard to see
$$
\tr((A+2B+C)\#(A-2B+C))\geq (a_1+2b_1+c_1)(a-2b+c).
$$
Hence we obtain \be
\frac{d}{dt}(a-2b+c)\geq (a_1+2b_1+c_1)(a-2b+c). 
\ee

We now state and prove the following pinching estimates of
Hamilton for the associated ODE (5.2.3).

\begin{proposition}[{Hamilton \cite{Ha86}}]
If we choose successively positive constants $G$ large enough, $H$
large enough, $\delta$ small enough, $J$ large enough,
$\varepsilon$ small enough, $K$ large enough, $\theta$ small
enough, and $L$ large enough, with each depending on those chosen
before, then the closed convex subset $X$ of
$\{M_{\alpha\beta}\geq 0 \}$ defined by the inequalities
\begin{itemize}
\item[(1)] $(b_2+b_3)^2\leq Ga_1c_1,$ \item[(2)] $a_3\leq Ha_1 \ \
\text{and}\ \ c_3\leq Hc_1,$ \item[(3)] $(b_2+b_3)^{2+\delta}\leq
Ja_1c_1(a-2b+c)^{\delta},$ \item[(4)]
$(b_2+b_3)^{2+\varepsilon}\leq Ka_1c_1,$ \item[(5)] $a_3\leq
a_1+La_1^{1-\theta}\ \ \text{and} \ \ c_3\leq
c_1+Lc_1^{1-\theta},$
\end{itemize}
is preserved by {\rm ODE} $(5.2.3)$. Moreover every compact subset
of $\{M_{\alpha\beta}>0\}$ lies in some such set $X$.
\end{proposition}

\begin{pf}
Clearly the subset $X$ is closed and convex.  We first note that
we may assume $b_2+b_3>0$ because if $b_2+b_3=0$, then from
(5.2.5), $b_2+b_3$ will remain zero and then the inequalities (1),
(3) and (4) concerning $b_2+b_3$ are automatically satisfied.
Likewise we may assume $a_3>0$ and $c_3>0$ from (5.2.4).

Let $G$ be a fixed positive constant. To prove the inequality (1)
we only need to check \be
\frac{d}{dt}\log\frac{a_1c_1}{(b_2+b_3)^2}\geq 0 
\ee whenever $(b_2+b_3)^2=Ga_1c_1$ and $b_2+b_3>0$. Indeed, it
follows from (5.2.4) and (5.2.5) that
\begin{align}
\frac{d}{dt}\log a_1 &\geq
2b_1+2a_3+\frac{(a_1-b_1)^2}{a_1}+2\frac{a_3}{a_1}(a_2-a_1), \\
\frac{d}{dt}\log c_1 &\geq
2b_1+2c_3+\frac{(c_1-b_1)^2}{c_1}+2\frac{c_3}{c_1}(c_2-c_1), \\
\intertext{and} \frac{d}{dt}\log(b_2+b_3)&\leq
2b_1+a_3+c_3-\frac{b_2}{b_2+b_3}[(a_3-a_2)+(c_3-c_2)],
\end{align}
which immediately give the desired inequality (5.2.7).

By (5.2.4), we have \be \frac{d}{dt}\log a_3\leq
a_3+2a_1+\frac{b_3^2}{a_3}-\frac{2a_1}{a_3}(a_3-a_2).
\ee {}From the inequality (1) there holds $b_3^2\leq Ga_1c_1$.
Since $\tr A=\tr C$, $c_1\leq c_1+c_2+c_3=a_1+a_2+a_3\leq 3a_3$
\noindent which shows
$$
\frac{b_3^2}{a_3}\leq 3Ga_1.
$$
Thus by (5.2.8) and (5.2.11),
$$
\frac{d}{dt}\log\frac{a_3}{a_1}\leq (3G+2)a_1-a_3.
$$
So if $H\geq (3G+2)$, then the inequalities $a_3\leq Ha_1$ and
likewise $c_3\leq Hc_1$ are preserved.

For the inequality (3), we compute from (5.2.8)-(5.2.10)
\begin{align*}
\frac{d}{dt}\log\frac{a_1c_1}{(b_2+b_3)^2}&  \geq
\frac{(a_1-b_1)^2}{a_1}+\frac{(c_1-b_1)^2}{c_1}
+2\frac{a_3}{a_1}(a_2-a_1)+2\frac{c_3}{c_1}(c_2-c_1)\\
&\quad +\frac{2b_2}{b_2+b_3}[(a_3-a_2)+(c_3-c_2)].
\end{align*}
If $b_1\leq a_1/2$, then
$$
\frac{(a_1-b_1)^2}{a_1}\geq \frac{a_1}{4}\geq\frac{1}{4H}a_3,
$$
and if $b_1\geq a_1/2$, then
$$
\frac{2b_2}{b_2+b_3}\geq \frac{2b_2}{\sqrt{Ga_1c_1}}\geq
\frac{2b_2}{\sqrt{3Ga_1a_3}}\geq \frac{2b_2}{\sqrt{3GH}\cdot a_1}
\geq \frac{1}{\sqrt{3GH}}.
$$
Thus by combining with $ 3a_3 \geq c_3$, we have
$$
\frac{d}{dt}\log\frac{a_1c_1}{(b_2+b_3)^2}\geq \delta
(a_3-a_1)+\delta(c_3-c_1)
$$
provided $\delta\leq \min(\frac{1}{24H},\frac{1}{\sqrt{3GH}})$. On
the other hand, it follows from (5.2.6) and (5.2.10) that
$$
\frac{d}{dt}\log\frac{b_2+b_3}{a-2b+c}\leq(a_3-a_1)+(c_3-c_1).
$$
Therefore the inequality (3)
$$
(b_2+b_3)^{2+\delta}\le Ja_1c_1(a-2b+c)^\delta
$$
will be preserved by any positive constant $J$.

To verify the inequality (4), we first note that there is a small
$\eta>0$ such that
$$
b\le (1-\eta)a,
$$
on the set defined by the inequality (3). Indeed, if $b\le
\frac{a}{2}$, this is trivial and if $b\ge\frac{a}{2}$, then
$$
\(\frac{a}{3}\)^{2+\delta} \le(b_2+b_3)^{2+\delta}\le 2^\delta
Ja^2(a-b)^\delta
$$
which makes $b\le(1-\eta)a$ for some $\eta>0$ small enough.
Consequently, $$\eta a\le a-b\le 3(a_3-b_1)$$ which implies either
$$
a_3-a_1\ge \frac{1}{6}\eta a,
$$
or
$$
a_1-b_1\ge \frac{1}{6}\eta a.
$$
Thus as in the proof of the inequality (3), we have
$$
\frac{d}{dt}\log\frac{a_1c_1}{(b_2+b_3)^2}\ge \delta(a_3-a_1)
$$
and
$$
\frac{d}{dt}\log\frac{a_1c_1}{(b_2+b_3)^2}\ge
\frac{(a_1-b_1)^2}{a_1},
$$
which in turn implies
$$
\frac{d}{dt}\log\frac{a_1c_1}{(b_2+b_3)^2} \ge
\(\max\left\{\frac{1}{6}\eta\delta,\;
\frac{1}{36}\eta^2\right\}\)\cdot a.
$$
On the other hand, it follows from (5.2.10) that
$$
\frac{d}{dt}\log(b_2+b_3)\le 2b_1+a_3+c_3\le 4a
$$
since $M_{\alpha\beta}\ge 0$. Then if $\varepsilon>0$ is small
enough
$$
\frac{d}{dt}\log\frac{a_1c_1}{(b_2+b_3)^{2+\varepsilon}}\ge 0
$$
and it follows that the inequality (4) is preserved by any
positive $K$.

Finally we consider the inequality (5). From (5.2.8) we have
$$
\frac{d}{dt}\log a_1\ge a_1+2a_3
$$
and then for $\theta\in (0,1)$,
$$
\frac{d}{dt}\log (a_1+La_1^{1-\theta})\ge
\frac{a_1+(1-\theta)La_1^{1-\theta}}{a_1+La_1^{1-\theta}}(a_1+2a_3).
$$
On the other hand, the inequality (4) tells us
$$
b_3^2\le \tilde{K}a_1^{1-\theta}a_3
$$
for some positive constant $\tilde{K}$ large enough with $\theta$
to be fixed small enough. And then
$$
\frac{d}{dt}\log a_3\le a_3+2a_1+\tilde{K}a_1^{1-\theta},
$$
by combining with (5.2.11). Thus by choosing
$\theta\le\frac{1}{6H}$ and $L\ge 2\tilde{K}$,
\begin{align*}
\frac{d}{dt}\log \frac{a_1+La_1^{1-\theta}}{a_3}& \ge
(a_3-a_1)-\theta\frac{La_1^{1-\theta}}{a_1+La_1^{1-\theta}}(a_1+2a_3)
-\tilde{K}a_1^{1-\theta}\\
&
\ge(a_3-a_1)-\theta\frac{La_1^{1-\theta}}{a_1+La_1^{1-\theta}}\cdot
3Ha_1-\tilde{K}a_1^{1-\theta}\\
&  \ge (a_3-a_1)-(3\theta HL+\tilde{K})a_1^{1-\theta}\\
&  = [L-(3\theta HL+\tilde{K})]a_1^{1-\theta}\\
&  \ge 0
\end{align*}
whenever $a_1+La_1^{1-\theta}=a_3$. Consequently the set
$\{a_1+La_1^{1-\theta}\ge a_3\}$ is preserved. A similar argument
works for the inequality in $C$. This completes the proof of
Proposition 5.2.6.
\end{pf}

The combination of the advanced maximum principle Theorem 2.3.1
and the pinching estimates of the ODE (5.2.3) in Proposition 5.2.6
immediately gives the following pinching estimate for the Ricci
flow on a compact four-manifold.

\begin{corollary} [{Hamilton \cite{Ha86}}]
Suppose that the initial metric of the solution to the Ricci flow
on a compact four-manifold has positive curvature operator. Then
for any $\varepsilon>0$ we can find positive constant
$C_\varepsilon<+\infty$ such that
$$
|{\mathop{Rm}^ \circ}|\le\varepsilon R+C_\varepsilon
$$
for all $t\ge 0$ as long as the solution exists, where ${\mathop
{Rm}\limits^ \circ}$ is the traceless part of the curvature
operator.
\end{corollary}

\smallskip
{\bf\em Proof of Theorem} {\bf 5.2.2.} \ Let $M$ be a compact
four-manifold with positive curvature operator and let us evolve
the metric by the Ricci flow. Again the evolution equation of the
scalar curvature tells us that the scalar curvature remains
positive and becomes unbounded in finite time. Pick a sequence of
points $x_j\in M$ and times $\tau_j$ where the curvature at $x_j$
is as large as it has been anywhere for $0\le t\le \tau_j$. Dilate
the metrics so that the maximum curvature
$|Rm(x_j,\tau_j)|=\max_{t\le\tau_j}\max_{x\in M}|Rm(x,t)|$ becomes
1 and shift the time so that the time $\tau_j$ becomes the new
time 0. The Klingenberg injectivity radius estimate in Section 4.2
tells us that the injectivity radii of the rescaled metrics at the
origins $x_j$ and at the new time 0 are uniformly bounded from
below. Then we can apply the Hamilton's compactness theorem
(Theorem 4.1.5) to take a limit. By the pinching estimate in
Corollary 5.2.7, we know that the limit metric has positive
constant curvature which is either $\mathbb{S}^4$ or
$\mathbb{R}\mathbb{P}^4$. Therefore the compact four-manifold $M$
is diffeomorphic to the sphere $\mathbb{S}^4$ or the real
projective space $\mathbb{R}\mathbb{P}^4$.
\endproof

\begin{remark}
The proofs of Theorem 5.2.1 and Theorem 5.2.2 also show that the
Ricci flow on a compact three-manifold with positive Ricci
curvature or a compact four-manifold with positive curvature
operator is subsequentially converging (up to scalings) in the
$C^{\infty}$ topology to the same underlying compact manifold with
a metric of positive constant curvature. Of course, this
subsequential convergence is in the sense of Hamilton's
compactness theorem (Theorem 4.1.5) which is also up to the
pullbacks of diffeomorphisms. Actually in \cite{Ha82} and
\cite{Ha86}, Hamilton obtained the convergence in the stronger
sense that the (rescaled) metrics converge (in the $C^{\infty}$
topology) to a constant (positive) curvature metric.
\end{remark}

In the following we use Hamilton's strong maximum principle
(Theorem 2.2.1) to prove Theorem 5.2.3.

\medskip
{\bf \em Proof of Theorem} {\bf 5.2.3.} \ In views of Theorem
5.2.1 and Theorem 5.2.2, we may assume the Ricci curvature (in
dimension 3) and the curvature operator (in dimension 4) always
have nontrivial kernels somewhere along the Ricci flow.

\smallskip
(i) In the case of dimension 3, we consider the evolution equation
(1.3.5) of the Ricci curvature
$$
\frac{\partial R_{ab}}{\partial t} =\triangle
R_{ab}+2R_{acbd}R_{cd}
$$
in an orthonormal frame coordinate. At each point, we diagonalize
$R_{ab}$ with eigenvectors $e_{1}, e_{2}, e_{3}$ and eigenvalues
$\lambda_1\leq\lambda_2\leq\lambda_3$. Since
\begin{align*}
R_{1c1d}R_{cd}&  = R_{1212}R_{22}+R_{1313}R_{33}\\
&  =\frac{1}{2}((\lambda_{3}-\lambda_{2})^{2}
+\lambda_{1}(\lambda_2+\lambda_3)),
\end{align*}
we know that if $R_{ab}\geq 0$, then $R_{acbd}R_{cd}\geq0$. By
Hamilton's strong maximum principle (Theorem 2.2.1), there exists
an interval $0<t<\delta$ on which the rank of $R_{ab}$ is constant
and the null space of $R_{ab}$ is invariant under parallel
translation and invariant in time and also lies in the null space
of $R_{acbd}R_{cd}$. If the null space of $R_{ab}$ has rank one,
then $\lambda_1=0$ and $\lambda_2=\lambda_{3}>0$.  In this case,
by De Rham decomposition theorem, the universal cover $\tilde{M}$
of the compact $M$ splits isometrically as $\mathbb{R}\times
\Sigma^{2}$ and the curvature of $\Sigma^{2}$ has a positive lower
bound. Hence $\Sigma^{2}$ is diffeomorphic to $\mathbb{S}^{2}$.
Assume $M=\mathbb{R}\times \Sigma^{2}/\Gamma$, for some isometric
subgroup $\Gamma$ of $\mathbb{R}\times \Sigma^{2}$. Note that
$\Gamma$ remains to be an isometric subgroup of $\mathbb{R}\times
\Sigma^{2}$ during the Ricci flow by the uniqueness (Theorem
1.2.4). Since the Ricci flow on $\mathbb{R}\times
\Sigma^{2}/\Gamma$ converges to the standard metric by Theorem
5.1.11, $\Gamma$ must be an isometric subgroup of
$\mathbb{R}\times\mathbb{S}^{2}$ in the standard metric. If the
null space of $R_{ab}$ has rank greater than one, then $R_{ab}=0$
and the manifold is flat. This proves Theorem 5.2.3 part (i).

\smallskip
(ii) In the case of dimension 4, we classify the manifolds
according to the (restricted) holonomy algebra $\mathcal{G}$. Note
that the curvature operator has nontrivial kernel and
$\mathcal{G}$ is the image of the the curvature operator, we see
that $\mathcal{G}$ is a proper subalgebra of $so(4)$. We divide
the argument into two cases.

\smallskip
{\it Case} 1. $\mathcal{G}$ is reducible.

\smallskip
In this case the universal cover $\tilde{M}$ splits isometrically
as $\tilde{M_{1}}\times \tilde{M_{2}}$.  By the above results on
two and three dimensional Ricci flow, we see that ${M}$ is
diffeomorphic to a quotient of one of the spaces $\mathbb{R}^{4},$
$\mathbb{R}\times \mathbb{S}^{3}$, $\mathbb{R}^{2}\times
\mathbb{S}^{2}$, $\mathbb{S}^{2}\times \mathbb{S}^{2}$ by a group
of fixed point free isometries. As before by running the Ricci
flow until it converges and using the uniqueness (Theorem 1.2.4),
we see that this group is actually a subgroup of the isometries in
the standard metrics.

\smallskip
{\it Case} 2. $\mathcal{G}$ is not reducible (i.e., irreducible).

\smallskip
If the manifold is not Einstein, then by Berger's classification
theorem for holonomy groups \cite{Ber}, $\mathcal{G}=so(4)$ or
$u(2)$. Since the curvature operator is not strictly positive,
$\mathcal{G}=u(2)$, and the universal cover $\tilde{M}$ of $M$ is
K$\ddot{a}$hler and has positive bisectional curvature. In this
case $\tilde{M}$ is biholomorphic to $\mathbb{C}\mathbb{P}^{2}$ by
the result of Andreotti-Frankel \cite{Ft} (also cf. Mori
\cite{Mori} and Siu-Yau \cite{SiY}).

If the manifold is Einstein, then by the block decomposition of
the curvature operator matrix in four-manifolds (see the third
section of Chapter 1),
$$
Rm(\Lambda^{2}_{+},\Lambda^{2}_{-})=0.
$$
Let $\varphi\neq0$, and
$$
\varphi=\varphi_{+}+\varphi_{-}\in
  \Lambda^{2}_{+}\oplus\Lambda^{2}_{-},
$$
lies in the kernel of the curvature operator, then
$$
   0=Rm(\varphi_{+},\varphi_{+})+Rm(\varphi_{-},\varphi_{-}).
$$
It follows that
$$
   Rm(\varphi_{+},\varphi_{+})=0, \mbox{ and }
   Rm(\varphi_{-},\varphi_{-})=0.
$$
We may assume $\varphi_{+}\neq 0$ (the argument for the other case
is similar).  We consider the restriction of $Rm$ to
$\Lambda^{2}_{+}$, since $\Lambda^{2}_{+}$ is an invariant
subspace of $R{m}$ and the intersection of $\Lambda^{2}_{+}$ with
the null space of $Rm$ is nontrivial. By considering the null
space of $Rm$ and its orthogonal complement in $\Lambda^{2}_{+}$,
we obtain a parallel distribution of rank one in
$\Lambda^{2}_{+}$. This parallel distribution gives a parallel
translation invariant two-form $\omega \in \Lambda^{2}_{+}$ on the
universal cover $\tilde{M}$ of $M$. This two-form is
nondegenerate, so it induces a K$\ddot{a}$hler structure of
$\tilde{M}$. Since the K$\ddot{a}$hler metric is parallel with
respect to the original metric and the manifold is irreducible,
the K$\ddot{a}$hler metric is proportional to the original metric.
Hence the manifold $\tilde{M}$ is K$\ddot{a}$hler-Einstein with
nonnegative curvature operator. Taking into account the
irreducibility of $\mathcal{G}$, it follows that $\tilde{M}$ is
biholomorphic to $\mathbb{C}\mathbb{P}^{2}$.  Therefore the proof
of Theorem 5.2.3 is completed.
\endproof

To end this section, we mention some generalizations of Hamilton's
differential sphere theorem (Theorem 5.2.1 and Theorem 5.2.2) to
higher dimensions.

It is well-known that the curvature tensor $Rm=\{R_{ijkl}\}$ of a
Riemannian manifold can be decomposed into three orthogonal
components which have the same symmetries as $Rm$:
$$
Rm=W+V+U.
$$
Here $W=\{W_{ijkl}\}$ is the Weyl conformal curvature tensor,
whereas $V=\{V_{ijkl}\}$ and $U=\{U_{ijkl}\}$ denote the traceless
Ricci part and the scalar curvature part respectively. The
following pointwisely pinching sphere theorem under the additional
assumption that the manifold is compact was first obtained by
Huisken \cite{Hu}, Margerin \cite{Mar94}, \cite{Mar98} and
Nishikawa \cite{Ns} by using the Ricci flow. The compactness
assumption was later removed by Chen and the second author in
\cite{CZ00}.

\begin{theorem}
Let $n\ge4$. Suppose $M$ is a complete n-dimensional manifold with
positive and bounded scalar curvature and satisfies the
pointwisely pinching condition
$$
|W|^2+|V|^2\le\delta_n(1-\varepsilon)^2|U|^2,
$$
where $\varepsilon>0,\delta_4=\frac{1}{5},\delta_5=\frac{1}{10}$,
and
$$
\delta_n=\frac{2}{(n-2)(n+1)},n\ge 6.
$$
Then $M$ is diffeomorphic to the sphere $\mathbb{S}^n$ or a
quotient of it by a finite group of fixed point free isometries in
the standard metric.
\end{theorem}

Also, using the minimal surface theory, Micallef and Moore
\cite{MiMo} proved any compact simply connected $n$-dimensional
manifold with positive curvature operator is
homeomorphic\footnote{Very recently, B$\ddot{o}$hm and Wilking
\cite{BW} have proved, by using the Ricci flow, that a compact
simply connected $n$-dimensional manifold with positive (or
2-positive) curvature operator is {\it diffeomorphic} to
$\mathbb{S}^n$. This gives an affirmative answer to a
long-standing conjecture of Hamilton.} to the $n$-sphere
$\mathbb{S}^n$.

Finally, in \cite{CZ00}, Chen and the second author also used the
Ricci flow to obtain the following flatness theorem for noncompact
three-manifolds.

\begin{theorem}
Let $M$ be a three-dimensional complete noncompact Riemannian
manifold with bounded and nonnegative sectional curvature. Suppose
$M$ satisfies the following Ricci pinching condition
$$
R_{ij}\ge \varepsilon Rg_{ij},\quad \text{on }\; M,
$$
for some $\varepsilon>0$. Then $M$ is flat.
\end{theorem}

\vskip 1cm

\section{Nonsingular Solutions on Three-manifolds}

We have seen in the previous section that a good understanding of
the long time behaviors for solutions to the Ricci flow could lead
to remarkable topological or geometric consequences for the
underlying manifolds. Since one of the central themes of the Ricci
flow is to study the geometry and topology of three-manifolds, we
will start to analyze the long time behavior of the Ricci flow on
a compact three-manifold. Here, we shall first consider a special
class of solutions, the nonsingular solutions (see the definition
below). The main purpose of this section is to present Hamilton's
important result in \cite{Ha99} that any compact three-manifold
admitting a nonsingular solution is geometrizable in the sense of
Thurston\cite{Th82}. Most of the presentation is based on Hamilton
\cite{Ha99}.

Let $M$ be a compact three-manifold. We will consider the
(unnormalized) Ricci flow
$$
\frac{\partial}{\partial t}g_{ij}=-2R_{ij},
$$
and the normalized Ricci flow
$$
\frac{\partial}{\partial t}g_{ij}=\frac{2}{3}rg_{ij}-2R_{ij}
$$
where $r=r(t)$ is the function of the average of the scalar
curvature. Recall that the normalized flow differs from the
unnormalized flow only by rescaling in space and time so that the
total volume $V=\int_Md\mu$ remains constant. As we mentioned
before, in this section we only consider a special class of
solutions that we now define.

\begin{definition}
A {\bf nonsingular solution}\index{solution!nonsingular} of the
Ricci flow is one where the solution of the normalized flow exists
for all time $0\leq t<\infty$, and the curvature remains bounded
$|Rm|\leq C<+\infty$ for all time with some constant $C$
independent of $t$.
\end{definition}

Clearly any solution to the Ricci flow on a compact three-manifold
with nonnegative Ricci curvature is nonsingular. Currently there
are few conditions which guarantee a solution will remain
nonsingular. Nevertheless, the ideas and arguments of Hamilton
\cite{Ha99} as described below is extremely important. One will
see in Chapter 7 that these arguments will be modified to analyze
the long-time behavior of arbitrary solutions, or even the
solutions with surgery, to the Ricci flow on three-manifolds.

We begin with an improvement of Hamilton-Ivey pinching result
\index{Hamilton-Ivey pinching estimate} (Theorem 2.4.1).

\begin{theorem}[{Hamilton \cite{Ha99}}]
Suppose we have a solution to the $($unnormalized$\,)$\; Ricci\;
flow\; on\; a\; three--manifold\; which\; is\; complete\; with
bounded curvature for each $t\geq0$. Assume at $t=0$ the
eigenvalues $\lambda\geq\mu\geq\nu$ of the curvature operator at
each point are bounded below by $\nu\geq-1$. Then at all points
and all times $t\geq0$ we have the pinching estimate
$$
R\geq(-\nu)[\log(-\nu)+\log(1+t)-3]
$$
whenever $\nu<0$.
\end{theorem}

\begin{pf}
As before, we study the ODE system
\begin{equation*}
\begin{cases}{\displaystyle
    \frac{d\lambda}{dt} =\lambda^2+\mu\nu,}\\[3mm]
{\displaystyle \frac{d\mu}{dt} =\mu^2+\lambda\nu,}\\[3mm]
{\displaystyle \frac{d\nu}{dt} =\nu^2+\lambda\mu.}
\end{cases}
\end{equation*}

Consider again the function $$y=f(x)=x(\log x-3)$$ for $e^2\leq
x<+\infty$, which is increasing and convex with range $-e^2\leq
y<+\infty$. Its inverse function $x=f^{-1}(y)$ is increasing and
concave on $-e^2\leq y<+\infty$. For each $t\geq0$, we consider
the set $K(t)$ of $3\times3$ symmetric matrices defined by the
inequalities: \be
\lambda+\mu+\nu \geq-\frac{3}{1+t}, 
\ee and \be
\nu(1+t)+f^{-1}((\lambda+\mu+\nu)(1+t))\geq0, 
\ee which is closed and convex (as we saw in the proof of Theorem
2.4.1). By the assumptions at $t=0$ and the advanced maximum
principle Theorem 2.3.5, we only need to check that the set $K(t)$
is preserved by the ODE system.

Since $R=\lambda+\mu+\nu$, we get from the ODE that
$$
\frac{dR}{dt}\geq\frac{2}{3}R^2\geq\frac{1}{3}R^2
$$
which implies that
$$
R\geq-\frac{3}{1+t},\ \ \ \text{for\ all }\ t\geq0.
$$
Thus the first inequality (5.3.1) is preserved. Note that the
second inequality (5.3.2) is automatically satisfied when $(-\nu)
\leq 3/(1+t)$. Now we compute from the ODE system,
\begin{displaymath}
\begin{split}
\frac{d}{dt}(\frac{R}{(-\nu)}-\log(-\nu))&  =\frac{1}{(-\nu)^2}
\left[(-\nu)\cdot\frac{dR}{dt}-(R+(-\nu))\frac{d(-\nu)}{dt}\right]\\[1mm]
&  =\frac{1}{(-\nu)^2}[(-\nu)^3+(-\nu)\mu^2+\lambda^2((-\nu)+\mu)
-\lambda\mu(\nu-\mu)]\\[1mm]
&  \geq(-\nu)\\[3mm]
&  \geq\frac{3}{(1+t)}\\[1mm]
&  \geq\frac{d}{dt}[\log(1+t)-3]
\end{split}
\end{displaymath}
whenever $R=(-\nu)[\log(-\nu)+\log(1+t)-3]$ and
$(-\nu)\geq3/(1+t)$. Thus the second inequality (5.3.2) is also
preserved under the system of ODE.

Therefore we have proved the theorem.
\end{pf}

Denote by
$$
\hat{\rho}(t)=\max\{\inj(x,g_{ij}(t))\ |\ x\in M \}
$$
where $\inj(x,g_{ij}(t))$ is the injectivity radius of the
manifold $M$ at $x$ with respect to the metric $g_{ij}(t)$.


\begin{definition}
We say a solution to the normalized Ricci flow is {\bf
collapsed}\index{collapsed} if there is a sequence of times
$t_k\rightarrow+\infty$ such that $\hat{\rho}(t_k)\rightarrow0$ as
$k\rightarrow+\infty$.
\end{definition}

When a nonsingular solution of the Ricci flow on $M$ is collapsed,
it follows from the work of Cheeger-Gromov \cite{ChG86, ChG90} or
Cheeger-Gromov-Fukaya \cite{CGF92} that the manifold $M$ has an
$\mathcal{F}$-structure and then its topology is completely
understood. Thus, in the following, we always assume our
nonsingular solutions are not collapsed.

Now suppose that we have a nonsingular solution which does not
collapse. Then for arbitrary sequence of times
$t_j\rightarrow\infty$, we can find a sequence of points $x_j$ and
some $\delta>0$ so that the injectivity radius of $M$ at $x_j$ in
the metric at time $t_j$ is at least $\delta$. Clearly the
Hamilton's compactness theorem (Theorem 4.1.5) also holds for the
normalized Ricci flow. Then by taking the $x_j$ as origins and the
$t_j$ as initial times, we can extract a convergent subsequence.
We call such a limit a {\bf noncollapsing
limit}\index{noncollapsing limit}. Of course the limit has also
finite volume. However the volume of the limit may be smaller than
the original one if the diameter goes to infinity.

The main result of this section is the following theorem of
Hamilton \cite{Ha99}.

\begin{theorem}[{Hamilton \cite{Ha99}}]
Let $g_{ij}(t)$, $0\leq t<+\infty$, be a noncollapsing nonsingular
solution of the normalized Ricci flow on a compact three-manifold
$M$. Then either
\begin{itemize}
\item[(i)] there exist a sequence of times $t_k\rightarrow+\infty$
and a sequence of diffeomorphisms $\varphi_k:\ M\rightarrow M$ so
that the pull-back of the metric $g_{ij}(t_k)$ by $\varphi_k$
converges in the $C^\infty$ topology to a metric on $M$ with
constant sectional curvature; or \item[(ii)] we can find a finite
collection of complete noncompact hyperbolic three-manifolds
$\mathcal{H}_1,\ldots,\mathcal{H}_m$ with finite volume, and for
all $t$ beyond some time $T<+\infty$ we can find compact subsets
$K_1,\ldots,K_m$ of $\mathcal{H}_1,\ldots,\mathcal{H}_m$
respectively obtained by truncating each cusp of the hyperbolic
manifolds along constant mean curvature torus of small area, and
diffeomorphisms $\varphi_l(t)$, $1\leq l\leq m$, of $K_l$ into $M$
so that as long as $t$ sufficiently large, the pull-back of the
solution metric $g_{ij}(t)$ by $\varphi_l(t)$ is as close as to
the hyperbolic metric as we like on the compact sets
$K_1,\ldots,K_m$; and moreover if we call the \textbf{exceptional
part} of $M$ those points where they are not in the image of any
$\varphi_l$, we can take the injectivity radii of the exceptional
part at everywhere as small as we like and the boundary tori of
each $K_l$ are \textbf{incompressible} in the sense that each
$\varphi_l$ injects $\pi_1(\partial K_l)$ into $\pi_1(M)$.
\index{exceptional part} \index{incompressible}
\end{itemize}
\end{theorem}

\begin{remark}
The exceptional part has bounded curvature and arbitrarily small
injectivity radii everywhere as $t$ large enough. Moreover the
boundary of the exceptional part consists of a finite disjoint
union of tori with sufficiently small area and is convex. Then by
the work of Cheeger-Gromov \cite{ChG86}, \cite{ChG90} or
Cheeger-Gromov-Fukaya \cite{CGF92}, there exists an
$\mathcal{F}$-structure on the exceptional part. In particular,
the exceptional part is a graph manifold, which have been
topologically classified. Hence \textbf{any nonsingular solution
to the normalized Ricci flow is geometrizable} in the sense of
Thurston (see the last section of Chapter 7 for details).
\end{remark}

The rest of this section is devoted to the proof of Theorem 5.3.4.
We now present the proof given by Hamilton \cite{Ha99} and will
divide his arguments in \cite{Ha99} into the following three
parts.

\vskip 0.5cm {\bf Part I: Subsequence Convergence} \vskip 0.2cm

According to Lemma 5.1.1, the scalar curvature of the normalized
flow evolves by the equation
\begin{align}
\frac{\partial }{\partial t}R
&  =\Delta R+2|Ric|^2-\frac{2}{3}rR\\
&  = \Delta R+2|\stackrel{\circ}{Ric}|^2+\frac{2}{3}R(R-r)\nn
\end{align}  
where $\stackrel{\circ}{Ric}$ is the traceless part of the Ricci
tensor. As before, we denote by $R_{\min}(t)=\min_{x\in M}R(x,t)$.
It then follows from the maximum principle that \be
\frac{d}{dt}R_{\min}\geq\frac{2}{3}R_{\min}(R_{\min}-r),
\ee which implies that if $R_{\min}\leq0$ it must be
nondecreasing, and if $R_{\min}\geq0$ it cannot go negative again.
We can then divide the noncollapsing solutions of the normalized
Ricci flow into three cases.

\smallskip
{\it Case} (1): $R_{\min}(t)>0$ for some $t>0$;

\smallskip
{\it Case} (2): $R_{\min}(t)\leq0$ for all $t\in[0,+\infty)$ and
$\lim\limits_{t\rightarrow+\infty}R_{\min}(t)=0$;

\smallskip
{\it Case} (3): $R_{\min}(t)\leq0$ for all $t\in[0,+\infty)$ and
$\lim\limits_{t\rightarrow+\infty}R_{\min}(t)<0$.

\smallskip
Let us first consider Case (1). In this case the maximal time
interval $[0,T)$ of the corresponding solution of the unnormalized
flow is finite, since the unnormalized scalar curvature
$\tilde{R}$ satisfies
\begin{displaymath}
\begin{split}
\frac{\partial}{\partial t}\tilde{R}
&  =\Delta\tilde{R}+2|\tilde{R}ic|^2\\
&  \geq\Delta\tilde{R}+\frac{2}{3}\tilde{R}^2
\end{split}
\end{displaymath}
which implies that the curvature of the unnormalized solution
blows up in finite time. Without loss of generality, we may assume
that for the initial metric at $t=0$, the eigenvalues
$\tilde{\lambda}\geq\tilde{\mu}\geq\tilde{\nu}$ of the curvature
operator are bounded below by $\tilde{\nu}\geq-1$. It follows from
Theorem 5.3.2 that the pinching estimate
$$
\tilde{R}\geq(-\tilde{\nu})[\log(-\tilde{\nu})+\log(1+t)-3]
$$
holds whenever $\tilde{\nu}<0$. This shows that when the
unnormalized curvature big, the negative ones are not nearly as
large as the positive ones. Note that the unnormalized curvature
becomes unbounded in finite time. Thus when we rescale the
unnormalized flow to the normalized flow, the scaling factor must
go to infinity. In the nonsingular case the rescaled positive
curvature stay finite, so the rescaled negative curvature (if any)
go to zero. Hence we can take a noncollapsing limit for the
nonsingular solution of the normalized flow so that it has
nonnegative sectional curvature.

Since the volume of the limit is finite, it follows from a result
of Calabi and Yau (cf. \cite{ScY}) that the limit must be compact
and the limiting manifold is the original one. Then by the strong
maximum principle as in the proof of Theorem 5.2.3 (i), either the
limit is flat, or it is a compact metric quotient of the product
of a positively curved surface $\Sigma^2$ with $\mathbb{R}$, or it
has strictly positive curvature. By the work of Schoen-Yau
\cite{ScY79}, a flat three-manifold cannot have a metric of
positive scalar curvature, but our manifold does in Case (1). This
rules out the possibility of a flat limit. Clearly the limit is
also a nonsingular solution to the normalized Ricci flow. Note
that the curvature of the surface $\Sigma^2$ has a positive lower
bound and is compact since it comes from the lifting of the
compact limiting manifold. From Theorem 5.1.11, we see the metric
of the two-dimensional factor $\Sigma^2$ converges to the round
two-sphere $\mathbb{S}^2$ in the normalized Ricci flow. Note also
that the normalized factors in two-dimension and three-dimension
are different. This implies that the compact quotient of the
product $\Sigma^2 \times \mathbb{R}$ cannot be nonsingular, which
is also ruled out for the limit. Thus the limit must have strictly
positive sectional curvature. Since the convergence takes place
everywhere for the compact limit, it follows that as $t$ large
enough the original nonsingular solution has strictly positive
sectional curvature. This in turn shows that the corresponding
unnormalized flow has strictly positive sectional curvature after
some finite time. Then in views of the proof of Theorem 5.2.1, in
particular the pinching estimate in Proposition 5.2.5, the limit
has constant Ricci curvature and then constant sectional curvature
for three-manifolds. This finishes the proof in Case (1).

We next consider Case (2). In this case we only need to show that
we can take a noncollapsing limit which has nonnegative sectional
curvature. Indeed, if this is true, then as in the previous case,
the limit is compact and either it is flat, or it splits as a
product (or a quotient of a product) of a positively curved $S^2$
with a circle $S^1$, or it has strictly positive curvature. But
the assumption $R_{\min}(t)\leq0$ for all times $t\geq0$ in this
case implies the limit must be flat.

Let us consider the corresponding unnormalized flow
$\tilde{g}_{ij}(t)$ associated to the noncollapsing nonsingular
solution. The pinching estimate in Theorem 5.3.2 tells us that we
may assume the unnormalized flow $\tilde{g}_{ij}(t)$ exists for
all times $0\leq t<+\infty$, for otherwise, the scaling factor
approaches infinity as in the previous case which implies the
limit has nonnegative sectional curvature. The volume
$\tilde{V}(t)$ of the unnormalized solution $\tilde{g}_{ij}(t)$
now changes. We divide the discussion into three subcases.

Subcase (2.1): there is a sequence of times
$\tilde{t}_k\rightarrow+\infty$ such that
$\tilde{V}(\tilde{t}_k)\rightarrow+\infty$;

Subcase (2.2): there is a sequence of times
$\tilde{t}_k\rightarrow+\infty$ such that
$\tilde{V}(\tilde{t}_k)\rightarrow0$;

Subcase (2.3): there exist two positive constants $C_1$, $C_2$
such that $C_1\leq\tilde{V}(t)\leq C_2$ for all $0\leq t<+\infty$.

For Subcase (2.1), because
$$
\frac{d\tilde{V}}{dt}=-r\tilde{V}
$$
we have
$$
\log\frac{\tilde{V}(\tilde{t}_k)}{\tilde{V}(0)}
=-\int_0^{\tilde{t}_k}r(t)dt\rightarrow+\infty,\; \text{ as }\;
k\rightarrow+\infty,
$$
which implies that there exists another sequence of times, still
denoted by $\tilde{t}_k$, such that
$\tilde{t}_k\rightarrow+\infty$ and $r(\tilde{t}_k)\leq0$. Let
$t_k$ be the corresponding times for the normalized flow. Thus
there holds for the normalized flow
$$
r(t_k)\rightarrow0,\quad\text{as }\; k\rightarrow\infty,
$$
since $0\geq r(t_k)\geq R_{\min}(t_k)\rightarrow0$ as
$k\rightarrow+\infty$. Then
$$
\int_M(R-R_{\min})d\mu(t_k)=(r(t_k)-R_{\min}(t_k))V\rightarrow0,
\quad \text{as }\; k\rightarrow\infty.
$$
As we take a noncollapsing limit along the time sequence $t_k$, we
get
$$
\int_{M^\infty}Rd\mu^\infty=0
$$
for the limit of the normalized solutions at the new time $t=0$.
But $R\geq0$ for the limit because
$\lim\limits_{t\rightarrow+\infty}R_{\min}(t)=0$ for the
nonsingular solution. So $R=0$ at $t=0$ for the limit. Since the
limit flow exists for $-\infty<t<+\infty$ and the scalar curvature
of the limit flow evolves by
$$
\frac{\partial}{\partial t}R=\Delta
R+2|\Ric|^2-\frac{2}{3}r^\infty R,\quad t\in(-\infty,+\infty)
$$
where $r^\infty$ is the limit of the function $r(t)$ by
translating the times $t_k$ as the new time $t=0$. It follows from
the strong maximum principle that
$$
R\equiv0, \quad \text{on }\; M^\infty\times(-\infty,+\infty).
$$
This in turn implies, in view of the above evolution equation,
that
$$
\Ric\equiv0, \quad \text{on }\; M^\infty\times(-\infty,+\infty).
$$
Hence this limit must be flat. Since the limit $M^\infty$ is
complete and has finite volume, the flat manifold $M^\infty$ must
be compact. Thus the underlying manifold $M^\infty$ must agree
with the original $M$ (as a topological manifold). This says that
the limit was taken on $M$.

For Subcase (2.2), we may assume as before that for the initial
metric at $t=0$ of the unnormalized flow $\tilde{g}_{ij}(t)$, the
eigenvalues $\tilde{\lambda}\geq\tilde{\mu}\geq\tilde{\nu}$ of the
curvature operator satisfy $\tilde{\nu}\geq-1$. It then follows
from Theorem 5.3.2 that
$$
\tilde{R}\geq(-\tilde{\nu})[\log(-\tilde{\nu})+\log(1+t)-3], \quad
\text{for all }\; t\geq0
$$
whenever $\tilde{\nu}<0$.

Let $t_k$ be the sequence of times in the normalized flow which
corresponds to the sequence of times $\tilde{t}_k$. Take a
noncollapsing limit for the normalized flow along the times $t_k$.
Since $\tilde{V}(\tilde{t}_k)\rightarrow 0$, the normalized
curvatures at the times $t_k$ are reduced by multiplying the
factor $(\tilde{V}(\tilde{t}_k))^\frac{2}{3}$. We claim the
noncollapsing limit has nonnegative sectional curvature. Indeed if
the maximum value of $(-\tilde{\nu})$ at the time $\tilde{t}_k$
does not go to infinity, the normalized eigenvalue $-\nu$ at the
corresponding time $t_k$ must get rescaled to tend to zero; while
if the maximum value of $(-\tilde{\nu})$ at the time $\tilde{t}_k$
does go to infinity, the maximum value of $\tilde{R}$ at
$\tilde{t}_k$ will go to infinity even faster from the pinching
estimate, and when we normalize to keep the normalized scalar
curvature $R$ bounded at the time $t_k$ so the normalized
$(-{\nu})$ at the time $t_k$ will go to zero. Thus in either case
the noncollapsing limit has nonnegative sectional curvature at the
initial time $t=0$ and then has nonnegative sectional curvature
for all times $t\ge 0$.

For Subcase (2.3), normalizing the flow only changes quantities in
a bounded way. As before we have the pinching estimate
$$R\ge (-\nu)[\log(-\nu)+\log(1+t)-C]$$
for the normalized Ricci flow, where $C$ is a positive constant
depending only on the constants $C_1$, $C_2$ in the assumption of
Subcase (2.3). If $$(-\nu)\le \frac{A}{1+t}$$ for any fixed
positive constant $A$, then $(-\nu)\rightarrow 0$ as $t\rightarrow
+\infty$ and we can take a noncollapsing limit which has
nonnegative sectional curvature. On the other hand if we can pick
a sequence of times $t_k\rightarrow\infty$ and points $x_k$ where
$(-\nu)(x_k,t_k)=\max\limits_{x\in M}(-\nu)(x,t_k)$ satisfies
$$
(-\nu)(x_k,t_k)(1+t_k)\rightarrow +\infty,\quad \text{as } \;
k\rightarrow +\infty,
$$
then from the pinching estimate, we have
$$
\frac{R(x_k,t_k)}{(-\nu)(x_k,t_k)}\rightarrow +\infty,\quad
\text{as}\quad k\rightarrow +\infty.
$$
But $R(x_k,t_k)$ are uniformly bounded since normalizing the flow
only changes quantities in bounded way. This shows $\sup
(-\nu)(\cdot,t_k)\rightarrow 0$ as $k\rightarrow+\infty$. Thus we
can take a noncollapsing limit along $t_k$ which has nonnegative
sectional curvature. Hence we have completed the proof of Case
(2).

We now come to the most interesting Case (3) where $R_{min}$
increases monotonically to a limit strictly less than zero. By
scaling we can assume $R_{min}(t)\rightarrow -6$ as $t\rightarrow
+\infty$.

\begin{lemma} [{Hamilton \cite{Ha99}}]
In Case $(3)$ where $R_{\min}\rightarrow-6$ as $t\rightarrow
+\infty$, all noncollapsing limit are hyperbolic with constant
sectional curvature $-1.$
\end{lemma}

\begin{pf}
By (5.3.4) and the fact $R_{\min}(t)\le-6$, we have
$$
\frac{d}{dt}R_{\min}(t)\ge 4(r(t)-R_{\min}(t))
$$
and
$$
\int_0^\infty(r(t)-R_{\min}(t))dt<+\infty.
$$
Since $r(t)-R_{\min}(t)\ge 0$ and $R_{min}(t)\rightarrow-6$ as
$t\rightarrow+\infty$, it follows that the function $r(t)$ has the
limit
$$
r=-6,
$$
for any convergent subsequence. And since
$$
\int_M(R-R_{\min}(t))d\mu=(r(t)-R_{\min}(t))\cdot V,
$$
it then follows that
$$
R\equiv-6\quad \text{for the limit}.
$$
The limit still has the following evolution equation for the
limiting scalar curvature
$$
\frac{\partial}{\partial t}R=\Delta R+2| {\mathop {\rm
Ric}\limits^\circ}|^2+\frac{2}{3}R(R-r).
$$
Since $R\equiv r\equiv -6$ in space and time for the limit, it
follows directly that $|{\mathop {\rm Ric}\limits^ \circ}|\equiv
0$ for the limit. Thus the limit metric has $\lambda=\mu=\nu=-2$,
so it has constant sectional curvature $-1$ as desired.
\end{pf}

If in the discussion above there exists a compact noncollapsing
limit, then we know that the underlying manifold $M$ is compact
and we fall into the conclusion of Theorem 5.3.4(i) for the
constant negative sectional curvature limit. Thus it remains to
show when every noncollapsing limit is a complete noncompact
hyperbolic manifold with finite volume, we have conclusion (ii) in
Theorem 5.3.4.

Now we first want to find a finite collection of persistent
complete noncompact hyperbolic manifolds as stated in Theorem
5.3.4 (ii).

\bigskip
{\bf Part II: Persistence of Hyperbolic Pieces}

\smallskip
We begin with the definition of the topology of $C^\infty$
convergence on compact sets for maps $F:M\rightarrow N$ of one
Riemannian manifold to another. For any compact set
$K\subset\subset M$ and any two maps $F,G:M\rightarrow N$, we
define
$$
d_K(F,G)=\sup_{x\in K}d(F(x),G(x))
$$
where $d(y,z)$ is the geodesic distance from $y$ to $z$ on $N$.
This gives the $C^0_{loc}$ topology for maps between $M$ and $N$.
To define $C^k_{loc}$ topology for any positive integer $k\ge1$,
we consider the \mbox{\boldmath{$k$}}-{\bf jet
space}\index{$k$-jet!space} $J^kM$ of a manifold $M$ which is the
collection of all
$$
(x,J^1,J^2,\ldots,J^k)
$$
where $x$ is a point on $M$ and $J^i$ is a tangent vector for
$1\le i\le k$ defined by the $i^{th}$ covariant derivative
$J^i=\nabla^i_\frac{\partial}{\partial t}\gamma(0)$ for a path
$\gamma$ passing through the point $x$ with $\gamma(0)=x$. A
smooth map $F:M\rightarrow N$ induces a map $$J^kF:\quad
J^kM\rightarrow J^kN$$ defined by
$$
J^kF(x,J^1,\ldots,J^k)=(F(x),\nabla_\frac{\partial}{\partial t}
(F(\gamma))(0),\ldots, \nabla^k_\frac{\partial}{\partial
t}(F(\gamma))(0))
$$
where $\gamma$ is a path passing through the point $x$ with
$J^i=\nabla^i_\frac{\partial}{\partial t}\gamma(0),1\le i\le k.$

Define the {\bf ${\mathbf k}$-jet
distance}\index{$k$-jet!distance} between $F$ and $G$ on a compact
set $K\subset\subset M$ by
$$
d_{C^k(K)}(F,G)=d_{BJ^kK}(J^kF,J^kG)
$$
where $BJ^kK$ consists of all $k$-jets $(x,J^1,\ldots,J^k)$ with
$x\in K$ and
$$
|J^1|^2+|J^2|^2+\cdots+|J^k|^2\le 1.
$$
Then the convergence in the metric $d_{C^k(K)}$ for all positive
integers $k$  and all compact sets $K$ defines the topology of
$C^\infty$ convergence on compact sets for the space of maps.

We will need the following \textbf{Mostow type rigidity}
\index{Mostow type rigidity} result (cf. Corollary 8.3 of
\cite{Ha99}).

\begin{lemma} For any complete noncompact
hyperbolic three-manifold $\mathcal{H}$ with finite volume with
metric $h$, we can find a compact set $\mathcal{K}$ of
$\mathcal{H}$ such that for every integer $k$ and every
$\varepsilon>0$, there exist an integer $q$ and a $\delta>0$ with
the following property: if $F$ is a diffeomorphism of
$\mathcal{K}$ into another complete noncompact hyperbolic
three-manifold $\tilde{\mathcal{H}}$ with no fewer cusps $($than
$\mathcal{H}),$ finite volume with metric $\tilde{h}$ such that
$$
\|F^*\tilde{h}-h\|_{C^q(\mathcal{K})}<\delta
$$
then there exists an isometry $I$ of $\mathcal{H}$ to
$\tilde{\mathcal{H}}$ such that
$$
d_{C^k(\mathcal{K})}(F,I)<\varepsilon.
$$
\end{lemma}

\begin{pf} This version of Mostow rigidity is given by Hamilton
(cf. section 8 of \cite{Ha99}). The following argument is in part
based on the Editors' notes in p.323-324 of \cite{CCCY}.

First we claim that $\mathcal{H}$ is isometric to
$\mathcal{\tilde{H}}$ for an appropriate choice of compact set
$\mathcal{K}$, positive integer $q$ and positive number $\delta$.
Let $l: \mathcal{H}\rightarrow \mathbb{R}$ be a function defined
at each point by the length of the shortest non-contractible loop
starting and ending at this point. Denote the Margulis constant by
$\mu$. Then by Margulis lemma (see for example \cite{Gro79} or
\cite{KM}), for any $0<\varepsilon_0<\frac{1}{2}\mu$, the set
$l^{-1}([0,\varepsilon_0])\subset \mathcal{H}$ consists of
finitely many components and each of these components is isometric
to a cusp or to a tube. Topologically, a tube is just a solid
torus. Let $\varepsilon_0$ be even smaller than one half of the
minimum of the lengths of the all closed geodesics on the tubes.
Then $l^{-1}([0,\varepsilon_0])$ consists of finite number of
cusps. Set $\mathcal{K}_0 = l^{-1}([\varepsilon_0,\infty))$. The
boundary of $\mathcal{K}_0$ consists of flat tori with constant
mean curvatures. Note that each embedded torus in a complete
hyperbolic three-manifold with finite volume either bounds a solid
torus or is isotopic to a standard torus in a cusp. The
diffeomorphism $F$ implies the boundary $F(\partial\mathcal{K}_0)$
are embedded tori. If one of components bounds a solid torus, then
as $\delta$ sufficiently small and $q$ sufficiently large,
$\mathcal{\tilde{H}}$ would have fewer cusps than $\mathcal{H}$,
which contradicts with our assumption. Consequently,
$\mathcal{\tilde{H}}$ is diffeomorphic to
$F({\overset{o}{\mathcal{K}_0}})$. Here
${\overset{o}{\mathcal{K}_0}}$ is the interior of the set
$\mathcal{K}_0$.  Since $\mathcal{H}$ is diffeomorphic to
${\overset{o}{\mathcal{K}_0}}$, $\mathcal{H}$ is diffeomorphic to
$\mathcal{\tilde{H}}$.  Hence by Mostow's rigidity theorem (see
\cite{Mos} and \cite{Pr}), $\mathcal{H}$ is isometric to
$\mathcal{\tilde{H}}$.

So we can assume $\mathcal{\tilde{H}}=\mathcal{H}$.  For
$\mathcal{K} = \mathcal{K}_0$, we argue by contradiction. Suppose
there is some $k>0$ and $\varepsilon>0$ so that there exist
sequences of integers $q_j\rightarrow\infty$, $\delta_j\rightarrow
0^+$ and diffeomorphisms $F_j$ mapping $\mathcal{K}$ into
$\mathcal{H}$ with
$$
  \|F_j^{*}h-h\|_{C^{q_j}(\mathcal{K})}<\delta_j
$$
 and
$$
 d_{C^{k}(\mathcal{K})}(F_j,I)\geq \varepsilon
$$
for all isometries $I$ of $\mathcal{H}$ to itself. We can extract
a subsequence of $F_j$ convergent to a map $F_{\infty}$ with
$F_{\infty}^{*}h=h$ on $\mathcal{K}$.

We need to check that $F_{\infty}$ is still a diffeomorphism on
$\mathcal{K}$. Since $F_{\infty}$ is a local diffeomorphism and is
the limit of diffeomorphisms, we can find an inverse of
$F_{\infty}$ on $F_{\infty}(\overset{o}{\mathcal{K}})$. So
$F_{\infty}$ is a diffeomorphism on $\overset{o}{\mathcal{K}}$. We
claim the image of the boundary can not touch the image of the
interior. Indeed, if $F_{\infty}(x_{1})=F_{\infty}(x_{2})$ with
$x_{1}\in \partial \mathcal{K} $ and $x_{2}\in
\overset{o}{\mathcal{K}} $, then we can find $x_{3}\in
\overset{o}{\mathcal{K}}$ near $x_{1}$ and $x_{4}\in
\overset{o}{\mathcal{K}}$ near $x_{2}$ with
$F_{\infty}(x_{3})=F_{\infty}(x_{4})$, since $F_{\infty}$ is a
local diffeomorphism. This contradicts with the fact that
$F_{\infty}$ is a diffeomorphism on $\overset{o}{\mathcal{K}}$.
This proves our claim. Hence, the only possible overlap is at the
boundary. But the image $F_{\infty}(\partial \mathcal{K})$ is
strictly concave, this prevents the boundary from touching itself.
We conclude that the mapping $F_{\infty}$ is a diffeomorphism on
$\mathcal{K}$, hence an isometry.

To extend $F_{\infty}$ to a global isometry, we argue as follows.
For each truncated cusp end of $\mathcal{K}$,  the area of
constant mean curvature flat torus is strictly decreasing. Since
$F_{\infty}$ takes each such torus to another of the same area, we
see that $F_{\infty}$ takes the foliation of an end by constant
mean curvature flat tori to another such foliation. So
$F_{\infty}$ takes cusps to cusps and preserves their foliations.
Note that the isometric type of a cusp is just the isometric type
of the torus, more precisely, let
$(N,dr^{2}+e^{-2r}g_{\mathcal{V}})$ be a cusp (where
$g_{\mathcal{V}}$ is the flat metric on the torus $\mathcal{V}$),
$0<a<b$ are two constants, any isometry of $N\cap l^{-1}[a,b]$ to
itself is just an isometry of $\mathcal{V}$. Hence the isometry
$F_{\infty}$ can be extended to the whole cusps. This gives a
global isometry $I$ contradicting our assumption when $j$ large
enough.

The proof of the Lemma 5.3.7 is completed.
\end{pf}

In order to obtain the persistent hyperbolic pieces stated in
Theorem 5.3.4 (ii), we will need to use a special parametrization
given by harmonic maps.

\begin{lemma} [{Hamilton \cite{Ha99}}]
Let\; $(X,\,g)$\; be\; a\; compact\; Riemannian\; manifold\; with
strictly negative Ricci curvature and with strictly concave
boundary. Then there are positive integer $l_0$ and small number
$\varepsilon_0>0$ such that for each positive integer $l\ge l_0$ and
positive number $\varepsilon\le\varepsilon_0$ we can find positive
integer $q$ and positive number $\delta>0$ such that for every
metric $\tilde{g}$ on $X$ with $||\tilde{g}-g||_{C^q(X)}\le \delta$
we can find a unique diffeomorphism $F$ of $X$ to itself so that
\begin{itemize}
\item[(a)] $F: (X,g)\rightarrow (X,\tilde{g})$ is harmonic,
\item[(b)] $F$ takes the boundary $\partial X$ to itself and
satisfies the {\bf free boundary condition}\index{free boundary
condition} that the normal derivative $\nabla_NF$ of $F$ at the
boundary is normal to the boundary, \item[(c)]
$d_{C^l(X)}(F,Id)<\varepsilon$, where $Id$ is the identity map.
\end{itemize}
\end{lemma}

\begin{pf}
The following argument is adapted from Hamilton's paper
\cite{Ha99} and the Editors' note on p.325 of \cite{CCCY}. Let
$\Phi(X,\partial X)$ be the space of maps of $X$ to itself which
take $\partial X$ to itself. Then $\Phi(X,\partial X)$ is a Banach
manifold and the tangent space to $\Phi(X,\partial X)$ at the
identity is the space of vector fields
$V=V^i\frac{\partial}{\partial x^i}$ tangent to the boundary.
Consider the map sending $F\in \Phi(X,\partial X)$ to the pair
$\{\Delta F,(\nabla_NF)_{//}\}$ consisting of the harmonic map
Laplacian and the tangential component (in the target) of the
normal derivative of $F$ at the boundary. By using the inverse
function theorem, we only need to check that the derivative of
this map is an isomorphism at the identity with $\tilde{g}=g$.

Let $\{x^i\}_{i=1,\ldots,n}$ be a local coordinates of $(X,g)$ and
$\{y^\alpha\}_{\alpha=1,\ldots,n}$ be a local coordinates of
$(X,\tilde{g})$. The harmonic map Laplacian of $F:
(X,g)\rightarrow (X,\tilde{g})$ is given in local coordinates by
$$
(\Delta F)^\alpha
=\Delta(F^\alpha)+g^{ij}(\tilde{\Gamma}^\alpha_{\beta\gamma}\circ
F) \frac{\partial F^\beta}{\partial x^i}\frac{\partial
F^\gamma}{\partial x^j}
$$
where $\Delta(F^\alpha)$ is the Laplacian of the function
$F^\alpha$ on $X$ and $\tilde{\Gamma}^\alpha_{\beta\gamma}$ is the
connection of $\tilde{g}$. Let $F$ be a one-parameter family with
$F|_{s=0}=Id$ and $\frac{dF}{ds}|_{s=0}=V$, a smooth vector field
on $X$ tangent to the boundary (with respect to $g$). At an
arbitrary given point $x\in X$, we choose the coordinates
$\{x^i\}_{i=1,\ldots,n}$ so that $\Gamma^i_{jk}(x)=0.$ We compute
at the point $x$ with $\tilde{g}=g$,
$$
\frac{d}{ds}{\Big|_{s=0}}(\Delta F)^\alpha
=\Delta(V^\alpha)+g^{ij}\(\frac{\partial}{\partial
x^k}\Gamma^\alpha_{ij}\)V^k.
$$
Since
$$(\nabla_iV)^{\alpha} =
\nabla_iV^{\alpha} + (\Gamma^{\alpha}_{i\beta}\circ F)V^{\beta},
$$
we have, at $s=0$ and the point $x$,
$$
(\Delta V)^\alpha=\Delta(V^\alpha)+g^{ij}\frac{\partial}{\partial
x^i}\Gamma^\alpha_{jk}V^k.
$$
Thus we obtain
\begin{align}
\frac{d}{ds}{|_{s=0}}(\Delta F)^\alpha&  = (\Delta
V)^\alpha+g^{ij}\(\frac{\partial}{\partial
x^k}\Gamma^\alpha_{ij}-\frac{\partial}{\partial
x^i}\Gamma^\alpha_{jk}\)V^k \\
&  = (\Delta V)^\alpha+g^{\alpha i}R_{ik}V^k.\nn
\end{align} 
Since
$$
(\nabla_NF)(F^{-1}(x)) =N^i(F^{-1}(x))\frac{\partial F^j}{\partial
x^i}(F^{-1}(x)) \frac{\partial}{\partial x^j}(x)\quad\text{on
}\;\partial X,
$$
we have
\begin{align}
\frac{d}{ds}|_{s=0}\{(\nabla_NF)_{//}\}
&  = \frac{d}{ds}|_{s=0}(\nabla_NF -\langle \nabla_NF,N \rangle N)\\
&  = \frac{d}{ds}|_{s=0}(\nabla_NF)-\left\langle
\frac{d}{ds}|_{s=0}\nabla_NF,N \right\rangle N  \nn\\
&\quad - \langle \nabla_NF,\nabla_VN \rangle N |_{s=0} - \langle
\nabla_NF,N
\rangle \nabla_VN |_{s=0} \nn\\
&  = \(\frac{d}{ds}|_{s=0}\nabla_NF\)_{//} -\nabla_VN \nn\\
&  = \(-V(N^i)\frac{\partial}{\partial
x^i}+N(V^j)\frac{\partial}{\partial x^j}\)_{//} - II(V) \nn\\
&  = [N,V]_{//} - II(V) \nn\\
&  = (\nabla_NV)_{//} - 2II(V) \nn
\end{align} 
where $II$ is the second fundamental form of the boundary (as an
automorphism of $T(\partial X)$). Thus by (5.3.5) and (5.3.6), the
kernel of the map sending $F\in \Phi(X,\partial X)$ to the pair
$\{\Delta F,(\nabla_NF)_{//}\}$ is the space of solutions of
elliptic boundary value problem \be
\begin{cases}
\Delta V+\Ric(V)=0 &\quad \text{on }\;X\\
V_{\bot}=0,&\quad \text{at }\;\partial X,\\
(\nabla_NV)_{//}-2II(V)=0,&\quad \text{at } \;\partial X,
\end{cases} 
\ee where $V_\bot$ is the normal component of $V$.

Now using these equations and integrating by parts gives
$$
\int\int\limits_X|\nabla V|^2
=\int\int\limits_X\Ric(V,V)+2\int\limits_{\partial X}II(V,V).
$$
Since $Rc<0$ and $II<0$ we conclude that the kernel is trivial.
Clearly this elliptic boundary value is self-adjoint because of
the free boundary condition. Thus the cokernel is trivial also.
This proves the lemma.
\end{pf}

Now we can prove the persistence of hyperbolic pieces. Let
$g_{ij}(t),\;0\le t< +\infty$, be a noncollapsing nonsingular
solution of the normalized Ricci flow on a compact three-manifold
$M$. Assume that any noncollapsing limit of the nonsingular
solution is a complete noncompact hyperbolic three-manifold with
finite volume. Consider all the possible hyperbolic limits of the
given nonsingular solution, and among them choose one such
complete noncompact hyperbolic three-manifold $\mathcal{H}$ with
the least possible number of cusps. In particular, we can find a
sequence of times $t_k\rightarrow +\infty$ and a sequence of
points $P_k$ on $M$ such that the marked three-manifolds
$(M,g_{ij}(t_k),P_k)$ converge in the $C^\infty_{loc}$ topology to
$\mathcal{H}$ with hyperbolic metric $h_{ij}$ and marked point
$P\in \mathcal{H}$.

For any small enough $a>0$ we can truncate each cusp of
$\mathcal{H}$ along a constant mean curvature torus of area $a$
which is uniquely determined; the remainder we denote by
$\mathcal{H}_a$. Clearly as $a\rightarrow 0$ the $\mathcal{H}_a$
exhaust $\mathcal{H}$. Pick a sufficiently small number $a>0$ to
truncate cusps so that Lemma 5.3.7 is applicable for the compact
set $\mathcal{K} = \mathcal{H}_a$. Choose an integer $l_0$ large
enough and an $\varepsilon_0$ sufficiently small to guarantee from
Lemma 5.3.8 the uniqueness of the identity map $Id$ among maps
close to $Id$ as a harmonic map $F$ from $\mathcal{H}_a$ to itself
with taking $\partial\mathcal{H}_a$ to itself, with the normal
derivative of $F$ at the boundary of the domain normal to the
boundary of the target, and with $d_{C^{l_0}(\mathcal{H}_a)}(F,Id)
< \varepsilon_0$. Then choose positive integer $q_0$ and small
number $\delta_0>0$ from Lemma 5.3.7 such that if $\tilde{F}$ is a
diffeomorphism of $\mathcal{H}_a$ into another complete noncompact
hyperbolic three-manifold $\tilde{\mathcal{H}}$ with no fewer
cusps (than $\mathcal{H}$), finite volume with metric
$\tilde{h}_{ij}$ satisfying
$$
||\tilde{F}^*\tilde{h}_{ij}-h_{ij}||_{C^{q_0}(\mathcal{H}_a)} \le
\delta_0,
$$
then there exists an isometry $I$ of $\mathcal{H}$ to
$\tilde{\mathcal{H}}$ such that \be
d_{C^{l_0}(\mathcal{H}_a)}(\tilde{F},I)<\varepsilon_0. 
\ee And we further require $q_0$ and $\delta_0$ from Lemma 5.3.8
to guarantee the existence of harmonic diffeomorphism from
$(\mathcal{H}_a,\tilde{g}_{ij})$ to $(\mathcal{H}_a,h_{ij})$ for
any metric $\tilde{g}_{ij}$ on $\mathcal{H}_a$ with
$||\tilde{g}_{ij}-h_{ij}||_{C^{q_0}(\mathcal{H}_a)}\le\delta_0.$

By definition, there exist a sequence of exhausting compact sets
$U_k$ of $\mathcal{H}$ (each $U_k \supset \mathcal{H}_a$) and a
sequence of diffeomorphisms $F_k$ from $U_k$ into $M$ such that
$F_k(P)=P_k$ and
$||F^*_kg_{ij}(t_k)-h_{ij}||_{C^m(U_k)}\rightarrow 0$ as
$k\rightarrow+\infty$ for all positive integers $m$. Note that
$\partial \mathcal{H}_a$ is strictly concave and we can foliate a
neighborhood of $\partial \mathcal{H}_a$ with constant mean
curvature hypersurfaces where the area $a$ has a nonzero gradient.
As the approximating maps
$F_k:\;(U_k,h_{ij})\rightarrow(M,g_{ij}(t_k))$ are close enough to
isometries on this collar of $\partial \mathcal{H}_a$, the metrics
$g_{ij}(t_k)$ on $M$ will also admit a unique constant mean
curvature hypersurface with the same area $a$ near $F_k(\partial
\mathcal{H}_a)(\subset M)$ by the inverse function theorem. Thus
we can change the map $F_k$ by an amount which goes to zero as
$k\rightarrow\infty$ so that now $F_k(\partial \mathcal{H}_a)$ has
constant mean curvature with the area $a$. Furthermore, by
applying Lemma 5.3.8 we can again change $F_k$ by an amount which
goes to zero as $k\rightarrow\infty$ so as to make $F_k$ a
harmonic diffeomorphism and take $\partial \mathcal{H}_a$ to the
constant mean curvature hypersurface $F_k(\partial \mathcal{H}_a)$
and also satisfy the free boundary condition that the normal
derivative of $F_k$ at the boundary of the domain is normal to the
boundary of the target. Hence for arbitrarily given positive
integer $q\ge q_0$ and positive number $\delta<\delta_0$, there
exists a positive integer $k_0$ such that for the modified
harmonic diffeomorphism $F_k$, when $k\ge k_0$,
$$
||F^*_kg_{ij}(t_k)-h_{ij}||_{C^q(\mathcal{H}_a)}<\delta.
$$

For each fixed $k\ge k_0$, by the implicit function theorem we can
first find a constant mean curvature hypersurface near
$F_k(\partial \mathcal{H}_a)$ in $M$ with the metric $g_{ij}(t)$
for $t$ close to $t_k$ and with the same area for each component
since $\partial \mathcal{H}_a$ is strictly concave and a
neighborhood of $\partial\mathcal{H}_a$ is foliated by constant
mean curvature hypersurfaces where the area $a$ has a nonzero
gradient and
$F_k:(\mathcal{H}_a,h_{ij})\rightarrow(M,g_{ij}(t_k))$ is close
enough to an isometry and $g_{ij}(t)$ varies smoothly. Then by
applying Lemma 5.3.8 we can smoothly continue the harmonic
diffeomorphism $F_k$ forward in time a little to a family of
harmonic diffeomorphisms $F_k(t)$ from $\mathcal{H}_a$ into $M$
with the metric $g_{ij}(t)$, with $F_k(t_k)=F_k$, where each
$F_k(t)$ takes $\partial \mathcal{H}_a$ into the constant mean
curvature hypersurface we just found in $(M,g_{ij}(t))$ and
satisfies the free boundary condition, and also satisfies
$$
||F^*_k(t)g_{ij}(t)-h_{ij}||_{C^q(\mathcal{H}_a)}<\delta.
$$
We claim that for all sufficiently large $k$, we can smoothly
extend the harmonic diffeomorphism $F_k$ to the family harmonic
diffeomorphisms $F_k(t)$ with
$||F^*_k(t)g_{ij}(t)-h_{ij}||_{C^q(\mathcal{H}_a)}\le\delta$ on a
maximal time interval $t_k\le t\leq\omega_k$ (or $t_k\leq
t<\omega_k$ when $\omega_k=+\infty$); and if $\omega_k<+\infty$,
then \be
||F^*_k(\omega_k)g_{ij}(\omega_k)-h_{ij}||_{C^q(\mathcal{H}_a)}
=\delta.  
\ee

Clearly the above argument shows that the set of $t$ where we can
extend the harmonic diffeomorphisms as desired is open. To verify
claim (5.3.9), we thus only need to show that if we have a family
of harmonic diffeomorphisms $F_k(t)$ such as we desire for $t_k\le
t<\omega (<+\infty)$, we can take the limit of $F_k(t)$ as
$t\rightarrow \omega$ to get a harmonic diffeomorphism
$F_k(\omega)$ satisfying
$$
||F^*_k(\omega)g_{ij}(\omega)-h_{ij}||_{C^q(\mathcal{H}_a)}\le\delta,
$$
and if
$$
||F^*_k(\omega)g_{ij}(\omega)-h_{ij}||_{C^q(\mathcal{H}_a)} <
\delta,
$$
then we can extend $F_k(\omega)$ forward in time a little (i.e.,
we can find a constant mean curvature hypersurface near
$F_k(\omega)(\partial \mathcal{H}_a)$ in $M$ with the metric
$g_{ij}(t)$ for each $t$ close to $\omega$ and with the same area
$a$ for each component). Note that \be
||F^*_k(t)g_{ij}(t)-h_{ij}||_{C^q(\mathcal{H}_a)}<\delta.
\ee for $t_k\le t<\omega$ and the metrics $g_{ij}(t)$ for $t_k\le
t\le\omega$ are uniformly equivalent. We can find a subsequence
$t_n\rightarrow \omega$ for which $F_k(t_n)$ converge to
$F_k(\omega)$ in $C^{q-1}(\mathcal{H}_a)$ and the limit map has
$$
||F^*_k(\omega)g_{ij}(\omega)-h_{ij}||_{C^{q-1}(\mathcal{H}_a)}
\le\delta.
$$
We need to check that $F_k(\omega)$ is still a diffeomorphism. We
at least know $F_k(\omega)$ is a local diffeomorphism, and
$F_k(\omega)$ is the limit of diffeomorphisms, so the only
possibility of overlap is at the boundary. Hence we use the fact
that $F_k(\omega)(\partial\mathcal{H}_a)$ is still strictly
concave since $q$ is large and $\delta$ is small to prevent the
boundary from touching itself. Thus $F_k(\omega)$ is a
diffeomorphism. A limit of harmonic maps is harmonic, so
$F_k(\omega)$ is a harmonic diffeomorphism from $\mathcal{H}_a$
into $M$ with the metric $g_{ij}(\omega)$. Moreover $F_k(\omega)$
takes $\partial\mathcal{H}_a$ to the constant mean curvature
hypersurface $\partial (F_k(\omega)(\mathcal{H}_a))$ of the area
$a$ in $(M,g_{ij}(\omega))$ and continue to satisfy the free
boundary condition. As a consequence of the standard regularity
result of elliptic partial differential equations (see for example
\cite{GT}), the map $F_k(\omega)\in C^\infty(\mathcal{H}_a)$ and
then from (5.3.10) we have
$$
||F^*_k(\omega)g_{ij}(\omega)-h_{ij}||_{C^q(\mathcal{H}_a)}\le\delta.
$$
If
$||F^*_k(\omega)g_{ij}(\omega)-h_{ij}||_{C^q(\mathcal{H}_a)}=\delta$,
we then finish the proof of the claim. So we may assume that
$||F^*_k(\omega)g_{ij}(\omega)-h_{ij}||_{C^q(\mathcal{H}_a)}<\delta$.
We want to show that $F_k(\omega)$ can be extended forward in time
a little.

We argue by contradiction. Suppose not, then we consider the new
sequence of the manifolds $M$ with metric $g_{ij}(\omega)$ and the
origins $F_k(\omega)(P)$. Since $F_k(\omega)$ are close to
isometries, the injectivity radii of the metrics $g_{ij}(\omega)$
at $F_k(\omega)(P)$ do not go to zero, and we can extract a
subsequence which converges to a hyperbolic limit
$\widetilde{\mathcal{H}}$ with the metric $\widetilde{h}_{ij}$ and
the origin $\widetilde{P}$ and with finite volume. The new limit
$\widetilde{\mathcal{H}}$ has at least as many cusps as the old
limit $\mathcal{H}$, since we choose $\mathcal{H}$ with cusps as
few as possible. By the definition of convergence, we can find a
sequence of compact sets $\widetilde{B}_k$ exhausting
$\widetilde{\mathcal{H}}$ and containing $\widetilde{P}$, and a
sequence of diffeomorphisms $\widetilde{F}_k$ of neighborhoods of
$\widetilde{B}_k$ into $M$ with $\widetilde{F}_k(\widetilde{P}) =
F_k(\omega)(P)$ such that for each compact set $\widetilde{B}$ in
$\widetilde{\mathcal{H}}$ and each integer $m$
$$
||\widetilde{F}^*_k(g_{ij}(\omega)) -
\widetilde{h}_{ij}||_{C^m(\widetilde{B})} \rightarrow 0
$$
as $k \rightarrow +\infty$. For large enough $k$ the set
$\widetilde{F}_k(\widetilde{B}_k)$ will contain all points out to
any fixed distance we need from the point $F_k(\omega)(P)$, and
then
$$
\widetilde{F}_k(\widetilde{B}_k)\supset F_k(\omega)(\mathcal{H}_a)
$$
since the points of $\mathcal{H}_a$ have a bounded distance from
$P$ and $F_k(\omega)$ are reasonably close to preserving the
metrics. Hence we can form the composition
$$
G_k=\widetilde{F}^{-1}_k\circ F_k(\omega): \mathcal{H}_a
\rightarrow \widetilde{\mathcal{H}}.
$$
Arbitrarily fix $\delta' \in (\delta, \delta_0)$. Since the
$\widetilde{F}_k$ are as close to preserving the metric as we
like, we have
$$
||G^*_k\widetilde{h}_{ij} - h_{ij}||_{C^q(\mathcal{H}_a)} <
\delta'
$$
for all sufficiently large $k$. By Lemma 5.3.7, we deduce that
there exists an isometry $I$ of $\mathcal{H}$ to
$\widetilde{\mathcal{H}}$, and then
$(M,g_{ij}(\omega),F_k(\omega)(P))$ (on compact subsets) is very
close to $(\mathcal{H},h_{ij},P)$ as long as $\delta$ small enough
and $k$ large enough. Since $F_k(\omega)(\partial \mathcal{H}_a)$
is strictly concave and the foliation of a neighborhood of
$F_k(\omega)(\partial \mathcal{H}_a)$ by constant mean curvature
hypersurfaces has the area as a function with nonzero gradient, by
the implicit function theorem, there exists a unique constant mean
curvature hypersurface with the same area $a$ near
$F_k(\omega)(\partial \mathcal{H}_a)$ in $M$ with the metric
$g_{ij}(t)$ for $t$ close to $\omega$. Hence, when $k$
sufficiently large, $F_k(\omega)$ can be extended forward in time
a little. This is a contradiction and we have proved claim
(5.3.9).

We further claim that there must be some $k$ such that
$\omega_k=+\infty$ (i.e., we can smoothly continue the family of
harmonic diffeomorphisms $F_k(t)$ for all $t_k\le t<+\infty$, in
other words, there must be at least one hyperbolic piece
persisting). We argue by contradiction. Suppose for each $k$ large
enough, we can continue the family $F_k(t)$ for $t_k\le t\le
\omega_k<+\infty$ with
$$
||F^*_k(\omega_k)g_{ij}(\omega_k)-h_{ij}||_{C^q(\mathcal{H}_a)}
=\delta.
$$
Then as before, we consider the new sequence of the manifolds $M$
with metrics $g_{ij}(\omega_k)$ and origins $F_k(\omega_k)(P)$.
For sufficiently large $k$, we can obtain diffeomorphisms
$\widetilde{F}_k$ of neighborhoods of $\widetilde{B}_k$ into $M$
with $\widetilde{F}_k(\widetilde{P})=F_k(\omega_k)(P)$ which are
as close to preserving the metric as we like, where
$\widetilde{B}_k$ is a sequence of compact sets, exhausting some
hyperbolic three-manifold $\tilde{\mathcal{H}}$, of finite volume
and with no fewer cusps (than $\mathcal{H}$), and containing
$\tilde{P}$; moreover, the set $\widetilde{F}_k(\widetilde{B}_k)$
will contain all the points out to any fixed distance we need from
the point $F_k(\omega_k)(P);$ and hence
$$
\widetilde{F}_k(\widetilde{B}_k)\supseteq
F_k(\omega_k)(\mathcal{H}_a)
$$
since $\mathcal{H}_a$ is at bounded distance from $P$ and
$F_k(\omega_k)$ is reasonably close to preserving the metrics.
Then we can form the composition
$$
G_k=\widetilde{F}^{-1}_k\circ F_k(\omega_k):\quad
\mathcal{H}_a\rightarrow\tilde{\mathcal{H}}.
$$
Since the $\widetilde{F}_k$ are as close to preserving the metric
as we like, for any $\widetilde{\delta}>\delta$ we have
$$
||G^*_k\widetilde{h}_{ij}-h_{ij}||_{C^q(\mathcal{H}_a)}
<\widetilde{\delta}
$$
for large enough $k$. Then a subsequence of $G_k$ converges at
least in $C^{q-1}(\mathcal{H}_a)$ topology to a map $G_\infty$ of
$\mathcal{H}_a$ into $\widetilde{\mathcal{H}}$. By the same reason
as in the argument of previous two paragraphs, the limit map
$G_\infty$ is a smooth harmonic diffeomorphism from
$\mathcal{H}_a$ into $\widetilde{\mathcal{H}}$ with the metric
$\tilde{h}_{ij}$, and takes $\partial \mathcal{H}_a$ to a constant
mean curvature hypersurface $G_\infty(\partial \mathcal{H}_a)$ of
$(\tilde{\mathcal{H}},\widetilde{h}_{ij})$ with the area $a$, and
also satisfies the free boundary condition. Moreover we still have
\be ||G^*_\infty\widetilde{h}_{ij}-h_{ij}||_{C^q(\mathcal{H}_a)}
=\delta. 
\ee Now by Lemma 5.3.7 we deduce that there exists an isometry $I$
of $\mathcal{H}$ to $\widetilde{\mathcal{H}}$ with
$$
d_{C^{l_0}(\mathcal{H}_a)}(G_\infty,I)<\varepsilon_0.
$$
By using $I$ to identify $\widetilde{\mathcal{H}}_a$ and
$\mathcal{H}_a$, we see that the map $I^{-1}\circ G_\infty $ is a
harmonic diffeomorphism of $\mathcal{H}_a$ to itself which
satisfies the free boundary condition and
$$
d_{C^{l_0}(\mathcal{H}_a)}(I^{-1}\circ G_\infty,Id)<\varepsilon_0.
$$
{}From the uniqueness in Lemma 5.3.8 we conclude that $I^{-1}\circ
G_\infty =Id$ which contradicts with (5.3.11). This shows at least
one hyperbolic piece persists. Moreover the pull-back of the
solution metric $g_{ij}(t)$ by $F_k(t)$, for $t_k\leq t<+\infty$,
is as close to the hyperbolic metric $h_{ij}$ as we like.

We can continue to form other persistent hyperbolic pieces in the
same way as long as there are any points $P_k$ outside of the
chosen pieces where the injectivity radius at times
$t_k\rightarrow \infty$ are all at least some fixed positive
number $\rho >0$. The only modification in the proof is to take
the new limit $\mathcal{H}$ to have the least possible number of
cusps out of all remaining possible limits.

Note that the volume of the normalized Ricci flow is constant in
time. Therefore by combining with Margulis lemma (see for example
\cite{Gro79} \cite{KM}), we have proved that there exists a finite
collection of complete noncompact hyperbolic three-manifolds
${\mathcal{H}}_1,\ldots,{\mathcal{H}}_m $ with finite volume, a
small number $a>0$ and a time $T<+\infty$ such that for all $t$
beyond $T$ we can find diffeomorphisms $\varphi_l(t)$ of
$({\mathcal{H}}_{l})_a $ into $M$, $1\leq l\leq m$, so that the
pull-back of the solution metric $g_{ij}(t)$ by $\varphi_l(t)$ is
as close to the hyperbolic metrics as we like and the exceptional
part of $M$ where the points are not in the image of any
$\varphi_l$ has the injectivity radii everywhere as small as we
like.

\medskip
{\bf Part III: Incompressibility} \vskip 0.2cm

We remain to show that the boundary tori of any persistent
hyperbolic piece are incompressible, in the sense that the
fundamental group of the torus injects into that of the whole
manifold. The argument of this part is a parabolic version of
Schoen and Yau's minimal surface argument in \cite{ScY79m, ScY79,
ScY82}.

Let $B$ be a small positive number and assume the above positive
number $a$ is much smaller than $B$. Denote by $M_a$ a persistent
hyperbolic piece of the manifold $M$ truncated by boundary tori of
area $a$ with constant mean curvature and denote by $M_a^c=M
\setminus \stackrel{\circ}{M_a}$ the part of $M$ exterior to
$M_a$. Thus there is a persistent hyperbolic piece $M_B \subset
M_a$ of the manifold $M$ truncated by boundary tori of area $B$
with constant mean curvature. We also denote by $M_B^c=M \setminus
\stackrel{\circ}{M_B}$. By Van Kampen's Theorem, if
$\pi_1(\partial M_B)$ injects into $\pi_1(M_B^c)$ then it injects
into $\pi_1(M)$ also. Thus we only need to show $\pi_1(\partial
M_B)$ injects into $\pi_1(M_B^c)$.

We will argue by contradiction. Let $T$ be a torus in $\partial
M_B$. Suppose $\pi_1(T)$ does not inject into $\pi_1(M_B^c)$, then
by Dehn's Lemma the kernel is a cyclic subgroup of $\pi_1(T)$
generated by a primitive element. The work of Meeks-Yau \cite{MY}
or Meeks-Simon-Yau \cite{MSY} shows that among all disks in
$M_B^c$ whose boundary curve lies in $T$ and generates the kernel,
there is a smooth embedded disk normal to the boundary which has
the least possible area. Let $A=A(t)$ be the area of this disk.
This is defined for all $t$ sufficiently large. We will show that
$A(t)$ decreases at a rate bounded away from zero which will be a
contradiction.

Let us compute the rate at which $A(t)$ changes under the Ricci
flow. We need to show $A(t)$ decrease at least at a certain rate,
and since $A(t)$ is the minimum area to bound any disk in the
given homotopy class, it suffices to find some such disk whose
area decreases at least that fast. We choose this disk as follows.
Pick the minimal disk  at time $t_0$, and extend it smoothly a
little past the boundary torus since the minimal disk is normal to
the boundary. For times $t$ a little bigger than $t_0$, the
boundary torus may need to move a little to stay constant mean
curvature with area $B$ as the metrics change, but we leave the
surface alone and take the bounding disk to be the one cut off
from it by the new torus. The change of the area $\tilde{A}(t)$ of
such disk comes from the change in the metric and the change in
the boundary.

For the change in the metric, we choose an orthonormal frame
$X,Y,Z$ at a point $x$ in the disk so that $X$ and $Y$ are tangent
to the disk while $Z$ is normal and compute the rate of change of
the area element $d\sigma$ on the disk as
\begin{align*}
\frac{\partial}{\partial t}d\sigma&  =
\frac{1}{2}(g^{ij})^T\(\frac{2}{3}r(g_{ij}\)^T-2(R_{ij})^T)d\sigma\\[2mm]
&  = \left[\frac{2}{3}r-\Ric(X,X)-\Ric(Y,Y)\right]d\sigma,
\end{align*}
since the metric evolves by the normalized Ricci flow. Here
$(\cdot)^T$ denotes the tangential projections on the disk. Notice
the torus $T$ may move in time to preserve constant mean curvature
and constant area $B$. Suppose the boundary of the disk evolves
with a normal velocity $N$. The change of the area at boundary
along a piece of length $ds$ is given by $Nds$. Thus the total
change of the area $\tilde{A}(t)$ is given by
$$
\frac{d\tilde{A}}{dt}
=\int\int\(\frac{2}{3}r-\Ric(X,X)-\Ric(Y,Y)\)d\sigma
+\int_\partial Nds.
$$
Note that
\begin{align*}
\Ric(X,X)+\Ric(Y,Y)&  = R(X,Y,X,Y)+R(X,Z,X,Z)\\
&\quad +R(Y,X,Y,X)+R(Y,Z,Y,Z)\\
&  = \frac{1}{2}R+R(X,Y,X,Y).
\end{align*}
By the Gauss equation, the Gauss curvature $K$ of the disk is
given by
$$
K=R(X,Y,X,Y)+\det II
$$
where $II$ is the second fundamental form of the disk in $M_B^c$.
This gives at $t=t_0$,
$$
\frac{dA}{dt}\leq\int\int\(\frac{2}{3}r-\frac{1}{2}R\)d\sigma
-\int\int(K-\det II)d\sigma+\int_\partial N ds
$$
Since the bounding disk is a minimal surface, we have
$$
\det II\leq 0.
$$
The Gauss-Bonnet Theorem tells us that for a disk
$$
\int\int Kd\sigma+\int_\partial k ds=2\pi
$$
where $k$ is the geodesic curvature of the boundary. Thus we
obtain \be \frac{dA}{dt}\leq
\int\int\(\frac{2}{3}r-\frac{1}{2}R\)d\sigma+\int_\partial
kds+\int_\partial Nds -2\pi. 
\ee Recall that we are assuming $R_{min}(t)$ increases
monotonically to $-6$ as $t\rightarrow +\infty$. By the evolution
equation of the scalar curvature,
$$
\frac{d}{dt}R_{\min}(t)\geq 4(r(t)-R_{\min}(t))
$$
and then
$$
\int_0^\infty (r(t)-R_{\min}(t))dt<+\infty.
$$
This implies that $r(t)\rightarrow -6$ as $t\rightarrow +\infty$
by using the derivatives estimate for the curvatures. Thus for
every $\varepsilon >0$ we have
$$
\frac{2}{3}r-\frac{1}{2}R\leq -(1-\varepsilon)
$$
for $t$ sufficiently large. And then the first term on RHS of
(5.3.12) is bounded above by
$$
\int\int \(\frac{2}{3}r-\frac{1}{2}R\)d\sigma\leq
-(1-\varepsilon)A.
$$
The geodesic curvature $k$ of the boundary of the minimal disk is
the acceleration of a curve moving with unit speed along the
intersection of the disk with the torus; since the disk and torus
are normal, this is the same as the second fundamental form of the
torus in the direction of the curve of intersection. Now if the
metric were actually hyperbolic, the second fundamental form of
the torus would be exactly 1 in all directions. Note that the
persistent hyperbolic pieces are as close to the standard
hyperbolic as we like. This makes that the second term of RHS of
(5.3.12) is bounded above by
$$
\int_\partial kds\leq(1+\varepsilon_0)L
$$
for some sufficiently small positive number $\varepsilon_0>0$,
where $L$ is the length of the boundary curve. Also since the
metric on the persistent hyperbolic pieces are close to the
standard hyperbolic as we like, its change under the normalized
Ricci flow is as small as we like; So the motion of the constant
mean curvature torus of fixed area $B$ will have a normal velocity
$N$ as small as we like. This again makes the third term of RHS of
(5.3.12) bounded above by
$$
\int_\partial Nds\leq \varepsilon_0 L.
$$
Combining these estimates, we obtain \be
\frac{dA}{dt}\leq(1+2\varepsilon_0)L-(1-\varepsilon_0)A-2\pi
\ee on the persistent hyperbolic piece, where $\varepsilon_0$ is
some sufficiently small positive number.

We next need to bound the length $L$ in terms of the area $A$.
Since $a$ is much smaller than $B$, for large $t$ the metric is as
close as we like to the standard hyperbolic one; not just on the
persistent hyperbolic piece $M_B$ but as far beyond as we like.
Thus for a long distance into $M_B^c$ the metric will look nearly
like a standard hyperbolic cusplike collar.

Let us first recall a special coordinate system on the standard
hyperbolic cusp projecting beyond torus $T_1$ in $\partial
{\mathcal {H}}_1$ as follows. The universal cover of the flat
torus $T_1$ can be mapped conformally to the $x$-$y$ plane so that
the deck transformation of $T_1$ become translations in $x$ and
$y$, and so that the Euclidean area of the quotient is 1; then
these coordinates are unique up to a translation. The hyperbolic
cusp projecting beyond the torus $T_1$ in $\partial {\mathcal
{H}}_1$  can be parametrized by $\{(x,y,z)\in {\mathbb{R}}^3\ |\ z
> 0\}$ with the hyperbolic metric
\be
ds^2=\frac{dx^2+dy^2+dz^2}{z^2}. 
\ee Note that we can make the solution metric, in an arbitrarily
large neighborhood of the torus $T$ (of $\partial M_B$), as close
to hyperbolic as we wish (in the sense that there exists a
diffeomorphism from a large neighborhood of the torus $T_B$ (of
$\partial\mathcal{H}_B$) on the standard hyperbolic cusp to the
above neighborhood of the torus $T$ (of $\partial M_B$) such that
the pull-back of the solution metric by the diffeomorphism is as
close to the hyperbolic metric as we wish). Then by using this
diffeomorphism (up to a slight modification) we can parametrize
the cusplike tube of $M_B^c$ projecting beyond the torus $T$ in
$\partial M_B$ by $\{(x,y,z)\ |\  z\geq \zeta\}$ where the height
$\zeta$ is chosen so that the torus in the hyperbolic cusp at
height $\zeta$ has the area $B$.

Now consider our minimal disk, and let $L(z)$ be the length of the
curve of the intersection of the disk with the torus at height $z$
in the above coordinate system, and also let $A(z)$ be the area of
the part of the disk between $\zeta$ and $z$. We now want to
derive a monotonicity formula on the area $A(z)$ for the minimal
surface.

For almost every $z$ the intersection of the disk with the torus
at height $z$ is a smooth embedded curve or a finite union of them
by the standard transversality theorem. If there is more than one
curve, at least one of them is not homotopic to a point in $T$ and
represents the primitive generator in the kernel of $\pi_1(T)$
such that a part of the original disk beyond height $z$ continues
to a disk that bounds it. We extend this disk back to the initial
height $\zeta$ by dropping the curve straight down. Let
$\tilde{L}(z)$ be the length of the curve we picked at height $z$;
of course $\tilde{L}(z)\leq L(z)$ with equality if it is the only
piece. Let $\tilde{L}(w)$ denote the length of the same curve in
the $x$-$y$ plane dropped down to height $w$ for $\zeta\leq w \leq
z$. In the hyperbolic space we would have
$$
\tilde{L}(w)=\frac{z}{w}\tilde{L}(z)
$$
exactly. In our case there is a small error proportional to
$\tilde{L}(z)$ and we can also take it proportional to the
distance $z-w$ by which it drops since
$\tilde{L}(w)|_{w=z}=\tilde{L}(z)$ and the solution metric is
close to the hyperbolic in the $C_{loc}^{\infty}$ topology. Thus,
for arbitrarily given $\delta>0$ and $\zeta^{\ast}>\zeta$, as the
solution metric is sufficiently close to hyperbolic, we have
$$
|\tilde{L}(w)-\frac{z}{w}\tilde{L}(z)|\leq \delta(z-w)\tilde{L}(z)
$$
for all $z$ and $w$ in $\zeta\leq w\leq z\leq {\zeta}^{\ast}$. Now
given $\varepsilon$ and ${\zeta}^{\ast}$ pick
$\delta=2\varepsilon/\zeta^{\ast} $. Then \be \tilde{L}(w)\leq
\frac{z}{w}\tilde{L}(z)\left[1+\frac{2\varepsilon(z-w)}{w}\right].
\ee When we drop the curve vertically for the construction of the
new disk we get an area $\tilde{A}(z)$ between $\zeta$ and $z$
given by
$$
\tilde{A}(z)=(1+o(1))\int_{\zeta}^{z} \frac{\tilde{L}(w)}{w}dw.
$$
Here and in the following $o(1)$ denotes various small error
quantities as the solution metric close to hyperbolic. On the
other hand if we do not drop vertically we pick up even more area,
so the area $A(z)$ of the original disk between $\zeta$ and $z$
has \be A(z)\geq (1-o(1))\int_{\zeta}^{z} \frac{L(w)}{w} dw.
\ee Since the original disk minimized among all disks bounded a
curve in the primitive generator of the kernel of $\pi_1(T)$, and
the new disk beyond the height $z$ is part of the original disk,
we have
$$
A(z)\leq \tilde{A}(z)
$$
and then by combining with (5.3.15),
\begin{align*}
\int_{\zeta}^{z} \frac{L(w)}{w}dw &  \leq (1+o(1))
z\tilde{L}(z)\int_{\zeta}^{z}
\left[\frac{1-2\varepsilon}{w^2}+\frac{2\varepsilon z}{w^3}\right]dw \\
& \leq (1+o(1))L(z)\frac{(z-\zeta)}{\zeta}
\left[1+\varepsilon\(\frac{z-\zeta}{\zeta}\)\right].
\end{align*}
Here we used the fact that $\tilde{L}(z)\leq L(z)$. Since the
solution metric is sufficiently close to hyperbolic, we have
\begin{align*}
\frac{d}{dz}\(\int_{\zeta}^{z}\frac{L(w)}{w}dw\)
&  = \frac{L(z)}{z}\\
& \geq(1-o(1))\frac{\zeta}{z(z-\zeta)}
\left[1-\varepsilon\(\frac{z-\zeta}{\zeta}\)\right]
\int_{\zeta}^{z}\frac{L(w)}{w}dw \\
&  \geq\left[\frac{1}{z-\zeta}-\frac{1+2\varepsilon}{z}\right]
\int_{\zeta}^{z}\frac{L(w)}{w}dw,
\end{align*}
or equivalently \be \frac{d}{dz}\log\left\{
\frac{z^{1+2\varepsilon}}{(z-\zeta)}
\int_{\zeta}^{z}\frac{L(w)}{w}dw\right\}\geq 0.
\ee This is the desired monotonicity formula for the area $A(z)$.

It follows directly from (5.3.16) and (5.3.17) that
$$
\frac{z^{1+2\varepsilon}}{(z-\zeta)}A(z)\geq
(1-o(1))\zeta^{2\varepsilon}L(\zeta),
$$
or equivalently
$$
L(\zeta)\leq (1+o(1))
\(\frac{z}{\zeta}\)^{2\varepsilon}\frac{z}{z-\zeta} A(z)
$$
for all $z\in [\zeta,\zeta^{\ast}]$. Since the solution metric, in
an arbitrarily large neighborhood of the torus $T$ (of $\partial
M_B)$, as close to hyperbolic as we wish, we may assume that
${\zeta}^{\ast}$ is so large that $\sqrt{\zeta^{\ast}}>\zeta$ and
$\frac{\sqrt{\zeta^{\ast}}}{\sqrt{\zeta^{\ast}}-\zeta}$ is close
to 1, and also $\varepsilon >0$ is so small that
$(\frac{\sqrt{\zeta^{\ast}}}{\zeta})^{2\varepsilon}$ is close to
1. Thus for arbitrarily small $\delta_0 > 0$, we have \be
L(\zeta)\leq (1+\delta_0)A\(\sqrt{\zeta^{\ast}}\). 
\ee Now recall that (5.3.13) states
$$
\frac{dA}{dt}\leq (1+2\varepsilon_0)L-(1-\varepsilon_0)A-2\pi.
$$
We now claim that if
$$
(1+2\varepsilon_0)L-(1-\varepsilon_0)A\geq 0
$$
then $L=L(\zeta)$ is uniformly bounded from above.

Indeed by the assumption we have
$$
A(\zeta^{\ast})\leq
\frac{(1+2\varepsilon_0)}{(1-\varepsilon_0)}L(\zeta)
$$
since $A(\zeta^{\ast})\leq A$. By combining with (5.3.16) we have
some $z_0\in (\sqrt{\zeta^{\ast}},\zeta^{\ast})$ satisfying
\begin{align*}
\frac{L(z_0)}{z_0}\(\zeta^{\ast}-\sqrt{\zeta^{\ast}}\)&  \leq
\int_{\sqrt{\zeta^{\ast}}}^{\zeta^{\ast}} \frac{L(w)}{w}dw \\
&  \leq  (1+o(1))A(\zeta^{\ast}) \\
&  \leq
(1+o(1))\(\frac{1+2\varepsilon_0}{1-\varepsilon_0}\)L(\zeta).
\end{align*}
Thus for $\zeta^{\ast}$ suitably large, by noting that the
solution metric on a large neighborhood of $T$ (of $\partial
{M}_B)$ is sufficiently close to hyperbolic, we have \be
L(z_0)\leq (1+4\varepsilon_0)\frac{z_0}{\zeta^{\ast}}L(\zeta)
\ee for some $z_0\in (\sqrt{\zeta^{\ast}},\zeta^{\ast})$. It is
clear that we may assume the intersection curve between the
minimal disk with the torus at this height $z_0$ is smooth and
embedded. If the intersection curve at the height $z_0$ has more
than one piece, as before one of them will represent the primitive
generator in the kernel of $\pi_1(T)$, and we can ignore the
others. Let us move (the piece of) the intersection curve on the
torus at height $z_0$ through as small as possible area in the
same homotopy class of $\pi_1(T)$ to a curve which is a geodesic
circle in the flat torus coming from our special coordinates, and
then drop this geodesic circle vertically in the special
coordinates to obtain another new disk. We will compare the area
of this new disk with the original minimal disk as follows.

Denote by $G$ the length of the geodesic circle in the standard
hyperbolic cusp at height 1. Then the length of the geodesic
circle at height $z_0$ will be $G/z_0$. Observe that given an
embedded curve of length $l$ circling the cylinder $S^1\times
\mathbb{R}$ of circumference $w$ once, it is possible to deform
the curve through an area not bigger than $lw$ into a meridian
circle. Note that (the piece of) the intersection curve represents
the primitive generator in the kernel of $\pi_1(T)$. Note also
that the solution metric is sufficiently close to the hyperbolic
metric. Then the area of the deformation from (the piece of) the
intersection curve on the torus at height $z_0$ to the geodesic
circle at height $z_0$ is bounded by
$$
(1+o(1))\(\frac{G}{z_0}\)\cdot L(z_0).
$$
The area to drop the geodesic circle from height $z_0$ to height
$\zeta$ is bounded by
$$
(1+o(1))\int_{\zeta}^{z_0}\frac{G}{w^2}dw.
$$
Hence comparing the area of the original minimal disk to that of
this new disk gives
$$
A(z_0)\leq (1+o(1))G\left[\frac{L(z_0)}{z_0}
+\(\frac{1}{\zeta}-\frac{1}{z_0}\)\right].
$$
By (5.3.18), (5.3.19) and the fact that
$z_0\in(\sqrt{\zeta^{\ast}},\zeta^{\ast})$, this in turn gives
\begin{align*}
L(\zeta)&  \leq (1+\delta_0)A(z_0) \\
&  \leq(1+\delta_0)G
\left[(1+4\varepsilon_0)\frac{L(\zeta)}{\zeta^{\ast}}+\frac{1}{\zeta}\right].
\end{align*}
Since $\zeta^{\ast}$ is suitably large, we obtain
$$
L(\zeta)\leq 2G/\zeta
$$
This gives the desired assertion since $G$ is fixed from the
geometry of the limit hyperbolic manifold $\mathcal{H}$ and
$\zeta$ is very large as long as the area $B$ of $\partial M_B$
small enough.

Thus the combination of (5.3.13), (5.3.18) and the assertion
implies that either
$$
\frac{d}{dt}A\leq-2\pi,
$$
or
\begin{align*}
\frac{d}{dt}A&  \leq (1+2\varepsilon_0)L-(1-\varepsilon_0)A-2\pi\\
&  \leq (1+2\varepsilon_0)\frac{2G}{\zeta} - 2\pi\\
&  \leq  -\pi,
\end{align*}
since the solution metric on a very large neighborhood of the
torus $T$ (of $\partial {M}_B$) is sufficiently close to
hyperbolic and $\zeta$ is very large as the area $B$ of $\partial
M_B$ small enough. This is impossible because $A\geq 0$ and the
persistent hyperbolic pieces go on forever. The contradiction
shows that $\pi_1(T)$ in fact injects into $\pi_1(M_B^c)$. This
proves that $\pi_1(\partial{M}_B) $ injects into $\pi_1(M)$.

Therefore we have completed the proof of Theorem 5.3.4.
\endproof

\newpage
\part{{\Large Ancient $\kappa$-solutions}}

\bigskip
Let us consider a solution of the Ricci flow on a compact
manifold. If the solution blows up in finite time (i.e., the
maximal solution exists only on a finite time interval), then as
we saw in Chapter 4 a sequence of rescalings of the solution
around the singularities converge to a solution which exists at
least on the time interval $(-\infty,T)$ for some finite number
$T$. Furthermore, by Perelman's no local collapsing theorem I
(Theorem 3.3.2), we see that the limit is $\kappa$-noncollapsed on
all scales for some positive constant $\kappa$. In addition, if
the dimension $n=3$ then the Hamilton-Ivey pinching estimate
implies that the limiting solution must have nonnegative curvature
operator.

We call a solution to the Ricci flow an \textbf{ancient
$\kappa$-solution}\index{ancient!$\kappa$-solution} if it is
complete (either compact or noncompact) and defined on an ancient
time interval $(-\infty,T)$ with $T>0$, has nonnegative curvature
operator and bounded curvature, and is $\kappa$-noncollapsed on
all scales for some positive constant $\kappa$.

In this chapter we study ancient $\kappa$-solutions of the Ricci
flow. We will obtain crucial curvature estimates of such solutions
and determine their structures in lower dimensional cases. In
particular, Sections 6.2-6.4 give a detailed exposition of
Perelman's work in section 11 of \cite{P1} and section 1 of
\cite{P2}. We also remak that the earlier work on ancient
solutions can be found in Hamilton \cite{Ha95F}.

\section{Preliminaries}

We first present a useful geometric property (cf. Proposition 2.2 of
\cite{CZ05F}) for complete noncompact Riemannian manifolds with
nonnegative sectional curvature.

Let $(M,g_{ij})$ be an $n$-dimensional complete Riemannian
manifold and let $\varepsilon$ be a positive constant. We call an
open subset $N\subset M$ an \textbf{$\varepsilon$-neck of radius
$r$}\index{$\varepsilon$-neck of radius $r$} if $(N,
r^{-2}g_{ij})$ is $\varepsilon$-close, in the
$C^{[\varepsilon^{-1}]}$ topology, to a standard neck
$\mathbb{S}^{n-1}\times \mathbb{I}$, where $\mathbb{S}^{n-1}$ is
the round $(n-1)$-sphere with scalar curvature 1 and $\mathbb{I}$
is an interval of length $2\varepsilon^{-1}$. The following result
is, to some extent, in similar spirit of Yau's volume lower bound
estimate \cite{Y76}.

\begin{proposition}
There exists a positive constant $\varepsilon_0=\varepsilon_0(n)$
such that every complete noncompact $n$-dimensional Riemannian
manifold $(M, g_{ij})$ of nonnegative sectional curvature has a
positive constant $r_0$ such that any $\varepsilon$-neck of radius
$r$ on $(M, g_{ij})$ with $\varepsilon \leq \varepsilon_0$ must
have $r\geq r_0$.
\end{proposition}

\begin{pf} The following argument is taken from \cite{CZ05F}.
We argue by contradiction. Suppose there exist a sequence of
positive constants $\varepsilon^\alpha \rightarrow 0$ and a
sequence of $n$-dimensional complete noncompact Riemannian
manifolds $(M^\alpha, g^{\alpha}_{ij})$ such that for each fixed
$\alpha$, there exists a sequence of $\varepsilon^\alpha$-necks
$N_{k}$ of radius at most $1/k$ in $M^\alpha$ with centers $P_{k}$
divergent to infinity.

Fix a point $P$ on the manifold $M^\alpha$ and connect each $P_k$
to $P$ by a minimizing geodesic $\gamma_k$. By passing to a
subsequence we may assume the angle $\theta_{kl}$ between geodesic
$\gamma_k$ and $\gamma_l$ at $P$ is very small and tends to zero
as $k, l\rightarrow+\infty$, and the length of $\gamma_{k+1}$ is
much bigger than the length of $\gamma_k$. Let us connect $P_k$ to
$P_l$ by a minimizing geodesic $\eta_{kl}$. For each fixed $l>k$,
let $\tilde{P}_k$ be a point on the geodesic $\gamma_l$ such that
the geodesic segment from $P$ to $\tilde{P}_k$ has the same length
as $\gamma_k$ and consider the triangle $\Delta PP_k\tilde{P}_k$
in $M^\alpha$ with vertices $P$, $P_k$ and $\tilde{P}_k$. By
comparing with the corresponding triangle in the Euclidean plane
$\mathbb{R}^2$ whose sides have the same corresponding lengths,
Toponogov's comparison theorem implies
$$
d(P_k,\tilde{P}_k)\le 2\sin\(\frac{1}{2}\theta_{kl}\)\cdot
d(P_k,P).
$$
Since $\theta_{kl}$ is very small, the distance from $P_k$ to the
geodesic $\gamma_l$ can be realized by a geodesic $\zeta_{kl}$
which connects $P_k$ to a point $P_k'$ on the interior of the
geodesic $\gamma_l$ and has length at most
$2\sin(\frac{1}{2}\theta_{kl})\cdot d(P_k,P).$ Clearly the angle
between $\zeta_{kl}$ and $\gamma_l$ at the intersection point
$P_k'$ is $\frac{\pi}{2}$. Consider $\alpha$ to be fixed and
sufficiently large. We claim that for large enough $k$, each
minimizing geodesic $\gamma_l$ with $l>k$, connecting $P$ to
$P_l$, goes through the neck $N_k$.

Suppose not; then the angle between $\gamma_k$ and $\zeta_{kl}$ at
$P_k$ is close to either zero or $\pi$ since $P_k$ is in the
center of an $\varepsilon^\alpha$-neck and $\alpha$ is
sufficiently large. If the angle between $\gamma_k$ and
$\zeta_{kl}$ at $P_k$ is close to zero, we consider the triangle
$\Delta PP_kP_k'$ in $M^\alpha$ with vertices $P$, $P_k$, and
$P_k'$. Note that the length between $P_k$ and $P_k'$ is much
smaller than the lengths from $P_k$ or $P_k'$ to $P$. By comparing
the angles of this triangle with those of the corresponding
triangle in the Euclidean plane with the same corresponding
lengths and using Toponogov's comparison theorem, we find that it
is impossible. Thus the angle between $\gamma_k$ and $\zeta_{kl}$
at $P_k$ is close to $\pi$. We now consider the triangle $\Delta
P_kP_k'P_l$ in $M^\alpha$ with the three sides $\zeta_{kl}$,
$\eta_{kl}$ and the geodesic segment from $P_k'$ to $P_l$ on
$\gamma_l$. We have seen that the angle of $\Delta P_k P_k'P_l$ at
$P_k$ is close to zero and the angle at $P_k'$ is $\frac{\pi}{2}$.
By comparing with corresponding triangle
$\bar{\Delta}\bar{P_k}\bar{P_k'}\bar{P_l}$ in the Euclidean plane
$\mathbb{R}^2$ whose sides have the same corresponding lengths,
Toponogov's comparison theorem implies
$$
\angle \bar{P_l}\bar{P_k}\bar{P_k'} +\angle
\bar{P_l}\bar{P_k'}\bar{P_k} \le\angle P_lP_kP_k'+\angle
P_lP_k'P_k<\frac{3}{4}\pi.
$$
This is impossible since the length between $\bar{P_k}$ and
$\bar{P_k'}$ is much smaller than the length from $\bar{P_l}$ to
either $\bar{P_k}$ or $\bar{P_k'}$. So we have proved each
$\gamma_l$ with $l>k$ passes through the neck $N_k$.

Hence by taking a limit, we get a geodesic ray $\gamma$ emanating
from $P$ which passes through all the necks $N_k$, $k = 1, 2,
\ldots,$  except a finite number of them. Throwing these finite
number of necks away, we may assume $\gamma$ passes through all
necks $N_k$, $k=1,2,\ldots.$ Denote the center sphere of $N_k$ by
$S_k$, and their intersection points with $\gamma$ by  $p_{k}\in
S_k\cap \gamma$, $k=1,2,\ldots.$

Take a sequence of points  $\gamma(m)$ with $m=1,2,\ldots.$ For
each fixed neck $N_k$, arbitrarily choose a point $q_{k}\in N_k$
near the center sphere $S_k$ and draw a geodesic segment
$\gamma^{km}$ from $q_{k}$ to $\gamma(m)$. Now we claim that for
any neck $N_l$ with $l>k$, $\gamma^{km}$ will pass through $N_l$
for all sufficiently large $m$.

We argue by contradiction. Let us place all the necks $N_i$
horizontally so that the geodesic $\gamma$ passes through each
$N_i$ from the left to the right. We observe that the geodesic
segment $\gamma^{km}$ must pass through the right half of $N_k$;
otherwise $\gamma^{km}$ cannot be minimal. Then for large enough
$m$, the distance from $p_{l}$ to the geodesic segment
$\gamma^{km}$ must be achieved by the distance from $p_l$ to some
interior point ${p_k}'$ of $\gamma^{km}$. Let us draw a minimal
geodesic $\eta$ from $p_{l}$ to the interior point ${p_k}'$ with
the angle at the intersection point ${p_k}'\in \eta\cap
\gamma^{km}$ to be $\frac{\pi}{2}.$ Suppose the claim is false.
Then the angle between $\eta$ and $\gamma$ at $p_{l}$ is close to
$0$ or $\pi$ since $\varepsilon^\alpha$ is small.

If the angle between $\eta$ and $\gamma$ at $p_{l}$ is close to
$0$, we consider the triangle ${\Delta} {p}_{l}{p_k}'{\gamma}(m)$
and construct a comparison triangle $ \bar{\Delta}
\bar{p}_{l}\bar{p_k}'\bar{\gamma}(m)$ in the plane with the same
corresponding length. Then by Toponogov's comparison theorem, we
see the sum of the inner angles of the comparison triangle $
\bar{\Delta} \bar{p}_{l}\bar{p_k}'\bar{\gamma}(m)$ is less than
$3\pi/4$, which is impossible.

If the angle between $\eta$ and $\gamma$ at $p_{l}$ is close to
$\pi$, by drawing a minimal geodesic $\xi$ from $q_k$ to $p_{l}$,
we see that $\xi$ must pass through the right half of $N_k$ and
the left half of $N_l$; otherwise $\xi$ cannot be minimal. Thus
the three inner angles of the triangle $\Delta p_{l}{p_k}'q_k$ are
almost $0,\pi/2$, and $0$ respectively. This is also impossible by
the Toponogov comparison theorem.

Hence we have proved that the geodesic segment $\gamma^{km}$
passes through $N_l$ for $m$ large enough.

Consider the triangle $\Delta p_{k}q_k\gamma(m)$ with two long
sides $\overline{p_{k}\gamma(m)}(\subset\gamma)$ and
$\overline{q_{k}\gamma(m)}(= \gamma^{km})$. For any $s>0$, choose
points ${\tilde{p}_{k}}$ on $\overline{p_{k}\gamma(m)}$ and
${\tilde{q}_{k}}$ on $\overline{q_{k}\gamma(m)}$ with
$d(p_{k},{\tilde{p}_{k}})=d(q_{k},{\tilde{q}_{k}})=s$. By
Toponogov's comparison theorem, we have {\small
\begin{align*}
&  \(\frac{d({\tilde{p}_{k}},{\tilde{q}_{k}})}{d(p_{k},q_{k})}\)^{2}\\
&  =\frac{d({\tilde{p}_{k}},\gamma(m))^{2}
+d({\tilde{q}_{k}},\gamma(m))^{2}-2d({\tilde{p}_{k}},
\gamma(m))d({\tilde{q}_{k}}, \gamma(m))\cos\bar{\measuredangle}
({\tilde{p}_{k}}\gamma(m){\tilde{q}_{k}})}
{d({p_{k}},\gamma(m))^{2}+d({q_{k}},\gamma(m))^{2}-2d({p_{k}},
\gamma(m))d({q_{k}},\gamma(m))\cos\bar{\measuredangle}
({p_{k}}\gamma(m){q_{k}})}\\[1mm]
&  \geq \frac{d({\tilde{p}_{k}},\gamma(m))^{2}+d({\tilde{q}_{k}},
\gamma(m))^{2}-2d({\tilde{p}_{k}},\gamma(m))d({\tilde{q}_{k}},
\gamma(m))\cos\bar{\measuredangle}
({\tilde{p}_{k}}\gamma(m){\tilde{q}_{k}})}
{d({p_{k}},\gamma(m))^{2}+d({q_{k}},\gamma(m))^{2}
-2d({p_{k}},\gamma(m))d({q_{k}},
\gamma(m))\cos\bar{\measuredangle}
({\tilde{p}_{k}}\gamma(m){\tilde{q}_{k}})}\\[1mm]
&  = \frac{(d({\tilde{p}_{k}},\gamma(m)) -
d({\tilde{q}_{k}},\gamma(m)))^{2}+2d({\tilde{p}_{k}},
\gamma(m))d({\tilde{q}_{k}}, \gamma(m))(1 -
\cos\bar{\measuredangle}
({\tilde{p}_{k}}\gamma(m){\tilde{q}_{k}}))}
{(d({\tilde{p}_{k}},\gamma(m)) -
d({\tilde{q}_{k}},\gamma(m)))^{2}+2d({{p}_{k}},\gamma(m))d({{q}_{k}},
\gamma(m))(1 - \cos\bar{\measuredangle}
({\tilde{p}_{k}}\gamma(m){\tilde{q}_{k}}))}\\[1mm]
&  \geq \frac{d({\tilde{p}_{k}},\gamma(m))d({\tilde{q}_{k}},
\gamma(m))} {d({p_{k}},\gamma(m))d({q_{k}}, \gamma(m))}\\[1mm]
&   \rightarrow 1
\end{align*}
} as $m\rightarrow\infty$, where $\bar{\measuredangle}
({{p}_{k}}\gamma(m){{q}_{k}})$ and $\bar{\measuredangle}
({\tilde{p}_{k}}\gamma(m){\tilde{q}_{k}})$ are the corresponding
angles of the comparison triangles.

Letting $m\rightarrow \infty$, we see that $\gamma^{km}$ has a
convergent subsequence whose limit $\gamma^{k}$ is a geodesic ray
passing through all $N_l$ with $l > k$. Let us denote by
$p_{j}=\gamma(t_j), j=1, 2, \ldots$. From the above computation,
we deduce that
$$
d(p_{k},q_{k})\leq d(\gamma(t_k+s),\gamma^{k}(s))
$$
for all $s>0$.

Let $\varphi(x)=\lim_{t\rightarrow+\infty}(t-d(x,\gamma(t)))$ be
the Busemann function constructed from the ray $\gamma$. Note that
the level set $\varphi^{-1}(\varphi(p_{j}))\cap N_j$ is close to
the center sphere $S_j$ for any $j = 1, 2, \ldots$. Now let $q_k$
be any fixed point in $\varphi^{-1}(\varphi(p_{k}))\cap N_k$. By
the definition of Busemann function $\varphi$ associated to the
ray $\gamma$, we see that
$\varphi(\gamma^{k}(s_1))-\varphi(\gamma^{k}(s_2))=s_1-s_2$ for
any $s_1$, $s_2\geq0$. Consequently, for each $l>k$, by choosing
$s=t_l-t_k$, we see $\gamma^{k}(t_l-t_k)\in
\varphi^{-1}(\varphi(p_{l}))\cap N_l.$ Since
$\gamma(t_k+t_l-t_k)=p_{l}$, it follows that
$$
d(p_{k},q_{k})\leq d(p_{l},\gamma^{k}(s)).
$$
with $s =t_l-t_k >0$. This implies that the diameter of
$\varphi^{-1}(\varphi(p_{k}))\cap N_k$ is not greater than the
diameter of $\varphi^{-1}(\varphi(p_{l}))\cap N_l$ for any $l>k$,
which is a contradiction for $l$ much larger than $k$.

Therefore we have proved the proposition.
\end{pf}

In \cite{Ha95F}, Hamilton discovered an important repulsion
principle (cf. Theorem 21.4 of \cite{Ha95F}) about the influence
of a bump of strictly positive curvature in a complete noncompact
manifold of nonnegative sectional curvature. Namely minimal
geodesic paths that go past the bump have to avoid it. As a
consequence he obtained a finite bump theorem (cf. Theorem 21.5 of
\cite{Ha95F}) that gives a bound on the number of bumps of
curvature.

Let $M$ be a complete noncompact Riemannian manifold with
nonnegative sectional curvature $K\ge 0$. A geodesic ball $B(p,r)$
of radius $r$ centered at a point $p\in M$ is called a {\bf
curvature $\beta$-bump}\index{curvature $\beta$-bump} if sectional
curvature $K\ge \beta/r^2$ at all points in the ball. The ball
$B(p,r)$ is called {\bf $\lambda$-remote}\index{$\lambda$-remote}
from an origin $O$ if $d(p, O)\ge \lambda r$.

\medskip
{\bf Finite Bump Theorem}\index{finite bump theorem} (Hamilton
\cite{Ha95F}){\bf .} \emph{ For every $\beta >0$ there exists
$\lambda<\infty$ such that in any complete manifold of nonnegative
sectional curvature there are at most a finite number of disjoint
balls which are $\lambda$-remote curvature $\beta$-bumps.}

\medskip
This finite bump theorem played an important role in Hamilton's
study of the behavior of singularity models at infinity and in the
dimension reduction argument he developed for the Ricci flow (cf.
Section 22 of \cite{Ha95F}, see also \cite{CTZ04} for application
to the K\"ahler-Ricci flow and uniformization problem in complex
dimension two). A special consequence of the finite bump theorem
is that if we have a complete noncompact solution to the Ricci
flow on an ancient time interval $-\infty<t<T$ with $T>0$
satisfying certain local injectivity radius bound, with curvature
bounded at each time and with asymptotic scalar curvature ratio
$A=\limsup Rs^2=\infty$, then we can find a sequence of points
$p_j$ going to $\infty$ (as in the following Lemma 6.1.3) such
that a cover of the limit of dilations around these points at time
$t=0$ splits as a product with a flat factor. The following
result, somewhat known in Alexandrov space theory (e.g., Perelman
wrote in \cite{P1} (page 29, line 3) ``this follows from a
standard argument based on the Aleksandrov-Toponogov concavity"
and the essentially same statement also appeared in \cite{KL},
\cite{CZ05F} and \cite{CLN}), is in similar spirit as Hamilton's
finite bumps theorem and its consequence. The advantage is that we
will get in the limit of dilations a product of the line with a
lower dimensional manifold, instead of a quotient of such a
product.

\begin{proposition}
Suppose $(M,g_{ij})$ is a complete $n$-dimensional Riemannian
manifold with nonnegative sectional curvature. Let $P\in M$ be
fixed, and $P_k\in M$ a sequence of points and $\lambda_k$ a
sequence of positive numbers with $d(P,P_k)\rightarrow +\infty$
and $\lambda_kd(P,P_k)\rightarrow +\infty$. Suppose also that the
marked manifolds $(M,\lambda_k^2g_{ij},P_k)$ converge in the
$C^{\infty}_{loc}$ topology to a Riemannian manifold
$\widetilde{M}$. Then the limit $\widetilde{M}$ splits
isometrically as the metric product of the form $\mathbb{R}\times
N$, where $N$ is a Riemannian manifold with nonnegative sectional
curvature.
\end{proposition}

\begin{pf}
We now follow an argument given in \cite{CZ05F}. Let us denote by
$|OQ|=d(O,Q)$ the distance between two points $O,Q\in M$. Without
loss of generality, we may assume that for each $k$, \be
1+2|PP_k|\leq|PP_{k+1}|.
\ee Draw a minimal geodesic $\gamma_k$ from $P$ to $P_k$ and a
minimal geodesic $\sigma_k$ from $P_k$ to $P_{k+1}$, both
parametrized by arclength. We may further assume \be
\theta_k=|\measuredangle(\dot{\gamma}_k(0),\dot{\gamma}_{k+1}(0))
|<\frac{1}{k}.
\ee

By assumption, the sequence $(M,\lambda^2_kg_{ij},P_k)$ converges
(in the $C^{\infty}_{loc}$ topology) to a Riemannian manifold
$(\widetilde{M},\widetilde{g}_{ij},\widetilde{P})$ with
nonnegative sectional curvature. By a further choice of
subsequences, we may also assume $\gamma_k$ and $\sigma_k$
converge to geodesic rays $\widetilde{\gamma}$ and
$\widetilde{\sigma}$ starting at $\widetilde{P}$ respectively. We
will prove that $\tilde{\gamma}\cup\tilde{\sigma}$ forms a line in
$\widetilde{M}$, and then by the Toponogov splitting theorem
\cite{Mi} the limit $\widetilde{M}$ must be splitted as
$\mathbb{R}\times N$.

We argue by contradiction. Suppose
$\widetilde{\gamma}\cup\widetilde{\sigma}$ is not a line; then for
each $k$, there exist two points $A_k\in \gamma_k$ and $B_k\in
\sigma_k$ such that as $k\rightarrow +\infty$, \be \left\{
\begin{array}{llll}
\lambda_kd(P_k,A_k)\rightarrow A>0,\\[2mm]
\lambda_kd(P_k,B_k)\rightarrow B>0,\\[2mm]
\lambda_kd(A_k,B_k)\rightarrow C>0,\\[2mm]
\mbox{but  }A+B>C.\\
\end{array}\right. 
\ee

\begin{center}
\setlength{\unitlength}{2mm}
\begin{picture}(60,20)
\linethickness{1pt} \put(0,0){\line(6,1){60}}

\put(0,0){\line(2,1){25}} \thicklines

\qbezier(25,12.5)(40,11.5)(60,10)

\put(20,10){\line(6,1){12}}

\put(0,-2){\makebox(2,1)[c]{$P$}}

\put(25,15){\makebox(2,1)[c]{$P_k$}}

\put(60,8){\makebox(2,1)[c]{$P_{k+1}$}}

\put(20,8){\makebox(2,1)[c]{$A_k$}}

\put(25,9){\makebox(2,1)[c]{$\delta_k$}}

\put(32,10){\makebox(2,1)[c]{$B_k$}}

\put(42,12){\makebox(2,1)[c]{$\sigma_k$}}

\put(12,7){\makebox(2,1)[c]{$\gamma_k$}}

\end{picture}

\end{center}

\bigskip\bigskip
Now draw a minimal geodesic $\delta_k$ from $A_k$ to $B_k$.
Consider comparison triangles
$\bar{\triangle}{\bar{P}_k\bar{P}\bar{P}_{k+1}}$ and
$\bar{\triangle}{\bar{P}_k\bar{A}_k\bar{B}_k}$ in $\mathbb{R}^2$
with
$$
|\bar{P}_k\bar{P}|=|P_kP|, |\bar{P}_k\bar{P}_{k+1}|=|P_kP_{k+1}|,
|\bar{P}\bar{P}_{k+1}|=|PP_{k+1}|,
$$
$$\mbox{ and }\;|\bar{P}_k\bar{A}_k|
=|P_kA_k|, |\bar{P}_k\bar{B}_k|=|P_kB_k|,
|\bar{A}_k\bar{B}_k|=|A_kB_k|.
$$
By Toponogov's comparison theorem \cite{BGP}, we have \be
\measuredangle\bar{A}_k\bar{P}_k\bar{B}_k\geq
\measuredangle \bar{P} \bar{{P}_k}\bar{P}_{k+1}. 
\ee On the other hand, by (6.1.2) and using Toponogov's comparison
theorem again, we have \be
\measuredangle\bar{P}_k\bar{P}\bar{P}_{k+1} \leq \measuredangle
P_kPP_{k+1}<\frac{1}{k},
\ee and since $|\bar{P}_k\bar{P}_{k+1}|>|\bar{P}\bar{P}_k|$ by
(6.1.1), we further have \be
\measuredangle\bar{P}_k\bar{P}_{k+1}\bar{P} \leq
\measuredangle\bar{P}_k\bar{P}\bar{P}_{k+1}<\frac{1}{k}.
\ee Thus the above inequalities (6.1.4)-(6.1.6) imply that
$$
\measuredangle\bar{A}_k\bar{P}_k\bar{B}_k>\pi-\frac{2}{k}.
$$
Hence \be |\bar{A}_k\bar{B}_k|^2\geq
|\bar{A}_k\bar{P}_k|^2+|\bar{P}_k\bar{B}_k|^2
-2|\bar{A}_k\bar{P}_k|\cdot|\bar{P}_k\bar{B}_k|\cos\(\pi-\frac{2}{k}\).
\ee

Multiplying the above inequality by $\lambda^2_k$ and letting
$k\rightarrow+\infty$, we get $$C\geq A+B$$ which contradicts
(6.1.3).

Therefore we have proved the proposition.
\end{pf}

Let $M$ be an $n$-dimensional complete noncompact Riemannian
manifold. Pick an origin $O\in M$. Let $s$ be the geodesic
distance to the origin $O$ of $M$, and $R$ the scalar curvature.
Recall that in Chapter 4 we have defined the {\bf asymptotic
scalar curvature ratio}\index{asymptotic scalar curvature ratio}
$$
A=\limsup_{s\rightarrow+\infty}Rs^2.
$$

We now state a useful lemma of Hamilton (cf. Lemma 22.2 of
\cite{Ha95F}) about picking local (almost) maximum curvature
points at infinity.

\begin{lemma} [{Hamilton \cite{Ha95F}}]
 Given a complete noncompact Riemannian manifold with
bounded curvature and with asymptotic scalar curvature ratio
$$
A=\limsup_{s\rightarrow+\infty}Rs^2=+\infty,
$$
we can find a sequence of points $x_j$ divergent to infinity, a
sequence of radii $r_j$ and a sequence of positive numbers
$\delta_j\rightarrow0$ such that
\begin{itemize}
\item[(i)] $R(x)\leq(1+\delta_j)R(x_j)$ for all $x$ in the ball
$B(x_j,r_j)$ of radius $r_j$ around $x_j$, \item[(ii)]
$r^2_jR(x_j)\rightarrow+\infty$, \item[(iii)]
$\lambda_j=d(x_j,O)/r_j\rightarrow+\infty$, \item[(iv)]
the balls $B(x_j,r_j)$ are disjoint, \\
where $d(x_j,O)$ is the distance of $x_j$ from the origin $O$.
\end{itemize}
\end{lemma}

\begin{pf} The proof is essentially from Hamilton \cite{Ha95F}.
Pick a sequence of positive numbers $\epsilon_j\rightarrow0$, then
choose $A_j\rightarrow+\infty$ so that
$A_j\epsilon^2_j\rightarrow+\infty$. Let $\sigma_j$ be the largest
number such that
$$
\sup\{ R(x)d(x,O)^2\ |\ d(x,O)\leq\sigma_j\}\leq A_j.
$$
Then there exists some $y_j\in M$ such that
$$
R(y_j)d(y_j,O)^2=A_j\ \ \ \mbox{and}\ \ \ d(y_j,O)=\sigma_j.
$$
Now pick $x_j\in M$ so that $d(x_j,O)\geq\sigma_j$ and
$$
R(x_j)\geq\frac{1}{1+\epsilon_j}\sup\{ R(x)\ |\
d(x,O)\geq\sigma_j\}.
$$
Finally pick $r_j=\epsilon_j\sigma_j$. We check the properties
(i)-(iv) as follows.

\smallskip
(i)\ \ \ If $x\in B(x_j,r_j)\cap\{d(\cdot,O)\geq\sigma_j\}$, we
have
$$
R(x)\leq(1+\epsilon_j)R(x_j)
$$
by the choice of the point $x_j$; while if $x\in
B(x_j,r_j)\cap\{d(\cdot,O)\leq\sigma_j\}$, we have
\begin{displaymath}
\begin{split}
   R(x)&\leq A_j/d(x,O)^2\\[1mm]
       &\leq \frac{1}{(1-\epsilon_j)^2}(A_j/\sigma_j^2)\\[1mm]
       &=  \frac{1}{(1-\epsilon_j)^2}R(y_j)\\[1mm]
       &\leq \frac{(1+\epsilon_j)}{(1-\epsilon_j)^2}R(x_j),
\end{split}
\end{displaymath}
since $d(x,O)\geq
d(x_j,O)-d(x,x_j)\geq\sigma_j-r_j=(1-\epsilon_j)\sigma_j$. Thus we
have obtained
$$
R(x)\leq(1+\delta_j)R(x_j),\ \ \ \forall\ x\in B(x_j,r_j),
$$
where
$\delta_j=\frac{(1+\epsilon_j)}{(1-\epsilon_j)^2}-1\rightarrow0$
as $j\rightarrow+\infty$.

\smallskip
(ii)\ \ \ By the choices of $r_j$, $x_j$ and $y_j$, we have
\begin{displaymath}
\begin{split}
r^2_jR(x_j)&= \epsilon_j^2\sigma_j^2R(x_j)\\[1mm]
      &\geq \epsilon_j^2\sigma_j^2
\left[\frac{1}{1+\epsilon_j}R(y_j)\right]\\[1mm]
      &= \frac{\epsilon_j^2}{1+\epsilon_j}A_j\rightarrow+\infty,\
      \ \ \mbox{as }\; j\rightarrow+\infty.
\end{split}
\end{displaymath}

\smallskip
(iii)\ \ \ Since $d(x_j,O)\geq\sigma_j={r_j}/{\epsilon_j}$, it
follows that $\lambda_j={d(x_j,O)}/{r_j}\rightarrow+\infty$ as
$j\rightarrow+\infty$.

\smallskip
(iv)\ \ \ For any $x\in B(x_j,r_j)$, the distance from the origin
\begin{displaymath}
\begin{split}
   d(x,O)&\geq d(x_j,O)-d(x,x_j)\\[1mm]
         &\geq \sigma_j-r_j\\[1mm]
         &= (1-\epsilon_j)\sigma_j\rightarrow+\infty,\ \ \mbox{as }\;
         j\rightarrow+\infty.
\end{split}
\end{displaymath}
Thus any fixed compact set does not meet the balls $B(x_j,r_j)$
for large enough $j$. If we pass to a subsequence, the balls will
all avoid each other.
\end{pf}

The above point picking lemma of Hamilton, as written down in
Lemma 22.2 of \cite{Ha95F}, requires the curvature of the manifold
to be bounded. When the manifold has unbounded curvature, we will
appeal to the following simple lemma.

\begin{lemma}
Given a complete noncompact Riemannian manifold with unbounded
curvature, we can find a sequence of points $x_j$ divergent to
infinity such that for each positive integer $j$, we have
$|Rm(x_{j})|\geq j$, and
$$
|Rm(x)|\leq4|Rm(x_{j})|
$$
for $x\in B(x_j,\frac{j}{\sqrt{|Rm(x_{j})|}})$.
\end{lemma}

\begin{pf}
Each $x_j$ can be constructed as a limit of a finite sequence
$\{y_i\}$, defined as follows. Let $y_0$ be any fixed point with
$|Rm(y_{0})| \geq j$. Inductively, if $y_i$ cannot be taken as
$x_j$, then there is a $y_{i+1}$ such that
\begin{equation*}
 \left\{
\begin{aligned}
|Rm(y_{i+1})|&  > 4|Rm(y_{i})|,\\
d(y_i,y_{i+1}) &  \leqslant \frac{j}{\sqrt{|Rm(y_{i})|}}.
\end{aligned}
\right.
\end{equation*}
Thus we have
$$
|Rm(y_{i})|> 4^i|Rm(y_{0})|\geq 4^ij,$$
$$d(y_i,y_{0})\leq
j\sum^i_{k=1}\frac{1}{\sqrt{4^{k-1}j}}<2\sqrt{j}.$$ Since the
manifold is smooth, the sequence $\{y_i\}$ must be finite. The
last element fits.
\end{pf}

\section{Asymptotic Shrinking Solitons}
The main purpose of this section is to prove a result of Perelman
(cf. Proposition 11.2 of \cite{P1}) on the asymptotic shapes of
ancient $\kappa$-solutions as time $t \rightarrow -\infty$.

We begin with the study of the asymptotic behavior of an ancient
$\kappa$-solution $g_{ij}(x,t)$, on $M\times (-\infty,T)$ with
$T>0$, to the Ricci flow as $t\rightarrow-\infty$.

Pick an arbitrary point $(p,t_0)\in M\times(-\infty,0]$ and recall
from Chapter 3 that
$$
\tau=t_0-t,\ \ \ \mbox{for }\; t<t_0,
$$
\begin{multline*}
l(q,\tau)=\frac{1}{2\sqrt{\tau}}\inf \bigg\{
\int_0^\tau\sqrt{s}\Big(R(\gamma(s),t_0-s) \\
+|\dot{\gamma}(s)|^2_{g_{ij}(t_0-s)}\Big)ds \left|
\begin{array}{c}
\gamma:[0,\tau]\rightarrow M\ \mbox{with}\\[1mm]
\gamma(0)=p,\ \gamma(\tau)=q\
\end{array}\bigg\}\right.
\end{multline*}
and
$$
\tilde{V}(\tau)=\int_M(4\pi\tau)^{-\frac{n}{2}}
\exp(-l(q,\tau))dV_{t_0-\tau}(q).
$$
We first observe that Corollary 3.2.6 also holds for the general
complete manifold $M$. Indeed, since the scalar curvature is
nonnegative, the function $\bar{L}(\cdot,\tau)=4\tau
l(\cdot,\tau)$ achieves its minimum on $M$ for each fixed
$\tau>0$. Thus the same argument in the proof of Corollary 3.2.6
shows there exists $q=q(\tau)$ such that \be
l(q(\tau),\tau)\leq\frac{n}{2} 
\ee for each $\tau>0$.

Recall from (3.2.11)-(3.2.13), the Li-Yau-Perelman distance $l$
satisfies the following
\begin{align}
\frac{\partial}{\partial\tau}l
&=-\frac{l}{\tau}+R+\frac{1}{2\tau^{3/2}}K, \\  
|\nabla l|^2&=-R+\frac{l}{\tau}-\frac{1}{\tau^{3/2}}K, \\  
\Delta l& \leq-R+\frac{n}{2\tau}-\frac{1}{2\tau^{3/2}}K,  
\end{align}
and the equality in (6.2.4) holds everywhere if and only if we are
on a gradient shrinking soliton. Here $K=\int_0^\tau
s^{3/2}Q(X)ds,$ $Q(X)$ is the trace Li-Yau-Hamilton quadratic
given by
$$
Q(X)=-R_{\tau} -\frac{R}{\tau}-2 \langle\nabla
R,X\rangle+2\Ric(X,X)
$$
and $X$ is the tangential (velocity) vector field of an
$\mathcal{L}$-shortest curve $\gamma: [0,\tau]\rightarrow M$
connecting $p$ to $q$.

By applying the trace Li-Yau-Hamilton inequality (Corollary 2.5.5)
to the ancient $\kappa$-solution, we have
\begin{align*}
Q(X)&  = -R_{\tau} -\frac{R}{\tau}-2 \langle\nabla
R,X\rangle+2\Ric(X,X) \\
&  \geq -\frac{R}{\tau}
\end{align*}
and hence
\begin{align*}
K&  = \int_0^{\tau}s^{3/2}Q(X)ds   \\
&  \geq -\int_0^{\tau}\sqrt{s}Rds \\
&  \geq -L(q,\tau).
\end{align*}
Thus by (6.2.3) we get \be |\nabla l|^2+R\leq \frac{3l}{\tau}.
\ee 

\medskip
We are now state and prove the following

\begin{theorem}[{Perelman \cite{P1}}]
Let $g_{ij}(\cdot,t), -\infty<t<T$ with some $T>0$, be a nonflat
ancient $\kappa$-solution for some $\kappa >0$. Then there exist a
sequence of points $q_k$ and a sequence of times $t_k \rightarrow
-\infty$ such that the scalings of $g_{ij}(\cdot,t)$ around $q_k$
with factor $|t_k|^{-1}$ and with the times $t_k$ shifting to the
new time zero converge to a nonflat gradient shrinking soliton in
$C_{\rm loc}^\infty$ topology.
\end{theorem}

\begin{pf} The proof basically follows the argument of Perelman
(11.2 of \cite{P1}).  Clearly, we may assume that the nonflat
ancient $\kappa$-solution is not a gradient shrinking soliton. For
the arbitrarily fixed $(p,t_0)$, let $q(\tau) (\tau=t_0 -t)$ be
chosen as in (6.2.1) with $l(q(\tau),\tau) \leq \frac{n}{2}$. We
only need to show that the scalings of $g_{ij}(\cdot,t)$ around
$q(\tau)$ with factor $\tau^{-1}$ converge along a subsequence of
$\tau \rightarrow +\infty$ to a nonflat gradient shrinking soliton
in the $C^{\infty}_{\rm loc}$ topology.

We first claim that for any $A \geq 1$, one can find
$B=B(A)<+\infty$ such that for every $\bar{\tau}>1$ there holds
\be l(q,\tau)\leq B \; \mbox{ and }\;  \tau R(q,t_0-\tau)\leq B,
\ee whenever $\frac{1}{2}\bar{\tau}\leq \tau \leq A\bar{\tau}$ and
$d_{t_0-\frac{\bar{\tau}}{2}}^2 (q,q(\frac{\bar{\tau}}{2}))\leq
A\bar{\tau}.$

Indeed, by using (6.2.5) at $\tau=\frac{\bar{\tau}}{2}$, we have
\begin{align}
\sqrt{l(q,\frac{\bar{\tau}}{2})} &  \leq \sqrt{\frac{n}{2}}+
\sup\{|\nabla \sqrt{l}|\}\cdot d_{t_0-\frac{\bar{\tau}}{2}}
\(q,q\(\frac{\bar{\tau}}{2}\)\) \\
&  \leq \sqrt{\frac{n}{2}}+ \sqrt{\frac{3}{2\bar{\tau}}}\cdot
\sqrt{A\bar{\tau}} \nn\\
&  = \sqrt{\frac{n}{2}}+\sqrt{\frac{3A}{2}},\nn
\end{align}  
and
\begin{align}
R\(q,t_0-\frac{\bar{\tau}}{2}\)&  \leq
\frac{3l(q,\frac{\bar{\tau}}{2})}{(\frac{\bar{\tau}}{2})}\\
&  \leq
\frac{6}{\bar{\tau}}\(\sqrt{\frac{n}{2}}+\sqrt{\frac{3A}{2}}\)^2,\nn
\end{align}  
for $q\in
B_{t_0-\frac{\bar{\tau}}{2}}(q(\frac{\bar{\tau}}{2}),\sqrt{A\bar{\tau}})$.
Recall that the Li-Yau-Hamilton inequality implies that the scalar
curvature of the ancient solution is pointwise nondecreasing in
time. Thus we know from (6.2.8) that \be \tau R(q,t_0-\tau)\leq
6A\(\sqrt{\frac{n}{2}}+\sqrt{\frac{3A}{2}}\)^2 
\ee whenever $\frac{1}{2}\bar{\tau}\leq \tau \leq A\bar{\tau}$ and
$d_{t_0-\frac{\bar{\tau}}{2}}^2 (q,q(\frac{\bar{\tau}}{2}))\leq
A\bar{\tau}.$

By (6.2.2) and (6.2.3) we have
$$
\frac{\partial l}{\partial \tau}+\frac{1}{2}|\nabla
l|^2=-\frac{l}{2\tau}+\frac{R}{2}.
$$
This together with (6.2.9) implies that
\begin{align*}
\frac{\partial l}{\partial \tau} &\leq
-\frac{l}{2\tau}+\frac{3A}{\tau}
\(\sqrt{\frac{n}{2}}+\sqrt{\frac{3A}{2}}\)^2 \\
\intertext{i.e.,} \frac{\partial}{\partial\tau}(\sqrt{\tau}l)
&\leq\frac{3A}{\sqrt{\tau}}\(\sqrt{\frac{n}{2}}+\sqrt{\frac{3A}{2}}\)^2
\end{align*}
whenever $\frac{1}{2}\bar{\tau}\leq \tau \leq A\bar{\tau}$ and
$d_{t_0-\frac{\bar{\tau}}{2}}^2
\(q,q\(\frac{\bar{\tau}}{2}\)\)\leq A\bar{\tau}.$ Hence by
integrating this differential inequality, we obtain
$$
\sqrt{\tau}l(q,\tau)-\sqrt{\frac{\bar{\tau}}{2}}l\(q,\frac{\bar{\tau}}{2}\)
\leq 6A\(\sqrt{\frac{n}{2}}+\sqrt{\frac{3A}{2}}\)^2\sqrt{\tau}
$$
and then by (6.2.7),
\begin{align}
l(q,\tau)&  \leq
l\(q,\frac{\bar{\tau}}{2}\)+6A\(\sqrt{\frac{n}{2}}+\sqrt{\frac{3A}{2}}\)^2\\
&  \leq 7A\(\sqrt{\frac{n}{2}}+\sqrt{\frac{3A}{2}}\)^2 \nn
\end{align}  
whenever $\frac{1}{2}\bar{\tau}\leq \tau \leq A\bar{\tau}$ and
$d_{t_0-\frac{\bar{\tau}}{2}}^2 (q,q(\frac{\bar{\tau}}{2}))\leq
A\bar{\tau}.$  So we have proved claim (6.2.6).

Recall that $g_{ij}(\tau)=g_{ij}(\cdot,t_0-\tau) $ satisfies
$(g_{ij})_\tau=2R_{ij}$. Let us take the scaling of the ancient
$\kappa$-solution around $q(\frac{\bar{\tau}}{2})$ with factor
$(\frac{\bar{\tau}}{2})^{-1}$, i.e.,
$$
\tilde{g}_{ij}(s)=\frac{2}{\bar{\tau}}g_{ij}
\(\cdot,t_0-s\frac{\bar{\tau}}{2}\)
$$
where $s\in [0,+\infty)$. Claim (6.2.6) says that for all
$s\in[1,2A]$ and all $q$ such that
dist$^2_{\tilde{g}_{ij}(1)}(q,q(\frac{\bar{\tau}}{2}))\leq A,$ we
have $\tilde{R}(q,s)=\frac{\bar
\tau}{2}R(q,t_0-s\bar{\frac{\tau}{2}})\leq B$. Now taking into
account the $\kappa$-noncollapsing assumption and Theorem 4.2.2,
we can use Hamilton's compactness theorem (Theorem 4.1.5) to
obtain a sequence $\bar{\tau}_k\rightarrow +\infty$ such that the
marked evolving manifolds
$(M,\tilde{g}_{ij}^{(k)}(s),q(\frac{\bar{\tau}_k}{2})),$ with
$\tilde{g}_{ij}^{(k)}(s)=\frac{2}{\bar{\tau}_k}g_{ij}
(\cdot,t_0-s\frac{\bar{\tau}_k}{2})$ and $s \in [1,+\infty)$,
converge to a manifold $(\bar{M},\bar{g}_{ij}(s),\bar{q})$ with $s
\in [1,+\infty)$, where $\bar{g}_{ij}(s)$ is also a solution to
the Ricci flow on $\bar{M}$.

Denote by $\tilde{l}_k$ the corresponding Li-Yau-Perelman distance
of $\tilde{g}_{ij}^{(k)}(s)$. It is easy to see that
$\tilde{l}_k(q,s)=l(q,\frac{\bar{\tau}_k}{2}s),$ for $s\in
[1,+\infty).$ From (6.2.5), we also have \be |\nabla
\tilde{l}_k|^2_{\tilde{g}_{ij}^{(k)}}+\tilde{R}^{(k)} \leq
6\tilde{l}_k,
\ee where $\tilde{R}^{(k)}$ is the scalar curvature of the metric
$\tilde{g}_{ij}^{(k)}$. Claim (6.2.6) says that $\tilde{l}_k$ are
uniformly bounded on compact subsets of $M \times [1,+\infty)$
(with the corresponding origins $q(\frac{\bar{\tau}_k}{2})$). Thus
the above gradient estimate (6.2.11) implies that the functions
$\tilde{l}_k$ tend (up to a subsequence) to a function $\bar{l}$
which is a locally Lipschitz function on $\bar{M}$.

We know from (6.2.2)-(6.2.4) that the Li-Yau-Perelman distance
$\tilde{l}_k$ satisfies the following inequalities: \be
(\tilde{l}_k)_s-\Delta
\tilde{l}_k+|\nabla\tilde{l}_k|^2-\tilde{R}^{(k)}
+\frac{n}{2s}\geq 0, 
\ee \be 2\Delta\tilde{l}_k-|\nabla\tilde{l}_k|^2
+\tilde{R}^{(k)}+\frac{\tilde{l}_k-n}{s}\leq 0. 
\ee We next show that the limit $\bar{l}$ also satisfies the above
two inequalities in the sense of distributions. Indeed the above
two inequalities can be rewritten as \be
\(\frac{\partial}{\partial s}-\triangle+\widetilde{R}^{(k)}\)
\((4\pi s)^{-\frac{n}{2}}\exp(-\tilde{l}_k)\)\leq 0,
\ee \be -(4\triangle
-\widetilde{R}^{(k)})e^{-\frac{\tilde{l}_k}{2}}
+\frac{\tilde{l}_{k}-n}{s}e^{-\frac{\tilde{l}_k}{2}}\leq 0,
\ee in the sense of distributions. Note that the estimate (6.2.11)
implies that $\tilde{l}_k\rightarrow \bar{l}$ in the $C_{\rm
loc}^{0,\alpha}$ norm for any $0<\alpha<1.$ Thus the inequalities
(6.2.14) and (6.2.15) imply that the limit $\overline{l}$
satisfies \be \(\frac{\partial}{\partial s}-\triangle
+\overline{R}\) \((4\pi s)^{-\frac{n}{2}}\exp(-\overline{l})\)\leq
0,
\ee \be -(4\triangle-\overline{R})e^{-\frac{\overline{l}}{2}}
+\frac{\bar{l}-n}{s}e^{-\frac{\overline{l}}{2}}\leq 0,
\ee in the sense of distributions.

Denote\; by\; $\tilde{V}^{(k)}\,(s)$\; Perelman's\; reduced\;
volume\; of\; the\; scaled\; metric $\tilde{g}_{ij}^{(k)}(s)$.
Since $\tilde{l}_k(q,s)=l(q,\frac{\bar{\tau}_k}{2}s)$, we see that
$\tilde{V}^{(k)}(s)=\tilde{V}(\frac{\bar{\tau}_k}{2}s)$ where
$\tilde{V}$ is Perelman's reduced volume of the ancient
$\kappa$-solution. The monotonicity of Perelman's reduced volume
(Theorem 3.2.8) then implies that \be
\lim_{k\rightarrow\infty}\tilde{V}^{(k)}(s)=\bar{V}, \;
\mbox{ for }\; s\in[1,2],  
\ee for some nonnegative constant $\bar{V}$.

(We remark that by the Jacobian comparison theorem (Theorem
3.2.7), (3.2.18) and (3.2.19), the integrand of
$\tilde{V}^{(k)}(s)$ is bounded by
$$
(4\pi s)^{-\frac{n}{2}}\exp(-\tilde{l}_k(X,s))
\tilde{\mathcal{J}}^{(k)}(s) \leq (4\pi)^{-\frac{n}{2}}
\exp(-|X|^2)
$$
on $T_pM$, where $\tilde{\mathcal{J}}^{(k)}(s)$ is the
$\mathcal{L}$-Jacobian of the $\mathcal{L}$-exponential map of the
metric $\tilde{g}_{ij}^{(k)}(s)$ at $T_pM$. Thus we can apply the
dominant convergence theorem to get the convergence in (6.2.18).
But we are not sure whether the limiting $\bar{V}$ is exactly
Perelman's reduced volume of the limiting manifold
($\bar{M},\bar{g}_{ij}(s)$), because the points
$q(\frac{\bar{\tau}_k}{2})$ may diverge to infinity. Nevertheless,
we can ensure that $\bar{V}$ is not less than Perelman's reduced
volume of the limit.)

Note by (6.2.5) that
\begin{align}
& \tilde{V}^{(k)}(2)-\tilde{V}^{(k)}(1)\\[1mm]
& = \int_{1}^2\frac{d}{ds}(\tilde{V}^{(k)}(s))ds \nn\\[1mm]
& = \int_{1}^2ds\int_M\(\frac{\partial}{\partial s}
-\Delta+\tilde{R}^{(k)}\) \((4\pi s)^{-\frac{n}{2}}
\exp(-\tilde{l}_k)\)dV_{\tilde{g}_{ij}^{(k)}(s)}. \nn
\end{align}
Thus we deduce that in the sense of distributions, \be
\(\frac{\partial}{\partial s}-\Delta+\bar{R}\) \((4\pi
s)^{-\frac{n}{2}}\exp(-\bar{l})\)= 0,
\ee and
$$
(4\Delta-\bar{R})e^{-\frac{\bar{l}}{2}}
=\frac{\bar{l}-n}{s}e^{-\frac{\bar{l}}{2}}
$$
or equivalently, \be 2\Delta \bar{l} - |\nabla \bar{l}|^2 +
\bar{R}
+ \frac{\bar{l} -n}{s} = 0, 
\ee on $\bar{M} \times [1,2]$. Thus by applying standard parabolic
equation theory to (6.2.20) we find that $\bar{l}$ is actually
smooth. Here we used (6.2.2)-(6.2.4) to show that the equality in
(6.2.16) implies the equality in (6.2.17).

Set
$$
v=[s(2\Delta \bar{l}-|\nabla \bar{l}|^2+\bar{R})+\bar{l}-n]\cdot
(4\pi s)^{-\frac{n}{2}}e^{-\bar{l}}.
$$
Then by applying Lemma 2.6.1, we have \be
\(\frac{\partial}{\partial s}-\Delta+\bar{R}\)v
=-2s|\bar{R}_{ij}+\nabla_i\nabla_j\bar{l}
-\frac{1}{2s}\bar{g}_{ij}|^2\cdot
(4\pi s)^{-\frac{n}{2}}e^{-\bar{l}}. 
\ee We see from (6.2.21) that the LHS of the equation (6.2.22) is
identically zero. Thus the limit metric $\bar{g}_{ij}$ satisfies
\be
\bar{R}_{ij}+\nabla_i\nabla_j\bar{l}-\frac{1}{2s}\bar{g}_{ij}=0,
\ee so we have shown the limit is a gradient shrinking soliton.

To show that the limiting gradient shrinking soliton is nonflat,
we first show that the constant function $\bar{V}(s)$ is strictly
less than 1. Consider Perelman's reduced volume $\tilde{V}(\tau)$
of the ancient $\kappa$-solution. By using Perelman's Jacobian
comparison theorem (Theorem 3.2.7), (3.2.18) and (3.2.19) as
before, we have
\begin{align*}
\tilde{V}(\tau)&  = \int (4\pi
\tau)^{-\frac{n}{2}}e^{-l(X,\tau)}{\mathcal{J}}(\tau)dX \\
&  \leq \int_{T_pM}(4\pi)^{-\frac{n}{2}}e^{-|X|^2}dX \\
&  = 1.
\end{align*}

Recall that we have assumed the nonflat ancient $\kappa$-solution
is not a gradient shrinking soliton. Thus for $\tau>0,$ we must
have $\tilde{V}(\tau)<1$. Since the limiting function $\bar{V}(s)$
is the limit of $\tilde{V}(\frac{\bar{\tau}_k}{2}s)$ with
$\bar{\tau}_k\rightarrow +\infty$, we deduce that the constant
$\bar{V}(s)$ is strictly less than $1$, for $s\in[1,2].$

We now argue by contradiction. Suppose the limiting gradient
shrinking soliton $\bar{g}_{ij}(s)$ is flat. Then by (6.2.23),
$$
\nabla_i\nabla_j\bar{l} =\frac{1}{2s}\bar{g}_{ij} \; \mbox{ and }
\;\Delta\bar{l} =\frac{n}{2s}.
$$
Putting these into the identity (6.2.21), we get
$$
|\nabla\bar{l}|^2=\frac{\bar{l}}{s}
$$
Since the function $\bar{l}$ is strictly convex, it follows that
$\sqrt{4s\bar{l}}$ is a distance function (from some point) on the
complete flat manifold $\bar{M}$. From the smoothness of the
function $\bar{l}$, we conclude that the flat manifold $\bar{M}$
must be $\mathbb{R}^n$. In this case we would have its reduced
distance to be $\bar{l}$ and its reduced volume to be $1$. Since
$\bar{V}$ is not less than the reduced volume of the limit, this
is a contradiction. Therefore the limiting gradient shrinking
soliton $\bar{g}_{ij}$ is not flat.
\end{pf}

To conclude this section, we use the above theorem to derive the
classification of all two-dimensional ancient $\kappa$-solutions
which was obtained earlier by Hamilton in Section 26 of
\cite{Ha95F}.

\begin{theorem} [{Hamilton \cite{Ha95F}}]
The only nonflat ancient $\kappa$-solutions to Ricci flow on
two-dimensional manifolds are the round sphere $\mathbb{S}^2$ and
the round real projective plane $\mathbb{R}\mathbb{P}^2$.
\end{theorem}

\begin{pf}
Let $g_{ij}(x,t)$ be a nonflat ancient $\kappa$-solution defined
on $M\times (-\infty,T)$ (for some $T>0$), where $M$ is a
two-dimensional manifold. Note that the ancient $\kappa$-solution
satisfies the Li-Yau-Hamilton inequality (Corollary 2.5.5). In
particular by Corollary 2.5.8, the scalar curvature of the ancient
$\kappa$-solution is pointwise nondecreasing in time. Moreover by
the strong maximum principle, the ancient $\kappa$-solution has
strictly positive curvature everywhere.

By the above Theorem 6.2.1, we know that the scalings of the
ancient $\kappa$-solution along a sequence of points $q_k$ in $M$
and a sequence of times $t_k \rightarrow -\infty$ converge to a
nonflat gradient shrinking soliton $(\bar{M},\bar{g}_{ij}(x,t))$
with $-\infty < t \leq 0$.

We first show that the limiting gradient shrinking soliton
$(\bar{M},\bar{g}_{ij}(x,t))$ has uniformly bounded curvature.
Clearly, the limiting soliton has nonnegative curvature and is
$\kappa$-noncollapsed on all scales, and its scalar curvature is
still pointwise nondecreasing in time. Thus we only need to show
that the limiting soliton has bounded curvature at $t=0$. We argue
by contradiction. Suppose the curvature of the limiting soliton is
unbounded at $t=0$. Of course in this case the limiting soliton
$\bar{M}$ is noncompact. Then by applying Lemma 6.1.4, we can
choose a sequence of points $x_j, j = 1, 2, \ldots,$ divergent to
infinity such that the scalar curvature $\bar{R}$ of the limit
satisfies
$$
\bar{R}(x_j,0) \geq j\; \mbox{ and }\; \bar{R}(x,0) \leq
4\bar{R}(x_j,0)
$$
for all $j = 1, 2,\, \ldots,$ and $x\in
B_0(x_j,{j}/{\sqrt{\bar{R}(x_{j},0)}})$. And then by the
nondecreasing (in time) of the scalar curvature, we have
$$
\bar{R}(x,t) \leq 4\bar{R}(x_j,0),
$$
for all $j = 1, 2,\, \ldots$, $x\in
B_0(x_j,{j}/{\sqrt{\bar{R}(x_{j},0)}})$ and $t\leq0$. By combining
with Hamilton's compactness theorem (Theorem 4.1.5) and the
$\kappa$-noncollapsing, we know that a subsequence of the
rescaling solutions
$$
(\bar{M},\bar{R}(x_j,0)\bar{g}_{ij}(x,t/\bar{R}(x_j,0)),x_j),
\quad j=1, 2, \ldots,
$$
converges in the $C^{\infty}_{loc}$ topology to a nonflat smooth
solution of the Ricci flow. Then Proposition 6.1.2 implies that
the new (two-dimensional) limit must be flat. This arrives at a
contradiction. So we have proved that the limiting gradient
shrinking soliton has uniformly bounded curvature.

We next adapt an argument of Perelman (cf. the proof of Lemma 1.2
in \cite{P2}) to show that the limiting soliton is compact.
Suppose the limiting soliton is (complete and) noncompact. By the
strong maximum principle we know that the limiting soliton also
has strictly positive curvature everywhere. After a shift of the
time, we may assume that the limiting soliton satisfies the
following equation \be \nabla_i\nabla_j
f+\bar{R}_{ij}+\frac{1}{2t}\bar{g}_{ij}=0,\quad
\mbox{on }\;  -\infty<t<0, 
\ee everywhere for some function $f$. Differentiating the equation
(6.2.24) and switching the order of differentiations, as in the
derivation of (1.1.14), we get \be
\nabla_i\bar{R}=2\bar{R}_{ij}\nabla_jf. 
\ee

Fix some $t<0$, say $t=-1$, and consider a long shortest geodesic
$\gamma(s)$, $0\leq s\leq \overline{s}$. Let $x_0=\gamma(0)$ and
$X(s)=\dot{\gamma}(s)$. Let $V(0)$ be any unit vector orthogonal
to $\dot{\gamma}(0)$ and translate $V(0)$ along $\gamma(s)$ to get
a parallel vector field $V(s)$, $0\leq s\leq \overline{s}$ on
$\gamma$. Set
$$
\widehat{V}(s)=
\begin{cases}
sV(s), & \mbox{for  $0\leq s\leq 1$},\\
V(s), & \mbox{for  $1\leq s\leq \overline{s}-1$},\\
(\overline{s}-s)V(s), & \mbox{for  $\overline{s}-1\leq s\leq
\overline{s}.$}
\end{cases}
$$
It follows from the second variation formula of arclength that
$$
\int^{\overline{s}}_0(|\dot{\widehat{V}}(s)|^2
-\bar{R}(X,\widehat{V},X,\widehat{V}))ds\geq 0.
$$
Thus we clearly have
$$
\int^{\overline{s}}_0\bar{R}(X,\widehat{V},X,\widehat{V})ds\leq
{\rm const.},
$$
and then \be \int^{\overline{s}}_0\bar{\rm Ric}(X,X)ds\leq
{\rm const.} . 
\ee By integrating the equation (6.2.24) we get
$$
X(f(\gamma(\overline{s})))-X(f(\gamma(0)))
+\int^{\overline{s}}_0\bar{\rm
Ric}(X,X)ds-\frac{1}{2}\overline{s}=0
$$
and then by (6.2.26), we deduce
$$
\frac{d}{ds}(f\circ\gamma(s))\geq \frac{s}{2}-{\rm const.},
$$
$$
\mbox{and } f\circ\gamma(s)\geq \frac{s^2}{4}-{\rm const.}\cdot
s-{\rm const.}
$$
for $s>0$ large enough. Thus at large distances from the fixed
point $x_0$ the function $f$ has no critical points and is proper.
It then follows from the Morse theory that any two high level sets
of $f$ are diffeomorphic via the gradient curves of $f$. Since by
(6.2.25),
\begin{align*}
\frac{d}{ds}\bar{R}(\eta(s),-1)&  = \<\nabla \bar{R},\dot{\eta}(s)\>\\
&  = 2\bar{R}_{ij}\nabla_if\nabla_jf \\
&  \geq 0
\end{align*}
for any integral curve $\eta(s)$ of $\nabla f$, we conclude that
the scalar curvature $\bar{R}(x,-1)$ has a positive lower bound on
$\bar{M}$, which contradicts the Bonnet-Myers Theorem. So we have
proved that the limiting gradient shrinking soliton is compact.

By Proposition 5.1.10, the compact limiting gradient shrinking
soliton has constant curvature. This says that the scalings of the
ancient $\kappa$-solution $(M,g_{ij}(x,t))$ along a sequence of
points $q_k \in M$ and a sequence of times $t_k \rightarrow
-\infty$ converge in the $C^{\infty}$ topology to the round
$\mathbb{S}^2$ or the round $\mathbb{R}\mathbb{P}^2$. In
particular, by looking at the time derivative of the volume and
the Gauss-Bonnet theorem, we know that the ancient
$\kappa$-solution $(M,g_{ij}(x,t))$ exists on a maximal time
interval $(-\infty,T)$ with $T < +\infty$.

Consider the scaled entropy of Hamilton \cite{Ha88}
$$
E(t)=\int_MR\log[R(T-t)]dV_t.
$$
We compute
\begin{align}
\frac{d}{dt}E(t)&  = \int_M \log[R(T-t)]\Delta
RdV_t+\int_M\left[\Delta R+R^2-\frac{R}{(T-t)}\right]dV_t\\
&  = \int_M\left[-\frac{|\nabla R|^2}{R}+R^2-rR\right]dV_t \nn\\
&  = \int_M\left[-\frac{|\nabla R|^2}{R}+(R-r)^2\right]dV_t \nn
\end{align}  
where $r=\int_MRdV_t/Vol_t(M)$ and we have used
$\Vol_t(M)=(\int_MRdV_t)\cdot(T-t)$ (by the Gauss-Bonnet theorem).

We now need an inequality of Chow in \cite{ChowE}. For sake of
completeness, we present his proof as follows. For a smooth function
$f$ on the surface $M$, one can readily check
\begin{align*}
\int_M(\Delta f)^2 &=2\int_M\left|\nabla_i\nabla_j
f-\frac{1}{2}(\Delta f)g_{ij}\right|^2
+\int_MR|\nabla f|^2, \\
\int_M\frac{|\nabla R+R\nabla f|^2}{R} &=\int_M\frac{|\nabla
R|^2}{R}-2\int_M R(\Delta f)+\int_MR|\nabla f|^2,
\end{align*}
and then
\begin{align*}
&\int_M\frac{|\nabla R|^2}{R}+\int_M(\Delta f)(\Delta f-2R)\\
&=2\int_M\left|\nabla_i\nabla_jf-\frac{1}{2}(\Delta
f)g_{ij}\right|^2 +\int_M\frac{|\nabla R+R\nabla f|^2}{R}.
\end{align*}
By choosing $\Delta f=R-r$, we get
\begin{align*}
&\int_M\frac{|\nabla R|^2}{R}-\int_M(R-r)^2 \\
&=2\int_M\left|\nabla_i\nabla_jf-\frac{1}{2}(\Delta
f)g_{ij}\right|^2 +\int_M\frac{|\nabla R+R\nabla f|^2}{R}\ge 0.
\end{align*}
If the equality holds, then we have
$$
\nabla_i\nabla_jf-\frac{1}{2}(\Delta f)g_{ij}=0
$$
i.e., $\nabla f$ is conformal. By the Kazdan-Warner identity
\cite{KW}, it follows that
$$
\int_M\nabla R\cdot \nabla f=0,
$$
so
\begin{align*}
0&  = -\int_MR\Delta f\\
&  = -\int_M(R-r)^2.
\end{align*}
Hence we have proved the following inequality due to Chow
\cite{ChowE} \be \int_M\frac{|\nabla R|^2}{R}\ge\int_M(R-r)^2,
\ee and the equality holds if and only if $R\equiv r$.

The combination (6.2.27) and (6.2.28) shows that the scaled
entropy $E(t)$ is strictly decreasing along the Ricci flow unless
we are on the round sphere $\mathbb{S}^2$ or its quotient
$\mathbb{R}\mathbb{P}^2$. Moreover the convergence result in
Theorem 5.1.11 shows that the scaled entropy $E$ has its minimum
value at the constant curvature metric (round $\mathbb{S}^2$ or
round $\mathbb{R}\mathbb{P}^2$). We had shown that the scalings of
the nonflat ancient $\kappa$-solution along a sequence of times
$t_k\rightarrow-\infty$ converge to the constant curvature metric.
Then $E(t)$ has its minimal value at $t=-\infty$, so it was
constant all along, hence the ancient $\kappa$-solution must have
constant curvature for each $t\in (-\infty,T)$. This proves the
theorem.
\end{pf}

\section{Curvature Estimates via Volume Growth}

For solutions to the Ricci flow, Perelman's no local collapsing
theorems tell us that the local curvature upper bounds imply the
local volume lower bounds. Conversely, one would expect to get
local curvature upper bounds from local volume lower bounds. If
this is the case, one will be able to establish an elliptic type
estimate for the curvatures of solutions to the Ricci flow. This
will provide the key estimate for the canonical neighborhood
structure and thick-thin decomposition of the Ricci flow on
three-manifolds. In this section we derive such curvature
estimates for nonnegatively curved solutions. In the next chapter
we will derive similar estimates for all smooth solutions, as well
as surgically modified solutions, of the Ricci flow on
three-manifolds.

Let $M$ be an $n$-dimensional complete noncompact Riemannian
manifold with nonnegative Ricci curvature. Pick an origin $O\in
M$. The well-known Bishop-Gromov volume comparison theorem tells
us the ratio $Vol(B(O,r))/r^n$ is monotone nonincreasing in
$r\in[0,+\infty)$. Thus there exists a limit
$$
\nu_M=\lim_{r\rightarrow+\infty}\frac{\Vol(B(O,r))}{r^n}.
$$
Clearly the number $\nu_M$ is invariant under dilation and is
independent of the choice of the origin. $\nu_M$ is called the
{\bf asymptotic volume ratio}\index{asymptotic volume ratio} of
the Riemannian manifold $M$.

The following result obtained by Perelman (cf. Proposition 11.4 of
\cite{P1}) shows that any ancient $\kappa$-solution must have zero
asymptotic volume ratio. This result for the Ricci flow on
K\"ahler manifolds was implicitly and independently given by Chen
and the second author in the proof of Theorem 5.1 of \cite{CZ04}.
Moreover in the K\"ahler case, as shown by Chen, Tang and the
second author (implicitly in the proof of Theorem 4.1 of
\cite{CTZ04} for complex two dimension) and by Ni (in \cite{Ni}
for all dimensions), the condition of nonnegative curvature
operator can be replaced by the weaker condition of nonnegative
bisectional curvature.

\begin{lemma} [{Perelman \cite{P1}}]
Let $M$ be an $n$-dimensional complete noncompact Riemannian
manifold. Suppose $g_{ij}(x,t)$, $x\in M$ and $t\in(-\infty,T)$
with $T>0$, is a nonflat ancient solution of the Ricci flow with
nonnegative curvature operator and bounded curvature. Then the
asymptotic volume ratio of the solution metric satisfies
$$
\nu_M(t)=\lim_{r\rightarrow+\infty}\frac{\Vol_t(B_t(O,r))}{r^n}=0
$$
for each $t\in(-\infty,T)$.
\end{lemma}

\begin{pf}
The proof is by induction on the dimension. When the dimension is
two, the lemma is valid by Theorem 6.2.2. For dimension $\geq 3$,
we argue by contradiction.

Suppose the lemma is valid for dimensions $\leq n-1$ and suppose
$\nu_{M}(t_0) > 0$ for some $n$-dimensional nonflat ancient
solution with nonnegative curvature operator and bounded curvature
at some time $t_0 \leq 0$. Fixing a point $x_0 \in M$, we consider
the asymptotic scalar curvature ratio
$$
A=\limsup_{d_{t_0}(x,x_0)\rightarrow+\infty}
R(x,t_0)d^2_{t_0}(x,x_0).
$$
We divide the proof into three cases.

\medskip
{\it Case} 1: $A=+\infty$.

\smallskip
By Lemma 6.1.3, there exist sequences of points $x_k \in M$
divergent to infinity, of radii $r_k \rightarrow +\infty$, and of
positive constants $\delta_k \rightarrow 0$ such that
\begin{itemize}
\item[(i)] $R(x,t_0)\leq(1+\delta_k)R(x_k,t_0)$ for all $x$ in the
ball $B_{t_0}(x_k,r_k)$ of radius $r_k$ around $x_k$, \item[(ii)]
$r^2_kR(x_k,t_0)\rightarrow+\infty$ as $k\rightarrow+\infty$,
\item[(iii)] $d_{t_0}(x_k,x_0)/r_k\rightarrow+\infty$.
\end{itemize}

By scaling the solution around the points $x_k$ with factor
$R(x_k,t_0)$, and shifting the time $t_0$ to the new time zero, we
get a sequence of rescaled solutions
$$
{g}_k(s) = R(x_k,t_0)g\(\cdot,t_0 + \frac{s}{R(x_k,t_0)}\)
$$
to the Ricci flow. Since the ancient solution has nonnegative
curvature operator and bounded curvature, there holds the
Li-Yau-Hamilton inequality (Corollary 2.5.5). Thus the rescaled
solutions satisfy
$$
{R}_k(x,s) \leq (1+\delta_k)
$$
for all $s\leq 0$ and $x\in
B_{{g}_k(0)}(x_k,r_k\sqrt{R(x_k,t_0)}).$ Since $\nu_{M}(t_0)
> 0$, it follows from the standard volume comparison and Theorem
4.2.2 that the injectivity radii of the rescaled solutions ${g}_k$
at the points $x_k$ and the new time zero is uniformly bounded
below by a positive number. Then by Hamilton's compactness theorem
(Theorem 4.1.5), after passing to a subsequence,
$(M,{g}_k(s),x_k)$ will converge to a solution
$(\tilde{M},\tilde{g}(s),O)$ to the Ricci flow with
$$
\tilde{R}(y,s) \leq 1, \mbox{ for all } s \leq 0 \mbox{ and } y\in
\tilde{M},
$$
and
$$
\tilde{R}(O,0) = 1.
$$
Since the metric is shrinking, by (ii) and (iii), we get
$$
R(x_k,t_0)d^2_{g(\cdot,t_0 + \frac{s}{R(x_k,t_0)})}(x_0,x_k) \geq
R(x_k,t_0)d^2_{g(\cdot,t_0)}(x_0,x_k)
$$
which tends to $+\infty$, as $k\rightarrow +\infty$, for all
$s\leq 0$. Thus by Proposition 6.1.2, for each $s\leq 0$,
$(\tilde{M},\tilde{g}(s))$ splits off a line. We now consider the
lifting of the solution $(\tilde{M},\tilde{g}(s)), s \leq 0,$ to
its universal cover and denote it by
$(\tilde{\tilde{M}},\tilde{\tilde{g}}(s)), s \leq 0.$ Clearly we
still have
$$
\nu_{\tilde{M}}(0) > 0 \mbox{ and }\nu_{\tilde{\tilde{M}}}(0) > 0.
$$
By applying Hamilton's strong maximum principle and the de Rham
decomposition theorem, the universal cover $\tilde{\tilde{M}}$
splits isometrically as $X \times \mathbb{R}$ for some
$(n-1)$-dimensional nonflat (complete) ancient solution $X$ with
nonnegative curvature operator and bounded curvature. These imply
that $\nu_X(0) > 0$, which contradicts the induction hypothesis.

\medskip
{\it Case} 2: $0 < A < +\infty$.

\smallskip
Take a sequence of points $x_k$ divergent to infinity such that
$$
R(x_k,t_0)d^2_{t_0}(x_k,x_0) \rightarrow A,  \mbox{ as } k
\rightarrow +\infty.
$$
Consider the rescaled solutions $(M,g_k(s))$ (around the fixed
point $x_0$), where
$$
g_k(s) = R(x_k,t_0)g\(\cdot,t_0 + \frac{s}{R(x_k,t_0)}\), s\in
(-\infty,0].
$$
Then there is a constant $C > 0$ such that \be
   \left\{
   \begin{array}{lll}
    R_k(x,0) \leq {C}/{d^2_k(x,x_0,0)},
         \\[2mm]
    R_k(x_k,0) = 1,
         \\[2mm]
    d_k(x_k,x_0,0) \rightarrow \sqrt{A} >0,
\end{array}
\right.
\ee where $d_k(\cdot,x_0,0)$ is the distance function from the
point $x_0$ in the metric $g_k(0)$.

Since $\nu_M(t_0) >0$, it is a basic result in Alexandrov space
theory (see for example Theorem 7.6 of \cite{CC96}) that a
subsequence of $(M,g_k(s),x_0)$ converges in the Gromov-Hausdorff
sense to an $n$-dimensional metric cone
$(\tilde{M},\tilde{g}(0),x_0)$ with vertex $x_0$.

By (6.3.1), the standard volume comparison and Theorem 4.2.2, we
know that the injectivity radius of $(M,g_k(0))$ at $x_k$ is
uniformly bounded from below by a positive number $\rho_0$. After
taking a subsequence, we may assume $x_k$ converges to a point
$x_{\infty}$ in $\tilde{M}\setminus \{x_0\}$. Then by Hamilton's
compactness theorem (Theorem 4.1.5), we can take a subsequence
such that the metrics $g_k(s)$ on the metric balls
$B_0(x_k,\frac{1}{2}\rho_0) (\subset M $ with respect to the
metric $g_k(0))$ converge in the $C^{\infty}_{loc}$ topology to a
solution of the Ricci flow on a ball
$B_0(x_{\infty},\frac{1}{2}\rho_0)$. Clearly the
$C^{\infty}_{loc}$ limit has nonnegative curvature operator and it
is a piece of the metric cone at the time $s=0$. By (6.3.1), we
have \be
\tilde{R}(x_{\infty},0) = 1. 
\ee

Let $x$ be any point in the limiting ball
$B_0(x_{\infty},\frac{1}{2}\rho_0)$ and $e_1$ be any radial
direction at $x$. Clearly $\tilde{Ric}(e_1,e_1) = 0$. Recall that
the evolution equation of the Ricci tensor in frame coordinates is
$$
\frac{\partial}{\partial t}{\tilde{R}_{ab}} =
\tilde{\triangle}{\tilde{R}_{ab}} +
2{\tilde{R}_{acbd}}{\tilde{R}_{cd}}.
$$
Since the curvature operator is nonnegative, by applying
Hamilton's strong maximum principle (Theorem 2.2.1) to the above
equation, we deduce that the null space of $\tilde{Ric}$ is
invariant under parallel translation. In particular, all radial
directions split off locally and isometrically. While by (6.3.2),
the piece of the metric cone is nonflat. This gives a
contradiction.

\medskip
{\it Case} 3: $A = 0$.

\smallskip
The gap theorem as was initiated by Mok-Siu-Yau \cite{MoSY} and
established by Greene-Wu \cite{GW82, GW},
Eschenberg-Shrader-Strake \cite{ESS}, and Drees \cite{Dr} shows
that a complete noncompact $n$-dimensional (except $n=4$ or $8$)
Riemannian manifold with nonnegative sectional curvature and the
asymptotic scalar curvature ratio $A=0$ must be flat. So the
present case is ruled out except in dimension $n=4$ or $8$. Since
in our situation the asymptotic volume ratio is positive and the
manifold is the solution of the Ricci flow, we can give an
alternative proof for all dimensions as follows.

We claim the sectional curvature of $(M,g_{ij}(x,t_0))$ is
positive everywhere. Indeed, by Theorem 2.2.2, the image of the
curvature operator is just the restricted holonomy algebra
$\mathcal{G}$ of the manifold. If the sectional curvature vanishes
for some two-plane, then the holonomy algebra $\mathcal{G}$ cannot
be $so(n)$. We observe the manifold is not Einstein since it is
noncompact, nonflat and has nonnegative curvature operator. If
$\mathcal{G}$ is irreducible, then by Berger's Theorem \cite{Ber},
$\mathcal{G} = u(\frac{n}{2})$. Thus the manifold is  K\"ahler
with bounded and nonnegative bisectional curvature and with
curvature decay faster than quadratic. Then by the gap theorem
obtained by Chen and the second author in \cite{CZ03}, this
K\"ahler manifold must be flat. This contradicts the assumption.
Hence the holonomy algebra $\mathcal{G}$ is reducible and the
universal cover of $M$ splits isometrically as $\tilde{M}_1 \times
\tilde{M}_2$ nontrivially. Clearly the universal cover of $M$ has
positive asymptotic volume ratio. So $\tilde{M}_1$ and
$\tilde{M}_2$ still have positive asymptotic volume ratio and at
least one of them is nonflat. By the induction hypothesis, this is
also impossible. Thus our claim is proved.

Now we know that the sectional curvature of $(M,g_{ij}(x,t_0))$ is
positive everywhere. Choose a sequence of points $x_k$ divergent
to infinity such that
$$
   \left\{
   \begin{array}{lll}
    R(x_k,t_0)d^2_{t_0}(x_k,x_0) = \sup \{ R(x,t_0)d^2_{t_0}(x,x_0)\
    |\
    d_{t_0}(x,x_0)\geq d_{t_0}(x_k,x_0)\},
         \\[2mm]
    d_{t_0}(x_k,x_0) \geq k,
         \\[2mm]
    R(x_k,t_0)d^2_{t_0}(x_k,x_0)=\varepsilon_k \rightarrow 0.
\end{array}
\right.
$$
Consider the rescaled metric
$$
g_k(0) = R(x_k,t_0)g(\cdot,t_0)
$$
as before. Then \be
   \left\{
   \begin{array}{lll}
    R_k(x,0) \leq {\varepsilon_k}/{d^2_k(x,x_0,0)},\quad \mbox{ for }
    d_k(x,x_0,0) \geq \sqrt{\varepsilon_k},
         \\[2mm]
    d_{k}(x_k,x_0,0) = \sqrt{\varepsilon_k} \rightarrow 0.
\end{array}
\right. 
\ee As in Case 2, the rescaled marked solutions $(M,g_k(0),x_0)$
will converge in the Gromov-Hausdorff sense to a metric cone
$(\tilde{M},\tilde{g}(0),x_0)$. And by the virtue of Hamilton's
compactness theorem (Theorem 4.1.5), up to a subsequence, the
convergence is in the $C^{\infty}_{loc}$ topology in
$\tilde{M}\setminus \{x_0\}$. We next claim the metric cone
$(\tilde{M},\tilde{g}(0),x_0)$ is isometric to $\mathbb{R}^n$.

Indeed, let us write the metric cone $\tilde{M}$ as a warped
product $\mathbb{R}_{+}\times_r X^{n-1}$ for some
$(n-1)$-dimensional manifold $X^{n-1}$. By (6.3.3), the metric
cone must be flat and $X^{n-1}$ is isometric to a quotient of the
round sphere $\mathbb{S}^{n-1}$ by fixed point free isometries in
the standard metric. To show $\tilde{M}$ is isometric to
$\mathbb{R}^n$, we only need to verify that $X^{n-1}$ is simply
connected.

Let $\varphi$ be the Busemann function of $(M,g_{ij}(\cdot,t_0))$
with respect to the point $x_0$. Since $(M,g_{ij}(\cdot,t_0))$ has
nonnegative sectional curvature, it is easy to see that for any
small $\varepsilon >0$, there is a $r_0 >0$ such that
$$(1-\varepsilon)d_{t_0}(x,x_0) \leq \varphi(x) \leq
d_{t_0}(x,x_0)$$ for all $x \in M\setminus B_{t_0}(x_0,r_0)$. The
strict positivity of the sectional curvature of the manifold
$(M,g_{ij}(\cdot,t_0))$ implies that the square of the Busemann
function is strictly convex (and exhausting). Thus every level set
$\varphi^{-1}(a)$, with $a > \inf \{\varphi(x)\ |\ x \in M \}$, of
the Busemann function $\varphi$ is diffeomorphic to the
$(n-1)$-sphere $\mathbb{S}^{n-1}$. In particular,
$\varphi^{-1}([a,\frac{3}{2}a])$ is simply connected for $a > \inf
\{\varphi(x)\ |\ x \in M \}$ since $n \geq 3$.

Consider an annulus portion $[1,2]\times X^{n-1}$ of the metric
cone $\tilde{M} = \mathbb{R}_+\times_r X^{n-1}$. It is the limit
of $(M_k,g_k(0))$, where
$$
M_k = \left\{ x \in M \ \Big|\  \frac{1}{\sqrt{R(x_k,t_0)}} \leq
d_{t_0}(x,x_0) \leq \frac{2}{\sqrt{R(x_k,t_0)}} \right\}.
$$
It is clear that
\begin{multline*}
\varphi^{-1}\(\left[\frac{1}{\sqrt{R(x_k,t_0)}},
\frac{2(1-\varepsilon)}{\sqrt{R(x_k,t_0)}}\right]\)
\subset M_k \\
\subset
\varphi^{-1}\(\left[\frac{1-\varepsilon}{\sqrt{R(x_k,t_0)}},
\frac{2}{\sqrt{R(x_k,t_0)}}\right]\)\!\!\!
\end{multline*}
for $k$ large enough. Thus any closed loop in $\{\frac{3}{2}\}
\times X^{n-1}$ can be shrunk to a point by a homotopy in
$[1,2]\times X^{n-1}$. This shows that $X^{n-1}$ is simply
connected. Hence we have proven that the metric cone
$(\tilde{M},\tilde{g}(0),x_0)$ is isometric to $\mathbb{R}^n$.
Consequently,
$$
\lim_{k\rightarrow+\infty}\!\frac{\Vol_{g(t_0)}\!\(\!\!B_{g(t_0)}
\(\!x_0,\frac{r}{\sqrt{R(x_k,t_0)}}\)\setminus
B_{g(t_0)}\!\(\!x_0,\frac{\sigma r}{\sqrt{R(x_k,t_0)}}\)\)}{
\(\frac{r}{\sqrt{R(x_k,t_0)}}\)^n} = \alpha_n (1\!-\!\sigma^n)
$$
for any $r>0$ and $0<\sigma<1$, where $\alpha_n$ is the volume of
the unit ball in the Euclidean space $\mathbb{R}^n$. Finally, by
combining with the monotonicity of the Bishop-Gromov volume
comparison, we conclude that the manifold $(M,g_{ij}(\cdot,t_0))$
is flat and isometric to $\mathbb{R}^n$. This contradicts the
assumption.

Therefore, we have proved the lemma.
\end{pf}

Finally we would like to include an alternative simpler argument,
inspired by Ni \cite{Ni}, for the above Case 2 and Case 3 to avoid
the use of the gap theorem, holonomy groups, and asymptotic cone
structure.

\medskip
{\bf \em Alternative Proof for Case {\bf 2} and Case {\bf 3}.}

Let us consider the situation of $0 \leq A < +\infty$ in the above
proof. Observe that $\nu_M(t)$ is nonincreasing in time $t$ by
using Lemma 3.4.1(ii) and the fact that the metric is shrinking in
$t$. Suppose $\nu_M(t_0)>0$, then the solution $g_{ij}(\cdot,t)$
is $\kappa$-noncollapsed for $t\leq t_0$ for some uniform
$\kappa>0$. By combining with Theorem 6.2.1, there exist a
sequence of points $q_k$ and a sequence of times $t_k \rightarrow
-\infty$ such that the scalings of $g_{ij}(\cdot,t)$ around $q_k$
with factor $|t_k|^{-1}$ and with the times $t_k$ shifting to the
new time zero converge to a nonflat gradient shrinking soliton
$\bar{M}$ in the $C_{loc}^\infty$ topology. This gradient soliton
also has maximal volume growth (i.e. $\nu_{\bar{M}}(t)>0$) and
satisfies the Li-Yau-Hamilton estimate (Corollary 2.5.5). If the
curvature of the shrinking soliton $\bar{M}$ at the time $-1$ is
bounded, then we see from the proof of Theorem 6.2.2 that by using
the equations (6.2.24)-(6.2.26), the scalar curvature has a
positive lower bound everywhere on $\bar{M}$ at the time $-1$. In
particular, this implies the asymptotic scalar curvature ratio
$A=\infty$ for the soliton at the time $-1$, which reduces to Case
1 and arrives at a contradiction by the dimension reduction
argument. On the other hand, if the scalar curvature is unbounded,
then by Lemma 6.1.4, the Li-Yau-Hamilton estimate (Corollary
2.5.5) and Lemma 6.1.2, we can do the same dimension reduction as
in Case 1 to arrive at a contradiction also.
\endproof

The following lemma is a local and space-time version of Lemma
6.1.4 on picking local (almost) maximum curvature points. We
formulate it from Perelman's arguments in section 10 of \cite{P1}.

\begin{lemma}
For any positive constants $B, C$ with $B>4$ and $C>1000$, there
exists $1\leq A<\min\{\frac{1}{4}B,\frac{1}{1000}C\}$ which tends
to infinity as $B$ and $C$ tend to infinity and satisfies the
following property. Suppose we have a (not necessarily complete)
solution $g_{ij}(t)$ to the Ricci flow, defined on $M\times
[-t_0,0]$, so that at each time $t\in [-t_0,0]$ the metric ball
$B_t(x_0,1)$ is compactly contained in $M$. Suppose there exists a
point $(x',t') \in M\times (-t_0,0]$ such that
$$
d_{t'}(x',x_0)\leq \frac{1}{4} \; \mbox{ and }\;
|Rm(x',t')|>C+B(t'+t_0)^{-1}.
$$
Then we can find a point $(\bar{x},\bar{t}) \in M\times (-t_0,0]$
such that
$$
d_{\bar{t}}(\bar{x},x_0)<\frac{1}{3}\; \mbox{ with }\; Q
=|Rm(\bar{x},\bar{t})|>C+ B(\bar{t}+t_0)^{-1},
$$
and
$$
|Rm(x,t)|\leq 4Q
$$
for all $(-t_0<) \ \bar{t}-AQ^{-1}\leq t\leq \bar{t}$ and
$d_t(x,\bar{x})\leq \frac{1}{10}A^{\frac{1}{2}}Q^{-\frac{1}{2}}.$
\end{lemma}

\begin{pf}
We first claim that there exists a point $(\bar{x},\bar{t})$ with
$-t_0<\bar{t}\leq 0$ and $d_{\bar{t}}(\bar{x},x_0)<\frac{1}{3}$
such that
$$
Q=|Rm(\bar{x},\bar{t})|>C+ B(\bar{t}+t_0)^{-1},
$$
and \be
|Rm(x,t)|\leq 4Q 
\ee wherever $\bar{t}-AQ^{-1}\leq t\leq \bar{t},\ \ d_t(x,x_0)\leq
d_{\bar{t}}(\bar{x},x_0) +(AQ^{-1})^{\frac{1}{2}}.$

We will construct such $(\bar{x},\bar{t})$ as a limit of a finite
sequence of points.  Take an arbitrary $(x_1,t_1)$ such that
$$
d_{t_1}(x_1,x_0)\leq \frac{1}{4},\ \ \ -t_0<t_1\leq 0 \; \mbox{
and } \;|Rm(x_1,t_1)|>C+B(t_1+t_0)^{-1}.
$$
Such a point clearly exists by our assumption. Assume we have
already constructed $(x_k,t_k).$ If we cannot take the point
$(x_k,t_k)$ to be the desired point $(\bar{x},\bar{t})$, then
there exists a point $(x_{k+1},t_{k+1})$ such that
$$
t_k-A|Rm(x_k,t_k)|^{-1}\leq t_{k+1}\leq t_k,
$$
and
$$
d_{t_{k+1}}(x_{k+1},x_0) \leq
d_{t_{k}}(x_{k},x_0)+(A|Rm(x_k,t_k)|^{-1})^{\frac{1}{2}},
$$
but
$$
|Rm(x_{k+1},t_{k+1})|>4|Rm(x_k,t_k)|.$$

\vskip 0.1cm\noindent It then follows that
\begin{align*}
d_{t_{k+1}}(x_{k+1},x_0)& \leq
d_{t_{1}}(x_{1},x_0)+A^{\frac{1}{2}}
\(\sum\limits_{i=1}^k|Rm(x_i,t_i)|^{-\frac{1}{2}}\)\\
&  \leq \frac{1}{4}+A^{\frac{1}{2}}
\(\sum\limits_{i=1}^k2^{-(i-1)}|Rm(x_1,t_1)|^{-\frac{1}{2}}\)\\
&  \leq \frac{1}{4}+2(AC^{-1})^{\frac{1}{2}} \\
&  < \frac{1}{3},
\end{align*}
\begin{displaymath}
\begin{split}
t_{k+1}-(-t_0)&  =\sum_{i=1}^k(t_{i+1}-t_i)+(t_1-(-t_0))\\[1mm]
&  \geq-\sum_{i=1}^kA|Rm(x_i,t_i)|^{-1}+(t_1-(-t_0))\\[1mm]
&  \geq-A\sum_{i=1}^k4^{-(i-1)}|Rm(x_1,t_1)|^{-1}+(t_1-(-t_0))\\[1mm]
&  \geq-\frac{2A}{B}(t_1+t_0)+(t_1+t_0)\\[1mm]
&  \geq\frac{1}{2}(t_1+t_0),
\end{split}
\end{displaymath}
and
\begin{align*}
|Rm(x_{k+1},t_{k+1})|&  > 4^k|Rm(x_1,t_1)| \\
&  \geq 4^kC\rightarrow +\infty\; \mbox{ as } \;
k\rightarrow+\infty.
\end{align*}
Since the solution is smooth, the sequence $\{(x_k, t_k)\}$ is
finite and its last element fits. Thus we have proved  assertion
(6.3.4).

{}From the above construction we also see that the chosen point
$(\bar{x},\bar{t})$ satisfies
$$
d_{\bar{t}}(\bar{x},x_0) < \frac{1}{3}
$$
and
$$
Q=|Rm(\bar{x},\bar{t})|>C+B(\bar{t}+t_0)^{-1}.
$$
Clearly, up to some adjustment of the constant $A$, we only need
to show that
$$
|Rm(x,t)|\leq 4Q\leqno(6.3.4)' 
$$
wherever $\bar{t}-\frac{1}{200n}A^{\frac{1}{2}}Q^{-1}\leq t\leq
\bar{t} \ \ \mbox{and} \ \ d_t(x,\bar{x})\leq
\frac{1}{10}A^{\frac{1}{2}}Q^{-\frac{1}{2}}$.

For any point $(x,\bar{t})$ with $d_{\bar{t}}(x,\bar{x})\leq
\frac{1}{10}A^{\frac{1}{2}}Q^{-\frac{1}{2}}$, we have
\begin{align*}
d_{\bar{t}}(x,x_0)&  \leq
d_{\bar{t}}(\bar{x},x_0)+d_{\bar{t}}(x,\bar{x})\\
&  \leq d_{\bar{t}}(\bar{x},x_0)+(AQ^{-1})^{\frac{1}{2}}
\end{align*}
and then by (6.3.4)
$$
|Rm(x,\bar{t})|\leq 4Q.
$$
Thus by continuity, there is a minimal
$\bar{t}'\in[\bar{t}-\frac{1}{200n}A^{\frac{1}{2}}Q^{-1},\bar{t}]$
such that \be \sup\left\{|Rm(x,t)|\  | \  \bar{t}'\leq t\leq
\bar{t},\quad d_t(x,\bar{x})\leq
\frac{1}{10}A^{\frac{1}{2}}Q^{-\frac{1}{2}}\right\}\leq 5Q.
\ee

For any point $(x,t)$ with $\bar{t}'\leq t\leq \bar{t}$ and
$d_t(x,\bar{x})\leq \frac{1}{10}(AQ^{-1})^{\frac{1}{2}}$, we
divide the discussion into two cases.

\medskip
{\it Case} (1): $d_t(\bar{x},x_0)\leq
\frac{3}{10}(AQ^{-1})^{\frac{1}{2}}$.

\smallskip
{}From assertion (6.3.4) we see that
$$
\sup\{|Rm(x,t)|\ \ | \ \ \bar{t}'\leq t\leq \bar{t},\ \
d_t(x,x_0)\leq (AQ^{-1})^{\frac{1}{2}}\}\leq 4Q. \leqno (6.3.5)'
$$
Since $d_t(\bar{x},x_0)\leq \frac{3}{10}(AQ^{-1})^{\frac{1}{2}}$,
we have
\begin{align*}
d_t({x},x_0) &\leq d_t(x,\bar{x})+d_t(\bar{x},x_0)\\
&  \leq \frac{1}{10}(AQ^{-1})^{\frac{1}{2}}
+\frac{3}{10}(AQ^{-1})^{\frac{1}{2}}\\
&  \leq (AQ^{-1})^{\frac{1}{2}}
\end{align*}
which implies the estimate $|Rm(x,t)|\leq 4Q $ from (6.3.5)$'$.

\medskip
{\it Case} (2): $d_t(\bar{x},x_0)>
\frac{3}{10}(AQ^{-1})^{\frac{1}{2}}$.

\smallskip
{}From the curvature bounds in (6.3.5) and (6.3.5)$'$, we can
apply Lemma 3.4.1 (ii) with $r_0=\frac{1}{10}Q^{-\frac{1}{2}}$ to
get
$$
\frac{d}{dt}(d_t(\bar{x},x_0))\geq -40(n-1)Q^{\frac{1}{2}}
$$
and then
\begin{align*}
d_t(\bar{x},x_0)&  \leq d_{\hat{t}}(\bar{x},x_0)
+40n(Q^{\frac{1}{2}})\(\frac{1}{200n}A^{\frac{1}{2}}Q^{-1}\)\\
& \leq d_{\hat{t}}(\bar{x},x_0)+\frac{1}{5}(AQ^{-1})^{\frac{1}{2}}
\end{align*}
where $\hat{t} \in (t,\bar{t}]$ satisfies the property that
$d_s(\bar{x},x_0)\geq \frac{3}{10}(AQ^{-1})^{\frac{1}{2}}$
whenever $s \in [t,\hat{t}]$. So we have either
\begin{align*}
d_t({x},x_0)&  \leq d_t(x,\bar{x})+d_t(\bar{x},x_0)\\
&  \leq \frac{1}{10}(AQ^{-1})^{\frac{1}{2}}
+\frac{3}{10}(AQ^{-1})^{\frac{1}{2}}+\frac{1}{5}(AQ^{-1})^{\frac{1}{2}}\\
&  \leq (AQ^{-1})^{\frac{1}{2}},
\end{align*}
or
\begin{align*}
d_t(x,x_0)&  \leq d_t(x,\bar{x})+d_t(\bar{x},x_0)\\
&  \leq \frac{1}{10}(AQ^{-1})^{\frac{1}{2}}
+d_{\bar{t}}(\bar{x},x_0)+\frac{1}{5}(AQ^{-1})^{\frac{1}{2}} \\
&  \leq d_{\bar{t}}(\bar{x},x_0)+(AQ^{-1})^{\frac{1}{2}}.
\end{align*}
It then follows from (6.3.4) that $|Rm(x,t)|\leq 4Q $.

Hence we have proved
$$
|Rm(x,t)|\leq 4Q
$$
for any point $(x,t)$ with $\bar{t}'\leq t\leq \bar{t}$ and
$d_t(x,\bar{x})\leq \frac{1}{10}(AQ^{-1})^{\frac{1}{2}}$. By
combining with the choice of $\bar{t}'$ in (6.3.5), we must have
$\bar{t}'= \bar{t}-\frac{1}{200n}A^{\frac{1}{2}}Q^{-1}$. This
proves assertion (6.3.4)$'$.

Therefore we have completed the proof of the lemma.
\end{pf}

We now use the volume lower bound assumption to establish the
crucial curvature upper bound estimate of Perelman \cite{P1} for
the Ricci flow. For the Ricci flow on K\"ahler manifolds, a global
version of this estimate (i.e., curvature decaying linear in time
and quadratic in space) was independently obtained in \cite{CTZ04}
and \cite{CZ04}. Note that the volume estimate conclusion in the
following Theorem 6.3.3 (ii) was not stated in Corollary 11.6 (b)
of Perelman \cite{P1}. The estimate will be used later in the
proof of Theorem 7.2.2 and Theorem 7.5.2.

\begin{theorem}[{Perelman \cite{P1}}]
For every $w>0$ there exist $B=B(w)<+\infty,\ \ C=C(w)<+\infty, \
\ \tau_0=\tau_0(w)>0,$ and $\xi=\xi(w)>0$ $($depending also on the
dimension$)$ with the following properties. Suppose we have a
$($not necessarily complete$)$ solution $g_{ij}(t)$ to the Ricci
flow, defined on $M\times[-t_0r_0^2,0],$ so that at each time
$t\in [-t_0r_0^2,0]$ the metric ball $B_t(x_0,r_0)$ is compactly
contained in $M.$
\begin{itemize}
\item[(i)] If at each time $t \in [-t_0r_0^2,0]$,
$$
Rm(.,t)\geq -r_0^{-2} \; \mbox{ on }\; B_t(x_0,r_0)
$$
$$
\mbox{and } \; \Vol_t(B_t(x_0,r_0))\geq wr_0^n,
$$
then we have the estimate
$$
|Rm(x,t)|\leq Cr_0^{-2}+B(t+t_0r_0^2)^{-1}
$$
whenever $-t_0r_0^2 < t\leq 0$ and $d_t(x,x_0)\leq
\frac{1}{4}r_0.$ \item[(ii)] If for some $0 < \bar{\tau} \leq
t_0$,
$$
Rm(x,t)\geq -r_0^{-2} \; \mbox{ for }\; t\in [-\bar{\tau}
r_0^2,0], x\in B_t(x_0,r_0),
$$
$$
\mbox{and } \; \Vol_0(B_0(x_0,r_0))\geq wr_0^n,
$$
then we have the estimates
$$
\Vol_t(B_t(x_0,r_0))\geq \xi r_0^n\; \mbox{ for all }\;
\max\{-\bar{\tau}r_0^2,-\tau_0r_0^2\}   \leq t \leq 0,
$$
and
$$
|Rm(x,t)|\leq Cr_0^{-2}+B(t -\max\{-\bar{\tau}r_0^2,-\tau_0r_0^2\}
)^{-1}
$$
whenever $\max\{-\bar{\tau}r_0^2,-\tau_0r_0^2\} < t\leq 0$ and
$d_t(x,x_0)\leq \frac{1}{4}r_0.$
\end{itemize}
\end{theorem}

\begin{pf}
By scaling we may assume $r_0=1.$

\medskip
(i) By the standard (relative) volume comparison, we know that
there exists some $w'>0$, with $w'\leq w $, depending only on $w$,
such that for each point $(x,t)$ with $-t_0\leq t\leq 0$ and
$d_t(x,x_0)\leq \frac{1}{3},$ and for each $r\leq \frac{1}{3}$,
there holds \be
\Vol_t(B_t(x,r))\geq w'r^n. 
\ee

We argue by contradiction. Suppose there are sequences $B,
C\rightarrow +\infty,$ of solutions $g_{ij}(t)$ and points
$(x',t')$ such that
$$
d_{t'}(x',x_0)\leq \frac{1}{4},\quad -t_0<t'\leq 0\; \mbox{ and }
\; |Rm(x',t')|>C+B(t'+t_0)^{-1}.
$$
Then by Lemma 6.3.2, we can find a sequence of points
$(\bar{x},\bar{t})$ such that
$$
d_{\bar{t}}(\bar{x},x_0)< \frac{1}{3},
$$
$$
Q=|Rm(\bar{x},\bar{t})|> C+B(\bar{t}+t_0)^{-1},
$$
and
$$
|Rm(x,t)|\leq 4Q
$$
wherever $(-t_0<) \ \bar{t}-AQ^{-1}\leq t\leq \bar{t},\ \
d_{t}(x,\bar{x})\leq \frac{1}{10}A^{\frac{1}{2}}Q^{-\frac{1}{2}}$,
where $A$ tends to infinity with $B,C$. Thus we may take a blow-up
limit along the points $(\bar{x},\bar{t})$ with factors $Q$ and
get a non-flat ancient solution
$(M_{\infty},g_{ij}^{(\infty)}(t))$ with nonnegative curvature
operator and with the asymptotic volume ratio
$\nu_{M_{\infty}}(t)\geq w'>0$ for each $t\in(-\infty,0] \ \ $(by
(6.3.6)). This contradicts  Lemma 6.3.1.

\medskip
(ii) Let $B(w),\ C(w)$ be good for the first part of the theorem.
By the volume assumption at $t=0$ and the standard (relative)
volume comparison, we still have the estimate
$$
\Vol_0(B_0(x,r))\geq w'r^n \leqno{(6.3.6)'}  
$$
for each $x \in M$ with $d_0(x,x_0)\leq \frac{1}{3}$ and $r\leq
\frac{1}{3}$. We will show that $\xi=5^{-n}w'$, $B=B(5^{-n}w')$
and $C=C(5^{-n}w')$ are good for the second part of the theorem.

By continuity and the volume assumption at $t=0$, there is a
maximal subinterval $[-\tau,0]$ of the time interval
$[-\bar{\tau},0]$ such that
$$
\Vol_t(B_t(x_0,1))\geq w \geq 5^{-n}w'\quad \mbox{for all } \;
t\in[-\tau,0].
$$
This says that the assumptions of (i) hold with $5^{-n}w'$ in
place of $w$ and with $\tau$ in place of $t_0$. Thus the
conclusion of the part (i) gives us the estimate \be
|Rm(x,t)|\leq C+B(t+\tau)^{-1} 
\ee whenever $t\in (-\tau,0]$ and $d_t(x,x_0)\leq \frac{1}{4}$.

We need to show that one can choose a positive $\tau_0$ depending
only on $w$ and the dimension such that the maximal $\tau \geq
\min \{\bar{\tau},\tau_0\}$.

For $t\in(-\tau,0]$ and $\frac{1}{8}\leq d_t(x,x_0)\leq
\frac{1}{4}$, we use (6.3.7) and Lemma 3.4.1(ii) to get
$$
\frac{d}{dt}d_t(x,x_0)\geq
-10(n-1)(\sqrt{C}+(\sqrt{B}/\sqrt{t+\tau}))
$$
which further gives
$$
d_0(x,x_0)\geq
d_{-\tau}(x,x_0)-10(n-1)(\tau\sqrt{C}+2\sqrt{B\tau}).
$$
This means \be B_{(-\tau)}(x_0,\frac{1}{4})\supset
B_0\(x_0,\frac{1}{4}-10(n-1)(\tau\sqrt{C}+2\sqrt{B\tau})\).
\ee

Note that the scalar curvature $R\geq -C(n)$ for some constant
$C(n)$ depending only on the dimension since $Rm\geq -1.$ We have
\begin{align*}
&\frac{d}{dt}\Vol_t\(B_0\(x_0,\frac{1}{4}
-10(n-1)(\tau\sqrt{C}+2\sqrt{B\tau})\)\)\\
&  = \int_{B_0(x_0,\frac{1}{4}-10(n-1)(\tau\sqrt{C}
+2\sqrt{B\tau}))}(-R)dV_t\\
&  \leq C(n)\Vol_t\(B_0\(x_0,\frac{1}{4}
-10(n-1)(\tau\sqrt{C}+2\sqrt{B\tau})\)\)
\end{align*}
and then
\begin{align}
&  \Vol_t\(B_0\(x_0,\frac{1}{4}-10(n-1)(\tau\sqrt{C}
+2\sqrt{B\tau})\)\)\\
& \leq e^{C(n)\tau}\Vol_{(-\tau)}\(B_0\(x_0,\frac{1}{4}
-10(n-1)(\tau\sqrt{C}+2\sqrt{B\tau})\)\).\nn
\end{align}  
Thus by (6.3.6)$'$, (6.3.8) and (6.3.9),
\begin{align*}
& \Vol_{(-\tau)}(B_{(-\tau)})(x_0,1) \\
&  \geq \Vol_{(-\tau)}(B_{(-\tau)})\(x_0,\frac{1}{4}\)\\
&  \geq \Vol_{(-\tau)}\(B_0\(x_0,\frac{1}{4}-10(n-1)(\tau\sqrt{C}
+2\sqrt{B\tau})\)\)\\
&  \geq e^{-C(n)\tau}\Vol_{0}\(B_0\(x_0,\frac{1}{4}
-10(n-1)(\tau\sqrt{C}+2\sqrt{B\tau})\)\)\\
& \geq e^{-C(n)\tau}w'\(\frac{1}{4}
-10(n-1)(\tau\sqrt{C}+2\sqrt{B\tau})\)^n.
\end{align*}
So it suffices to choose $\tau_0=\tau_0(w)$ small enough so that
$$
e^{-C(n)\tau_0}\(\frac{1}{4}-10(n-1)(\tau_0\sqrt{C}
+2\sqrt{B\tau_0})\)^n\geq \(\frac{1}{5}\)^n.
$$
Therefore we have proved the theorem.
\end{pf}

\section{Ancient $\kappa$-solutions on Three-manifolds}

In this section we will determine the structures of ancient
$\kappa$-solutions on three-manifolds.

First of all, we consider a special class of ancient solutions ---
gradient shrinking Ricci solitons. Recall that a solution
$g_{ij}(t)$ to the Ricci flow is said to be a \textbf{gradient
shrinking Ricci soliton}\index{gradient shrinking Ricci soliton}
if there exists a smooth function $f$ such that \be
\nabla_i\nabla_j f+R_{ij}+\frac{1}{2t}g_{ij}=0\; \mbox{ for }\;
-\infty<t<0. 
\ee A gradient shrinking Ricci soliton moves by the one parameter
group of diffeomorphisms generated by $\nabla f$ and shrinks by a
factor at the same time.

The following result of Perelman \cite{P2} gives a complete
classification for all three-dimensional complete
$\kappa$-noncollapsed gradient shrinking solitons with bounded and
nonnegative sectional curvature.

\begin{lemma}[Classification of three-dimensional shrinking
solitons]\index{classification of three-dimensional shrinking
solitons} Let $(M,g_{ij}(t))$ be a nonflat gradient shrinking
soliton on a three-manifold. Suppose $(M,g_{ij}(t))$ has bounded
and nonnegative sectional curvature and is $\kappa$-noncollapsed
on all scales for some $\kappa>0$. Then $(M,g_{ij}(t))$ is one of
the following:
\begin{itemize}
\item[(i)] the round three-sphere $\mathbb{S}^3$, or a metric
quotient of $\mathbb{S}^3$; \item[(ii)] the round infinite
cylinder $\mathbb{S}^2\times \mathbb{R}$, or one of its
$\mathbb{Z}_2$ quotients.
\end{itemize}
\end{lemma}

\begin{pf}
We first consider the case that the sectional curvature of the
nonflat gradient shrinking soliton is not strictly positive. Let
us pull back the soliton to its universal cover. Then the
pull-back metric is again a nonflat ancient $\kappa$-solution. By
Hamilton's strong maximum principle (Theorem 2.2.1), we know that
the pull-back solution splits as the metric product of a
two-dimensional nonflat ancient $\kappa$-solution and
$\mathbb{R}$. Since the two-dimensional nonflat ancient
$\kappa$-solution is simply connected, it follows from Theorem
6.2.2 that it must be the round sphere $\mathbb{S}^2$. Thus, the
gradient shrinking soliton must be $\mathbb{S}^2\times\mathbb{R}
/\Gamma$, a metric quotient of the round cylinder.

For each $\sigma \in \Gamma$ and $(x,s) \in
\mathbb{S}^2\times\mathbb{R}$, we write $\sigma(x,s) =
(\sigma_1(x,s),\sigma_2(x,s))$ $\in \mathbb{S}^2\times\mathbb{R}$.
Since $\sigma$ sends lines to lines, and sends cross spheres to
cross spheres, we have $\sigma_2(x,s) = \sigma_2(y,s)$, for all
$x, y \in \mathbb{S}^2$. This says that $\sigma_2$ reduces to a
function of $s$ alone on $\mathbb{R}$. Moreover, for any $(x,s),
(x',s') \in \mathbb{S}^2\times\mathbb{R}$, since $\sigma$
preserves the distances between cross spheres $\mathbb{S}^2 \times
\{s\}$ and $\mathbb{S}^2 \times \{s'\}$, we have $|\sigma_2(x,s) -
\sigma_2(x',s')| = |s-s'|$. So the projection $\Gamma_2$ of
$\Gamma$ to the second factor $\mathbb{R}$ is an isometry subgroup
of $\mathbb{R}$. If the metric quotient
$\mathbb{S}^2\times\mathbb{R} /\Gamma$
 were compact, it would not be $\kappa$-noncollapsed on sufficiently
large scales as $t \rightarrow -\infty$. Thus the metric quotient
$\mathbb{S}^2\times\mathbb{R} /\Gamma$ is noncompact. It follows
that $\Gamma_2 = \{1\} \mbox{ or } \mathbb{Z}_2$. In particular,
there is a $\Gamma$-invariant cross sphere $\mathbb{S}^2$ in the
round cylinder $\mathbb{S}^2\times\mathbb{R}$. Denote it by
$\mathbb{S}^2\times \{0\}$. Then $\Gamma$ acts on the round
two-sphere $\mathbb{S}^2\times \{0\}$ isometrically without fixed
points. This implies $\Gamma$ is either $\{1\}$ or $\mathbb{Z}_2$.
Hence we conclude that the gradient shrinking soliton is either
the round cylinder $\mathbb{S}^2\times\mathbb{R}$, or
$\mathbb{R}\mathbb{P}^2 \times \mathbb{R}$, or the twisted product
$\mathbb{S}^2 \tilde{\times} \mathbb{R}$ where $\mathbb{Z}_2$
flips both $\mathbb{S}^2$ and $\mathbb{R}$.

We next consider the case that the gradient shrinking soliton is
compact and has strictly positive sectional curvature everywhere.
By the proof of Theorem 5.2.1 (see also Remark 5.2.8) we see that
the compact gradient shrinking soliton is getting round and tends
to a space form (with positive constant curvature) as the time
approaches the maximal time $t=0$. Since the shape of a gradient
shrinking Ricci soliton does not change up to reparametrizations
and homothetical scalings, the gradient shrinking soliton has to
be the round three-sphere $\mathbb{S}^3$ or a metric quotient of
$\mathbb{S}^3$.

Finally we want to exclude the case that the gradient shrinking
soliton is noncompact and has strictly positive sectional
curvature everywhere. The following argument follows the proof of
Lemma 1.2 of Perelman \cite{P2}.

Suppose there is a (complete three-dimensional) noncompact
$\kappa$-non\-collapsed gradient shrinking soliton $g_{ij}(t)$,
$-\infty<t<0$, with bounded and positive sectional curvature at
each $t\in (-\infty,0)$ and satisfying the equation (6.4.1). Then
as in (6.2.25), we have \be
\nabla_iR=2R_{ij}\nabla_jf. 
\ee

Fix some $t<0$, say $t=-1$, and consider a long shortest geodesic
$\gamma(s)$, $0\le s\le\bar{s}.$ Let $x_0=\gamma(0)$ and
$X(s)=\dot{\gamma}(s).$ Let $U(0)$ be any unit vector orthogonal
to $\dot{\gamma}(0)$ and translate $U(0)$ along $\gamma(s)$ to get
a parallel vector field $U(s)$, $0\le s\le\bar{s}$, on $\gamma$.
Set
$$
\widetilde{U}(s)=\left\{\arraycolsep=1.5pt\begin{array}{lll}
&  sU(s),\quad &  \mbox{for}\quad 0\le s\le 1,\\[2mm]
&  U(s),\quad &  \mbox{for}\quad 1\le s\le\bar{s}-1\\[2mm]
&  (\bar{s}-s)U(s),\quad &  \mbox{for}\quad \bar{s}-1\le s\le
\bar{s}.\end{array}\right.
$$
It follows from the second variation formula of arclength that
$$
\int_0^{\bar{s}}(|\dot{\widetilde{U}}(s)|^2
-R(X,\widetilde{U},X,\widetilde{U}))ds\ge 0.
$$
Since the curvature of the metric $g_{ij}(-1)$ is bounded, we
clearly have
$$
\int_0^{\bar{s}}R(X,U,X,U)ds\le {\rm const.}
$$
and then \be
\int_0^{\bar{s}}\Ric(X,X)ds\le {\rm const.}. 
\ee Moreover, since the curvature of the metric $g_{ij}(-1)$ is
positive, it follows from the Cauchy-Schwarz inequality that for
any unit vector field $Y$ along $\gamma$ and orthogonal to
$X(=\dot{\gamma}(s)),$ we have
\begin{align*}
\int_0^{\bar{s}}|\Ric(X,Y)|^2ds
&  \le\int_0^{\bar{s}}\Ric(X,X)\Ric(Y,Y)ds\\
&  \le {\rm const.}\,\cdot\int_0^{\bar{s}}\Ric(X,X)ds\\
&  \le {\rm const.}
\end{align*}
and then \be \int_0^{\bar{s}}|\Ric(X,Y)|ds
\le {\rm const.}\,\cdot(\sqrt{\bar{s}}+1). 
\ee {}From (6.4.1) we know
$$
\nabla_X\nabla_Xf+\Ric(X,X)-\frac{1}{2}=0
$$
and by integrating this equation we get
$$
X(f(\gamma(\bar{s})))-X(f(\gamma(0)))
+\int_0^{\bar{s}}\Ric(X,X)ds-\frac{1}{2}\bar{s}=0.
$$
Thus by (6.4.3) we deduce \be \frac{\bar{s}}{2}-{\rm const.}\,\le
\langle X,\nabla f(\gamma(\bar{s}))\rangle
\le\frac{\bar{s}}{2}+{\rm const.}.  
\ee Similarly by integrating (6.4.1) and using (6.4.4) we can
deduce \be |\langle Y,\nabla f(\gamma(\bar{s}))\rangle|
\le {\rm const.}\,\cdot (\sqrt{\bar{s}}+1). 
\ee These two inequalities tell us that at large distances from
the fixed point $x_0$ the function $f$ has no critical point, and
its gradient makes a small angle with the gradient of the distance
function from $x_0$.

Now from (6.4.2) we see that at large distances from $x_0$, $R$ is
strictly increasing along the gradient curves of $f$, in
particular
$$
\bar{R}=\limsup_{d_{(-1)}(x,x_0)\rightarrow+\infty}R(x,-1)>0.
$$
Let us choose a sequence of points $(x_k,-1)$ where
$R(x_k,-1)\rightarrow\bar{R}$. By the noncollapsing assumption we
can take a limit along this sequence of points of the gradient
soliton and get an ancient $\kappa$-solution defined on
$-\infty<t<0$. By Proposition 6.1.2, we deduce that the limiting
ancient $\kappa$-solution splits off a line. Since the soliton has
positive sectional curvature, we know from Gromoll-Meyer \cite{GM}
that it is orientable. Then it follows from Theorem 6.2.2 that the
limiting solution is the shrinking round infinite cylinder with
scalar curvature $\bar{R}$ at time $t=-1$. Since the limiting
solution exists on $(-\infty,0)$, we conclude that $\bar{R}\le 1$.
Hence
$$
R(x,-1)< 1
$$
when the distance from $x$ to the fixed $x_0$ is large enough on
the gradient shrinking soliton.

Let us consider the level surface $\{f=a\}$ of $f$. The second
fundamental form of the level surface is given by
\begin{align*}
h_{ij}&  = \left\langle\nabla_i\(\frac{\nabla f}{|\nabla
f|}\),e_j\right\rangle\\
&  = \nabla_i\nabla_j f/|\nabla f|,\qquad i,j=1,2,
\end{align*}
where $\{e_1,e_2\}$ is an orthonormal basis of the level surface.
By (6.4.1), we have
$$
\nabla_{e_i}\nabla_{e_i}f =\frac{1}{2}-\Ric(e_i,e_i)
\ge\frac{1}{2}-\frac{R}{2}>0,\qquad i=1,2,
$$
since for a three-manifold the positivity of sectional curvature
is equivalent to $R\ge 2\Ric$. It then follows from the first
variation formula that
\begin{align}
\frac{d}{da}{\rm Area}\,\{f=a\}
& = \int_{\{f=a\}}{\rm div}\,\(\frac{\nabla f}{|\nabla f|}\)\\
&  \ge \int_{\{f=a\}}\frac{1}{|\nabla f|}(1-R)\nn\\
&  > \int_{\{f=a\}}\frac{1}{|\nabla f|}(1-\bar{R})\nn\\
&  \ge 0\nn
\end{align}     
for $a$ large enough. We conclude that Area $\{f=a\}$ strictly
increases as $a$ increases. From (6.4.5) we see that for $s$ large
enough
$$
\left|\frac{df}{ds}-\frac{s}{2}\right|\le{\rm const.},
$$
and then
$$
\left|f-\frac{s^2}{4}\right|\le {\rm const.}\,\cdot(s+1).
$$
Thus we get from (6.4.7)
$$
\frac{d}{da}{\rm Area}\,\{f=a\}
>\frac{1-\bar{R}}{2\sqrt{a}}{\rm Area}\,\{f=a\}
$$
for $a$ large enough. This implies that
$$
\log {\rm Area}\,\{f=a\}>(1-\bar{R})\sqrt{a}-{\rm const.}
$$
for $a$ large enough. But it is clear from (6.4.7) that
Area$\,\{f=a\}$ is uniformly bounded from above by the area of the
round sphere of scalar curvature $\bar{R}$ for all large $a$. Thus
we deduce that $\bar{R}=1$. So \be
{\rm Area}\,\{f=a\}<8\pi  
\ee for $a$ large enough.

Denote by $X$ the unit normal vector to the level surface
$\{f=a\}$. By using the Gauss equation and (6.4.1), the intrinsic
curvature of the level surface $\{f=a\}$ can be computed as
\begin{align}
&\text{intrinsic curvature} \\
&  = R_{1212}+\det(h_{ij})  \nn\\
&  = \frac{1}{2}(R-2\Ric(X,X))
+\frac{\det({\rm Hess}\,(f))}{|\nabla f|^2} \nn\\
&  \le \frac{1}{2}(R-2\Ric(X,X))
+\frac{1}{4|\nabla f|^2}(\tr({\rm Hess}\,(f)))^2 \nn\\
&  = \frac{1}{2}(R-2\Ric(X,X))
+\frac{1}{4|\nabla f|^2}(1-(R-\Ric(X,X)))^2 \nn\\
& =\frac{1}{2}\left[1-\Ric(X,X)-(1-R+\Ric(X,X))
+\frac{(1-R+\Ric(X,X))^2}{2|\nabla f|^2}\right] \nn\\
&  < \frac{1}{2} \nn
\end{align} 
for sufficiently large $a$, since $(1-R+Ric(X,X))>0$ and $|\nabla
f|$ is large when $a$ is large. Thus the combination of (6.4.8)
and (6.4.9) gives a contradiction to the Gauss-Bonnet formula.

Therefore we have proved the lemma.
\end{pf}

As a direct consequence, there is a universal positive constant
$\kappa_0$ such that any nonflat three-dimensional gradient
shrinking soliton, which is also an ancient $\kappa$-solution, to
the Ricci flow must be $\kappa_0$-noncollapsed on all scales
unless it is a metric quotient of  round three-sphere. The
following result, claimed by Perelman in the section 1.5 of
\cite{P2}, shows that this property actually holds for all nonflat
three-dimensional ancient $\kappa$-solutions.

\begin{proposition}[{Universal noncollapsing}\index{universal
noncollapsing}] There exists a positive constant $\kappa_0$ with
the following property. Suppose we have a nonflat
three-dimensional ancient $\kappa$-solution for some $\kappa >0$.
Then either the solution is $\kappa_0$-noncollapsed on all scales,
or it is a metric quotient of the round three-sphere.
\end{proposition}

\begin{pf}
Let $g_{ij}(x,t), x\in M$ and $t\in(-\infty, 0]$, be a nonflat
ancient $\kappa$-solution for some $\kappa>0$. For an arbitrary
point $(p,t_0)\in M\times (-\infty,0]$, we define as in Chapter 3
that
\begin{align*}
& \tau=t_0-t,\quad \mbox{for }\;t<t_0,\\
& l(q,\tau)=\frac{1}{2\sqrt{\tau}}
\inf\bigg\{\int_0^\tau\sqrt{s}(R(\gamma(s),t_0-s)
+|\dot{\gamma}(s)|^2_{g_{ij}(t_0-s)})ds|\\
& \qquad \gamma:[0,\tau]\rightarrow M\;\mbox{ with }
\;\gamma(0)=p,\gamma(\tau)=q\bigg\},\\
&  \mbox{and }\widetilde{V}(\tau)
=\int_M(4\pi\tau)^{-\frac{3}{2}}\exp(-l(q,\tau))dV_{t_0-\tau}(q).
\end{align*}
Recall from (6.2.1) that for each $\tau>0$ we can find $q=q(\tau)$
such that $l(q,\tau)\le\frac{3}{2}$. In view of Lemma 6.4.1, we
may assume that the ancient $\kappa$-solution is not a gradient
shrinking Ricci soliton. Thus by (the proof of) Theorem 6.2.1, the
scalings of $g_{ij}(t_0-\tau)$ at $q(\tau)$ with factor
$\tau^{-1}$ converge along a subsequence of
$\tau\rightarrow+\infty$ to a nonflat gradient shrinking soliton
with nonnegative curvature operator which is $\kappa$-noncollapsed
on all scales. We now show that the limit has bounded curvature.

Denote the limiting nonflat gradient shrinking soliton by
$(\bar{M},\bar{g}_{ij}(x,t))$ with $-\infty < t \leq 0$. Note that
there holds the Li-Yau-Hamilton inequality (Theorem 2.5.4) on any
ancient $\kappa$-solution and in particular, the scalar curvature
of the ancient $\kappa$-solution is pointwise nondecreasing in
time. This implies that the scalar curvature of the limiting
soliton $(\bar{M},\bar{g}_{ij}(x,t))$ is still pointwise
nondecreasing in time. Thus we only need to show that the limiting
soliton has bounded curvature at $t=0$.

We argue by contradiction. By lifting to its orientable cover, we
may assume that $\bar{M}$ is orientable. Suppose the curvature of
the limiting soliton is unbounded at $t=0$. Of course in this case
the limiting soliton $\bar{M}$ is noncompact. Then by applying
Lemma 6.1.4, we can choose a sequence of points $x_j, j = 1, 2,
\ldots,$ divergent to infinity such that the scalar curvature
$\bar{R}$ of the limit satisfies
$$
\bar{R}(x_j,0) \geq j\; \mbox{ and }\; \bar{R}(x,0) \leq
4\bar{R}(x_j,0)
$$
for all $x\in B_{0}(x_j,{j}/{\sqrt{\bar{R}(x_{j},0)}})$ and $j =
1, 2,\, \ldots$.  Since the scalar curvature is nondecreasing in
time, we have \be
\bar{R}(x,t) \leq 4\bar{R}(x_j,0), 
\ee for all $x\in B_{0}(x_j,{j}/{\sqrt{\bar{R}(x_{j},0)}})$, all
$t\leq 0$ and $j = 1, 2,\, \ldots$. By combining with Hamilton's
compactness theorem (Theorem 4.1.5) and the
$\kappa$-noncollapsing, we know that a subsequence of the rescaled
solutions
$$
(\bar{M},\bar{R}(x_j,0)\bar{g}_{ij}(x, t/\bar{R}(x_j,0)),x_j), \ \
j=1, 2, \ldots,
$$
converges in the $C^{\infty}_{\rm loc}$ topology to a nonflat
smooth solution of the Ricci flow. Then Proposition 6.1.2 implies
that the new limit at the new time $\{t = 0\}$ must split off a
line. By pulling back the new limit to its universal cover and
applying Hamilton's strong maximum principle, we deduce that the
pull-back of the new limit on the universal cover splits off a
line for all time $t \leq 0 $. Thus by combining with Theorem
6.2.2 and the argument in the proof of Lemma 6.4.1, we further
deduce that the new limit is either the round cylinder
$\mathbb{S}^2 \times \mathbb{R}$ or the round
$\mathbb{R}\mathbb{P}^2 \times \mathbb{R}$. Since $\bar{M}$ is
orientable, the new limit must be $\mathbb{S}^2 \times
\mathbb{R}$. Since $(\bar{M},\bar{g}_{ij}(x,0))$ has nonnegative
curvature operator and the points $\{ x_j \}$ going to infinity
and $ \bar{R}(x_j,0)\rightarrow +\infty $, this gives a
contradiction to Proposition 6.1.1. So we have proved that the
limiting gradient shrinking soliton has bounded curvature at each
time.

Hence by Lemma 6.4.1, the limiting gradient shrinking soliton is
either the round three-sphere $\mathbb{S}^3$ or its metric
quotients, or the infinite cylinder $\mathbb{S}^2\times
\mathbb{R}$ or one of its $\mathbb{Z}_2$ quotients. If the
asymptotic gradient shrinking soliton is the round three-sphere
$\mathbb{S}^3$ or its metric quotients, it follows from Lemma
5.2.4 and Proposition 5.2.5 that the ancient $\kappa$-solution
must be round. Thus in the following we may assume the asymptotic
gradient shrinking soliton is the infinite cylinder
$\mathbb{S}^2\times\mathbb{R}$ or a $\mathbb{Z}_2$ quotient of
$\mathbb{S}^2\times\mathbb{R}$.

We now come back to consider the original ancient
$\kappa$-solution $(M,g_{ij}(x,$ $t))$. By rescaling, we can
assume that $R(x,t)\le 1$ for all $(x,t)$ satisfying
$d_{t_0}(x,p)\le 2$ and $t\in[t_0-1,t_0]$. We will argue as in the
proof of Theorem 3.3.2 (Perelman's no local collapsing theorem I)
to obtain a positive lower bound for $\Vol_{t_0}(B_{t_0}(p,1))$.

Denote by $\xi=\Vol_{t_0}(B_{t_0}(p,1))^\frac{1}{3}$. For any
$v\in T_pM$ we can find an $\mathcal{L}$-geodesic $\gamma(\tau)$,
starting at $p$, with $\lim_{\tau\rightarrow
0^+}\sqrt{\tau}\dot{\gamma}(\tau)=v.$ It follows from the
$\mathcal{L}$-geodesic equation (3.2.1) that
$$
\frac{d}{d\tau}(\sqrt{\tau}\dot{\gamma})
-\frac{1}{2}\sqrt{\tau}\nabla
R+2\Ric(\sqrt{\tau}\dot{\gamma},\cdot)=0.
$$
By integrating as before we see that for $\tau\le\xi$ with the
property $\gamma(\sigma)\in B_{t_0}(p,1)$ as long as
$\sigma<\tau$, there holds
$$
|\sqrt{\tau}\dot{\gamma}(\tau)-v|\le C\xi (|v|+1)
$$
where $C$ is some positive constant depending only on the
dimension. Without loss of generality, we may assume
$C\xi\le\frac{1}{4}$ and $\xi\le\frac{1}{100}$. Then for $v\in
T_pM$ with $|v|\le\frac{1}{4}\xi^{-\frac{1}{2}}$ and for
$\tau\le\xi$ with the property $\gamma(\sigma)\in B_{t_0}(p,1)$ as
long as $\sigma<\tau$, we have
\begin{align*}
d_{t_0}(p,\gamma(\tau))
&  \le \int_0^\tau|\dot{\gamma}(\sigma)|d\sigma\\
&  < \frac{1}{2}\xi^{-\frac{1}{2}}
\int_0^\tau\frac{d\sigma}{\sqrt{\sigma}}\\
&  = 1.
\end{align*}
This shows \be
\mathcal{L}\exp\left\{|v|\le\frac{1}{4}\xi^{-\frac{1}{2}}\right\}(\xi)
\subset B_{t_0}(p,1). 
\ee

We decompose Perelman's reduced volume $\widetilde{V}(\xi)$ as
\begin{align}
\widetilde{V}(\xi) &=\int_{\mathcal{L}\exp\left\{|v|
\le\frac{1}{4}\xi^{-\frac{1}{2}}\right\}(\xi)} \\
&\quad+\int_{M\setminus\mathcal{L}\exp\left\{|v|
\leq\frac{1}{4}\xi^{-\frac{1}{2}}\right\}(\xi)}(4\pi\xi)^{-\frac{3}{2}}
\exp(-l(q,\xi))dV_{t_0-\xi}(q). \nn
\end{align}
By using (6.4.11) and the metric evolution equation of the Ricci
flow, the first term on the RHS of (6.4.12) can be estimated by
\begin{align*}
&
\int_{\mathcal{L}\exp\{|v|\le\frac{1}{4}\xi^{-\frac{1}{2}}\}(\xi)}
(4\pi\xi)^{-\frac{3}{2}}\exp(-l(q,\xi))dV_{t_0-\xi}(q)\\
&  \le\int_{B_{t_0}(p,1)}(4\pi\xi)^{-\frac{3}{2}}e^{3\xi}dV_{t_0}(q)\\
&  = (4\pi)^{-\frac{3}{2}}e^{3\xi}{\xi}^{\frac{3}{2}}\\
&  < \xi^\frac{3}{2},
\end{align*}
while by using Theorem 3.2.7 (Perelman's Jacobian comparison
theorem), the second term on the RHS of (6.4.12) can be estimated
by
\begin{align}
& \int_{M\setminus\mathcal{L}\exp
\left\{|v|\leq\frac{1}{4}\xi^{-\frac{1}{2}}\right\}(\xi)}
(4\pi\xi)^{-\frac{3}{2}}
\exp(-l(q,\xi))dV_{t_0-\xi}(q)\\
& \le\int_{\{|v|>\frac{1}{4}\xi^{-\frac{1}{2}}\}}
(4\pi\tau)^{-\frac{3}{2}}
\exp(-l(\tau))\mathcal{J}(\tau)|_{\tau=0}dv \nn\\
&  =
(4\pi)^{-\frac{3}{2}}\int_{\{|v|>\frac{1}{4}\xi^{-\frac{1}{2}}\}}
\exp(-|v|^2)dv \nn\\
&  < \xi^\frac{3}{2} \nn
\end{align}  
since $\lim_{\tau\rightarrow
0^+}\tau^{-\frac{3}{2}}\mathcal{J}(\tau)=1$ and
$\lim_{\tau\rightarrow 0^+}l(\tau)=|v|^2$ by (3.2.18) and (3.2.19)
respectively. Thus we obtain \be
\widetilde{V}(\xi)<2\xi^\frac{3}{2}. 
\ee

On the other hand, we recall that there exist a sequence
$\tau_k\rightarrow+\infty$ and a sequence of points $q(\tau_k)\in
M$ with $l(q(\tau_k),\tau_k)\le\frac{3}{2}$ so that the scalings
of the ancient $\kappa$-solution at $q(\tau_k)$ with factor
$\tau_k^{-1}$ converge to either round $\mathbb{S}^2\times
\mathbb{R}$ or one of its $\mathbb{Z}_2$ quotients. For
sufficiently large $k$, we construct a path
$\gamma:\quad[0,2\tau_k]\rightarrow M$, connecting $p$ to any
given point $q\in M$, as follows: the first half path
$\gamma|_{[0,\tau_k]}$ connects $p$ to $q(\tau_k)$ such that
$$
l(q(\tau_k),\tau_k)
=\frac{1}{2\sqrt{\tau_k}}\int_0^{\tau_k}\sqrt{\tau}
(R+|\dot{\gamma}(\tau)|^2)d\tau\le2,
$$
and the second half path $\gamma|_{[\tau_k,2\tau_k]}$ is a
shortest geodesic connecting $q(\tau_k)$ to $q$ with respect to
the metric $g_{ij}(t_0-\tau_k)$. Note that the rescaled metric
$\tau_k^{-1}g_{ij}(t_0-\tau)$ over the domain
$B_{t_0-\tau_k}(q(\tau_k),\sqrt{\tau_k})\times[t_0-2\tau_k,t_0-\tau_k]$
is sufficiently close to the round $\mathbb{S}^2\times \mathbb{R}$
or its $\mathbb{Z}_2$ quotients. Then there is a universal
positive constant $\beta$ such that
\begin{align*}
l(q,2\tau_k)& \le \frac{1}{2\sqrt{2\tau_k}}\(\int_0^{\tau_k}
+\int_{\tau_k}^{2\tau_k}\)\sqrt{\tau}
(R+|\dot{\gamma}(\tau)|^2)d\tau\\
&  \le \sqrt{2}+\frac{1}{2\sqrt{2\tau_k}}
\int_{\tau_k}^{2\tau_k}\sqrt{\tau}(R+|\dot{\gamma}(\tau)|^2)d\tau\\
&  \le \beta
\end{align*}
for all $q\in B_{t_0-\tau_k}(q(\tau_k),\sqrt{\tau_k})$. Thus
\begin{align*}
\widetilde{V}(2\tau_k)&  = \int_M (4\pi(2\tau_k))^{-\frac{3}{2}}
\exp(-l(q,2\tau_k))dV_{t_0-2\tau_k}(q)\\
&  \ge e^{-\beta}\int_{B_{t_0-\tau_k}(q(\tau_k),
\sqrt{\tau_k})}(4\pi(2\tau_k))^{-\frac{3}{2}}dV_{t_0-2\tau_k}(q)\\
&  \ge \tilde{\beta}
\end{align*}
for some universal positive constant $\tilde{\beta}$. Here we have
used the curvature estimate (6.2.6). By combining with the
monotonicity of Perelman's reduced volume (Theorem 3.2.8) and
(6.4.14), we deduce that
$$
\tilde{\beta}\le\widetilde{V}(2\tau_k)
\le\widetilde{V}(\xi)<2\xi^\frac{3}{2}.
$$
This proves
$$
\Vol_{t_0}(B_{t_0}(p,1))\ge \kappa_0>0
$$
for some universal positive constant $\kappa_0$. So we have proved
that the ancient $\kappa$-solution is also an ancient
$\kappa_0$-solution.
\end{pf}

The important Li-Yau-Hamilton inequality gives rise to a parabolic
Harnack estimate (Corollary 2.5.7) for solutions of the Ricci flow
with bounded and nonnegative curvature operator. As explained in
the previous section, the no local collapsing theorem of Perelman
implies a volume lower bound from a curvature upper bound, while
the estimate in the previous section implies a curvature upper
bound from a volume lower bound. The combination of these two
estimates as well as the Li-Yau-Hamilton inequality will give an
important elliptic type property for three-dimensional ancient
$\kappa$-solutions. This elliptic type property was first
implicitly given by Perelman in \cite{P1} and it will play a
crucial role in the analysis of singularities.

\begin{theorem}[{Elliptic type estimate}\index{elliptic
type estimate}] There exist a positive constant $\eta$ and a
positive increasing function $\omega:\quad [0,+\infty)\rightarrow
(0,+\infty)$ with the following properties. Suppose we have a
three-dimensional ancient $\kappa$-solution
$(M,g_{ij}(t)),-\infty<t\le 0,$ for some $\kappa>0$. Then
\begin{itemize}
\item[(i)]for every $x,y\in M$ and $t\in(-\infty,0]$, there holds
$$
R(x,t)\le R(y,t)\cdot \omega(R(y,t)d_t^2(x,y));
$$
\item[(ii)] for all $x\in M$ and $t\in (-\infty,0]$, there hold
$$
|\nabla R|(x,t)\le\eta R^\frac{3}{2}(x,t)\;  \mbox{ and }\;
\left|\frac{\partial R}{\partial t}\right|(x,t)\le \eta R^2(x,t).
$$
\end{itemize}
\end{theorem}

\begin{pf} \

(i) Consider a three-dimensional nonflat ancient $\kappa$-solution
$g_{ij}(x,t)$ on $M \times (-\infty,0]$. In view of Proposition
6.4.2, we may assume that the ancient solution is universal
$\kappa_0$-noncollapsed. Obviously we only need to establish the
estimate at $t=0$. Let $y$ be an arbitrarily fixed point in $M$.
By rescaling, we can assume $R(y,0)=1$.

Let us first consider the case that $\sup\{R(x,0)d_0^2(x,y)\ |\
x\in M\}>1$. Define $z$ to be the closest point to $y$ (at time
$t=0$) satisfying $R(z,0)d_0^2(z,y)=1$. We want to bound
$R(x,0)/R(z,0)$ from above for $x\in
B_0(z,2R(z,0)^{-\frac{1}{2}}).$

Connect $y$ and $z$ by a shortest geodesic and choose a point
$\tilde{z}$ lying on the geodesic satisfying
$d_0(\tilde{z},z)=\frac{1}{4}R(z,0)^{-\frac{1}{2}}$. Denote by $B$
the ball centered at $\tilde{z}$ and with radius
$\frac{1}{4}R(z,0)^{-\frac{1}{2}}$ (with respect to the metric at
$t=0$). Clearly the ball $B$ lies in
$B_0(y,R(z,0)^{-\frac{1}{2}})$ and lies outside
$B_0(y,\frac{1}{2}R(z,0)^{-\frac{1}{2}})$. Thus for $x\in B$, we
have
$$
R(x,0)d_0^2(x,y)\le 1\quad \mbox{and} \quad
d_0(x,y)\ge\frac{1}{2}R(z,0)^{-\frac{1}{2}}
$$
and hence
$$
R(x,0)\le\frac{1}{(\frac{1}{2}R(z,0)^{-\frac{1}{2}})^2} \quad
\mbox{for all } \;x\in B.
$$
Then by the Li-Yau-Hamilton inequality and the
$\kappa_0$-noncollapsing, we have
$$
\Vol_0(B)\ge\kappa_0\(\frac{1}{4}R(z,0)^{-\frac{1}{2}}\)^3,
$$
and then
$$
\Vol_0(B_0(z,8R(z,0)^{-\frac{1}{2}}))
\ge\frac{\kappa_0}{2^{15}}(8R(z,0)^{-\frac{1}{2}})^3.
$$
So by Theorem 6.3.3(ii), there exist positive constants
$B(\kappa_0), C(\kappa_0),$ and $\tau_0(\kappa_0)$ such that \be
R(x,0)\le (C(\kappa_0)+\frac{B(\kappa_0)}{\tau_0(\kappa_0)})R(z,0)
\ee for all $x\in B_0(z,2R(z,0)^{-\frac{1}{2}})$.

We now consider the remaining case. If $R(x,0)d_0^2(x,y)\le 1$
everywhere, we choose a point $z$ satisfying $\sup\{R(x,0)\ |\
x\in M\}\leq2R(z,0)$. Obviously we also have the estimate (6.4.15)
in this case.

We next want to bound $R(z,0)$ for the chosen $z\in M$. By
(6.4.15) and the Li-Yau-Hamilton inequality, we have
$$
R(x,t)\le (C(\kappa_0)+\frac{B(\kappa_0)}{\tau_0(\kappa_0)})R(z,0)
$$
for all $x\in B_0(z,2R(z,0)^{-\frac{1}{2}})$ and all $t\le 0$. It
then follows from the local derivative estimates of Shi that
$$
\frac{\partial R}{\partial t}(z,t)
\le\widetilde{C}(\kappa_0)R(z,0)^2,\quad \mbox{for all }\;
-R^{-1}(z,0)\leq t\leq 0
$$
which implies \be
R(z,-cR^{-1}(z,0))\ge cR(z,0)  
\ee for some small positive constant $c$ depending only on
$\kappa_0$. On the other hand, by using the Harnack estimate in
Corollary 2.5.7, we have \be
1=R(y,0)\ge\widetilde{c}R(z,-cR^{-1}(z,0)) 
\ee for some small positive constant $\widetilde{c}$ depending
only on $\kappa_0$, since $d_0(y,z)\le R(z,0)^{-\frac{1}{2}}$ and
the metric $g_{ij}(t)$ is equivalent on
$$
B_0(z,2R(z,0)^{-\frac{1}{2}})\times[-cR^{-1}(z,0),0]
$$
with $c>0$ small enough. Thus we get from (6.4.16) and (6.4.17)
that \be
R(z,0)\le\widetilde{A} 
\ee for some positive constant $\widetilde{A}$ depending only on
$\kappa_0$.

Since $B_0(z,2R(z,0)^{-\frac{1}{2}})\supset
B_0(y,R(z,0)^{-\frac{1}{2}})$ and $R(z,0)^{-\frac{1}{2}}\ge
(\widetilde{A})^{-\frac{1}{2}}$, the combination of (6.4.15) and
(6.4.18) gives \be R(x,0)\le
(C(\kappa_0)+\frac{B(\kappa_0)}{\tau_0(\kappa_0)})
\widetilde{A}  
\ee whenever $x\in B_0(y,(\widetilde{A})^{-\frac{1}{2}})$. Then by
the $\kappa_0$-noncollapsing there exists a positive constant
$r_0$ depending only on $\kappa_0$ such that
$$
\Vol_0(B_0(y,r_0))\ge \kappa_0r_0^3.
$$
For any fixed $R_0\ge r_0$, we then have
$$
\Vol_0(B_0(y,R_0))\ge \kappa_0r_0^3
=\kappa_0(\frac{r_0}{R_0})^3\cdot R_0^3.
$$
By applying Theorem 6.3.3 (ii) again and noting that the constant
$\kappa_0$ is universal, there exists a positive constant
$\omega(R_0)$ depending only on $R_0$ such that
$$
R(x,0)\le\omega(R_0^2)\qquad \mbox{for all }\; x\in
B_0(y,\frac{1}{4}R_0).
$$
This gives the desired estimate.

\medskip
(ii) This follows immediately from  conclusion (i), the
Li-Yau-Hamilton inequality and the local derivative estimate of
Shi.
\end{pf}

As a consequence, we have the following compactness result due to
Perelman \cite{P1}.

\begin{corollary}[{Compactness of ancient $\kappa_0$-solutions}%
\index{compactness of ancient $\kappa_0$-solutions}] The set of
nonflat three-dimensional ancient $\kappa_0$-solutions is compact
modulo scaling in the sense that for any sequence of such
solutions and marking points $(x_k,0)$ with $R(x_k,0)=1$, we can
extract a $C^\infty_{loc}$ converging subsequence whose limit is
also an ancient $\kappa_0$-solution.
\end{corollary}

\begin{pf}
Consider any sequence of three-dimensional ancient
$\kappa_0$-solutions and marking points $(x_k,0)$ with
$R(x_k,0)=1$. By Theorem 6.4.3(i), the Li-Yau-Hamilton inequality
and Hamilton's compactness theorem (Theorem 4.1.5), we can extract
a $C^\infty_{loc}$ converging subsequence such that the limit
$(\bar{M},\bar{g}_{ij}(x,t))$, with $-\infty <t \leq 0$, is an
ancient solution to the Ricci flow with nonnegative curvature
operator and $\kappa_0$-noncollapsed on all scales. Since any
ancient $\kappa_0$-solution satisfies the Li-Yau-Hamilton
inequality, it implies that the scalar curvature $\bar{R}(x,t)$ of
the limit $(\bar{M},\bar{g}_{ij}(x,t))$ is pointwise nondecreasing
in time. Thus it remains to show that the limit solution has
bounded curvature at $t=0$.

Obviously we may assume the limiting manifold $\bar{M}$ is
noncompact. By pulling back the limiting solution to its
orientable cover, we can assume that the limiting manifold
$\bar{M}$ is orientable. We now argue by contradiction. Suppose
the scalar curvature $\bar{R}$ of the limit at $t=0$ is unbounded.

By applying Lemma 6.1.4, we can choose a sequence of points $x_j
\in \bar{M}, j = 1, 2, \ldots,$ divergent to infinity such that
the scalar curvature $\bar{R}$ of the limit satisfies
$$
\bar{R}(x_j,0) \geq j\; \mbox{ and }\; \bar{R}(x,0) \leq
4\bar{R}(x_j,0)
$$
for all $j = 1, 2,\, \ldots,$ and $x\in
B_{0}(x_j,{j}/{\sqrt{\bar{R}(x_{j},0)}})$. Then from the fact that
the limiting scalar curvature $\bar{R}(x,t)$ is pointwise
nondecreasing in time, we have \be
\bar{R}(x,t) \leq 4\bar{R}(x_j,0) 
\ee for all $j = 1, 2,\, \ldots$, $x\in
B_{0}(x_j,{j}/{\sqrt{\bar{R}(x_{j},0)}})$ and $t\leq 0$. By
combining with Hamilton's compactness theorem (Theorem 4.1.5) and
the $\kappa_0$-noncollapsing, we know that a subsequence of the
rescaled solutions
$$
(\bar{M},\bar{R}(x_j,0)\bar{g}_{ij}(x, t/\bar{R}(x_j,0)),x_j), \ \
j=1, 2, \ldots,
$$
converges in the $C^{\infty}_{loc}$ topology to a nonflat smooth
solution of the Ricci flow. Then Proposition 6.1.2 implies that
the new limit at the new time $\{t = 0\}$ must split off a line.
By pulling back the new limit to its universal cover and applying
Hamilton's strong maximum principle, we deduce that the pull-back
of the new limit on the universal cover splits off a line for all
time $t \leq 0$. Thus by combining with Theorem 6.2.2 and the
argument in the proof of Lemma 6.4.1, we further deduce that the
new limit is either the round cylinder $\mathbb{S}^2 \times
\mathbb{R}$ or the round $\mathbb{R}\mathbb{P}^2 \times
\mathbb{R}$. Since $\bar{M}$ is orientable, the new limit must be
$\mathbb{S}^2 \times \mathbb{R}$. Moreover, since
$(\bar{M},\bar{g}_{ij}(x,0))$ has nonnegative curvature operator
and the points $\{ x_j \}$ are going to infinity and $
\bar{R}(x_j,0)\rightarrow +\infty $, this gives a contradiction to
Proposition 6.1.1. So we have proved that the limit
$(\bar{M},\bar{g}_{ij}(x,t))$ has uniformly bounded curvature.
\end{pf}

Arbitrarily fix $\varepsilon>0$. Let $g_{ij}(x,t)$ be a nonflat
ancient $\kappa$-solution on a three-manifold $M$ for some $\kappa
>0$. We say that a point $x_0\in M$ is the {\bf center of an
evolving $\varepsilon$-neck}\index{center of an evolving
$\varepsilon$-neck} at $t=0$, if the solution $g_{ij}(x,t)$ in the
set $\{(x,t)\ |\ -\varepsilon^{-2} Q^{-1}<t\le 0,
d_t^2(x,x_0)<\varepsilon^{-2} Q^{-1}\}$, where $Q=R(x_0,0),$ is,
after scaling with factor $Q$, $\varepsilon$-close (in the
$C^{[\varepsilon^{-1}]}$ topology) to the corresponding set of the
evolving round cylinder having scalar curvature one at $t=0$.

As another consequence of the elliptic type estimate, we have the
following global structure result obtained by Perelman in
\cite{P1} for noncompact ancient $\kappa$-solutions.

\begin{corollary} [{Perelman \cite{P1}}]
For any $\varepsilon>0$ there exists $C=C(\varepsilon)>0$, such
that if $g_{ij}(t)$ is a nonflat ancient $\kappa$-solution on a
noncompact three-manifold $M$ for some $\kappa>0$, and
$M_\varepsilon$ denotes the set of points in $M$ which are not
centers of evolving $\varepsilon$-necks at $t=0$, then at $t=0$,
either the whole manifold $M$ is the round cylinder $\mathbb{S}^2
\times \mathbb{R}$ or its $\mathbb{Z}_2$ metric quotients, or
$M_\varepsilon$ satisfies the following
\begin{itemize}
\item[(i)] $M_\varepsilon$ is compact, \item[(ii)] ${\rm diam}\,
M_\varepsilon \le CQ^{-\frac{1}{2}}$ and $C^{-1}Q\le R(x,0)\le
CQ$, whenever $x\in M_\varepsilon$, where $Q=R(x_0,0)$ for some
$x_0\in\partial M_\varepsilon$.
\end{itemize}
\end{corollary}

\begin{pf}
We first consider the easy case that the curvature operator of the
ancient $\kappa$-solution has a nontrivial null vector somewhere
at some time. Let us pull back the solution to its universal
cover. By applying Hamilton's strong maximum principle and Theorem
6.2.2, we see that the universal cover is the evolving round
cylinder $\mathbb{S}^2 \times \mathbb{R}$. Thus in this case, by
the argument in the proof of Lemma 6.4.1, we conclude that the
ancient $\kappa$-solution is either isometric to the round
cylinder $\mathbb{S}^2 \times \mathbb{R}$ or one of its
$\mathbb{Z}_2$ metric quotients (i.e., $\mathbb{R}\mathbb{P}^2
\times \mathbb{R}$, or the twisted product $\mathbb{S}^2
\tilde{\times} \mathbb{R}$ where $\mathbb{Z}_2$ flips both
$\mathbb{S}^2$, or $\mathbb{R}$).

We then assume that the curvature operator of the nonflat ancient
$\kappa$-solution is positive everywhere. Firstly we want to show
$M_\varepsilon$ is compact. We argue by contradiction. Suppose
there exists a sequence of points $z_k,\ k=1,2,\ldots$, going to
infinity (with respect to the metric $g_{ij}(0)$) such that each
$z_k$ is not the center of any evolving $\varepsilon$-neck. For an
arbitrarily fixed point $z_0\in M$, it follows from Theorem
6.4.3(i) that
$$
0<R(z_0,0)\leq R(z_k,0)\cdot \omega(R(z_k,0)d^2_0(z_k,z_0))
$$
which implies that
$$
\lim_{k\rightarrow\infty}R(z_k,0)d^2_0(z_k,z_0)=+\infty.
$$
Since the sectional curvature of the ancient $\kappa$-solution is
positive everywhere, the underlying manifold is diffeomorphic to
$\mathbb{R}^3$, and in particular, orientable. Then as before, by
Proposition 6.1.2, Theorem 6.2.2 and Corollary 6.4.4, we conclude
that $z_k$ is the center of an evolving $\varepsilon$-neck for $k$
sufficiently large. This is a contradiction, so we have proved
that $M_{\varepsilon}$ is compact.

Again, we notice that $M$ is diffeomorphic to $\mathbb{R}^3$ since
the curvature operator is positive. According to the resolution of
the Schoenflies conjecture in three-dimensions, every
approximately round two-sphere cross-section through the center of
an evolving $\varepsilon$-neck divides $M$ into two parts such
that one of them is diffeomorphic to the three-ball
$\mathbb{B}^3$. Let $\varphi$ be the Busemann function on $M$, it
is a standard fact that $\varphi$ is convex and proper. Since
$M_{\varepsilon}$ is compact, $M_\varepsilon$ is contained in a
compact set $K=\varphi^{-1}((-\infty,A])$ for some large $A$. We
note that each point $x\in M\setminus M_\varepsilon$ is the center
of an $\varepsilon$-neck. It is clear that there is an
$\varepsilon$-neck $N$ lying entirely outside $K$.  Consider a
point $x$ on one of the boundary components of the
$\varepsilon$-neck $N$. Since $x \in M\setminus M_{\varepsilon}$,
there is an $\varepsilon$-neck adjacent to the initial
$\varepsilon$-neck, producing a longer neck. We then take a point
on the boundary of the second $\varepsilon$-neck and continue.
This procedure can either terminate when we get into
$M_{\varepsilon}$ or go on infinitely to produce a semi-infinite
(topological) cylinder. The same procedure can be repeated for the
other boundary component of the initial $\varepsilon$-neck. This
procedure will give a maximal extended neck $\tilde{N}$. If
$\tilde{N}$ never touches $M_\varepsilon$, the manifold will be
diffeomorphic to the standard infinite cylinder, which is a
contradiction. If both ends of $\tilde{N}$ touch $M_\varepsilon$,
then there is a geodesic connecting two points of
$M_{\varepsilon}$ and passing through $N$. This is impossible
since the function $\varphi$ is convex. So we conclude that one
end of $\tilde{N}$ will touch $M_\varepsilon$ and the other end
will tend to infinity to produce a semi-infinite (topological)
cylinder. Thus we can find an approximately round two-sphere
cross-section which encloses the whole set $M_{\varepsilon}$ and
touches some point $x_{0}\in
\partial M_{\varepsilon}$. We next want to show that
$R(x_0,0)^{\frac{1}{2}}\cdot diam(M_{\varepsilon})$ is bounded
from above by some positive constant $C=C(\varepsilon)$ depending
only on $\varepsilon$.

Suppose not; then there exists a sequence of nonflat noncompact
three-dimensional ancient $\kappa$-solutions with positive
curvature operator such that for the above chosen points $x_0\in
\partial M_{\varepsilon}$ there would hold
\be R(x_0,0)^{\frac{1}{2}}\cdot
{\rm diam}\,(M_{\varepsilon})\rightarrow +\infty. 
\ee By Proposition 6.4.2, we know that the ancient solutions are
$\kappa_0$-noncollapsed on all scales for some universal positive
constant $\kappa_0$. Let us dilate the ancient solutions around
the points $x_0$ with the factors $R(x_0,0)$. By Corollary 6.4.4,
we can extract a convergent subsequence. From the choice of the
points $x_0$ and (6.4.21), the limit has at least two ends. Then
by Toponogov's splitting theorem the limit is isometric to
$X\times\mathbb{R} $ for some nonflat two-dimensional ancient
$\kappa_0$-solution $X$. Since $M$ is orientable, we conclude from
Theorem 6.2.2 that limit must be the evolving round cylinder
$\mathbb{S}^2\times\mathbb{R}$. This contradicts the fact that
each chosen point $x_0$ is not the center of any evolving
$\varepsilon$-neck. Therefore we have proved
$$
{\rm diam}\,(M_{\varepsilon})\leq CQ^{-\frac{1}{2}}
$$
for some positive constant $C=C(\varepsilon)$ depending only on
$\varepsilon$, where $Q=R(x_0,0)$.

Finally by combining this diameter estimate with Theorem 6.4.3(i),
we immediately deduce
$$
\widetilde{C}^{-1}Q\leq R(x,0)\leq \widetilde{C}Q,\;\mbox{
whenever }\;x\in M_{\varepsilon},
$$
for some positive constant $\widetilde{C}$ depending only on
$\varepsilon$.
\end{pf}

We now can describe the canonical structures for three-dimensional
nonflat (compact or noncompact) ancient $\kappa$-solutions. The
following theorem was given by Perelman in the section 1.5 of
\cite{P2}. Recently in \cite{CZ05F}, this canonical neighborhood
result has been extended to four-dimensional ancient
$\kappa$-solutions with isotropic curvature pinching.

\begin{theorem}[{Canonical neighborhood theorem}\index{canonical
neighborhood theorem}] For any $\varepsilon>0$ one can find
positive constants $C_1=C_1(\varepsilon)$ and
$C_2=C_2(\varepsilon)$ with the following property. Suppose we
have a three-dimensional nonflat $($compact or noncompact$)$
ancient $\kappa$-solution $(M,g_{ij}(x,t))$. Then either the
ancient solution is the round $\mathbb{R}\mathbb{P}^2 \times
\mathbb{R}$, or every point $(x,t)$ has an open neighborhood $B$,
with $B_t(x,r)\subset B \subset B_t(x,2r)$ for some
$0<r<C_1R(x,t)^{-\frac{1}{2}}$, which falls into one of the
following three categories:
\begin{itemize}
\item[(a)] $B$ is an {\bf evolving $\varepsilon$-neck} $($in the
sense that it is the slice at the time $t$ of the parabolic region
$\{(x',t')\ |\ x' \in B, t' \in [t - \varepsilon^{-2}R(x,t)^{-1},
t] \}$ which is, after scaling with factor $R(x,t)$ and shifting
the time $t$ to zero, $\varepsilon$-close $($in the
$C^{[\varepsilon^{-1}]}$ topology$)$ to the subset
$(\mathbb{S}^2\times \mathbb{I})\times [-\varepsilon^{-2},0]$ of
the evolving standard round cylinder with scalar curvature $1$ and
length $2\varepsilon^{-1}$ to $\mathbb{I}$ at the time zero$),$ or
\index{evolving $\varepsilon$-neck} \item[(b)] $B$ is an {\bf
evolving $\varepsilon$-cap} $($in the sense that it is the time
slice at the time $t$ of an evolving metric on $\mathbb{B}^3$ or
$\mathbb{R}\mathbb{P}^3\setminus\bar{\mathbb{B}^3}$ such that the
region outside some suitable compact subset of $\mathbb{B}^3$ or
$\mathbb{R}\mathbb{P}^3\setminus\bar{\mathbb{B}^3}$ is an evolving
$\varepsilon$-neck$),$ or \index{evolving $\varepsilon$-cap}
\item[(c)] $B$ is a compact manifold $($without boundary$)$ with
positive sectional curvature $($thus it is diffeomorphic to the
round three-sphere $\mathbb{S}^3$ or a metric quotient of
$\mathbb{S}^3);$
\end{itemize}
furthermore, the scalar curvature of the ancient $\kappa$-solution
on $B$ at time $t$ is between $C_2^{-1}R(x,t)$ and $C_2R(x,t)$,
and the volume of $B$ in case {\rm (a)} and case {\rm (b)}
satisfies
$$
(C_2R(x,t))^{-\frac{3}{2}} \leq \Vol_t(B) \leq \varepsilon r^3.
$$
\end{theorem}

\begin{pf}
As before, we first consider the easy case that the curvature
operator has a nontrivial null vector somewhere at some time. By
pulling back the solution to its universal cover and applying
Hamilton's strong maximum principle and Theorem 6.2.2, we deduce
that the universal cover is the evolving round cylinder
$\mathbb{S}^2 \times \mathbb{R}$. Then exactly as before, by the
argument in the proof of Lemma 6.4.1, we conclude that the ancient
$\kappa$-solution is isometric to the round $\mathbb{S}^2 \times
\mathbb{R}$, $\mathbb{R}\mathbb{P}^2 \times \mathbb{R}$, or the
twisted product $\mathbb{S}^2 \tilde{\times} \mathbb{R}$ where
$\mathbb{Z}_2$ flips both $\mathbb{S}^2$ and $\mathbb{R}$. Clearly
each point of the round cylinder $\mathbb{S}^2 \times \mathbb{R}$
or the twisted product $\mathbb{S}^2 \tilde{\times} \mathbb{R}$
has a neighborhood falling into the category (a) or (b) (over
$\mathbb{R}\mathbb{P}^3\setminus\bar{\mathbb{B}^3}$).

We now assume that the curvature operator of the nonflat ancient
$\kappa$-solution is positive everywhere. Then the manifold is
orientable by the Cheeger-Gromoll theorem \cite{CG} for the
noncompact case or the Synge theorem \cite{CE} for the compact
case.

Without loss of generality, we may assume $\varepsilon$ is
suitably small, say $0<\varepsilon < \frac{1}{100}$. If the
nonflat ancient $\kappa$-solution is noncompact, the conclusions
follow immediately from the combination of Corollary 6.4.5 and
Theorem 6.4.3(i). Thus we may assume the nonflat ancient
$\kappa$-solution is compact. By Proposition 6.4.2, either the
compact ancient $\kappa$-solution is isometric to a metric
quotient of the round $\mathbb{S}^3$, or it is
$\kappa_0$-noncollapsed on all scales for the universal positive
constant $\kappa_0$. Clearly each point of a metric quotient of
the round $\mathbb{S}^3$ has a neighborhood falling into category
(c). Thus we may further assume the ancient $\kappa$-solution is
also $\kappa_0$-noncollapsing.

We argue by contradiction. Suppose that for some $\varepsilon \in
(0, \frac{1}{100})$, there exist a sequence of compact orientable
ancient $\kappa_0$-solutions $(M_k,g_k)$ with positive curvature
operator, a sequence of points $(x_k,0)$ with $x_k\in M_k$ and
sequences of positive constants $C_{1k}\rightarrow\infty$ and
$C_{2k}=\omega(4C^2_{1k})$, with the function $\omega$ given in
Theorem 6.4.3, such that for every radius $r$,
$0<r<C_{1k}R(x_k,0)^{-\frac{1}{2}}$, any open neighborhood $B$,
with $B_0(x_k,r) \subset B \subset B_0(x_k,2r)$, does not fall
into one of the three categories (a), (b) and (c), where in the
case (a) and case (b), we require the neighborhood $B$ to satisfy
the volume estimate
$$
(C_{2k}R(x_k,0))^{-\frac{3}{2}} \leq \Vol_0(B) \leq \varepsilon
r^3.
$$

By Theorem 6.4.3(i) and the choice of the constants $C_{2k}$ we
see that the diameter of each $M_k$ at $t=0$ is at least
$C_{1k}R(x_k,0)^{-\frac{1}{2}}$; otherwise we can choose suitable
$r \in (0,C_{1k}R(x_k,0)^{-\frac{1}{2}})$ and $B=M_k$, which falls
into the category (c) with the scalar curvature between
$C_{2k}^{-1}R(x,0)$ and $C_{2k}R(x,0)$ on $B$. Now by scaling the
ancient $\kappa_0$-solutions along the points $(x_k,0)$ with
factors $R(x_k,0)$, it follows from Corollary 6.4.4 that a
sequence of the ancient $\kappa_0$-solutions converge in the
$C^{\infty}_{loc}$ topology to a noncompact orientable ancient
$\kappa_0$-solution.

If the curvature operator of the noncompact limit has a nontrivial
null vector somewhere at some time, it follows exactly as before
by using the argument in the proof of Lemma 6.4.1 that the
orientable limit is isometric to the round $\mathbb{S}^2 \times
\mathbb{R}$, or the twisted product $\mathbb{S}^2 \tilde{\times}
\mathbb{R}$ where $\mathbb{Z}_2$ flips both $\mathbb{S}^2$ and
$\mathbb{R}$. Then for $k$ large enough, a suitable neighborhood
$B$ (for suitable $r$) of the point $(x_k,0)$ would fall into the
category (a) or (b) (over
$\mathbb{R}\mathbb{P}^3\setminus\bar{\mathbb{B}^3}$) with the
desired volume estimate. This is a contradiction.

If the noncompact limit has positive sectional curvature
everywhere, then by using Corollary 6.4.5 and Theorem 6.4.3(i) for
the noncompact limit we see that for $k$ large enough, a suitable
neighborhood $B$ (for suitable $r$) of the point $(x_k,0)$ would
fall into category (a) or (b) (over $\mathbb{B}^3$) with the
desired volume estimate. This is also a contradiction.

Finally, the statement on the curvature estimate in the
neighborhood $B$ follows directly from Theorem 6.4.3(i).
\end{pf}

\newpage
\part{{\Large Ricci Flow on Three-manifolds}}

\bigskip
We will use the Ricci flow to study the topology of compact
orientable three-manifolds. Let $M$ be a compact three-dimensional
orientable manifold. Arbitrarily given a Riemannian metric on the
manifold, we evolve it by the Ricci flow. The basic idea is to
understand the topology of the underlying manifold by studying
long-time behavior of the solution of the Ricci flow. We have seen
in Chapter 5 that for a compact three-manifold with positive Ricci
curvature as initial data, the solution to the Ricci flow tends,
up to scalings, to a metric of positive constant curvature.
Consequently, a compact three-manifold with positive Ricci
curvature is diffeomorphic to the round three-sphere or a metric
quotient of it.

However, for general initial metrics, the Ricci flow may develop
singularities in some parts while it keeps smooth in other parts.
Naturally one would like to cut off the singularities and continue
to run the Ricci flow. If the Ricci flow still develops
singularities after a while, one can do the surgeries and run the
Ricci flow again. By repeating this procedure, one will get a kind
of ``weak" solution to the Ricci flow. Furthermore, if the ``weak"
solution has only a finite number of surgeries at any finite time
interval and one can remember what had been cut during the
surgeries, and if the ``weak" solution has a well-understood
long-time behavior, then one will also get the topology structure
of the initial manifold. This theory of surgically modified Ricci
flow was first developed by Hamilton \cite{Ha97} for compact
four-manifolds and further developed more recently by Perelman
\cite{P2} for compact orientable three-manifolds.

The main purpose of this chapter is to give a complete and
detailed discussion of Perelman's work on the Ricci flow with
surgery on three-manifolds. More specifically, Sections 7.1-7.2
give a detailed exposition of section 12 of Perelman's first paper
\cite{P1}; Sections 7.3-7.6 give a detailed exposition of sections
2-7 of Perelman's second paper \cite{P2}, except the general
collapsing result, Theorem 7.4 of \cite{P2}, claimed by Perelman
(a special case of which has been proved by Shioya-Yamaguchi
\cite{ShY}). In Section 7.7, we combine the long time behavior
result in Section 7.6 and the collapsing result of
Shioya-Yamaguchi \cite{ShY} to present a proof of Thurston's
geometrization conjecture.

\section{Canonical Neighborhood Structures}

Let us call a Riemannian metric on a compact orientable
three-dimen\-sional manifold \textbf{normalized}\index{normalized}
if the eigenvalues of its curvature operator at every point are
bounded by $\frac{1}{10}\geq \lambda \geq \mu \geq \nu \geq
-\frac{1}{10}$, and every geodesic ball of radius one has volume
at least one. By the evolution equation of the curvature and the
maximum principle, it is easy to see that any solution to the
Ricci flow with (compact and three-dimensional) normalized initial
metric exists on a maximal time interval $[0,t_{\max})$ with
$t_{\max}>1$.

Consider a smooth solution $g_{ij}(x,t)$ to the Ricci flow on
$M\times [0,T)$, where $M$ is a compact orientable three-manifold
and $T<+\infty$. After rescaling, we may always assume the initial
metric $g_{ij}(\cdot,0)$ is normalized. By Theorem 5.3.2, the
solution $g_{ij}(\cdot,t)$ then satisfies the pinching estimate
\be
R \geq (-\nu)[\log(-\nu) + \log(1+t) -3] 
\ee whenever $\nu < 0$ on $M\times [0,T)$. Recall the function
$$
y= f(x) = x(\log x - 3), \ \mbox{ for } e^2 \leq x < +\infty,
$$
is increasing and convex with range $-e^2 \leq y < +\infty$, and
its inverse function is also increasing and satisfies
$$
\lim_{y\rightarrow +\infty} f^{-1}(y)/y = 0.
$$
We can rewrite the pinching estimate (7.1.1) as \be
Rm(x,t) \geq-[f^{-1}(R(x,t)(1+t))/(R(x,t)(1+t))]R(x,t) 
\ee on $M\times [0,T)$.

Suppose that the solution $g_{ij}(\cdot,t)$ becomes singular as
$t\rightarrow T$. Let us take a sequence of times $t_k \rightarrow
T$, and a sequence of points $p_k \in M$ such that for some
positive constant $C$, $|Rm|(x,t) \leq CQ_k$ with
$Q_k=|Rm(p_k,t_k)|$ for all $x \in M$ and $t \in [0,t_k]$. Thus,
$(p_k,t_k)$ is a sequence of (almost) maximum points. By applying
Hamilton's compactness theorem and Perelman's no local collapsing
theorem I as well as the pinching estimate (7.1.2), a sequence of
the scalings of the solution $g_{ij}(x,t)$ around the points $p_k$
with factors $Q_k$ converges to a nonflat complete
three-dimensional orientable ancient $\kappa$-solution (for some
$\kappa > 0$). For an arbitrarily given $\varepsilon>0$, the
canonical neighborhood theorem (Theorem 6.4.6) in the previous
chapter implies that each point in the ancient $\kappa$-solution
has a neighborhood which is either an evolving $\varepsilon$-neck,
or an evolving $\varepsilon$-cap, or a compact (without boundary)
positively curved manifold. This gives the structure of
singularities coming from a sequence of (almost) maximum points.

However the above argument does not work for singularities coming
from a sequence of points $(y_k,\tau_k)$ with $\tau_k \rightarrow T$
and $|Rm(y_k,\tau_k)| \rightarrow +\infty$ when $|Rm(y_k,\tau_k)|$
is not comparable with the maximum of the curvature at the time
$\tau_k$, since we cannot take a limit directly. In \cite{P1},
Perelman developed a refined rescaling argument to obtain the
following singularity structure theorem. We remark that our
statement of the singularity structure theorem below is slightly
different from Perelman's original statement (cf. Theorem 12.1 of
\cite{P1}). While Perelman assumed the condition of
$\kappa$-noncollapsing on scales less than $r_0$, we assume that the
initial metric is normalized so that from the rescaling argument one
can get the $\kappa$-noncollapsing \emph{on all scales} for the
limit solutions.

\begin{theorem}[{Singularity structure theorem}\index{singularity
structure theorem}] Given $\varepsilon>0$ and $T_0 >1$, one can
find $r_0>0$ with the following property. If $g_{ij}(x,t), x\in M$
and $t\in [0,T)$ with $1 < T \leq T_0$, is a solution to the Ricci
flow on a compact orientable three-manifold $M$ with normalized
initial metric, then for any point $(x_0,t_0)$ with $t_0\geq 1$
and $Q=R(x_0,t_0)\geq r_0^{-2}$, the solution in $\{(x,t)\ |\
d_{t_0}^2(x,x_0)<\varepsilon^{-2} Q^{-1},\ \ t_0-\varepsilon^{-2}
Q^{-1}\leq t\leq t_0\}$ is, after scaling by the factor $Q$,
$\varepsilon$-close $($in the $C^{[\varepsilon^{-1}]}$-topology$)$
to the corresponding subset of some orientable ancient
$\kappa$-solution $($for some $\kappa >0)$.
\end{theorem}

\begin{pf} The proof is basically along the line sketched by Perelman
(cf. section 12 of \cite{P1}). However, the proof of Step 2
follows, with some modifications, the argument given in the
initial notes of Kleiner-Lott \cite{KL} on Perelman's first paper
\cite{P1}. Also, we give a proof of Step 4 which is different from
both Perelman \cite{P1} and Kleiner-Lott \cite{KL}.

Since the initial metric is normalized, it follows from the no
local collapsing theorem I or I' (and their proofs) that there is
a positive constant $\kappa$, depending only on $T_0$, such that
the solution in Theorem 7.1.1 is $\kappa$-noncollapsed on all
scales less than $\sqrt{T_0}$. Let $C(\varepsilon)$ be a positive
constant larger than or equal to $\varepsilon^{-2}$. It suffices
to prove that there exists $r_0
>0$ such that for any point $(x_0,t_0)$ with $t_0 \geq 1$ and
$Q=R(x_0,t_0) \geq r_0^{-2}$, the solution in the parabolic region
$\{ (x,t) \in M \times [0,T)\ | \ d^2_{t_0}(x,x_0)<
C(\varepsilon)Q^{-1}, t_0-C(\varepsilon)Q^{-1}\leq t \leq t_0\}$
is, after scaling by the factor $Q$, $\varepsilon$-close to the
corresponding subset of some orientable ancient $\kappa$-solution.
The constant $C(\varepsilon)$ will be determined later.

We argue by contradiction. Suppose for some $\varepsilon>0$, there
exist a sequence of solutions $(M_k,g_k(\cdot,t))$ to the Ricci
flow on compact orientable three-manifolds with normalized initial
metrics, defined on the time intervals [0,$T_k)$ with $1<T_k \leq
T_0$, a sequence of positive numbers $r_k\rightarrow 0$, and a
sequence of points $x_k\in M_k$ and times $t_k\geq 1$ with
$Q_k=R_k(x_k,t_k)\geq r_k^{-2}$ such that each solution
$(M_k,g_k(\cdot,t))$ in the parabolic region $\{(x,t) \in
M_k\times [0,T_k)\ |\ d^2_{t_k}(x,x_k)<C(\varepsilon)Q_k^{-1},
t_k-C(\varepsilon)Q_k^{-1}\leq t\leq t_k\}$ is not, after scaling
by the factor $Q_k$, $\varepsilon$-close to the corresponding
subset of any orientable ancient $\kappa$-solution, where $R_k$
denotes the scalar curvature of $(M_k,g_k)$.

For each solution $(M_k,g_k(\cdot,t))$, we may adjust the point
$(x_k,t_k)$ with $t_k\geq \frac{1}{2}$ and with $Q_k=R_k(x_k,t_k)$
to be as large as possible so that the conclusion of the theorem
fails at $(x_k,t_k)$, but holds for any $(x,t)\in M_k\times
[t_k-H_kQ_k^{-1},t_k]$ satisfying $R_k(x,t)\geq 2Q_k$, where
$H_k=\frac{1}{4}r_k^{-2}\rightarrow +\infty$ as $k\rightarrow
+\infty.$  Indeed, suppose not, by setting
$(x_{k_1},t_{k_1})=(x_k,t_k)$, we can choose a sequence of points
$(x_{k_l},t_{k_l})\in M_k\times
[t_{k_{(l-1)}}-H_kR_k(x_{k_{(l-1)}},t_{k_{(l-1)}})^{-1},t_{k_{(l-1)}}]$
such that $R_k(x_{k_l},t_{k_l})\geq
2R_k(x_{k_{(l-1)}},t_{k_{(l-1)}})$ and the conclusion of the
theorem fails at $(x_{k_l},t_{k_l})$ for each $l =2, 3, \ldots.$
 Since the solution is smooth, but
$$
R_k(x_{k_l},t_{k_l})\geq 2R_k(x_{k_{(l-1)}},t_{k_{(l-1)}})\geq
\cdots \geq 2^{l-1}R_k(x_{k},t_{k}),
$$
and
\begin{align*}
t_{k_l}&  \geq
t_{k_{(l-1)}}-H_kR_k(x_{k_{(l-1)}},t_{k_{(l-1)}})^{-1}\\
&  \geq
t_{k}-H_k\sum\limits_{i=1}^{l-1}\frac{1}{2^{i-1}}R_k(x_{k},t_{k})^{-1}\\
&  \geq  \frac{1}{2},
\end{align*}
this process must terminate after a finite number of steps and the
last element fits.

Let\; $(M_k,\,\tilde{g}_k(\cdot,t),\,x_k)$\, be\, the\, rescaled\,
solutions\, obtained\, by\, rescaling $(M_k,g_k(\cdot,t))$ around
$x_k$ with the factors $Q_k=R_k(x_k,t_k)$ and shifting the time
$t_k$ to the new time zero. Denote by $\tilde{R}_k$ the rescaled
scalar curvature.  We will show that a subsequence of the
orientable rescaled solutions $(M_k,\tilde{g}_k(\cdot,t),x_k)$
converges in the $C^{\infty}_{loc}$ topology to an orientable
ancient $\kappa$-solution, which is a contradiction. In the
following we divide the argument into four steps.

\medskip
{\it Step} 1. \ First of all, we need a local bound on curvatures.
The following lemma is the Claim 1 of Perelman in his proof Theorem
12.1 of \cite{P1}.

\begin{lemma}
For each $(\bar{x},\bar{t})$ with $t_k-\frac{1}{2}H_kQ_k^{-1}\leq
\bar{t}\leq t_k$, we have $R_k(x,t)\leq 4\bar{Q}_k$ whenever
$\bar{t}-c\bar{Q}_k^{-1}\leq t\leq \bar{t}$ and
$d^2_{\bar{t}}(x,\bar{x})\leq c\bar{Q}_k^{-1}$, where
$\bar{Q}_k=Q_k+R_k(\bar{x},\bar{t})$ and $c>0$ is a small
universal constant.
\end{lemma}

\begin{pf}
Consider any point $(x,t)\in B_{\bar{t}}(\bar{x},
(c\bar{Q}_k^{-1})^{\frac{1}{2}})\times
[\bar{t}-c\bar{Q}_k^{-1},\bar{t}]$ with $c>0$ to be determined. If
$R_k(x,t)\leq 2Q_k$, there is nothing to show. If $R_k(x,t)>
2Q_k$, consider a space-time curve $\gamma$ from $(x,t)$ to
$(\bar{x},\bar{t})$ that goes straight from $(x,t)$ to
$(x,\bar{t})$ and goes from $(x,\bar{t})$ to $(\bar{x},\bar{t})$
along a minimizing geodesic (with respect to the metric
$g_k(\cdot,\bar{t}))$. If there is a point on $\gamma$ with the
scalar curvature $2Q_k$, let $y_0$ be the nearest such point to
$(x,t)$. If not, put $y_0=(\bar{x},\bar{t}).$ On the segment of
$\gamma$ from $(x,t)$ to $y_0$, the scalar curvature is at least
$2Q_k$. According to the choice of the point $(x_k,t_k)$, the
solution along the segment is $\varepsilon$-close to that of some
ancient $\kappa$-solution. It follows from Theorem 6.4.3 (ii) that
$$
|\nabla(R_k^{-\frac{1}{2}})|\leq 2\eta \; \mbox{ and } \;
\left|\frac{\partial}{\partial t}(R_k^{-1})\right|\leq 2\eta
$$
on the segment. (Here, without loss of generality, we may assume
$\varepsilon$ is suitably small). Then by choosing $c>0$
(depending only on $\eta$) small enough we get the desired
curvature bound by integrating the above derivative estimates
along the segment. This proves the lemma.
\end{pf}

\smallskip
{\it Step} 2. \  Next we want to show that for each $A<+\infty$,
there exist a positive constant $C(A)$ (independent of $k$) such
that the curvatures of the rescaled solutions
$\tilde{g}_k(\cdot,t)$ at the new time $t =0$ (corresponding to
the original times $t_k$) satisfy the estimate
$$
|\widetilde{Rm}_{k}|(y,0) \leq C(A)
$$
whenever $d_{\tilde{g}_k(\cdot,0)}(y,x_k) \leq A$ and $k \geq 1$.
(This is just a weak version of the Claim 2 of Perelman in his
proof of Theorem 12.1 of \cite{P1}. The first detailed exposition
on this weak version was given by Kleiner-Lott in their June 2003
notes \cite{KL}. We now follow their argument as in \cite{KL} here
with some modifications.)

For each $\rho \geq 0$, set
$$
M(\rho)=\sup\{\tilde{R}_k(x,0)\ |\ k\geq 1, x\in M_k \; \mbox{
with } \;d_{0}(x,x_k)\leq \rho\}
$$
and
$$
\rho_0=\sup\{\rho \geq 0\  |\ \ M(\rho)<+\infty \}.
$$
By the pinching estimate (7.1.1), it suffices to show
$\rho_0=+\infty.$

Note that $\rho_0>0$ by applying Lemma 7.1.2 with
$(\bar{x},\bar{t})=(x_k,t_k).$ We now argue by contradiction to
show $\rho_0=+\infty.$ Suppose not, we may find (after passing to
a subsequence if necessary) a sequence of points $y_k\in M_k$ with
$d_{0}(x_k,y_k)\rightarrow \rho_0<+\infty$ and
$\tilde{R}_k(y_k,0)\rightarrow +\infty$. Let $\gamma_k(\subset
M_k)$ be a minimizing geodesic segment from $x_k$ to $y_k$. Let
$z_k\in \gamma_k$ be the point on $\gamma_k$ closest to $y_k$ with
$\tilde{R}_k(z_k,0)=2$, and let $\beta_k$ be the subsegment of
$\gamma_k$ running from $z_k$ to $y_k$. By Lemma 7.1.2 the length
of $\beta_k$ is bounded away from zero independent of $k$. By the
pinching estimate (7.1.1), for each $\rho<\rho_0$, we have a
uniform bound on the curvatures on the open balls
$B_{0}(x_k,\rho)\subset (M_k,\tilde{g}_k)$. The injectivity radii
of the rescaled solutions $\tilde{g}_k$ at the points $x_k$ and
the time $t=0$ are also uniformly bounded from below by the
$\kappa$-noncollapsing property. Therefore by Lemma 7.1.2 and
Hamilton's compactness theorem (Theorem 4.1.5), after passing to a
subsequence, we can assume that the marked sequence
$(B_{0}(x_k,\rho_0),\tilde{g}_k(\cdot,0),x_k)$ converges in the
$C_{\rm loc}^{\infty}$ topology to a marked (noncomplete) manifold
$(B_{\infty},\tilde{g}_{\infty},x_{\infty})$, the segments
$\gamma_k$ converge to a geodesic segment (missing an endpoint)
$\gamma_{\infty}\subset B_{\infty}$ emanating from $x_{\infty}$,
and $\beta_k$ converges to a subsegment $\beta_{\infty}$ of
$\gamma_{\infty}$. Let $\bar{B}_{\infty}$ denote the completion of
$(B_{\infty},\tilde{g}_{\infty})$, and $y_{\infty}\in
\bar{B}_{\infty}$ the limit point of $\gamma_{\infty}$.

Denote by $\tilde{R}_{\infty}$ the scalar curvature of
$(B_{\infty},\tilde{g}_{\infty})$. Since the rescaled scalar
curvatures $\tilde{R}_k$ along $\beta_k$ are at least 2, it
follows from the choice of the points $(x_k,0)$ that for any
$q_0\in \beta_{\infty}$, the manifold
$(B_{\infty},\tilde{g}_{\infty})$ in $\{q\in B_{\infty}| \ \ {\rm
dist}^2_{\tilde{g}_{\infty}}(q,q_0)<
C(\varepsilon)(\tilde{R}_{\infty}(q_0))^{-1}\}$ is
$2\varepsilon$-close to the corresponding subset of (a time slice
of) some orientable ancient $\kappa$-solution. Then by Theorem
6.4.6, we know that the orientable ancient $\kappa$-solution at
each point $(x,t)$ has a radius $r$, $0 < r <
C_1(2\varepsilon)R(x,t)^{-\frac{1}{2}}$, such that its canonical
neighborhood $B$, with $B_t(x,r) \subset  B \subset B_t(x,2r)$, is
either an evolving $2\varepsilon$-neck, or an evolving
$2\varepsilon$-cap, or a compact manifold (without boundary)
diffeomorphic to a metric quotient of the round three-sphere
$\mathbb{S}^3$, and moreover the scalar curvature is between
$(C_2(2\varepsilon))^{-1}R(x,t)$ and $C_2(2\varepsilon)R(x,t)$,
where $C_1(2\varepsilon)$ and $C_2(2\varepsilon)$ are the positive
constants in Theorem 6.4.6.

We now choose $C(\varepsilon) = \max \{
2C_1^2(2\varepsilon),\varepsilon^{-2}\}$. By the local curvature
estimate in Lemma 7.1.2, we see that the scalar curvature
$\tilde{R}_{\infty}$ becomes unbounded when approaching
$y_{\infty}$ along $\gamma_{\infty}$. This implies that the
canonical neighborhood around $q_0$ cannot be a compact manifold
(without boundary) diffeomorphic to a metric quotient of the round
three-sphere $\mathbb{S}^3$. Note that $\gamma_{\infty}$ is
shortest since it is the limit of a sequence of shortest
geodesics. Without loss of generality, we may assume $\varepsilon$
is suitably small (say, $\varepsilon \leq \frac{1}{100}$). These
imply that as $q_0$ gets sufficiently close to $y_{\infty}$, the
canonical neighborhood around $q_0$ cannot be an evolving
$2\varepsilon$-cap. Thus we conclude that each $q_0 \in
\gamma_{\infty}$ sufficiently close to $y_{\infty}$ is the center
of an evolving $2\varepsilon$-neck.

Let
$$
U=\bigcup\limits_{q_0\in
\gamma_{\infty}}B(q_0,4\pi(\tilde{R}_{\infty}(q_0))^{-\frac{1}{2}})\
\ (\subset (B_{\infty}, \tilde{g}_{\infty})),
$$
where $B(q_0,4\pi(\tilde{R}_{\infty}(q_0))^{-\frac{1}{2}})$ is the
ball centered at $q_0\in B_{\infty}$ with the radius
$4\pi(\tilde{R}_{\infty}(q_0))^{-\frac{1}{2}}$. Clearly $U$ has
nonnegative sectional curvature by the pinching estimate (7.1.1).
Since the metric $\tilde{g}_{\infty}$ is cylindrical at any point
$q_0\in \gamma_{\infty}$ which is sufficiently close to
$y_{\infty}$, we see that the metric space $\overline{U}=U\cup
\{y_{\infty}\}$ by adding in the point $y_{\infty}$, is locally
complete and intrinsic near $y_{\infty}$. Furthermore $y_{\infty}$
cannot be an interior point of any geodesic segment in
$\overline{U}$. This implies the curvature of $\overline{U}$ at
$y_{\infty}$ is nonnegative in the Alexandrov sense. It is a basic
result in Alexandrov space theory (see for example Theorem 10.9.3
and Corollary 10.9.5 of \cite{BBI}) that there exists a
three-dimensional tangent cone $C_{y_{\infty}}\overline{U}$ at
$y_{\infty}$ which is a metric cone. It is clear that its aperture
is $\leq 10\varepsilon$, thus the tangent cone is nonflat.

Pick a point $p\in C_{y_{\infty}}\overline{U}$ such that the
distance from the vertex $y_{\infty}$ to $p$ is one and it is
nonflat around $p$. Then the ball $B(p,\frac{1}{2})\subset
C_{y_{\infty}}\overline{U}$ is the Gromov-Hausdorff limit of the
scalings of a sequence of balls
$B_0(p_k,s_k)\subset(M_k,\tilde{g}_{k}(\cdot,0))$ by some factors
$a_k$, where $s_k\rightarrow 0^+$. Since the tangent cone is
three-dimensional and nonflat around $p$, the factors $a_k$ must
be comparable with $\tilde{R}_k(p_k,0)$. By using the local
curvature estimate in Lemma 7.1.2, we actually have the
convergence in the $C^{\infty}_{\rm loc}$ topology for the
solutions $\tilde{g}_{k}(\cdot,t)$ on the balls $B_0(p_k,s_k)$ and
over some time interval $t\in [-\delta, 0]$ for some sufficiently
small $\delta>0$. The limiting ball $B(p,\frac{1}{2})\subset
C_{y_{\infty}}\overline{U}$ is a piece of the nonnegative curved
and nonflat metric cone whose radial directions are all Ricci
flat. On the other hand, by applying Hamilton's strong maximum
principle to the evolution equation of the Ricci curvature tensor
as in the proof of Lemma 6.3.1, the limiting ball
$B(p,\frac{1}{2})$ would split off all radial directions
isometrically (and locally). Since the limit is nonflat around
$p$, this is impossible. Therefore we have proved that the
curvatures of the rescaled solutions $\tilde{g}_k(\cdot,t)$ at the
new times $t=0$ (corresponding to the original times $t_k$) stay
uniformly bounded at bounded distances from $x_k$ for all $k$.

We have proved that for each $A<+\infty$, the curvature of the
marked manifold $(M_k,\tilde{g}_k(\cdot,0),x_k)$ at each point
$y\in M_k$ with distance from $x_k$ at most $A$ is bounded by
$C(A)$. Lemma 7.1.2 extends this curvature control to a backward
parabolic neighborhood centered at $y$ whose radius depends only
on the distance from $y$ to $x_k$. Thus by Shi's local derivative
estimates (Theorem 1.4.2) we can control all derivatives of the
curvature in such backward parabolic neighborhoods. Then by using
the $\kappa$-noncollapsing and Hamilton's compactness theorem
(Theorem 4.1.5), we can take a $C_{\rm loc}^{\infty}$ subsequent
limit to obtain $(M_{\infty},
\tilde{g}_{\infty}(\cdot,t),x_{\infty})$, which is
$\kappa$-noncollapsed on all scales and is defined on a space-time
open subset of $M_{\infty}\times(-\infty,0]$ containing the time
slice $M_{\infty}\times\{0\}$. Clearly it follows from the
pinching estimate (7.1.1) that the limit $(M_{\infty},
\tilde{g}_{\infty}(\cdot,0),x_{\infty})$ has nonnegative curvature
operator (and hence nonnegative sectional curvature).

\medskip
{\it Step} 3. \ We further claim that the limit $(M_{\infty},
\tilde{g}_{\infty}(\cdot,0),x_{\infty})$ at the time slice $\{ t=0
\}$ has bounded curvature.

We know that the sectional curvature of the limit $(M_{\infty},
\tilde{g}_{\infty}(\cdot,0),x_{\infty})$ is nonnegative
everywhere. Argue by contradiction. Suppose the curvature of
$(M_{\infty}, \tilde{g}_{\infty}(\cdot,0),x_{\infty})$ is not
bounded, then by Lemma 6.1.4, there exists a sequence of points
$q_j \in M_{\infty}$ diverging to infinity such that their scalar
curvatures $\tilde{R}_{\infty}(q_j,0) \rightarrow +\infty$ as $j
\rightarrow +\infty$ and
$$
\tilde{R}_{\infty}(x,0) \leq 4\tilde{R}_{\infty}(q_j,0)
$$
for $x \in B(q_j, j/\sqrt{\tilde{R}_{\infty}(q_j,0)})\subset
(M_{\infty},\tilde{g}_{\infty}(\cdot,0))$. By combining with Lemma
7.1.2 and the $\kappa$-noncollapsing, a subsequence of the
rescaled and marked manifolds
$(M_{\infty},\tilde{R}_{\infty}(q_j,0)\tilde{g}_{\infty}(\cdot,0),q_j)$
converges in the $C^{\infty}_{\rm loc}$ topology to a smooth
nonflat limit $Y$. By Proposition 6.1.2, the new limit $Y$ is
isometric to a metric product $N \times \mathbb{R}$ for some
two-dimensional manifold $N$. On the other hand, in view of the
choice of the points $(x_k,t_k)$, the original limit $(M_{\infty},
\tilde{g}_{\infty}(\cdot,0),x_{\infty})$ at the point $q_j$ has a
canonical neighborhood which is either a $2\varepsilon$-neck, a
$2\varepsilon$-cap, or a compact manifold (without boundary)
diffeomorphic to a metric quotient of the round $\mathbb{S}^3$. It
follows that for $j$ large enough, $q_j$ is the center of a
$2\varepsilon$-neck of radius
$(\tilde{R}_{\infty}(q_j,0))^{-\frac{1}{2}}$. Without loss of
generality, we may further assume that $2\varepsilon <
\varepsilon_0$, where $\varepsilon_0$ is the positive constant
given in Proposition 6.1.1. Since
$(\tilde{R}_{\infty}(q_j,0))^{-\frac{1}{2}}\rightarrow 0$ as
$j\rightarrow +\infty$, this contradicts Proposition 6.1.1. So the
curvature of $(M_{\infty}, \tilde{g}_{\infty}(\cdot,0))$ is
bounded.

\medskip
{\it Step} 4. \  Finally we want to extend the limit backwards in
time to $-\infty$.

By Lemma 7.1.2 again, we now know that the limiting solution
$(M_{\infty},$ $\tilde{g}_{\infty}(\cdot,t))$ is defined on a
backward time interval $[-a,0]$ for some $a>0$.

Denote by
\begin{align*}
t' = \inf \{ \tilde{t}&  | \mbox{ we can take a
smooth limit on }\; (\tilde{t},0]  \mbox { (with bounded} \\
& \mbox {curvature at each time slice) from a subsequence}\\
&  \mbox {of the convergent rescaled solutions }\; \tilde{g}_k\}.
\end{align*}
We first claim that there is a subsequence of the rescaled
solutions $\tilde{g}_k$ which converges in the $C^{\infty}_{\rm
loc}$ topology to a smooth limit
$(M_{\infty},\tilde{g}_{\infty}(\cdot,t))$ on the maximal time
interval $(t',0]$.

Indeed, let $t'_k$ be a sequence of negative numbers such that
$t'_k\rightarrow t'$ and there exist smooth limits
$(M_{\infty},\tilde{g}_{\infty}^{k}(\cdot,t))$ defined on
$(t'_k,0]$. For each $k$, the limit has nonnegative sectional
curvature and has bounded curvature at each time slice. Moreover
by Lemma 7.1.2, the limit has bounded curvature on each
subinterval $[-b,0]\subset(t'_k,0]$. Denote by $\tilde{Q}$ the
scalar curvature upper bound of the limit at time zero (where
$\tilde{Q}$ is the same for all $k$). Then we can apply the
Li-Yau-Hamilton estimate (Corollary 2.5.7) to get
$$
\tilde{R}_{\infty}^{k}(x,t)\leq\tilde{Q} \(\frac{-t'_k}{t-t'_k}\),
$$
where $\tilde{R}_{\infty}^{k}(x,t)$ are the scalar curvatures of
the limits $(M_{\infty},\tilde{g}_{\infty}^{k}(\cdot,t))$.  Hence
by the definition of convergence and the above curvature
estimates, we can find a subsequence of the above convergent
rescaled solutions $\tilde{g}_k$ which converges in the
$C^{\infty}_{\rm loc}$ topology to a smooth limit
$(M_{\infty},\tilde{g}_{\infty}(\cdot,t))$ on the maximal time
interval $(t',0]$.

We next claim that $t'=-\infty$.

Suppose not, then by Lemma 7.1.2, the curvature of the limit
$(M_{\infty},$ $\tilde{g}_{\infty}(\cdot,t))$ becomes unbounded as
$t\rightarrow t'>-\infty$. By applying the maximum principle to
the evolution equation of the scalar curvature, we see that the
infimum of the scalar curvature is nondecreasing in time. Note
that $\tilde{R}_{\infty}(x_{\infty},0) =1$. Thus there exists some
point $y_{\infty}\in M_{\infty}$ such that
$$
\tilde{R}_{\infty}\(y_{\infty},t'+\frac{c}{10}\)<\frac{3}{2}
$$
where $c>0$ is the universal constant in Lemma 7.1.2. By using
Lemma 7.1.2 again we see that the limit
$(M_{\infty},\tilde{g}_{\infty}(\cdot,t))$ in a small neighborhood
of the point $(y_{\infty},t'+\frac{c}{10})$ extends backwards to
the time interval $[t'-\frac{c}{10},t'+\frac{c}{10}]$. We remark
that the distances at time $t$ and time $0$ are roughly equivalent
in the following sense \be
d_t(x,y) \geq d_0(x,y) \geq d_t(x,y)-{\rm const.}  
\ee for any $x,y \in M_{\infty}$ and $t \in (t',0]$. Indeed from
the Li-Yau-Hamilton inequality (Corollary 2.5.7) we have the
estimate
$$
\tilde{R}_\infty(x,t) \le \tilde{Q}\(\frac{-t'}{t-t'}\),\quad
\mbox{on }\; M_\infty\times(t',0].
$$
By applying Lemma 3.4.1 (ii), we have
$$
d_t(x,y)\le d_0(x,y)+30(-t')\sqrt{\tilde{Q}}
$$
for any $x,y\in M_\infty$ and $t\in(t',0]$. On the other hand,
since the curvature of the limit metric
$\tilde{g}_\infty(\cdot,t)$ is nonnegative, we have
$$
d_t(x,y)\ge d_0(x,y)
$$
for any $x,y\in M_\infty$ and $t\in(t',0]$. Thus we obtain the
estimate (7.1.3).

Let us still denote by $(M_k,\tilde{g}_k(\cdot,t))$ the
subsequence which converges on the maximal time interval
$(t'\!,0]$. Consider the rescaled sequence
$(\!M_k,\tilde{g}_k(\cdot,t))$ with the marked points $x_k$
replaced by the associated sequence of points $y_k \rightarrow
y_{\infty}$ and the (original unshifted) times $t_k$ replaced by
any $s_k \in
[t_k+(t'-\frac{c}{20})Q_k^{-1},t_k+(t'+\frac{c}{20})Q_k^{-1}]$. It
follows from Lemma 7.1.2 that for $k$ large enough, the rescaled
solutions $(M_k,\tilde{g}_k(\cdot,t))$ at $y_k$ satisfy
$$
\tilde{R}_k(y_k,t) \leq 10
$$
for all $t \in [t'-\frac{c}{10},t'+\frac{c}{10}]$. By applying the
same arguments as in the above Step 2, we conclude that for any $A
> 0$, there is a positive constant $C(A)<+\infty$ such that
$$
\tilde{R}_k(x,t) \leq C(A)
$$
for all $(x,t)$ with ${d}_t(x,y_k) \leq A$ and $t \in
[t'-\frac{c}{20},t'+\frac{c}{20}]$. The estimate (7.1.3) implies
that there is a positive constant $A_0$ such that for arbitrarily
given small $\epsilon' \in(0,\frac{c}{100})$, for $k$ large
enough, there hold
$$
{d}_t(x_k,y_k) \leq A_0
$$
for all $t \in [t'+ \epsilon',0]$. By combining with Lemma 7.1.2,
we then conclude that for any $A>0$, there is a positive constant
$\tilde{C}(A)$ such that for $k$ large enough, the rescaled
solutions $(M_k,\tilde{g}_k(\cdot,t))$ satisfy
$$
\tilde{R}_k(x,t) \leq \tilde{C}(A)
$$
for all $x \in \tilde{B}_0(x_k,A)$ and $t \in [t' -
\frac{c}{100}(C(A))^{-1},0]$.

Now, by taking convergent subsequences from the (original)
rescaled solutions $(M_k,\tilde{g}_k(\cdot,t),x_k)$, we see that
the limiting solution $(M_{\infty},\tilde{g}_{\infty}(\cdot,t))$
is defined on a space-time open subset of $M_{\infty} \times
(-\infty,0]$ containing $M_{\infty} \times [t',0]$. By repeating
the argument of Step 3 and using Lemma 7.1.2, we further conclude
the limit $(M_{\infty},\tilde{g}_{\infty}(\cdot,t))$ has uniformly
bounded curvature on $M_{\infty} \times [t',0]$. This is a
contradiction.

Therefore we have proved a subsequence of the rescaled solutions
$(M_k,$ $\tilde{g}_k(\cdot,t),x_k)$ converges to an orientable
ancient $\kappa$-solution, which gives the desired contradiction.
This completes the proof of the theorem.
\end{pf}

We remark that this singularity structure theorem has been
extended by Chen and the second author in \cite{CZ05F} to the
Ricci flow on compact four-manifolds with positive isotropic
curvature.

\section{Curvature Estimates for Smooth Solutions}

Let us consider solutions to the Ricci flow on compact orientable
three-manifolds with normalized initial metrics. The above
singularity structure theorem of Perelman (Theorem 7.1.1) tells us
that the solutions around high curvature points are sufficiently
close to ancient $\kappa$-solutions. It is thus reasonable to
expect that the elliptic type estimate (Theorem 6.4.3) and the
curvature estimate via volume growth (Theorem 6.3.3) for ancient
$\kappa$-solutions are heritable to general solutions of the Ricci
flow on three-manifolds. The main purpose of this section is to
establish such curvature estimates. In the fifth section of this
chapter, we will further extend these estimates to surgically
modified solutions.

The first result of this section is an extension of the elliptic
type estimate (Theorem 6.4.3) by Perelman (cf. Theorem 12.2 of
\cite{P1}). This result is reminiscent of the second step in the
proof of Theorem 7.1.1.

\begin{theorem}[{Perelman \cite{P1}}]
For any $A<+\infty$, there exist $K=K(A)<+\infty$ and
$\alpha=\alpha(A)>0$ with the following property. Suppose we have
a solution to the Ricci flow on a three-dimensional, compact and
orientable manifold $M$ with normalized initial metric. Suppose
that for some $x_0\in M$ and some $r_0>0$ with $r_0<\alpha$, the
solution is defined for $0\le t\le r_0^2$ and satisfies
$$
|Rm|(x,t)\le r_0^{-2},\quad \text{for }\; 0\le t\le
r_0^2,\;d_0(x,x_0)\le r_0,
$$
and
$$
\Vol_0(B_0(x_0,r_0))\ge A^{-1}r_0^3.
$$
Then $R(x,r_0^2)\le Kr_0^{-2}$ whenever $d_{r_0^2}(x,x_0)<Ar_0$.
\end{theorem}

\begin{pf}
Given any large $A>0$ and letting $\alpha>0$ be chosen later,  by
Perelman's no local collapsing theorem II (Theorem 3.4.2), there
exists a positive constant $\kappa=\kappa(A)$ (independent of
$\alpha$) such that any complete solution satisfying the
assumptions of the theorem is $\kappa$-noncollapsed on scales $\le
r_0$ over the region $\{(x,t)\ |\ \frac{1}{5}r_0^2\le t\le
r_0^2,\;d_t(x,x_0)\le 5Ar_0\}.$ Set
$$
\varepsilon = \min \left\{ \frac{1}{4}\varepsilon_0,
\frac{1}{100}\right\},
$$
where $\varepsilon_0$ is the positive constant in Proposition
6.1.1. We first prove the following assertion.

\medskip
{\bf Claim.} \ For the above fixed $\varepsilon>0$, one can find
$K=K(A,\varepsilon)<+\infty$ such that if we have a
three-dimensional complete orientable solution with normalized
initial metric and satisfying
$$
|Rm|(x,t)\le r_0^{-2}\quad \mbox{for }\; 0\le t\le
r_0^2,\;d_0(x,x_0)\le r_0,
$$
and
$$
\Vol_0(B_0(x_0,r_0))\ge A^{-1}r_0^3
$$
for some $x_0 \in M$ and some $r_0 > 0$, then for any point $x\in
M$ with $d_{r_0^2}(x,x_0)<3Ar_0$, either $$R(x,r_0^2)< Kr_0^{-2}$$
or the subset $\{(y,t)\ |\ d_{r_0^2}^2(y,x)\leq \varepsilon^{-2}
R(x,r_0^2)^{-1},\;r_0^2-\varepsilon^{-2} R(x,r_0^2)^{-1}\le t\le
r_0^2\}$ around the point $(x,r_0^2)$ is $\varepsilon$-close to
the corresponding subset of an orientable ancient
$\kappa$-solution.
\medskip

Notice that in this assertion we don't impose the restriction of
$r_0<\alpha$, so we can consider for the moment $r_0>0$ to be
arbitrary in proving the above claim. Note that the assumption on
the normalization of the initial metric is just to ensure the
pinching estimate. By scaling, we may assume $r_0=1$. The proof of
the claim is essentially adapted from that of Theorem 7.1.1. But
we will meet the difficulties of adjusting points and verifying a
local curvature estimate.

Suppose that the claim is not true. Then there exist a sequence of
solutions $(M_k,g_k(\cdot,t))$ to the Ricci flow satisfying the
assumptions of the claim with the origins $x_{0_k}$, and a
sequence of positive numbers $K_k\rightarrow\infty$, times $t_k=1$
and points $x_k\in M_k$ with $d_{t_k}(x_k,x_{0_k})<3A$ such that
$Q_k=R_k(x_k,t_k)\ge K_k$ and the solution in $\{(x,t)\ |\
t_k-C(\varepsilon)Q_k^{-1}\le t\le t_k,\;d_{t_k}^2(x,x_k)\le
C(\varepsilon)Q_k^{-1}\}$ is not, after scaling by the factor
$Q_k$, $\varepsilon$-close to the corresponding subset of any
orientable ancient $\kappa$-solution, where $R_k$ denotes the
scalar curvature of $(M_k,g_k(\cdot,t))$ and $C(\varepsilon)(\geq
\varepsilon^{-2})$ is the constant defined in the proof of Theorem
7.1.1. As before we need to first adjust the point $(x_k,t_k)$
with $t_k\ge\frac{1}{2}$ and $d_{t_k}(x_k,x_{0_k})<4A$ so that
$Q_k=R_k(x_k,t_k)\ge K_k$ and the conclusion of the claim fails at
$(x_k,t_k)$, but holds for any $(x,t)$ satisfying $R_k(x,t)\ge
2Q_k,$ $t_k-H_kQ_k^{-1}\le t\le t_k$ and
$d_t(x,x_{0_k})<d_{t_k}(x_k,x_{0_k})
+H_k^{\frac{1}{2}}Q_k^{-\frac{1}{2}}$, where
$H_k=\frac{1}{4}K_k\rightarrow\infty,$ as $k\rightarrow+\infty$.

Indeed, by starting with $(x_{k_1},t_{k_1})=(x_k,1)$ we can choose
$(x_{k_2},t_{k_2})\in M_k\times (0,1]$ with
$t_{k_1}-H_kR_k(x_{k_1},t_{k_1})^{-1}\le t_{k_2}\le t_{k_1}$, and
$d_{t_{k_2}}(x_{k_2},x_{0_k})<d_{t_{k_1}}(x_{k_1},x_{0_k})
+H_k^\frac{1}{2}R_k(x_{k_1},t_{k_1})^{-\frac{1}{2}}$ such that
$R_k(x_{k_2},t_{k_2})\ge 2R_k(x_{k_1},t_{k_1})$ and the conclusion
of the claim fails at $(x_{k_2},t_{k_2})$; otherwise we have the
desired point. Repeating this process, we can choose points
$(x_{k_i},t_{k_i})$, $i=2, \ldots, j$, such that
$$
R_k(x_{k_i},t_{k_i})\ge 2R_k(x_{k_{i-1}},t_{k_{i-1}}),
$$
$$
t_{k_{i-1}}-H_kR_k(x_{k_{i-1}},t_{k_{i-1}})^{-1} \le t_{k_{i}}\le
t_{k_{i-1}},
$$
$$
d_{t_{k_{i}}}(x_{k_{i}},x_{0_{k}})
<d_{t_{k_{i-1}}}(x_{k_{i-1}},x_{0_k})+H_k^\frac{1}{2}
R_k(x_{k_{i-1}},t_{k_{i-1}})^{-\frac{1}{2}},
$$
and the conclusion of the claim fails at the points
$(x_{k_i},t_{k_i})$, $i=2, \ldots, j$. These inequalities imply
$$
R_k(x_{k_{j}},t_{k_{j}})\ge 2^{j-1}R_k(x_{k_1},t_{k_1})\ge
2^{j-1}K_k,
$$
$$
1\ge t_{k_{j}}\ge t_{k_{1}}
-H_k\sum_{i=0}^{j-2}\frac{1}{2^i}R_k(x_{k_1},t_{k_1})^{-1}
\ge\frac{1}{2},
$$
and
$$
d_{t_{k_{j}}}(x_{k_j},x_{0_k}) <d_{t_{k_1}}(x_{k_1},x_{0_k})
+H_k^\frac{1}{2}\sum_{i=0}^{j-2}\frac{1}{(\sqrt{2})^i}
R_k(x_{k_1},t_{k_1})^{-\frac{1}{2}}<4A.
$$
Since the solutions are smooth, this process must terminate after
a finite number of steps to give the desired point, still denoted
by $(x_k,t_k).$

For each adjusted $(x_k,t_k)$, let $[t',t_k]$ be the maximal
subinterval of $[t_k-\frac{1}{2}\varepsilon^{-2}Q_k^{-1},t_k]$ so
that the conclusion of the claim with $K=2Q_k$ holds on
\begin{align*}
& P\(x_k,t_k,\frac{1}{10}H_k^{\frac{1}{2}}Q_k^{-\frac{1}{2}},t'-t_k\)\\
&=\left\{ ({x},{t}) \ |\ {x} \in
B_{{t}}\(x_k,\frac{1}{10}H_k^{\frac{1}{2}}Q_k^{-\frac{1}{2}}\),
{t} \in [t',t_k]\right\}
\end{align*}
for all sufficiently large $k$. We now want to show $t' =
t_k-\frac{1}{2}\varepsilon^{-2}Q_k^{-1}$.

Consider the scalar curvature $R_k$ at the point $x_k$ over the
time interval $[t',t_k]$. If there is a time $\tilde{t} \in
[t',t_k]$ satisfying $R_k(x_k,\tilde{t}) \geq 2Q_k$, we let
$\tilde{t}$ be the first of such time from $t_k$. Then the
solution $(M_k,g_k(\cdot,t))$ around the point $x_k$ over the time
interval $[\tilde{t}-
\frac{1}{2}\varepsilon^{-2}Q_k^{-1},{\tilde{t}}]$ is
$\varepsilon$-close to some orientable ancient $\kappa$-solution.
Note from the Li-Yau-Hamilton inequality that the scalar curvature
of any ancient $\kappa$-solution is pointwise nondecreasing in
time. Consequently, we have the following curvature estimate
$$
R_k(x_k,t) \leq 2(1+\varepsilon)Q_k
$$
for $ t \in [\tilde{t}- \frac{1}{2}\varepsilon^{-2}Q_k^{-1},t_k]$
(or $t \in [t',t_k]$ if there is no such time $\tilde{t}$).  By
combining with the elliptic type estimate for ancient
$\kappa$-solutions (Theorem 6.4.3) and the Hamilton-Ivey pinching
estimate, we further have \be
|Rm(x,t)| \leq 5\omega(1)Q_k  
\ee for all $x \in B_t(x_k,(3Q_k)^{-\frac{1}{2}})$ and $t \in
[\tilde{t}- \frac{1}{2}\varepsilon^{-2}Q_k^{-1},t_k]$ (or $t \in
[t',t_k]$) and all sufficiently large $k$, where $\omega$ is the
positive function in Theorem 6.4.3.

For any point $(x,t)$ with $ \tilde{t}-
\frac{1}{2}\varepsilon^{-2}Q_k^{-1} \leq t \leq t_k$ (or $t \in
[t',t_k]$) and $d_t(x,x_k) \leq
\frac{1}{10}H_k^{\frac{1}{2}}Q_k^{-\frac{1}{2}}$, we divide the
discussion into two cases.

\medskip
{\it Case} (1): $d_t(x_k,x_{0_k}) \leq
\frac{3}{10}H_k^{\frac{1}{2}}Q_k^{-\frac{1}{2}}.$
\begin{align}
d_t(x,x_{0_k})&  \leq
d_t(x,x_k) + d_t(x_k,x_{0_k})\\
&  \leq\frac{1}{10}H_k^{\frac{1}{2}}Q_k^{-\frac{1}{2}}
+ \frac{3}{10}H_k^{\frac{1}{2}}Q_k^{-\frac{1}{2}}\nn\\
&  \leq \frac{1}{2}H_k^{\frac{1}{2}}Q_k^{-\frac{1}{2}}.\nn
\end{align} 

\medskip
{\it Case} (2): $d_t(x_k,x_{0_k}) >
\frac{3}{10}H_k^{\frac{1}{2}}Q_k^{-\frac{1}{2}}.$

\smallskip
{}From the curvature bound $(7.2.1)$ and the assumption, we apply
Lemma 3.4.1(ii) with $r_0 = Q_k^{-\frac{1}{2}}$ to get
$$
\frac{d}{dt}(d_t(x_k,x_{0_k})) \geq -20(\omega(1) +
1)Q_k^{\frac{1}{2}},
$$
and then for $k$ large enough,
\begin{align*}
d_t(x_k,x_{0_k})&  \leq d_{\hat{t}}(x_k,x_{0_k})
+ 20(\omega(1)+1)\varepsilon^{-2}Q_k^{-\frac{1}{2}}\\
&  \leq d_{\hat{t}}(x_k,x_{0_k}) +
\frac{1}{10}H_k^{\frac{1}{2}}Q_k^{-\frac{1}{2}},
\end{align*}
where $\hat{t} \in (t,t_k]$ satisfies the property that
$d_s(x_k,x_{0_k}) \geq
\frac{3}{10}H_k^{\frac{1}{2}}Q_k^{-\frac{1}{2}}$ whenever $s \in
[t,\hat{t}]$. So we have
\begin{align}
d_t(x,x_{0_k})&  \leq d_t(x,x_k) + d_t(x_k,x_{0_k})\\
&  \leq\frac{1}{10}H_k^{\frac{1}{2}}Q_k^{-\frac{1}{2}} +
d_{\hat{t}}(x_k,x_{0_k})
+ \frac{1}{10}H_k^{\frac{1}{2}}Q_k^{-\frac{1}{2}} \nn\\
&  \leq d_{t_k}(x_k,x_{0_k}) +
\frac{1}{2}H_k^{\frac{1}{2}}Q_k^{-\frac{1}{2}},\nn
\end{align}   
for all sufficiently large $k$. Then the combination of (7.2.2),
(7.2.3) and the choice of the points $(x_k,t_k)$ implies $t'=
t_k-\frac{1}{2}\varepsilon^{-2}Q_k^{-1}$ for all sufficiently
large $k$. (Here we also used the maximality of the subinterval
$[t',t_k]$ in the case that there is no time in $ [t',t_k]$ with
$R_k(x_k,\cdot) \geq 2Q_k$.)

Now we rescale the solutions $(M_k,g_k(\cdot,t))$ into
$(M_k,\tilde{g}_k(\cdot,t))$ around the points $x_k$ by the
factors $Q_k=R_k(x_k,t_k)$ and shift the times $t_k$ to the new
times zero. Then the same arguments from Step 1 to Step 3 in the
proof of Theorem 7.1.1 prove that a subsequence of the rescaled
solutions $(M_k,\tilde{g}_k(\cdot,t))$ converges in the
$C^\infty_{\rm loc}$ topology to a limiting (complete) solution
$(M_{\infty},\tilde{g}_{\infty}(\cdot,t))$, which is defined on a
backward time interval $[-a,0]$ for some $a>0$. (The only
modification is in Lemma 7.1.2 of Step 1 by further requiring
$t_k- \frac{1}{4}\varepsilon^{-2}Q_k^{-1} \leq \bar{t} \leq t_k$).

We next study how to adapt the argument of Step 4 in the proof of
Theorem 7.1.1. As before, we have a maximal time interval
$(t_{\infty},0]$ for which we can take a smooth limit
$(M_{\infty},\tilde{g}_{\infty}(\cdot,t),x_{\infty})$ from a
subsequence of the rescaled solutions
$(M_k,\tilde{g}_k(\cdot,t),x_k)$. We want to show
$t_{\infty}=-\infty$.

Suppose not; then $t_{\infty} > -\infty.$ Let $c>0$ be a positive
constant much smaller than $\frac{1}{10}\varepsilon^{-2}$. Note
that the infimum of the scalar curvature is nondecreasing in time.
Then we can find some point $y_{\infty} \in M_{\infty}$ and some
time $t=t_{\infty} + \theta$ with $0<\theta<\frac{c}{3}$ such that
$\tilde{R}_{\infty}(y_{\infty},t_{\infty}+\theta) \leq
\frac{3}{2}$.

Consider the (unrescaled) scalar curvature $R_k$ of
$(M_k,{g}_k(\cdot,t))$ at the point $x_k$ over the time interval
$[t_k+(t_{\infty}+\frac{\theta}{2})Q_k^{-1},t_k]$. Since the
scalar curvature $\tilde{R}_{\infty}$ of the limit on
$M_{\infty}\times [t_{\infty}+\frac{\theta}{3},0]$ is uniformly
bounded by some positive constant $C$, we have the curvature
estimate $$R_k(x_k,t) \leq 2CQ_k$$ for all $t \in
[t_k+(t_{\infty}+\frac{\theta}{2})Q_k^{-1},t_k]$ and all
sufficiently large $k$. Then by repeating the same arguments as in
deriving (7.2.1), (7.2.2) and (7.2.3), we deduce that the
conclusion of the claim with $K=2Q_k$ holds on the parabolic
neighborhood
$P(x_k,t_k,\frac{1}{10}H_k^{\frac{1}{2}}Q_k^{-\frac{1}{2}},
(t_{\infty}+\frac{\theta}{2})Q_k^{-1})$ for all sufficiently large
$k$.

Let $(y_k,t_k+(t_{\infty}+\theta_k)Q_k^{-1})$ be a sequence of
associated points and times in the (unrescaled) solutions
$(M_k,g_k(\cdot,t))$ so that after rescaling, the sequence
converges to the $(y_{\infty},t_{\infty}+\theta)$ in the limit.
Clearly $\frac{\theta}{2}\leq \theta_k \leq 2\theta$ for all
sufficiently large $k$. Then, by considering the scalar curvature
$R_k$ at the point $y_k$ over the time interval
$[t_k+(t_{\infty}-\frac{c}{3})Q_k^{-1},t_k+(t_{\infty}+\theta_k)Q_k^{-1}]$,
the above argument (as in deriving the similar estimates
(7.2.1)-(7.2.3)) implies that the conclusion of the claim with
$K=2Q_k$ holds on the parabolic neighborhood
$P(y_k,t_k,\frac{1}{10}H_k^{\frac{1}{2}}Q_k^{-\frac{1}{2}},
(t_{\infty}-\frac{c}{3})Q_k^{-1})$ for all sufficiently large $k$.
In particular, we have the curvature estimate
$$
R_k(y_k,t) \leq 4(1+\varepsilon)Q_k
$$
for $t \in [t_k+(t_{\infty}-\frac{c}{3})Q_k^{-1},
t_k+(t_{\infty}+\theta_k)Q_k^{-1}]$ for all sufficiently large
$k$.

We now consider the rescaled sequence $(M_k,\tilde{g}_k(\cdot,t))$
with the marked points replaced by $y_k$ and the times replaced by
$s_k \in [t_k+(t_{\infty}-\frac{c}{4})Q_k^{-1},
t_k+(t_{\infty}+\frac{c}{4})Q_k^{-1}]$. By applying the same
arguments from Step 1 to Step 3 in the proof of Theorem 7.1.1 and
the Li-Yau-Hamilton inequality as in Step 4 of Theorem 7.1.1, we
conclude that there is some small constant $a'>0$ such that the
original limit $(M_{\infty},\tilde{g}_{\infty}(\cdot,t))$ is
actually well defined on $M_{\infty} \times [t_{\infty}-a',0]$
with uniformly bounded curvature. This is a contradiction.
Therefore we have checked the claim.

To finish the proof, we next argue by contradiction. Suppose there
exist sequences of positive numbers $K_k\rightarrow+\infty$,
$\alpha_k\rightarrow0$, as $k\rightarrow+\infty$, and a sequence
of solutions $(M_k,g_k(\cdot,t))$ to the Ricci flow satisfying the
assumptions of the theorem with origins $x_{0_k}$ and with radii
$r_{0_k}$ satisfying $r_{0_k}<\alpha_k$ such that for some points
$x_k\in M_k$ with $d_{r^2_{0_k}}(x_k,x_{0_k})<Ar_{0_k}$ we have
\be
R(x_k,r^2_{0_k})>K_kr^{-2}_{0_k} 
\ee for all $k$. Let $(M_k,\hat{g}_k(\cdot,t),x_{0_k})$ be the
rescaled solutions of $(M_k,g_k(\cdot,t))$ around the origins
$x_{0_k}$ by the factors $r^{-2}_{0_k}$ and shifting the times
$r^2_{0_k}$ to the new times zero. The above claim tells us that
for $k$ large, any point $(y,0) \in
(M_k,\hat{g}_k(\cdot,0),x_{0_k})$ with
$d_{\hat{g}_k(\cdot,0)}(y,x_{0_k}) < 3A$ and with the rescaled
scalar curvature $\hat{R}_k(y,0) > K_k$ has a canonical
neighborhood which is either a $2\varepsilon$-neck, or a
$2\varepsilon$-cap, or a compact manifold (without boundary)
diffeomorphic to a metric quotient of the round three-sphere. Note
that the pinching estimate (7.1.1) and the condition
$\alpha_k\rightarrow0$ imply any subsequential limit of the
rescaled solutions $(M_k,\hat{g}_k(\cdot,t),x_{0_k})$ must have
nonnegative sectional curvature. Thus the same argument as in Step
2 of the proof of Theorem 7.1.1 shows that for all sufficiently
large $k$, the curvatures of the rescaled solutions at the time
zero stay uniformly bounded at those points whose distances from
the origins $x_{0_k}$ do not exceed $2A$. This contradicts (7.2.4)
for $k$ large enough.

Therefore we have completed the proof of the theorem.
\end{pf}

The next result is a generalization of the curvature estimate via
volume growth in Theorem 6.3.3 (ii) where the condition on the
curvature lower bound over a time interval is replaced by that at
a time slice only.

\begin{theorem}[{Perelman \cite{P1}}]
For any $w>0$ there exist $\tau=\tau(w)>0$, $K=K(w)<+\infty$,
$\alpha=\alpha(w)>0$ with the following property. Suppose we have
a three-dimensional, compact and orientable solution to the Ricci
flow defined on $M\times[0,T)$ with normalized initial metric.
Suppose that for some radius $r_0>0$ with $r_0<\alpha$ and a point
$(x_0,t_0)\in M\times[0,T)$ with $T
> t_0\geq4\tau r^2_0$, the solution on the ball $B_{t_0}(x_0,r_0)$
satisfies
$$
Rm(x,t_0)\geq-r_0^{-2}\ \ \ on\ B_{t_0}(x_0,r_0),
$$
$$
and\ \ \ \Vol_{t_0}(B_{t_0}(x_0,r_0))\geq wr_0^3.
$$
Then $R(x,t)\leq Kr_0^{-2}$ whenever $t\in[t_0-\tau r^2_0,t_0]$
and $d_t(x,x_0)\leq\frac{1}{4}r_0.$
\end{theorem}

\begin{pf} The following argument basically follows
the proof of Theorem 12.3 of Perelman \cite{P1}.

 If we knew that
$$
Rm(x,t)\geq-r_0^{-2}
$$
for all $t\in[0,t_0]$ and $d_t(x,x_0)\leq r_0$, then we could just
apply Theorem 6.3.3 (ii) and take $\tau(w)=\tau_0(w)/2$,
$K(w)=C(w)+2B(w)/\tau_0(w)$. Now fix these values of $\tau$ and
$K$.

We argue by contradiction. Consider a three-dimensional, compact
and orientable solution $g_{ij}(t)$ to the Ricci flow with
normalized initial metric, a point $(x_0,t_0)$ and some radius
$r_0>0$ with $r_0<\alpha$, for $\alpha>0$ a sufficiently small
constant to be determined later, such that the assumptions of the
theorem do hold whereas the conclusion does not. We first claim
that we may assume that any other point $(x',t')$ and radius
$r'>0$ with the same property has either $t'>t_0$ or $t'<t_0-2\tau
r_0^2$, or $2r'>r_0$. Indeed, suppose otherwise. Then there exist
$(x_0',t_0')$ and $r_0'$ with $t_0'\in[t_0-2\tau r_0^2,t_0]$ and
$r_0'\leq\frac{1}{2}r_0$, for which the assumptions of the theorem
hold but the conclusion does not. Thus, there is a point $(x,t)$
such that
$$
t\in[t_0'-\tau(r_0')^2,t_0']\subset\left[t_0-2\tau
r^2_0-\frac{\tau}{4}r_0^2,t_0\right]
$$
$$
\mbox{and}\ \ \ R(x,t)> K(r_0')^{-2}\geq4Kr_0^{-2}.
$$
If the point $(x_0',t_0')$ and the radius $r_0'$ satisfy the claim
then we stop, and otherwise we iterate the procedure. Since
$t_0\geq4\tau r^2_0$ and the solution is smooth, the iteration
must terminate in a finite number of steps, which provides the
desired point and the desired radius.

Let $\tau'\geq0$ be the largest number such that \be
Rm(x,t)\geq-r_0^{-2}  
\ee whenever $t\in[t_0-\tau'r_0^2,t_0]$ and $d_t(x,x_0)\leq r_0$.
If $\tau'\geq2\tau$, we are done by Theorem 6.3.3 (ii). Thus we
may assume $\tau'<2\tau$.  By applying Theorem 6.3.3(ii), we know
that at time $t'=t_0-\tau'r_0^2$, the ball $B_{t'}(x_0,r_0)$ has
\be
\Vol_{t'}(B_{t'}(x_0,r_0))\geq \xi(w) r_0^3 
\ee for some positive constant $\xi(w)$ depending only on $w$. We
next claim that there exists a ball (at time $t'=t_0-\tau'r_0^2$)
$B_{t'}(x',r')\subset B_{t'}(x_0,r_0)$ with \be
\Vol_{t'}(B_{t'}(x',r'))\geq\frac{1}{2}\alpha_3(r')^3 
\ee and with \be
\frac{r_0}{2} > r'\geq c(w)r_0 
\ee for some small positive constant $c(w)$ depending only on $w$,
where $\alpha_3$ is  the volume of the unit ball $\mathbb{B}^3$ in
the Euclidean space $\mathbb{R}^3$. (The following argument in
deriving (7.2.7) and (7.2.8) is a standard one in the Alexandrov
space theory and has nothing to do with the Ricci flow. Our
presentation here is inspired from Lemma 53.1 of Kleiner-Lott notes
\cite{KL}.)

Indeed, suppose that it is not true. Then after rescaling, there
is a sequence of Riemannian manifolds $M_i,\ i=1,2,\ldots,$ with
balls $B(x_i,1)\subset M_i$ so that
$$
Rm\geq-1\ \ \ \mbox{on}\ B(x_i,1) \leqno(7.2.5)'
$$
and
$$
\Vol(B(x_i,1))\geq \xi(w) \leqno(7.2.6)'
$$
for all $i$, but all balls $B(x_i',r_i')\subset B(x_i,1)$ with
$\frac{1}{2} > r_i'\geq\frac{1}{i}$ satisfy \be
\Vol(B(x_i',r_i'))<\frac{1}{2}\alpha_3(r_i')^3. 
\ee It follows from basic results in Alexandrov space theory (see
for example Theorem 10.7.2 and Theorem 10.10.10 of \cite{BBI})
that, after taking a subsequence, the marked balls
$(B(x_i,1),x_i)$ converge in the Gromov-Hausdorff topology to a
marked length space $(B_\infty,x_\infty)$ with curvature bounded
from below by $-1$ in the Alexandrov space sense, and the
associated Riemannian volume forms $d\Vol_{M_i}$ over
$(B(x_i,1),x_i)$ converge weakly to the Hausdorff measure $\mu$ of
$B_\infty$. It is well-known that the Hausdorff dimension of any
Alexandrov space is either an integer or infinity (see for example
Theorem 10.8.2 of \cite{BBI}). Then by (7.2.6)$'$, we know the
limit $(B_\infty,x_\infty)$ is a three-dimensional Alexandrov
space of curvature $\geq -1$. In the Alexandrov space theory, a
point $p \in B_\infty$ is said to be
\textbf{regular}\index{regular} if the tangent cone of $B_\infty$
at $p$ is isometric to $\mathbb{R}^3$. It is also a basic result
in Alexandrov space theory (see for example Corollary 10.9.13 of
\cite{BBI}) that the set of regular points in $B_\infty$ is dense
and for each regular point there is a small neighborhood which is
almost isometric to an open set of the Euclidean space
$\mathbb{R}^3$. Thus for any $\varepsilon>0$, there are balls
$B(x'_\infty,r'_\infty)\subset B_\infty$ with $0 < r'_\infty <
\frac{1}{3}$ and satisfying
$$
\mu(B(x'_\infty,r'_\infty))\geq(1-\varepsilon)\alpha_3(r'_\infty)^3.
$$
This is a contradiction with (7.2.9).

Without loss of generality, we may assume
$w\leq\frac{1}{4}\alpha_3$. Since $\tau'<2\tau$, it follows from
the choice of the point $(x_0,t_0)$ and the radius $r_0$ and
(7.2.5), (7.2.7), (7.2.8) that the conclusion of the theorem holds
for $(x',t')$ and $r'$. Thus we have the estimate
$$
R(x,t)\leq K(r')^{-2}
$$
whenever $t\in[t'-\tau(r')^2,t']$ and
$d_t(x,x')\leq\frac{1}{4}r'$. For $\alpha > 0$ small, by combining
with the pinching estimate (7.1.1), we have
$$
|Rm(x,t)|\leq K'(r')^{-2}
$$
whenever $t\in[t'-\tau(r')^2,t']$ and
$d_t(x,x')\leq\frac{1}{4}r'$, where $K'$ is some positive constant
depending only on $K$. Note that this curvature estimate implies
the evolving metrics are equivalent over a suitable subregion of
$\{(x,t)\ |\ t\in[t'-\tau(r')^2,t'] \ \ \mbox{and} \ \
d_t(x,x')\leq\frac{1}{4}r'\}$. Now we can apply Theorem 7.2.1 to
choose $\alpha=\alpha(w)>0$ so small that \be
R(x,t)\leq\tilde{K}(w)(r')^{-2}
\leq\tilde{K}(w)c(w)^{-2}r_0^{-2} 
\ee whenever $t\in[t'-\frac{\tau}{2}(r')^2,t']$ and
$d_t(x,x')\leq10r_0$. Then the combination of (7.2.10) with the
pinching estimate (7.1.2) would imply
\begin{displaymath}
\begin{split}
Rm(x,t)&  \geq-[f^{-1}(R(x,t)(1+t))/(R(x,t)(1+t))]R(x,t)\\[1mm]
    &  \geq-\frac{1}{2}r_0^{-2}
\end{split}
\end{displaymath}
on the region $\{(x,t)\ |\ t\in[t'-\frac{\tau}{2}(r')^2,t']$ and
$d_t(x,x_0)\leq r_0\}$ when $\alpha=\alpha(w)>r_0$ small enough.
This contradicts the choice of $\tau'$. Therefore we have proved
the theorem.
\end{pf}

The combination of the above two theorems immediately gives the
following consequence.

\begin{corollary}
For any $w>0$ and $A<+\infty$, there exist $\tau=\tau(w,A)>0$,
$K=K(w,A)<+\infty$, and $\alpha=\alpha(w,A)>0$ with the following
property. Suppose we have a three-dimensional, compact and
orientable solution to the Ricci flow defined on $M\times[0,T)$
with normalized initial metric. Suppose that for some radius
$r_0>0$ with $r_0 <\alpha$ and a point $(x_0,t_0)\in M\times[0,T)$
with $T > t_0\geq 4\tau r^2_0$, the solution on the ball
$B_{t_0}(x_0,r_0)$ satisfies
$$
Rm(x,t_0)\geq-r_0^{-2}\ \ \ on\ B_{t_0}(x_0,r_0),
$$
$$
and\ \ \ \Vol_{t_0}(B_{t_0}(x_0,r_0))\geq wr_0^3.
$$
Then $R(x,t)\leq Kr_0^{-2}$ whenever $t\in[t_0-\tau r^2_0,t_0]$
and $d_t(x,x_0)\leq Ar_0.$
\end{corollary}

We can also state the previous corollary in the following version.

\begin{corollary}[{Perelman \cite{P1}}]
For any $w>0$ one can find $\rho = \rho(w)>0$ such that if
$g_{ij}(t)$ is a complete solution to the Ricci flow defined on
$M\times[0,T)$ with $T>1$ and with normalized initial metric,
where $M$ is a three-dimensional, compact and orientable manifold,
and if $B_{t_0}(x_0,r_0)$ is a metric ball at time $t_0\geq1$,
with $r_0<\rho$, such that$$\min\{ Rm(x,t_0)\ |\ x\in
B_{t_0}(x_0,r_0)\}=-r_0^{-2},$$ then
$$
\Vol_{t_0}(B_{t_0}(x_0,r_0))\leq wr_0^3.
$$
\end{corollary}

\begin{pf}
We argue by contradiction. Suppose for any $\rho>0$, there is a
solution and a ball $B_{t_0}(x_0,r_0)$ satisfying the assumption
of the corollary with $r_0<\rho$, $t_0\geq1$, and with
$$
\min\{ Rm(x,t_0)\ |\ x\in B_{t_0}(x_0,r_0)\}=-r_0^{-2},
$$
but
$$
\Vol_{t_0}(B_{t_0}(x_0,r_0))>wr_0^3.
$$
We can apply Corollary 7.2.3 to get
$$
R(x,t)\leq Kr_0^{-2}
$$
whenever $t\in[t_0-\tau r_0^2,t_0]$ and $d_t(x,x_0)\leq2r_0$,
provided $\rho>0$ is so small that $4\tau\rho^2\leq1$ and
$\rho<\alpha$, where $\tau,\ \alpha$ and $K$ are the positive
constants obtained in Corollary 7.2.3. Then for $r_0<\rho$ and
$\rho>0$ sufficiently small, it follows from the pinching estimate
(7.1.2) that
\begin{displaymath}
\begin{split}
Rm(x,t)&  \geq -[f^{-1}(R(x,t)(1+t))/(R(x,t)(1+t))]R(x,t)\\[1mm]
           &  \geq-\frac{1}{2}r_0^{-2}
\end{split}
\end{displaymath}
in the region $\{(x,t)\ |\ t\in[t_0-\tau(r_0)^2,t_0]$ and
$d_t(x,x_0)\leq 2r_0\}$. In particular, this would imply
$$
\min\{ Rm(x,t_0)|x\in B_{t_0}(x_0,r_0)\}>-r_0^{-2}.
$$
This contradicts the assumption.
\end{pf}

\section{Ricci Flow with Surgery}

One of the central themes of the Ricci flow theory is to give a
classification of all compact orientable three-manifolds. As we
mentioned before, the basic idea is to obtain long-time behavior
of solutions to the Ricci flow. However the solutions will in
general become singular in finite time. Fortunately, we now
understand the precise structures of the solutions around the
singularities, thanks to Theorem 7.1.1. When a solution develops
singularities, one can perform geometric surgeries by cutting off
the canonical neighborhoods around the singularities and gluing
back some known pieces, and then continue running the Ricci flow.
By repeating this procedure, one hopes to get a kind of ``weak"
solution. In this section we will give a detailed description of
this surgery procedure (cf. \cite{Ha97, P2}) and define a global
``weak" solution to the Ricci flow. This section is a detailed
exposition of sections 3 and 4 of Perelman \cite{P2}.

Given any $\varepsilon>0$, based on the singularity structure
theorem (Theorem 7.1.1), we can get a clear picture of the
solution near the singular time as follows.

Let $(M,g_{ij}(\cdot,t))$ be a maximal solution to the Ricci flow
on the maximal time interval $[0,T)$ with $T<+\infty$, where $M$
is a connected compact orientable three-manifold and the initial
metric is normalized. For the given $\varepsilon>0$ and the
solution $(M,g_{ij}(\cdot,t))$, we can find $r_0>0$ such that each
point $(x,t)$ with $R(x,t)\ge r_0^{-2}$ satisfies the derivative
estimates \be |\nabla R(x,t)|<\eta R^{\frac{3}{2}}(x,t)\quad
\mbox{and }\; \left|\frac{\partial}{\partial t}R(x,t)\right|<\eta
R^2(x,t),
\ee where $\eta>0$ is a universal constant, and has a canonical
neighborhood which is either an evolving $\varepsilon$-neck, or an
evolving $\varepsilon$-cap, or a compact positively curved
manifold (without boundary). In the last case the solution becomes
extinct at the maximal time $T$ and the manifold $M$ is
diffeomorphic to the round three-sphere $\mathbb{S}^3$ or a metric
quotient of $\mathbb{S}^3$ by Theorem 5.2.1.

Let $\Omega$ denote the set of all points in $M$ where the
curvature stays bounded as $t\rightarrow T$. The gradient
estimates in (7.3.1) imply that $\Omega$ is open and that
$R(x,t)\rightarrow\infty$ as $t\rightarrow T$ for each $x\in
M\setminus \Omega$.

If $\Omega$ is empty, then the solution becomes extinct at time
$T$. In this case, either the manifold $M$ is compact and
positively curved, or it is entirely covered by evolving
$\varepsilon$-necks and evolving $\varepsilon$-caps shortly before
the maximal time $T$. So the manifold $M$ is diffeomorphic to
either $\mathbb{S}^3$, or a metric quotient of the round
$\mathbb{S}^3$, or $\mathbb{S}^2\times \mathbb{S}^1$, or
$\mathbb{RP}^3\#\mathbb{RP}^3$. The reason is as follows. Clearly,
we only need to consider the situation that the manifold $M$ is
entirely covered by evolving $\varepsilon$-necks and evolving
$\varepsilon$-caps shortly before the maximal time $T$. If $M$
contains a cap $C$, then there is a cap or a neck adjacent to the
neck-like end of $C$. The former case implies that $M$ is
diffeomorphic to $\mathbb{S}^3$, $\mathbb{RP}^3$, or
$\mathbb{RP}^3\#\mathbb{RP}^3$. In the latter case, we get a new
longer cap and continue. Finally, we must end up with a cap,
producing a $\mathbb{S}^3$, $\mathbb{RP}^3$, or
$\mathbb{RP}^3\#\mathbb{RP}^3$. If $M$ contains no caps, we start
with a neck $N$. By connecting with the necks that are adjacent to
the boundary of $N$, we get a longer neck and continue. After a
finite number of steps, the resulting neck must repeat itself.
Since $M$ is orientable, we conclude that $M$ is diffeomorphic to
$\mathbb{S}^2\times \mathbb{S}^1$.

We can now assume that $\Omega$ is nonempty. By using the local
derivative estimates of Shi (Theorem 1.4.2), we see that as
$t\rightarrow T$ the solution $g_{ij}(\cdot,t)$ has a smooth limit
$\bar{g}_{ij}(\cdot)$ on $\Omega$. Let $\bar{R}(x)$ denote the
scalar curvature of $\bar{g}_{ij}$. For any $\rho<r_0$, let us
consider the set
$$
\Omega_\rho=\{x\in\Omega\ |\ \bar{R}(x)\le\rho^{-2}\}.
$$
By the evolution equation of the Ricci flow, we see that the
initial metric $g_{ij}(\cdot,0)$ and the limit metric
$\overline{g}_{ij}(\cdot)$ are equivalent over any fixed region
where the curvature remains uniformly bounded. Note that for any
fixed $x\in
\partial \Omega$, and any sequence of points $x_j\in \Omega$ with
$x_j\rightarrow x$ with respect to the initial metric
$g_{ij}(\cdot,0)$, we have $\overline{R}(x_j)\rightarrow +\infty$.
In fact, if there were a subsequence $x_{j_k}$ so that
$\lim_{k\rightarrow \infty}\overline{R}(x_{j_k})$ exists and is
finite, then it would follow from the gradient estimates (7.3.1)
that $\overline{R}$ is uniformly bounded in some small
neighborhood of $x\in \partial \Omega$ (with respect to the
induced topology of the initial metric $g_{ij}(\cdot,0)$); this is
a contradiction. From this observation and the compactness of the
initial manifold, we see that $\Omega_{\rho}$ is compact (with
respect to the metric $\overline{g}_{ij}(\cdot)$).

For further discussions, let us introduce the following
terminologies. Denote by $\mathbb{I}$ an interval.

Recall that an {\bf $\varepsilon$-neck}\index{$\varepsilon$-neck}
(of radius $r$) is an open set with a Riemannian metric, which is,
after scaling the metric with factor $r^{-2}$, $\varepsilon$-close
to the standard neck $\mathbb{S}^2\times \mathbb{I}$ with the
product metric, where $\mathbb{S}^2$ has constant scalar curvature
one and $\mathbb{I}$ has length $2\varepsilon^{-1}$ and the
$\varepsilon$-closeness refers to the $C^{[\varepsilon^{-1}]}$
topology.

A metric on $\mathbb{S}^2\times \mathbb{I}$, such that each point
is contained in some $\varepsilon$-neck, is called an {\bf
$\varepsilon$-tube}\index{$\varepsilon$-tube}, or an {\bf
$\varepsilon$-horn}\index{$\varepsilon$-horn}, or a {\bf double
$\varepsilon$-horn}\index{double $\varepsilon$-horn}, if the
scalar curvature stays bounded on both ends, or stays bounded on
one end and tends to infinity on the other, or tends to infinity
on both ends, respectively.

A metric on $\mathbb{B}^3$ or $\mathbb{RP}^3\setminus
\bar{\mathbb{B}}^3$ is called a {\bf
$\varepsilon$-cap}\index{$\varepsilon$-cap} if the region outside
some suitable compact subset is an $\varepsilon$-neck. A metric on
$\mathbb{B}^3$ or $\mathbb{RP}^3\setminus \bar{\mathbb{B}}^3$ is
called an {\bf capped $\varepsilon$-horn}\index{capped
$\varepsilon$-horn} if each point outside some compact subset is
contained in an $\varepsilon$-neck and the scalar curvature tends
to infinity on the end.

Now take any $\varepsilon$-neck in $(\Omega,\bar{g}_{ij})$ and
consider a point $x$ on one of its boundary components. If $x\in
\Omega\setminus\Omega_\rho$, then there is either an
$\varepsilon$-cap or an $\varepsilon$-neck, adjacent to the
initial $\varepsilon$-neck. In the latter case we can take a point
on the boundary of the second $\varepsilon$-neck and continue.
This procedure can either terminate when we get into $\Omega_\rho$
or an $\varepsilon$-cap, or go on indefinitely, producing an
$\varepsilon$-horn. The same procedure can be repeated for the
other boundary component of the initial $\varepsilon$-neck.
Therefore, taking into account that $\Omega$ has no compact
components, we conclude that each $\varepsilon$-neck of
$(\Omega,\bar{g}_{ij})$ is contained in a subset of $\Omega$ of
one of the following types:
\begin{align}
\text{(a)}&\quad \text{an $\varepsilon$-tube with boundary
components in $\Omega_\rho$, or} \nn\\[1mm]
\text{(b)}&\quad \text{an $\varepsilon$-cap with boundary
in $\Omega_\rho$, or} \nn\\[1mm]
\text{(c)}&\quad \text{an $\varepsilon$-horn with boundary
in $\Omega_\rho$, or}\\[1mm]  
\text{(d)}&\quad \text{a capped $\varepsilon$-horn, or} \nn\\[1mm]
\text{(e)}&\quad \text{a double $\varepsilon$-horn.} \nn
\end{align}

Similarly, each $\varepsilon$-cap of $(\Omega,\bar{g}_{ij})$ is
contained in a subset of $\Omega$ of either type (b) or type (d).

It is clear that there is a definite lower bound (depending on
$\rho$) for the volume of subsets of types (a), (b) and (c), so
there can be only a finite number of them. Thus we conclude that
there is only a finite number of components of $\Omega$ containing
points of $\Omega_\rho$, and every such component has a finite
number of ends, each being an $\varepsilon$-horn. On the other
hand, every component of $\Omega$, containing no points of
$\Omega_\rho$, is either a capped $\varepsilon$-horn, or a double
$\varepsilon$-horn. If we look at the solution at a slightly
earlier time, the above argument shows each $\varepsilon$-neck or
$\varepsilon$-cap of $(M,g_{ij}(\cdot,t))$ is contained in a
subset of types (a) or (b), while the $\varepsilon$-horns, capped
$\varepsilon$-horns and double $\varepsilon$-horns (at the maximal
time T) are connected together to form $\varepsilon$-tubes and
$\varepsilon$-caps at the times $t$ shortly before $T$.

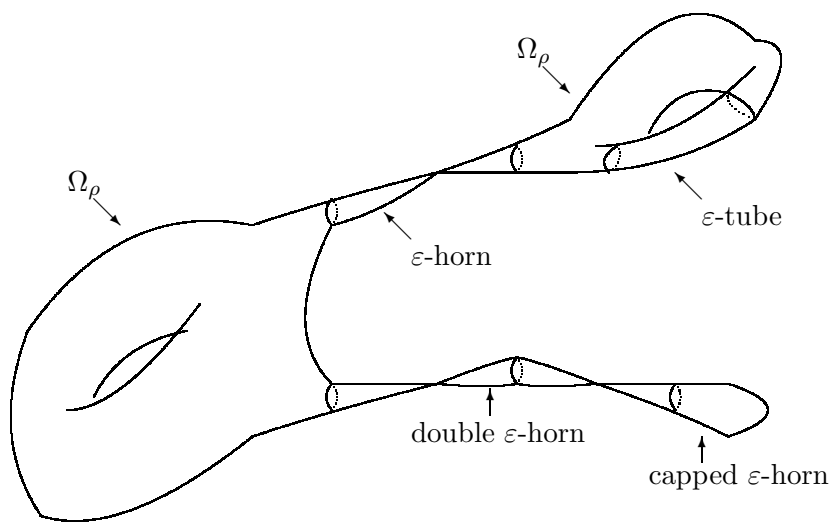
\begin{figure}[h]
\begin{picture}(320, 240)
\qbezier(15,90)(0,50)(20,20)\qbezier(20,20)(50,10)(100,50)\qbezier(100,50)(130,60)(170,70)
\qbezier(15,90)(50,140)(100,130)\qbezier(100,130)(130,140)(170,150)\qbezier(130,130)(110,90)(130,70)
\qbezier(130,140)(126,135)(130,130)\qbezier(130,70)(126,65)(130,60)
\qbezier[8](130,140)(134,135)(130,130)\qbezier[8](130,70)(134,65)(130,60)
\qbezier(130,130)(150,135)(170,150)\qbezier(130,70)(150,70)(170,70)
\qbezier(170,70)(190,78)(200,80)\qbezier(170,70)(190,68)(200,70)
\qbezier(200,80)(210,78)(230,70)\qbezier(200,70)(210,68)(230,70)
\qbezier(200,80)(196,75)(200,70)\qbezier[8](200,80)(204,75)(200,70)
\qbezier(230,70)(260,70)(280,70)\qbezier(230,70)(260,60)(280,50)\qbezier(280,70)(310,60)(280,50)
\qbezier(260,70)(256,65)(260,60)\qbezier[8](260,70)(264,65)(260,60)
\qbezier(170,150)(200,160)(220,170)\qbezier(170,150)(200,150)(220,150)
\qbezier[8](200,160)(204,155)(200,150)\qbezier(200,160)(196,155)(200,150)
\qbezier(220,170)(260,230)(290,200)\qbezier(220,150)(260,150)(290,170)\qbezier(290,200)(310,200)(290,170)
\qbezier(230,160)(260,160)(290,190)\qbezier(250,165)(260,185)(280,180)
\qbezier(30,60)(50,60)(80,100)\qbezier(40,65)(50,85)(75,90)
\qbezier(238,160)(230,155)(235,150)\qbezier[8](238,160)(241,155)(235,150)
\qbezier(280,180)(290,175)(290,170)\qbezier[10](280,180)(278,175)(290,170)

\put(30,140){\makebox(0,0)[bl]{$\Omega_{\rho}$}}
\put(40,140){\vector(1,-1){10}}
\put(200,190){\makebox(0,0)[bl]{$\Omega_{\rho}$}}
\put(210,190){\vector(1,-1){10}}
\put(160,115){\makebox(0,0)[bl]{$\varepsilon$-horn}}
\put(160,125){\vector(-1,1){10}}
\put(270,130){\makebox(0,0)[bl]{$\varepsilon$-tube}}
\put(270,140){\vector(-1,1){10}}
\put(160,48){\makebox(0,0)[bl]{double $\varepsilon$-horn}}
\put(190,58){\vector(0,1){10}}
 \put(250,30){\makebox(0,0)[bl]{capped $\varepsilon$-horn}}
 \put(270,40){\vector(0,1){10}}

\end{picture}
\end{figure}

Hence, by looking at the solution at times shortly before $T$, we
see that the topology of $M$ can be reconstructed as follows: take
the components $\Omega_j$, $1\le j\le k$, of $\Omega$ which
contain points of $\Omega_\rho$, truncate their
$\varepsilon$-horns, and glue to the boundary components of
truncated $\Omega_j$ a finite collection of tubes
$\mathbb{S}^2\times\mathbb{I}$ and caps $\mathbb{B}^3$ or
$\mathbb{RP}^3 \setminus \bar{\mathbb{B}}^3$. Thus, $M$ is
diffeomorphic to a connected sum of $\bar{\Omega}_j$, $1\le j\le
k$, with a finite number of copies of $\mathbb{S}^2\times
\mathbb{S}^1$ (which correspond to gluing a tube to two boundary
components of the same $\Omega_j$), and a finite number of copies
of $\mathbb{RP}^3$. Here $\bar{\Omega}_j$ denotes $\Omega_j$ with
each $\varepsilon$-horn one point compactified. More
geometrically, one can get $\bar{\Omega}_j$ in the following way:
in every $\varepsilon$-horn of $\Omega_j$ one can find an
$\varepsilon$-neck, cut it along the middle two-sphere, remove the
horn-shaped end, and glue back a cap (or more precisely, a
differentiable three-ball). Thus to understand the topology of
$M$, one only needs to understand the topologies of the compact
orientable three-manifolds $\bar{\Omega}_j$, $1\le j\le k$.

Naturally one can evolve each $\bar{\Omega}_j$ by the Ricci flow
again and, when singularities develop again, perform the above
surgery for each $\varepsilon$-horn to get new compact orientable
three-manifolds. By repeating this procedure indefinitely, it will
likely give a long-time ``weak" solution to the Ricci flow. The
following abstract definition for this kind of ``weak" solution
was introduced by Perelman in \cite{P2} .

\begin{definition}
Suppose we are given a (finite or countably infinite) collection
of three-dimensional smooth solutions $g^k_{ij}(t)$ to the Ricci
flow defined on $M_k\times [t^-_k,t^+_k)$ and go singular as
$t\rightarrow t^+_k$, where each manifold $M_k$ is compact and
orientable, possibly disconnected with only a finite number of
connected components. Let $(\Omega_k,\bar{g}^k_{ij})$ be the
limits of the corresponding solutions $g^k_{ij}(t)$ as
$t\rightarrow t^+_k$, as above. Suppose also that for each $k$ we
have $t^-_k=t^+_{k-1}$, and that
$(\Omega_{k-1},\bar{g}^{k-1}_{ij})$ and $(M_k,g^k_{ij}(t^-_k))$
contain compact (possibly disconnected) three-dimensional
submanifolds with smooth boundary which are isometric. Then by
identifying these isometric submanifolds, we say the collection of
solutions $g^k_{ij}(t)$ is a solution to the \textbf{Ricci flow
with surgery}\index{Ricci flow with surgery} (or a
\textbf{surgically modified solution} to the Ricci flow)
\index{surgically modified solution} on the time interval which is
the union of all $[t^-_k,t^+_k)$, and say the times $t^+_k$ are
\textbf{surgery
times}\index{surgery times}. 
\end{definition}

To get the topology of the initial manifold from the solution to
the Ricci flow with surgery, one has to overcome the following two
difficulties:
\begin{itemize}
\item[(i)] how to prevent the surgery times from accumulating?
\item[(ii)] how to obtain the long time behavior of the solution
to the Ricci flow with surgery?
\end{itemize}

Thus it is natural to consider those solutions having ``good"
properties. For any arbitrarily fixed positive number
$\varepsilon$, we will only consider those solutions to the Ricci
flow with surgery which satisfy the following \textbf{a priori
assumptions (with accuracy $\varepsilon$)}\index{a priori
assumptions (with accuracy $\varepsilon$)}.

\medskip
{\bf Pinching assumption.}\index{pinching assumption} \ The
eigenvalues $\lambda \geq \mu \geq \nu$ of the curvature operator
of the solution to the Ricci flow with surgery at each point and
each time satisfy \be
R \geq (-\nu)[\log(-\nu) + \log(1+t) -3] 
\ee whenever $\nu < 0$.

\medskip
{\bf Canonical neighborhood assumption (with accuracy $\boldsymbol
\varepsilon$).}\index{canonical neighborhood assumption (with
accuracy $\varepsilon$)} \ For any given $\varepsilon >0$, there
exist positive constants $C_1$ and $C_2$ depending only on
$\varepsilon$, and a nonincreasing positive function $r:
[0,+\infty) \rightarrow (0,+\infty)$ such that at each time $t>0$,
every point $x$ where scalar curvature $R(x,t)$ is at least
$r^{-2}(t)$ has a neighborhood $B$, with $B_t(x,\sigma)\subset B
\subset B_t(x,2\sigma)$ for some $ 0<\sigma <
C_1R^{-\frac{1}{2}}(x,t)$, which falls into one of the following
three categories:
\begin{itemize}
\item[(a)] $B$ is a \textbf{strong
$\varepsilon$-neck}\index{$\varepsilon$-neck!strong} (in the sense
$B$ is the slice at time $t$ of the parabolic neighborhood
$\{(x',t')\ |\ x'\in B, t'\in[t-R(x,t)^{-1},t]\}$, where the
solution is well defined on the whole parabolic neighborhood and
is, after scaling with factor $R(x,t)$ and shifting the time to
zero, $\varepsilon$-close (in the $C^{[\varepsilon^{-1}]}$
topology) to the subset $(\mathbb{S}^2 \times \mathbb{I}) \times
[-1,0]$ of the evolving standard round cylinder with scalar
curvature 1 to $\mathbb{S}^2$ and length $2\varepsilon^{-1}$ to
$\mathbb{I}$ at time zero), or \item[(b)] $B$ is an
$\varepsilon$-cap, or \item[(c)] $B$ is a compact manifold
(without boundary) of positive sectional curvature.
\end{itemize}
Furthermore, the scalar curvature in $B$ at time $t$ is between
$C^{-1}_2R(x,t)$ and $C_2R(x,t)$, satisfies the gradient estimates
\be |\nabla R|<\eta R^{\frac{3}{2}}\;  \mbox{ and }\;
\left|\frac{\partial R}{\partial t}\right|<\eta R^2, 
\ee and the volume of $B$ in case (a) and case (b) satisfies
$$
(C_2R(x,t))^{-\frac{3}{2}} \leq \Vol_t(B) \leq \varepsilon
\sigma^3.
$$
Here $\eta$ is a universal positive constant.

\medskip
Without loss of generality, we always assume the above constants
$C_1$ and $C_2$ are twice bigger than the corresponding constants
$C_1(\frac{\varepsilon}{2})$ and $C_2(\frac{\varepsilon}{2})$ in
Theorem 6.4.6 with the accuracy $\frac{\varepsilon}{2}$.

We remark that the above definition of the canonical neighborhood
assumption is slightly different from that of Perelman in
\cite{P2} in two aspects: (1) it allows the parameter $r$ to
depend on time; (2) it also includes an volume upper bound for the
canonical neighborhoods of types (a) and (b).

Arbitrarily given a compact orientable three-manifold with a
Riemannian metric, by scaling, we may assume the Riemannian metric
is normalized. In the rest of this section and the next section,
we will show the Ricci flow with surgery, with the normalized
metric as initial data, has a long-time solution which satisfies
the above a priori assumptions and has only a finite number of
surgery times at each finite time interval. The construction of
the long-time solution will be given by an induction argument.

First, for the arbitrarily given compact orientable normalized
three-dimensional Riemannian manifold $(M,g_{ij}(x))$, the Ricci
flow with it as initial data has a maximal solution $g_{ij}(x,t)$
on a maximal time interval $[0,T)$ with $T>1$. It follows from
Theorem 5.3.2 and Theorem 7.1.1 that the a priori assumptions
(with accuracy $\varepsilon$) hold for the smooth solution on
$[0,T)$. If $T=+\infty$, we have the desired long time solution.
Thus, without loss of generality, we may assume the maximal time
$T<+\infty$ so that the solution goes singular at time $T$.

Suppose that we have a solution to the Ricci flow with surgery,
with the normalized metric as initial data, satisfying the a
priori assumptions (with accuracy $\varepsilon$), defined on
$[0,T)$ with $T<+\infty$, going singular at time $T$ and having
only a finite number of surgery times on $[0,T)$. Let $\Omega$
denote the set of all points where the curvature stays bounded as
$t\rightarrow T$. As we have seen before, the canonical
neighborhood assumption implies that $\Omega$ is open and that
$R(x,t)\rightarrow\infty$ as $t\rightarrow T$ for all $x$ lying
outside $\Omega$. Moreover, as $t\rightarrow T$, the solution
$g_{ij}(x,t)$ has a smooth limit $\bar{g}_{ij}(x)$ on $\Omega$.

For some $\delta>0$ to be chosen much smaller than $\varepsilon$,
we let $\rho=\delta r(T)$ where $r(t)$ is the positive
nonincreasing function in the definition of the canonical
neighborhood assumption. We consider the corresponding compact set
$$
\Omega_\rho=\{x\in\Omega\ |\ \bar{R}(x)\leq\rho^{-2}\},
$$
where $\bar{R}(x)$ is the scalar curvature of $\bar{g}_{ij}$. If
$\Omega_\rho$ is empty, the manifold (near the maximal time $T$)
is entirely covered by $\varepsilon$-tubes, $\varepsilon$-caps and
compact components with positive curvature. Clearly, the number of
the compact components is finite. Then in this case the manifold
(near the maximal time $T$) is diffeomorphic to the union of a
finite number of copies of $\mathbb{S}^3$, or metric quotients of
the round $\mathbb{S}^3$, or $\mathbb{S}^2 \times \mathbb{S}^1$,
or a connected sum of them. Thus when $\Omega_\rho$ is empty, the
procedure stops here, and we say the \textbf{solution becomes
extinct}\index{solution becomes extinct}. We now assume
$\Omega_\rho$ is not empty. Then we know that every point
$x\in\Omega\setminus\Omega_\rho$ lies in one of the subsets of
$\Omega$ listed in (7.3.2), or in a compact component with
positive curvature, or in a compact component which is contained
in $\Omega\setminus\Omega_\rho$ and is diffeomorphic to either
$\mathbb{S}^3$, or $\mathbb{S}^2 \times \mathbb{S}^1$ or
$\mathbb{RP}^3\#\mathbb{RP}^3$. Note again that the number of the
compact components is finite. Let us throw away all the compact
components lying in $\Omega\setminus\Omega_\rho$ and all the
compact components with positive curvature, and then consider
those components $\Omega_j$, $1\leq j\leq k$, of $\Omega$ which
contain points of $\Omega_\rho$. (We will consider those
components of $\Omega\setminus\Omega_\rho$ consisting of capped
$\varepsilon$-horns and double $\varepsilon$-horns later). We will
perform surgical procedures, as we roughly described before, by
finding an $\varepsilon$-neck in every horn of $\Omega_j$, $1\leq
j\leq k$, and then cutting it along the middle two-sphere,
removing the horn-shaped end, and gluing back a cap.

In order to maintain the a priori assumptions \textbf{with the
same accuracy} after the surgery, we will need to find sufficient
``fine" necks in the $\varepsilon$-horns and to glue sufficient
``fine" caps. Note that $\delta>0$ will be chosen much smaller
than $\varepsilon>0$. The following lemma due to Perelman
\cite{P2} gives us the ``fine" necks in the $\varepsilon$-horns.
(At the first sight, we should also cut off all those
$\varepsilon$-tubes and $\varepsilon$-caps in the surgery
procedure. However, in general we are not able to find a ``fine"
neck in an $\varepsilon$-tube or in an $\varepsilon$-cap, and
surgeries at ``rough" $\varepsilon$-necks will certainly lose some
accuracy. If we perform surgeries at the necks with some fixed
accuracy $\varepsilon$ in the high curvature region at each
surgery time, then it is possible that the errors of surgeries may
accumulate to a certain amount so that at some later time we
cannot recognize the structure of very high curvature regions.
This prevents us from carrying out the whole process in finite
time with a finite number of steps. This is the reason why we will
only perform the surgeries at the $\varepsilon$-horns.)

\begin{lemma}[{Perelman \cite{P2}}]
Given $0<\varepsilon \leq \frac{1}{100}, 0<\delta<\varepsilon$ and
$0<T<+\infty$, there exists a radius $0<h<\delta\rho$, depending
only on $\delta$ and $r(T)$, such that if we have a solution to
the Ricci flow with surgery, with a normalized metric as initial
data, satisfying the a priori assumptions $($with accuracy
$\varepsilon),$ defined on $[0,T)$, going singular at time $T$ and
having only a finite number of surgery times on $[0,T)$, then for
each point $x$ with $h(x)=\bar{R}^{-\frac{1}{2}}(x)\leq h$ in an
$\varepsilon$-horn of $(\Omega,\bar{g}_{ij})$ with boundary in
$\Omega_\rho$, the neighborhood
$B_T(x,\delta^{-1}h(x))\triangleq\{y\in\Omega\ |\ d_T(y,x)\leq
\delta^{-1}h(x)\}$ is a strong $\delta$-neck $($i.e.,
$B_T(x,\delta^{-1}h(x)) \times [T-h^2(x),T]$ is, after scaling
with factor $h^{-2}(x)$, $\delta$-close $($in the
$C^{[\delta^{-1}]}$ topology$)$ to the corresponding subset of the
evolving standard round cylinder $\mathbb{S}^2 \times \mathbb{R}$
over the time interval $[-1,0]$ with scalar curvature $1$ at time
zero$)$.
\end{lemma}

\vskip 0.5cm
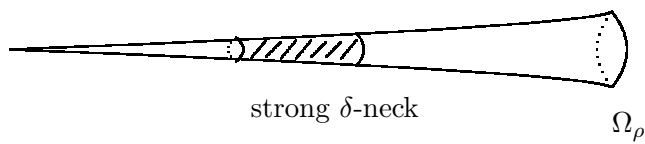
\begin{center}
\setlength{\unitlength}{2mm}
\begin{picture}(60,8)
\linethickness{0.5pt}

\qbezier[10](50,0)(48,2.5)(50,5)

\qbezier(50,0)(52,2.5)(50,5)

\qbezier[5](25,1.7)(24,2.5)(25,3.3)
\qbezier(25,1.7)(26,2.5)(25,3.3)

\qbezier(26,1.9)(26.5,2.5)(27,3.1) \qbezier(27,2)(27.5,2.5)(28,3)
\qbezier(28,1.9)(28.5,2.5)(29,3.1)
\qbezier(29,1.9)(29.5,2.5)(30,3.1)
\qbezier(30,2)(30.5,2.5)(31,3)\qbezier(32,1.9)(32.5,2.5)(33,3.1)
\qbezier(31,2)(31.5,2.5)(32,3)

\qbezier(33,1.4)(34,2.5)(33,3.6)

\qbezier(10,2.5)(48,4)(50,5)

\qbezier(10,2.5)(48,1)(50,0)

\put(26,-2){{strong $\delta$-neck}}
\put(50,-3){\textbf{$\Omega_{\rho}$}}
\end{picture}
\end{center}

\vskip 1.0cm

\begin{pf} The following proof is essentially given by Perelman
(cf. 4.3 of \cite{P2}).

We argue by contradiction. Suppose that there exists a sequence of
solutions $g^k_{ij}(\cdot,t)$, $k=1,2,\ldots$, to the Ricci flow
with surgery, satisfying the a priori assumptions (with accuracy
$\varepsilon$), defined on $[0,T)$ with limit metrics
$(\Omega^k,\bar{g}^k_{ij}), k=1, 2, \ldots$, and points $x^k$,
lying inside an $\varepsilon$-horn of $\Omega^k$ with boundary in
$\Omega_\rho^k$, and having $h(x^k)\rightarrow0$ such that the
neighborhoods $B_T(x^k,\delta^{-1}h(x^k))=\{y\in\Omega^k\ |\
d_T(y,x^k)\leq \delta^{-1}h(x^k)\}$ are not strong $\delta$-necks.

Let $\widetilde{g}^k_{ij}(\cdot,t)$ be the solutions obtained by
rescaling by the factor $\bar{R}(x^k)=h^{-2}(x^k)$ around $x^k$
and shifting the time $T$ to the new time zero. We now want to
show that a subsequence of
$\widetilde{g}^k_{ij}(\cdot,t),k=1,2,\ldots$, converges to the
evolving round cylinder, which will give a contradiction.

Note that $\widetilde{g}^k_{ij}(\cdot,t), k=1, 2, \ldots,$ are
solutions modified by surgery. So, we cannot apply Hamilton's
compactness theorem directly since it is stated only for smooth
solutions. For each (unrescaled) surgical solution
$\bar{g}^k_{ij}(\cdot,t)$, we pick a point $z^k$, with
$\bar{R}(z^k) = 2C_2^2(\varepsilon)\rho^{-2},$ in the
$\varepsilon$-horn of $(\Omega^k,\bar{g}^k_{ij})$ with boundary in
$\Omega_\rho^k$, where $C_2(\varepsilon)$ is the positive constant
in the canonical neighborhood assumption. From the definition of
$\varepsilon$-horn and the canonical neighborhood assumption, we
know that each point $x$ lying inside the $\varepsilon$-horn of
$(\Omega^k,\bar{g}^k_{ij})$ with
$d_{\bar{g}^k_{ij}}(x,\Omega_\rho^k) \geq
d_{\bar{g}^k_{ij}}(z^k,\Omega_\rho^k)$ has a strong
$\varepsilon$-neck as its canonical neighborhood. Since
$h(x^k)\rightarrow0$, each $x^k$ lies deeply inside an
$\varepsilon$-horn. Thus for each positive $A < +\infty$, the
rescaled (surgical) solutions $\widetilde{g}^k_{ij}(\cdot,t)$ with
the marked origins $x^k$ over the geodesic balls
$B_{\widetilde{g}^k_{ij}(\cdot,0)}(x^k,A)$, centered at $x^k$ of
radii $A$ (with respect to the metrics
$\widetilde{g}^k_{ij}(\cdot,0)$), will be smooth on some uniform
(size) small time intervals for all sufficiently large $k$
whenever the curvatures of the rescaled solutions
$\widetilde{g}^k_{ij}$ at $t=0$ in
$B_{\widetilde{g}^k_{ij}(\cdot,0)}(x^k,A)$ are uniformly bounded.
In such a situation, Hamilton's compactness theorem is applicable.
Then we can apply the same argument as in the second step of the
proof of Theorem 7.1.1 to conclude that for each $A<+\infty$,
there exists a positive constant $C(A)$ such that the curvatures
of the rescaled solutions $\widetilde{g}^k_{ij}(\cdot,t)$ at the
new time $0$ satisfy the estimate $$|\widetilde{R}m_k|(y,0) \leq
C(A)$$ whenever $d_{\widetilde{g}^k_{ij}(\cdot,0)}(y,x^k) \leq A$
and $k \geq 1$; otherwise we would get a piece of a non-flat
nonnegatively curved metric cone as a blow-up limit, which
contradicts Hamilton's strong maximum principle. Moreover, by
Hamilton's compactness theorem (Theorem 4.1.5), a subsequence of
the rescaled solutions $\widetilde{g}^k_{ij}(\cdot,t)$ converges
to a $C^\infty_{loc}$ limit $\widetilde{g}^\infty_{ij}(\cdot,t)$,
defined on a spacetime set which is relatively open in the half
spacetime $\{t\leq 0\}$ and contains the time slice $t=0$.

By the pinching assumption, the limit is a complete manifold with
nonnegative sectional curvature. Since $x^k$ was contained in an
$\varepsilon$-horn with boundary in $\Omega^k_\rho$ and
$h(x^k)/\rho\rightarrow0$, the limiting manifold has two ends.
Thus, by Toponogov's splitting theorem, the limiting manifold
admits a metric splitting $\Sigma^2\times \mathbb{R}$, where
$\Sigma^2$ is diffeomorphic to the two-sphere $\mathbb{S}^2$
because $x^k$ was the center of a strong $\varepsilon$-neck.

By combining with the canonical neighborhood assumption (with
accuracy $\varepsilon$), we see that the limit is defined on the
time interval $[-1,0]$ and is $\varepsilon$-close to the evolving
standard round cylinder. In particular, the scalar curvature of
the limit at time $t=-1$ is $\varepsilon$-close to $1/2$.

Since $h(x^k)/\rho\rightarrow0$, each point in the limiting
manifold at time $t=-1$ also has a strong $\varepsilon$-neck as
its canonical neighborhood. Thus the limit is defined at least on
the time interval $[-2,0]$ and the limiting manifold at time
$t=-2$ is, after rescaling, $\varepsilon$-close to the standard
round cylinder.

By using the canonical neighborhood assumption again, every point
in the limiting manifold at time $t=-2$ still has a strong
$\varepsilon$-neck as its canonical neighborhood. Also note that
the scalar curvature of the limit at $t=-2$ is not bigger than
$1/2+\varepsilon$. Thus the limit is defined at least on the time
interval $[-3,0]$ and the limiting manifold at time $t=-3$ is,
after rescaling, $\varepsilon$-close to the standard round
cylinder. By repeating this argument we prove that the limit
exists on the ancient time interval $(-\infty, 0]$.

The above argument also shows that at every time, each point of
the limit has a strong $\varepsilon$-neck as its canonical
neighborhood. This implies that the limit is $\kappa$-noncollaped
on all scales for some $\kappa>0$. Therefore, by Theorem 6.2.2,
the limit is the evolving round cylinder
$\mathbb{S}^2\times\mathbb{R}$, which gives the desired
contradiction.
\end{pf}

In the above lemma, the property that the radius $h$ depends only
on $\delta$ and the time $T$ but is independent of the surgical
solution is crucial; otherwise we will not be able to cut off
enough volume at each surgery to guarantee the number of surgeries
being finite in each finite time interval. We also remark that the
above proof actually proves a stronger result: the parabolic
region $\{(y,t)\ |\ y\in
B_T(x,\delta^{-1}h(x)),t\in[T-\delta^{-2}h^2(x),T]\}$ is, after
scaling with factor $h^{-2}(x)$, $\delta$-close (in the
$C^{[\delta^{-1}]}$ topology) to the corresponding subset of the
evolving standard round cylinder $\mathbb{S}^2 \times \mathbb{R}$
over the time interval $[-\delta^{-2},0]$ with scalar curvature
$1$ at the time zero. This fact will be used later in the proof of
Proposition 7.4.1.

We next want to construct ``fine" caps. Take a rotationally
symmetric metric on $\mathbb{R}^3$ with nonnegative sectional
curvature and positive scalar curvature such that outside some
compact set it is a semi-infinite standard round cylinder (i.e.
the metric product of a ray with the round two-sphere of scalar
curvature 1). We call such a metric on $\mathbb{R}^3$  a
\textbf{standard capped infinite cylinder}\index{standard capped
infinite cylinder}. By the short-time existence theorem of Shi
(Theorem 1.2.3), the Ricci flow with a standard capped infinite
cylinder as initial data has a complete solution on a maximal time
interval $[0,T)$ such that the curvature of the solution is
bounded on $\mathbb{R}^3 \times [0,T']$ for each $0<T'<T$.  Such a
solution is called a \textbf{standard
solution}\index{solution!standard} by Perelman \cite{P2}.

The following result, proved by Chen and the second author in
\cite{CZ05F}, gives the curvature estimate for standard solutions.
This curvature estimate in the special case, when the dimension is
three and the initial metric is rotationally symmetric, was first
claimed by Perelman in \cite{P2}.

\begin{proposition}
Let $g_{ij}$ be a complete Riemannian metric on $\mathbb{R}^{n}$
($n \geq 3$) with nonnegative curvature operator and positive
scalar curvature which is asymptotic to a round cylinder of scalar
curvature $1$ at infinity.  Then there is a complete solution
$g_{ij}({\cdot,t})$ to the Ricci flow, with $g_{ij}$ as initial
metric, which exists on the time interval $[0,\frac{n-1}{2})$, has
bounded curvature at each time $t \in [0,\frac{n-1}{2})$, and
satisfies the estimate
$$
  R(x,t)\geq \frac{C^{-1}}{\frac{n-1}{2}-t}
$$
for some $C$ depending only on the initial metric $g_{ij}$.
\end{proposition}

\begin{pf}
Since the initial metric has bounded curvature operator and a
positive lower bound on its scalar curvature, the Ricci flow has a
solution $g_{ij}(\cdot,t)$ defined on a maximal time interval
$[0,T)$ with $T<\infty$ which has bounded curvature on
$\mathbb{R}^n \times [0,T']$ for each $0<T'<T$. By Proposition
2.1.4, the solution $g_{ij}(\cdot,t)$ has nonnegative curvature
operator for all $t \in [0, T)$.  Note that the injectivity radius
of the initial metric has a positive lower bound. As we remarked
at the beginning of Section 3.4, the same proof of Perelman's no
local collapsing theorem I concludes that $g_{ij}({\cdot,t})$ is
$\kappa$-noncollapsed on all scales less than $\sqrt{T}$ for some
$\kappa>0$ depending only on the initial metric.

We will first prove the following assertion.

\medskip
{\bf Claim 1.} \ There is a positive function
$\omega:[0,\infty)\rightarrow[0,\infty)$ depending only on the
initial metric and $\kappa$ such that
$$
  R(x,t)\leq R(y,t)\omega(R(y,t)d_{t}^{2}(x,y))
$$
for all $x,y\in \mathbb{R}^{n}$, $t\in[0,T)$.

\medskip
The proof is similar to that of Theorem 6.4.3. Notice that the
initial metric has nonnegative curvature operator and its scalar
curvature satisfies the bounds \be
C^{-1}\leqslant R(x)\leqslant C 
\ee for some positive constant $C$. By the maximum principle, we
know $T\geq\frac{1}{2nC}$ and $R(x,t)\leq 2C$ for
$t\in[0,\frac{1}{4nC}]$. The assertion is clearly true for
$t\in[0,\frac{1}{4nC}]$.

Now fix $(y,t_{0})\in \mathbb{R}^{n}\times[0,T)$ with
$t_{0}\geq\frac{1}{4nC}$. Let $z$ be the closest point to $y$ with
the property $R(z,t_{0})d^{2}_{t_{0}}(z,y)=1$ (at time $t_{0}$).
Draw a shortest geodesic from $y$ to $z$ and choose a point
$\tilde{z}$ on the geodesic satisfying $d_{t_{0}}(z,\tilde{z})
=\frac{1}{4}R(z,t_{0})^{-\frac{1}{2}}$, then we have
$$
 R(x,t_{0})\leq\frac{1}{(\frac{1}{2}R(z,t_{0})^{-\frac{1}{2}})^{2}}
    \qquad \mbox{on}\ \
 B_{t_{0}}\(\tilde{z},\frac{1}{4}R(z,t_{0}\)^{-\frac{1}{2}}).
$$

Note that $R(x,t)\geqslant C^{-1}$ everywhere by the evolution
equation of the scalar curvature. Then by the Li-Yau-Hamilton
inequality (Corollary 2.5.5), for all $(x,t)\in
B_{t_{0}}(\tilde{z},\frac{1}{8nC}R(z,t_{0})^{-\frac{1}{2}})
\times[t_{0}-(\frac{1}{8nC}R(z,t_{0})^{-\frac{1}{2}})^{2},t_{0}]
$, we have
\begin{equation*}
\begin{split}
R(x,t)&  \leq \(\frac{t_{0}}{t_{0}-\(\frac{1}{8n\sqrt{C}}\)^{2}}\)
\frac{1}{\(\frac{1}{2}R(z,t_{0})^{-\frac{1}{2}}\)^{2}}\\
&  \leq \left[\frac{1}{8nC}R(z,t_{0})^{-\frac{1}{2}}\right]^{-2}.
\end{split}
\end{equation*}
Combining this with the $\kappa$-noncollapsing property, we have
$$
\Vol
\(B_{t_{0}}\(\tilde{z},\frac{1}{8nC}R(z,t_{0})^{-\frac{1}{2}}\)\)
\geq \kappa \(\frac{1}{8nC}R(z,t_{0})^{-\frac{1}{2}}\)^{n}
$$
and then
$$
\Vol \(B_{t_{0}}\(z,8R(z,t_{0})^{-\frac{1}{2}}\)\) \geq \kappa
\(\frac{1}{64nC}\)^{n}\(8R(z,t_{0})^{-\frac{1}{2}}\)^{n}.
$$
So by Theorem 6.3.3 (ii), we have
$$
R(x,t_{0})\leq C(\kappa)R(z,t_{0})\quad \text{for all }\; x\in
B_{t_{0}}\(z,4R(z,t_{0})^{-\frac{1}{2}}\).
$$
Here and in the following we denote by $C(\kappa)$ various
positive constants depending only on $\kappa,n$ and the initial
metric.

Now by the Li-Yau-Hamilton inequality (Corollary 2.5.5) and local
gradient estimate of Shi (Theorem 1.4.2), we obtain
\begin{align*}
R(x,t)\leq C(\kappa)R(z,t_{0})\; \text{ and }\;
\left|\frac{\partial}{\partial t}R\right|(x,t) \leq
C(\kappa)(R(z,t_{0}))^{2}
\end{align*}
for all $ (x,t)\in
B_{t_{0}}(z,2R(z,t_{0})^{-\frac{1}{2}}))\times[t_{0}-
(\frac{1}{8nC}R(z,t_{0})^{-\frac{1}{2}})^{2},t_{0}]$. Therefore by
combining with the Harnack estimate (Corollary 2.5.7), we obtain
\begin{align*}
R(y,t_{0})&  \geq C(\kappa)^{-1}R(z,t_{0}-C(\kappa)^{-1}R(z,t_{0})^{-1})\\
&  \geq C(\kappa)^{-2} R(z,t_{0})
\end{align*}

Consequently, we have showed that there is a constant $C(\kappa)$
such that
$$
\Vol \(B_{t_{0}}\(y,R(y,t_{0})^{-\frac{1}{2}}\)\) \geq
C(\kappa)^{-1}\(R(y,t_{0})^{-\frac{1}{2}}\)^{n}
$$
and
$$
R(x,t_{0})\leq C(\kappa)R(y,t_{0}) \quad \text{for all }\; x\in
B_{t_{0}}\(y,R(y,t_{0})^{-\frac{1}{2}}\).
$$
In general, for any $r\geq R(y,t_{0})^{-\frac{1}{2}}$, we have
$$
\Vol (B_{t_{0}}(y,r)) \geq
C(\kappa)^{-1}(r^{2}R(y,t_{0}))^{-\frac{n}{2}}r^{n}.
$$
By\; applying\; Theorem\; 6.3.3(ii)\; again,\; there\; exists\;
a\; positive\; constant $\omega(r^{2}R(y,t_{0}))$ depending only
on the constant $r^{2}R(y,t_{0})$ and $\kappa$ such that
$$
R(x,t_{0})\leq R(y,t_{0})\omega(r^{2}R(y,t_{0}))\quad \mbox{for
all }\;  x\in B_{t_{0}}\(y,\frac{1}{4}r\).
$$
This proves the desired Claim 1.

Now we study the asymptotic behavior of the solution at infinity.
For any $0<t_{0}<T$, we know that the metrics $g_{ij}(x,t)$ with
$t\in [0,t_{0}]$ has uniformly bounded curvature.  Let $x_{k}$ be
a sequence of points with $d_{0}(x_{0},x_{k})\rightarrow \infty$.
After passing to a subsequence, $g_{ij}(x,t)$ around $x_{k}$ will
converge to a solution to the Ricci flow on $\mathbb{R}\times
\mathbb{S}^{n-1}$ with round cylinder metric of scalar curvature 1
as initial data.  Denote the limit by $\tilde{g}_{ij}$. Then by
the uniqueness theorem (Theorem 1.2.4), we have
$$
\tilde{R}(x,t) =\frac{\frac{n-1}{2}}{\frac{n-1}{2}-t}\quad
\text{for all }\; t\in[0,t_{0}].
$$
It follows that $T\leq \frac{n-1}{2}$. In order to show
$T=\frac{n-1}{2}$, it suffices to prove the following assertion.

\medskip
{\bf Claim 2.} \ Suppose $T<\frac{n-1}{2}$. Fix a point $x_{0}\in
\mathbb{R}^{n} $, then there is a $\delta>0$, such that for any
$x\in M$ with $d_{0}(x,x_{0})\geq \delta^{-1}$, we have
$$
R(x,t)\leq 2C+\frac{n-1}{\frac{n-1}{2}-t}\quad \text{ for all }\;
t\in[0,T),
$$
where $C$ is the constant in (7.3.5).

\medskip
In view of Claim 1, if Claim 2 holds, then
\begin{align*}
\sup_{M^{n}\times[0,T)}R(y,t) &  \leq
\omega\(\delta^{-2}\(2C+\frac{n-1}{\frac{n-1}{2}-T}\)\)
\(2C+\frac{n-1}{\frac{n-1}{2}-T}\)\\
&  < \infty
\end{align*}
which will contradict the definition of $T$.

To show Claim 2, we argue by contradiction. Suppose for each
$\delta>0$, there is a point $(x_{\delta},t_{\delta})$ with
$0<t_{\delta}<T $ such that
$$
R(x_{\delta},t_{\delta})>2C+\frac{n-1}{\frac{n-1}{2}-t_{\delta}}\;
\mbox{ and }\; d_{0}(x_{\delta},x_{0})\geq \delta^{-1}.
$$
Let
$$
\bar{t}_{\delta}=\sup\left\{t\; \Big|\; \sup_{M^{n} \setminus
B_{0}(x_{0},\delta^{-1})}
R(y,t)<2C+\frac{n-1}{\frac{n-1}{2}-t}\right\}.
$$
Since $\lim\limits_{d_{0}(y,x_{0})\rightarrow\infty}R(y,t)
=\frac{\frac{n-1}{2}}{\frac{n-1}{2}-t}$ and
$\sup_{M\times[0,\frac{1}{4nC}]}R(y,t)\leq 2C$, we know
$\frac{1}{4nC}\leq \bar{t}_{\delta}\leq t_{\delta}$ and there is a
$\bar{x}_{\delta}$ such that
$d_{0}(x_{0},\bar{x}_{\delta})\geq\delta^{-1}$ and
$R(\bar{x}_{\delta},\bar{t}_{\delta})
=2C+\frac{n-1}{\frac{n-1}{2}-\bar{t}_{\delta}}.$ By Claim 1 and
Hamilton's compactness theorem (Theorem 4.1.5), for $\delta\
\rightarrow 0$ and after taking a subsequence, the metrics
$g_{ij}(x,t)$ on $B_{0}(\bar{x}_{\delta},\frac{\delta^{-1}}{2})$
over the time interval $[0,\bar{t}_{\delta}]$ will converge to a
solution $\tilde{g}_{ij}$ on $\mathbb{R}\times \mathbb{S}^{n-1}$
with the standard metric of scalar curvature 1 as initial data
over the time interval $[0,\bar{t}_{\infty}]$, and its scalar
curvature satisfies
\begin{align*}
\tilde{R}(\bar{x}_{\infty},\bar{t}_{\infty})
&  = 2C+\frac{n-1}{\frac{n-1}{2}-\bar{t}_{\infty}},\\
\tilde{R}(x,t)&  \leqslant
2C+\frac{n-1}{\frac{n-1}{2}-\bar{t}_{\infty}},
 \quad \mbox{for all }\;  t\in [0,\bar{t}_{\infty}],
\end{align*}
where $(\bar{x}_{\infty},\bar{t}_{\infty})$ is the limit of
$(\bar{x}_{\delta},\bar{t}_{\delta})$. On the other hand, by the
uniqueness theorem (Theorem 1.2.4) again, we know
$$
\tilde{R}(\bar{x}_{\infty},\bar{t}_{\infty})
=\frac{\frac{n-1}{2}}{\frac{n-1}{2}-\bar{t}_{\infty}}
$$
which is a contradiction.  Hence we have proved Claim 2 and then
have verified $T=\frac{n-1}{2}$.

Now we are ready to show \be R(x,t)\geq
\frac{\tilde{C}^{-1}}{\frac{n-1}{2}-t}, \quad \text{for all }\;
(x,t)\in \mathbb{R}^{n}\times\Big[0,\frac{n-1}{2}\Big),
\ee for some positive constant $\tilde{C}$ depending only on the
initial metric.

For any $(x,t)\in \mathbb{R}^{n}\times[0,\frac{n-1}{2})$, by Claim
1 and $\kappa$-noncollapsing, there is a constant $C(\kappa)>0$
such that
$$
\Vol_{t}(B_{t}(x,{R(x,t)}^{-\frac{1}{2}}))\geq
C(\kappa)^{-1}({R(x,t)}^{-\frac{1}{2}})^{n}.
$$
Then by the well-known volume estimate of Calabi-Yau (see for
example \cite{Y76} or \cite{ScY}) for complete manifolds with
$\Ric\geq 0$, for any $a\geq 1$, we have
$$
\Vol_{t}(B_{t}(x,a{R(x,t)}^{-\frac{1}{2}}))\geq
C(\kappa)^{-1}\frac{a}{8n}({R(x,t)}^{-\frac{1}{2}})^{n}.
$$
On the other hand, since $(\mathbb{R}^{n},g_{ij}(\cdot,t))$ is
asymptotic to a cylinder of scalar curvature
${(\frac{n-1}{2})}/{(\frac{n-1}{2}-t)}$, for sufficiently large
$a>0$, we have
$$
\Vol_{t}\(B_{t}\(x,a\sqrt{\frac{n-1}{2}-t}\)\)\leq
C(n)a\(\frac{n-1}{2}-t\)^{\frac{n}{2}}.
$$
Combining the two inequalities, for  all sufficiently large $a$,
we have:
\begin{align*}
C(n)a\(\frac{n-1}{2}-t\)^{\frac{n}{2}}&  \geq
\Vol_{t}\(B_{t}\(x,a\(\frac{\sqrt{\frac{n-1}{2}
-t}}{R(x,t)^{-\frac{1}{2}}}\){R(x,t)}^{-\frac{1}{2}}\)\)\\
&  \geq C(\kappa)^{-1}\frac{a}{8n}\(\frac{\sqrt{\frac{n-1}{2}
-t}}{{R(x,t)}^{-\frac{1}{2}}}\) \({R(x,t)}^{-\frac{1}{2}}\)^{n},
\end{align*}
which gives the desired estimate (7.3.6). Therefore the proof of
the proposition is complete.
\end{pf}

We now fix a standard capped infinite cylinder for dimension $n=3$
as follows.  Consider the semi-infinite standard round cylinder
$N_0 = \mathbb{S}^2 \times (-\infty,4)$ with the metric $g_0$ of
scalar curvature 1. Denote by $z$ the coordinate of the second
factor $(-\infty,4)$. Let $f$ be a smooth nondecreasing convex
function on $(-\infty,4)$ defined by \be
   \left\{
   \begin{array}{lll}
    f(z) = 0, \ \ \ z\leq 0,
         \\[3mm]
    f(z) = ce^{-\frac{P}{z}}, \ \ \ z \in (0,3], \\[3mm]
    f(z) \mbox{ is strictly convex on } z \in [3,3.9], \\[3mm]
    f(z) = -\frac{1}{2}\log(16-z^2), \ \ \ z \in [3.9,4),
\end{array}
\right.         
\ee where the (small) constant $c>0$ and (big) constant $P>0$ will
be determined later. Let us replace the standard metric $g_0$ on
the portion $\mathbb{S}^2 \times [0,4)$ of the semi-infinite
cylinder by $\hat{g} = e^{-2f}g_0$. Then the resulting metric
$\hat{g}$ will be smooth on $\mathbb{R}^3$ obtained by adding a
point to $\mathbb{S}^2 \times (-\infty,4)$ at $z=4$. We denote by
$C(c,P) = (\mathbb{R}^3,\hat{g})$. Clearly, $C(c,P)$ is a standard
capped infinite cylinder.

We next use a compact portion of the standard capped infinite
cylinder $C(c,P)$ and the $\delta$-neck obtained in Lemma 7.3.2 to
perform the following surgery procedure due to Hamilton
\cite{Ha97}.

Consider the metric $\bar{g}$ at the maximal time $T<+\infty$.
Take an $\varepsilon$-horn with boundary in $\Omega_\rho$. By
Lemma 7.3.2, there exists a $\delta$-neck $N$ of radius
$0<h<\delta \rho$ in the $\varepsilon$-horn. By definition,
$(N,h^{-2}\bar{g})$ is $\delta$-close (in the $C^{[\delta^{-1}]}$
topology) to the standard round neck $\mathbb{S}^2\times
\mathbb{I}$ of scalar curvature 1 with
$\mathbb{I}=(-\delta^{-1},\delta^{-1})$. Using the parameter $z
\in \mathbb{I}$, we see the above function $f$ is defined on the
$\delta$-neck $N$.

Let\; us\; cut\; the\; $\delta$-neck\; $N$\; along\; the\;
middle\; (topological)\; two-sphere $N\bigcap\{z=0\}$. Without
loss of generality, we may assume that the right hand half portion
$N\bigcap\{z\geq 0\}$ is contained in the horn-shaped end. Let
$\varphi$ be a smooth bump function with $\varphi = 1$ for
$z\leq2$, and $\varphi = 0$ for $z\geq3$. Construct a new metric
$\tilde{g}$ on a (topological) three-ball $\mathbb{B}^3$ as
follows \be
   \tilde{g} =  \begin{cases}
    \bar{g}, &\quad z= 0, \\[1mm]
    e^{-2f}\bar{g},&\quad  z \in [0,2], \\[1mm]
    \varphi e^{-2f}\bar{g} + (1-\varphi)e^{-2f}h^2g_0,
&\quad z \in [2,3], \\[1mm]
    h^2e^{-2f}g_0, &\quad z\in [3,4].
\end{cases}  
\ee The surgery is to replace the horn-shaped end by the cap
$(\mathbb{B}^3,\tilde{g})$. We call such surgery procedure a
\textbf{$\delta$-cutoff surgery}\index{$\delta$-cutoff surgery}.

The following lemma determines the constants $c$ and $P$ in the
$\delta$-cutoff surgery so that the pinching assumption is
preserved under the surgery (cf 4.4 of Perelman \cite{P2}).

\begin{lemma}[{Justification of the pinching
assumption}\index{justification of the pinching assumption}] There
are universal positive constants $\delta_0$, $c_0$ and $P_0$ such
that if we take a $\delta$-cutoff surgery at a $\delta$-neck of
radius $h$ at time $T$ with $\delta \leq \delta_0$ and $h^{-2}
\geq 2e^2\log(1+T)$, then we can choose $c=c_0$ and $P=P_0$ in the
definition of $f(z)$ such that after the surgery, there still
holds the pinching condition \be \tilde{R} \geq
(-\tilde{\nu})[\log(-\tilde{\nu}) + \log(1+T) -3]
\ee whenever $\tilde{\nu} < 0$, where $\tilde{R}$ is the scalar
curvature of the metric $\tilde{g}$ and $\tilde{\nu}$ is the least
eigenvalue of the curvature operator of $\tilde{g}$. Moreover,
after the surgery, any metric ball of radius
$\delta^{-\frac{1}{2}}h$ with center near the tip $($i.e., the
origin of the attached cap$)$ is, after scaling with factor
$h^{-2}$, $\delta^{\frac{1}{2}}$-close $($in the
$C^{[\delta^{-\frac{1}{2}}]}$ topology$)$ to the corresponding
ball of the standard capped infinite cylinder $C(c_0,P_0)$.
\end{lemma}

\begin{pf}
First, we consider the metric $\tilde{g}$ on the portion $\{0 \leq
z \leq 2\}$. Under the conformal change $\tilde{g}=
e^{-2f}\bar{g}$, the curvature tensor $\tilde{R}_{ijkl}$ is given
by
\begin{align*}
\tilde{R}_{ijkl} &= e^{-2f}\Big[\bar{R}_{ijkl} + |\bar{\nabla}
f|^2(\bar{g}_{il}\bar{g}_{jk}-\bar{g}_{ik}\bar{g}_{jl}) + (f_{ik}
+ f_if_k)\bar{g}_{jl} \\
&\qquad + (f_{jl} + f_jf_l)\bar{g}_{ik} - (f_{il} +
f_if_l)\bar{g}_{jk} - (f_{jk} + f_jf_k)\bar{g}_{il}\Big].
\end{align*}
If $\{\bar{F}_a = \bar{F}^i_a\frac{\partial}{\partial x^i}\}$ is
an orthonormal frame for $\bar{g}_{ij}$, then $\{\tilde{F}_a =
e^f\bar{F}_a = \tilde{F}^i_a\frac{\partial}{\partial x^i}\}$ is an
orthonormal frame for $\tilde{g}_{ij}$. Let
\begin{align*}
\bar{R}_{abcd}
&= \bar{R}_{ijkl}\bar{F}^i_a\bar{F}^j_b\bar{F}^k_c\bar{F}^l_d,\\
\tilde{R}_{abcd}
&=\tilde{R}_{ijkl}\tilde{F}^i_a\tilde{F}^j_b\tilde{F}^k_c\tilde{F}^l_d,
\end{align*}
then
\begin{align}
\tilde{R}_{abcd} &= e^{2f}\Big[\bar{R}_{abcd} + |\bar{\nabla}
f|^2(\delta_{ad}\delta_{bc}-\delta_{ac}\delta_{bd}) + (f_{ac} +
f_af_c)\delta_{bd} \\
&\qquad + (f_{bd} + f_bf_d)\delta_{ac} - (f_{ad} +
f_af_d)\delta_{bc} - (f_{bc} + f_bf_c)\delta_{ad}\Big], \nn
\end{align}
and \be \tilde{R} = e^{2f}(\bar{R} + 4\bar{\triangle} f
-2|\bar{\nabla}f|^2). 
\ee Since
$$
\frac{df}{dz}= ce^{-\frac{P}{z}}\frac{P}{z^2}, \ \
\frac{d^2f}{dz^2} =
ce^{-\frac{P}{z}}\(\frac{P^2}{z^4}-\frac{2P}{z^3}\),
$$
then for any small $\theta >0$, we may choose $c>0$ small and
$P>0$ large such that for $z \in [0,3]$, we have \be |e^{2f} - 1|
+ \left|\frac{df}{dz}\right| + \left|\(\frac{df}{dz}\)^2\right| <
\theta \frac{d^2f}{dz^2}, \ \ \frac{d^2f}{dz^2} < \theta.
\ee

On the other hand, by the definition of $\delta$-neck of radius
$h$, we have
$$
|\bar{g} - h^2g_0|_{g_0} < \delta h^2,
$$
$$
|\overset{o}{\nabla^j}\bar{g}|_{g_0} < \delta h^2, \; \mbox{ for
}\; 1 \leq j\leq [\delta^{-1}],
$$
where $g_0$ is the standard metric of the round cylinder
$\mathbb{S}^2 \times \mathbb{R}$. Note that in three dimensions,
we can choose the orthonormal frame
$\{\bar{F}_1,\bar{F}_2,\bar{F}_3 \}$ for the metric $\bar{g}$ so
that its curvature operator is diagonal in the orthonormal frame
$\{ \sqrt{2}\bar{F}_2\wedge \bar{F}_3, \sqrt{2}\bar{F}_3\wedge
\bar{F}_1, \sqrt{2}\bar{F}_1\wedge \bar{F}_2 \}$ with eigenvalues
$\bar{\nu} \leq \bar{\mu} \leq \bar{\lambda}$ and
$$
\bar{\nu}=2\bar{R}_{2323}, \ \ \bar{\mu}=2\bar{R}_{3131}, \ \
\bar{\lambda}=2\bar{R}_{1212}.
$$
Since $h^{-2}\bar{g}$ is $\delta$-close to the standard round
cylinder metric $g_0$ on the $\delta$-neck, we have \be \left\{
\begin{array}{lll}
|\bar{R}_{3131}| + |\bar{R}_{2323}| < \delta^{\frac{7}{8}}h^{-2}, \\[1mm]
|\bar{R}_{1212} - \frac{1}{2}h^{-2}| < \delta^{\frac{7}{8}}h^{-2}, \\[1mm]
|\bar{F}_3 - h^{-1}\frac{\partial}{\partial z}|_{g_0} <
\delta^{\frac{7}{8}}h^{-1},\\[1mm]
\end{array}
\right. 
\ee for suitably small $\delta > 0$. Since $\bar{\nabla}_az =
\bar{\nabla} z(\bar{F}_a)$ and $\bar{\nabla}_a\bar{\nabla}_b z =
\bar{\nabla}^2 z(\bar{F}_a,\bar{F}_b)$, it follows that
\begin{align*}
|\bar{\nabla}_3z - h^{-1}| &< \delta^{\frac{7}{8}}h^{-1},\\
|\bar{\nabla}_1z| + |\bar{\nabla}_2z|
&<\delta^{\frac{7}{8}}h^{-1},
\end{align*}
and
$$
|\bar{\nabla}_a\bar{\nabla}_b z| < \delta^{\frac{7}{8}}h^{-2},\;
\mbox{ for }\; 1 \leq a, b \leq 3.
$$
By combining with
$$
\bar{\nabla}_af = \frac{df}{dz}\bar{\nabla}_az, \ \
\bar{\nabla}_a\bar{\nabla}_bf =
\frac{df}{dz}\bar{\nabla}_a\bar{\nabla}_bz +
\frac{d^2f}{dz^2}\bar{\nabla}_az\bar{\nabla}_bz
$$
and (7.3.12), we get \be
\begin{cases}
|\bar{\nabla}_af| < 2\theta h^{-1}\frac{d^2f}{dz^2}, &
\mbox{for }\; 1 \leq a \leq 3, \\[1mm]
   |\bar{\nabla}_a\bar{\nabla}_b f| <
\delta^{\frac{3}{4}}h^{-2}\frac{d^2f}{dz^2},
& \mbox{unless }\;  a=b=3, \\[1mm]
|\bar{\nabla}_3\bar{\nabla}_3f - h^{-2}\frac{d^2f}{dz^2}| <
\delta^{\frac{3}{4}}h^{-2}\frac{d^2f}{dz^2}.
\end{cases} 
\ee By combining (7.3.10) and (7.3.14), we have \be \left\{
\begin{array}{lll}
\tilde{R}_{1212} \geq \bar{R}_{1212} - (\theta^{\frac{1}{2}}
+ \delta^{\frac{5}{8}})h^{-2}\frac{d^2f}{dz^2}, \\[1mm]
\tilde{R}_{3131} \geq \bar{R}_{3131} + (1 - \theta^{\frac{1}{2}}
- \delta^{\frac{5}{8}})h^{-2}\frac{d^2f}{dz^2},  \\[1mm]
\tilde{R}_{2323} \geq \bar{R}_{2323} + (1 - \theta^{\frac{1}{2}}
- \delta^{\frac{5}{8}})h^{-2}\frac{d^2f}{dz^2}, \\[1mm]
|\tilde{R}_{abcd}| \leq (\theta^{\frac{1}{2}} +
\delta^{\frac{5}{8}})h^{-2}\frac{d^2f}{dz^2}, \qquad
\mbox{otherwise},
\end{array}
\right. 
\ee where $\theta$ and $\delta$ are suitably small. Then it
follows that
$$
\tilde{R} \geq \bar{R} + [4-6(\theta^{\frac{1}{3}} +
\delta^{\frac{1}{2}})]h^{-2}\frac{d^2f}{dz^2},
$$
$$
-\tilde{\nu} \leq -\bar{\nu} - [2-2(\theta^{\frac{1}{3}} +
\delta^{\frac{1}{2}})]h^{-2}\frac{d^2f}{dz^2},
$$
for suitably small $\theta$ and $\delta$.

If $0<-\tilde{\nu} \leq e^2$, then by the assumption that $h^{-2}
\geq 2e^2\log(1+T)$, we have
\begin{align*}
\tilde{R} &\geq \bar{R}\\
&  \geq \frac{1}{2}h^{-2} \\
&  \geq e^2\log(1+T)\\[2mm]
&  \geq (-\tilde{\nu})[\log(-\tilde{\nu}) + \log(1+T) -3].
\end{align*}
While if $-\tilde{\nu} > e^2$, then by the pinching estimate of
$\bar{g}$, we have
\begin{align*}
\tilde{R}&  \geq \bar{R}\\
&  \geq (-\bar{\nu})[\log(-\bar{\nu}) + \log(1+T) -3] \\
&  \geq (-\tilde{\nu})[\log(-\tilde{\nu}) + \log(1+T) -3].
\end{align*}
So we have verified the pinching condition on the portion $\{0
\leq z \leq 2\}$.

Next, we consider the metric $\tilde{g}$ on the portion $\{2 \leq
z \leq 4\}$. Let $\theta$ be a fixed suitably small positive
number. Then the constant $c=c_0$ and $P=P_0$ are fixed. So $\zeta
= \min_{z \in [1,4]} \frac{d^2f}{dz^2}
> 0$ is also fixed. By the same argument as in the derivation of (7.3.15)
from (7.3.10), we see that the curvature of the metric $\hat{g} =
e^{-2f}g_0$ of the standard capped infinite cylinder $C(c_0,P_0)$
on the portion $\{ 1\leq z \leq 4\}$ is bounded from below by
$\frac{2}{3}\zeta >0$.  Since $h^{-2}\bar{g}$ is $\delta$-close to
the standard round metric $g_0$, the metric $h^{-2}\tilde{g}$
defined by (7.3.8) is clearly $\delta^{\frac{3}{4}}$-close to the
metric $\hat{g} = e^{-2f}g_0$ of the standard capped infinite
cylinder on the portion $\{ 1\leq z \leq 4 \}$. Thus as $\delta$
is sufficiently small, the curvature operator of $\tilde{g}$ on
the portion $\{2 \leq z \leq 4\}$ is positive. Hence the pinching
condition (7.3.9) holds trivially on the portion $\{2 \leq z \leq
4\}$.

The last statement in Lemma 7.3.4 is obvious from the definition
(7.3.8).
\mbox{ \ \ }\end{pf}

Recall from Lemma 7.3.2 that the $\delta$-necks at a time $t>0$,
where we performed Hamilton's surgeries, have their radii $0 < h <
\delta\rho =\delta^2r(t)$. Without loss of generality, we may
assume the positive nonincreasing function $r(t)$ in the
definition of the canonical neighborhood assumption is less than
$1$ and the universal constant $\delta_0$ in Lemma 7.3.4 is also
less than $1$. We define a positive function $\bar{\delta}(t)$ by
\be \bar{\delta}(t) = \min\left\{\frac{1}{{2e^2\log(1+t)}},
\delta_0\right\} \quad \mbox{for } \; t\in [0,+\infty). 
\ee

{}From now on, we always assume $0 < \delta < \bar{\delta}(t)$ for
any $\delta$-cutoff surgery at time $t>0$ and assume $c=c_0$ and
$P=P_0$. As a result, the standard capped infinite cylinder and
the standard solution are also fixed. The following lemma, which
was claimed by Perelman in \cite{P2}, gives the canonical
neighborhood structure for the fixed standard solution.

\begin{lemma}
Let $g_{ij}(x,t)$ be the above fixed standard solution to the
Ricci flow on $\mathbb{R}^3\times [0,1)$.  Then for any
$\varepsilon>0$, there is a positive constant $C(\varepsilon)$
such that each point $(x,t) \in \mathbb{R}^3\times [0,1)$ has an
open neighborhood $B$, with $B_t(x,r)\subset B \subset B_t(x,2r)$
for some $0<r<C(\varepsilon)R(x,t)^{-\frac{1}{2}}$, which falls
into one of the following two categories: either
\begin{itemize}
\item[(a)] $B$ is  an $\varepsilon$-cap, or \item[(b)] $B$ is an
$\varepsilon$-neck and it is the slice at the time $t$ of the
parabolic neighborhood
$B_t(x,\varepsilon^{-1}R(x,t)^{-\frac{1}{2}}) \times [t-\min
\{R(x,t)^{-1},t\},t]$, on which the standard solution is, after
scaling with the factor $R(x,t)$ and shifting the time $t$ to
zero, $\varepsilon$-close $($in the $C^{[\varepsilon^{-1}]}$
topology$)$ to the corresponding subset of the evolving standard
cylinder $\mathbb{S}^2 \times \mathbb{R}$ over the time interval
$[-\min \{tR(x,t),1\},0]$ with scalar curvature $1$ at the time
zero.
\end{itemize}
\end{lemma}

\begin{pf}
The proof of the lemma is reduced to two assertions. We now state
and prove the first assertion which takes care of those points
with times close to $1$.

\begin{assertion}
For any $\varepsilon>0$, there is a positive number
$\theta=\theta(\varepsilon)$ with $0<\theta<1$ such that for any
$(x_0,t_0)\in \mathbb{R}^3\times [\theta,1)$, the standard
solution on the parabolic neighborhood
$B_{t_0}(x,\varepsilon^{-1}R(x_0,t_0)^{-\frac{1}{2}}) \times
[t_0-\varepsilon^{-2}R(x_0,t_0)^{-1},t_0]$ is well-defined and is,
after scaling with the factor $R(x_0,t_0)$, $\varepsilon$-close
(in the $C^{[\varepsilon^{-1}]}$ topology) to the corresponding
subset of some orientable ancient $\kappa$-solution.
\end{assertion}

We argue by contradiction. Suppose Assertion 1 is not true, then
there exist $\bar\varepsilon>0$ and a sequence of points
$(x_{k},t_k)$ with $t_{k}\rightarrow 1$, such that for each $k$,
the standard solution on the parabolic neighborhood
$$
B_{t_k}(x_k,\bar\varepsilon^{-1}R(x_k,t_k)^{-\frac{1}{2}}) \times
[t_k-\bar\varepsilon^{-2}R(x_k,t_k)^{-1},t_k]
$$
is not, after scaling by the factor $R(x_k,t_k)$,
$\bar\varepsilon$-close to the corresponding subset of any ancient
$\kappa$-solution. Note that by Proposition 7.3.3, there is a
constant $C>0$ (depending only on the initial metric, hence it is
universal ) such that $R(x,t)\geq C^{-1}/(1-t)$. This implies
$$
\bar\varepsilon^{-2}R(x_k,t_k)^{-1} \leq
C\bar\varepsilon^{-2}(1-t_k)<t_k ,
$$
and\; then\; the\; standard\; solution\; on\; the\; parabolic\;
neighborhood\; $B_{t_k}(x_k,$
$\bar\varepsilon^{-1}R(x_k,t_k)^{-\frac{1}{2}}) \times
[t_k-\bar\varepsilon^{-2}R(x_k,t_k)^{-1},t_k]$ is well-defined for
$k$ large. By Claim 1 in the proof of Proposition 7.3.3, there is
a positive function $\omega:[0,\infty)\rightarrow[0,\infty)$ such
that
$$
  R(x,t_k)\leq R(x_k,t_k)\omega(R(x_k,t_k)d_{t_k}^{2}(x,x_k))
$$
for all $x\in \mathbb{R}^3$.  Now by scaling the standard solution
$g_{ij}(\cdot,t)$ around $x_k$ with the factor $R(x_k,t_k)$ and
shifting the time $t_k$ to zero, we get a sequence of the rescaled
solutions $
\tilde{g}^{k}_{ij}(x,\tilde{t})=R(x_k,t_k)g_{ij}(x,t_k+\tilde{t}/R(x_k,t_k))$
to the Ricci flow defined on $\mathbb{R}^3$ with $\tilde{t}\in
[-R(x_k,t_k)t_k,0]$. We denote the scalar curvature and the
distance of the rescaled metric $\tilde{g}^{k}_{ij}$ by
$\tilde{R}^{k}$ and $\tilde{d}$. By combining with Claim 1 in the
proof of Proposition 7.3.3 and the Li-Yau-Hamilton inequality, we
get
\begin{align*}
\tilde{R}^{k}(x,0)& \leq \omega(\tilde{d}_{0}^{2}(x,x_k))\\
\tilde{R}^{k}(x,\tilde{t})& \leq
\frac{R(x_k,t_k)t_k}{\tilde{t}+R(x_k,t_k)t_k}
\omega(\tilde{d}_{0}^{2}(x,x_k))
\end{align*}
for any $x\in \mathbb{R}^3$ and $\tilde{t}\in (-R(x_k,t_k)t_k,0]$.
Note that $R(x_k,t_k)t_k\rightarrow \infty$ by Proposition 7.3.3.
We have shown in the proof of Proposition 7.3.3 that the standard
solution is $\kappa$-noncollapsed on all scales less than $1$ for
some $\kappa>0$. Then from the $\kappa$-noncollapsing property,
the above curvature estimates and Hamilton's compactness theorem,
we know $\tilde{g}^{k}_{ij}(x,\tilde{t})$ has a convergent
subsequence (as $k\rightarrow\infty$) whose limit is an ancient,
$\kappa$-noncollapsed, complete and orientable solution with
nonnegative curvature operator. This limit must have bounded
curvature by the same proof of Step 3 in the proof of Theorem
7.1.1. This gives a contradiction. Hence Assertion 1 is proved.

We now fix the constant $\theta(\varepsilon)$ obtained in
Assertion 1. Let $O$ be the tip of the standard capped infinite
cylinder $\mathbb{R}^3$ (it is rotationally symmetric about $O$ at
time $0$, and it remains so as $t>0$ by the uniqueness Theorem
1.2.4).

\begin{assertion}
There are constants $B_1(\varepsilon)$ and $B_2(\varepsilon)$
depending only on $\varepsilon$, such that if $(x_0,t_0)\in
\mathbb{R}^3\times [0,\theta)$ with $d_{t_0}(x_0,O)\leq
B_{1}(\varepsilon)$, then there is a $0<r<B_2(\varepsilon)$ such
that $B_{t_0}(x_0,r)$ is an $\varepsilon$-cap; if $(x_0,t_0)\in
\mathbb{R}^3\times [0,\theta)$ with $d_{t_0}(x_0,O)\geq
B_1(\varepsilon)$, then the parabolic neighborhood
$B_{t_0}(x_0,\varepsilon^{-1}R(x_0,t_0)^{-\frac{1}{2}})$ $\times
[t_0-\min \{R(x_0,t_0)^{-1},t_0\},t_0]$ is after scaling with the
factor $R(x_0,t_0)$ and shifting the time $t_0$ to zero,
$\varepsilon$-close (in the $C^{[\varepsilon^{-1}]}$ topology) to
the corresponding subset of the evolving standard cylinder
$\mathbb{S}^2 \times \mathbb{R}$ over the time interval $[-\min
\{t_0R(x_0,t_0),1\},0]$ with scalar curvature $1$ at time zero.
\end{assertion}

Since the standard solution exists on the time interval $[0,1)$,
there is a constant $B_{0}(\varepsilon)$ such that the curvatures
on $[0,\theta(\varepsilon)]$ are uniformly bounded by
$B_{0}(\varepsilon)$. This implies that the metrics in
$[0,\theta(\varepsilon)]$ are equivalent.  Note that the initial
metric is asymptotic to the standard capped infinite cylinder. For
any sequence of points $x_{k}$ with $d_{0}(O,x_{k})\rightarrow
\infty$, after passing to a subsequence, $g_{ij}(x,t)$ around
$x_{k}$ will converge to a solution to the Ricci flow on
$\mathbb{R}\times \mathbb{S}^{2}$ with round cylinder metric of
scalar curvature 1 as initial data. By the uniqueness theorem
(Theorem 1.2.4), the limit solution must be the standard evolving
round cylinder.  This implies that there is a constant
$B_{1}(\varepsilon)>0$ depending on $\varepsilon$ such that for
any $(x_0,t_0)$ with $t_0\leq \theta(\varepsilon)$ and
$d_{t_0}(x_0,O)\geq B_{1}(\varepsilon)$, the standard solution on
the parabolic neighborhood
$B_{t_0}(x_0,\varepsilon^{-1}R(x_0,t_0)^{-\frac{1}{2}}) \times
[t_0-\min \{R(x_0,t_0)^{-1},t_0\},t_0]$ is, after scaling with the
factor $R(x_0,t_0)$, $\varepsilon$-close to the corresponding
subset of the evolving round cylinder. Since the solution is
rotationally symmetric around $O$, the cap neighborhood structures
of those points $x_0$ with $d_{t_0}(x_0,O)\leq B_{1}(\varepsilon)$
follow directly. Hence Assertion 2 is proved.

The combination of these two assertions proves the lemma.
\end{pf}

Since there are only a finite number of horns with the other end
connected to $\Omega_\rho$, \emph{we perform only a finite number
of such $\delta$-cutoff surgeries at time $T$}. Besides those
horns, there could be capped horns and double horns which lie in
$\Omega \setminus \Omega_\rho$. As explained before, they are
connected to form tubes or capped tubes at any time slightly
before $T$. So \emph{we can regard the capped horns and double
horns (of $\Omega \setminus \Omega_\rho$) to be extinct and throw
them away at time $T$}. We only need to remember that the
connected sums were broken there. Remember that \emph{we have
thrown away all compact components, either lying in $\Omega
\setminus\Omega_\rho$ or with positive sectional curvature}, each
of which is diffeomorphic to either $\mathbb{S}^3$, or a metric
quotient of $\mathbb{S}^3$, or $\mathbb{S}^2 \times \mathbb{S}^1$
or $\mathbb{RP}^3\#\mathbb{RP}^3$. So we have also removed a
finite number of copies of $\mathbb{S}^3$, or metric quotients of
$\mathbb{S}^3$, or $\mathbb{S}^2 \times \mathbb{S}^1$ or
$\mathbb{RP}^3\#\mathbb{RP}^3$ at the time $T$. Let us agree to
\textbf{declare extinct every compact component either with
positive sectional curvature or lying in $\Omega
\setminus\Omega_\rho$}; in particular, this allows us to exclude
the components with positive sectional curvature from the list of
canonical neighborhoods.

\medskip
\emph{In summary, our surgery at time $T$ consists of the
following four procedures:}

(1) \emph{perform $\delta$-cutoff surgeries for all
$\varepsilon$-horns, whose other ends are  connected to
$\Omega_\rho$;}

(2) \emph{declare extinct every compact component which has
positive sectional curvature;}

(3) \emph{throw away all capped horns and double horns lying in
$\Omega \setminus \Omega_\rho$;}

(4) \emph{declare extinct all compact components lying in $\Omega
\setminus\Omega_\rho$.}

\medskip \noindent
({\it In Sections $7.6$ and $7.7,$ we will add one more procedure
by declaring extinct every compact component which has nonnegative
scalar curvature.})

\medskip
 By Lemma 7.3.4, after performing
surgeries at time $T$, the pinching assumption (7.3.3) still holds
for the surgically modified manifold. With this surgically
modified manifold (possibly disconnected) as initial data, we now
continue our solution under the Ricci flow until it becomes
singular again at some time $T'(>T)$. Therefore, we have extended
the solution to the Ricci flow with surgery, originally defined on
$[0,T)$ with $T<+\infty$, to the new time interval $[0,T')$ with
$T'>T$. By the proof of Theorem 5.3.2, we see that the solution to
the Ricci flow with surgery also satisfies the pinching assumption
on $[0,T')$. It remains to verify the canonical neighborhood
assumption (with accuracy $\varepsilon$) for the solution on the
time interval $[T,T')$ and to prove that this extension procedure
works indefinitely (unless it becomes extinct at some finite time)
and that there exists at most a finite number of surgeries at
every finite time interval. We leave these arguments to the next
section.

Before we end this section, we check the following two results of
Perelman in \cite{P2} which will be used in the next section to
estimate the Li-Yau-Perelman distance of space-time curves which
stretch to surgery regions. The proofs are basically given by
Perelman (cf. 4.5 and 4.6 of \cite{P2}).

\begin{lemma}[{Perelman \cite{P2}}]
For any $0<\varepsilon \leq 1/100$, $1<A<+\infty$ and
$0<\theta<1$, one can find
$\bar{\delta}=\bar{\delta}(A,\theta,\varepsilon)$ with the
following property. Suppose we have a solution to the Ricci flow
with surgery which satisfies the a priori assumptions $($with
accuracy $\varepsilon)$ on $[0,T]$ and is obtained from a compact
orientable three-manifold by a finite number of $\delta$-cutoff
surgeries with each $\delta<\bar{\delta}$. Suppose we have a
cutoff surgery at time $T_0\in(0,T)$, let $x_0$ be any fixed point
on the gluing caps $($i.e., the regions affected by the cutoff
surgeries at time $T_0),$ and let $T_1=min\{T,T_0+\theta h^2\}$,
where $h$ is the cutoff radius around $x_0$ obtained in Lemma
$7.3.2.$ Then either
\begin{itemize}
\item[(i)] the solution is defined on
$P(x_0,T_0,Ah,T_1-T_0)\triangleq\{(x,t)\ |\ x\in B_t(x_0,Ah),\\
 t\in[T_0,T_1]\}$ and is, after scaling with factor $h^{-2}$ and
shifting time $T_0$ to zero, $A^{-1}$-close to a corresponding
subset of the standard solution, or \item[(ii)] the assertion (i)
holds with $T_1$ replaced by some time $t^+\in(T_0,T_1)$, where
$t^+$ is a surgery time; moreover, for each point in
$B_{T_0}(x_0,Ah)$, the solution is defined for $t\in[T_0,t^+)$ but
is not defined past $t^+$ $($i.e., the whole ball
$B_{T_0}(x_0,Ah)$ is cut off at the time $t^+).$
\end{itemize}
\end{lemma}

\begin{pf}
Let $Q$ be the maximum of the scalar curvature of the standard
solution in the time interval $[0,\theta]$ and choose a large
positive integer $N$ so that $\Delta
t=\frac{(T_1-T_0)}{N}<\varepsilon\eta^{-1}Q^{-1}h^2$, where the
positive constant $\eta$ is given in the canonical neighborhood
assumption. Set $t_k=T_0+k\Delta t$, $k=0, 1, \ldots, N$.

{}From Lemma 7.3.4, the geodesic ball $B_{T_0}(x_0,A_0h)$ at time
$T_0$, with $A_0=\delta^{-\frac{1}{2}}$ is, after scaling with
factor $h^{-2}$, $\delta^{\frac{1}{2}}$-close to the corresponding
ball in the standard capped infinite cylinder with the center near
the tip. Assume first that for each point in $B_{T_0}(x_0,A_0h)$,
the solution is defined on $[T_0,t_1]$. By the gradient estimates
(7.3.4) in the canonical neighborhood assumption and the choice of
$\Delta t$ we have a uniform curvature bound on this set for
$h^{-2}$-scaled metric. Then by the uniqueness theorem (Theorem
1.2.4), if $\delta^{\frac{1}{2}}\rightarrow0$ (i.e.
$A_0=\delta^{-\frac{1}{2}} \rightarrow +\infty$), the solution
with $h^{-2}$-scaled metric will converge to the standard solution
in the $C^\infty_{\rm loc}$ topology. Therefore we can define
$A_1$, depending only on $A_0$ and tending to infinity with $A_0$,
such that the solution in the parabolic region
$P(x_0,T_0,A_1h,t_1-T_0)\triangleq\{(x,t)\ |\ x\in
B_t(x_0,A_1h),t\in[T_0,T_0+(t_1-T_0)]\}$ is, after scaling with
factor $h^{-2}$ and shifting time $T_0$ to zero, $A^{-1}_1$-close
to the corresponding subset in the standard solution. In
particular, the scalar curvature on this subset does not exceed
$2Qh^{-2}$. Now if for each point in $B_{T_0}(x_0,A_1h)$ the
solution is defined on $[T_0,t_2]$, then we can repeat the
procedure, defining $A_2$, such that the solution in the parabolic
region $P(x_0,T_0,A_2h,t_2-T_0)\triangleq\{(x,t)\ |\ x\in
B_t(x_0,A_2h),t\in[T_0,T_0+(t_2-T_0)]\}$ is, after scaling with
factor $h^{-2}$ and shifting time $T_0$ to zero, $A^{-1}_2$-close
to the corresponding subset in the standard solution. Again, the
scalar curvature on this subset still does not exceed $2Qh^{-2}$.
Continuing this way, we eventually define $A_N$. Note that $N$ is
depends only on $\theta$ and $\varepsilon$. Thus there exists a
positive $\bar{\delta}=\bar{\delta}(A,\theta,\varepsilon)$ such
that for $\delta<\bar{\delta}$, we have $A_0>A_1>\cdots>A_N>A$,
and assertion (i) holds when the solution is defined on
$B_{T_0}(x_0,A_{(N-1)}h)\times [T_0,T_1]$.

The above argument shows that either assertion (i) holds, or there
exists some $k$ ($0\leq k\leq N-1$) and a surgery time
$t^+\in(t_k,t_{k+1}]$ such that the solution on
$B_{T_0}(x_0,A_kh)$ is defined on $[T_0,t^+)$, but for some point
of this set it is not defined past $t^+$. Now we consider the
latter case. Clearly the above argument also shows that the
parabolic region $P(x_0,T_0,A_{k+1}h,t^+-T_0)\triangleq\{(x,t)\ |\
x\in B_t(x,A_{k+1}h),t\in[T_0,t^+)\}$ is, after scaling with
factor $h^{-2}$ and shifting time $T_0$ to zero,
$A^{-1}_{k+1}$-close to the corresponding subset in the standard
solution. In particular, as time tends to $t^+$, the ball
$B_{T_0}(x_0,A_{k+1}h)$ keeps on looking like a cap. Since the
scalar curvature on $B_{T_0}(x_0,A_kh) \times [T_0,t_k]$ does not
exceed $2Qh^{-2}$, it follows from the pinching assumption, the
gradient estimates in the canonical neighborhood assumption and
the evolution equation of the metric that the diameter of the set
$B_{T_0}(x_0,A_kh)$ at any time $t \in [T_0,t^+)$ is bounded from
above by $4 \delta^{-\frac{1}{2}}h$. These imply that no point of
the ball $B_{T_0}(x_0,A_kh)$ at any time near $t^+$ can be the
center of a $\delta$-neck for any $0 < \delta <
\bar{\delta}(A,\theta,\varepsilon)$ with
$\bar{\delta}(A,\theta,\varepsilon)
> 0$ small enough, since $4 \delta^{-\frac{1}{2}}h$ is much
smaller than $\delta^{-1}h$. However the solution disappears
somewhere in the set $B_{T_0}(x_0,A_kh)$ at time $t^+$ by a cutoff
surgery and the surgery is always done along the middle two-sphere
of a $\delta$-neck. So the set $B_{T_0}(x_0,A_kh)$  at  time $t^+$
is a part of a capped horn. (Recall that we have declared extinct
every compact component with positive curvature and every compact
component lying in $\Omega \setminus \Omega_{\rho}$). Hence for
each point of $B_{T_0}(x_0,A_kh)$ the solution terminates at
$t^+$. This proves assertion (ii).
\end{pf}

\begin{corollary}[{Perelman \cite{P2}}]
For any $l<\infty$ one can find $A=A(l)<\infty$ and
$\theta=\theta(l)$, $0<\theta<1$, with the following property.
Suppose we are in the situation of the lemma above, with
$\delta<\bar{\delta}(A,\theta,\varepsilon)$. Consider smooth
curves $\gamma$ in the set $B_{T_0}(x_0,Ah)$, parametrized by
$t\in[T_0,T_{\gamma}]$, such that $\gamma(T_0)\in
B_{T_0}(x_0,\frac{Ah}{2})$ and either $T_{\gamma}=T_1<T$, or
$T_{\gamma}<T_1$ and $\gamma(T_{\gamma})\in \partial
B_{T_0}(x_0,Ah)$, where $x_0$ is any fixed point on a gluing cap
at $T_0$ and $T_1=min\{T,T_0+\theta h^2\}$. Then
$$
\int^{T_{\gamma}}_{T_0}(R(\gamma(t),t)+|\dot{\gamma}(t)|^2)dt>l.
$$
\end{corollary}

\begin{pf}
We know from Proposition 7.3.3 that on the standard solution,
\begin{align*}
\int^{\theta}_{0}Rdt
&  \geq{\rm  const.}\,\int^{\theta}_0(1-t)^{-1}dt\\
&  = -{\rm  const.}\,\cdot\log(1-\theta).
\end{align*}
By choosing $\theta=\theta(l)$ sufficiently close to 1 we have the
desired estimate for the standard solution.

Let us consider the first case: $T_{\gamma}=T_1<T$. For
$\theta=\theta(l)$ fixed above, by Lemma 7.3.6, our solution in
the subset $B_{T_0}(x_0,Ah)$ and in the time interval
$[T_0,T_{\gamma}]$ is, after scaling with factor $h^{-2}$ and
shifting time $T_0$ to zero, $A^{-1}$-close to the corresponding
subset in the standard solution for any sufficiently large $A$. So
we have
\begin{align*}
\int^{T_{\gamma}}_{T_0}(R(\gamma(t),t) + |\dot{\gamma}(t)|^2)dt
&  \geq {\rm const.}\,\int^{\theta}_0(1-t)^{-1}dt\\
&  = -{\rm const.}\,\cdot\log(1-\theta).
\end{align*}
Hence we have obtained the desired estimate in the first case.

We now consider the second case: $T_{\gamma}<T_1$ and
$\gamma(T_{\gamma})\in \partial B_{T_0}(x_0,Ah)$. Let
$\theta=\theta(l)$ be chosen above and let $Q=Q(l)$ be the maximum
of the scalar curvature on the standard solution in the time
interval $[0,\theta]$.

On the standard solution, we can choose $A=A(l)$ so large that for
each $t\in[0,\theta]$,
\begin{align*}
{\rm dist}_t(x_0,\partial B_0(x_0,A))
&  \geq {\rm dist}_0(x_0,\partial B_0(x_0,A))-4(Q+1)t\\
&  \geq A-4(Q+1)\theta\\
&  \geq \frac{4}{5}A
\end{align*}
and
$$
{\rm dist}_t\(x_0,\partial
B_0\(x_0,\frac{A}{2}\)\)\leq\frac{A}{2},
$$
where we have used Lemma 3.4.1(ii) in the first inequality. Now
our solution in the subset $B_{T_0}(x_0,Ah)$ and in the time
interval $[T_0,T_{\gamma}]$ is, after scaling with factor $h^{-2}$
and shifting time $T_0$ to zero, $A^{-1}$-close to the
corresponding subset in the standard solution. This implies that
for $A=A(l)$ large enough
$$
\frac{1}{5}Ah\leq\int^{T_{\gamma}}_{T_0}|\dot{\gamma}(t)|dt
\leq\(\int^{T_{\gamma}}_{T_0}|
\dot{\gamma}(t)|^2dt\)^{\frac{1}{2}}\cdot(T_{\gamma}-T_0)^{\frac{1}{2}},
$$
Hence
$$
\int^{T_{\gamma}}_{T_0}(R(\gamma(t),t) +
|\dot{\gamma}(t)|^2)dt\geq \frac{A^2}{25\theta}>l.
$$
This proves the desired estimate.
\end{pf}

\section{Justification of the Canonical Neighborhood Assumptions}

We continue the induction argument for the construction of a
long-time solution to the Ricci flow with surgery. Let us recall
what we have done in the previous section. Let $\varepsilon$ be an
arbitrarily given positive constant satisfying $0<\varepsilon \leq
1/100$. For an arbitrarily given compact orientable normalized
three-manifold, we evolve it by the Ricci flow. We may assume that
the solution goes singular at some time $0<t^+_1<+\infty$ and know
that the solution satisfies the a priori assumptions (with
accuracy $\varepsilon$) on $[0,t^+_1)$ for a nonincreasing
positive function $r=r_1(t)$ (defined on $[0,+\infty)$). Suppose
that we have a solution to the Ricci flow with surgery, defined on
$[0,t^+_k)$ with $0<t^+_1<t^+_2<\cdots<t^+_k<+\infty$, satisfying
the a priori assumptions (with accuracy $\varepsilon$) for some
nonincreasing positive function $r=r_k(t)$ (defined on
$[0,+\infty)$), going singular at time $t^+_k$ and having
$\delta_i$-cutoff surgeries at each time $t^+_i$, $1\leq i\leq
k-1$, where $\delta_i < \bar{\delta}(t^+_i)$ for each $1\leq i\leq
k-1$. Then for any $0<\delta_k < \bar{\delta}(t^+_k)$, we can
perform $\delta_k$-cutoff surgeries at the time $t^+_k$ and extend
the solution to the interval $[0,t^+_{k+1})$ with
$t^+_{k+1}>t^+_k$. Here $\bar{\delta}(t)$ is the positive function
defined in (7.3.16). We have already shown in Lemma 7.3.4 that the
extended solution still satisfies the pinching assumption on
$[0,t^+_{k+1})$.

In view of Theorem 7.1.1, there always is a nonincreasing positive
function $r=r_{k+1}(t)$, defined on $[0,+\infty)$, such that the
canonical neighborhood assumption (with accuracy $\varepsilon$)
holds on the extended time interval $[0,t^+_{k+1})$ with the
positive function $r=r_{k+1}(t)$. Nevertheless, in order to
prevent the surgery times from accumulating, the key is to choose
the nonincreasing positive functions $r=r_i(t), i= 1, 2, \ldots$,
uniformly. That is, to justify the canonical neighborhood
assumption (with accuracy $\varepsilon$) for the indefinitely
extending solution, we need to show that there exists a
nonincreasing positive function $\widetilde{r}(t)$, defined on
$[0,+\infty)$, which is independent of $k$, such that the above
chosen nonincreasing positive functions satisfy
$$r_i(t) \geq \widetilde{r}(t), \ \
\mbox{ on } \ \ [0,+\infty),$$ for all $i = 1, 2, \ldots, k+1$.

By a further restriction on the positive function
$\bar{\delta}(t)$, we can verify this after proving the following
assertion which was stated by Perelman in \cite{P2}.

\begin{proposition}[{Justification of the canonical neighborhood
assumption} \index{justification of the canonical neighborhood
assumption}] Given any small $\varepsilon>0$, there exist
decreasing sequences $0<\widetilde{r}_j<\varepsilon$,
$\kappa_j>0$, and $0<\widetilde{\delta}_j<\varepsilon^2 $,
$j=1,2,\ldots$, with the following property. Define the positive
function $\widetilde{\delta}(t)$ on $[0,+\infty)$ by
$\widetilde{\delta}(t) = \widetilde{\delta}_j$ for $t \in
[(j-1)\varepsilon^2,j\varepsilon^2)$. Suppose there is a
surgically modified solution, defined on $[0,T)$ with $T
<+\infty$, to the Ricci flow which satisfies the following:
\begin{itemize}
\item[(1)] it starts on a compact orientable three-manifold with
normalized initial metric, and \item[(2)] it has only a finite
number of surgeries such that each surgery at a time $t \in (0,T)$
is a $\delta(t)$-cutoff surgery with
$$
0 < \delta(t) \leq \min \{\widetilde{\delta}(t),\bar{\delta}(t)\}.
$$
\end{itemize}
Then on each time interval
$[(j-1)\varepsilon^2,j\varepsilon^2]\bigcap [0,T)$, $ j = 1, 2,
\cdots $, the solution satisfies the $\kappa_j$-noncollapsing
condition on all scales less than $\varepsilon$ and the canonical
neighborhood assumption $($with accuracy $\varepsilon)$ with
$r=\widetilde{r}_j$.
\end{proposition}

Here and in the following, we call a (three-dimensional)
surgically modified solution $g_{ij}(t), 0 \leq t <T$,
\textbf{$\kappa$-noncollapsed} \index{$\kappa$-noncollapsed} at
$(x_0,t_0)$ on the scales less than $\rho$ (for some
$\kappa>0,\rho>0$) if it satisfies the following property:
whenever $r < \rho$ and
$$
|Rm(x,t)| \leq r^{-2}
$$
for all those $(x,t) \in P(x_0,t_0,r,-r^2)=\{ (x',t') \ |\ x'\in
B_{t'}(x_0,r), t' \in [t_0 -r^2,t_0] \}$, for which the solution
is defined, we have
$$
\Vol_{t_0}(B_{t_0}(x_0,r)) \geq \kappa r^3.
$$
Before we give the proof of the proposition, we need to verify a
$\kappa$-non\-collapsing estimate which was given by Perelman in
\cite{P2}.

\begin{lemma} [{Perelman \cite{P2}}]
Given any $0<\varepsilon \leq \bar{\varepsilon}_0 $ $($for some
sufficiently small universal constant $\bar{\varepsilon}_0),$
suppose we have constructed the sequences satisfying the
proposition for $1\leq j\leq m$ $($for some positive integer $m)$.
Then there exists $\kappa>0$, such that for any $r$,
$0<r<\varepsilon$, one can find
$\widetilde{\delta}=\widetilde{\delta}(r,\varepsilon)$,
$0<\widetilde{\delta}<\varepsilon^2$, which may also depend on the
already constructed sequences, with the following property.
Suppose we have a solution with a compact orientable normalized
three-manifold as initial data, to the Ricci flow with finite
number of surgeries on a time interval $[0,\bar{T}]$ with
$m\varepsilon^2\leq \bar{T}<(m+1)\varepsilon^2$, satisfying the
assumptions and the conclusions of Proposition $7.4.1$ on
$[0,m\varepsilon^2)$, and the canonical neighborhood assumption
$($with accuracy $\varepsilon)$ with $r$ on
$[m\varepsilon^2,\bar{T}]$, as well as $0<\delta(t) \leq \min
\{\widetilde{\delta},\bar{\delta}(t)\}$ for any $\delta$-cutoff
surgery with $\delta =\delta(t)$ at a time
$t\in[(m-1)\varepsilon^2,\bar{T}]$. Then the solution is
$\kappa$-noncollapsed on $[0,\bar{T}]$ for all scales less than
$\varepsilon$.
\end{lemma}

\begin{pf}
Consider a parabolic neighborhood
$$
P(x_0,t_0,r_0,-r^2_0)\triangleq\{(x,t)\ |\ x\in
B_t(x_0,r_0),t\in[t_0-r^2_0,t_0]\}
$$
with $m\varepsilon^2\leq t_0\leq \bar{T}$ and $0<r_0<\varepsilon$,
where the solution satisfies $|Rm|\leq r^{-2}_0$, whenever it is
defined. We will use an argument analogous to the proof of Theorem
3.3.2 (no local collapsing theorem I) to prove \be
\Vol_{t_0}(B_{t_0}(x_0,r_0)) \geq \kappa
r^3_0. 
\ee

Let $\eta$ be the universal positive constant in the definition of
the canonical neighborhood assumption. Without loss of generality,
we always assume $\eta \geq 10$. Firstly, we want to show that one
may assume $r_0 \geq \frac{1}{2\eta}r$.

Obviously, the curvature satisfies the estimate
$$
|Rm(x,t)|\leq 20r^{-2}_0,
$$
for those $(x,t) \in
P(x_0,t_0,\frac{1}{2\eta}r_0,-\frac{1}{8\eta}r^2_0) = \{(x,t)\ |\
x\in
B_t(x_0,\frac{1}{2\eta}r_0),t\in[t_0-\frac{1}{8\eta}r^2_0,t_0]\}$,
for which the solution is defined. When $r_0< \frac{1}{2\eta}r$,
we can enlarge $r_0$ to some $r'_0\in[r_0,r]$ so that
$$
|Rm|\leq 20r'^{-2}_0
$$
on $P(x_0,t_0,\frac{1}{2\eta}r'_0,-\frac{1}{8\eta}r'^2_0)$
(whenever it is defined), and either the equality holds somewhere
in $P(x_0,t_0,\frac
{1}{2\eta}r'_0,-(\frac{1}{8\eta}{r'}^2_0+\epsilon'))$ for
arbitrarily small $\epsilon'>0$, or $r'_0=r$.

In the case that the equality holds somewhere, it follows from the
pinching assumption that we have
$$
R > 10r'^{-2}_0
$$
somewhere in $P(x_0,t_0,\frac{1}{2\eta}r'_0,-(\frac{1}{8\eta}{r'}
^2_0+\epsilon'))$ for arbitrarily small $\epsilon'>0$. Here,
without loss of generality, we have assumed $r$ is suitably small.
Then by the gradient estimates in the definition of the canonical
neighborhood assumption, we know $$R(x_0,t_0) > r'^{-2}_0 \geq
r^{-2}.$$ Hence the desired noncollapsing estimate (7.4.1) in this
case follows directly from the canonical neighborhood assumption.
(Recall that we have excluded every compact component which has
positive sectional curvature in the surgery procedure and then we
have excluded them from the list of canonical neighborhoods. Here
we also used the standard volume comparison when the canonical
neighborhood is an $\varepsilon$-cap.)

While in the case that $r'_0 = r$, we have the curvature bound
$$
|Rm(x,t)|\leq \(\frac{1}{2\eta}r\)^{-2},
$$
for those $(x,t) \in
P(x_0,t_0,\frac{1}{2\eta}r,-(\frac{1}{2\eta}r)^{2}) = \{(x,t)\ |\
x\in
B_t(x_0,\frac{1}{2\eta}r),t\in[t_0-(\frac{1}{2\eta}r)^2,t_0]\}$,
for which the solution is defined. It follows from the standard
volume comparison that we only need to verify the noncollapsing
estimate (7.4.1) for $r_0 = \frac{1}{2\eta}r$. Thus we have
reduced the proof to the case $r_0\geq \frac{1}{2\eta}r$.

Recall from Theorem 3.3.2 that if a solution is smooth everywhere,
we can get a lower bound for the volume of the ball
$B_{t_0}(x_0,r_0)$ as follows: define $\tau(t)=t_0-t$ and consider
Perelman's reduced volume function and the Li-Yau-Perelman
distance associated to the point $x_0$; take a point $\bar{x}$ at
the time $ t=\varepsilon^2$ so that the Li-Yau-Perelman distance
$l$ attains its minimum
$l_{\min}(\tau)=l(\bar{x},\tau)\leq\frac{3}{2}$ for $\tau
=t_0-\varepsilon^2$; use it to obtain an upper bound for the
Li-Yau-Perelman distance from $x_0$ to each point of
$B_0(\bar{x},1)$, thus getting a lower bound for Perelman's
reduced volume at $\tau=t_0$; apply the monotonicity of Perelman's
reduced volume to deduce a lower bound for Perelman's reduced
volume at $\tau$ near $0$, and then get the desired estimate for
the volume of the ball $B_{t_0}(x_0,r_0)$. Now since our solution
has undergone surgeries, we need to localize this argument to the
region which is unaffected by surgery.

We call a space-time curve in the solution track
\textbf{admissible}\index{admissible curve} if it stays in the
space-time region unaffected by surgery, and we call a space-time
curve in the solution track a \textbf{barely admissible
curve}\index{barely admissible curve} if it is on the boundary of
the set of admissible curves.

 First of all, we want to estimate the $\mathcal{L}$-length of a
barely admissible curve.

\medskip
{\bf Claim.} \ For any $L<\infty$ one can find
$\bar{\delta}=\bar{\delta}(L,r,\widetilde{r}_m,\varepsilon)>0$
with the following property. Suppose that we have a curve
$\gamma$, parametrized by $t\in[T_0,t_0]$, $(m-1)\varepsilon^2\leq
T_0<t_0$, such that $\gamma(t_0)=x_0$, $T_0$ is a surgery time,
and $\gamma(T_0)$ lies in the gluing cap. Suppose also each
$\delta$-cutoff surgery at a time in
$[(m-1)\varepsilon^2,\bar{T}]$ has $\delta \leq \bar{\delta}$.
Then we have an estimate \be
\int^{t_0}_{T_0}\sqrt{t_0-t}(R_+(\gamma(t),t)
+|\dot{\gamma}(t)|^2)dt\geq L  
\ee where $R_+=\max\{R,0\}$.

\medskip
Since $r_0\geq \frac{1}{2\eta}r$ and $|Rm|\leq r^{-2}_0$ on
$P(x_0,t_0,r_0,-r^2_0)$ (whenever it is defined), we can require
$\bar{\delta}>0$, depending on $r$ and $\widetilde{r}_m$, to be so
small that $\gamma(T_0)$ does not lie in the region
$P(x_0,t_0,r_0,-r^2_0)$. Let $\Delta t$ be the maximal number such
that $\gamma|_{[t_0-\Delta t,t_0]}\subset P(x_0,t_0,r_0,-\Delta
t)$ (i.e., $t_0-\Delta t$ is the first time when $\gamma$ escapes
the parabolic region $P(x_0,t_0,r_0,-r^2_0)).$ Obviously we only
need to consider the case:
$$
\int^{t_0}_{t_0 -\Delta t}\sqrt{t_0 -
t}(R_+(\gamma(t),t)+|\dot{\gamma}(t)|^2)dt<L.
$$

We observe that $\Delta t$ can be bounded from below in terms of
$L$ and $r_0$. Indeed, if $\Delta t \geq r_0^2$, there is nothing
to prove. Thus we may assume $\Delta t<r_0^2$. By the curvature
bound $|Rm|\leq r^{-2}_0$ on $P(x_0,t_0,r_0,-r^2_0)$ and the Ricci
flow equation we see
$$
\int^{t_0}_{t_0-\Delta t}|\dot{\gamma}(t)|dt\geq cr_0
$$
for some universal positive constant $c$. On the other hand, by
the Cauchy-Schwarz inequality, we have
\begin{align*}
\int^{t_0}_{t_0-\Delta t}|\dot{\gamma}(t)|dt & \leq
\(\int^{t_0}_{t_0-\Delta t}\sqrt{t_0-t}(R_++|
\dot{\gamma}|^2)dt\)^{\frac{1}{2}} \cdot\(\int^{t_0}_{t_0-\Delta
t}\frac{1}{\sqrt{t_0-t}}dt\)^{\frac{1}{2}}\\
&  \leq (2L)^{\frac{1}{2}}(\Delta t)^{\frac{1}{4}},
\end{align*}
which implies \be
(\Delta t)^{\frac{1}{2}}\geq\frac{c^2r^2_0}{2L}. 
\ee Thus
\begin{align*}
\int^{t_0}_{T_0}\sqrt{t_0-t}(R_++|\dot{\gamma}|^2)dt& \geq
\int^{t_0-\Delta
t}_{T_0}\sqrt{t_0-t}(R_++|\dot{\gamma}|^2)dt\\
&  \geq (\Delta t)^{\frac{1}{2}}\int^{t_0-\Delta
t}_{T_0}(R_++|\dot{\gamma}|^2)dt\\
&  \geq \(\min \left\{\frac{c^2r^2_0}{2L},r_0 \right\}\)
\int^{t_0-\Delta t}_{T_0}(R_++|\dot{\gamma}|^2)dt,
\end{align*}
while by Corollary 7.3.7, we can find
$\bar{\delta}=\bar{\delta}(L,r,\widetilde{r}_m,\varepsilon)>0$ so
small that
$$
\int^{t_0-\Delta t}_{T_0}(R_++|\dot{\gamma}|^2)dt \geq
L\(\min\left\{\frac{c^2r^2_0}{2L},r_0 \right\}\)^{-1}.
$$
Then we have proved the desired assertion (7.4.2).

Recall that for a curve $\gamma$, parametrized by
$\tau=t_0-t\in[0,\bar{\tau}],$ with $\gamma(0)=x_0$ and
$\bar{\tau}\leq t_0-(m-1)\varepsilon^2$, we have
$L(\gamma)=\int^{\bar{\tau}}_0\sqrt{\tau}(R+|\dot{\gamma}|^2)d\tau$.
We can also define $L_+(\gamma)$ by replacing $R$ with $R_+$ in
the previous formula. Recall that $R\geq-1$ at the initial time
$t=0$ for the normalized initial manifold. Recall that the
surgeries occur at the parts where the scalar curvatures are very
large. Thus we can apply the maximum principle to conclude that
the solution with surgery still satisfies $R\geq -1$  everywhere
in space-time. This implies \be L_+(\gamma)\leq
L(\gamma)+(2\varepsilon^2)^{\frac{3}{2}}.
\ee By applying the assertion (7.4.2), we now choose
$\tilde{\delta}>0$ (depending on $r$, $\varepsilon$ and
$\widetilde{r}_m$) such that as each $\delta$-cutoff surgery at
the time interval $t\in[(m-1)\varepsilon^2,T]$ has $\delta \leq
\tilde{\delta}$, every barely admissible curve $\gamma$ from
$(x_0,t_0)$ to a point $(x,t)$ (with
$t\in[(m-1)\varepsilon^2,t_0)$) has
$$
L_+(\gamma)\geq22\sqrt{2}.
$$
Thus if the Li-Yau-Perelman distance from $(x_0,t_0)$ to a point
$(x,t)$ (with $t\in[(m-1)\varepsilon^2,t_0)$) is achieved by a
space-time curve which is not admissible, then its Li-Yau-Perelman
distance has \be
l\geq\frac{L_+-(2\varepsilon^2)^\frac{3}{2}}{2\sqrt{2}
\varepsilon}>10\varepsilon^{-1}. 
\ee We also observe that the absolute value of $l(x_0,\tau)$ is
very small as $\tau$ closes to zere. Thus the maximum principle
argument in Corollary 3.2.6 still works for our solutions with
surgery because barely admissible curves do not attain the
minimum. So we conclude that
$$
l_{\min}(\bar{\tau})=\min\{l(x,\bar{\tau})\ |\ x \mbox{ lies in
the solution manifold at }t_0-\bar{\tau}\} \leq\frac{3}{2}
$$
for $\bar{\tau}\in(0,t_0-(m-1)\varepsilon^2]$. In particular,
there exists a minimizing curve $\gamma$ of
$l_{\min}(t_0-(m-1)\varepsilon^2)$, defined on
$\tau\in[0,t_0-(m-1)\varepsilon^2]$ with $\gamma(0)=x_0$, such
that
\begin{align}
L_+(\gamma)
& \leq \frac{3}{2}\cdot2\sqrt{2}\varepsilon+2\sqrt{2}\varepsilon^3\\
&  \leq 5\varepsilon, \nn
\end{align}  
since $0<\varepsilon \leq \bar{\varepsilon}_0$ with
$\bar{\varepsilon}_0$ sufficiently small (to be further
determined). Consequently, there exists a point
$(\bar{x},\bar{t})$ on the minimizing curve $\gamma$ with
$\bar{t}\in[(m-1)\varepsilon^2+\frac{1}{4}\varepsilon^2,
(m-1)\varepsilon^2+\frac{3}{4}\varepsilon^2]$ (i.e.,
$\tau\in[t_0-(m-1)\varepsilon^2-\frac{3}{4}\varepsilon^2,
t_0-(m-1)\varepsilon^2-\frac{1}{4}\varepsilon^2])$ such that \be
R(\bar{x},\bar{t})\leq25\widetilde{r}^{-2}_m. 
\ee Otherwise, we have
\begin{align*}
L_+(\gamma)& \geq
\int^{t_0-(m-1)\varepsilon^2-\frac{1}{4}\varepsilon^2}_{t_0-(m-1)
\varepsilon^2-\frac{3}{4}\varepsilon^2}
\sqrt{\tau}R(\gamma(\tau),t_0-\tau)d\tau\\
&  > 25\widetilde{r}^{-2}_m\sqrt{\frac{1}{4}\varepsilon^2}
\(\frac{1}{2}\varepsilon^2\)\\
&  > 5\varepsilon,
\end{align*}
since $0<\widetilde{r}_m<\varepsilon$. This contradicts (7.4.6).

Next we want to get a lower bound for Perelman's reduced volume of
a ball around $\bar{x}$ of radius about $\widetilde{r}_m$ at some
time slightly before $\bar{t}$.

Denote by $\theta_1 = \frac{1}{16}\eta^{-1}$ and $\theta_2 =
\frac{1}{64}\eta^{-1}$, where $\eta$ is the universal positive
constant in the gradient estimates (7.3.4). Since the solution
satisfies the canonical neighborhood assumption on the time
interval $[(m-1)\varepsilon^2,m\varepsilon^2)$, it follows from
the gradient estimates (7.3.4) that \be
R(x,t)\leq 400\widetilde{r}^{-2}_m  
\ee for those $(x,t)\in P(\bar{x},\bar{t},\theta_1\widetilde{r}_m,
-\theta_2\widetilde{r}^2_m)\triangleq\{(x',t')\ |\ x'\in
B_{t'}(\bar{x},\theta_1\widetilde{r}_m),t'\in[\bar{t}
-\theta_2\widetilde{r}^2_m,\bar{t}]\}$, for which the solution is
defined. And since the scalar curvature at the points where the
$\delta$-cutoff surgeries occur in the time interval
$[(m-1)\varepsilon^2,m\varepsilon^2)$ is at least
$(\widetilde{\delta})^{-2}\widetilde{r}^{-2}_m$, the solution is
well-defined on the whole parabolic region
$P(\bar{x},\bar{t},\theta_1\widetilde{r}_m,-\theta_2\widetilde{r}^2_m)$
(i.e., this parabolic region is unaffected by surgery). Thus by
combining (7.4.6) and (7.4.8), we know that the Li-Yau-Perelman
distance from $(x_0,t_0)$ to each point of the ball
$B_{\bar{t}-\theta_2\widetilde{r}^2_m}(\bar{x},
\theta_1\widetilde{r}_m)$ is uniformly bounded by some universal
constant. Let us define Perelman's reduced volume of the ball
$B_{\bar{t}-\theta_2\widetilde{r}^2_m}(\bar{x},
\theta_1\widetilde{r}_m)$, by
\begin{align*}
&\widetilde{V}_{t_0-\bar{t}+\theta_2\widetilde{r}^2_m}
(B_{\bar{t}-\theta_2\widetilde{r}^2_m}(\bar{x},\theta_1\widetilde{r}_m))\\
&=\int_{B_{\bar{t}-\theta_2\widetilde{r}^2_m}(\bar{x},
\theta_1\widetilde{r}_m)}(4\pi(t_0-\bar{t}+\theta_2\widetilde{r}^2_m))^
{-\frac{3}{2}} \\
&\qquad\cdot\exp(-l(q,t_0-\bar{t}+\theta_2\widetilde{r}^2_m))
dV_{\bar{t}-\theta_2\widetilde{r}^2_m}(q),
\end{align*}
where $l(q,\tau)$ is the Li-Yau-Perelman distance from
$(x_0,t_0)$. Hence by the $\kappa_m$-noncollapsing assumption on
the time interval $[(m-1)\varepsilon^2,m\varepsilon^2)$, we
conclude that Perelman's reduced volume of the ball
$B_{\bar{t}-\theta_2\widetilde{r}^2_m}(\bar{x},\theta_1\widetilde{r}_m)$
is bounded from below by a positive constant depending only on
$\kappa_m$ and $\widetilde{r}_m$.

Finally we want to apply a local version of the monotonicity of
Perelman's reduced volume to get a lower bound estimate for the
volume of the ball $B_{t_0}(x_0,r_0)$.

We have seen that the Li-Yau-Perelman distance from $(x_0,t_0)$ to
each point of the ball
$B_{\bar{t}-\theta_2\widetilde{r}^2_m}(\bar{x},\theta_1\widetilde{r}_m)$
is uniformly bounded by some universal constant. Now we can choose
a sufficiently small (universal) positive constant
$\bar{\varepsilon}_0$ such that when $0<\varepsilon \leq
\bar{\varepsilon}_0$, by (7.4.5), all the points in the ball
$B_{\bar{t}-\theta_2\widetilde{r}^2_m}(\bar{x},\theta_1\widetilde{r}_m)$
can be connected to $(x_0,t_0)$ by shortest
$\mathcal{L}$-geodesics, and all of these $\mathcal{L}$-geodesics
are admissible (i.e., they stay in the region unaffected by
surgery). The union of all shortest $\mathcal{L}$-geodesics from
$(x_0,t_0)$ to the ball
$B_{\bar{t}-\theta_2\widetilde{r}^2_m}(\bar{x},
\theta_1\widetilde{r}_m)$ defined by
$CB_{\bar{t}-\theta_2\widetilde{r}^2_m}(\bar{x},
\theta_1\widetilde{r}_m)=\{(x,t)\ |\ (x,t)$ lies in a shortest
$\mathcal{L}$-geodesic from $(x_0,t_0)$ to a point in
$B_{\bar{t}-\theta_2\widetilde{r}^2_m}(\bar{x},
\theta_1\widetilde{r}_m)\}$, forms a cone-like subset in
space-time with the vertex $(x_0,t_0)$. Denote $B(t)$ by the
intersection of the cone-like subset
$CB_{\bar{t}-\theta_2\widetilde{r}^2_m}(\bar{x},\theta_1\widetilde{r}_m)$
with the time-slice at $t$. Perelman's reduced volume of the
subset $B(t)$ is given by
$$
\widetilde{V}_{t_0-t}(B(t))
=\int_{B(t)}(4\pi(t_0-t))^{-\frac{3}{2}}exp(-l(q,t_0-t))dV_t(q).
$$
Since the cone-like subset
$CB_{\bar{t}-\theta_2\widetilde{r}^2_m}(\bar{x},\theta_1\widetilde{r}_m)$
lies entirely in the region unaffected by surgery, we can apply
Perelman's Jacobian comparison theorem (Theorem 3.2.7) to conclude
that
\begin{align}
\widetilde{V}_{t_0-t}(B(t))& \geq
\widetilde{V}_{t_0-\bar{t}+\theta_2\widetilde{r}^2_m}
(B_{\bar{t}-\theta_2\widetilde{r}^2_m}(\bar{x},
\theta_1\widetilde{r}_m))\\
&  \geq c(\kappa_m,\widetilde{r}_m), \nn
\end{align}  
for all $t\in[\bar{t}-\theta_2\widetilde{r}^2_m,t_0]$, where
$c(\kappa_m,\widetilde{r}_m)$ is some positive constant depending
only on $\kappa_m$ and $\widetilde{r}_m$.

Set $\xi=r^{-1}_0\Vol_{t_0}(B_{t_0}(x_0,r_0))^{\frac{1}{3}}$. Our
purpose is to give a positive lower bound for $\xi$. Without loss
of generality, we may assume $\xi<\frac{1}{4}$, thus $0<\xi
r^2_0<t_0-\bar{t}+\theta_2\widetilde{r}^2_m$. Denote by
$\widetilde{B}(t_0-\xi r^2_0)$ the subset of the time-slice
$\{t=t_0-\xi r^2_0\}$ of which every point can be connected to
$(x_0,t_0)$ by an admissible shortest $\mathcal{L}$-geodesic.
Clearly, $B(t_0-\xi r^2_0)\subset\widetilde{B}(t_0-\xi r^2_0)$. We
now argue as in the proof of Theorem 3.3.2 to bound Perelman's
reduced volume of $\widetilde{B}(t_0-\xi r^2_0)$ from above.

Since $r_0\geq \frac{1}{2\eta}r$ and
$\tilde{\delta}=\tilde{\delta}(r,\varepsilon,\widetilde{r}_m)$
sufficiently small, the whole region $P(x_0,t_0,r_0,$ $-r^2_0)$ is
unaffected by surgery. Then by exactly the same argument as in
deriving (3.3.5), we see that there exists a universal positive
constant $\xi_0$ such that when $0<\xi\leq\xi_0$, there holds \be
\mathcal{L} \exp_{\{|\upsilon|\leq
\frac{1}{4}\xi^{-\frac{1}{2}}\}}(\xi r^2_0)\subset
B_{t_0}(x_0,r_0). 
\ee Perelman's reduced volume of $\widetilde{B}(t_0-\xi r^2_0)$ is
given by
\begin{align}
& \widetilde{V}_{\xi r^2_0}(\widetilde{B}(t_0-\xi r^2_0))\\
&  = \int_{\widetilde{B}(t_0-\xi r^2_0)}(4\pi\xi
r^2_0)^{-\frac{3}{2}}\exp(-l(q,\xi r^2_0))dV_{t_0-\xi
r^2_0}(q)\nn\\
&  = \int_{\widetilde{B}(t_0-\xi r^2_0)\cap \mathcal{L}
\exp_{\{|\upsilon|\leq \frac{1}{4}\xi^{-\frac{1}{2}}\}}(\xi
r^2_0)}(4\pi\xi r^2_0)^{-\frac{3}{2}}\exp(-l(q,\xi
r^2_0))dV_{t_0-\xi r^2_0}(q)\nn\\
& \quad +\int_{\widetilde{B}(t_0-\xi r^2_0)\setminus
\mathcal{L}\exp_{\{|\upsilon|\leq
\frac{1}{4}\xi^{-\frac{1}{2}}\}}(\xi r^2_0)}(4\pi\xi
r^2_0)^{-\frac{3}{2}}\exp(-l(q,\xi r^2_0))dV_{t_0-\xi
r^2_0}(q).\nn
\end{align}  
The first term on the RHS of (7.4.11) can be estimated by
\begin{align}
&\int_{\widetilde{B}(t_0-\xi r^2_0)\cap \mathcal{L}
\exp_{\{|\upsilon| \leq\frac{1}{4}\xi^{-\frac{1}{2}}\}}(\xi
r^2_0)}(4\pi\xi
r^2_0)^{-\frac{3}{2}}\exp(-l(q,\xi r^2_0))dV_{t_0-\xi r^2_0}(q)\\
&\leq e^{C\xi}(4\pi)^{-\frac{3}{2}}\cdot\xi^{\frac{3}{2}} \nn
\end{align}  
for some universal constant $C$, as in deriving (3.3.7). While as
in deriving (3.3.8), the second term on the RHS of (7.4.11) can be
estimated by
\begin{align}
&  \int_{\widetilde{B}(t_0-\xi r^2_0)\setminus\mathcal{L}
\exp_{\{|\upsilon|\leq\frac{1}{4}\xi^{-\frac{1}{2}}\}}{(\xi
r^2_0)}}(4\pi \xi r^2_0)^{-\frac{3}{2}}\exp(-l(q,\xi
r^2_0))dV_{t_0-\xi r^2_0}(q)\\
&  \leq \int_{\{|\upsilon|>\frac{1}{4}\xi^{-\frac{1}{2}}\}}
(4\pi\tau)^{-\frac{3}{2}}\exp(-l(\tau))\mathcal{J}
(\tau)|_{\tau=0}d\upsilon\nn\\
& =(4\pi)^{-\frac{3}{2}}\int_{\{|\upsilon|>\frac{1}{4}
\xi^{-\frac{1}{2}}\}}\exp(-|\upsilon|^2)d\upsilon,\nn
\end{align}  
where we have used Perelman's Jacobian comparison theorem (Theorem
3.2.7) in the first inequality. Hence the combination of (7.4.9),
(7.4.11), (7.4.12) and (7.4.13) bounds $\xi$ from below by a
positive constant depending only on $\kappa_m$ and
$\widetilde{r}_m$. Therefore we have completed the proof of the
lemma.
\mbox{ \ \ }\end{pf}

We are ready to prove the proposition.

\medskip
{\bf\em Proof of Proposition} {\bf 7.4.1.} \ We now follow
Perelman \cite{P2} to prove the proposition by induction: having
constructed our sequences for $1\leq j\leq m$, we make one more
step, defining $\widetilde{r}_{m+1}$, $\kappa_{m+1}$,
$\widetilde{\delta}_{m+1}$, and redefining
$\widetilde{\delta}_m=\widetilde{\delta}_{m+1}$. In view of the
previous lemma, we only need to define $\widetilde{r}_{m+1}$ and
$\widetilde{\delta}_{m+1}$.

In Theorem 7.1.1 we have obtained the canonical neighborhood
structure for smooth solutions. When adapting the arguments in the
proof of Theorem 7.1.1 to the present surgical solutions, we will
encounter the new difficulty of how to take a limit for the
surgically modified solutions. The idea to overcome the difficulty
consists of two parts. The first part, due to Perelman \cite{P2},
is to choose $\widetilde{\delta}_{m}$ and
$\widetilde{\delta}_{m+1}$ small enough to push the surgical
regions to infinity in space. (This is the reason why we need to
redefine $\widetilde{\delta}_{m} = \widetilde{\delta}_{m+1}$.) The
second part (see Assertions 1-3 proved in step 2 below), due to
the authors and Bing-Long Chen, is to show that solutions are
smooth on some small, but uniform, time intervals (on compact
subsets) so that we can apply Hamilton's compactness theorem,
since we only have curvature bounds; otherwise Shi's interior
derivative estimate may not be applicable.

We now argue by contradiction. Suppose for sequence of positive
numbers $r^{\alpha}$ and $\widetilde{\delta}^{\alpha\beta}$,
satisfying $r^{\alpha}\rightarrow0$ as $\alpha\rightarrow\infty$ and
$\widetilde{\delta}^{\alpha\beta}\leq\frac{1}{\alpha \cdot
\beta}(\rightarrow0)$, there exist sequences of solutions
$g^{\alpha\beta}_{ij}$ to the Ricci flow with surgery, where each of
them has only a finite number of cutoff surgeries and has a compact
orientable normalized three-manifold as initial data, so that the
following two assertions hold:
\begin{itemize}
\item[(i)] each $\delta$-cutoff at a time
$t\in[(m-1)\varepsilon^2,(m+1)\varepsilon^2]$ satisfies $\delta
\leq \widetilde{\delta}^{\alpha\beta}$; and \item[(ii)] the
solutions satisfy the statement of the proposition on
$[0,m\varepsilon^2]$, but violate the canonical neighborhood
assumption (with accuracy $\varepsilon$) with $r=r^{\alpha}$ on
$[m\varepsilon^2,(m+1)\varepsilon^2]$.
\end{itemize}

For each solution $g^{\alpha\beta}_{ij}$, we choose $\bar{t}$
(depending on $\alpha$ and $\beta$) to be the nearly first time
for which the canonical neighborhood assumption (with accuracy
$\varepsilon$) is violated. More precisely, we choose
$\bar{t}\in[m\varepsilon^2,(m+1)\varepsilon^2]$ so that the
canonical neighborhood assumption with $r=r^{\alpha}$ and with
accuracy parameter $\varepsilon$ is violated at some
$(\bar{x},\bar{t})$, however the canonical neighborhood assumption
with accuracy parameter $2\varepsilon$ holds on
$t\in[m\varepsilon^2,\bar{t}]$. After passing to subsequences, we
may assume each $\widetilde{\delta}^{\alpha\beta}$ is less than
the $\widetilde{\delta}$ in Lemma 7.4.2 with $r=r^{\alpha}$ when
$\alpha$ is fixed. Then by Lemma 7.4.2 we have uniform
$\kappa$-noncollapsing on all scales less than $\varepsilon$ on
$[0,\bar{t}]$ with some $\kappa>0$ independent of $\alpha,\beta$.

Slightly abusing notation, we will often drop the indices $\alpha$
and $\beta$.

Let $\widetilde{g}^{\alpha\beta}_{ij}$ be the rescaled solutions
around $(\bar{x},\bar{t})$ with factors $R(\bar{x},\bar{t})(\geq
r^{-2}\rightarrow+\infty)$ and shift the times $\bar{t}$ to zero.
We hope to take a limit of the rescaled solutions for subsequences
of $\alpha,\beta\rightarrow\infty$ and show the limit is an
orientable ancient $\kappa$-solution, which will give the desired
contradiction. We divide our arguments into the following six
steps.

\medskip
{\it Step} 1. \  Let $(y,\hat{t})$ be a point on the rescaled
solution $\widetilde{g}^{\alpha\beta}_{ij}$ with
$\widetilde{R}(y,\hat{t})\leq A$ (for some $A\geq 1$) and
$\hat{t}\in [-(\bar{t} -
(m-1)\varepsilon^{2})R(\bar{x},\bar{t}),0]$. Then we have estimate
\be
\widetilde{R}(x,t)\leq10A 
\ee for those $(x,t)$ in the parabolic neighborhood
$P(y,\hat{t},\frac{1}{2}\eta^{-1}A^{-\frac{1}{2}},
-\frac{1}{8}\eta^{-1}A^{-1})$ $\triangleq\{(x',t')\ |\ x'\in
\widetilde{B}_{t'}(y,\frac{1}{2}\eta^{-1}A^{-\frac{1}{2}}),
t'\in[\hat{t}-\frac{1}{8}\eta^{-1}A^{-1},\hat{t}]\}$, for which
the rescaled solution is defined.

Indeed, as in the first step of the proof of Theorem 7.1.1, this
follows directly from the gradient estimates (7.3.4) in the
canonical neighborhood assumption with parameter $2\varepsilon$.

\medskip
{\it Step} 2. \ In this step, we will prove three time extension
results.

\medskip
{\bf Assertion 1.} \ For arbitrarily fixed $\alpha$,
$0<A<+\infty$,  $1 \leq C<+\infty$ and $0 \leq B <
\frac{1}{2}\varepsilon^2(r^{\alpha})^{-2}-
\frac{1}{8}\eta^{-1}C^{-1}$, there is a
$\beta_0=\beta_0(\varepsilon, A, B, C)$ (independent of $\alpha$)
such that if $\beta\geq\beta_0$ and the rescaled solution
$\widetilde{g}^{\alpha\beta}_{ij}$ on the ball
$\widetilde{B}_{0}(\bar{x},A)$ is defined on a time interval
$[-b,0]$ with $0 \leq b \leq B$ and the scalar curvature satisfies
$$\widetilde{R}(x,t)\leq C , \ \ \mbox{ on }
\widetilde{B}_{0}(\bar{x},A) \times [-b,0],$$ then the rescaled
solution $\widetilde{g}^{\alpha\beta}_{ij}$ on the ball
$\widetilde{B}_{0}(\bar{x},A)$ is also defined on the extended
time interval $ [-b-\frac{1}{8}\eta^{-1}C^{-1},0]$.

\medskip
Before giving the proof, we make a simple \textbf{observation}:
once a space point in the Ricci flow with surgery is removed by
surgery at some time, then it never appears for later time; if a
space point at some time $t$ cannot be defined before the time $t$
, then either the point lies in a gluing cap of the surgery at
time $t$ or the time $t$ is the initial time of the Ricci flow.

\medskip
{\bf \em Proof of Assertion} {\bf 1.} \ Firstly we claim that
there exists $\beta_0=\beta_0( \varepsilon, A, B, C)$ such that
when $\beta\geq\beta_0$, the rescaled solution
$\widetilde{g}^{\alpha\beta}_{ij}$ on the ball
$\widetilde{B}_{0}(\bar{x},A)$ can be defined before the time $-b$
(i.e., there are no surgeries interfering in
$\widetilde{B}_{0}(\bar{x},A)\times [-b-\epsilon',-b]$ for some
$\epsilon'>0$).

We argue by contradiction. Suppose not, then there is some point
$\tilde{x}\in \widetilde{B}_{0}(\bar{x},A)$ such that the rescaled
solution $\widetilde{g}^{\alpha\beta}_{ij}$ at $\tilde{x}$ cannot
be defined before the time $-b$. By the above observation, there
is a surgery at the time $-b$ such that the point $\tilde{x}$ lies
in the instant gluing cap.

Let\; $\tilde{h}$\; $(=R(\bar{x},\,\bar{t})^{\frac{1}{2}}h$)\;
be\; the\; cut-off\; radius at the\; time\; $-b$\; for the
rescaled solution.  Clearly, there is a universal constant $D$
such that
${D}^{-1}\tilde{h}\leq\widetilde{R}(\tilde{x},-b)^{-\frac{1}{2}}\leq
D\tilde{h}$.

By Lemma 7.3.4 and looking at the rescaled solution at the time
$-b$, the gluing cap and the adjacent $\delta$-neck, of radius
$\tilde{h}$, constitute a
${(\widetilde{\delta}^{\alpha\beta})}^{\frac{1}{2}}$-cap
$\mathcal{K}$. For any fixed small positive constant
$\delta^{\prime}$ (much smaller than $\varepsilon$), we see that
$$
\widetilde{B}_{(-b)}(\tilde{x},{(\delta^{\prime})}^{-1}
\widetilde{R}(\tilde{x}, -b)^{-\frac{1}{2}})\subset\mathcal{K}
$$
when $\beta$ large enough.  We first verify the following

\medskip
{\bf Claim 1.} \ For any small constants $0<\tilde{\theta}<1$,
$\delta'>0$, there exists a $\beta(\delta',\varepsilon,
\tilde{\theta})> 0$ such that when $\beta \geq
\beta(\delta',\varepsilon, \tilde{\theta})$, we have
\begin{itemize}
\item[(i)] the rescaled solution
$\widetilde{g}^{\alpha\beta}_{ij}$ over
$\widetilde{B}_{(-b)}(\tilde{x},{(\delta^{\prime})}^{-1}
\tilde{h})$  is defined on the time interval
$[-b,0]\cap[-b,-b+(1-\tilde{\theta})\tilde{h}^{2}]$; \item[(ii)]
the ball
$\widetilde{B}_{(-b)}(\tilde{x},{(\delta^{\prime})}^{-1}\tilde{h})$
in the ${(\widetilde{\delta}^{\alpha\beta})}^{\frac{1}{2}}$-cap
$\mathcal{K}$ evolved by the Ricci flow  on the time interval
$[-b,0]\cap[-b,-b+(1-\tilde{\theta})\tilde{h}^{2}]$ is, after
scaling with factor $\tilde{h}^{-2}$, ${\delta}^{\prime}$-close
(in the $C^{[\delta'^{-1}]}$ topology) to the corresponding subset
of the standard solution.
\end{itemize}

\medskip
This claim essentially follows from Lemma 7.3.6. Indeed, suppose
there is a surgery at some time $\tilde{\tilde{t}}\in
[-b,0]\cap(-b,-b+(1-\tilde{\theta})\tilde{h}^{2}]$ which removes
some point  $\tilde{\tilde{x}}\in \widetilde
B_{(-b)}(\tilde{x},{(\delta^{\prime})}^{-1}\tilde{h})$. We assume
$\tilde{\tilde{t}}\in (-b,0]$ is the first time with that
property.

Then by Lemma 7.3.6, there is a
$\bar{\delta}=\bar{\delta}(\delta',\varepsilon,\tilde{\theta})$
such that if $\widetilde{\delta}^{\alpha\beta}<\bar{\delta}$, then
the ball
$\widetilde{B}_{(-b)}(\tilde{x},{(\delta^{\prime})}^{-1}\tilde{h})$
in the ${(\widetilde{\delta}^{\alpha\beta})}^{\frac{1}{2}}$-cap
$\mathcal{K}$ evolved by the Ricci flow  on the time interval
$[-b,\tilde{\tilde{t}})$ is, after scaling with factor
$\tilde{h}^{-2}$, ${\delta}^{\prime}$-close to the corresponding
subset of the standard solution. Note that the metrics for times
in $[-b,\tilde{\tilde{t}})$ on $\widetilde
B_{(-b)}(\tilde{x},{(\delta^{\prime})}^{-1}\tilde{h})$  are
equivalent. By Lemma 7.3.6, the solution on $\widetilde
B_{(-b)}(\tilde{x},{(\delta^{\prime})}^{-1}\tilde{h})$ keeps
looking like a cap for $t\in [-b,\tilde{\tilde{t}})$. On the other
hand, by the definition, the surgery is always done along the
middle two-sphere of a $\delta$-neck with $\delta <
\widetilde{\delta}^{\alpha\beta}$. Then for $\beta$ large, all the
points in $\widetilde
B_{(-b)}(\tilde{x},{(\delta^{\prime})}^{-1}\tilde{h})$ are removed
(as a part of a capped horn) at the time $\tilde{\tilde{t}}$. But
$\tilde{x}$ (near the tip of the cap) exists past the time
$\tilde{\tilde{t}}$. This is a contradiction. Hence we have proved
that
$\widetilde{B}_{(-b)}(\tilde{x},{(\delta^{\prime}})^{-1}\tilde{h})$
is defined on the time interval
$[-b,0]\cap[-b,-b+(1-\tilde{\theta})\tilde{h}^{2}]$.

The $\delta'$-closeness of the solution on $\widetilde{B}_{(-b)}
(\tilde{x},{(\delta^{\prime})}^{-1}h)\times ([-b,0]\cap[-b,
-b+(1-\tilde{\theta})\tilde{h}^{2}])$ with the corresponding
subset of the standard solution follows from Lemma 7.3.6. Then we
have proved Claim 1.

\smallskip
We next verify the following

\medskip
{\bf Claim 2.} \ There is $\tilde{\theta}=\tilde{\theta}(CB)$,
$0<\tilde{\theta}<1$, such that
$b\leq(1-\tilde{\theta})\tilde{h}^{2}$ when $\beta$ large.

\medskip
Note from Proposition 7.3.3, there is a universal constant $D'>0$
such that the standard solution satisfies the following curvature
estimate
$$
R(y,s) \geq \frac{2D'}{1-s}.
$$
We choose $\tilde{\theta}= {D'}/{2(D'+CB)}$. Then for $\beta$
large enough, the rescaled solution satisfies \be
\widetilde{R}(x,t)\geq
\frac{D'}{1-(t+b)\tilde{h}^{-2}}\tilde{h}^{-2} 
\ee on
$\widetilde{B}_{(-b)}(\tilde{x},{(\delta^{\prime})}^{-1}\tilde{h})\times
([-b,0]\cap[-b,-b+(1-\tilde{\theta})\tilde{h}^{2}])$.

Suppose $b\geq (1-\tilde{\theta})\tilde{h}^{2}$. Then by combining
with the assumption $\widetilde{R}(\tilde{x},t)\leq {C}$ for
$t=(1-\tilde{\theta})\tilde{h}^{2}-b$, we have
$$
C\geq \frac{D'}{1-(t+b)\tilde{h}^{-2}}\tilde{h}^{-2},
$$
and then
$$
1\geq(1-\tilde{\theta})\(1+\frac{D'}{CB}\).
$$
This is a contradiction. Hence we have proved Claim 2.

The combination of the above two claims shows that there is a
positive constant $0 <\tilde{\theta}=\tilde{\theta}(CB)<1$ such
that for any small $\delta'>0$, there is a positive
$\beta(\delta',\varepsilon, \tilde{\theta})$ such that when $\beta
\geq \beta(\delta',\varepsilon, \tilde{\theta})$, we have $b\leq
(1-\tilde{\theta})\tilde{h}^{2}$ and the rescaled solution in the
ball
$\widetilde{B}_{(-b)}(\tilde{x},{(\delta^{\prime}})^{-1}\tilde{h})$
on the time interval $[-b,0]$ is, after scaling with factor
$\tilde{h}^{-2}$, ${\delta}^{\prime}$-close ( in the
$C^{[(\delta')^{-1}]}$ topology) to the corresponding subset of
the standard solution.

By (7.4.15) and the assumption $\widetilde{R}\leq {C}$ on
$\widetilde{B}_{0}(\bar{x},A) \times [-b,0],$ we know that the
cut-off radius $\tilde{h}$ at the time $-b$ for the rescaled
solution satisfies
$$\tilde{h} \geq \sqrt{\frac{D'}{C}}.$$

Let $\delta'>0$ be much smaller than $\varepsilon$ and
$\min\{A^{-1},A\}$. Since $\tilde{d}_0(\tilde{x},\bar{x})\leq A$,
it follows that there is constant $C(\tilde{\theta})$ depending
only on $\tilde{\theta}$ such that
$\tilde{d}_{(-b)}(\tilde{x},\bar{x})\leq C(\tilde{\theta})A \ll
(\delta')^{-1}\tilde{h}$. We now apply Lemma 7.3.5 with the
accuracy parameter ${\varepsilon}/{2}$. Let $C({\varepsilon}/{2})$
be the positive constant in Lemma 7.3.5. Without loss of
generality, we may assume the positive constant $C_1(\varepsilon)$
in the canonical neighborhood assumption is larger than
$4C({\varepsilon}/{2})$.  When $\delta'(>0)$ is much smaller than
$\varepsilon$ and $\min\{A^{-1},A\}$, the point $\bar{x}$ at the
time $\bar{t}$ has a neighborhood which is either a
$\frac{3}{4}\varepsilon$-cap or a $\frac{3}{4}\varepsilon$-neck.

Since the canonical neighborhood assumption with accuracy
parameter $\varepsilon$ is violated at $(\bar{x},\bar{t})$, the
neighborhood of the point $\bar{x}$ at the new time zero for the
rescaled solution must be a $\frac{3}{4}\varepsilon$-neck. By
Lemma 7.3.5 (b), we know the neighborhood is the slice at the time
zero of the parabolic neighborhood
$$
P(\bar{x},0,\frac{4}{3}\varepsilon^{-1}\widetilde{R}
(\bar{x},0)^{-\frac{1}{2}},
-\min\{\widetilde{R}(\bar{x},0)^{-1},b\})
$$
(with $\widetilde{R}(\bar{x},0)=1$) which is
$\frac{3}{4}\varepsilon$-close (in the
$C^{[\frac{4}{3}\varepsilon^{-1}]}$ topology) to the corresponding
subset of the evolving standard cylinder $\mathbb{S}^2 \times
\mathbb{R}$ over the time interval $[-\min \{b,1\},0]$ with scalar
curvature $1$ at the time zero. If $b \geq 1$, the
$\frac{3}{4}\varepsilon$-neck is strong, which is a contradiction.
While if $b<1$, the $\frac{3}{4}\varepsilon$-neck at time $-b$ is
contained in the union of the gluing cap and the adjacent
$\delta$-neck where the $\delta$-cutoff surgery took place. Since
$\varepsilon$ is small (say $\varepsilon< 1/100$), it is clear
that the point $\bar{x}$ at time $-b$ is the center of an
$\varepsilon$-neck which is entirely contained in the adjacent
$\delta$-neck. By the proof of Lemma 7.3.2, the adjacent
$\delta$-neck approximates an ancient $\kappa$-solution. This
implies the point $\bar{x}$ at the time $\bar{t}$ has a strong
$\varepsilon$-neck, which is also a contradiction.

Hence we have proved that there exists $\beta_0=\beta_0(
\varepsilon, A, B, C)$ such that when $\beta\geq\beta_0$, the
rescaled solution on the ball $\widetilde{B}_{0}(\bar{x},A)$ can
be defined before the time $-b$.

Let $[t_{A}^{\alpha\beta},0]\supset[-b,0]$ be the largest time
interval so that the rescaled solution
$\widetilde{g}^{\alpha\beta}_{ij}$ can be defined on
$\widetilde{B}_{0}(\bar{x},A)\times[t_{A}^{\alpha\beta},0]$. We
finally claim that $t_{A}^{\alpha\beta}\le
-b-\frac{1}{8}\eta^{-1}C^{-1}$ for $\beta$ large enough.

Indeed, suppose not, by the gradient estimates as in Step 1, we
have the curvature estimate
$$
\widetilde{R}(x,t)\leq {10C}
$$
on $\widetilde{B}_{0}(\bar{x},A)\times [t_{A}^{\alpha\beta},-b]$.
Hence we have the curvature estimate
$$
\widetilde{R}(x,t)\leq {10C}
$$
on $\widetilde{B}_{0}(\bar{x},A)\times [t_{A}^{\alpha\beta},0]$.
By the above argument there is a $\beta_0=\beta_0( \varepsilon, A,
B + \frac{1}{8}\eta^{-1}C^{-1}, 10C)$ such that for
$\beta\geq\beta_0$, the solution in the ball
$\widetilde{B}_{0}(\bar{x},A)$ can be defined before the time
$t_{A}^{\alpha\beta}$. This is a contradiction.

Therefore we have proved Assertion 1.

\medskip
{\bf Assertion 2.} \ For arbitrarily fixed $\alpha$,
$0<A<+\infty$, $1 \leq C<+\infty$ and $0 < B
<\frac{1}{2}\varepsilon^2 (r^{\alpha})^{-2}-
\frac{1}{50}\eta^{-1}$, there is a $\beta_0=\beta_0(\varepsilon,
A,B,C)$ (independent of $\alpha$) such that if $\beta\geq\beta_0$
and the rescaled solution $\widetilde{g}^{\alpha\beta}_{ij}$ on
the ball $\widetilde{B}_{0}(\bar{x},A)$ is defined on a time
interval $[-b+\epsilon',0]$ with $0 < b \leq B$ and $0<\epsilon'
<\frac{1}{50}\eta^{-1}$ and the scalar curvature satisfies
$$
\widetilde{R}(x,t)\leq {C} \ \ \mbox{ on } \
\widetilde{B}_{0}(\bar{x},A)\times [-b+\epsilon',0],
$$
 and there is a point $y\in \widetilde{B}_{0}(\bar{x},A)$ such
that $\widetilde{R}(y,-b+\epsilon')\leq \frac{3}{2}$, then the
rescaled solution $\widetilde{g}^{\alpha\beta}_{ij}$ at $y$ is
also defined on the extended time interval
$[-b-\frac{1}{50}\eta^{-1},0]$ and satisfies the estimate
$$
\widetilde{R}(y,t)\leq 15
$$
for $t\in[-b-\frac{1}{50}\eta^{-1},-b+\epsilon']$.

\medskip
{\bf \em Proof of Assertion} {\bf 2.} \ We imitate the proof of
Assertion 1. If the rescaled solution
$\widetilde{g}^{\alpha\beta}_{ij}$ at $y$ cannot be defined for
some time in $[-b-\frac{1}{50}\eta^{-1},-b+\epsilon')$, then there
is a surgery at some time $\tilde{\tilde{t}} \in
[-b-\frac{1}{50}\eta^{-1},-b+\epsilon']$ such that $y$ lies in the
instant gluing cap.  Let $\tilde{h}$
$(=R(\bar{x},\bar{t})^{\frac{1}{2}}h$) be the cutoff radius at the
time $\tilde{\tilde{t}}$ for the rescaled solution. Clearly, there
is a universal constant $D>1$ such that
${D}^{-1}\tilde{h}\leq\widetilde{R}(y,\tilde{\tilde{t}})^{-\frac{1}{2}}
\leq D\tilde{h}$. By the gradient estimates as in Step 1, the
cutoff radius satisfies
$$
\tilde{h} \geq D^{-1}15^{-\frac{1}{2}}.
$$

As in Claim 1 (i) in the proof of Assertion 1, for any small
constants $0<\tilde{\theta}<\frac{1}{2}$, $\delta'>0$, there
exists a $\beta(\delta',\varepsilon, \tilde{\theta})> 0$ such that
for $\beta \geq \beta(\delta',\varepsilon, \tilde{\theta})$, there
is no surgery interfering in
$\widetilde{B}_{\tilde{\tilde{t}}}(y,(\delta')^{-1}\tilde{h})\times
([\tilde{\tilde{t}},(1-\tilde{\theta})\tilde{h}^{2}
+\tilde{\tilde{t}}]\cap(\tilde{\tilde{t}},0])$. Without loss of
generality, we may assume that the universal constant $\eta$ is
much larger than $D$. Then we have
$(1-\tilde{\theta})\tilde{h}^{2}
+\tilde{\tilde{t}}>-b+\frac{1}{50}\eta^{-1}$. As in Claim 2 in the
proof of Assertion 1, we can use the curvature bound assumption to
choose $\tilde{\theta}=\tilde{\theta}(B,C)$ such that
$(1-\tilde{\theta})\tilde{h}^{2} +\tilde{\tilde{t}}\geq 0$;
otherwise
$$
C\geq\frac{D'}{\tilde{\theta}\tilde{h}^{2}}
$$
for some universal constant $D'>1$, and
$$
|\tilde{\tilde{t}}+b|\leq \frac{1}{50}\eta^{-1},
$$
which implies
$$
1\geq(1-\tilde{\theta})\(1+\frac{D'}{C\(B+\frac{1}{50}\eta^{-1}\)}\).
$$
This is a contradiction if we choose
$\tilde{\theta}={D'}/{2(D'+C(B+\frac{1}{50}\eta^{-1}))}$.

So there is a positive constant $0
<\tilde{\theta}=\tilde{\theta}(B,C)<1$ such that for any
$\delta'>0$, there is a positive $\beta(\delta',\varepsilon,
\tilde{\theta})$ such that when $\beta \geq
\beta(\delta',\varepsilon, \tilde{\theta})$, we have
$-\tilde{\tilde{t}}\leq (1-\tilde{\theta})\tilde{h}^{2}$  and the
solution in the ball $\widetilde{B}_{\tilde{\tilde{t}}}(\tilde{x},
{(\delta^{\prime})}^{-1}\tilde{h})$ on the time interval
$[\tilde{\tilde{t}},0]$ is, after scaling with factor
$\tilde{h}^{-2}$, ${\delta}^{\prime}$-close (in the
$C^{[\delta'^{-1}]}$ topology) to the corresponding subset of the
standard solution.

Then exactly as in the proof of Assertion 1, by using the
canonical neighborhood structure of the standard solution in Lemma
7.3.5, this gives the desired contradiction with the hypothesis
that the canonical neighborhood assumption with accuracy parameter
$\varepsilon$ is violated at $(\bar{x},\bar{t})$, for $\beta$
sufficiently large.

The curvature estimate at the point $y$ follows from Step 1.
Therefore the proof of Assertion 2 is complete.

Note that the standard solution satisfies $R(x_1,t)\leq D''
R(x_2,t)$ for any $t\in [0,\frac{1}{2}]$ and any two points
$x_1,x_2$, where $D''\geq 1$ is a universal constant.

\medskip
{\bf Assertion 3.} \ For arbitrarily fixed $\alpha$,
$0<A<+\infty$, $1 \leq C<+\infty$, there is a
$\beta_0=\beta_0(\varepsilon, AC^{\frac{1}{2}})$ such that if any
point $(y_0,t_0)$  with  $0 \leq -t_0 <\frac{1}{2}\varepsilon^2
(r^{\alpha})^{-2}- \frac{1}{8}\eta^{-1}C^{-1}$ of the rescaled
solution $\widetilde{g}^{\alpha\beta}_{ij}$ for $\beta\geq\beta_0$
satisfies $\widetilde{R}(y_0,t_0)\leq C$ , then either the
rescaled solution at $y_0$ can be defined at least on
$[t_0-\frac{1}{16}\eta^{-1}C^{-1},t_0]$ and the rescaled scalar
curvature satisfies
$$
\widetilde{R}(y_0,t)\leq 10 {C} \; \mbox{ for } t\in
\Big[t_0-\frac{1}{16}\eta^{-1}C^{-1},t_0\Big],
$$
or we have
$$
\widetilde{R}(x_1,t_0)\leq 2D'' \widetilde{R}(x_2,t_0)
$$
for any two points $x_1, x_2\in \widetilde{B}_{t_0}(y_0,A)$, where
$D''$ is the above universal constant.

\medskip
{\bf \em Proof of Assertion} {\bf 3.} \ Suppose the rescaled
solution $\widetilde{g}^{\alpha\beta}_{ij}$ at $y_0$ cannot be
defined for some $t\in [t_0-\frac{1}{16}\eta^{-1}C^{-1},t_0)$;
then there is a surgery at some time $\tilde{t} \in
[t_0-\frac{1}{16}\eta^{-1}C^{-1},t_0]$ such that $y_0 $ lies in
the instant gluing cap. Let $\tilde{h}$
$(=R(\bar{x},\bar{t})^{\frac{1}{2}}h$) be the cutoff radius at the
time $\tilde{t}$ for the rescaled solution
$\widetilde{g}^{\alpha\beta}_{ij}$.  By the gradient estimates as
in Step 1, the cutoff radius satisfies
$$
\tilde{h} \geq D^{-1}10^{-\frac{1}{2}}C^{-\frac{1}{2}},
$$
where $D$ is the universal constant in the proof of the Assertion
1. Since we assume $\eta$ is suitably larger than $D$ as before,
we have $\frac{1}{2}\tilde{h}^{2} +\tilde{t}>t_0$. As in Claim 1
(ii) in the proof of Assertion 1, for arbitrarily small
$\delta'>0$, we know that for $\beta$ large enough the rescaled
solution on
$\widetilde{B}_{\tilde{t}}(y_0,{(\delta^{\prime})}^{-1}\tilde{h})\times
[\tilde{t},t_0]$ is, after scaling with factor $\tilde{h}^{-2}$,
${\delta}^{\prime}$-close (in the $C^{[(\delta')^{-1}]}$ topology)
to the corresponding subset of the standard solution. Since
$(\delta')^{-1}\tilde{h}\gg A$ for $\beta$ large enough, Assertion
3 follows from the curvature estimate of standard solution in the
time interval $ [0,\frac{1}{2}]$.

\medskip
{\it Step} 3. \ For any subsequence $(\alpha_k,\beta_k)$ of
$(\alpha,\beta)$ with $r^{\alpha_k}\rightarrow 0$ and
$\delta^{\alpha_k\beta_k}\rightarrow 0$ as $k\rightarrow \infty$,
we next argue as in the second step of the proof of Theorem 7.1.1
to show that the curvatures of the rescaled solutions
$\tilde{g}^{\alpha_k\beta_k}_{ij}$ at the new times zero (after
shifting) stay uniformly bounded at bounded distances from
$\bar{x}$ for all sufficiently large $k$. More precisely, we will
prove the following assertion:

\medskip
{\bf Assertion 4.} \ Given any subsequence of the rescaled
solutions $\tilde{g}^{\alpha_k\beta_k}_{ij}$ with
$r^{\alpha_k}\rightarrow 0$ and
$\delta^{\alpha_k\beta_k}\rightarrow 0$ as $k\rightarrow \infty$,
then for any $L>0$, there are constants $C(L)>0$ and $k(L)$ such
that the rescaled solutions $\tilde{g}^{\alpha_k\beta_k}_{ij}$
satisfy
\begin{itemize}
\item[(i)] $\tilde{R}(x,0)\leq C(L)$ for all points $x$ with
$\tilde{d}_{0}(x,\bar{x})\leq L$ and all $k \geq 1$; \item[(ii)]
the rescaled solutions over the ball $\tilde{B}_{0}(\bar{x},L)$
are defined at least on the time interval
$[-\frac{1}{16}\eta^{-1}C(L)^{-1},0]$ for all $k\geq k(L)$.
\end{itemize}

\medskip
{\bf \em Proof of Assertion} {\bf 4.} \ For each $\rho>0$, set
\begin{multline*}
M(\rho)=\sup\Big\{\tilde{R}(x,0)\ |\ k \geq 1\; \mbox{ and }\;
\tilde{d}_0(x,\bar{x})\leq\rho \\
\mbox{in the rescaled solutions } \;
\tilde{g}^{\alpha_k\beta_k}_{ij}\Big\}
\end{multline*}
and
$$
\rho_0=\sup\{\rho > 0\ |\ M(\rho)<+\infty\}.
$$
Note that the estimate (7.4.14) implies that $\rho_0>0$. For (i),
it suffices to prove $\rho_0=+\infty$.

We argue by contradiction. Suppose $\rho_0<+\infty$. Then there is
a sequence of points $y$ in the rescaled solutions
$\tilde{g}^{\alpha_k\beta_k}_{ij}$ with
$\tilde{d}_0(\bar{x},y)\rightarrow\rho_0<+\infty$ and
$\tilde{R}(y,0)\rightarrow +\infty$. Denote by $\gamma$ a
minimizing geodesic segment from $\bar{x}$ to $y$ and denote by
$\tilde{B}_0(\bar{x},\rho_0)$ the open geodesic ball centered at
$\bar{x}$ of radius $\rho_0$ on the rescaled solution
$\tilde{g}^{\alpha_k\beta_k}_{ij}$.

First, we claim that for any $0<\rho<\rho_0$ with $\rho$ near
$\rho_0$, the rescaled solutions on the balls
$\tilde{B}_{0}(\bar{x},\rho)$ are defined on the time interval
$[-\frac{1}{16}\eta^{-1}M(\rho)^{-1},0]$ for all large $k$.
Indeed, this follows from Assertion 3 or Assertion 1. For the
later purpose in Step 6, we now present an argument by using
Assertion 3.  If the claim is not true, then there is a surgery at
some time $\tilde{t} \in [-\frac{1}{16}\eta^{-1}M(\rho)^{-1},0]$
such that some point $\tilde{y}\in \tilde{B}_{0}(\bar{x},\rho) $
lies in the instant gluing cap.  We can choose sufficiently small
$\delta'>0$ such that $2\rho_0<(\delta')^{-\frac{1}{2}}\tilde{h}$,
where $\tilde{h} \geq
D^{-1}20^{-\frac{1}{2}}M(\rho)^{-\frac{1}{2}}$ is the cutoff
radius of the rescaled solutions at $\tilde{t}$. By applying
Assertion 3 with $(\tilde{y},0)=(y_0,t_0)$, we see that there is a
$k(\rho_0,M(\rho))>0$ such that when $k\geq k(\rho_0,M(\rho))$,
$$
\widetilde{R}(x,0)\leq 2D''
$$
for all $x\in \widetilde{B}_{0}(\bar{x},\rho)$. This is a
contradiction as $\rho\rightarrow\rho_0$.

Since for each fixed $0<\rho<\rho_0$ with $\rho$ near $\rho_0$,
the rescaled solutions are defined on
$\tilde{B}_{0}(\bar{x},\rho)\times
[-\frac{1}{16}\eta^{-1}M(\rho)^{-1},0]$ for all large $k$, by Step
1 and Shi's derivative estimate, we know that the covariant
derivatives and higher order derivatives of the curvatures on
$\tilde{B}_{0}(\bar{x},\rho -
\frac{(\rho_0-\rho)}{2})\times[-\frac{1}{32}\eta^{-1}M(\rho)^{-1},0]$
are also uniformly bounded.

By the uniform $\kappa$-noncollapsing property and Hamilton's
compactness theorem (Theorem 4.1.5), after passing to a
subsequence, we can assume that the marked sequence
$(\tilde{B}_0(\bar{x},\rho_0),\widetilde{g}^{\alpha_k\beta_k}_{ij},
\bar{x})$ converges in the $C^{\infty}_{loc}$ topology to a marked
(noncomplete) manifold
($B_{\infty},\widetilde{g}^{\infty}_{ij},\bar{x})$ and the
geodesic segments $\gamma$ converge to a geodesic segment (missing
an endpoint) $\gamma_{\infty}\subset B_{\infty}$ emanating from
$\bar{x}$.

Clearly, the limit has nonnegative sectional curvature by the
pinching assumption. Consider a tubular neighborhood along
$\gamma_{\infty}$ defined by
$$
V=\bigcup_{q_0\in\gamma_{\infty}}B_{\infty}
(q_0,4\pi(\widetilde{R}_{\infty}(q_0))^{-\frac{1}{2}}),
$$
where $\widetilde{R}_{\infty}$ denotes the scalar curvature of the
limit and
$$
B_{\infty}(q_0,4\pi(\widetilde{R}_{\infty}(q_0))^{-\frac{1}{2}})
$$
is the ball centered at $q_0\in B_{\infty}$ with the radius
$4\pi(\widetilde{R}_{\infty}(q_0))^{-\frac{1}{2}}$. Let
$\bar{B}_{\infty}$ denote the completion of
$(B_{\infty},\widetilde{g}^{\infty}_{ij})$, and
$y_{\infty}\in\bar{B}_{\infty}$ the limit point of
$\gamma_{\infty}$. Exactly as in the second step of the proof of
Theorem 7.1.1, it follows from the canonical neighborhood
assumption with accuracy parameter $2\varepsilon$ that the
limiting metric $\widetilde{g}^{\infty}_{ij}$ is cylindrical at
any point $q_0\in\gamma_{\infty}$ which is sufficiently close to
$y_{\infty}$ and then the metric space $\bar{V}=V\cup
\{y_{\infty}\}$ by adding the point $y_{\infty}$ has nonnegative
curvature in the Alexandrov sense. Consequently we have a
three-dimensional non-flat tangent cone $C_{y_{\infty}}\bar{V}$ at
$y_{\infty}$ which is a metric cone with aperture $\leq
20\varepsilon$.

On the other hand, note that by the canonical neighborhood
assumption, the canonical $2\varepsilon$-neck neighborhoods are
strong. Thus at each point $q\in V$ near $y_{\infty}$, the
limiting metric $\widetilde{g}^{\infty}_{ij}$ actually exists on
the whole parabolic neighborhood
$$
V \bigcap
P\(q,0,\frac{1}{3}\eta^{-1}(\widetilde{R}_{\infty}(q))^{-\frac{1}{2}},
-\frac{1}{10}\eta^{-1}(\widetilde{R}_{\infty}(q))^{-1}\),
$$
and is a smooth solution of the Ricci flow there. Pick $z\in
C_{y_{\infty}}\bar{V}$ with distance one from the vertex
$y_{\infty}$ and it is nonflat around $z$. By definition the ball
$B(z,\frac{1}{2})\subset C_{y_{\infty}}\bar{V}$ is the
Gromov-Hausdorff convergent limit of the scalings of a sequence of
balls $B_{\infty}(z_{\ell},\sigma_{\ell})
(\subset(V,\widetilde{g}^{\infty}_{ij}))$ where
$\sigma_{\ell}\rightarrow0$. Since the estimate (7.4.14) survives
on $(V,\widetilde{g}^{\infty}_{ij})$ for all $A < +\infty$, and
the tangent cone is three-dimensional and nonflat around $z$, we
see that this convergence is actually in the $C^{\infty}_{\rm
loc}$ topology and over some ancient time interval. Since the
limiting $B_{\infty}(z,\frac{1}{2})(\subset
C_{y_{\infty}}\bar{V})$ is a piece of nonnegatively curved nonflat
metric cone, we get a contradiction with Hamilton's strong maximum
principle (Theorem 2.2.1) as before. So we have proved
$\rho_0=\infty$. This proves (i).

By the same proof of Assertion 1 in Step 2, we can further show
that for any $L$, the rescaled solutions on the balls
$\tilde{B}_{0}(\bar{x},L)$ are defined at least on the time
interval $[-\frac{1}{16}\eta^{-1}C(L)^{-1},0]$ for all
sufficiently large $k$. This proves (ii).

\medskip
{\it Step} 4. \ For any subsequence $(\alpha_k,\beta_k)$ of
$(\alpha,\beta)$ with $r^{\alpha_k}\rightarrow 0$ and
$\widetilde{\delta}^{\alpha_k\beta_k}\rightarrow 0$ as
$k\rightarrow \infty$, by Step 3, the $\kappa$-noncollapsing
property and Hamilton's compactness theorem, we can extract a
$C^{\infty}_{loc}$ convergent subsequence of
$\tilde{g}^{\alpha_k\beta_k}_{ij}$ over some space-time open
subsets containing the slice $\{t=0\}$. We now want to show
\textbf{any} such limit has bounded curvature at $t=0$. We prove
by contradiction. Suppose not, then there is a sequence of points
$z_{\ell}$ divergent to infinity in the limiting metric at time
zero with curvature divergent to infinity. Since the curvature at
$z_{\ell}$ is large (comparable to one), $z_{\ell}$ has a
canonical neighborhood which is a $2\varepsilon$-cap or strong
$2\varepsilon$-neck. Note that the boundary of $2\varepsilon$-cap
lies in some $2\varepsilon$-neck. So we get a sequence of
$2\varepsilon$-necks with radius going to zero. Note also that the
limit has nonnegative sectional curvature. Without loss of
generality, we may assume $2\varepsilon < \varepsilon_0$, where
$\varepsilon_0$ is the positive constant in Proposition 6.1.1.
Thus this arrives at a contradiction with Proposition 6.1.1.

\medskip
{\it Step} 5. \ In this step, we will choose some subsequence
$(\alpha_k,\beta_k)$ of $(\alpha,\beta)$ so that we can extract a
complete smooth limit of the rescaled solutions
$\widetilde{g}^{\alpha_k \beta_k}_{ij}$ to the Ricci flow with
surgery on a time interval $[-a,0]$ for some $a>0$.

Choose $\alpha_k,\beta_k\rightarrow\infty$ so that
$r^{\alpha_k}\rightarrow0$,
$\widetilde{\delta}^{\alpha_{k}\beta_{k}}\rightarrow 0$, and
Assertion 1, 2, 3 hold with $\alpha=\alpha_k, \beta=\beta_k$ for
all $A \in \{{p}/{q} \ |\ p, q=1, 2, \ldots, k\}$, and $B,C \in
\{1,2,\ldots,k\}$. By Step 3, we may assume the rescaled solutions
$\widetilde{g}^{\alpha_k \beta_k}_{ij}$ converge in the
$C^{\infty}_{\rm loc}$ topology at the time $t=0$. Since the
curvature of the limit at $t=0$ is bounded by Step 4, it follows
from Assertion 1 in Step 2 and the choice of the sequence
$(\alpha_k,\beta_k)$ that the limiting
$(M_{\infty},\widetilde{g}^{\infty}_{ij}(\cdot,t))$ is defined at
least on a backward time interval $[-a,0]$ for some positive
constant $a$ and is a smooth solution to the Ricci flow there.

\medskip
{\it Step} 6. \ We further want to extend the limit in Step 5
backwards in time to infinity to get an ancient $\kappa$-solution.
Let $\widetilde{g}^{\alpha_k \beta_k}_{ij}$ be the convergent
sequence obtained in the above Step 5.

Denote by
\begin{eqnarray*}
t_{\max} = \sup \Big\{\ t'&  |&  \mbox {we can take a
smooth limit on } (-t',0] \mbox{ (with bounded}  \\
& &   \mbox{curvature at each time slice) from a subsequence
of}\\
&   &   \mbox{the rescaled solutions
}\widetilde{g}^{\alpha_k\beta_k}_{ij} \Big\}.
\end{eqnarray*}
We first claim that there is a subsequence of the rescaled
solutions $\widetilde{g}^{\alpha_k\beta_k}_{ij}$ which converges
in the $C^{\infty}_{\rm loc}$ topology to a smooth limit
$(M_{\infty},\widetilde{g}^{\infty}_{ij}(\cdot,t))$ on the maximal
time interval $(-t_{\max},0]$.

Indeed, let $t_{\ell}$ be a sequence of positive numbers such that
$t_{\ell}\rightarrow t_{\max}$ and there exist smooth limits
$(M_{\infty},\widetilde{g}^{\infty}_{{\ell}}(\cdot,t))$ defined on
$(-t_{\ell},0]$. For each ${\ell}$, the limit has nonnegative
sectional curvature and has bounded curvature at each time slice.
Moreover by the gradient estimate in canonical neighborhood
assumption with accuracy parameter $2\varepsilon$, the limit has
bounded curvature on each subinterval
$[-b,0]\subset(-t_{\ell},0]$. Denote by $\widetilde{Q}$ the scalar
curvature upper bound of the limit at time zero ($\widetilde{Q}$
is independent of ${\ell}$). Then we can apply Li-Yau-Hamilton
inequality (Corollary 2.5.5) to get
$$
 \widetilde{R}^{\infty}_{{\ell}}(x,t)\leq
 \frac{t_{\ell}}{t+t_{\ell}}\widetilde{Q},
$$
where $\widetilde{R}^{\infty}_{{\ell}}(x,t)$ are the scalar
curvatures of the limits
$(M_{\infty},\widetilde{g}^{\infty}_{{\ell}}(\cdot,t))$. Hence by
the definition of convergence and the above curvature estimates,
we can find a subsequence of the rescaled solutions
$\widetilde{g}^{\alpha_k\beta_k}_{ij}$ which converges in the
$C^{\infty}_{loc}$ topology to a smooth limit
$(M_{\infty},\widetilde{g}^{\infty}_{ij}(\cdot,t))$ on the maximal
time interval $(-t_{\max},0]$.

We need to show $-t_{\max}=-\infty$. Suppose $-t_{\max}>-\infty$,
there are only the following two possibilities: either
\begin{itemize}
\item[(1)] The curvature of the limiting solution
$(M_{\infty},\widetilde{g}^{\infty}_{ij}(\cdot,t))$ becomes
unbounded as $t\searrow -t_{\max}$; or \item[(2)] For each small
constant $\theta>0$ and each large integer $k_0>0$, there is some
$k\geq k_0$ such that the rescaled solution
$\widetilde{g}^{\alpha_k\beta_k}_{ij}$ has a surgery time
$T_k\in[-t_{\max}-\theta,0]$ and a surgery point $x_k$ lying in a
gluing cap at the times $T_k$ so that $d^2_{T_k}(x_k,\bar{x})$ is
uniformly bounded from above by a constant independent of $\theta$
and $k_0$.
\end{itemize}

We next claim that the possibility (1) always occurs. Suppose not;
then the curvature of the limiting solution
$(M_{\infty},\widetilde{g}^{\infty}_{ij}(\cdot,t))$ is bounded on
$M_{\infty}\times(-t_{\max},0]$ by some positive constant
$\hat{C}$. In particular, for any $A>0$, there is a sufficiently
large integer $k_1>0$ such that any rescaled solution
$\widetilde{g}^{\alpha_k\beta_k}_{ij}$ with $k\geq k_1$ on the
geodesic ball $\widetilde{B}_{0}(\bar{x},A)$ is defined on the
time interval $[-t_{\max}+\frac{1}{50}\eta^{-1}{\hat{C}}^{-1},0]$
and its scalar curvature is bounded by $2\hat{C}$ there. (Here,
without loss of generality, we may assume that the upper bound
$\hat{C}$ is so large that
$-t_{\max}+\frac{1}{50}\eta^{-1}{\hat{C}}^{-1} < 0$.) By Assertion
1 in Step 2, for $k$ large enough, the rescaled solution
$\widetilde{g}^{\alpha_k \beta_k}_{ij}$ over
$\widetilde{B}_{0}(\bar{x},A)$ can be defined on the extended time
interval $[-t_{\max}-\frac{1}{50}\eta^{-1}{\hat{C}}^{-1},0]$ and
has the scalar curvature $\widetilde{R}\leq 10 \hat{C}$ on
$\widetilde{B}_{0}(\bar{x},A)\times
[-t_{\max}-\frac{1}{50}\eta^{-1}{\hat{C}}^{-1},0]$. So we can
extract a smooth limit from the sequence to get the limiting
solution which is defined on a larger time interval
$[-t_{\max}-\frac{1}{50}\eta^{-1}{\hat{C}}^{-1},0]$. This
contradicts the definition of the maximal time $-t_{\max}$.

It remains to exclude the possibility (1).

By using Li-Yau-Hamilton inequality (Corollary 2.5.5) again, we
have
$$
\widetilde{R}_{\infty}(x,t)\leq
\frac{t_{\max}}{t+t_{\max}}\widetilde{Q}.
$$
So we only need to control the curvature near $-t_{\max}$. Exactly
as in Step 4 in the proof of Theorem 7.1.1, it follows from
Li-Yau-Hamilton inequality that \be \tilde{d}_0(x,y)\leq
\tilde{d}_t(x,y)\leq \tilde{d}_0(x,y)
+30t_{\max}\sqrt{\widetilde{Q}}
\ee for any $x,y\in M_{\infty}$ and $t\in(-t_{\max},0]$.

Since the infimum of the scalar curvature is nondecreasing in
time, we have some point $y_{\infty}\in M_{\infty}$  and some time
$-t_{\max} < t_{\infty} < -t_{\max}+\frac{1}{50}\eta^{-1}$ such
that $\widetilde{R}_{\infty}(y_{\infty},t_{\infty})<{5}/{4}$. By
(7.4.16), there is a constant $\widetilde{A}_0>0$ such that
$\tilde{d}_t(\bar{x},y_{\infty})\leq \widetilde{A}_0/2$ for all $t
\in (-t_{\max},0]$.

Now we come back to the rescaled solution
$\widetilde{g}^{\alpha_k\beta_k}_{ij}$. Clearly, for arbitrarily
given small $\epsilon' > 0$, when $k$ large enough, there is a
point $y_k$ in the underlying manifold of
$\widetilde{g}^{\alpha_k\beta_k}_{ij}$ at time $0$ satisfying the
following properties \be
\widetilde{R}(y_k,t_{\infty})<\frac{3}{2},\qquad
\widetilde{d}_{t}(\bar{x},y_k)\leq \widetilde{A}_0 
\ee for $t\in [-t_{\max} + \epsilon',0]$. By the definition of
convergence, we know that for any fixed $A_0 \geq
2\widetilde{A}_0$, for $k$ large enough, the rescaled solution
over $\widetilde{B}_{0}(\bar{x},A_0)$ is defined on the time
interval $[t_{\infty},0]$ and satisfies
$$
\widetilde{R}(x,t)\leq \frac{2t_{\max}}{t+t_{\max}}\widetilde{Q}
$$
on $\widetilde{B}_{0}(\bar{x},A_0)\times [t_{\infty},0]$. Then by
Assertion 2 of Step 2, we have proved that there is a sufficiently
large integer $\bar{k}_0$ such that when $k \geq \bar{k}_0$, the
rescaled solutions $\widetilde{g}^{\alpha_k\beta_k}_{ij}$ at $y_k$
can be defined on $[-t_{\max}-\frac{1}{50}\eta^{-1},0]$, and
satisfy
$$
\widetilde{R}(y_k,t)\leq 15
$$
for $t\in[-t_{\max}-\frac{1}{50}\eta^{-1},t_{\infty}]$.

We now prove a statement analogous to Assertion 4 (i) of Step 3.

\medskip
{\bf Assertion 5.} \ For the above rescaled solutions
$\widetilde{g}^{\alpha_k\beta_k}_{ij}$ and $\bar{k}_0$, we have
that for any $L>0$, there is a positive constant $\omega(L)$ such
that the rescaled solutions $\widetilde{g}^{\alpha_k\beta_k}_{ij}$
satisfy
$$
\widetilde{R}(x,t)\leq\omega(L)
$$
for all $(x,t)$ with $\tilde{d}_{t}(x,y_k)\leq L$ and $t\in
[-t_{\max}-\frac{1}{50}\eta^{-1},t_{\infty}]$, and for all $k \geq
\bar{k}_0$.

\medskip
{\bf \em Proof of Assertion} {\bf 5.} \ We slightly modify the
argument in the proof of Assertion 4 (i). Let

\begin{eqnarray*}
M(\rho)=\sup\Big\{\widetilde{R}(x,t) &  |&
\tilde{d}_t(x,y_k)\leq\rho\; \mbox{ and }\; t \in
[-t_{\max}-\frac{1}{50}\eta^{-1},t_{\infty}]\\
&   &   \mbox{in the rescaled solutions } \;
\widetilde{g}^{\alpha_k\beta_k}_{ij}, k\geq \bar{k}_0 \Big\}
\end{eqnarray*}
and
$$
\rho_0=\sup\{\rho > 0\ |\  M(\rho)<+\infty\}.
$$
Note that the estimate (7.4.14) implies that $\rho_0>0$. We only
need to show $\rho_0 = +\infty$.

We argue by contradiction. Suppose $\rho_0<+\infty$. Then, after
passing to a subsequence, there is a sequence $(\tilde{y}_k,t_k)$
in the rescaled solutions $\widetilde{g}^{\alpha_k\beta_k}_{ij}$
with $t_k \in [-t_{\max}-\frac{1}{50}\eta^{-1},t_{\infty}]$ and
$\tilde{d}_{t_k}(y_k,\tilde{y}_k)\rightarrow\rho_0<+\infty$ such
that $\widetilde{R}(\tilde{y}_k,t_k)\rightarrow +\infty$. Denote
by $\gamma_k$ a minimizing geodesic segment from $y_k$ to
$\tilde{y}_k$ at the time $t_k$ and denote by
$\widetilde{B}_{t_k}(y_k,\rho_0)$ the open geodesic ball centered
at $y_k$ of radius $\rho_0$ on the rescaled solution
$\widetilde{g}^{\alpha_k\beta_k}_{ij}(\cdot,t_k)$.

For any $0<\rho<\rho_0$ with $\rho$ near $\rho_0$, by applying
Assertion 3 as before, we get that  the rescaled solutions on the
balls $\widetilde{B}_{t_k}(y_k,\rho)$ are defined on the time
interval $[t_k-\frac{1}{16}\eta^{-1}M(\rho)^{-1},t_k]$ for all
large $k$. By Step 1 and Shi's derivative estimate, we further
know that the covariant derivatives and higher order derivatives
of the curvatures on
$\widetilde{B}_{t_k}(y_k,\rho-\frac{(\rho_0-\rho)}{2})
\times[t_k-\frac{1}{32}\eta^{-1}M(\rho)^{-1},t_k]$ are also
uniformly bounded. Then by the uniform $\kappa$-noncollapsing
property and Hamilton's compactness theorem (Theorem 4.1.5), after
passing to a subsequence, we can assume that the marked sequence
$(\tilde{B}_{t_k}(y_k,\rho_0),\widetilde{g}^{\alpha_k\beta_k}_{ij}
(\cdot,t_k),y_k)$ converges in the $C^{\infty}_{loc}$ topology to
a marked (noncomplete) manifold
($B_{\infty},\widetilde{g}^{\infty}_{ij},y_{\infty})$ and the
geodesic segments $\gamma_k$ converge to a geodesic segment
(missing an endpoint) $\gamma_{\infty}\subset B_{\infty}$
emanating from $y_{\infty}$.

Clearly, the limit also has nonnegative sectional curvature by the
pinching assumption. Then by repeating the same argument as in the
proof of Assertion 4 (i) in the rest, we derive a contradiction
with Hamilton's strong maximum principle. This proves Assertion 5.

We then apply the second estimate of (7.4.17) and Assertion 5 to
conclude that for any large constant $0<A<+\infty$, there is a
positive constant $C(A)$ such that for any small $\epsilon'>0$,
the rescaled solutions $\widetilde{g}^{\alpha_k\beta_k}_{ij}$
satisfy \be
\widetilde{R}(x,t)\leq C(A), 
\ee for all $x \in \widetilde{B}_{0}(\bar{x},A)$ and $t\in
[-t_{\max} + \epsilon',0]$, and for all sufficiently large $k$.
Then by applying Assertion 1 in Step 2, we conclude that the
rescaled solutions $\widetilde{g}^{\alpha_k\beta_k}_{ij}$ on the
geodesic balls $\widetilde{B}_{0}(\bar{x},A)$ are also defined on
the extended time interval $ [-t_{\max} +
\epsilon'-\frac{1}{8}\eta^{-1}C(A)^{-1},0]$ for all sufficiently
large $k$. Furthermore, by the gradient estimates as in Step 1, we
have
$$
\widetilde{R}(x,t)\leq 10C(A),
$$
for $x \in \widetilde{B}_{0}(\bar{x},A)$ and $t\in [-t_{\max} +
\epsilon'-\frac{1}{8}\eta^{-1}C(A)^{-1},0]$. Since $\epsilon'>0$
is arbitrarily small and the positive constant $C(A)$ is
independent of $\epsilon'$, we conclude that the rescaled
solutions $\widetilde{g}^{\alpha_k\beta_k}_{ij}$ on
$\widetilde{B}_{0}(\bar{x},A)$ are defined on the extended time
interval $ [-t_{\max} -\frac{1}{16}\eta^{-1}C(A)^{-1},0]$ and
satisfy \be
\widetilde{R}(x,t)\leq 10C(A), 
\ee for $x \in \widetilde{B}_{0}(\bar{x},A)$ and $t\in [-t_{\max}
-\frac{1}{16}\eta^{-1}C(A)^{-1},0]$, and for all sufficiently
large $k$.

Now, by taking convergent subsequences from the rescaled solutions
$\widetilde{g}^{\alpha_k\beta_k}_{ij}$, we see that the limit
solution is defined smoothly on a space-time open subset of
$M_{\infty}\times (-\infty,0]$ containing $M_{\infty}\times
[-t_{\max},0]$. By Step 4, we see that the limiting metric
$\widetilde{g}^{\infty}_{ij}(\cdot,-t_{\max})$ at time $-t_{\max}$
has bounded curvature. Then by combining with the canonical
neighborhood assumption of accuracy $2\varepsilon$, we conclude
that the curvature of the limit is uniformly bounded on the time
interval $[-t_{\max},0]$. So we have excluded the possibility (1).

Hence we have proved a subsequence of the rescaled solutions
converges to an orientable ancient $\kappa$-solution.

Finally by combining with the canonical neighborhood theorem
(Theorem 6.4.6), we see that $(\bar{x},\bar{t})$ has a canonical
neighborhood with parameter $\varepsilon$, which is a
contradiction. Therefore we have completed the proof of the
proposition.
\endproof

Summing up, we have proved that for any $\varepsilon>0$, (without
loss of generality, we may assume $\varepsilon \leq
\bar{\varepsilon}_0$), there exist nonincreasing (continuous)
positive functions $\widetilde{\delta}(t)$ and $\widetilde{r}(t)$,
defined on $[0,+\infty)$ with
$$
\widetilde{\delta}(t) \leq \bar{\delta}(t) = \min
\left\{\frac{1}{{2e^2\log(1+t)}}, \delta_0\right\},
$$
such that for arbitrarily given (continuous) positive function
$\delta(t)$ with $\delta(t)<\widetilde{\delta}(t)$ on
$[0,+\infty)$, and arbitrarily given a compact orientable
normalized three-manifold as initial data, the Ricci flow with
surgery has a solution on $[0,T)$ obtained by evolving the Ricci
flow and by performing $\delta$-cutoff surgeries at a sequence of
times $0<t_1 < t_2 < \cdots < t_i < \cdots<T$, with ${\delta}(t_i)
\leq \delta \leq \widetilde{\delta}(t_i)$ at each time $t_i$, so
that the pinching assumption and the canonical neighborhood
assumption (with accuracy $\varepsilon$) with $r=\widetilde{r}(t)$
are satisfied.  (At this moment we still do not know whether the
surgery times $t_i$ are discrete.)

Since the $\delta$-cutoff surgeries occur at the points lying
deeply in the $\varepsilon$-horns, the minimum of the scalar
curvature $R_{min}(t)$ of the solution to the Ricci flow with
surgery at each time-slice is achieved in the region unaffected by
the surgeries. Thus we know from the evolution equation of the
scalar curvature that \be
\frac{d}{dt}R_{\min}(t) \geq \frac{2}{3}R^2_{\min}(t).
\ee In particular, the minimum of the scalar curvature
$R_{\min}(t)$ is nondecreasing in time. Also note that each
$\delta$-cutoff surgery decreases volume. Then the upper
derivative of the volume in time satisfies
\begin{align*}
\bar{\(\frac{d}{dt}\)}V(t) &\triangleq \lim\sup_{\triangle
t\rightarrow
0}\frac{V(t+\triangle t)-V(t)}{\triangle t}\\
&\leq-R_{\min}(0)V(t)
\end{align*}
which implies that
$$
V(t)\leq V(0)e^{-R_{\min}(0)t}.
$$

On the other hand, by Lemma 7.3.2 and the $\delta$-cutoff procedure
given in the previous section, we know that at each time $t_i$, each
$\delta$-cutoff surgery cuts down the volume at least at an amount
of $h^3(t_i)$ with $h(t_i)$ depending only on ${\delta}(t_i)$ and
$\widetilde{r}(t_i)$. Thus the surgery times $t_i$ cannot accumulate
in any finite interval. When the solution becomes extinct at some
finite time $T$, the solution at time near $T$ is entirely covered
by canonical neighborhoods and then the initial manifold is
diffeomorphic to a connected sum of a finite copies of
$\mathbb{S}^2\times \mathbb{S}^1$ and $\mathbb{S}^3/\Gamma$ (the
metric quotients of round three-sphere). So we have proved the
following long-time existence result which was given by Perelman in
\cite{P2}.

\begin{theorem}[Long-time existence theorem]
For any fixed constant $\varepsilon >0$, there exist nonincreasing
$($continuous$)$ positive functions $\widetilde{\delta}(t)$ and
$\widetilde{r}(t)$, defined on $[0,+\infty)$, such that for an
arbitrarily given $($continuous$)$ positive function $\delta(t)$
with $\delta(t) \leq \widetilde{\delta}(t)$ on $[0,+\infty)$, and
arbitrarily given a compact orientable normalized three-manifold
as initial data, the Ricci flow with surgery has a solution with
the following properties: either
\begin{itemize}
\item[(i)] it is defined on a finite interval $[0,T)$ and obtained
by evolving the Ricci flow and by performing a finite number of
cutoff surgeries, with each $\delta$-cutoff at a time $t \in
(0,T)$ having ${\delta}= \delta(t)$, so that the solution becomes
extinct at the finite time $T$, and the initial manifold is
diffeomorphic to a connected sum of a finite copies of
$\mathbb{S}^2\times \mathbb{S}^1$ and $\mathbb{S}^3/\Gamma$ $($the
metric quotients of round three-sphere$)$ ; or \item[(ii)] it is
defined on $[0,+\infty)$ and obtained by evolving the Ricci flow
and by performing at most countably many cutoff surgeries, with
each $\delta$-cutoff at a time $t \in [0,+\infty)$ having
${\delta}=\delta(t)$, so that the pinching assumption and the
canonical neighborhood assumption $($with accuracy $\varepsilon)$
with $r=\widetilde{r}(t)$ are satisfied, and there exist at most a
finite number of surgeries on every finite time interval.
\end{itemize}
\end{theorem}

In particular, if the initial manifold has positive scalar
curvature, say $R\geq a>0$, then by (7.4.20), the solution becomes
extinct at $T\leq \frac{3}{2}a$. Hence we have the following
topological description of compact three-manifolds with
nonnegative scalar curvature which improves the well-known work of
Schoen-Yau \cite{ScY79m}, \cite{ScY79}.

\begin{corollary}[{Perelman \cite{P2}}]
 Let $M$ be a compact orientable three-manifold with nonnegative
scalar curvature. Then either $M$ is flat or it is diffeomorphic
to a connected sum of a finite copies of $\mathbb{S}^2\times
\mathbb{S}^1$ and $\mathbb{S}^3/\Gamma$ $($the metric quotients of
the round three-sphere$).$
\end{corollary}

The famous \textbf{Poincar\'{e} conjecture}\index{Poincar\'{e}
conjecture} states that every compact three-manifold with trivial
fundamental group is diffeomorphic to $\mathbb{S}^3$. Developing
tools to attack the conjecture formed the basis for much of the
works in three-dimensional topology over the last one hundred
years. Now we use the Ricci flow to discuss the Poincar$\acute{e}$
conjecture.

Let $M$ be a compact three-manifold with trivial fundamental
group. In particular, the three-manifold $M$ is orientable.
Arbitrarily given a Riemannian metric on $M$, by scaling we may
assume the metric is normalized. With this normalized metric as
initial data, we consider the solution to the Ricci flow with
surgery. If one can show the solution becomes extinct in finite
time, it will follow from Theorem 7.4.3 (i) that the
three-manifold $M$ is diffeomorphic to the three-sphere
$\mathbb{S}^3$. Such a finite extinction time result was first
proposed by Perelman in \cite{P3}. Recently, Colding-Minicozzi has
published a proof of it in \cite{CM}. So \textbf{the combination
of Theorem 7.4.3 (i) and the finite extinction result \cite{P3,
CM} gives a proof of the Poincar\'{e} conjecture}.

We also remark that the above long-time existence result of
Perelman has been extended to compact four-manifolds with positive
isotropic curvature by Chen and the second author in \cite{CZ05F}.
As a consequence it gave a complete proof of the following
classification theorem of compact four-manifolds, with no
essential incompressible space-form and with a metric of positive
isotropic curvature. The theorem was first proved by Hamilton in
(\cite{Ha97}), though it was later found that the proof contains
some gaps (see for example the comment of Perelman in Page 1, the
second paragraph, of \cite{P2}).

\begin{theorem}
A compact four-manifold with no essential incompressible
space-form and with a metric of positive isotropic curvature is
diffeomorphic to $\mathbb{S}^4$, or $\mathbb{RP}^4$, or
$\mathbb{S}^3\times \mathbb{S}^1$, or
$\mathbb{S}^3\widetilde{\times}\mathbb{S}^1$ $($the $\mathbb{Z}_2$
quotient of $\mathbb{S}^{3}\times \mathbb{S}^{1}$ where
$\mathbb{Z}_{2}$ flips $\mathbb{S}^{3}$ antipodally and rotates
$\mathbb{S}^{1}$ by $180^{0}),$ or a connected sum of them.
\end{theorem}

\section{Curvature Estimates for Surgically Modified Solutions}

This section is a detailed exposition of section 6 of Perelman
\cite{P2}. Here we will generalize the curvature estimates for
smooth solutions in Section 7.2 to that of solutions with cutoff
surgeries. We first state and prove a version of Theorem 7.2.1.

\begin{theorem}[{Perelman \cite{P2}}]
For any $\varepsilon > 0$ and  $1 \leq A<+\infty$, one can find
$\kappa=\kappa(A,\varepsilon)>0$,
$K_1=K_1(A,\varepsilon)<+\infty$, $K_2=K_2(A,\varepsilon)<+\infty$
and $\bar{r}=\bar{r}(A,\varepsilon)>0$ such that for any
$t_0<+\infty$ there exists $\bar{\delta}_A=\bar{\delta}_A(t_0)>0$
$($depending also on $\varepsilon),$ nonincreasing in $t_0$, with
the following property. Suppose we have a solution, constructed by
Theorem $7.4.3$ with the nonincreasing (continuous) positive
functions $\widetilde{\delta}(t)$ and $\widetilde{r}(t)$, to the
Ricci flow with $\delta$-cutoff surgeries on time interval $[0,T]$
and with a compact orientable normalized three-manifold as initial
data, where each $\delta$-cutoff at a time $t$ satisfies
$\delta=\delta(t) \leq \widetilde{\delta}(t)$ on $[0,T]$ and
$\delta=\delta(t) \leq \bar{\delta}_A$ on $[\frac{t_0}{2},t_0]$;
assume that the solution is defined on the whole parabolic
neighborhood $P(x_0,t_0,r_0,-r^2_0)\triangleq\{(x,t)\ |\ x\in
B_t(x_0,r_0),t\in[t_0-r^2_0,t_0]\}$, $2r^2_0<t_0$, and satisfies
$$
|Rm|\leq r^{-2}_0 \ \ on \ \ P(x_0,t_0,r_0,-r^2_0),
$$
$$
\text{and} \qquad \Vol_{t_0}(B_{t_0}(x_0,r_0))\geq A^{-1}r^3_0.
$$
Then
\begin{itemize}
\item[(i)] the solution is $\kappa$-noncollapsed on all scales
less than $r_0$ in the ball $B_{t_0}(x_0,Ar_0)$; \item[(ii)] every
point $x\in B_{t_0}(x_0,Ar_0)$ with $R(x,t_0)\geq K_1r^{-2}_0$ has
a canonical neighborhood $B$, with $ B_{t_0}(x,\sigma)\subset
B\subset B_{t_0}(x,2\sigma)$ for some
$0<\sigma<C_1(\varepsilon)R^{-\frac{1}{2}}(x,t_0), $ which is
either a strong $\varepsilon$-neck or an $\varepsilon$-cap;
\item[(iii)] if $r_0\leq \bar{r}\sqrt{t_0}$ then $R\leq
K_2r^{-2}_0$ in $B_{t_0}(x_0,Ar_0)$.
\end{itemize}
Here $C_1(\varepsilon)$ is the positive constant in the canonical
neighborhood assumption.

\end{theorem}

\begin{pf}
Without loss of generality, we may assume $0<\varepsilon \leq
\bar{\varepsilon}_0$, where $\bar{\varepsilon}_0$ is the
sufficiently small (universal) positive constant in Lemma 7.4.2.

\medskip
(i) This is analog of no local collapsing theorem II (Theorem
3.4.2). In comparison with the no local collapsing theorem II,
this statement gives $\kappa$-noncollapsing property no matter how
big the time is and it also allows the solution to be modified by
surgery.

Let $\eta (\geq 10)$ be the universal constant in the definition
of the canonical neighborhood assumption. Recall that we had
removed every component which has positive sectional curvature in
our surgery procedure. By the same argument as in the first part
of the proof of Lemma 7.4.2, the canonical neighborhood assumption
of the solution implies the $\kappa$-noncollapsing on the scales
less than $\frac{1}{2\eta}\widetilde{r}(t_0)$ for some positive
constant $\kappa$ depending only on $C_1(\varepsilon)$ and
$C_2(\varepsilon)$ (in the definition of the canonical
neighborhood assumption). So we may assume
$\frac{1}{2\eta}\widetilde{r}(t_0)\leq r_0\leq
\sqrt{\frac{t_0}{2}}$, and study the scales $\rho$,
$\frac{1}{2\eta}\widetilde{r}(t_0)\leq \rho\leq r_0$. Let $x\in
B_{t_0}(x_0,Ar_0)$ and assume that the solution satisfies
$$
|Rm|\leq \rho^{-2}
$$
for those points in $P(x,t_0,\rho,-\rho^{2})\triangleq\{(y,t)\ |\
y\in B_t(x,\rho),t\in[t_0-\rho^2,t_0]\}$ for which the solution is
defined. We want to bound the ratio
$\Vol_{t_0}(B_{t_0}(x,\rho))/\rho^3$ from below.

Recall that a space-time curve is called admissible if it stays in
the region unaffected by surgery, and a space-time curve on the
boundary of the set of admissible curves is called a barely
admissible curve. Consider any barely admissible curve $\gamma$,
parametrized by $t\in[t_{\gamma},t_0]$, $t_0-r^2_0\leq
t_{\gamma}\leq t_0$, with $\gamma(t_0)=x$. The same proof for the
assertion (7.4.2) (in the proof of Lemma 7.4.2) shows that for
arbitrarily large $L>0$ (to be determined later), one can find a
sufficiently small
$\bar{\delta}(L,t_0,\widetilde{r}(t_0),\widetilde{r}
(\frac{t_0}{2}),\varepsilon)>0$ such that when each
$\delta$-cutoff in $[\frac{t_0}{2},t_0]$ satisfies $\delta \leq
\bar{\delta}(L,t_0,\widetilde{r}(t_0),\widetilde{r}
(\frac{t_0}{2}),\varepsilon)$, there holds \be
\int^{t_0}_{t_{\gamma}}\sqrt{t_0-t}(R_+(\gamma(t),t)
+|\dot{\gamma}(t)|^2)dt\geq Lr_0. 
\ee

{}From now on, we assume that each $\delta$-cutoff of the solution
in the time interval $[\frac{t_0}{2},t_0]$ satisfies $\delta \leq
\bar{\delta}(L,t_0,\widetilde{r}(t_0),
\widetilde{r}(\frac{t_0}{2}),\varepsilon).$

Let us scale the solution, still denoted by $g_{ij}(\cdot,t)$, to
make $r_0=1$ and the time as $t_0=1$. By the maximum principle, it
is easy to see that the (rescaled) scalar curvature satisfies
$$
R \geq -\frac{3}{2t}
$$
on $(0,1]$. Let us consider the time interval $[\frac{1}{2},1]$
and define a function of the form
$$
h(y,t)=\phi(d_t(x_0,y)-A(2t-1))(\bar{L}(y,\tau)+2\sqrt{\tau})
$$
where $\tau=1-t$, $\phi$ is the function of one variable chosen in
the proof of Theorem 3.4.2 which is equal to one on
$(-\infty,\frac{1}{20})$, rapidly increasing to infinity on
$(\frac{1}{20},\frac{1}{10})$, and satisfies
$2\frac{(\phi')^2}{\phi}-\phi''\geq(2A+300)\phi'-C(A)\phi$ for
some constant $C(A)<+\infty$, and $\bar{L}$ is the function
defined by
\begin{align*}
\bar{L}(q,\tau) &=\inf\Big\{2\sqrt{\tau}\int^{\tau}_0\sqrt{s}
(R+|\dot{\gamma}|^2)ds\ |\ (\gamma(s),s),
s\in[0,\tau]\\
&\text{is a space-time curve with $\gamma(0)=x$ and
$\gamma(\tau)=q$}\Big\}.
\end{align*}
Note that
\begin{align}
\bar{L}(y,\tau)
& \geq 2\sqrt{\tau}\int^{\tau}_0\sqrt{s}Rds\\
& \geq -4\tau^2 \nn\\
&  >   -2\sqrt{\tau} \nn
\end{align}
since $R\geq -3$ and $0<\tau\leq \frac{1}{2}$. This says $h$ is
positive for $t\in[\frac{1}{2},1]$. Also note that \be
\frac{\partial}{\partial \tau}\bar{L}+\triangle \bar{L}\leq 6
\ee as long as $\bar{L}$ is achieved by admissible curves. Then as
long as the shortest $\mathcal{L}$-geodesics from $(x_0,0)$ to
$(y,\tau)$ are admissible, there holds at $y$ and $t=1-\tau$,
\begin{align*}
\(\frac{\partial}{\partial t}-\triangle\)h &\geq
\(\phi'\left[\(\frac{\partial}{\partial t}
-\triangle\)d_t-2A\right]-\phi''\)\cdot(\bar{L}+2\sqrt{\tau})\\
&\quad-\(6+\frac{1}{\sqrt{\tau}}\)\phi-2\<\nabla
\phi,\nabla\bar{L}\>.
\end{align*}

Firstly,\; we\; may\; assume\; the\; constant\; $L$\; in (7.5.1)\;
is\; not\; less\; than $2\exp(C(A)+100)$. We claim that Lemma
3.4.1(i) is applicable for $d=d_t(\cdot,x_0)$ at $y$ and $
t=1-\tau $ (with $\tau\in[0,\frac{1}{2}]$) whenever
$\bar{L}(y,\tau)$ is achieved by admissible curves and satisfies
the estimate
$$
\bar{L}(y,\tau) \leq 3\sqrt{\tau}\exp(C(A)+100).
$$
Indeed, since the solution is defined on the whole neighborhood
$P(x_0,t_0,r_0,$ $-r^2_0)$ with $r_0=1$ and $t_0=1$, the point
$x_0$ at the time $t=1-\tau$ lies on the region unaffected by
surgery. Note that $R \geq -3$ for $t \in [\frac{1}{2},1]$.  When
$\bar{L}(y,\tau)$ is achieved by admissible curves and satisfies
$\bar{L}(y,\tau) \leq 3\sqrt{\tau}\exp(C(A)+100)$, the estimate
(7.5.1) implies that the point $y$ at the time $t=1-\tau$ does not
lie in the collars of the gluing caps. Thus any minimal geodesic
(with respect to the metric $g_{ij}(\cdot,t)$ with $t=1-\tau$)
connecting $x_0$ and $y$ also lies in the region unaffected by
surgery; otherwise the geodesic is not minimal.  Then from the
proof of Lemma 3.4.1(i), we see that it is applicable.

Assuming the minimum of $h$ at a time, say $t=1-\tau$, is achieved
at a point, say $y$, and assuming $\bar{L}(y,\tau)$ is achieved by
admissible curves and satisfies $\bar{L}(y,\tau) \leq
3\sqrt{\tau}\exp(C(A)+100)$, we have
$$
(\bar{L}+2\sqrt{\tau})\nabla\phi=-\phi\nabla\bar{L},
$$
and then by the computations and estimates in the proof of Theorem
3.4.2,
\begin{align*}
&\(\frac{\partial}{\partial t}-\triangle\)h  \\
&  \geq \(\phi'\left[\(\frac{\partial}{\partial
t}-\triangle\)d_t-2A\right]-\phi''
+2\frac{(\phi')^2}{\phi}\)\cdot(\bar{L}+2\sqrt{\tau})
-\(6+\frac{1}{\sqrt{\tau}}\)\phi\\
& \geq -C(A)h-\(6+\frac{1}{\sqrt{\tau}}\)
\frac{h}{(2\sqrt{\tau}-4\tau^2)},
\end{align*}
at $y$ and $t=1-\tau$. Here we used (7.5.2) and Lemma 3.4.1(i).

As before, denoting by $h_{min}(\tau)=\min_{z} h(z,1-\tau)$, we
obtain
\begin{align}
\frac{d}{d\tau}\(\log\(\frac{h_{\min}(\tau)}{\sqrt{\tau}}\)\) &
\leq C(A)+\frac{6\sqrt{\tau}+1}{2\tau-4\tau^2\sqrt{\tau}}
-\frac{1}{2\tau}\\
& \leq C(A)+\frac{50}{\sqrt{\tau}}, \nn
\end{align} 
as long as the associated shortest $\mathcal{L}$-geodesics are
admissible with $\bar{L} \leq 3\sqrt{\tau}\exp(C(A)+100)$. On the
other hand, by definition, we have \be \lim_{\tau\rightarrow
0^+}\frac{h_{\min}(\tau)}{\sqrt{\tau}}\leq
\phi(d_1(x_0,x)-A)\cdot 2=2. 
\ee The combination of (7.5.4) and (7.5.5) gives the following
assertion:

\vskip 0.2cm \emph{Let $\tau\in[0,\frac{1}{2}]$. If for each
$s\in[0,\tau]$, $\inf\{\bar{L}(y,s)\ |\ \mbox{ } d_t(x_0,y)\leq
A(2t-1)+\frac{1}{10}\mbox{ with }s=1-t\}$ is achieved by
admissible curves, then we have
\begin{align}
&\inf\left\{\bar{L}(y,\tau)\ |\ \mbox{ } d_t(x_0,y)
\leq A(2t-1)+\frac{1}{10}\mbox{ with }\tau=1-t\right\}\\
&\leq 2\sqrt{\tau}\exp(C(A)+100). \nn  
\end{align}
}

\vskip 0.1cm Note again that $R \geq -3$ for $t \in
[\frac{1}{2},1]$. By combining with (7.5.1), we know that any
barely admissible curve $\gamma$, parametrized by $s\in[0,\tau]$,
$0\leq \tau \leq \frac{1}{2}$, with $\gamma(0)=x$, satisfies
$$
\int^{\tau}_0\sqrt{s}(R + |\dot{\gamma}|^2)ds \geq
\frac{7}{4}\exp(C(A)+100),
$$
by assuming $L\geq 2\exp(C(A)+100)$.

Since $|Rm|\leq \rho^{-2}$ on $P(x,t_0,\rho,-\rho^{2})$ with
$\rho\geq \frac{1}{2\eta}\widetilde{r}(t_0)$ (and $t_0=1$) and
$\bar{\delta}(L,t_0,\widetilde{r}(t_0),
\widetilde{r}(\frac{t_0}{2}),\varepsilon)>0$ is sufficiently
small, the parabolic neighborhood $P(x,1,\rho,-\rho^2)$ around the
point $(x,1)$ is contained in the region unaffected by the
surgery. Thus as $\tau=1-t$ is sufficiently close to zero,
$\frac{1}{2\sqrt{\tau}}\inf\bar{L}$ can be bounded from above by a
small positive constant and then the infimum
$\inf\{\bar{L}(y,\tau)\ |\ \mbox{ } d_t(x_0,y)\leq
A(2t-1)+\frac{1}{10}\mbox{ with }\tau=1-t\}$ is achieved by
admissible curves.

Hence we conclude that for each $\tau\in [0,\frac{1}{2}]$, any
minimizing curve $\gamma_{\tau}$ of $\inf\{\bar{L}(y,\tau)|\mbox{
} d_t(x_0,y)\leq A(2t-1)+\frac{1}{10}\mbox{ with }\tau=1-t\}$ is
admissible and satisfies
$$
\int^{\tau}_0\sqrt{s}(R+|\dot{\gamma}_{\tau}|^2)ds\leq
\exp(C(A)+100).
$$

Now we come back to the unrescaled solution. It then follows that
the Li-Yau-Perelman distance $l$ from $(x,t_0)$ satisfies the
following estimate \be \min\left\{l\(y,t_0-\frac{1}{2}r^2_0\)\
\Big|\ y\in
B_{t_0-\frac{1}{2}r^2_0}\(x_0,\frac{1}{10}r_0\)\right\}
\leq\exp(C(A)+100), 
\ee by noting the (parabolic) scaling invariance of the
Li-Yau-Perelman distance.

By the assumption that $|Rm|\leq r^{-2}_0$ on
$P(x_0,t_0,r_0,-r^2_0)$, exactly as before, for any $q\in
B_{t_0-r^2_0}(x_0,r_0)$, we can choose a path $\gamma$
parametrized by $\tau\in[0,r^2_0]$ with $\gamma(0)=x$,
$\gamma(r^2_0)=q$, and $\gamma(\frac{1}{2}r^2_0)=y\in
B_{t_0-\frac{1}{2}r^2_0}(x_0,\frac{1}{10}r_0)$, where
$\gamma|_{[0,\frac{1}{2}r^2_0]}$ achieves the minimum
$\min\{l(y,t_0-\frac{1}{2}r^2_0)\ |\ y\in
B_{t_0-\frac{1}{2}r^2_0}(x_0,\frac{1}{10}r_0)\}$ and
$\gamma|_{[\frac{1}{2}r^2_0,r^2_0]}$ is a suitable curve
satisfying $\gamma|_{[\frac{1}{2}r^2_0,r^2_0]}(\tau) \in B_{t_0 -
\tau}(x_0,r_0)$, for each $\tau \in [\frac{1}{2}r^2_0,r^2_0]$, so
that the $\mathcal{L}$-length of $\gamma$ is uniformly bounded
from above by a positive constant (depending only on $A$)
multiplying $r_0$. This implies that the Li-Yau-Perelman distance
from $(x,t_0)$ to the ball $B_{t_0-r^2_0}(x_0,r_0)$ is uniformly
bounded by a positive constant $L(A)$ (depending only on $A$). Now
we can choose the constant $L$ in (7.5.1) by
$$
L=\max\{2L(A),2\exp(C(A)+100)\}.
$$
Thus every shortest $\mathcal{L}$-geodesic from $(x,t_0)$ to the
ball $B_{t_0-r^2_0}(x_0,r_0)$ is necessarily admissible. By
combining with the assumption that $\Vol_{t_0}(B_{t_0}(x_0,r_0))$
$\geq\; \;A^{-1}\,r^3_0$,\; we\; conclude\; that\; Perelman's\;
reduced\; volume\; of\; the\; ball $B_{t_0-r^2_0}(x_0,r_0)$
satisfies the estimate
\begin{align}
\widetilde{V}_{r^2_0}(B_{t_0-r^2_0}(x_0,r_0)) &  =
\int_{B_{t_0-r^2_0}(x_0,r_0)}
(4\pi r^2_0)^{-\frac{3}{2}}\exp(-l(q,r^2_0))dV_{t_0-r^2_0}(q)\\
&  \geq c(A) \nn
\end{align}  
for some positive constant $c(A)$ depending only on $A$.

We can now argue as in the last part of the proof of Lemma 7.4.2
to get a lower bound estimate for the volume of the ball
$B_{t_0}(x,\rho)$. The union of all shortest
$\mathcal{L}$-geodesics from $(x,t_0)$ to the ball
$B_{t_0-r^2_0}(x_0,r_0)$, defined by
\begin{multline*}
CB_{t_0-r^2_0}(x_0,r_0)=\{(y,t)\ |\ (y,t)
\mbox{ lies in a shortest $\mathcal{L}$-geodesic from}\\
(x,t_0)\mbox{ to a point in }B_{t_0-r^2_0}(x_0,r_0)\},
\end{multline*}
forms a cone-like subset in space-time with vertex $(x,t_0)$.
Denote by $B(t)$ the intersection of the cone-like subset
$CB_{t_0-r^2_0}(x_0,r_0)$ with the time-slice at $t$. Perelman's
reduced volume of the subset $B(t)$ is given by
$$
\widetilde{V}_{t_0-t}(B(t))=\int_{B(t)}(4\pi(t_0-t))^{-\frac{3}{2}}
\exp(-l(q,t_0-t))dV_t(q).
$$
Since the cone-like subset $CB_{t_0-r^2_0}(x_0,r_0)$ lies entirely
in the region unaffected by surgery, we can apply Perelman's
Jacobian comparison Theorem 3.2.7 and the estimate (7.5.8) to
conclude that
\begin{align}
\widetilde{V}_{t_0-t}(B(t))&  \geq
\widetilde{V}_{r^2_0}(B_{t_0-r^2_0}(x_0,r_0))\\
&  \geq c(A) \nn
\end{align}         
for all $t\in[t_0-r^2_0,t_0]$.

As before, denoting by
$\xi=\rho^{-1}\Vol_{t_0}(B_{t_0}(x,\rho))^{\frac{1}{3}}$, we only
need to get a positive lower bound for $\xi$. Of course we may
assume $\xi<1$. Consider $\widetilde{B}(t_0-\xi\rho^2)$, the
subset at the time-slice $\{t=t_0-\xi\rho^2\}$ where every point
can be connected to $(x,t_0)$ by an admissible shortest
$\mathcal{L}$-geodesic. Perelman's reduced volume of
$\widetilde{B}(t_0-\xi\rho^2)$ is given by
\begin{align}
& \widetilde{V}_{\xi\rho^2}(\widetilde{B}(t_0-\xi\rho^2))\\
&  = \int_{\widetilde{B}(t_0-\xi\rho^2)}
(4\pi\xi\rho^2)^{-\frac{3}{2}}
\exp(-l(q,\xi\rho^2))dV_{t_0-\xi\rho^2}(q) \nn\\
& =\int_{\widetilde{B}(t_0-\xi\rho^2)\cap\mathcal{L}
\exp_{\{|\upsilon|\leq\frac{1}{4}
\xi^{-\frac{1}{2}}\}}(\xi\rho^2)} (4\pi\xi\rho^2)^{-\frac{3}{2}}
\exp(-l(q,\xi\rho^2))dV_{t_0-\xi\rho^2}(q) \nn\\
&\quad  +\int_{\widetilde{B}(t_0-\xi\rho^2)
\setminus\mathcal{L}\exp_{\{|\upsilon|\leq\frac{1}{4}
\xi^{-\frac{1}{2}}\}}(\xi\rho^2)} (4\pi\xi\rho^2)^{-\frac{3}{2}}
\exp(-l(q,\xi\rho^2))dV_{t_0-\xi\rho^2}(q). \nn
\end{align}  
Note that the whole region $P(x,t_0,\rho,-\rho^2)$ is unaffected
by surgery because $\rho\geq \frac{1}{2\eta}\widetilde{r}(t_0)$
and $\bar{\delta}(L,t_0,\widetilde{r}(t_0),
\widetilde{r}(\frac{t_0}{2}),\varepsilon)>0$ is sufficiently
small. Then exactly as before, there is a universal positive
constant $\xi_0$ such that when $0<\xi\leq\xi_0$, there holds
$$
\mathcal{L}\exp_{\{|\upsilon|\leq\frac{1}{4}
\xi^{-\frac{1}{2}}\}}(\xi\rho^2)\subset B_{t_0}(x,\rho)
$$
and the first term on RHS of (7.5.10) can be estimated by
\begin{align}
&\int_{\widetilde{B}(t_0-\xi\rho^2)\cap\mathcal{L}
\exp_{\{|\upsilon|\leq\frac{1}{4}\xi^{-\frac{1}{2}}\}}
(\xi\rho^2)}(4\pi\xi\rho^2)^{-\frac{3}{2}}
\exp(-l(q,\xi\rho^2))dV_{t_0-\xi\rho^2}(q)\\
&\leq e^{C\xi}(4\pi)^{-\frac{3}{2}}\xi^{\frac{3}{2}}\nn  
\end{align}
for some universal constant $C$; while the second term on RHS of
(7.5.10) can be estimated by
\begin{align}
&\int_{\widetilde{B}(t_0-\xi\rho^2)\setminus\mathcal{L}
\exp_{\{|\upsilon|\leq\frac{1}{4}\xi^{-\frac{1}{2}}\}}
(\xi\rho^2)}(4\pi\xi\rho^2)^{-\frac{3}{2}}
\exp(-l(q,\xi\rho^2))dV_{t_0-\xi\rho^2}(q)\\
&\leq(4\pi)^{-\frac{3}{2}}\int_{\{|\upsilon|>\frac{1}{4}
\xi^{-\frac{1}{2}}\}}\exp(-|\upsilon|^2)d\upsilon. \nn
\end{align}
Since $B(t_0-\xi\rho^2)\subset\widetilde{B}(t_0-\xi\rho^2)$, the
combination of (7.5.9)-(7.5.12) bounds $\xi$ from below by a
positive constant depending only on $A$. This proves the statement
(i).

\medskip
(ii) This is analogous to the claim in the proof of Theorem 7.2.1.
We argue by contradiction. Suppose that for some $A<+\infty$ and a
sequence $K^{\alpha}_1\rightarrow\infty$, there exists a sequence
$t^{\alpha}_0$ such that for any sequences
$\bar{\delta}^{\alpha\beta}>0$ with
$\bar{\delta}^{\alpha\beta}\rightarrow0$ for fixed $\alpha$, we
have sequences of solutions $g^{\alpha\beta}_{ij}$ to the Ricci
flow with surgery and sequences of points $x^{\alpha\beta}_0$, of
radii $r^{\alpha\beta}_0$, which satisfy the assumptions but
violate the statement (ii) at some $x^{\alpha\beta}\in
B_{t^{\alpha}_0}(x^{\alpha\beta}_0,Ar^{\alpha\beta}_0)$ with
$R(x^{\alpha\beta},t^{\alpha}_0)\geq
K^{\alpha}_1(r^{\alpha\beta}_0)^{-2}$. Slightly abusing notation,
we will often drop the indices $\alpha,\beta$ in the following
argument.

Exactly as in the proof of Theorem 7.2.1, we need to adjust the
point $(x,t_0)$. More precisely, we claim that there exists a
point $(\bar{x},\bar{t})\!\in\! B_{\bar{t}}(x_0,2Ar_0)$
$\times[t_0-\frac{r^2_0}{2},t_0]$ with $\bar{Q}\triangleq
R(\bar{x},\bar{t})\geq K_1r^{-2}_0$ such that the point
$(\bar{x},\bar{t})$ does not satisfy the canonical neighborhood
statement, but each point $(y,t)\in \bar{P}$ with $R(y,t)\geq
4\bar{Q}$ does, where $\bar{P}$ is the set of all $(x',t')$
satisfying $\bar{t}-\frac{1}{4}K_1\bar{Q}^{-1}\leq t'\leq
\bar{t}$, $d_{t'}(x_0,x')\leq
d_{\bar{t}}(x_0,\bar{x})+K^{\frac{1}{2}}_1\bar{Q}^{-\frac{1}{2}}$.
Indeed as before, the point $(\bar{x},\bar{t})$ is chosen by an
induction argument. We first choose $(x_1,t_1)=(x,t_0)$ which
satisfies $d_{t_1}(x_0,x_1)\leq Ar_0$ and $R(x_1,t_1)\geq
K_1r^{-2}_0$, but does not satisfy the canonical neighborhood
statement. Now if $(x_k,t_k)$ is already chosen and is not the
desired $(\bar{x},\bar{t})$, then some point $(x_{k+1},t_{k+1})$
satisfies $t_k-\frac{1}{4}K_1R(x_k,t_k)^{-1}\leq t_{k+1}\leq t_k$,
$d_{t_{k+1}}(x_0,x_{k+1})\leq
d_{t_k}(x_0,x_k)+K^{\frac{1}{2}}_1R(x_k,t_k)^{-\frac{1}{2}}$, and
$R(x_{k+1},t_{k+1})\geq 4R(x_k,t_k)$, but $(x_{k+1},t_{k+1})$ does
not satisfy the canonical neighborhood statement. Then we have
\begin{align*}
R(x_{k+1},t_{k+1})&\geq 4^kR(x_1,t_1)\geq 4^kK_1r^{-2}_0,\\
d_{t_{k+1}}(x_0,x_{k+1}) &\leq d_{t_1}(x_0,x_1)+K^{\frac{1}{2}}_1
\sum^k_{i=1}R(x_i,t_i)^{-\frac{1}{2}} \leq Ar_0+2r_0,
\end{align*}
and
$$
t_0\geq t_{k+1}\geq
t_0-\frac{1}{4}K_1\sum^k_{i=1}R(x_i,t_i)^{-1}\geq
t_0-\frac{1}{2}r^2_0.
$$
So the sequence must be finite and its last element is the desired
$(\bar{x},\bar{t})$.

Rescale the solutions along $(\bar{x},\bar{t})$ with factor
$R(\bar{x},\bar{t})(\geq K_1r^{-2}_0)$ and shift the times
$\bar{t}$ to zero. We will adapt both the proof of Proposition
7.4.1 and that of Theorem 7.2.1 to show that a sequence of the
rescaled solutions $\widetilde{g}_{ij}^{\alpha\beta}$ converges to
an ancient $\kappa$-solution, which will give the desired
contradiction. Since we only need to consider the scale of the
curvature less than $\widetilde{r}({\bar{t}})^{-2}$, the present
situation is much easier than that of Proposition 7.4.1.

Firstly as before, we need to get a local curvature estimate.

For each adjusted $(\bar x,\bar t)$, let $[t',\bar t]$ be the
maximal subinterval of $[\bar
t-\frac{1}{20}\eta^{-1}\bar{Q}^{-1},\bar t]$\; so\; that\; for\;
each\; sufficiently\; large $\alpha$ and then sufficiently large
$\beta$, the canonical neighborhood statement holds for any
$(y,t)$ in $P(\bar x,\bar t,$ $\frac{1}{10}K_1^{\frac{1}{2}}
\bar{Q}^{-\frac{1}{2}},t'-\bar t)=\{ ({x},{t}) \ |\ {x} \in
B_{{t}}(\bar
x,\frac{1}{10}K_1^{\frac{1}{2}}\bar{Q}^{-\frac{1}{2}}), {t} \in
[t',\bar{t}]\}$ with $R(y,t)\geq 4\bar{Q}$, where $\eta$ is the
universal positive constant in the definition of canonical
neighborhood assumption. We want to show \be
t' = \bar t-\frac{1}{20}\eta^{-1}\bar{Q}^{-1}. 
\ee

Consider the scalar curvature $R$ at the point $\bar x$ over the
time interval $[t',\bar t]$. If there is a time $\tilde{t} \in
[t',\bar t]$ satisfying $R(\bar x,\tilde{t}) \geq 4\bar Q$, we let
$\tilde{t}$ be the first of such time from $\bar t$. Since the
chosen point $(\bar{x},\bar{t})$ does not satisfy the canonical
neighborhood statement, we know $R(\bar{x},\bar{t})\leq
\widetilde{r}(\bar{t})^{-2}$. Recall from our designed surgery
procedure that if there is a cutoff surgery at a point $x$ at a
time $t$, the scalar curvature at $(x,t)$ is at least
$(\bar{\delta}^{\alpha\beta})^{-2}\widetilde{r}({t})^{-2}$. Then
for each fixed $\alpha$, for $\beta$ large enough, the solution
$g^{\alpha\beta}_{ij}(\cdot,t)$ around the point $\bar x$ over the
time interval $[\tilde{t}-
\frac{1}{20}\eta^{-1}\bar{Q}^{-1},{\tilde{t}}]$ is well defined
and satisfies the following curvature estimate
$$
R(\bar x,t) \leq 8\bar Q,
$$
for $ t \in [\tilde{t}- \frac{1}{20}\eta^{-1}\bar{Q}^{-1},\bar t]$
(or $t \in [t',\bar t]$ if there is no such time $\tilde{t}$). By
the assumption that $t_0>2r^2_0$, we have
\begin{align*}
\bar{t}R(\bar{x},\bar{t})&  \geq
\frac{t_0}{2}R(\bar{x},\bar{t})\\
&  \geq r^2_0(K_1r^{-2}_0)\\
&  = K_1\rightarrow +\infty.
\end{align*}
Thus by using the pinching assumption and the gradient estimates
in the canonical neighborhood assumption, we further have
$$
|Rm(x,t)| \leq 30\bar{Q},
$$
for all $x \in B_t(\bar
x,\frac{1}{10}\eta^{-1}\bar{Q}^{-\frac{1}{2}})$ and $t \in
[\tilde{t}- \frac{1}{20}\eta^{-1}\bar{Q}^{-1},\bar t]$ (or $t \in
[t',\bar t]$) and all sufficiently large $\alpha$ and $\beta$.
Observe that Lemma 3.4.1 (ii) is applicable for $d_t(x_0,\bar{x})$
with $t \in [\tilde{t}- \frac{1}{20}\eta^{-1}\bar{Q}^{-1}, \bar
t]$ (or $t \in [t',\bar t]$) since any minimal geodesic, with
respect to the metric $g_{ij}(\cdot,t)$, connecting $x_0$ and
$\bar{x}$ lies in the region unaffected by surgery; otherwise the
geodesic is not minimal. After having obtained the above curvature
estimate, we can argue as deriving (7.2.2) and (7.2.3) in the
proof of Theorem 7.2.1 to conclude that any point $(x,t)$, with $
\tilde{t}- \frac{1}{20}\eta^{-1}\bar{Q}^{-1} \leq t \leq \bar t$
(or $t \in [t',\bar t]$) and $d_t(x,\bar x) \leq
\frac{1}{10}K_1^{\frac{1}{2}}\bar{Q}^{-\frac{1}{2}}$, satisfies
$$
d_t(x,x_{0}) \leq d_{\bar t}(\bar x,x_{0}) +
\frac{1}{2}K_1^{\frac{1}{2}}\bar{Q}^{-\frac{1}{2}},
$$
for all sufficiently large $\alpha$ and $\beta$. Then by combining
with the choice of the points $(\bar x,\bar t)$, we prove $t'=
\bar t-\frac{1}{20}\eta^{-1}\bar{Q}^{-1}$ (i.e., the canonical
neighborhood statement holds for any point $(y,t)$ in the
parabolic neighborhood $P(\bar x,\bar
t,\frac{1}{10}K_1^{\frac{1}{2}}\bar{Q}^{-\frac{1}{2}}$,$
-\frac{1}{20}\eta^{-1}\bar{Q}^{-1})$ with $R(y,t) \geq 4\bar Q$)
for all sufficiently large $\alpha$ and then sufficiently large
$\beta$.

Now it follows from the gradient estimates in the canonical
neighborhood assumption that the scalar curvatures of the rescaled
solutions $\widetilde{g}^{\alpha\beta}_{ij}$ satisfy
$$
\widetilde{R}(x,t)\leq 40
$$
for those $(x,t)\in P(\bar{x},0,\frac{1}{10}\eta^{-1},
-\frac{1}{20}\eta^{-1})\triangleq\{(x',t')\ |\
x'\in\widetilde{B}_{t'}(\bar{x},\frac{1}{10}\eta^{-1}),
t'\in[-\frac{1}{20}\eta^{-1},0]\}$, for which the rescaled
solution is defined. (Here $\widetilde{B}_{t'}$ denotes the
geodesic ball in the rescaled solution at time $t'$). Note again
that $R(\bar{x},\bar{t})\leq \widetilde{r}(\bar{t})^{-2}$ and
recall from our designed surgery procedure that if there is a
cutoff surgery at a point $x$ at a time $t$, the scalar curvature
at $(x,t)$ is at least
$(\bar{\delta}^{\alpha\beta})^{-2}\widetilde{r}({t})^{-2}$. Then
for each fixed sufficiently large $\alpha$, for $\beta$ large
enough, the rescaled solution $\widetilde{g}^{\alpha\beta}_{ij}$
is defined on the whole parabolic neighborhood
$P(\bar{x},0,\frac{1}{10}\eta^{-1},-\frac{1}{20}\eta^{-1})$. More
generally, for arbitrarily fixed $0<\widetilde{K}<+\infty$, there
is a positive integer $\alpha_0$ so that for each $\alpha \geq
\alpha_0$ we can find $\beta_0 >0$ (depending on $\alpha$) such
that if $\beta\geq \beta_0$ and $(y,0)$ is a point on the rescaled
solution $\widetilde{g}^{\alpha\beta}_{ij}$ with
$\widetilde{R}(y,0)\leq\widetilde{K}$ and
$\widetilde{d}_0(y,\bar{x}) \leq \widetilde{K}$, we have estimate
\be
\widetilde{R}(x,t)\leq 40\widetilde{K} 
\ee for\; $(x,\,t)\;\in\; P(y,\,0,\,\frac{1}{10}\,\eta^{-1}\,
\widetilde{K}^{-\frac{1}{2}},\;
-\frac{1}{20}\,\eta^{-1}\widetilde{K}^{-1})\;\triangleq\;\{(x',\,t')\
|\ x'\;\in\; \widetilde{B}_{t'}(y,$ $\frac{1}{10}\eta^{-1}
\widetilde{K}^{-\frac{1}{2}}),
t'\in[-\frac{1}{20}\eta^{-1}\widetilde{K}^{-1},0]\}$. In
particular, the rescaled solution is defined on the whole
parabolic neighborhood
$P(y,0,\frac{1}{10}\eta^{-1}\widetilde{K}^{-\frac{1}{2}},
-\frac{1}{20}\eta^{-1}\widetilde{K}^{-1})$.

Next, we want to show the curvature of the rescaled solutions at
the new times zero (after shifting) stay uniformly bounded at
bounded distances from $\bar{x}$ for some subsequences of $\alpha$
and $\beta$. Let $\alpha_m,\beta_m\rightarrow+\infty$ be chosen so
that the estimate (7.5.14) holds with $\widetilde{K}=m$. For all
$\rho \geq 0$, set
$$
M(\rho)=\sup\{\widetilde{R}(x,0)\ |\ m \geq 1, d_0(x,\bar{x})\leq
\rho\mbox{ in the rescaled solutions }
\widetilde{g}^{\alpha_m\beta_m}_{ij}\}
$$
and
$$
\rho_0=\sup\{\rho\geq0\ |\ M(\rho)<+\infty\}.
$$
Clearly the estimate (7.5.14) yields $\rho_0>0$. As we consider
the unshifted time $\bar{t}$, by combining with the assumption
that $t_0>2r^2_0$, we have
\begin{align}
\bar{t}R(\bar{x},\bar{t})&  \geq
\frac{t_0}{2}R(\bar{x},\bar{t})\\
&  \geq r^2_0(K_1r^{-2}_0) \nn\\
&  = K_1\rightarrow +\infty. \nn
\end{align}  
It then follows from the pinching assumption that we only need to
show $\rho_0=+\infty$. As before, we argue by contradiction.
Suppose we have a sequence of points $y_m$ in the rescaled
solutions $\widetilde{g}^{\alpha_m\beta_m}_{ij}$ with
$\widetilde{d}_0(\bar{x},y_m)\rightarrow \rho_0<+\infty$ and
$\widetilde{R}(y_m,0)\rightarrow+\infty$. Denote by $\gamma_m$ a
minimizing geodesic segment from $\bar{x}$ to $y_m$ and denote by
$\widetilde{B}_0(\bar{x},\rho_0)$ the open geodesic balls centered
at $\bar{x}$ of radius $\rho_0$ of the rescaled solutions. By
applying the assertion in statement (i), we have uniform
$\kappa$-noncollapsing at the points $(\bar{x},\bar{t})$. By
combining with the local curvature estimate (7.5.14) and
Hamilton's compactness theorem, we can assume that, after passing
to a subsequence, the marked sequence
$(\widetilde{B}_0(\bar{x},\rho_0),
\widetilde{g}^{\alpha_m\beta_m}_{ij},\bar{x})$ converges in the
$C^{\infty}_{\rm loc}$ topology to a marked (noncomplete) manifold
($B_{\infty},\widetilde{g}^{\infty}_{ij},x_{\infty}$) and the
geodesic segments $\gamma_m$ converge to a geodesic segment
(missing an endpoint) $\gamma_{\infty}\subset B_{\infty}$
emanating from $x_{\infty}$. Moreover, by the pinching assumption
and the estimate (7.5.15), the limit has nonnegative sectional
curvature.

Then exactly as before, we consider the tubular neighborhood along
$\gamma_{\infty}$
$$
V=\bigcup_{q_0\in\gamma_{\infty}}B_{\infty}
(q_0,4\pi(\widetilde{R}_{\infty}(q_0))^{-\frac{1}{2}})
$$
and the completion $\bar{B}_{\infty}$ of
($B_{\infty},\widetilde{g}^{\infty}_{ij}$) with $y_{\infty}\in
\bar{B}_{\infty}$ the limit point of $\gamma_{\infty}$. As before,
by the choice of the points $(\bar{x},\bar{t})$, we know that the
limiting metric $\widetilde{g}^{\infty}_{ij}$ is cylindrical at
any point $q_0\in \gamma_{\infty}$ which is sufficiently close to
$y_{\infty}$. Then by the same reason as before the metric space
$\bar{V}=V\cup\{y_{\infty}\}$ has nonnegative curvature in
Alexandrov sense, and we have a three-dimensional nonflat tangent
cone $C_{y_{\infty}}\bar{V}$ at $y_{\infty}$. Pick $z\in
C_{y_{\infty}}\bar{V}$ with distance one from the vertex and it is
nonflat around $z$. By definition $B(z,\frac{1}{2}) (\subset
C_{y_{\infty}}\bar{V})$ is the Gromov-Hausdoff convergent limit of
the scalings of a sequence of balls $B_{\infty}(z_k,\sigma_k)
(\subset(V,\widetilde{g}^{\infty}_{ij}))$ with
$\sigma_k\rightarrow0$. Since the estimate (7.5.14) survives on
$(V,\widetilde{g}^{\infty}_{ij})$ for all $\widetilde{K}<
+\infty$, we know that this convergence is actually in the
$C^{\infty}_{\rm loc}$ topology and over some time interval. Since
the limit $B(z,\frac{1}{2}) (\subset C_{y_{\infty}}\bar{V})$ is a
piece of a nonnegatively curved nonflat metric cone, we get a
contradiction with Hamilton's strong maximum Principle (Theorem
2.2.1) as before. Hence we have proved that a subsequence of the
rescaled solution $\widetilde{g}^{\alpha_m\beta_m}_{ij}$ has
uniformly bounded curvatures at bounded distance from $\bar{x}$ at
the new times zero.

Further, by the uniform $\kappa$-noncollapsing at the points
$(\bar{x},\bar{t})$ and the estimate (7.5.14) again, we can take a
$C^{\infty}_{\rm loc}$ limit
($M_{\infty},\widetilde{g}^{\infty}_{ij},x_{\infty}$), defined on
a space-time subset which contains the time slice $\{t=0\}$ and is
relatively open in $M_{\infty}\times(-\infty,0]$, for the
subsequence of the rescaled solutions. The limit is a smooth
solution to the Ricci flow, and is complete at $t=0$, as well as
has nonnegative sectional curvature by the pinching assumption and
the estimate (7.5.15). Thus by repeating the same argument as in
the Step 4 of the proof Proposition 7.4.1, we conclude that the
curvature of the limit at $t=0$ is bounded.

Finally we try to extend the limit backwards in time to get an
ancient $\kappa$-solution. Since the curvature of the limit is
bounded at $t=0$, it follows from the estimate (7.5.14) that the
limit is a smooth solution to the Ricci flow defined at least on a
backward time interval $[-a,0]$ for some positive constant $a$.
Let $(t_{\infty},0]$ be the maximal time interval over which we
can extract a smooth limiting solution. It suffices to show
$t_{\infty}=-\infty$. If $t_{\infty}>-\infty$, there are only two
possibilities: either there exist surgeries in finite distance
around the time $t_{\infty}$ or the curvature of the limiting
solution becomes unbounded as $t\searrow t_{\infty}$.

Let $c>0$ be a positive constant much smaller than
$\frac{1}{100}\eta^{-1}$. Note again that the infimum of the
scalar curvature is nondecreasing in time. Then we can find some
point $y_{\infty} \in M_{\infty}$ and some time $t=t_{\infty} +
\theta$ with $0<\theta<\frac{c}{3}$ such that
$\widetilde{R}_{\infty}(y_{\infty},t_{\infty}+\theta) \leq 2$.

Consider the (unrescaled) scalar curvature $R$ of
${g}_{ij}^{\alpha_m \beta_m}(\cdot,t)$ at the point $\bar x$ over
the time interval $[\bar
t+(t_{\infty}+\frac{\theta}{2})\bar{Q}^{-1},\bar t]$. Since the
scalar curvature $R_{\infty}$ of the limit on $M_{\infty}\times
[t_{\infty}+\frac{\theta}{3},0]$ is uniformly bounded by some
positive constant $C$, we have the curvature estimate
$$
R(\bar x,t) \leq 2C\bar{Q}
$$
for all $t \in [\bar
t+(t_{\infty}+\frac{\theta}{2})\bar{Q}^{-1},\bar t]$ and all
sufficiently large $m$. For each fixed $m$ and $\alpha_m$, we may
require the chosen $\beta_m$ to satisfy
$$
(\bar{\delta}^{\alpha_m\beta_m})^{-2}\(\widetilde{r}
\(\frac{t^{\alpha_m}_0}{2}\)\)^{-2} \geq
m\widetilde{r}({t^{\alpha_m}_0})^{-2}\geq m\bar{Q}.
$$
When $m$ is large enough, we observe again that Lemma 3.4.1 (ii)
is applicable for $d_t(x_0,\bar{x})$ with $t \in [\bar
t+(t_{\infty}+\frac{\theta}{2})\bar{Q}^{-1},\bar t]$. Then by
repeating the argument as in the derivation of (7.2.1), (7.2.2)
and (7.2.3), we deduce that for all sufficiently large $m$, the
canonical neighborhood statement holds for any $(y,t)$ in the
parabolic neighborhood $P(\bar x,\bar
t,\frac{1}{10}K_1^{\frac{1}{2}}\bar{Q}^{-\frac{1}{2}},
(t_{\infty}+\frac{\theta}{2})\bar{Q}^{-1})$ with $R(y,t) \geq
4\bar{Q}$.

Let $(y_m,\bar t+(t_{\infty}+\theta_m)\bar{Q}^{-1})$ be a sequence
of associated points and times in the (unrescaled) solutions
$g_{ij}^{\alpha_m \beta_m}(\cdot,t)$ so that after rescaling, the
sequence converges to $(y_{\infty},t_{\infty}+\theta)$ in the
limit. Clearly $\frac{\theta}{2}\leq \theta_m \leq 2\theta$ for
all sufficiently large $m$. Then by the argument as in the
derivation of (7.5.13), we know that for all sufficiently large
$m$, the solution $g_{ij}^{\alpha_m \beta_m}(\cdot,t)$ at $y_m$ is
defined on the whole time interval $[\bar
t+(t_{\infty}+\theta_m-\frac{1}{20}\eta^{-1})\bar{Q}^{-1}, \bar
t+(t_{\infty}+\theta_m)\bar{Q}^{-1}]$ and satisfies the curvature
estimate
$$
R(y_m,t) \leq 8\bar{Q}
$$
there; moreover the canonical neighborhood statement holds for any
$(y,t)$ with $R(y,t) \geq 4\bar Q$ in the parabolic neighborhood
$P(y_m,\bar t,\frac{1}{10}K_1^{\frac{1}{2}}\bar{Q}^{-\frac{1}{2}},
(t_{\infty}-\frac{c}{3})\bar{Q}^{-1}).$

We now consider the rescaled sequence
$\widetilde{g}_{ij}^{\alpha_m\beta_m}(\cdot,t)$ with the marked
points replaced by $y_m$ and the times replaced by $s_m \in [\bar
t+(t_{\infty}-\frac{c}{4})\bar{Q}^{-1},
\bar{t}+(t_{\infty}+\frac{c}{4})\bar{Q}^{-1}]$. As before the
Li-Yau-Hamilton inequality implies the rescaling limit around
$(y_m,s_m)$ agrees with the original one. Then the arguments in
previous paragraphs imply the limit is well-defined and smooth on
a space-time open neighborhood of the maximal time slice
$\{t=t_{\infty}\}$. Particularly this excludes the possibility of
existing surgeries in finite distance around the time
$t_{\infty}$. Moreover, the limit at $t=t_{\infty}$ also has
bounded curvature. By using the gradient estimates in the
canonical neighborhood assumption on the parabolic neighborhood
$P(y_m,\bar t,\frac{1}{10}K_1^{\frac{1}{2}}\bar{Q}^{-\frac{1}{2}},
(t_{\infty}-\frac{c}{3})\bar{Q}^{-1})$, we see that the second
possibility is also impossible. Hence we have proved a subsequence
of the rescaled solutions converges to an ancient
$\kappa$-solution.

Therefore we have proved the canonical neighborhood statement
(ii).

\medskip
(iii) This is analogous to Theorem 7.2.1. We also argue by
contradiction. Suppose for some $A<+\infty$ and sequences of
positive numbers $K^{\alpha}_2\rightarrow+\infty$,
$\bar{r}^{\alpha}\rightarrow0$ there exists a sequence of times
$t^{\alpha}_0$ such that for any sequences
$\bar{\delta}^{\alpha\beta}>0$ with
$\bar{\delta}^{\alpha\beta}\rightarrow0$ for fixed $\alpha$, we
have sequences of solutions $g^{\alpha\beta}_{ij}$ to the Ricci
flow with surgery and sequences of points $x^{\alpha\beta}_0$, of
positive constants $r^{\alpha\beta}_0$ with $r^{\alpha\beta}_0\leq
\bar{r}^{\alpha}\sqrt{t^{\alpha}_0}$ which satisfy the
assumptions, but for all $\alpha,\beta$ there hold \be
R(x^{\alpha\beta},t^{\alpha}_0)>K^{\alpha}_2(r^{\alpha\beta}_0)^{-2},\
\mbox{ for some } \ x^{\alpha\beta}\in
B_{t^{\alpha}_0}(x^{\alpha\beta}_0,Ar^{\alpha\beta}_0).
\ee

We may assume that $\bar{\delta}^{\alpha\beta}\leq
\bar{\delta}_{4A}(t^{\alpha}_0)$ for all $\alpha,\beta$, where
$\bar{\delta}_{4A}(t^{\alpha}_0)$ is chosen so that the statements
(i) and (ii) hold on
$B_{t^{\alpha}_0}(x^{\alpha\beta}_0,4Ar^{\alpha\beta}_0)$. Let
$\hat{g}^{\alpha\beta}_{ij}$ be the rescaled solutions of
$g^{\alpha\beta}_{ij}$ around the origins $x^{\alpha\beta}_0$ with
factor $(r^{\alpha\beta}_0)^{-2}$ and shift the times
$t^{\alpha}_0$ to zero. Then by applying the statement (ii), we
know that the regions, where the scalar curvature of the rescaled
solutions $\hat{g}^{\alpha\beta}_{ij}$ is at least
$K_1(=K_1(4A))$, are canonical neighborhood regions. Note that
canonical $\varepsilon$-neck neighborhoods are strong. Also note
that the pinching assumption and the assertion
$$
t^{\alpha}_0(r^{\alpha\beta}_0)^{-2} \geq
(\bar{r}^{\alpha})^{-2}\rightarrow+\infty, \mbox{ as
}\alpha\rightarrow+\infty,
$$
imply that any subsequent limit of the rescaled solutions
$\hat{g}^{\alpha\beta}_{ij}$ must have nonnegative sectional
curvature. Thus by the above argument in the proof of the
statement (ii) (or the argument in Step 2 of the proof of Theorem
7.1.1), we conclude that there exist subsequences
$\alpha=\alpha_m,\beta=\beta_m$ such that the curvatures of the
rescaled solutions $\hat{g}^{\alpha_m\beta_m}_{ij}$ stay uniformly
bounded at distances from the origins $x^{\alpha_m\beta_m}_0$ not
exceeding $2A$. This contradicts (7.5.16) for $m$ sufficiently
large. This proves the statements (iii).

Clearly for fixed $A$, after defining the $\bar{\delta}_A(t_0)$
for each $t_0$, one can adjust the $\bar{\delta}_A(t_0)$ so that
it is nonincreasing in $t_0$.

Therefore we have completed the proof of the theorem.
\end{pf}

{}From now on we redefine the function $\widetilde{\delta}(t)$  so
that it is also less than $\bar{\delta}_{2(t+1)}(2t)$ and then the
above theorem always holds for $A \in [1,2(t_0+1)]$. Particularly,
we still have
$$
\widetilde{\delta}(t)\leq \bar{\delta}(t) = \min
\{\frac{1}{{2e^2\log(1+t)}}, \delta_0\},
$$
which tends to zero as $t\rightarrow+\infty$. We may also require
that $\widetilde{r}(t)$ tends to zero as $t\rightarrow+\infty$.

The next result is a version of Theorem 7.2.2 for solutions with
surgery.

\begin{theorem}[{Perelman \cite{P2}}]
For any $\varepsilon > 0$ and $w>0$, there exist
$\tau=\tau(w,\varepsilon)>0$, $K=K(w,\varepsilon)<+\infty$,
$\bar{r}=\bar{r}(w,\varepsilon)>0$,
$\theta=\theta(w,\varepsilon)>0$ and $T=T(w) < +\infty$ with the
following property. Suppose we have a solution, constructed by
Theorem $7.4.3$ with the nonincreasing $($continuous$)$ positive
functions $\widetilde{\delta}(t)$ and $\widetilde{r}(t)$, to the
Ricci flow with surgery on the time interval $[0,t_0]$ with a
compact orientable normalized three-manifold as initial data,
where each $\delta$-cutoff at a time $t \in [0,t_0]$ has
$\delta=\delta(t) \leq \min
\{\widetilde{\delta}(t),\widetilde{r}(2t)\}$. Let $r_0,t_0$
satisfy $\theta^{-1}h\leq r_0\leq \bar{r}\sqrt{t_0}$ and $t_0 \geq
T$, where $h$ is the maximal cutoff radius for surgeries in
$[\frac{t_0}{2},t_0]$ $($if there is no surgery in the time
interval $[\frac{t_0}{2},t_0]$, we take $h=0),$ and assume that
the solution on the ball $B_{t_0}(x_0,r_0)$ satisfies
$$
Rm(x,t_0) \geq -r^{-2}_0,\; \mbox{ on }\; B_{t_0}(x_0,r_0),
$$
$$
\mbox{and }\quad  \Vol_{t_0}(B_{t_0}(x_0,r_0)) \geq w r^3_0.
$$
Then the solution is well defined and satisfies
$$
R(x,t)<Kr^{-2}_0
$$
in the whole parabolic neighborhood
$$
P\(x_0,t_0,\frac{r_0}{4},-\tau r^2_0\) =\left\{(x,t)\ |\ x\in
B_t\(x_0,\frac{r_0}{4}\),t\in[t_0-\tau r^2_0,t_0]\right\}.
$$
\end{theorem}

\begin{pf}
We are given that $Rm(x,t_0)\geq -r^{-2}_0$ for $x \in
B_{t_0}(x_0,r_0)$, and $\Vol_{t_0}(B_{t_0}(x_0,r_0))\geq w r^3_0$.
The same argument in the derivation of (7.2.7) and (7.2.8) (by
using the Alexandrov space theory) implies that there exists a
ball $B_{t_0}(x',r')\subset B_{t_0}(x_0,r_0)$ with \be
\Vol_{t_0}(B_{t_0}(x',r'))\geq \frac{1}{2}\alpha_3(r')^3
\ee and with \be
r'\geq c(w)r_0 
\ee for some small positive constant $c(w)$ depending only on $w$,
where $\alpha_3$ is the volume of the unit ball in $\mathbb{R}^3$.

As in (7.1.2), we can rewrite the pinching assumption (7.3.3) as
$$
Rm \geq-[f^{-1}(R(1+t))/(R(1+t))]R,
$$
where
$$
y= f(x) = x(\log x - 3), \quad \mbox{ for }\; e^2 \leq x <
+\infty,
$$
is increasing and convex with range $-e^2 \leq y < +\infty$, and
its inverse function is also increasing and satisfies
$$
\lim_{y\rightarrow +\infty} f^{-1}(y)/y = 0.
$$

Note that $t_0 r^{-2}_0 \geq \bar{r}^{-2}$ by the hypotheses. We
may require $T(w) \geq 8c(w)^{-1}$. Then by applying Theorem 7.5.1
(iii) with $A=8c(w)^{-1}$ and combining with the pinching
assumption, we can reduce the proof of the theorem to the special
case $w=\frac{1}{2}\alpha_3$. In the following we simply assume
$w=\frac{1}{2}\alpha_3$.

Let us first consider the case $r_0 < \widetilde{r}(t_0)$. We
claim that $R(x,t_0)\leq C^2_0r^{-2}_0$ on
$B_{t_0}(x_0,\frac{r_0}{3})$, for some sufficiently large positive
constant $C_0$ depending only on $\varepsilon$. If not, then there
is a canonical neighborhood around $(x,t_0)$. Note that the type
(c) canonical neighborhood has already been ruled out by our
design of cutoff surgeries. Thus $(x,t_0)$ belongs to an
$\varepsilon$-neck or an $\varepsilon$-cap. This tells us that
there is a nearby point $y$, with $R(y,t_0)\geq C^{-1}_2R(x,t_0) >
C^{-1}_2C^2_0r^{-2}_0$ and $d_{t_0}(y,x)\leq
C_1R(x,t_0)^{-\frac{1}{2}}\leq C_1C^{-1}_0r_0$, which is the
center of the $\varepsilon$-neck
$B_{t_0}(y,\varepsilon^{-1}R(y,t_0)^{-\frac{1}{2}})$. (Here
$C_1,C_2$ are the positive constants in the definition of
canonical neighborhood assumption). Clearly, when we choose $C_0$
to be much larger than $C_1,C_2$ and $\varepsilon^{-1}$, the whole
$\varepsilon$-neck
$B_{t_0}(y,\varepsilon^{-1}R(y,t_0)^{-\frac{1}{2}})$ is contained
in $B_{t_0}(x_0,\frac{r_0}{2})$ and we have \be
\frac{\Vol_{t_0}(B_{t_0}(y,\varepsilon^{-1}
R(y,t_0)^{-\frac{1}{2}}))}{(\varepsilon^{-1}
R(y,t_0)^{-\frac{1}{2}})^3}\leq
8\pi\varepsilon^2. 
\ee Without loss of generality, we may assume $\varepsilon>0$ is
very small. Since we have assumed that $Rm\geq -r^{-2}_0$ on
$B_{t_0}(x_0,r_0)$ and $\Vol_{t_0}(B_{t_0}(x_0,r_0))\geq
\frac{1}{2}\alpha_3r^3_0$, we then get a contradiction by applying
the standard Bishop-Gromov volume comparison. Thus we have the
desired curvature estimate $R(x,t_0)\leq C^2_0r^{-2}_0$ on
$B_{t_0}(x_0,\frac{r_0}{3})$.

Furthermore, by using the gradient estimates in the definition of
canonical neighborhood assumption, we can take
$K=10C^2_0,\tau=\frac{1}{100}\eta^{-1}C^{-2}_0$ and
$\theta=\frac{1}{5}C^{-1}_0$ in this case. And since
$r_0\geq\theta^{-1}h$, we have $R<10C^2_0r^{-2}_0\leq
\frac{1}{2}h^{-2}$ and the surgeries do not interfere in
$P(x_0,t_0,\frac{r_0}{4},-\tau r^2_0)$.

We now consider the remaining case $\widetilde{r}(t_0)\leq r_0\leq
\bar{r}\sqrt{t_0}$. Let us redefine
\begin{align*}
\tau&=\min\left\{\frac{\bar{\tau}_0}{2},\frac{1}{100}
\eta^{-1}C^{-2}_0 \right\},\\
K&=\max\left\{2\(\bar{C}+\frac{2\bar{B}}{\bar{\tau}_0}\),25C^2_0\right\},
\end{align*}
and
$$
\theta=\frac{1}{2}K^{-\frac{1}{2}}
$$
where $\bar{\tau}_0=\tau_0(w),\bar{B}=B(w)$ and $\bar{C}=C(w)$ are
the positive constants in Theorem 6.3.3(ii) with
$w=\frac{1}{2}\alpha_3$, and $C_0$ is the positive constant chosen
above. We will show there is a sufficiently small $\bar{r}>0$ such
that the conclusion of the theorem for $w=\frac{1}{2}\alpha_3$
holds for the chosen $\tau,K$ and $\theta$.

Argue by contradiction. Suppose not, then there exist a sequence
of $\bar{r}^{\alpha}\rightarrow0$, and a sequence of solutions
$g^{\alpha}_{ij}$ with points $(x^{\alpha}_0,t^{\alpha}_0)$ and
radii $r^{\alpha}_0$ such that the assumptions of the theorem do
hold with $\widetilde{r}(t^{\alpha}_0)\leq r^{\alpha}_0\leq
\bar{r}^{\alpha}\sqrt{t^{\alpha}_0}$ whereas the conclusion does
not. Similarly as in the proof of Theorem 7.2.2, we claim that we
may assume that for all sufficiently large $\alpha$, any other
point $(x^{\alpha},t^{\alpha})$ and radius $r^{\alpha}>0$ with
that property has either $t^{\alpha}> t^{\alpha}_0$ or
$t^{\alpha}= t^{\alpha}_0$ with $r^{\alpha}\geq r^{\alpha}_0$;
moreover $t^{\alpha}$ tends to $+\infty$ as $\alpha \rightarrow
+\infty$. Indeed, for fixed $\alpha$ and the solution
$g^{\alpha}_{ij}$, let $t^{\alpha}_{\min}$ be the infimum of all
possible times $t^{\alpha}$ with some point $x^{\alpha}$ and some
radius $r^{\alpha}$ having that property. Since each such
$t^{\alpha}$ satisfies $\bar{r}^{\alpha}\sqrt{t^{\alpha}}\geq
r^{\alpha}\geq \widetilde{r}(t^{\alpha})$, it follows that when
$\alpha$ is large, $t^{\alpha}_{\min}$ must be positive and very
large. Clearly for each fixed sufficiently large $\alpha$, by
passing to a limit, there exist some point $x^{\alpha}_{\min}$ and
some radius $r^{\alpha}_{\min}(\geq
\widetilde{r}(t^{\alpha}_{\min})>0)$ so that all assumptions of
the theorem still hold for $(x^{\alpha}_{\min},t^{\alpha}_{\min})$
and $r^{\alpha}_{\min}$, whereas the conclusion of the theorem
does not hold with $R\geq K(r^{\alpha}_{\min})^{-2}$ somewhere in
$P(x^{\alpha}_{\min},t^{\alpha}_{\min},\frac{1}{4}
r^{\alpha}_{\min},-\tau(r^{\alpha}_{\min})^2)$ for all
sufficiently large $\alpha$. Here we used the fact that if $R <
K(r^{\alpha}_{\min})^{-2}$ on
$P(x^{\alpha}_{\min},t^{\alpha}_{\min},\frac{1}{4}
r^{\alpha}_{\min},-\tau(r^{\alpha}_{\min})^2)$, then there is no
$\delta$-cutoff surgery there; otherwise there must be a point
there with scalar curvature at least
$\frac{1}{2}\delta^{-2}(\frac{t^{\alpha}_{\min}}{2})
(\widetilde{r}(\frac{t^{\alpha}_{\min}}{2}))^{-2}
\geq\frac{1}{2}(\widetilde{r}(t^{\alpha}_{\min}))^{-2}
(\widetilde{r}(\frac{t^{\alpha}_{\min}}{2}))^{-2} \gg K
(r^{\alpha}_{\min})^{-2}$ since
$\bar{r}^{\alpha}\sqrt{t^{\alpha}_{\min}}\geq
r^{\alpha}_{\min}\geq \widetilde{r}(t^{\alpha}_{\min})$ and
$\bar{r}^{\alpha} \rightarrow 0$, which is a contradiction.

After choosing the first time $t^{\alpha}_{\min}$, by passing to a
limit again, we can then choose $r^{\alpha}_{\min}$ to be the
smallest radius for all possible
$(x^{\alpha}_{\min},t^{\alpha}_{\min})$'s and
$r^{\alpha}_{\min}$'s with that property. Thus we have verified
the claim.

For simplicity, we will drop the index $\alpha$ in the following
arguments. By the assumption and the standard volume comparison,
we have
$$
\Vol_{t_0}\(B_{t_0}\(x_0,\frac{1}{2}r_0\)\)\geq \xi_0r^3_0
$$
for some universal positive $\xi_0$. As in deriving (7.2.7) and
(7.2.8), we can find a ball $B_{t_0}(x'_0,r'_0)\subset
B_{t_0}(x_0,\frac{r_0}{2})$ with
$$
\Vol_{t_0}(B_{t_0}(x'_0,r'_0)) \geq \frac{1}{2}\alpha_3(r'_0)^3\;
\mbox{ and }\; \frac{1}{2}r_0\geq r'_0\geq \xi'_0r_0
$$
for some universal positive constant $\xi'_0$. Then by what we had
proved in the previous case and by the choice of the first time
$t_0$ and the smallest radius $r_0$, we know that the solution is
defined in $P(x'_0,t_0,\frac{r'_0}{4},-\tau (r'_0)^2)$ with the
curvature bound $$R<K(r'_0)^{-2}\leq K(\xi'_0)^{-2}r^{-2}_0.$$
Since $\bar{r}\sqrt{t_0}\geq r_0\geq \widetilde{r}(t_0)$ and
$\bar{r}\rightarrow 0$ as $\alpha\rightarrow\infty$, we see that
$t_0\rightarrow +\infty$ and $t_0r^{-2}_0\rightarrow+\infty$ as
$\alpha\rightarrow+\infty$. Define $T(w) = 8c(w)^{-1} +
\bar{\xi}$, for some suitable large universal positive constant
$\bar{\xi}$. Then for $\alpha$ sufficiently large, we can apply
Theorem 7.5.1(iii) and the pinching assumption to conclude that
\be R\leq K'r^{-2}_0,\quad \mbox{ on }\;
P\(x_0,t_0,4r_0,-\frac{\tau}{2}(\xi'_0)^2r^2_0\), 
\ee for some positive constant $K'$ depending only on $K$ and
$\xi'_0$.

Furthermore, by combining with the pinching assumption, we deduce
that when $\alpha$ sufficiently large,
\begin{align}
Rm &  \geq -[f^{-1}(R(1+t))/(R(1+t))]R\\
&  \geq -r^{-2}_0 \nn
\end{align}  
on $P(x_0,t_0,r_0,-\frac{\tau}{2}(\xi'_0)^2r^2_0)$.  So by
applying Theorem 6.3.3(ii) with $w=\frac{1}{2}\alpha_3$, we have
that when $\alpha$ sufficiently large, \be
\Vol_t(B_t(x_0,r_0)) \geq \xi_1 r^3_0,  
\ee for all $ t \in [t_0 -\frac{\tau}{2}(\xi'_0)^2r^2_0, t_0]$,
and \be R\leq\(\bar{C}+\frac{2\bar{B}}{\bar{\tau}_0}\)r^{-2}_0\leq
\frac{1}{2}Kr^{-2}_0 
\ee on $P(x_0,t_0,\frac{r_0}{4},-\frac{\tau}{2}(\xi'_0)^2r^2_0)$,
where $\xi_1$ is some universal positive constant.

Next we want to extend the estimate (7.5.23) backwards in time.
Denote by $t_1=t_0-\frac{\tau}{2}(\xi'_0)^2r^2_0$. The estimate
(7.5.22) gives
$$
\Vol_{t_1}(B_{t_1}(x_0,r_0))\geq \xi_1r^3_0.
$$
By the same argument in the derivation of (7.2.7) and (7.2.8)
again, we can find a ball $B_{t_1}(x_1,r_1)\subset
B_{t_1}(x_0,r_0)$ with
$$
\Vol_{t_1}(B_{t_1}(x_1,r_1)) \geq \frac{1}{2}\alpha_3r^3_1
$$
and with
$$
r_1\geq \xi'_1r_0
$$
for some universal positive constant $\xi'_1$. Then by what we had
proved in the previous case and by the lower bound (7.5.21) at
$t_1$ and the choice of the first time $t_0$, we know that the
solution is defined on $P(x_1,t_1,\frac{r_1}{4},-\tau r^2_1)$ with
the curvature bound $R<K r^{-2}_1$. By applying Theorem 7.5.1(iii)
and the pinching assumption again we get that for $\alpha$
sufficiently large,
$$
R\leq K''r^{-2}_0 \leqno(7.5.20)'
$$ on
$P(x_0,t_1,4r_0,-\frac{\tau}{2}(\xi'_1)^2r^2_0)$, for some
positive constant $K''$ depending only on $K$ and $\xi'_1$.
Moreover, by combining with the pinching assumption, we have
\begin{align*}
Rm&  \geq-[f^{-1}(R(1+t))/(R(1+t))]R \tag*{(7.5.21)$'$}\\
&  \geq -r^{-2}_0
\end{align*}  
on $P(x_0,t_1,4r_0,-\frac{\tau}{2}(\xi'_1)^2r^2_0)$, for $\alpha$
sufficiently large. So by applying Theorem 6.3.3 (ii) with
$w=\frac{1}{2}\alpha_3$ again, we have that for $\alpha$
sufficiently large,
$$
\Vol_t(B_t(x_0,r_0)) \geq \xi_1 r^3_0,\leqno (7.5.22)'
$$
for all $ t \in [t_0 -\frac{\tau}{2}(\xi'_0)^2r^2_0
-\frac{\tau}{2}(\xi'_1)^2r^2_0, t_0]$, and
$$
R\leq \(\bar{C}+\frac{2\bar{B}}{\bar{\tau}_0}\)r^{-2}_0\leq
\frac{1}{2}Kr^{-2}_0 \leqno (7.5.23)'
$$
on $P(x_0,t_0,\frac{r_0}{4},-\frac{\tau}{2}(\xi'_0)^2
r^2_0-\frac{\tau}{2}(\xi'_1)^2r^2_0)$.

Note that the constants $\xi_0$, $\xi'_0$, $\xi_1$ and $\xi'_1$
are universal, independent of the time $t_1$ and the choice of the
ball $B_{t_1}(x_1,r_1)$. Then we can repeat the above procedure as
many times as we like, until we reach the time $t_0-\tau r^2_0$.
Hence we obtain the estimate
$$
R\leq(\bar{C}+\frac{2\bar{B}}{\bar{\tau}_0})r^{-2}_0\leq
\frac{1}{2}Kr^{-2}_0 \leqno (7.5.23)''
$$
on $P(x_0,t_0,\frac{r_0}{4},-\tau r^2_0)$, for sufficiently large
$\alpha$. This contradicts the choice of the point $(x_0,t_0)$ and
the radius $r_0$ which make $R\geq Kr^{-2}_0$ somewhere in
$P(x_0,t_0,\frac{r_0}{4},-\tau r^2_0)$.

Therefore we have completed the proof of the theorem.
\end{pf}

Consequently, we have the following result which is analog of
Corollary 7.2.4. This result is a weak version of a claim in the
section 7.3 of {Perelman \cite{P2}}.

\begin{corollary}
For any $\varepsilon > 0$ and $w>0$, there exist
$\bar{r}=\bar{r}(w,\varepsilon)>0$,
$\theta=\theta(w,\varepsilon)>0$ and $T=T(w)$ with the following
property. Suppose we have a solution, constructed by Theorem
$7.4.3$ with the positive functions $\widetilde{\delta}(t)$ and
$\widetilde{r}(t)$, to the Ricci flow with surgery with a compact
orientable normalized three-manifold as initial data, where each
$\delta$-cutoff at a time $t$ has $\delta=\delta(t)\leq
\min\{\widetilde{\delta}(t), \widetilde{r}(2t)\}$. If
$B_{t_0}(x_0,r_0)$ is a geodesic ball at time $t_0$, with
$\theta^{-1}h\leq r_0\leq \bar{r}\sqrt{t_0}$ and $t_0 \geq T$,
where $h$ is the maximal cutoff radii for surgeries in
$[\frac{t_0}{2},t_0]$ (if there is no surgery in the time interval
$[\frac{t_0}{2},t_0]$, we take $h=0$), and satisfies
$$
\min\{Rm(x,t_0)\ |\ x\in B_{t_0}(x_0,r_0)\}=-r^{-2}_0,
$$
then
$$
\Vol_{t_0}(B_{t_0}(x_0,r_0))<w r^3_0.
$$
\end{corollary}

\begin{pf}
We argue by contradiction. Let $\theta=\theta(w,\varepsilon)$  and
$T=2T(w)$, where $\theta(w,\varepsilon)$ and $T(w)$ are the
positive constant in Theorem 7.5.2. Suppose for any $\bar{r}>0$
there is a solution and a geodesic ball $B_{t_0}(x_0,r_0)$
satisfying the assumptions of the corollary with $\theta^{-1}h\leq
r_0\leq \bar{r}\sqrt{t_0}$ and $t_0 \geq T$, and with
$$
\min\{Rm(x,t_0)\ |\ x\in B_{t_0}(x_0,r_0)\}=-r^{-2}_0,
$$
but
$$
\Vol_{t_0}(B_{t_0}(x_0,r_0))\geq w r^3_0.
$$
Without loss of generality, we may assume that $\bar{r}$ is less
than the corresponding constant in Theorem 7.5.2. We can then
apply Theorem 7.5.2 to get $$R(x,t)\leq Kr^{-2}_0 $$ whenever
$t\in[t_0-\tau r^2_0,t_0]$ and $d_t(x,x_0)\leq \frac{r_0}{4}$,
where $\tau$ and $K$ are the positive constants in Theorem 7.5.2.
Note that $t_0r^{-2}_0\geq \bar{r}^{-2}\rightarrow+\infty$ as
$\bar{r}\rightarrow 0$. By combining with the pinching assumption
we have
\begin{align*}
Rm&  \geq-[f^{-1}(R(1+t))/(R(1+t))]R\\
&  \geq -\frac{1}{2}r^{-2}_0
\end{align*}
in the region $P(x_0,t_0,\frac{r_0}{4},-\tau r^2_0)= \{(x,t)\ |\
x\in B_t(x_0,\frac{r_0}{4}),t\in[t_0-\tau r^2_0,t_0]\}$, provided
$\bar{r}>0$ is sufficiently small. Thus we get the estimate
$$|Rm|\leq K'r^{-2}_0
$$
in $P(x_0,t_0,\frac{r_0}{4},-\tau r^2_0)$, where $K'$ is a
positive constant depending only on $w$ and $\varepsilon$.

We can now apply Theorem 7.5.1 (iii) to conclude that
$$
R(x,t)\leq \widetilde{K}r^{-2}_0
$$
whenever $t\in[t_0-\frac{\tau}{2}r^2_0,t_0]$ and $d_t(x,x_0)\leq
r_0$, where $\widetilde{K}$ is a positive constant depending only
on $w$ and $\varepsilon$. By using the pinching assumption again
we further have
$$
Rm(x,t)\geq -\frac{1}{2}r^{-2}_0
$$
in the region $P(x_0,t_0,r_0,-\frac{\tau}{2}r^2_0)=\{(x,t)\ |\
x\in B_t(x_0,r_0),t\in[t_0-\frac{\tau}{2}r^2_0,t_0]\}$, as long as
$\bar{r}$ is sufficiently small. In particular, this would imply
$$
\min\{Rm(x,t_0)\ |\ x\in B_{t_0}(x_0,r_0)\}>-r^{-2}_0,
$$
which is a contradiction.  
\end{pf}

\begin{remark}
In section 7.3 of \cite{P2}, Perelman claimed a stronger statement
than the above Corollary 7.5.3 that allows $r_0 < \theta^{-1}h$ in
the assumptions. Nevertheless, the above weaker statement is
sufficient to deduce the geometrization result.
\end{remark}

\section{Long Time Behavior}

In Section 5.3, we obtained the long time behavior for smooth
(compact) solutions to the three-dimensional Ricci flow with
bounded normalized curvature. The purpose of this section is to
adapt Hamilton's arguments there to solutions of the Ricci flow
with surgery and to drop the bounded normalized curvature
assumption as sketched by Perelman \cite{P2}.

Recall from Corollary 7.4.4 that we have completely understood the
topological structure of a compact, orientable three-manifold with
nonnegative scalar curvature. From now on we assume that our
initial manifold does not admit any metric with nonnegative scalar
curvature, and that \textbf{once we get a compact component with
nonnegative scalar curvature, it is immediately removed}.
Furthermore, if a solution to the Ricci flow with surgery becomes
extinct in a finite time, we have also obtained the topological
structure of the initial manifold. So in the following we only
consider those solutions to the Ricci flow with surgery which
exist for all times $t \geq 0$.

Let $g_{ij}(t)$, $0 \leq t < +\infty$, be a solution to the Ricci
flow with $\delta$-cutoff surgeries, constructed by Theorem 7.4.3
with normalized initial data. Let $0<t_1<t_2<\cdots<t_k<\cdots$ be
the surgery times, where each $\delta$-cutoff at a time $t_k$ has
$\delta = \delta(t_k) \leq \min
\{\widetilde{\delta}(t_k),\widetilde{r}(2t_k)\}$. On each time
interval $(t_{k-1},t_k)$ (denote by $t_0=0$), the scalar curvature
satisfies the evolution equation \be \frac{\partial}{\partial
t}R=\Delta
R+2|\stackrel{\circ}{\Ric}|^2+\frac{2}{3}R^2 
\ee where $\stackrel{\circ}{\Ric}$ is the trace-free part of
$\Ric.$ Then $R_{\min}(t)$, the minimum of the scalar curvature at
the time $t$, satisfies
$$
\frac{d}{dt}R_{\min}(t)\geq \frac{2}{3}R^2_{\min}(t)
$$
for $t\in(t_{k-1},t_k)$, for each $k = 1, 2, \ldots$. Since our
surgery procedure had removed all components with nonnegative
scalar curvature, the minimum $R_{\min}(t)$ is negative for all $t
\in [0,+\infty)$. Also recall that the cutoff surgeries were
performed only on $\delta$-necks. Thus the surgeries do not occur
at the parts where $R_{\min}(t)$ are achieved. So the differential
inequality
$$
\frac{d}{dt}R_{\min}(t)\geq \frac{2}{3}R^2_{\min}(t)
$$
holds for all $t\geq 0$, and then by normalization,
$R_{\min}(0)\geq -1$, we have \be R_{\min}(t)\geq
-\frac{3}{2}\cdot\frac{1}{t+\frac{3}{2}}, \
\mbox{ for all } \ t\geq 0. 
\ee Meanwhile, on each time interval $(t_{k-1},t_k)$, the volume
satisfies the evolution equation
$$
\frac{d}{dt}V=-\int RdV
$$
and then by (7.6.2),
$$
\frac{d}{dt}V\leq \frac{3}{2}\cdot\frac{1}{(t+\frac{3}{2})}V.
$$
Since the cutoff surgeries do not increase volume, we thus have
\be \frac{d}{dt}\log\(V(t)\(t+\frac{3}{2}\)^{-\frac{3}{2}}\)\leq 0
\ee for all $t\geq 0$. Equivalently, the function
$V(t)(t+\frac{3}{2})^{-\frac{3}{2}}$ is nonincreasing on
$[0,+\infty)$.

We can now use the monotonicity of the function
$V(t)(t+\frac{3}{2})^{-\frac{3}{2}}$ to extract the information of
the solution at large times. On each time interval
$(t_{k-1},t_k)$, we have
\begin{align*}
&\frac{d}{dt}\log\(V(t)\(t+\frac{3}{2}\)^{-\frac{3}{2}}\) \\
&=-\(R_{\min}(t)+\frac{3}{2\(t+\frac{3}{2}\)}\)
+\frac{1}{V}\int_M(R_{\min}(t)-R)dV.
\end{align*}
Then by noting that the cutoff surgeries do not increase volume,
we get
\begin{align}
\frac{V(t)}{(t+\frac{3}{2})^{\frac{3}{2}}}
&\leq\frac{V(0)}{(\frac{3}{2})^{\frac{3}{2}}}
\exp\Bigg\{-\int^t_0\(R_{\min}(t)
+\frac{3}{2(t+\frac{3}{2})}\)dt \\
&\qquad-\int^t_0\frac{1}{V}
\int_M(R-R_{\min}(t))dVdt\Bigg\} \nn 
\end{align}
for all $t>0$. Now by this inequality and the equation (7.6.1), we
obtain the following consequence (cf. Lemma 7.1 of Hamilton
\cite{Ha99} and section 7.1 of Perelman \cite{P2}).

\begin{lemma}
Let $g_{ij}(t)$ be a solution to the Ricci flow with surgery,
constructed by Theorem $7.4.3$ with normalized initial data. If for
a fixed $0<r<1$ and a sequence of times
$t^{\alpha}\rightarrow\infty$, the rescalings of the solution on the
parabolic neighborhoods
$P(x^{\alpha},t^{\alpha},r\sqrt{t^{\alpha}},-r^2t^{\alpha})=\{(x,t)\
|\ x\in B_t(x^{\alpha},r\sqrt{t^{\alpha}}),
t\in[t^{\alpha}-r^2t^{\alpha},t^{\alpha}]\}$, with factor
$(t^{\alpha})^{-1}$ and shifting the times $t^{\alpha}$ to $1$,
converge in the $C^{\infty}$ topology to some smooth limiting
solution, defined in an abstract parabolic neighborhood
$P(\bar{x},1,r,-r^2)$, then this limiting solution has constant
sectional curvature $-1/4t$ at any time $t\in [1-r^2,1]$.
\end{lemma}

In the previous section we obtained several curvature estimates
for the solutions to the Ricci flow with surgery. Now we combine
the curvature estimates with the above lemma to derive the
following asymptotic result for the curvature.

\begin{lemma}[{Perelman \cite{P2}}]
For any $\varepsilon > 0$, let $g_{ij}(t)$, $0 \leq t < +\infty$,
be a solution to the Ricci flow with surgery, constructed by
Theorem $7.4.3$ with normalized initial data.
\begin{itemize}
\item[(i)] Given $w>0$, $r>0$, $\xi>0$, one can find
$T=T(w,r,\xi,\varepsilon)<+\infty$ such that if the geodesic ball
$B_{t_0}(x_0,r\sqrt{t_0})$ at some time $t_0\geq T$ has volume at
least $w r^3t^{\frac{3}{2}}_0$ and the sectional curvature at
least $-r^{-2}t^{-1}_0$, then the curvature at $x_0$ at time
$t=t_0$ satisfies \be
|2tR_{ij}+g_{ij}|<\xi. 
\ee \item[(ii)] Given in addition $1 \leq A<\infty$ and allowing
$T$ to depend on $A$, we can ensure $(7.6.5)$ for all points in
$B_{t_0}(x_0,Ar\sqrt{t_0})$. \item[(iii)] The same is true for all
points in the forward parabolic neighborhood
$P(x_0,t_0,Ar\sqrt{t_0},Ar^2t_0)\triangleq\{(x,t)\ |\ x\in
B_t(x_0,Ar\sqrt{t_0}),t\in[t_0,t_0+Ar^2t_0]\}.$
\end{itemize}
\end{lemma}

\begin{pf}
(i) By the assumptions and the standard volume comparison, we have
$$
\Vol_{t_0}(B_{t_0}(x_0,\rho))\geq  cw \rho^3
$$
for all $0<\rho\leq r\sqrt{t_0}$, where $c$ is a universal
positive constant. Let $\bar{r}=\bar{r}(cw,\varepsilon)$ be the
positive constant in Theorem 7.5.2 and set
$r_0=\min\{r,\bar{r}\}$. On $B_{t_0}(x_0,r_0\sqrt{t_0})(\subset
B_{t_0}(x_0,r\sqrt{t_0}))$, we have
\begin{gather}
Rm\geq -(r_0\sqrt{t_0})^{-2} \\
\mbox{and }\qquad \Vol_{t_0}(B_{t_0}(x_0,r_0\sqrt{t_0}))\geq
cw(r_0\sqrt{t_0})^3.\nn
\end{gather}
Obviously, there holds $\theta^{-1}h\leq r_0\sqrt{t_0}\leq
\bar{r}\sqrt{t_0}$ when $t_0$ is large enough, where
$\theta=\theta(cw,\varepsilon)$ is the positive constant in
Theorem 7.5.2 and $h$ is the maximal cutoff radius for surgeries
in $[\frac{t_0}{2},t_0]$ (if there is no surgery in the time
interval $[\frac{t_0}{2},t_0]$, we take $h=0$). Then it follows
from Theorem 7.5.2 that the solution is defined and satisfies
$$
R<K(r_0\sqrt{t_0})^{-2}
$$
on whole parabolic neighborhood
$P(x_0,t_0,\frac{r_0\sqrt{t_0}}{4},-\tau(r_0\sqrt{t_0})^2)$. Here
$\tau=\tau(cw,\varepsilon)$ and $K=K(cw,\varepsilon)$ are the
positive constants in Theorem 7.5.2. By combining with the
pinching assumption we have
\begin{align*}
Rm&  \geq-[f^{-1}(R(1+t))/(R(1+t))]R\\
&  \geq -{\rm const.}\,K(r_0\sqrt{t_0})^{-2}
\end{align*}
in the region $P(x_0,t_0, \frac{r_0\sqrt{t_0}}{4},
-\tau(r_0\sqrt{t_0})^2)$. Thus we get the estimate \be
|Rm|\leq K'(r_0\sqrt{t_0})^{-2} 
\ee on
$P(x_0,t_0,\frac{r_0\sqrt{t_0}}{4},-\tau(r_0\sqrt{t_0})^2)$, for
some positive constant $K'=K'(w,\varepsilon)$ depending only on
$w$ and $\varepsilon$.

The curvature estimate (7.6.7) and the volume estimate (7.6.6)
ensure that as $t_0\rightarrow+\infty$ we can take smooth
(subsequent) limits for the rescalings of the solution with factor
$(t_0)^{-1}$ on parabolic neighborhoods
$P(x_0,t_0,\frac{r_0\sqrt{t_0}}{4},$ $-\tau(r_0\sqrt{t_0})^2)$.
Then by applying Lemma 7.6.1, we can find
$T=T(w,r,\xi,\varepsilon)<+\infty$ such that when $t_0\geq T$,
there holds \be
|2tR_{ij}+g_{ij}|<\xi, 
\ee on
$P(x_0,t_0,\frac{r_0\sqrt{t_0}}{4},-\tau(r_0\sqrt{t_0})^2)$, in
particular,
$$
|2tR_{ij}+g_{ij}|(x_0,t_0)<\xi.
$$
This proves the assertion (i).

\medskip
(ii) In view of the above argument, to get the estimate (7.6.5)
for all points in $B_{t_0}(x_0,Ar\sqrt{t_0})$, the key point is to
get a upper bound for the scalar curvature on the parabolic
neighborhood $P(x_0,t_0,Ar\sqrt{t_0},-\tau(r_0\sqrt{t_0})^2)$.
After having the estimates (7.6.6) and (7.6.7), one would like to
use Theorem 7.5.1(iii) to obtain the desired scalar curvature
estimate. Unfortunately it does not work since our $r_0$ may be
much larger than the constant $\bar{r}(A,\varepsilon)$ there when
$A$ is very large. In the following we will use Theorem 7.5.1(ii)
to overcome the difficulty.

Given $1 \leq A<+\infty$, based on (7.6.6) and (7.6.7), we can use
Theorem 7.5.1(ii) to find a positive constant
$K_1=K_1(w,r,A,\varepsilon)$ such that each point in
$B_{t_0}(x_0,2Ar\sqrt{t_0})$ with its scalar curvature at least
$K_1(r\sqrt{t_0})^{-2}$ has a canonical neighborhood. We claim
that there exists $T=T(w,r,A,\varepsilon)<+\infty$ so that when
$t_0\geq T$, we have \be R<K_1(r\sqrt{t_0})^{-2},\; \mbox{ on } \;
B_{t_0}(x_0,2Ar\sqrt{t_0}). 
\ee

Argue by contradiction. Suppose not; then there exist a sequence
of times $t^{\alpha}_0\rightarrow+\infty$ and sequences of points
$x^{\alpha}_0$, $x^{\alpha}$ with $x^{\alpha}\in
B_{t^{\alpha}_0}(x^{\alpha}_0,2Ar\sqrt{t^{\alpha}_0})$ and
$R(x^{\alpha},t^{\alpha}_0)=K_1(r\sqrt{t^{\alpha}_0})^{-2}$. Since
there exist canonical neighborhoods ($\varepsilon$-necks or
$\varepsilon$-caps) around the points $(x^{\alpha},t^{\alpha}_0)$,
there exist positive constants $c_1$, $C_2$ depending only on
$\varepsilon$ such that
$$
\Vol_{t^{\alpha}_0}(B_{t^{\alpha}_0}(x^{\alpha},
K^{-\frac{1}{2}}_1(r\sqrt{t^{\alpha}_0})))\geq
c_1(K^{-\frac{1}{2}}_1(r\sqrt{t^{\alpha}_0}))^3
$$
and
$$
C^{-1}_2K_1(r\sqrt{t^{\alpha}_0})^{-2} \leq R(x,t^{\alpha}_0) \leq
C_2K_1(r\sqrt{t^{\alpha}_0})^{-2},
$$
on $B_{t^{\alpha}_0} (x^{\alpha}, K^{-\frac{1}{2}}_1
(r\sqrt{t^{\alpha}_0}))$, for all $\alpha$. By combining with the
pinching assumption we have
\begin{align*}
Rm&  \geq-[f^{-1}(R(1+t))/(R(1+t))]R\\
&  \geq -{\rm const.}\, C_2K_1(r\sqrt{t^{\alpha}_0})^{-2},
\end{align*}
on
$B_{t^{\alpha}_0}(x^{\alpha},K^{-\frac{1}{2}}_1(r\sqrt{t^{\alpha}_0}))$,
for all $\alpha$. It then follows from the assertion (i) we just
proved that
$$
\lim_{\alpha\rightarrow+\infty}|2tR_{ij}+g_{ij}|
(x^{\alpha},t^{\alpha}_0)=0.
$$
In particular, we have
$$
t^{\alpha}_0R(x^{\alpha},t^{\alpha}_0)<-1
$$
for $\alpha$ sufficiently large. This contradicts our assumption
that $R(x^{\alpha},t^{\alpha}_0)=K_1(r\sqrt{t^{\alpha}_0})^{-2}$.
So we have proved assertion (7.6.9).

Now by combining (7.6.9) with the pinching assumption as before,
we have \be
Rm\geq  -K_2(r\sqrt{t_0})^{-2} 
\ee on $B_{t_0}(x_0,2Ar\sqrt{t_0})$, where
$K_2=K_2(w,r,A,\varepsilon)$ is some positive constant depending
only on $w$, $r$, $A$ and $\varepsilon$. Thus by (7.6.9) and
(7.6.10) we have \be |Rm|\leq K'_1(r\sqrt{t_0})^{-2},\
 \ \mbox{on}\ \ B_{t_0}(x_0,2Ar\sqrt{t_0}), 
\ee for some positive constant $K'_1=K'_1(w,r,A,\varepsilon)$
depending only on $w$, $r$, $A$ and $\varepsilon$. This gives us
the curvature estimate on $B_{t_0}(x_0,2Ar\sqrt{t_0})$ for all
$t_0 \geq T(w,r,A,\varepsilon)$.

{}From the arguments in proving the above assertion (i), we have
the estimates (7.6.6) and (7.6.7) and the solution is well-defined
on the whole parabolic neighborhood
$P(x_0,t_0,\frac{r_0\sqrt{t_0}}{4},-\tau(r_0\sqrt{t_0})^2)$ for
all $t_0 \geq T(w,r,A,\varepsilon)$. Clearly we may assume that
$(K'_1)^{-\frac{1}{2}}r < \min \{\frac{r_0}{4}, \sqrt{\tau}r_0
\}$. Thus by combining with the curvature estimate (7.6.11), we
can apply Theorem 7.5.1(i) to get the following volume control \be
\Vol_{t_0}(B_{t_0}(x,(K'_1)^{-\frac{1}{2}}r\sqrt{t_0}))
\geq \kappa((K'_1)^{-\frac{1}{2}}r\sqrt{t_0})^3 
\ee for any $x\in B_{t_0}(x_0,Ar\sqrt{t_0})$, where
$\kappa=\kappa(w,r,A,\varepsilon)$ is some positive constant
depending only on $w$, $r$, $A$ and $\varepsilon$. So by using the
assertion (i), we see that for $t_0\geq T$ with
$T=T(w,r,\xi,A,\varepsilon)$ large enough, the curvature estimate
(7.6.5) holds for all points in $B_{t_0}(x_0,Ar\sqrt{t_0})$.

\medskip
(iii) We next want to extend the curvature estimate (7.6.5) to all
points in the forward parabolic neighborhood
$P(x_0,t_0,Ar\sqrt{t_0},Ar^2t_0)$. Consider the time interval
$[t_0,t_0+Ar^2t_0]$ in the parabolic neighborhood. In assertion
(ii), we have obtained the desired estimate (7.6.5) at $t=t_0$.
Suppose estimate (7.6.5) holds on a maximal time interval
$[t_0,t')$ with $t'\leq t_0+Ar^2t_0$. This says that we have \be
|2tR_{ij}+g_{ij}|<\xi 
\ee on $P(x_0,t_0,Ar\sqrt{t_0},t'-t_0)\triangleq\{(x,t)\ |\ x\in
B_t(x_0,Ar\sqrt{t_0}),t\in[t_0,t')\}$ so that either there exists
a surgery in the ball $B_{t'}(x_0,Ar\sqrt{t_0})$ at $t=t'$, or
there holds $|2tR_{ij}+g_{ij}|=\xi$ somewhere in
$B_{t'}(x_0,Ar\sqrt{t_0})$ at $t=t'$. Since the Ricci curvature is
near $-\frac{1}{2t'}$ in the geodesic ball, the surgeries cannot
occur there. Thus we only need to consider the latter possibility.

Recall that the evolution of the length of a curve $\gamma$ and
the volume of a domain $\Omega$ are given by
\begin{align*}
\frac{d}{dt}L_t(\gamma)
&=-\int_{\gamma}\Ric(\dot{\gamma},\dot{\gamma})ds_t \\
\mbox{and } \quad \frac{d}{dt}Vol_t(\Omega) &=-\int_{\Omega}RdV_t.
\end{align*}
By substituting the curvature estimate (7.6.13) into the above two
evolution equations and using the volume lower bound (7.6.6), it
is not hard to see \be \Vol_{t'}(B_{t'}(x_0,\sqrt{t'}))\geq
 \kappa'(t')^{\frac{3}{2}} 
\ee for some positive constant
$\kappa'=\kappa'(w,r,\xi,A,\varepsilon)$ depending only on $w$,
$r$, $\xi$, $A$ and $\varepsilon$. Then by the above assertion
(ii), the combination of the curvature estimate (7.6.13) and the
volume lower bound (7.6.14) implies that the curvature estimate
(7.6.5) still holds for all points in $B_{t'}(x_0,Ar\sqrt{t_0})$
provided $T=T(w,r,\xi,A,\varepsilon)$ is chosen large enough. This
is a contradiction.  Therefore we have proved assertion (iii).
\end{pf}

We now state and prove the following important \textbf{Thick-thin
decomposition theorem}\index{thick-thin decomposition theorem}. A
more general version (without the restriction on $\varepsilon$)
was implicitly claimed by Perelman in \cite{P1} and \cite{P2}.

\begin{theorem}[The Thick-thin decomposition theorem]
For any $w>0$ and $0 < \varepsilon \leq \frac{1}{2}w$, there
exists a positive constant $\rho = \rho(w,\varepsilon) \leq 1$
with the following property. Suppose $g_{ij}(t)$ $(t \in
[0,+\infty))$ is a solution, constructed by Theorem $7.4.3$ with
the nonincreasing $($continuous$)$ positive functions
$\widetilde{\delta}(t)$ and $\widetilde{r}(t)$, to the Ricci flow
with surgery and with a compact orientable normalized
three-manifold as initial data, where each $\delta$-cutoff at a
time $t$ has $\delta=\delta(t) \leq \min
\{\widetilde{\delta}(t),\widetilde{r}(2t)\}$. Then for any
arbitrarily fixed $\xi >0$, for $t$ large enough, the manifold
$M_t$ at time $t$ admits a decomposition $M_t=M_{\rm
thin}(w,t)\cup M_{\rm thick}(w,t)$ with the following properties:
\begin{itemize}
\item[(a)] For every $x\in M_{\rm thin}(w,t)$, there exists some
$r=r(x,t)>0$, with $0<r\sqrt{t}<\rho\sqrt{t}$, such that
$$
Rm\geq -(r\sqrt{t})^{-2} \ \ on \ \ B_t(x,r\sqrt{t}), \ \
\mbox{and }
$$
$$
\Vol_t(B_t(x,r\sqrt{t}))<w(r\sqrt{t})^3.
$$
\item[(b)] For every $x\in M_{\rm thick}(w,t)$, we have
$$
|2tR_{ij}+g_{ij}|<\xi\ \ on \ \ B_t(x,\rho\sqrt{t}), \ \ \mbox{and
}
$$
$$
\Vol_t(B_t(x,\rho\sqrt{t}))\geq \frac{1}{10}w(\rho\sqrt{t})^3.
$$
\end{itemize}
Moreover, if we take any sequence of points $x^{\alpha}\in M_{\rm
thick}(w,t^{\alpha}),t^{\alpha}\rightarrow+\infty$, then the
scalings of $g_{ij}(t^{\alpha})$ around $x^{\alpha}$ with factor
$(t^{\alpha})^{-1}$ converge smoothly, along a subsequence of
$\alpha\rightarrow+\infty$, to a complete hyperbolic manifold of
finite volume with constant sectional curvature $-\frac{1}{4}$.
\end{theorem}

\begin{pf}
Let $\bar{r}=\bar{r}(w,\varepsilon),\theta=\theta(w,\varepsilon)$
and $h$ be the positive constants obtained in Corollary 7.5.3. We
may assume $\rho\leq\bar{r} \leq e^{-3}$. For any point $x\in
M_t$, there are two cases: either
$$
{\rm (i)}\ \ \ \mbox{ } \min\{Rm\ |\ B_t(x,\rho\sqrt{t})\}\geq
-(\rho\sqrt{t})^{-2},
$$
or
$$
{\rm (ii)}\ \ \ \mbox{    } \min\{Rm\ |\
B_t(x,\rho\sqrt{t})\}<-(\rho\sqrt{t})^{-2}.$$

Let us first consider Case (i). If
$\Vol_t(B_t(x,\rho\sqrt{t}))<\frac{1}{10}w(\rho\sqrt{t})^3$, then
we can choose $r$ slightly less than $\rho$ so that
$$
Rm\geq -(\rho\sqrt{t})^{-2}\geq -(r\sqrt{t})^{-2}
$$
on $B_t(x,r\sqrt{t})(\subset B_t(x,\rho\sqrt{t}))$, and
$$
\Vol_t(B_t(x,r\sqrt{t}))<\frac{1}{10}w(\rho\sqrt{t})^3<w(r\sqrt{t})^3;
$$
thus $x\in M_{thin}(w,t)$. If $\Vol_t(B_t(x,\rho\sqrt{t}))\geq
\frac{1}{10}w(\rho\sqrt{t})^3$, we can apply Lemma 7.6.2(ii) to
conclude that for $t$ large enough,
$$
|2tR_{ij}+g_{ij}|<\xi \quad \mbox{ on }\; B_t(x,\rho\sqrt{t});
$$
thus $x\in M_{thick}(w,t)$.

Next we consider Case (ii). By continuity, there exists
$0<r=r(x,t)<\rho$ such that \be
\min\{Rm\ |\ B_t(x,r\sqrt{t})\}=-(r\sqrt{t})^{-2}. 
\ee If $\theta^{-1}h\leq r\sqrt{t}\mbox{  }(\leq
\bar{r}\sqrt{t})$, we can apply Corollary 7.5.3 to conclude
$$
\Vol_t(B_t(x,r\sqrt{t}))<w(r\sqrt{t})^3;
$$
thus $x\in M_{thin}(w,t)$.

We now consider the difficult subcase $\mbox{}r\sqrt{t}
<\theta^{-1}h$. By the pinching assumption, we have
\begin{align*}
R &  \geq (r\sqrt{t})^{-2}(\log[(r\sqrt{t})^{-2}(1+t)] -3)\\
&  \geq (\log r^{-2}-3)(r\sqrt{t})^{-2}\\
&  \geq 2(r\sqrt{t})^{-2}\\
&  \geq 2\theta^{2}h^{-2}
\end{align*}
somewhere in $B_t(x,r\sqrt{t})$. Since $h$ is the maximal cutoff
radius for surgeries in $[\frac{t}{2},t]$, by the design of the
$\delta$-cutoff surgery, we have
\begin{align*}
 h &  \leq \sup \left\{{\delta}^2(s)\widetilde{r}(s)\ |\ s\in
\left[\frac{t}{2},t\right]\right\}\\
& \leq\widetilde{\delta}\(\frac{t}{2}\)
\widetilde{r}(t)\widetilde{r}\(\frac{t}{2}\).
\end{align*}
Note also $\widetilde{\delta}(\frac{t}{2})\rightarrow0$ as
$t\rightarrow+\infty$. Thus from the canonical neighborhood
assumption, we see that for $t$ large enough, there exists a point
in the ball $B_t(x,r\sqrt{t})$ which has a canonical neighborhood.

We claim that for $t$ sufficiently large, the point $x$ satisfies
$$
R(x,t)\geq \frac{1}{2}(r\sqrt{t})^{-2},
$$
and then the above argument shows that the point $x$ also has a
canonical $\varepsilon$-neck or $\varepsilon$-cap neighborhood.
Otherwise, by continuity, we can choose a point $x^*\in
B_t(x,r\sqrt{t})$ with $R(x^*,t)=\frac{1}{2}(r\sqrt{t})^{-2}$.
Clearly the new point $x^*$ has a canonical neighborhood $B^*$ by
the above argument. In particular, there holds
$$
C^{-1}_2(\varepsilon)R \leq  \frac{1}{2}(r\sqrt{t})^{-2} \leq
C_2(\varepsilon)R
$$
on the canonical neighborhood $B^*$. By the definition of
canonical neighborhood assumption, we have
$$
B_t(x^*,\sigma^*)\subset B^* \subset B_t(x^*,2\sigma^*)
$$
for some $\sigma^* \in
(0,C_1(\varepsilon)R^{-\frac{1}{2}}(x^*,t))$. Clearly, without
loss of generality, we may assume (in the definition of canonical
neighborhood assumption) that $\sigma^* >
2R^{-\frac{1}{2}}(x^*,t).$ Then
\begin{align*}
R(1+t) &  \geq\frac{1}{2}C^{-1}_2(\varepsilon)r^{-2}\\
&  \geq \frac{1}{2}C^{-1}_2(\varepsilon)\rho^{-2}
\end{align*}
on $B_t(x^*,2r\sqrt{t}).$ Thus when we choose
$\rho=\rho(w,\varepsilon)$ small enough, it follows from the
pinching assumption that
\begin{align*}
Rm&  \geq-[f^{-1}(R(1+t))/(R(1+t))]R\\
&  \geq -\frac{1}{2}(r\sqrt{t})^{-2},
\end{align*}
on $B_t(x^*,2r\sqrt{t}).$ This is a contradiction with (7.6.15).

We have seen that $tR(x,t)\geq \frac{1}{2}r^{-2} (\geq
\frac{1}{2}\rho^{-2})$. Since $r^{-2}
> \theta^2 h^{-2} t$ in this subcase, we conclude that for
arbitrarily given $A<+\infty$, \be
tR(x,t)>A^2\rho^{-2} 
\ee as long as $t$ is large enough.

Let $B$, with $B_t(x,\sigma)\subset B \subset B_t(x,2\sigma)$, be
the canonical $\varepsilon$-neck or $\varepsilon$-cap neighborhood
of $(x,t)$. By the definition of the canonical neighborhood
assumption, we have
$$
0 < \sigma < C_1(\varepsilon)R^{-\frac{1}{2}}(x,t),
$$
$$
C^{-1}_2(\varepsilon)R \leq  R(x,t) \leq C_2(\varepsilon)R, \
\mbox{ on } \ B,
$$
and \be \Vol_t(B) \leq \varepsilon \sigma^3 \leq
\frac{1}{2}w \sigma^3. 
\ee Choose $0 < A < C_1(\varepsilon)$ so that $\sigma =
AR^{-\frac{1}{2}}(x,t)$. For sufficiently large $t$, since
\begin{align*}
R(1+t) &  \geq C^{-1}_2(\varepsilon)(tR(x,t))\\
&  \geq \frac{1}{2}C^{-1}_2(\varepsilon)\rho^{-2},
\end{align*}
on $B$, we can require $\rho=\rho(w,\varepsilon)$ to be smaller
still, and use the pinching assumption to conclude
\begin{align}
Rm&  \geq-[f^{-1}(R(1+t))/(R(1+t))]R\\
&  \geq -(AR^{-\frac{1}{2}}(x,t))^{-2} \nn\\
&  = -\sigma^{-2}, \nn
\end{align} 
on $B$. For sufficiently large $t$, we adjust
\begin{align}
r&  = \sigma (\sqrt{t})^{-1}\\
&  = (AR^{-\frac{1}{2}}(x,t))(\sqrt{t})^{-1} \nn\\
&  < \rho, \nn
\end{align} 
by (7.6.16). Then the combination of (7.6.17), (7.6.18) and
(7.6.19) implies that $ x\in M_{thin}(w,t).$

The last statement in (b) follows directly from Lemma 7.6.2. (Here
we also used Bishop-Gromov volume comparison, Theorem 7.5.2 and
Hamilton's compactness theorem to take a subsequent limit.)

Therefore we have completed the proof of the theorem.
\end{pf}

To state the long-time behavior of a solution to the Ricci flow
with surgery, we first recall some basic terminology in
three-dimensional topology. A three-manifold $M$ is called
\textbf{irreducible}\index{irreducible} if every smooth two-sphere
embedded in $M$ bounds a three-ball in $M$. If we have a solution
$(M_t,g_{ij}(t))$ obtained by Theorem 7.4.3 with a compact,
orientable and irreducible three-manifold $(M,g_{ij})$ as initial
data, then at each time $t>0$, by the cutoff surgery procedure,
the solution manifold $M_t$ consists of a finite number of
components where one of the components, called the
\textbf{essential component} \index{essential component} and
denoted by $M_t^{(1)}$, is diffeomorphic to the initial manifold
$M$ while the rest are diffeomorphic to the three-sphere
$\mathbb{S}^3$.

The main result of this section is the following generalization of
Theorem 5.3.4. A more general version of the result (without the
restriction on $\varepsilon$) was implicitly claimed by Perelman
in \cite{P2}.

\begin{theorem}[Long-time behavior of the Ricci flow with surgery]
Let $w\!>0$ and $0<\!\varepsilon\leq\!\frac{1}{2}w$ be any small
positive constants and let $(M_t,g_{ij}(t)),$ $0<t<+\infty,$ be a
solution to the Ricci flow with surgery, constructed by Theorem
$7.4.3$ with the nonincreasing $($continuous$)$ positive functions
$\widetilde{\delta}(t)$ and $\widetilde{r}(t)$ and with a compact,
orientable, irreducible and normalized three-manifold $M$ as
initial data, where each $\delta$-cutoff at a time $t$ has
$\delta=\delta(t) \leq \min
\{\widetilde{\delta}(t),\widetilde{r}(2t)\}$. Then one of the
following holds: either
\begin{itemize}
\item[(i)] for all sufficiently large $t$, we have $M_t = M_{\rm
thin}(w,t)$; or \item[(ii)] there exists a sequence of times
$t^{\alpha}\rightarrow+\infty$ such that the scalings of
$g_{ij}(t^{\alpha})$ on the essential component
$M_{t^\alpha}^{(1)}$, with factor $(t^{\alpha})^{-1}$, converge in
the $C^{\infty}$ topology to a hyperbolic metric on the initial
compact manifold $M$ with constant sectional curvature
$-\frac{1}{4}$; or \item[(iii)] we can find a finite collection of
complete noncompact hyperbolic three-manifolds
$\mathcal{H}_1,\ldots,\mathcal{H}_m$, with finite volume, and
compact subsets $K_1,\ldots,K_m$ of
$\mathcal{H}_1,\ldots,\mathcal{H}_m$ respectively obtained by
truncating each cusp of the hyperbolic manifolds along constant
mean curvature torus of small area, and for all $t$ beyond some
time $T<+\infty$ we can find diffeomorphisms $\varphi_l,1\leq
l\leq m$, of $K_l$ into $M_t$ so that as long as $t$ is
sufficiently large, the metric $t^{-1}\varphi^*_l(t)g_{ij}(t)$ is
as close to the hyperbolic metric as we like on the compact sets
$K_1,\ldots,K_m$; moreover, the complement $M_t\backslash
(\varphi_1(K_1)\cup\cdots\cup\varphi_m(K_m))$ is contained in the
thin part $M_{\rm thin}(w,t)$, and the boundary tori of each $K_l$
are incompressible in the sense that each $\varphi_l$ injects
$\pi_1(\partial K_l)$ into $\pi_1(M_t)$.
\end{itemize}
\end{theorem}

\begin{pf}
The proof of the theorem follows, with some modifications, the
same argument of Hamilton \cite{Ha99} as in the proof of Theorem
5.3.4.

Clearly we may assume that the thick part $M_{\rm thick}(w,t)$ is
not empty for a sequence $t^{\alpha}\rightarrow+\infty$, since
otherwise we have case (i). If we take a sequence of points
$x^{\alpha}\in M_{\rm thick}(w,t^{\alpha})$, then by Theorem
7.6.3(b) the scalings of $g_{ij}(t^{\alpha})$ around $x^{\alpha}$
with factor $(t^{\alpha})^{-1}$ converge smoothly, along a
subsequence of $\alpha\rightarrow+\infty$, to a complete
hyperbolic manifold of finite volume with constant sectional
curvature $-\frac{1}{4}$. The limits may be different for
different choices of $(x^{\alpha},t^{\alpha})$. If a limit is
compact, we have case (ii). Thus we assume that all limits are
noncompact.

Consider all the possible hyperbolic limits of the solution, and
among them choose one such complete noncompact hyperbolic
three-manifold $\mathcal{H}$ with the least possible number of
cusps. Denote by $h_{ij}$ the hyperbolic metric of $\mathcal{H}$.
For all small $a>0$ we can truncate each cusp of $\mathcal{H}$
along a constant mean curvature torus of area $a$ which is
uniquely determined; we denote the remainder by $\mathcal{H}_a$.
Fix $a>0$ so small that Lemma 5.3.7 is applicable for the compact
set $\mathcal{K}=\mathcal{H}_a$. Pick an integer $l_0$
sufficiently large and an $\epsilon_0$ sufficiently small to
guarantee from Lemma 5.3.8 that the identity map $Id$ is the only
harmonic map $F$ from $\mathcal{H}_a$ to itself with taking
$\partial \mathcal{H}_a$ to itself, with the normal derivative of
$F$ at the boundary of the domain normal to the boundary of the
target, and with $d_{C^{l_0}(\mathcal{H}_a)}(F,Id) < \epsilon_0$.
Then choose a positive integer $q_0$ and a small number
$\delta_0>0$ from Lemma 5.3.7 such that if $\widetilde{F}$ is a
diffeomorphism of $\mathcal{H}_a$ into another complete noncompact
hyperbolic three-manifold
$(\widetilde{\mathcal{H}},\widetilde{h}_{ij})$ with no fewer cusps
(than $\mathcal{H}$), of finite volume and satisfying
$$
\|\widetilde{F}^*\widetilde{h}_{ij}-h_{ij}\|_{C^{q_0}
(\mathcal{H}_a)}\leq\delta_0,
$$
then there exists an isometry $I$ of $\mathcal{H}$ to
$\widetilde{\mathcal{H}}$ such that
$$
d_{C^{l_0}(\mathcal{H}_a)}(\widetilde{F},I)<\epsilon_0.
$$
By Lemma 5.3.8 we further require $q_0$ and $\delta_0$ to
guarantee the existence of a harmonic diffeomorphism from
$(\mathcal{H}_a,\widetilde{g}_{ij})$ to $(\mathcal{H}_a,h_{ij})$
for any metric $\widetilde{g}_{ij}$ on $\mathcal{H}_a$ with
$\|\widetilde{g}_{ij}-h_{ij}\|_{C^{q_0}(\mathcal{H}_a)}\leq
\delta_0$.

Let $x^{\alpha}\in M_{\rm
thick}(w,t^{\alpha}),t^{\alpha}\rightarrow+\infty$, be a sequence
of points such that the scalings of $g_{ij}(t^{\alpha})$ around
$x^{\alpha}$ with factor $(t^{\alpha})^{-1}$ converge to $h_{ij}$.
Then there exist a marked point $x^{\infty}\in \mathcal{H}_a$ and
a sequence of diffeomorphisms $F_{\alpha}$ from $\mathcal{H}_a$
into $M_{t^{\alpha}}$ such that
$F_{\alpha}(x^{\infty})=x^{\alpha}$ and
$$
\|(t^{\alpha})^{-1}F^*_{\alpha}g_{ij}(t^{\alpha})
-h_{ij}\|_{C^m(\mathcal{H}_a)}\rightarrow 0
$$
as $\alpha\rightarrow\infty$ for all positive integers $m$. By
applying Lemma 5.3.8 and the implicit function theorem, we can
change $F_{\alpha}$ by an amount which goes to zero as
$\alpha\rightarrow\infty$ so as to make $F_{\alpha}$ a harmonic
diffeomorphism taking $\partial \mathcal{H}_a$ to a constant mean
curvature hypersurface $F_{\alpha}(\partial \mathcal{H}_a)$ of
$(M_{t^{\alpha}},(t^{\alpha})^{-1}g_{ij}(t^{\alpha}))$ with the
area $a$ and satisfying the free boundary condition that the
normal derivative of $F_{\alpha}$ at the boundary of the domain is
normal to the boundary of the target; and by combining with Lemma
7.6.2 (iii), we can smoothly continue each harmonic diffeomorphism
$F_{\alpha}$ forward in time a little to a family of harmonic
diffeomorphisms $F_{\alpha}(t)$ from $\mathcal{H}_a$ into $M_t$
with the metric $t^{-1}g_{ij}(t)$, with
$F_{\alpha}(t^{\alpha})=F_{\alpha}$ and with the time $t$ slightly
larger than $t^{\alpha}$, where $F_{\alpha}(t)$ takes $\partial
\mathcal{H}_a$ into a constant mean curvature hypersurface of
$(M_t,t^{-1}g_{ij}(t))$ with the area $a$ and also satisfies the
free boundary condition. Moreover, since the surgeries do not take
place at the points where the scalar curvature is negative, by the
same argument as in Theorem 5.3.4 for an arbitrarily given
positive integer $q\geq q_0$, positive number $\delta<\delta_0$,
and sufficiently large $\alpha$, we can ensure the extension
$F_{\alpha}(t)$ satisfies
$\|t^{-1}F^*_{\alpha}(t)g_{ij}(t)-h_{ij}\|_{C^q(\mathcal{H}_a)}\leq
\delta$ on a maximal time interval $t^{\alpha}\leq t\leq
\omega^{\alpha}$ (or $t^{\alpha}\leq t< \omega^{\alpha}$ when
$\omega^{\alpha}=+\infty$), and with
$\|(\omega^{\alpha})^{-1}F^*_{\alpha}
(\omega^{\alpha})g_{ij}(\omega^{\alpha})-h_{ij}\|_{C^q(\mathcal{H}_a)}=
\delta$, when $\omega^{\alpha}<+\infty$. Here we have implicitly
used the fact that $F_{\alpha}(\omega^{\alpha})(\partial
\mathcal{H}_a)$ is still strictly concave to ensure the map
$F_{\alpha}(\omega^{\alpha})$ is diffeomorphic.

We further claim that there must be some $\alpha$ such that
$\omega^{\alpha}=+\infty$; in other words, at least one hyperbolic
piece persists. Indeed, suppose that for each large enough
$\alpha$ we can only continue the family $F_{\alpha}(t)$ on a
finite interval $t^{\alpha}\leq t\leq \omega^{\alpha}<+\infty$
with
$$
\|(\omega^{\alpha})^{-1}F^*_{\alpha}(\omega^{\alpha})g_{ij}
(\omega^{\alpha})-h_{ij}\|_{C^q(\mathcal{H}_a)}=\delta.
$$
Consider the new sequence of manifolds
($M_{\omega^{\alpha}},g_{ij}(\omega^{\alpha})$). Clearly by Lemma
7.6.1, the scalings of $g_{ij}(\omega^{\alpha})$ around the new
origins $F_{\alpha}(\omega^{\alpha})(x^{\infty})$ with factor
$(\omega^{\alpha})^{-1}$ converge smoothly (by passing to a
subsequence) to a complete noncompact hyperbolic three-manifold
$\widetilde{\mathcal{H}}$ with the metric $\widetilde{h}_{ij}$ and
the origin $\widetilde{x}^{\infty}$ and with finite volume. By the
choice of the old limit $\mathcal{H}$, the new limit
$\widetilde{\mathcal{H}}$ has at least as many cusps as
$\mathcal{H}$. By the definition of convergence, we can find a
sequence of compact subsets $\widetilde{U}_{\alpha}$ exhausting
$\widetilde{\mathcal{H}}$ and containing $\widetilde{x}^{\infty}$,
and a sequence of diffeomorphisms $\widetilde{F}_{\alpha}$ of
neighborhood of $\widetilde{U}_{\alpha}$ into
$M_{\omega^{\alpha}}$ with
$\widetilde{F}_{\alpha}(\widetilde{x}^{\infty})
=F_{\alpha}(\omega^{\alpha})(x^{\infty})$ such that for each
compact subset $\widetilde{U}$ of $\widetilde{\mathcal{H}}$ and
each integer $m$,
$$
\|(\omega^{\alpha})^{-1}\widetilde{F}^*_{\alpha}(g_{ij}
(\omega^{\alpha}))-\widetilde{h}_{ij}
\|_{C^m(\widetilde{U})}\rightarrow 0
$$
as $\alpha\rightarrow+\infty$. Thus for sufficiently large
$\alpha$, we have the map
$$
G_{\alpha}=\widetilde{F}^{-1}_{\alpha}\circ
F_{\alpha}(\omega^{\alpha}):\ \ \mathcal{H}_a\rightarrow
\widetilde{\mathcal{H}}
$$
such that
$$
\|G^*_{\alpha}\widetilde{h}_{ij}-h_{ij}
\|_{C^q(\mathcal{H}_a)}<\widetilde{\delta}
$$
for any fixed $\widetilde{\delta}>\delta$. Then a subsequence of
$G_{\alpha}$ converges at least in the $C^{q-1}(\mathcal{H}_a)$
topology to a map $G_{\infty}$ of $\mathcal{H}_a$ into
$\widetilde{\mathcal{H}}$ which is a harmonic map from
$\mathcal{H}_a$ into $\widetilde{\mathcal{H}}$ and takes $\partial
\mathcal{H}_a$ to a constant mean curvature hypersurface
$G_{\infty}(\partial \mathcal{H}_a)$ of
$(\widetilde{\mathcal{H}},\widetilde{h}_{ij})$ with the area $a$,
as well as satisfies the free boundary condition. Clearly,
$G_{\infty}$ is at least a local diffeomorphism. Since
$G_{\infty}$ is the limit of diffeomorphisms, the only possibility
of overlap is at the boundary. Note that $G_{\infty}(\partial
\mathcal{H}_a)$ is still strictly concave. So $G_{\infty}$ is
still a diffeomorphism. Moreover by using the standard regularity
result of elliptic partial differential equations (see for example
\cite{GT}), we also have \be
\|G^*_{\infty}\widetilde{h}_{ij}-h_{ij}
\|_{C^q(\mathcal{H}_a)}=\delta. 
\ee Now by Lemma 5.3.7 we deduce that there exists an isometry $I$
of $\mathcal{H}$ to $\widetilde{\mathcal{H}}$ with
$$
d_{C^{l_0}(\mathcal{H}_a)}(G_{\infty},I)<\epsilon_0.
$$
Thus $I^{-1}\circ G_{\infty} $ is a harmonic diffeomorphism of
$\mathcal{H}_a$ to itself which satisfies the free boundary
condition and
$$
d_{C^{l_0}(\mathcal{H}_a)}(I^{-1}\circ G_{\infty},Id)<\epsilon_0.
$$
However the uniqueness in Lemma 5.3.8 concludes that $I^{-1}\circ
G_{\infty}=Id$ which contradicts (7.6.20). So we have shown that
at least one hyperbolic piece persists and the metric
$t^{-1}F^*_{\alpha}(t)g_{ij}(t)$, for $\omega^{\alpha}\leq
t<\infty$, is as close to the hyperbolic metric $h_{ij}$ as we
like.

We can continue to form other persistent hyperbolic pieces in the
same way as long as there is a sequence of points $y^{\beta}\in
M_{\rm thick}(w,t^{\beta}),t^{\beta}\rightarrow+\infty$, lying
outside the chosen pieces. Note that
$V(t)(t+\frac{3}{2})^{-\frac{3}{2}}$ is nonincreasing on
$[0,+\infty)$. Therefore by combining with Margulis lemma (see for
example \cite{Gro79} or \cite{KM}), we have proved that there
exists a finite collection of complete noncompact hyperbolic
three-manifolds $\mathcal{H}_1,\ldots,\mathcal{H}_m$ with finite
volume, a small number $a>0$ and a time $T<+\infty$ such that for
all $t$ beyond $T$ we can find diffeomorphisms $\varphi_l(t)$ of
$(\mathcal{H}_l)_a$ into $M_t$, $1\leq l\leq m$, so that as long
as $t$ is sufficiently large, the metric
$t^{-1}\varphi^*_l(t)g_{ij}(t)$ is as close to the hyperbolic
metrics as we like and the complement $M_t\backslash
(\varphi_1(t)((\mathcal{H}_1)_a)\cup\cdots\cup
\varphi_m(t)((\mathcal{H}_m)_a))$ is contained in the thin part
$M_{\rm thin}(w,t)$.

It remains to show the boundary tori of any persistent hyperbolic
piece are incompressible. Let $B$ be a small positive number and
assume the above positive number $a$ is much smaller than $B$. Let
$M_a(t)=\varphi_l(t)((\mathcal{H}_l)_a)$ ($1\leq l\leq m$) be such
a persistent hyperbolic piece of the manifold $M_t$ truncated by
boundary tori of area $at$ with constant mean curvature, and
denote by $M^c_a(t)=M_t\setminus \stackrel{\circ}{M}_a(t)$ the
part of $M_t$ exterior to $M_a(t)$. Thus there exists a family of
subsets $M_B(t)\subset M_a(t)$ which is a persistent hyperbolic
piece of the manifold $M_t$ truncated by boundary tori of area
$Bt$ with constant mean curvature. We also denote by
$M^c_B(t)=M_t\setminus \stackrel{\circ}{M}_B(t)$. By Van Kampen's
Theorem, if $\pi_1(\partial M_B(t))$ injects into
$\pi_1(M_B^c(t))$ then it injects into $\pi_1(M_t)$ also. Thus we
only need to show $\pi_1(\partial M_B(t))$ injects into
$\pi_1(M_B^c(t))$.

As before we will use a contradiction argument to show
$\pi_1(\partial M_B(t))$ injects into $\pi_1(M^c_B(t))$. Let $T$
be a torus in $\partial M_B(t)$ and suppose $\pi_1(T)$ does not
inject into $\pi_1(M^c_B(t))$. By Dehn's Lemma we know that the
kernel is a cyclic subgroup of $\pi_1(T)$ generated by a primitive
element. Consider the normalized metric
$\widetilde{g}_{ij}(t)=t^{-1}g_{ij}(t)$ on $M_t$. Then by the work
of Meeks-Yau \cite{MY} or Meeks-Simon-Yau \cite{MSY}, we know that
among all disks in $M^c_B(t)$ whose boundary curve lies in $T$ and
generates the kernel of $\pi_1(T)$, there is a smooth embedded
disk normal to the boundary which has the least possible area
(with respect to the normalized metric $\widetilde{g}_{ij}(t)$).
Denote by $D$ the minimal disk and
$\widetilde{A}=\widetilde{A}(t)$ its area. We will show that
$\widetilde{A}(t)$ decreases at a certain rate which will arrive
at a contradiction.

We first consider the case that there exist no surgeries at the
time $t$. Exactly as in Part III of the proof of Theorem 5.3.4,
the change of the area $\widetilde{A}(t)$ comes from the change in
the metric and the change in the boundary. For the change in the
metric, we choose an orthonormal frame $X,Y,Z$ at a point $x$ in
the disk $D$ so that $X$ and $Y$ are tangent to the disk $D$ while
$Z$ is normal. Since the normalized metric $\widetilde{g}_{ij}$
evolves by
$$
\frac{\partial}{\partial t}\widetilde{g}_{ij}
=-t^{-1}(\widetilde{g}_{ij}+2\widetilde{R}_{ij}),
$$
the (normalized) area element $d\widetilde{\sigma}$ of the disk
$D$ around $x$ satisfies
$$\frac{\partial}{\partial t}d\widetilde{\sigma}
=-t^{-1}(1+\widetilde{\Ric}(X,X)+\widetilde{\Ric}(Y,Y))
d\widetilde{\sigma}.
$$
For the change in the boundary, we notice that the tensor
$\widetilde{g}_{ij}+2\widetilde{R}_{ij}$ is very small for the
persistent hyperbolic piece. Then by using the Gauss-Bonnet
theorem as before, we obtain the rate of change of the area \be
\frac{d\widetilde{A}}{dt}\leq
-\int_D\(\frac{1}{t}+\frac{\widetilde{R}}{2t}\)
d\widetilde{\sigma}+\frac{1}{t}\int_{\partial
D}\widetilde{k}d\widetilde{s}-\frac{2\pi}{t}+
o\(\frac{1}{t}\)\widetilde{L}, 
\ee where $\widetilde{k}$ is the geodesic curvature of the
boundary and $\widetilde{L}$ is the length of the boundary curve
$\partial D$ (with respect to the normalized metric
$\widetilde{g}_{ij}(t)$).  Since $\widetilde{R}\geq
-{3t}/{2(t+\frac{3}{2})}$ for all $t\geq 0$ by (7.6.2), the first
term on the RHS of (7.6.21) is bounded above by
$$
-\int_D\(\frac{1}{t}+\frac{\widetilde{R}}{2t}\)d\widetilde{\sigma}\leq
-\frac{1}{t}\(\frac{1}{4}- o(1)\)\widetilde{A};
$$
while the second term on the RHS of (7.6.21) can be estimated
exactly as before by
$$
\frac{1}{t}\int_{\partial D}\widetilde{k}d\widetilde{s}\leq
\frac{1}{t}\(\frac{1}{4}+o(1)\)\widetilde{L}.
$$
Thus we obtain \be \frac{d\widetilde{A}}{dt}\leq
\frac{1}{t}\left[\(\frac{1}{4}+o(1)\)
\widetilde{L}-\(\frac{1}{4}-o(1)\)\widetilde{A}-2\pi\right].
\ee

Next we show that these arguments also work for the case that
there exist surgeries at the time $t$. To this end, we only need
to check that the embedded minimal disk $D$ lies in the region
which is unaffected by surgery. Our surgeries for the irreducible
three-manifold took place on $\delta$-necks in
$\varepsilon$-horns, where the scalar curvatures are at least
$\delta^{-2} (\widetilde{r}(t))^{-1}$, and the components with
nonnegative scalar curvature have been removed. So the hyperbolic
piece is not affected by the surgeries. In particular, the
boundary $\partial D$ is unaffected by the surgeries. Thus if
surgeries occur on the minimal disk, the minimal disk has to pass
through a long thin neck before it reaches the surgery regions.
Look at the intersections of the embedding minimal disk with a
generic center two-sphere $\mathbb{S}^2$ of the long thin neck;
these are circles. Since the two-sphere $\mathbb{S}^2$ is simply
connected, we can replace the components of the minimal disk $D$
outside the center two-sphere $\mathbb{S}^2$ by some corresponding
components on the center two-sphere $\mathbb{S}^2$ to form a new
disk which also has $\partial D$ as its boundary. Since the metric
on the long thin neck is nearly a product metric, we could choose
the generic center two-sphere $\mathbb{S}^2$ properly so that the
area of the new disk is strictly less than the area of the
original disk $D$. This contradiction proves the minimal disk lies
entirely in the region unaffected by surgery.

Since $a$ is much smaller than $B$, the region within a long
distance from $\partial M_B(t)$ into $M^c_B(t)$ will look nearly
like a hyperbolic cusplike collar and is unaffected by the
surgeries. So we can repeat the arguments in the last part of the
proof of Theorem 5.3.4 to bound the length $\widetilde{L}$ by the
area $\widetilde{A}$ and to conclude
$$\frac{d\widetilde{A}}{dt}\leq -\frac{\pi}{t}$$ for all
sufficiently large times $t$, which is impossible because the RHS
is not integrable. This proves that the boundary tori of any
persistent hyperbolic piece are incompressible.

Therefore we have proved the theorem.
\end{pf}

\section{Geometrization of Three-manifolds}

In the late 70's and early 80's, Thurston \cite{Th82, Th86, Th88}
proved a number of remarkable results on the existence of
geometric structures on a class of three-manifolds:
\textbf{Haken}\index{Haken} manifolds (i.e. each of them contains
an incompressible surface of genus $\geq$ 1). These results
motivated him to formulate a profound conjecture which roughly
says every compact three-manifold admits a canonical decomposition
into domains, each of which has a canonical geometric structure.
To give a detailed description of the conjecture, we recall some
terminology as follows.

An $n$-dimensional complete Riemannian manifold $(M,g)$ is called
a \textbf{homogeneous manifold}\index{homogeneous manifold} if its
group of isometries acts transitively on the manifold. This means
that the homogeneous manifold looks the same metrically at
everypoint. For example, the round $n$-sphere $\mathbb{S}^n$, the
Euclidean space $\mathbb{R}^n$ and the standard hyperbolic space
$\mathbb{H}^n$ are homogeneous manifolds.  A Riemannian manifold
is said to \textbf{be modeled}\index{be modeled} on a given
homogeneous manifold $(M,g)$ if every point of the manifold has a
neighborhood isometric to an open set of $(M,g)$. And an
$n$-dimensional Riemannian manifold is called a \textbf{locally
homogeneous manifold}\index{homogeneous manifold!locally} if it is
complete and is modeled on a homogeneous manifold. By a theorem of
Singer \cite{Sin}, the universal cover of a locally homogeneous
manifold (with the pull-back metric) is a homogeneous manifold.

In dimension three, every locally homogeneous manifold with finite
volume is modeled on one of the following eight homogeneous
manifolds (see for example Theorem 3.8.4 of \cite{Th97}):
\begin{itemize}
\item[(1)] $\mathbb{S}^3$, the round three-sphere; \item[(2)]
$\mathbb{R}^3$, the Euclidean space ; \item[(3)] $\mathbb{H}^3$,
the standard hyperbolic space; \item[(4)] $\mathbb{S}^2 \times
\mathbb{R}$; \item[(5)] $\mathbb{H}^2 \times \mathbb{R}$;
\item[(6)] $Nil$, the three-dimensional nilpotent Heisenberg group
(consisting of upper triangular $3\times 3$ matrices with diagonal
entries 1); \item[(7)] $\widetilde{PSL}(2,\mathbb{R})$, the
universal cover of the unit sphere bundle of $\mathbb{H}^2$;
\item[(8)] $Sol$, the three-dimensional solvable Lie group.
\end{itemize}

\medskip
A three-manifold $M$ is called \textbf{prime}\index{prime} if it
is not diffeomorphic to $\mathbb{S}^3$ and if every (topological)
$\mathbb{S}^2 \subset M$, which separates $M$ into two pieces, has
the property that one of the two pieces is diffeomorphic to a
three-ball. Recall that a three-manifold is irreducible if every
embedded two-sphere bounds a three-ball in the manifold. Clearly
an irreducible three-manifold is either prime or is diffeomorphic
to $\mathbb{S}^3$. Conversely, an orientable prime three-manifold
is either irreducible or is diffeomorphic to $\mathbb{S}^2 \times
\mathbb{S}^1$ (see for example \cite{Hempel}). One of the first
results in three-manifold topology is the following prime
decomposition theorem obtained by Kneser \cite{Kn} in 1929 (see
also Theorem 3.15 of \cite{Hempel}).

\medskip
{\bf Prime Decomposition Theorem. }\index{prime! decomposition
theorem} \emph{ Every compact orientable three-manifold admits a
decomposition as a finite connected sum of orientable prime
three-manifolds. }

\medskip
In \cite{Mil62}, Milnor showed that the factors involved in the
above Prime Decomposition are unique.  Based on the prime
decomposition, the question about topology of compact orientable
three-manifolds is reduced to the question about prime
three-manifolds. Thurston's Geometrization Conjecture is about
prime three-manifolds.

\medskip
{\bf Thurston's Geometrization Conjecture. }\index{Thurston's
geometrization conjecture} Let $M$ be a compact, orientable and
prime three-manifold. Then there is an embedding of a finite
number of disjoint unions, possibly empty, of incompressible
two-tori $\coprod_{i} T_i^2 \subset M$ such that every component
of the complement admits a locally homogeneous Riemannian metric
of finite volume.

\medskip
We remark that the existence of a torus decomposition, also called
JSJ-decomposition, was already obtained by Jaco-Shalen \cite{Jaco}
and Johannsen \cite{Joha}. The JSJ-decomposition states that any
compact, orientable, and prime three-manifold has a finite
collection, possibly empty, of disjoint incompressible embedding
two-tori $\{T_i^2\}$ which separate the manifold into a finite
collection of compact three-manifolds (with toral boundary), each
of which is either a graph manifold or is \textbf{atoroidal}
\index{atoroidal} in the sense that any subgroup of its
fundamental group isomorphic to $\mathbb{Z}\times\mathbb{Z}$ is
conjugate into the fundamental group of some component of its
boundary. A compact three-manifold $X$, possibly with boundary, is
called a \textbf{graph manifold}\index{graph manifold} if there is
a finite collection of disjoint embedded tori $T_i \subset X$ such
that each component of $X \setminus \bigcup T_i$ is an
$\mathbb{S}^1$ bundle over a surface. Thus the point of the
conjecture is that the components should all be geometric.

The geometrization conjecture for a general compact orientable
3-man\-ifold is the statement that each of its prime factors
satisfies the above conjecture. We say a compact orientable
three-manifold is \textbf{geometrizable}\index{geometrizable} if
it satisfies the geometric conjecture.

We also remark, as is well-known, that the Poincar\'{e} conjecture
can be deduced from Thur\-ston's geometrization conjecture.
Indeed, suppose that we have a compact simply connected
three-manifold that satisfies the conclusion of the geometrization
conjecture. If it were not diffeomorphic to the three-sphere
$\mathbb{S}^3$, there would be a prime factor in the prime
decomposition of the manifold. Since the prime factor still has
vanishing fundamental group, the (torus) decomposition of the
prime factor in the geometrization conjecture must be trivial.
Thus the prime factor is a compact homogeneous manifold model.
>From the list of above eight models, we see that the only compact
three-dimensional model is $\mathbb{S}^3$. This is a
contradiction. Consequently, the compact simply connected
three-manifold is diffeomorphic to $\mathbb{S}^3$.

Now, based on the long-time behavior result (Theorem 7.6.4) we
apply the Ricci flow to discuss Thurston's geometrization
conjecture. Let $M$ be a compact, orientable and prime
three-manifold. Since a prime orientable three-manifold is either
irreducible or is diffeomorphic to $\mathbb{S}^2 \times
\mathbb{S}^1$, we may thus assume the manifold $M$ is irreducible
also. Arbitrarily given a (normalized) Riemannian metric for the
manifold $M$, we use it as initial data to evolve the metric by
the Ricci flow with surgery. From Theorem 7.4.3, we know that the
Ricci flow with surgery has a long-time solution on a maximal time
interval $[0,T)$ which satisfies the a priori assumptions and has
a finite number of surgeries on each finite time interval.
Furthermore, from the long-time behavior theorem (Theorem 7.6.4),
we have well-understood geometric structures on the thick part.
Whereas, to understand the thin part, Perelman announced the
following result in \cite{P2}.

\medskip
{\bf Perelman's Claim}\index{Perelman's claim} (cf. Theorem 7.4 of
\cite{P2}){\bf .} \ Suppose $(M^{\alpha},g^{\alpha}_{ij})$ is a
sequence of compact orientable three-manifolds, closed or with
convex boundary, and $w^{\alpha}\rightarrow  0$. Assume that
\begin{itemize}
\item[(1)] for each point $x \in M^{\alpha}$ there exists a radius
$\rho = \rho^{\alpha}(x), 0<\rho<1,$ not exceeding the diameter of
the manifold, such that the ball $B(x,\rho)$ in the metric
$g^{\alpha}_{ij}$ has volume at most $w^{\alpha}\rho^3$ and
sectional curvatures at least $-\rho^{-2}$; \item[(2)] each
component of the boundary of $M^{\alpha}$ has diameter at most
$w^{\alpha}$, and has a (topologically trivial) collar of length
one, where the sectional curvatures are between $-1/4 - \epsilon$
and $-1/4 + \epsilon$.
\end{itemize}

Then $M^{\alpha}$ for sufficiently large $\alpha$ are
diffeomorphic to graph manifolds.

\medskip
The topology of graph manifolds is well understood; in particular,
every graph manifold is geometrizable (see \cite{W}).

The proof of Perelman's Claim promised in \cite{P2} is still not
available in literature. Nevertheless, recently in \cite{ShY},
Shioya and Yamaguchi provided a proof of Perelman's Claim for the
special case when all the manifolds $(M^{\alpha},g^{\alpha}_{ij})$
are closed. That is, they proved the following weaker assertion.

\medskip
{\bf Weaker Assertion} (Theorem 8.1 of Shioya-Yamaguchi
\cite{ShY}){\bf .} \index{Weaker Assertion} \ Suppose
$(M^{\alpha},g^{\alpha}_{ij})$ is a sequence of compact orientable
three-manifolds without boundary, and $w^{\alpha}\rightarrow 0$.
Assume that for each point $x \in M^{\alpha}$ there exists a
radius $\rho = \rho^{\alpha}(x),$ not exceeding the diameter of
the manifold, such that the ball $B(x,\rho)$ in the metric
$g^{\alpha}_{ij}$ has volume at most $w^{\alpha}\rho^3$ and
sectional curvatures at least $-\rho^{-2}$. Then $M^{\alpha}$ for
sufficiently large $\alpha$ are diffeomorphic to graph manifolds.

\medskip
Based on the the long-time behavior theorem (Theorem 7.6.4) and
using the above Weaker Assertion, we can now give a proof for
Thurston's geometrization conjecture \index{Thurston's
geometrization conjecture}. We remark that if we assume the above
Perelman's Claim, then we does not need to use Thurston's theorem
for Haken manifolds in the proof of Theorem 7.7.1.

\begin{theorem}
Thurston's geometrization conjecture is true.
\end{theorem}

\begin{pf}
Let $M$ be a compact, orientable, and prime three-manifold
(without boundary). Without loss of generality, we may assume that
the manifold $M$ is irreducible also.

Recall that the theorem of Thurston (see for example Theorem A and
Theorem B in the third section of \cite{Morgan}, see also
\cite{Mc} and \cite{Ot})) says that any compact, orientable, and
irreducible Haken three-manifold (with or without boundary) is
geometrizable. Thus in the following, we may assume that the
compact three-manifold $M$ (without boundary) is atoroidal, and
then the fundamental group $\pi_1(M)$ contains no noncyclic,
abelian subgroup.

Arbitrarily given a (normalized) Riemannian metric on the manifold
$M$, we use it as initial data for the Ricci flow. Arbitrarily
take a sequence of small positive constants $w^{\alpha}\rightarrow
0$ as $\alpha \rightarrow +\infty$. For each fixed $\alpha$, we
set $\varepsilon={w^{\alpha}}/{2}>0$. Then by Theorem 7.4.3, the
Ricci flow with surgery has a long-time solution
$(M^{\alpha}_t,g^{\alpha}_{ij}(t))$ on a maximal time interval
$[0,T^{\alpha})$, which satisfies the a priori assumptions (with
the accuracy parameter $\varepsilon={w^{\alpha}}/{2}$) and has a
finite number of surgeries on each finite time interval. Since the
initial manifold is irreducible, by the surgery procedure, we know
that for each $\alpha$ and each $t>0$ the solution manifold
$M^{\alpha}_t$ consists of a finite number of components where the
essential component $(M^{\alpha}_t)^{(1)}$ is diffeomorphic to the
initial manifold $M$ and the others are diffeomorphic to the
three-sphere $\mathbb{S}^3$.

If for some $\alpha=\alpha_0$ the maximal time $T^{\alpha_0}$ is
finite, then the solution $(M^{\alpha_0}_t,g^{\alpha_0}_{ij}(t))$
becomes extinct at $T^{\alpha_0}$ and the (irreducible) initial
manifold $M$ is diffeomorphic to $\mathbb{S}^3/\Gamma$ (the metric
quotients of round three-sphere); in particular, the manifold $M$
is geometrizable. Thus we may assume that the maximal time
$T^{\alpha} = +\infty$ for all $\alpha$.

We now apply the long-time behavior theorem (Theorem 7.6.4). If
there is some $\alpha$ such that case (ii) of Theorem 7.6.4
occurs, then for some sufficiently large time $t$, the essential
component $(M^{\alpha}_t)^{(1)}$ of the solution manifold
$M^{\alpha}_t$ is diffeomorphic to a compact hyperbolic space, so
the initial manifold $M$ is geometrizable. Whereas if there is
some sufficiently large $\alpha$ such that case (iii) of Theorem
7.6.4  occurs, then it follows that for all sufficiently large
$t$, there is an embedding of a (nonempty) finite number of
disjoint unions of incompressible two-tori $\coprod_{i} T_i^2 $ in
the essential component $(M^{\alpha}_t)^{(1)}$ of $M^{\alpha}_t$.
This is a contradiction since we have assumed the initial manifold
$M$ is atoroidal.

It remains to deal with the situation that there is a sequence of
positive $\alpha_k \rightarrow +\infty$ such that the solutions
$(M^{\alpha_k}_t,g^{\alpha_k}_{ij}(t))$ always satisfy case (i) of
Theorem 7.6.4. That is, for each $\alpha_k$, $M^{\alpha_k}_t =
M_{\rm thin}(w^{\alpha_k},t)$ when the time $t$ is sufficiently
large. By the Thick-thin decomposition theorem (Theorem 7.6.3),
there is a positive constant, $0<\rho(w^{\alpha_k}) \leq 1$, such
that as long as $t$ is sufficiently large, for every $x\in
M^{\alpha_k}_t = M_{\rm thin}(w^{\alpha_k},t)$, we have some
$r=r(x,t)$, with $0<r\sqrt{t}<\rho(w^{\alpha_k})\sqrt{t}$, such
that \be Rm\geq -(r\sqrt{t})^{-2} \quad \mbox{ on } \;
B_t(x,r\sqrt{t}),
\ee and \be
\Vol_t(B_t(x,r\sqrt{t}))<w^{\alpha_k}(r\sqrt{t})^3. 
\ee

Clearly we only need to consider the essential component
$(M^{\alpha_k}_{t})^{(1)}$.  We divide the discussion into the
following two cases:

\smallskip
(1) there is a positive constant $1<C<+\infty$ such that for each
$\alpha_k$ there is a sufficiently large time $t_k>0$ such that
\be r(x,t_k)\sqrt{t_k} < C \cdot {\rm
diam}\,\((M^{\alpha_k}_{t_k})^{(1)}\)
\ee for all $x\in (M^{\alpha_k}_{t_k})^{(1)} \subset M_{\rm
thin}(w^{\alpha_k},{t_k})$;

\smallskip
(2) there are a subsequence $\alpha_k$ (still denoted by
$\alpha_k$), and sequences of positive constants $C_k \rightarrow
+\infty$ and times $T_k < +\infty$ such that for each $t \geq
T_k$, we have \be r(x(t),t)\sqrt{t} \geq C_k \cdot {\rm diam}\,
\((M^{\alpha_k}_{t})^{(1)}\)
\ee for some $x(t) \in (M^{\alpha_k}_{t})^{(1)}$, $k=1,2,\ldots.$
Here we denote by ${\rm diam}\,((M^{\alpha}_t)^{(1)})$ the
diameter of the essential component $(M^{\alpha}_t)^{(1)}$ with
the metric $g^{\alpha}_{ij}(t)$ at the time $t$.

Let us first consider case (1). For each point $x \in
(M^{\alpha_k}_{t_k})^{(1)} \subset M_{\rm thin}(w^{\alpha_k},$
${t_k})$, we denote by $\rho_k(x) = C^{-1}r(x,t_k)\sqrt{t_k}$.
Then by (7.7.1), (7.7.2) and (7.7.3), we have
$$
\rho_k(x) <{\rm diam}\,\((M^{\alpha_k}_{t_k})^{(1)}\),
$$
$$
\Vol_{t_k}(B_{t_k}(x,\rho_k(x))) \leq
\Vol_{t_k}(B_{t_k}(x,r(x,t_k)\sqrt{t_k}))<C^3w^{\alpha_k}(\rho_k(x))^3,
$$
and
$$
Rm\geq -(r(x,t_k)\sqrt{t_k})^{-2} \geq -(\rho_k(x))^{-2}
$$
on  $B_{t_k}(x,\rho_k(x)).$ Then it follows from the above Weaker
Assertion that $(M^{\alpha_k}_{t_k})^{(1)}$, for sufficiently
large $k$, are diffeomorphic to graph manifolds. This implies that
the (irreducible) initial manifold $M$ is diffeomorphic to a graph
manifold. So the manifold $M$ is geometrizable in case (1).

We next consider case (2). Clearly, for each $\alpha_k$ and the
chosen $T_k$, we may assume that the estimates (7.7.1) and (7.7.2)
hold for all $t \geq T_k$ and $x\in (M^{\alpha_k}_t)^{(1)}$. The
combination of (7.7.1) and (7.7.4) gives \be Rm \geq
-C^{-2}_k({\rm diam}\,((M^{\alpha_k}_{t})^{(1)}))^{-2}
\quad \mbox{ on } \; (M^{\alpha_k}_{t})^{(1)}, 
\ee for all $t \geq T_k$. If there are a subsequence $\alpha_k$
(still denoted by $\alpha_k$) and a sequence of times $t_k \in
(T_k,+\infty)$ such that \be
\Vol_{t_k}((M^{\alpha_k}_{t_k})^{(1)}) <
w'_k({\rm diam}\,((M^{\alpha_k}_{t_k})^{(1)}))^3 
\ee for some sequence $w'_k \rightarrow 0$, then it follows from
the Weaker Assertion that $(M^{\alpha_k}_{t_k})^{(1)}$, for
sufficiently large $k$, are diffeomorphic to graph manifolds which
implies the initial manifold $M$ is geometrizable. Thus we may
assume that there is a positive constant $w'$ such that \be
\Vol_t((M^{\alpha_k}_t)^{(1)})
\geq w'({\rm diam}\,((M^{\alpha_k}_t)^{(1)}))^3 
\ee for each $k$ and all $t \geq T_k$.

In view of the estimates (7.7.5) and (7.7.7), we now want to use
Theorem 7.5.2 to get a uniform upper bound  for the curvatures of
the essential components
$((M^{\alpha_k}_t)^{(1)},g^{\alpha_k}_{ij}(t))$ with sufficiently
large time $t$. Note that the estimate in Theorem 7.5.2 depends on
the parameter $\varepsilon$ and our $\varepsilon$'s depend on
$w^{\alpha_k}$ with $0<\varepsilon = {w^{\alpha_k}}/{2}$; so it
does not work in the present situation. Fortunately we notice that
the curvature estimate for smooth solutions in Corollary 7.2.3 is
independent of $\varepsilon$. In the following we try to use
Corollary 7.2.3 to obtain the desired curvature estimate.

We first claim that for each $k$, there is a sufficiently large
$T'_k \in (T_k,+\infty)$ such that the solution, when restricted
to the essential component $((M^{\alpha_k}_t)^{(1)},$
$g^{\alpha_k}_{ij}(t))$, has no surgery for all $t\geq T'_k$.
Indeed, for each fixed $k$, if there is a $\delta(t)$-cutoff
surgery at a sufficiently large time $t$, then the manifold
$((M^{\alpha_k}_t)^{(1)},g^{\alpha_k}_{ij}(t))$ would contain a
$\delta(t)$-neck $B_t(y,{\delta(t)}^{-1}R(y,t)^{-\frac{1}{2}})$
for some $y \in (M^{\alpha_k}_t)^{(1)}$ with the volume ratio \be
\frac{\Vol_t(B_t(y,{\delta(t)}^{-1}
R(y,t)^{-\frac{1}{2}}))}{({\delta(t)}^{-1}R(y,t)^{-\frac{1}{2}})^3}
\leq 8\pi{\delta(t)}^2. 
\ee On the other hand, by (7.7.5) and (7.7.7), the standard
Bishop-Gromov volume comparison implies that
$$
\frac{\Vol_t(B_t(y,{\delta(t)}^{-1}
R(y,t)^{-\frac{1}{2}}))}{({\delta(t)}^{-1}R(y,t)^{-\frac{1}{2}})^3}
\geq c(w')
$$
for some positive constant $c(w')$ depending only on $w'$. Since
$\delta(t)$ is very small when $t$ is large, this arrives at a
contradiction with (7.7.8). So for each $k$, the essential
component $((M^{\alpha_k}_t)^{(1)},g^{\alpha_k}_{ij}(t))$ has no
surgery for all sufficiently large $t$.

For each $k$, we consider any fixed time $\tilde{t}_k > 3 T'_k$.
Let us scale the solution $g^{\alpha_k}_{ij}(t)$ on the essential
component $(M^{\alpha_k}_t)^{(1)}$ by
$$
\tilde{g}^{\alpha_k}_{ij}(\cdot,s) =
(\tilde{t}_{k})^{-1}g^{\alpha_k}_{ij}(\cdot,\tilde{t}_{k}s).
$$
Note that $(M^{\alpha_k}_t)^{(1)}$ is diffeomorphic to $M$ for all
$t$. By the above claim, we see that the rescaled solution
$(M,\tilde{g}^{\alpha_k}_{ij}(\cdot,s))$ is a smooth solution to
the Ricci flow on the time interval $s \in [\frac{1}{2},1]$. Set
$$
\tilde{r}_k = \(\sqrt{\tilde{t}_k}\)^{-1} {\rm
diam}\,\((M^{\alpha_k}_{\tilde{t}_k})^{(1)}\).
$$
Then by (7.7.4), (7.7.5) and (7.7.7), we have
$$
\tilde{r}_k \leq C^{-1}_k \rightarrow 0, \quad \mbox{ as }\; k
\rightarrow +\infty,
$$
$$
\widetilde{Rm} \geq -C^{-2}_k(\tilde{r}_k)^{-2}, \quad \mbox{on }
\; B_1(x(\tilde{t}_k),\tilde{r}_k),
$$
and
$$
\Vol_1(B_1(x(\tilde{t}_k),\tilde{r}_k)) \geq w'(\tilde{r}_k)^3,
$$
where $\widetilde{Rm}$ is the rescaled curvature, $x(\tilde{t}_k)$
is the point given by (7.7.4) and
$B_1(x(\tilde{t}_k),\tilde{r}_k)$ is the geodesic ball of rescaled
solution at the time $s=1$. Moreover, the closure of
$B_1(x(\tilde{t}_k),\tilde{r}_k)$ is the whole manifold
$(M,\tilde{g}^{\alpha_k}_{ij}(\cdot,1)).$

Note that in Theorem 7.2.1, Theorem 7.2.2 and Corollary 7.2.3, the
condition about normalized initial metrics is just to ensure that
the solutions satisfy the Hamilton-Ivey pinching estimate. Since
our solutions $(M^{\alpha_k}_t,g^{\alpha_k}_{ij}(t))$ have already
satisfied the pinching assumption, we can then apply Corollary
7.2.3 to conclude
$$
|\widetilde{Rm}(x,s)| \leq K(w')(\tilde{r}_k)^{-2},
$$
whenever $s \in [1- \tau(w')(\tilde{r}_k)^{2}, 1]$, $x \in
(M,\tilde{g}^{\alpha_k}_{ij}(\cdot,s))$ and $k$ is sufficiently
large. Here $K(w')$ and $\tau(w')$ are positive constants
depending only on $w'$. Equivalently, we have the curvature
estimates \be |Rm(\cdot,t)| \leq K(w') ({\rm
diam}\,((M^{\alpha_k}_{\tilde{t}_k})^{(1)}))^{-2},
\quad \mbox{ on }\;  M, 
\ee whenever $t \in [\tilde{t}_k - \tau(w')({\rm
diam}\,((M^{\alpha_k}_{\tilde{t}_k})^{(1)}))^{2},\tilde{t}_k]$ and
$k$ is sufficiently large.

\smallskip\noindent
For each $k$, let us scale
$((M^{\alpha_k}_t)^{(1)},g^{\alpha_k}_{ij}(t))$ with the factor
$({\rm diam}((M^{\alpha_k}_{\tilde{t}_k})^{(1)}))^{-2}$ and shift
the time $\tilde{t}_k$ to the new time zero. By the curvature
estimate (7.7.9) and Hamilton's compactness theorem (Theorem
4.1.5), we can take a subsequential limit (in the $C^{\infty}$
topology) and get a smooth solution to the Ricci flow on $M \times
(-\tau(w'),0]$.  Moreover, by (7.7.5), the limit has nonnegative
sectional curvature on $M \times (-\tau(w'),0]$. Recall that we
have removed all compact components with nonnegative scalar
curvature. By combining this with the strong maximum principle, we
conclude that the limit is a flat metric. Hence in case (2), $M$
is diffeomorphic to a flat manifold and then it is also
geometrizable.

Therefore we have completed the proof of the theorem.
\end{pf}


\printindex
\end{document}